# THE LC METHOD

## A PARALLELIZABLE NUMERICAL METHOD FOR APPROXIMATING THE ROOTS OF SINGLE-VARIABLE POLYNOMIALS

*by*

*Daniel Alba Cuellar*

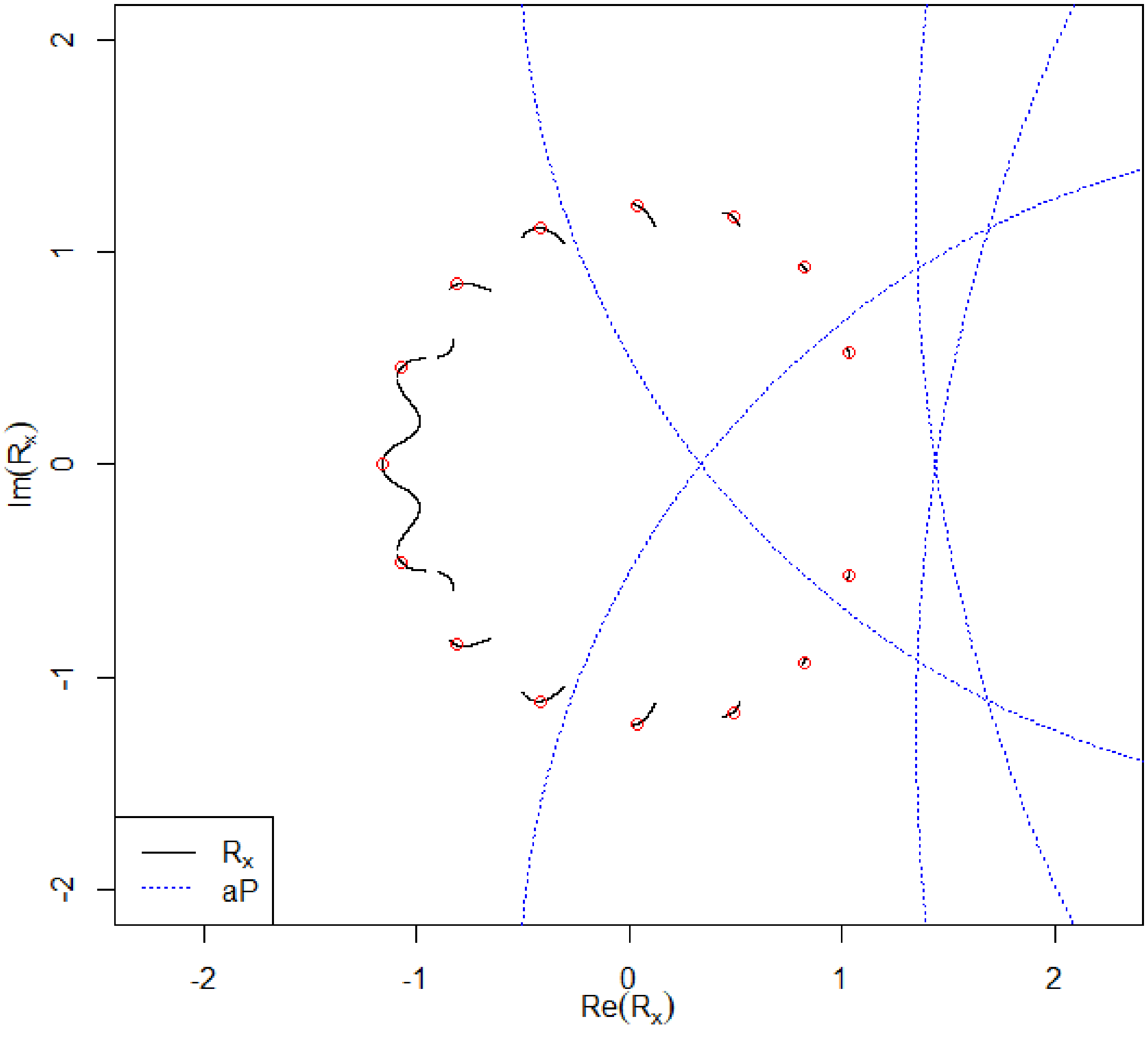

# Table of Contents

# Preface

In this work, our problem of interest consists in finding numerical approximations to a set of $n$ values $z = R_i$, $i = 1,2, \dots, n$, called roots, that satisfy the polynomial equation of the form

$$z^n + C_1 z^{n-1} + C_2 z^{n-2} + \cdots + C_{n-2} z^2 + C_{n-1} z + C_n = 0. \quad (1)$$

This is one of the oldest mathematical problems[1] (in ancient Babylon circa 2,000 BC, quadratic equations were used to solve practical geometric problems[2]), which has been, throughout the centuries, an important source for new ideas, new concepts, and even new branches in mathematics; this has had a profound impact on science (physics, chemistry, biology, astronomy, computer science, medicine, etc.), technology, engineering and even the arts. There are literally thousands of published works that address various aspects of this problem, as well as hundreds of well-established numerical methods that seek to directly solve equation (1); in references [1] and [2] we can find a historical and exhaustive review of the bibliography and relevant methods related to this topic. Even so, we believe that it is worth to continue searching for new ways to find solutions to this old and fundamental problem, as this could still give rise to new ideas and discoveries useful in the development of knowledge, especially in the fields of computer science and mathematics, which are, now more than ever before, very closely related.

The objective of this work is to describe a numerical method that helps to find the roots of equation (1), **focusing on obtaining reasonable initial approximations** that can later serve as inputs to fast local converging iterative methods, such as Newton-Raphson. This method, developed by the author of this work, and called *the LC method* (L is for line and C is for circumference), seeks to take advantage of the parallel processing capabilities from current computing equipment, to simultaneously find initial approximations to all roots of equation (1).

On how the idea to develop the LC method came about to be, we could say that with the help of a programmable graphing calculator (an old Texas Instruments TI-85), the author of this work realized that, by means of Möbius mappings from $\mathbb{C}$ to $\mathbb{C}$ ($\mathbb{C}$ being the set of complex numbers), it is possible to transform lines $\ell$ into circumferences $1/\ell$. The same author then observed that, for a quadratic equation $z^2 + C_1 z + C_2 = 0$, it is possible to construct a line $\ell_1\colon -C_1/2 + te^{i\theta^*}$ ($t, \theta^* \in \mathbb{R}$, $t$ variable) and a circumference $zC\colon C_2/\ell_1$ (by means of a Möbius transformation), whose intersections are the roots $R_1$, $R_2$ of this quadratic equation; the inclination angle $\theta^*$ for $\ell_1$ that causes both $\ell_1$ and $zC$ to contain $R_1$, $R_2$, is found by bisecting angle $\angle 0 P_1 P_2$, where $P_1 = -C_1/2$ is a fixed point in $\ell_1$, and $P_2 = C_2/P_1$ is a fixed point in $zC$. This is like drawing geometric constructions with ruler and compass, but with the help of a computer. The personal discovery of this result served as a motivation to devise a way to extend these geometric constructions of lines and circumferences towards univariate polynomials of higher degree, aiming to approximate all roots of equation (1), for $n > 2$.

The description of the LC method in this work focuses exclusively on proof of concept; it is not concerned on the design of advanced computational implementations, or the analysis of computational complexity, or the comparison against other numerical methods; as for the latter, we will only say that the vast majority of numerical methods that seek to find the roots of equation (1) are iterative in nature; i.e., they seek to improve the quality of a root's estimate from a previous estimate, repeating the refining process several times until a stopping criterion is met; on the other hand, **the LC method seeks to obtain, simultaneously, initial approximations to all roots of equation (1) by means of a two-step process; each step in this process is itself a parallel process**. Thus, rather than competing against other numerical methods that try to find the roots of equation (1), the LC method seeks to be a complementary preprocessing module, providing good quality initial estimates to other well-established numerical methods, so that these can compute, in fewer iterations, precise approximations to the roots of equation (1).

In the first step of the LC method, we construct on the plane of complex numbers $\mathbb{C}$, for an arbitrary value $\theta$ in the interval $[-\pi, \pi)$, the following geometric elements:

- a line $\ell_1(C_1, \theta) = -C_1/2 + te^{i\theta}$ ($t \in \mathbb{R}$), which seeks to contain one of the roots of equation (1), let's say $R_1$;
- a circumference $zC(C_1, C_n, \theta) = (-1)^n C_n/\ell_1$, which seeks to contain the product of $n-1$ roots $R_2 \cdot R_3 \cdot \ldots \cdot R_n$;
- a *terminal semi-line* $tL(C_1, C_2, \ldots, C_{n-1}, \theta)$, which is derived from $\ell_1$; $tL$ depends on coefficients $C_1, C_2, C_3, \ldots, C_{n-1}$ from equation (1), and, like $zC$, it also seeks to contain the product of $n-1$ roots $R_2 \cdot R_3 \cdot \ldots \cdot R_n$.

The construction of $\ell_1(C_1, \theta)$, $tL(C_1, C_2, \ldots, C_{n-1}, \theta)$ and $zC(C_1, C_n, \theta)$ is based on Vieta's relations, which link the roots $R_1, R_2, \ldots, R_n$ of equation (1) with its coefficients $C_1, C_2, \ldots, C_n$. Possible intersections between $tL$ and $zC$ provide information on the degree of angular proximity of line $\ell_1$ with respect to $R_1$, without knowing a priori where $R_1$ is located. This degree of angular proximity is called *weighted error* $e(\theta)$; we define $e(\theta)$ in such a way that it takes values in the interval $(-1,1]$; the closer $e(\theta)$ is to 0, the closer $\ell_1(C_1, \theta)$ is to $R_1$. In theory, it is perfectly possible to process simultaneously, on a computer with $N$ processors, $N$ triplets of geometric elements $\ell_1(C_1, \theta_k)$, $tL(C_1, C_2, \ldots, C_{n-1}, \theta_k)$ and $zC(C_1, C_n, \theta_k)$, with $\theta_k = -\pi + 2\pi k/N$, $k = 0,1,2, \ldots, N-1$, in order to construct a discrete *proximity map* $\{e(\theta_k)\}$.

In the second step of the LC method, we approximate the *theta roots* $\theta_i^*$ $i = 1,2, \ldots, n$, for the discrete proximity map $\{e(\theta_k)\}$; i.e., we compute crossings of $\{e(\theta_k)\}$ with the horizontal axis $y = 0$ via linear interpolation, which can be done efficiently on computers with parallel processing capabilities. The estimated theta roots $\hat{\theta}_i^*$ allow us to obtain initial approximations $\hat{R}_i$ to the roots $R_i$ of equation (1), **directly** from geometric structures $\ell_1(C_1, \hat{\theta}_i^*)$, $tL(C_1, C_2, \ldots, C_{n-1}, \hat{\theta}_i^*)$, $zC(C_1, C_n, \hat{\theta}_i^*)$.

To illustrate the operation of the LC method on a typical personal computer, the author of this work has developed a series of programs written in the R language that use vectorization of operations, which can be seen as a rudimentary form of parallel processing. These programs avoid the use of control structures that explicitly carry out iterations of blocks of instructions, such as

`for`, `while`, and `repeat`. These R programs (or scripts) are listed, as supplementary material, in the annexes at the end of this work, so that the reader can have access to basic programming tools that allow him/her not only to reproduce the results of the numerical examples included in this work, but to experiment with his/her own new ideas and seek answers to questions that may arise as the reading progresses. It is worth mentioning that to properly understand the contents of this work, it is enough to have basic working knowledge of algebra, computer programming, differential calculus, and numerical methods. In the end, what is sought is to describe and verify the operation of the LC method, using simple and direct language so that the contents of this book are accessible to a wide audience.

This book is divided into five chapters and contains supplementary material to support the expositions in the main chapters. In chapter 1, some preliminary concepts are presented, such as operations and properties of complex numbers, vectors in the field of complex numbers $\mathbb{C}$, scalars, parametric expressions of lines and circumferences in $\mathbb{C}$, and Möbius transformations, which will serve as the basis for the implementation of the LC method. Chapter 2 describes the LC method for univariate polynomials of degree 2; here the initial idea that motivated this whole study is explained and developed in detail. Chapter 3 extends the LC method to univariate polynomials of degree 3; it is here that the concept of discrete proximity map is described in detail. Chapter 4 extends the LC method to univariate polynomials of degree 4; here it is necessary to introduce two additional concepts that give definitive support to the LC method: the *terminal curve* $t\mathfrak{C}$ (concept analogous to that of the terminal semi-line $tL$) and the *dynamic squared distance* between terminal curve $t\mathfrak{C}$ and circumference $zC$; this second concept extends the notion of intersections between $tL$ and $zC$. Finally, chapter 5 generalizes the concepts described in chapter 4, without introducing new additional concepts, so it is possible to extend the LC method to polynomials of a single variable of degree $n = 5$ and above. At the end of chapter 5, we include some comments that serve as conclusions, and indicate possible lines of future work.

Each chapter, starting from chapter 2, begins with a theoretical-conceptual section that describes the basic theoretical ideas and concepts that support the LC method, and ends with a series of numerical examples that illustrate in a practical way the ideas developed in the theoretical part. During the development of the numerical examples, we sometimes encounter some additional theoretical issues not covered in the opening sections of the chapter; we will try to address these issues in detail during the development of the numerical examples. The supplemental annexes in this book list the R routines and programs used to produce the results in the numerical examples, and also include descriptions and comments that help to clarify how do these routines and programs work.

# PART I: Description of the LC Method

# Chapter One: Fundamentals of the LC Method

## Work Plan for this Chapter

First, we will describe some preliminary concepts: coefficients and roots of univariate polynomials, the fundamental theorem of algebra, and Vieta's relations between roots of a univariate polynomial and its coefficients. Next, we will talk about fundamental concepts used by the method of lines and circles (LC method) used to find initial approximations for the complex roots of univariate polynomials: arithmetic of complex numbers (addition, multiplication, multiplicative inverse, norm, and argument of a complex number), and lines and circles in the plane of complex numbers (definitions, Möbius transformation, properties).

## Univariate Polynomials

A univariate polynomial of degree $n$ is an algebraic expression of the form

$$p(x) \coloneqq x^n + c_1 x^{n-1} + c_2 x^{n-2} + \cdots + c_{n-2} x^2 + c_{n-1} x + c_n. \tag{1.1}$$

$x$ is the variable quantity in polynomial $p(x)$, while the fixed values $c_1$, $c_2$, ..., $c_n$ are the coefficients of $p(x)$; without loss of generality, in this work we will always assume that the coefficient of $x^n$ in expression (1.1) is 1. Our goal is to numerically approximate all the values $x$ that satisfy the equation $p(x) = 0$. These values are called *the roots of polynomial* $p(x)$. In this work we will assume that the variable $x$ and the coefficients $c_1$, $c_2$, ..., $c_n$ are complex numbers.

## Fundamental Theorem of Algebra and Factor Theorem

The Fundamental Theorem of Algebra states that the polynomial $p(x)$ of degree $n$ in expression (1.1) has exactly $n$ roots (some of these roots may be repeated).

The Factor Theorem allows us to express polynomial (1.1) as a product of factors, each of which is directly associated with one of its roots; it states that a polynomial $p(x)$ has a factor $(x - r)$ if and only if $p(r) = 0$; that is, $p(x)$ has a factor $(x - r)$ if and only if $r$ is a root of $p(x)$.

Denote the roots of $p(x)$ in expression (1.1) as $r_1$, $r_2$, ..., $r_n$ ($\in \mathbb{C}$). The Fundamental Theorem of Algebra, together with the Factor Theorem, allow us to rewrite polynomial $p(x)$ from expression (1.1) as

$$p(x) \coloneqq (x - r_1) \cdot (x - r_2) \cdot \cdots \cdot (x - r_n). \tag{1.2}$$

## Vieta's Relations

These relationships combine expressions (1.1) and (1.2); they allow us to express coefficients $c_1$, $c_2$, ..., $c_n$ from $p(x)$ in terms of its roots $r_1$, $r_2$, ..., $r_n$. Vieta's relations tell us that

$$-c_1 = r_1 + r_2 + \cdots + r_n,$$

$$\begin{aligned} c_2 = r_1r_2 + r_1r_3 + \cdots + r_1r_n \\ +r_2r_3 + \cdots + r_2r_n \\ \vdots \\ +r_{n-1}r_n, \end{aligned}$$

$$\begin{aligned} -c_3 = r_1r_2r_3 + r_1r_2r_4 + \cdots + r_1r_2r_n + r_1r_3r_4 + \cdots + r_1r_{n-1}r_n \\ +r_2r_3r_4 + \cdots + r_2r_{n-1}r_n \\ \vdots \\ +r_{n-2}r_{n-1}r_n, \end{aligned}$$

$$\vdots$$

$$(-1)^n c_n = r_1 \cdot r_2 \cdot r_3 \cdot \cdots \cdot r_n.$$

In words, Vieta's relations tell us that:

$c_1$ is the sum of all roots of $p(x)$, multiplied by $-1$.

$c_2$ is the sum of all possible "double products" of the roots of $p(x)$. A double product is a term of the form $r_ir_j$, with $i < j$; $i \in \{1,2,\dots,n-1\}$, $j \in \{2,3,\dots,n\}$.

$c_3$ is the sum of all possible "triple products" of the roots of $p(x)$, multiplied by $-1$. A triple product is a term of the form $r_ir_jr_k$, with $i < j < k$; $i \in \{1,2,\dots,n-2\}$, $j \in \{2,3,\dots,n-1\}$, $k \in \{3,4,\dots,n\}$.

$\vdots$

$c_n$ is the product of all roots of $p(x)$, multiplied by $(-1)^n$.

In general, $c_k$, with $k \in \{2,3,\dots,n\}$, is the sum of all possible $k$-products of the roots of $p(x)$, multiplied by $(-1)^k$. A $k$-product is the product of $k$ distinct roots of $p(x)$; the roots of $p(x)$ are distinguished from each other by the subscripts in the notation $r_i$ we are employing.

We will illustrate the validity of Vieta's relations by means of the general expression for a univariate polynomial of degree 4:

$$p(x) = x^4 + c_1x^3 + c_2x^2 + c_3x + c_4. \tag{1.3}$$

The fundamental theorem of algebra, together with the factor theorem, allow us to write the polynomial $p(x)$ of expression (1.3) as

$$p(x) = (x - r_1)(x - r_2)(x - r_3)(x - r_4), \tag{1.4}$$

where $r_1$, $r_2$, $r_3$, $r_4$ are the roots of $p(x)$ in (1.3). If we carry out the multiplications on the right side of expression (1.4), proceeding from left to right, we obtain the following results:

$$\begin{array}{lll}
(x - r_1) & & \\
\times\,(x - r_2) & & \\
\hline
x^2 - r_1x & & \\
\quad -r_2x & & + r_1r_2 \\
\hline
x^2 - (r_1 + r_2)x + r_1r_2 & &
\end{array}$$

$$\begin{array}{llll}
[x^2 - (r_1 + r_2)x + r_1r_2] & & & \\
\times\,(x - r_3) & & & \\
\hline
x^3 - (r_1 + r_2)x^2 & + r_1r_2x & & \\
\quad -r_3x^2 & + (r_1r_3 + r_2r_3)x & & - r_1r_2r_3 \\
\hline
\multicolumn{4}{l}{x^3 - (r_1 + r_2 + r_3)x^2 + (r_1r_2 + r_1r_3 + r_2r_3)x - r_1r_2r_3}
\end{array}$$

$$\begin{array}{l}
[x^3 - (r_1 + r_2 + r_3)x^2 + (r_1r_2 + r_1r_3 + r_2r_3)x - r_1r_2r_3] \\
\times\,(x - r_4) \\
\hline
x^4 - (r_1 + r_2 + r_3)x^3 + (r_1r_2 + r_1r_3 + r_2r_3)x^2 \qquad\qquad - r_1r_2r_3x \\
\qquad -r_4x^3 + (r_1r_4 + r_2r_4 + r_3r_4)x^2 - (r_1r_2r_4 + r_1r_3r_4 + r_2r_3r_4)x + r_1r_2r_3r_4 \\
\hline
x^4 - (r_1 + r_2 + r_3 + r_4)x^3 \\
\qquad +(r_1r_2 + r_1r_3 + r_1r_4 + r_2r_3 + r_2r_4 + r_3r_4)x^2 \\
\qquad\qquad -(r_1r_2r_3 + r_1r_2r_4 + r_1r_3r_4 + r_2r_3r_4)x \\
\qquad\qquad\qquad +r_1r_2r_3r_4
\end{array} \tag{1.5}$$

By equating the coefficients in (1.3) to those at the bottom of (1.5), according to like terms between the two expressions, we immediately note that

$$-c_1 = r_1 + r_2 + r_3 + r_4,$$

$$c_2 = r_1r_2 + r_1r_3 + r_1r_4 + r_2r_3 + r_2r_4 + r_3r_4,$$

$$-c_3 = r_1r_2r_3 + r_1r_2r_4 + r_1r_3r_4 + r_2r_3r_4,$$

$$c_4 = r_1r_2r_3r_4.$$

By mathematical induction, it is possible to prove that Vieta's relations are valid for any univariate polynomial $p(x)$ of degree $n \in \mathbb{N}\backslash\{1\}$.

## Complex Numbers

Since both the roots and the coefficients of the univariate polynomials with which we'll be working belong to the field of complex numbers, it may be convenient to review some basic properties of these objects. A complex number $z$ has the form

$z = x + yi$, where $x, y \in \mathbb{R}$, $i = \sqrt{-1}$, so that $i^2 = -1$.

Geometrically, we can visualize a complex number as a point located on a two-dimensional plane, which we call $\mathbb{C}$:

$\mathbb{C} = \{(x, y): x, y \in \mathbb{R}\}$ (Field of complex numbers).

$x \coloneqq \text{Re}(z)$ is the real part of $z$, while $y \coloneqq \text{Im}(z)$ is the imaginary part of $z$.

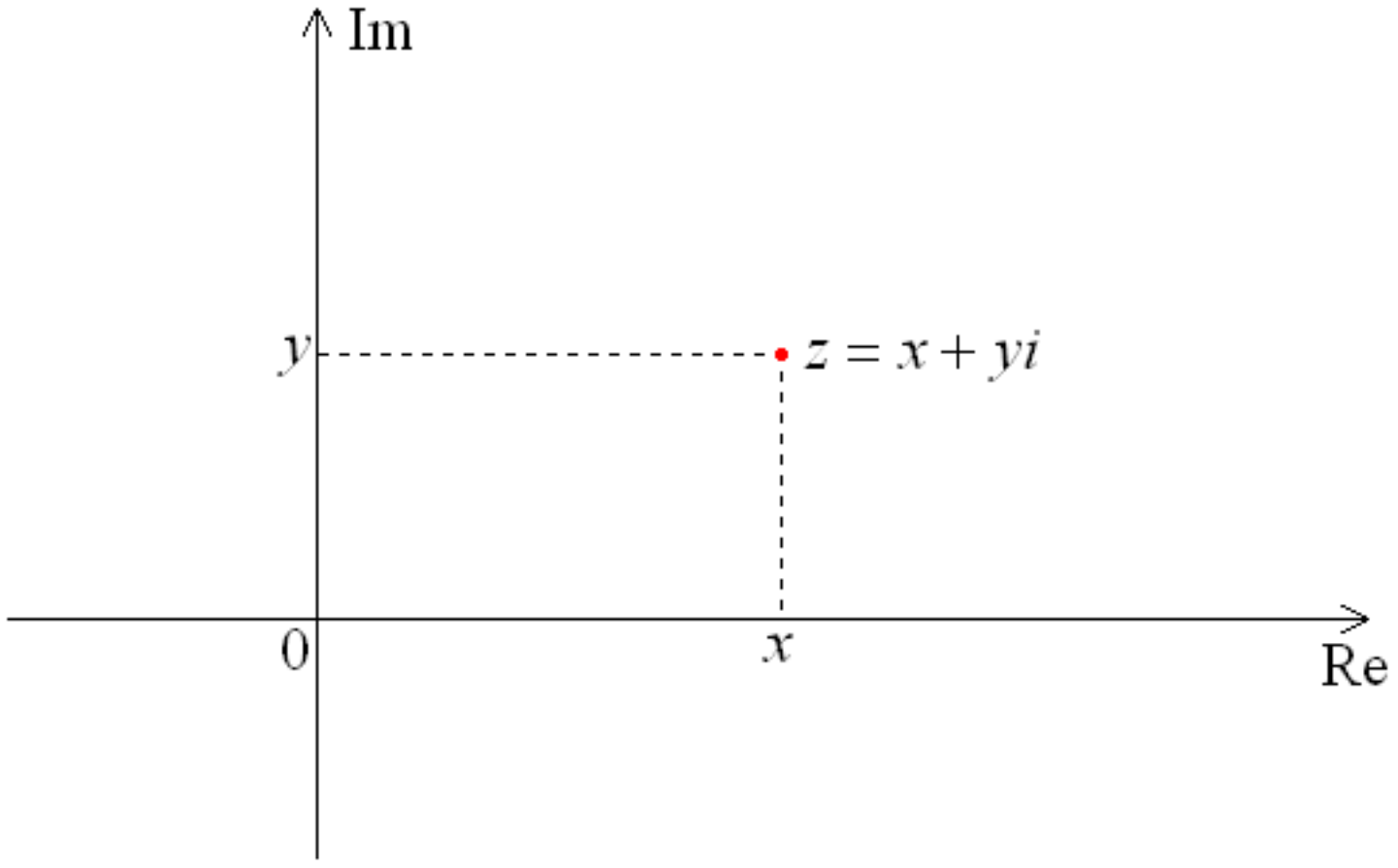

**Fig. 1.1**. Field of complex numbers $\mathbb{C}$, seen as a two-dimensional plane where the horizontal axis is called the *real axis*, and the vertical axis is called the *imaginary axis*. A complex number $z = x + yi$ can be seen as a point (or as a vector) contained in this plane.

In a complex number $z$, the imaginary part $\text{Im}(z)$ is multiplied by the *imaginary unit* $i = \sqrt{-1}$, which can also be written in terms of two-dimensional coordinates like this: $i = (0,1)$; moreover, any complex number $z$ can be written in terms of two-dimensional coordinates: thus, $x + yi$ and $(x, y)$ are two equivalent notations representing the same complex number.

In what follows, we will express a complex number either with the notation of two-dimensional coordinates $z = (x, y)$, or by means of algebraic notation $z = x + yi$, according to what is most convenient. Later, we will see that it is also advantageous to use two-dimensional coordinate notation to manipulate a complex number as if it were a vector in $\mathbb{C}$.

## Basic Operations of Complex Numbers

If $z_1 = (x_1, y_1) \in \mathbb{C}$ and $z_2 = (x_2, y_2) \in \mathbb{C}$,

Sum. $z_1 + z_2 = (x_1 + x_2, y_1 + y_2) \in \mathbb{C}$

Multiplication. $z_1 \cdot z_2 = (x_1 + y_1 i) \cdot (x_2 + y_2 i)$

$$= x_1x_2 + (x_1y_2 + x_2y_1)i + y_1y_2i^2$$

$$= (x_1x_2 - y_1y_2) + (x_1y_2 + x_2y_1)i$$

$$= (x_1x_2 - y_1y_2, x_1y_2 + x_2y_1) \in \mathbb{C}$$

Multiplication by a scalar α. If $\alpha \in \mathbb{R}$, $z = (x, y) \in \mathbb{C}$,

$$\alpha z = (\alpha x, \alpha y) \in \mathbb{C}$$

Multiplicative inverse or reciprocal of a complex number.

$$z^{-1} = \frac{1}{z} = \frac{1 + 0i}{x + yi} = \frac{(1 + 0i)(x - yi)}{(x + yi)(x - yi)} = \frac{x - yi}{x^2 + y^2} = \frac{1}{x^2 + y^2}(x, -y) \in \mathbb{C}$$

Norm (or modulus) of a complex number.

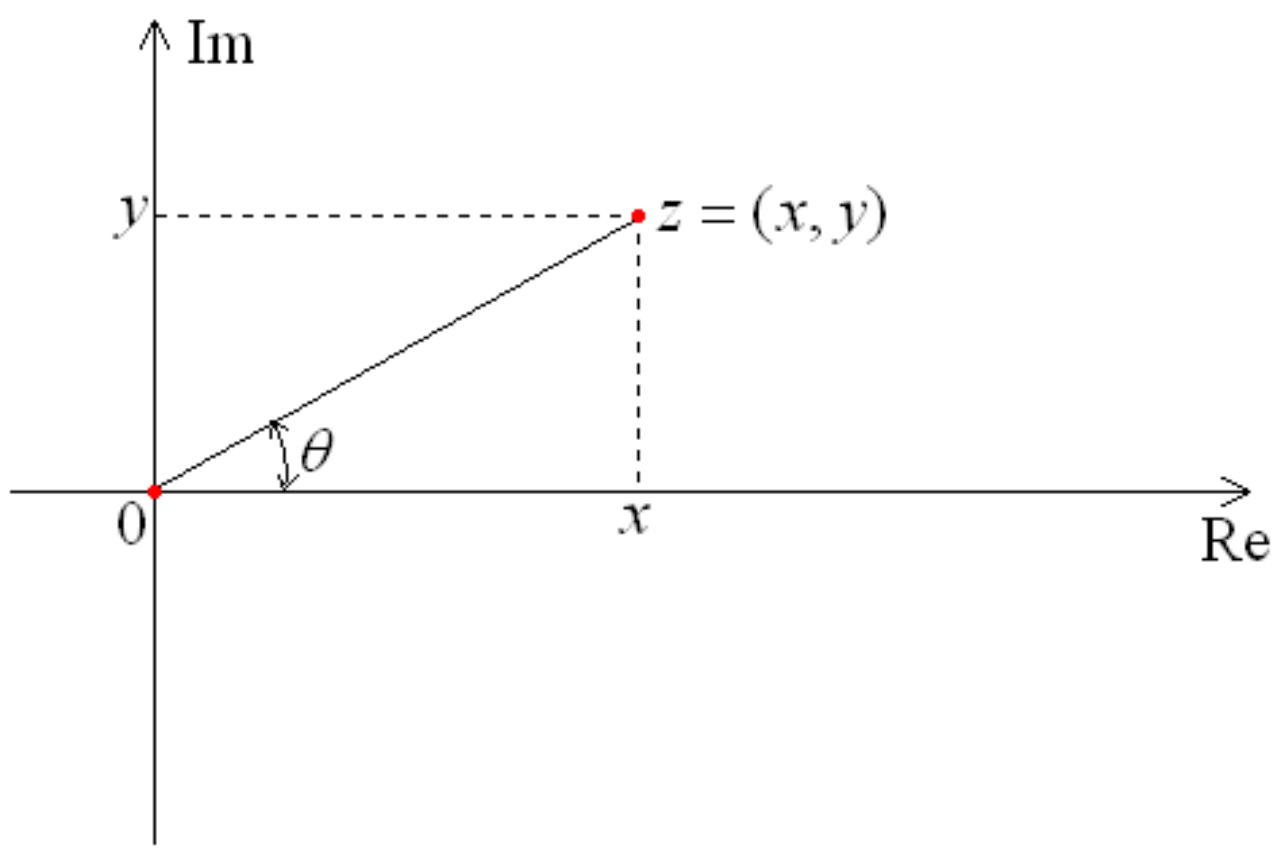

**Fig. 1.2**. A complex number $z = x + yi$, seen as a vector in $\mathbb{C}$. The length of the line segment from point 0 to point $z$ is called the norm or modulus of $z$. The angle $\theta$ that vector $z$ forms with the positive real axis is called the argument of $z$.

We denote the norm of $z$ as $|z|$, and we define it as:

$|z| :=$ length of line segment $\overline{0z}$.

Thus, $|z| = \sqrt{x^2 + y^2}$ (according to the Pythagorean theorem). From here, we see that

$z^{-1} = \frac{1}{|z|^2}(x, -y).$

Conjugate of a complex number.

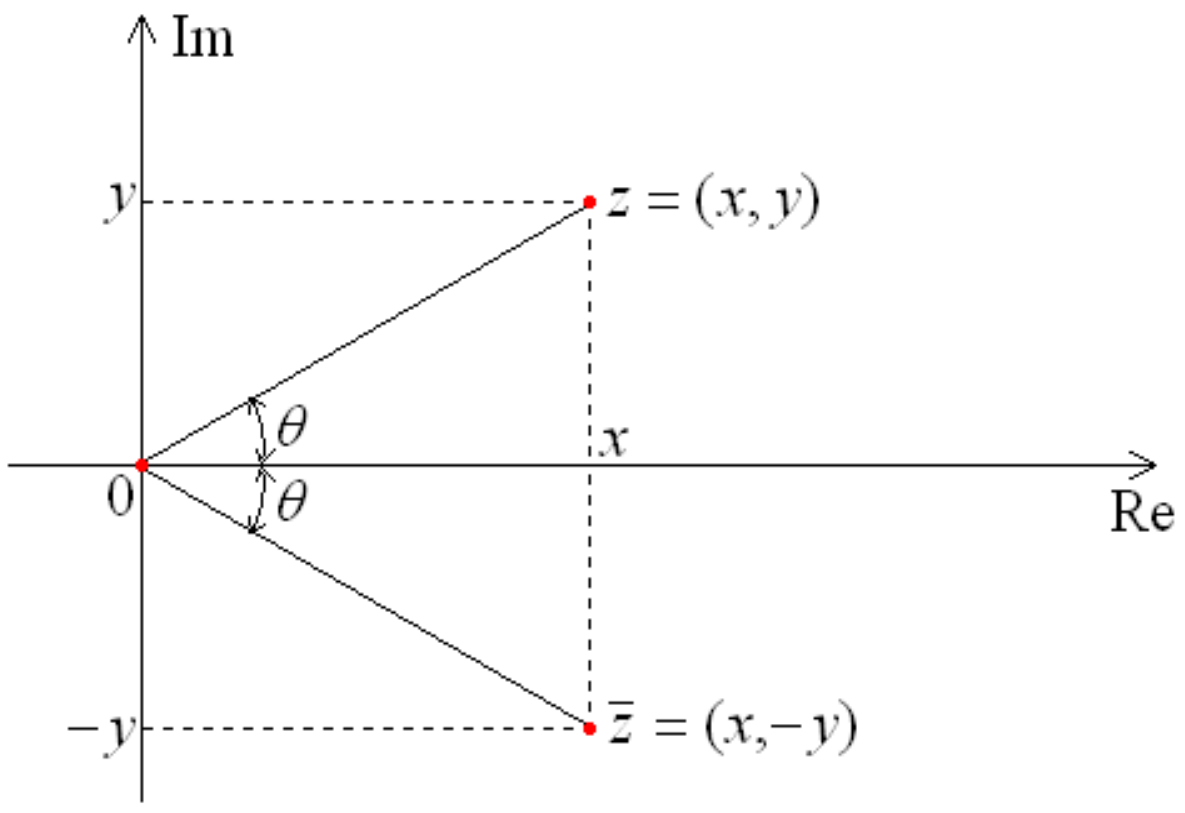

**Fig. 1.3**. The conjugate $\bar{z}$ of a complex number $z = x + yi$ is its reflection with respect to the real axis; thus, $\bar{z} = x - yi$.

The conjugate of $z = (x, y)$ is defined as $\bar{z} = (x, -y)$; when conjugating a complex number, the only thing that changes is the sign of the imaginary part.

Then, we see that

$$z^{-1} = \frac{\bar{z}}{|z|^2} \, . \qquad (1.6)$$

Argument of a complex number.

The argument of $z$, denoted as $\arg(z)$, is the angle $\theta$ shown in figure 1.2. Thus,

$$\arg(z) = \theta = \arctan\left(\frac{y}{x}\right).$$

There is an alternative way of expressing a complex number, in terms of its norm and its argument:

$z = |z|(\cos\theta, \sin\theta)$; this is the *polar form* of a complex number.

$\text{Re}(z) = |z|\cos\theta, \quad \text{Im}(z) = |z|\sin\theta.$

**REMARKS:**

1. A complex number can be thought of as a point, or as a vector in complex plane $\mathbb{C}$.

2. The sum of complex numbers is analogous to the sum of vectors in the Cartesian plane $\mathbb{R}^2$ (see figure 1.4).

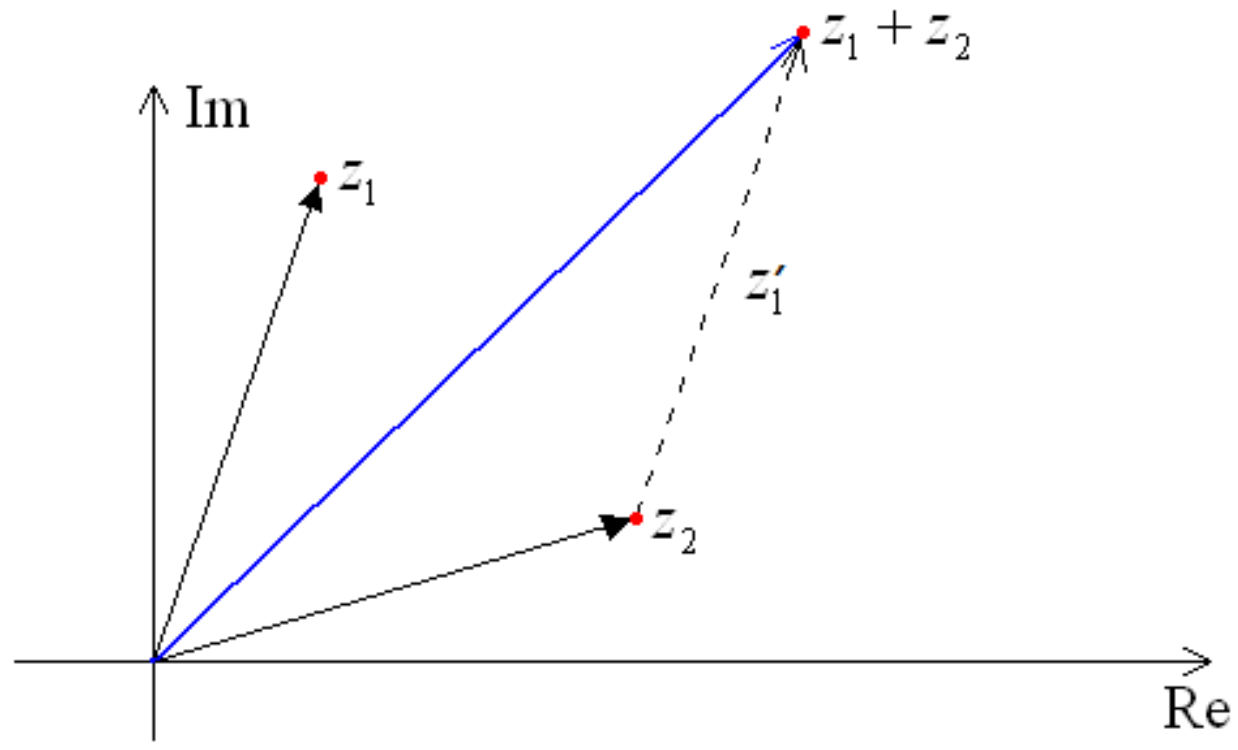

**Fig. 1.4**. The sum of two complex numbers $z_1$, $z_2$ is analogous to the sum of two-dimensional vectors in $\mathbb{R}^2$, using the parallelogram rule of vector addition; $z_1'$ is the translation of vector $z_1$.

3. A graphical procedure for multiplying two complex numbers $z_1 = (x_1, y_1)$ and $z_2 = (x_2, y_2)$ is as follows:

   a. On complex plane $\mathbb{C}$, superimpose a new plane $\mathbb{C}_{z_1}$, which has as "real" unit the vector $z_1$. The origin coincides in both planes.

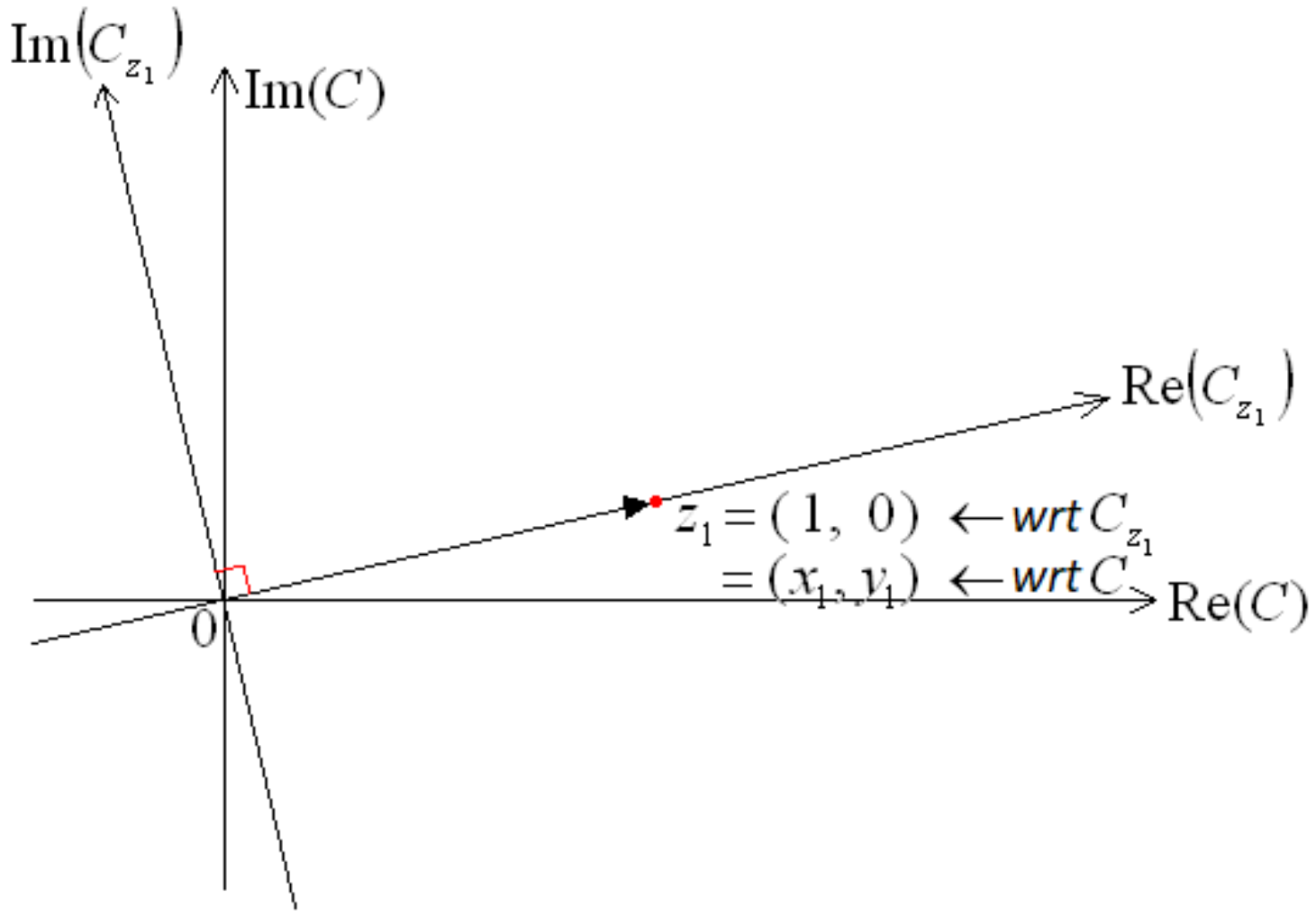

**Fig. 1.5**. Superimposition of plane $\mathbb{C}_{z_1}$ on plane $\mathbb{C}$. The vector $z_1 = (x_1, y_1)$, with respect to plane $\mathbb{C}$, is also the vector (1,0) with respect to plane $\mathbb{C}_{z_1}$.

   b. Locate the coordinates of $z_2$ in plane $\mathbb{C}_{z_1}$ in the usual way, taking vector $z_1$ as if it were the real unit (1,0).

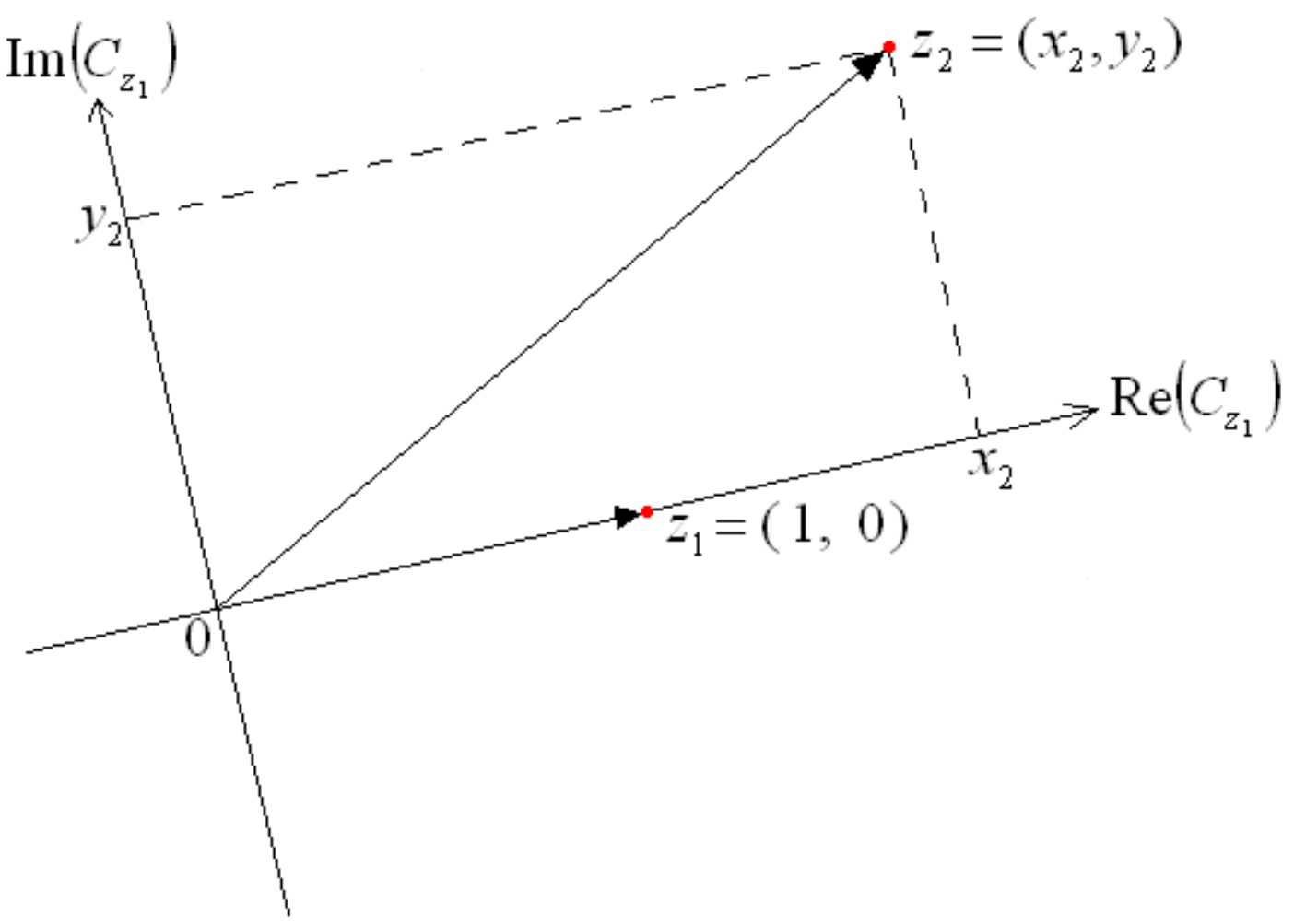

**Fig. 1.6**. Vector $z_2$ in the plane $\mathbb{C}_{z_1}$. The coordinates $x_2$, $y_2$ shown in this figure are relative to the "real unit" $z_1$.

c. Convert the coordinates of point $z_2$ located in plane $\mathbb{C}_{z_1}$ into coordinates with respect to the original plane $\mathbb{C}$. In this way we obtain $z_1 \cdot z_2$.

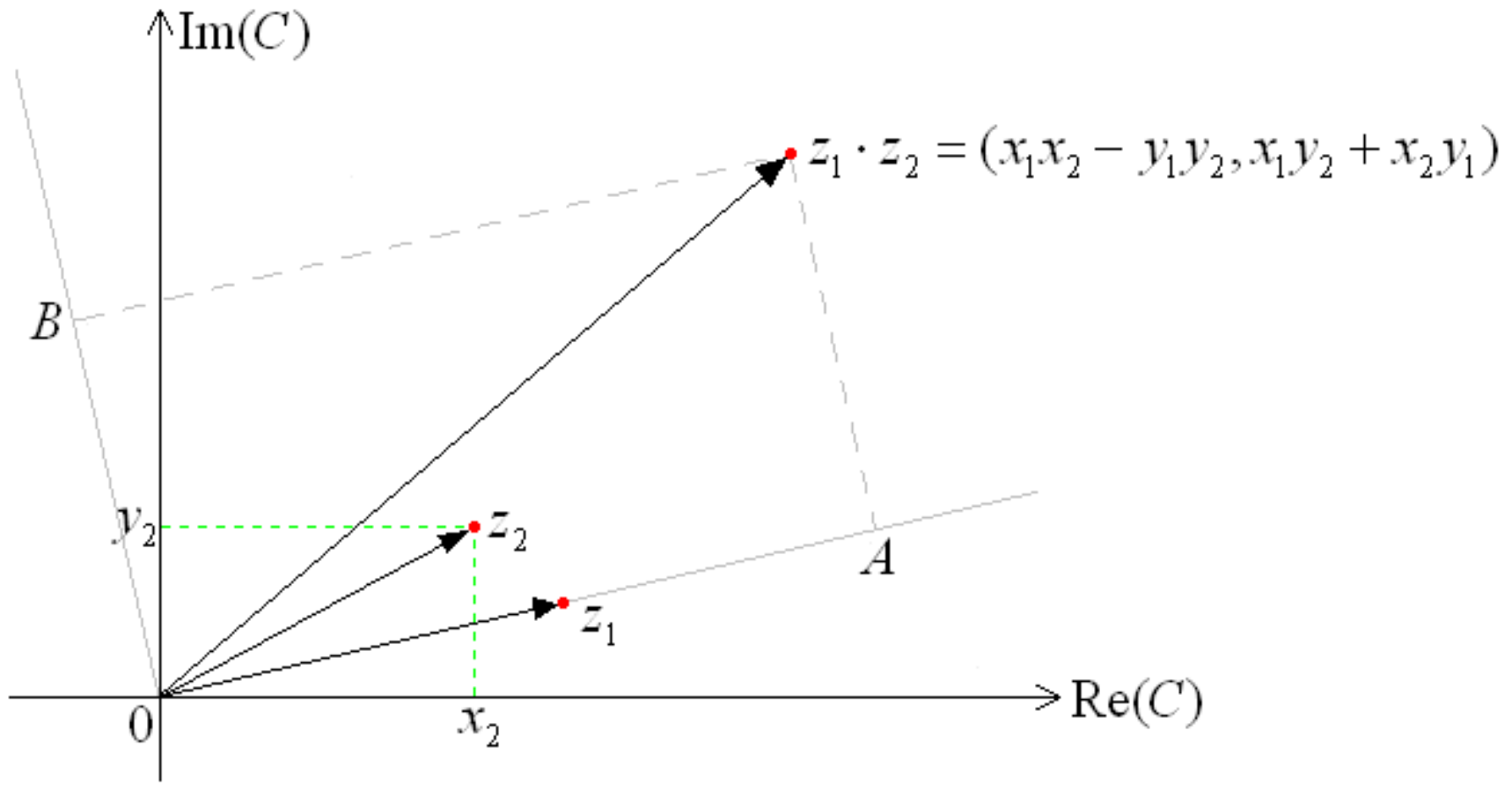

**Fig. 1.7**. The coordinates of $z_2$ in plane $\mathbb{C}_{z_1}$ from figure 1.6 are converted into coordinates relative to the complex plane $\mathbb{C}$. Thus $z_2$ in $\mathbb{C}_{z_1}$ becomes $z_1 \cdot z_2$ in $\mathbb{C}$. Rectangles $0x_2z_2y_2$ and $0A(z_1 \cdot z_2)B$ are similar.

4. If $z_1 = |z_1|(\cos\theta_1, \sin\theta_1)$ and $z_2 = |z_2|(\cos\theta_2, \sin\theta_2)$, then

$$z_1 \cdot z_2 = |z_1||z_2|[\cos(\theta_1 + \theta_2), \sin(\theta_1 + \theta_2)],$$
$$z_1/z_2 = (|z_1|/|z_2|)[\cos(\theta_1 - \theta_2), \sin(\theta_1 - \theta_2)].$$

Note that $|(\cos\alpha, \sin\alpha)| = 1 \qquad \forall\alpha \in \mathbb{R}$.

**5. De Moivre's theorem**

$$[|z|(\cos\theta, \sin\theta)]^n = |z|^n(\cos n\theta, \sin n\theta), n \in \mathbb{Z}. \quad (1.7)$$

**6. *n*-th roots of a complex number**

If $z = |z|(\cos\theta, \sin\theta) \in \mathbb{C}$, an *n-th root of* $z$ is

$$|z|^{1/n}\left[\cos\left(\frac{\theta+2\pi k}{n}\right), \sin\left(\frac{\theta+2\pi k}{n}\right)\right], \quad k = 0,1,2,\dots,n-1. \quad (1.8)$$

From remark 4, we see that if $|z_2| = 1$, $z_1 \cdot z_2$ is the same as rotating vector $z_1$ by an angle $\theta_2$ (counterclockwise if $\theta_2 > 0$, clockwise if $\theta_2 < 0$), starting from its original position. The same can be said of $z_1/z_2$, but the direction of rotation of $z_1$ is opposite with respect to that of $z_1 \cdot z_2$.

## Lines and Circumferences in the Plane of Complex Numbers

**Lines**. A line is a set of complex numbers of the form

$L = \{p + tv \mid p, v \in \mathbb{C} \text{ fixed}, \quad t \in \mathbb{R} \text{ variable}\}$.

In figure 1.8 we can see a geometric illustration of this concept.

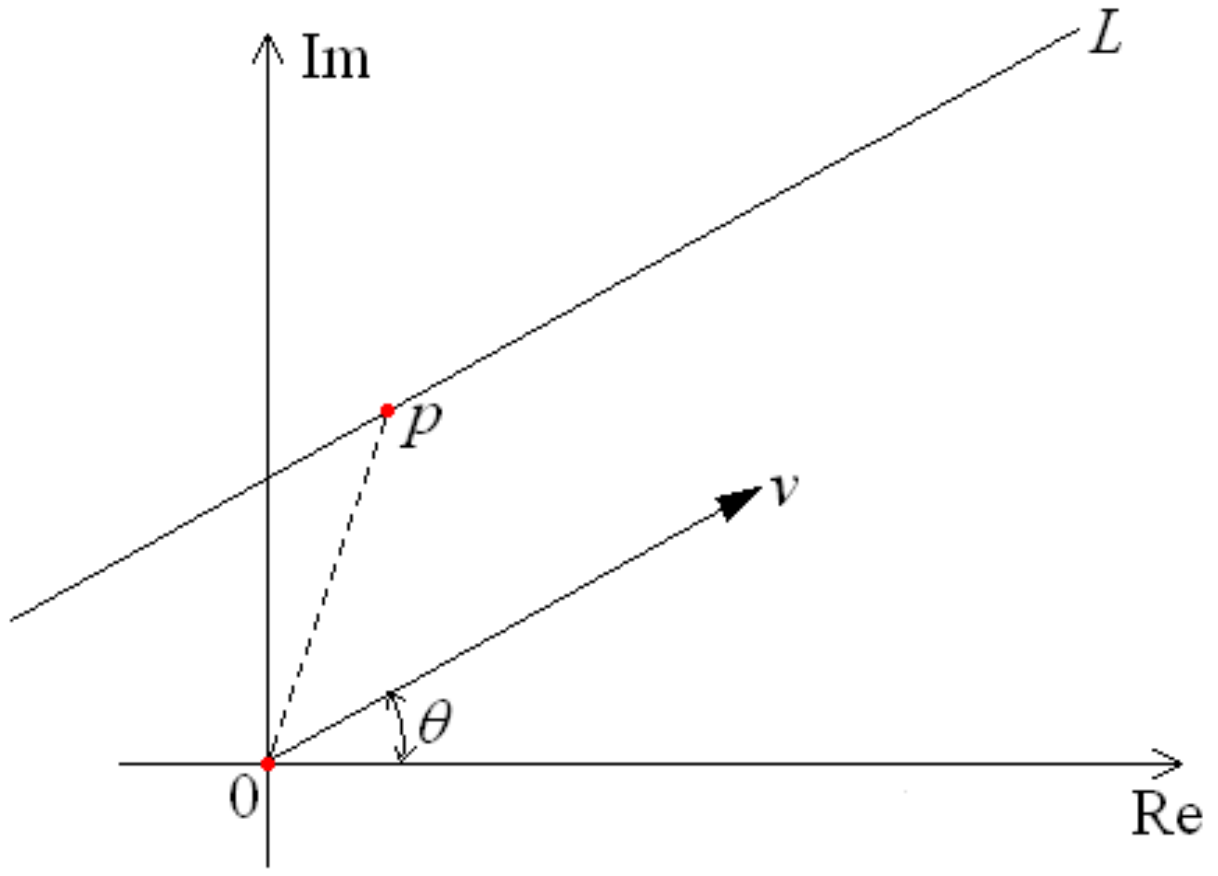

**Fig. 1.8**. A line $L$ in the complex plane $\mathbb{C}$, constructed from a fixed point $p$ and a direction vector $v$.

By convention, we will call $p$ *the fixed point of line* $L$, and $v$ *the direction vector of line* $L$. We will usually take $|v| = 1$, and say that the direction of $L$ is $\theta = \arg(v)$. The real parameter $t$ allows us to entirely construct line $L$ taking as references $p$ and $v$.

Next, we will see a concept that will be very useful in the construction of the LC method for the approximation of the roots of polynomial (1.1):

**Semi-line derived from $L$.** Suppose we have a line $L = \{p + tv\}$ with direction $\theta = \arg(v)$. We define the semi-line derived from $L$ as the set

$$L_d = \{(p + tv)(p - tv)\} = \{p^2 - t^2v^2\} = \{p^2 + t^2(-v^2)\}.$$

$L_d$ is a semi-line, since parameter $t^2 \in \mathbb{R}^+ \cup \{0\}$, while $p^2$ and $-v^2$ are fixed complex values, which serve respectively as the fixed point and the direction vector of $L_d$. The direction of $L_d$ is $2\theta - \pi$, since

$\arg(v^2) \quad = 2\theta$ (By remark 5, De Moivre's theorem),

$\arg(-v^2) = 2\theta - \pi$ (By remark 4, quotient $z_1/z_2$, taking $z_1 = v^2$ and $z_2 = -1$).

In figure 1.9 we can see, schematically, the spatial direction of a semi-line $L_d$ derived from $L$, in relation to the direction of the corresponding line $L$ that generates it:

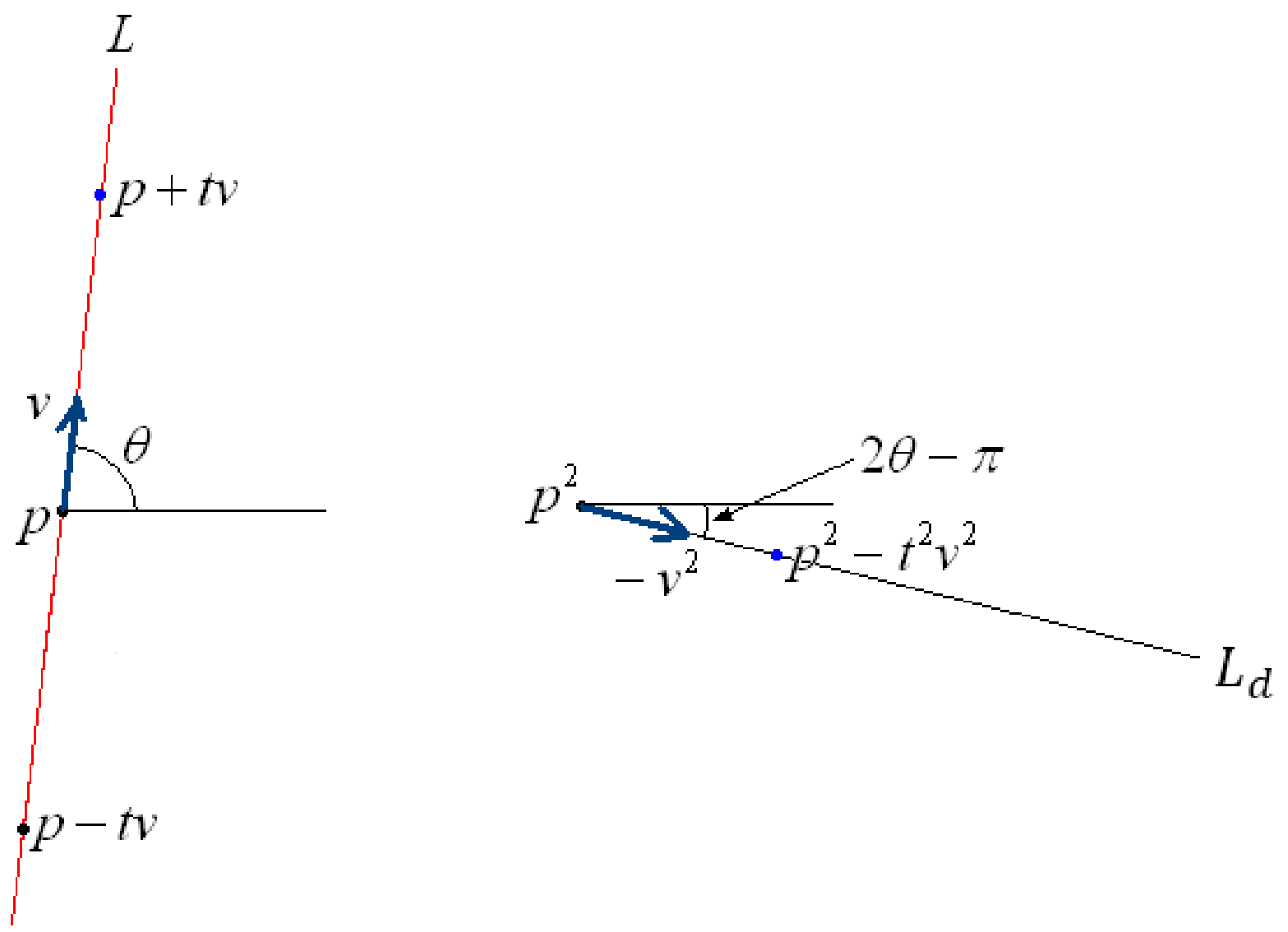

**Fig. 1.9**. A line $L$, and its corresponding derived semi-line, $L_d$.

Note that it is also true to state that $\arg(-v^2) = 2\theta + \pi$, since we could also have considered the product $z_1 \cdot z_2$ from remark 4 instead of the quotient $z_1/z_2$, taking, as we said before, $z_1 = v^2$ and $z_2 = -1$.

To simplify both the calculations and the presentation of results, we will see in coming chapters that it is useful to limit the angles of direction vectors to the interval $[-\pi, \pi)$, since this range of values covers minimally all possible directions of a line, moving as little away as possible from 0.

**Möbius transformations**. A Möbius transformation is a bijective mapping (also called one-to-one mapping or invertible mapping) from $\mathbb{C}$ to $\mathbb{C}$, having the form

$$z \to \frac{az+b}{cz+d}\,, \tag{1.9}$$

where $a, b, c, d$ are fixed complex numbers, which must meet the requirement $ad - bc \neq 0$; otherwise, the quotient $\frac{az+b}{cz+d}$ would be reduced to a constant value $\forall z \in \mathbb{C}$.

This transformation has several interesting properties. A property that will be very useful to us is the following:

In equation (1.9), let's take $a = 0, b = 1, c = 1, d = 0$. With this choice of parameters, the Möbius transformation acts as a multiplicative inversion

$$z \to \frac{1}{z}\,.$$

**Property 1.1**. If $z$ is a point on a line $L$ in the complex plane $\mathbb{C}$ that does not pass through the origin, then $1/z$ is a point on a circumference $zC$ in the complex plane $\mathbb{C}$ that passes through the origin without containing it (see figure 1.10). In other words, the mapping $z \to 1/z$ transforms line $L$ into circumference $zC$.

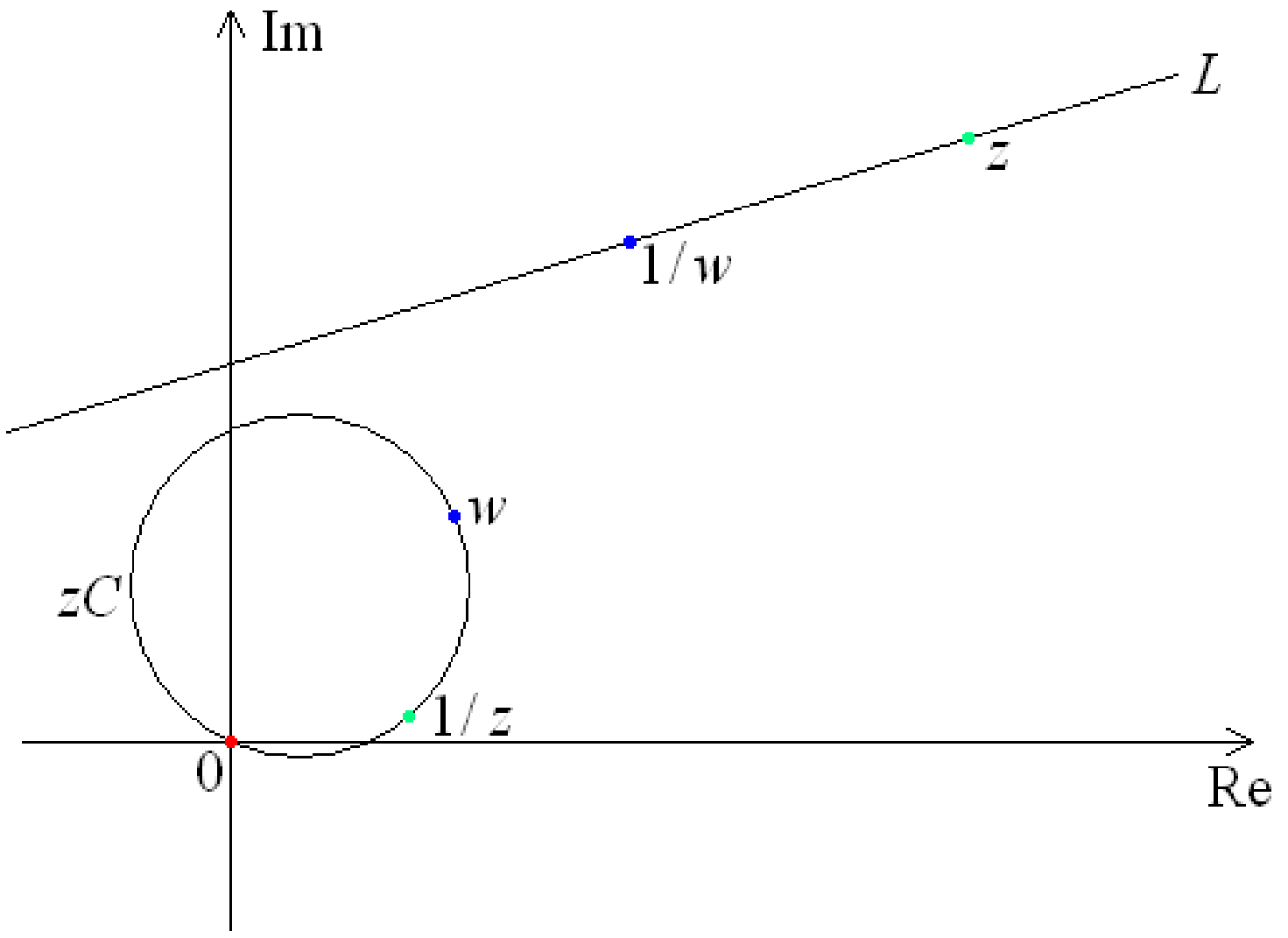

**Fig. 1.10**. A line $L$ is mapped into a circumference $zC$ by means of the Möbius transformation $z \to 1/z$. Conversely, given the bijectivity of Möbius transformations, a circumference $zC$ is mapped into a line $L$ by the Möbius transformation $w \to 1/w$.

**Property 1.1A**. If $w$ is a point on a circumference $zC$ in the complex plane $\mathbb{C}$ that passes through the origin without containing it, then $1/w$ is a point on a line $L$ in the complex plane $\mathbb{C}$ that does not pass through the origin. In other words, the mapping $w \to 1/w$ transforms circumference $zC$ into line $L$.

Property 1.1A is the converse of property 1.1, as can also be seen in figure 1.10. In summary, a Möbius transformation of the form $z \to 1/z$ converts a line that does not pass through the origin into a circumference that does pass through the origin and vice versa, given the bijectivity of the Möbius transformation. If a line passes through the origin, the Möbius transformation $z \to 1/z$ converts that line into another line that also passes through the origin (although, in general, with a different direction with respect to that of the original line).

Below, we demonstrate properties 1.1 and 1.1A; these demonstrations are based on ideas contained in notes by Dr. Carl Eberhart (University of Kentucky), which can be consulted in the following link:

http://www.ms.uky.edu/~carl/ma502/html/moebsln1.html

**Demonstration of property 1.1A**. Let $w \neq 0$ be a point on the circumference $zC$ which passes through the origin. There is a vector $u \in \mathbb{C}$ with $|u| = 1$ such that $u \cdot w$ is a point on the circumference $zC'$ that passes through the origin and has as its center $c$ a point on the real axis Re($\mathbb{C}$) (without loss of generality, we will assume that $c > 0$). Both circumferences $zC$ and $zC'$ have the same radius $r$; $u$ helps us to rotate $zC$ around the origin until its center coincides with the real axis Re($\mathbb{C}$). The point $u \cdot w$ has the form $c + r(\cos\theta, \sin\theta)$; $r$ is the radius of $zC'$ and is equal to $c$, so $u \cdot w = c + c(\cos\theta, \sin\theta)$.

Now,

$$\frac{1}{u \cdot w} = \frac{1}{c + c(\cos\theta, \sin\theta)} = \frac{1}{(c + c\cos\theta, c\sin\theta)}$$

$$\frac{1}{u \cdot w} = \frac{(c + c\cos\theta, -c\sin\theta)}{(c + c\cos\theta, c\sin\theta)(c + c\cos\theta, -c\sin\theta)}$$

$$\frac{1}{u \cdot w} = \frac{1}{(c + c\cos\theta)^2 + (c\sin\theta)^2}(c + c\cos\theta, -c\sin\theta)$$

Let's see in more detail the real part of $\frac{1}{u \cdot w}$:

$$\mathrm{Re}\left(\frac{1}{u \cdot w}\right) = \frac{c + c\cos\theta}{(c + c\cos\theta)^2 + (c\sin\theta)^2} = \frac{c + c\cos\theta}{c^2 + 2c^2\cos\theta + c^2\cos^2\theta + c^2\sin^2\theta}$$

$$\mathrm{Re}\left(\frac{1}{u \cdot w}\right) = \frac{c + c\cos\theta}{c^2 + 2c^2\cos\theta + c^2} = \frac{c + c\cos\theta}{2c^2 + 2c^2\cos\theta} = \frac{c + c\cos\theta}{2c(c + c\cos\theta)}$$

$$\mathrm{Re}\left(\frac{1}{u \cdot w}\right) = \frac{1}{2c}$$

The real part of $\frac{1}{u \cdot w}$ is constant, so $\frac{1}{u \cdot w}$ is on a vertical line whose abscissa is $\frac{1}{2c}$. On the other hand, the imaginary part of $\frac{1}{u \cdot w}$ is:

$$\mathrm{Im}\left(\frac{1}{u\cdot w}\right)=\frac{-c\sin\theta}{(c+c\cos\theta)^2+(c\sin\theta)^2}\qquad \theta\in(-\pi,\pi)$$

Let's see what happens to the imaginary part of $\frac{1}{u\cdot w}$ if we consider different points of $zC'$:

$$\lim_{\theta\to\pi^-}\mathrm{Im}\left(\frac{1}{u\cdot w}\right)=-\infty$$

$$\mathrm{Im}\left(\frac{1}{u\cdot w}\right)=-\frac{1}{2c}\qquad \text{if } \theta=\frac{\pi}{2}$$

$$\mathrm{Im}\left(\frac{1}{u\cdot w}\right)=0\qquad \text{if } \theta=0$$

$$\mathrm{Im}\left(\frac{1}{u\cdot w}\right)=\frac{1}{2c}\qquad \text{if } \theta=-\frac{\pi}{2}$$

$$\lim_{\theta\to-\pi^+}\mathrm{Im}\left(\frac{1}{u\cdot w}\right)=+\infty$$

This behavior indicates that $\mathrm{Im}\left(\frac{1}{u\cdot w}\right)$ covers the entire vertical line $x=\frac{1}{2c}$. Therefore, $\frac{1}{u\cdot w}$ is a vertical line, and $u\cdot\frac{1}{u\cdot w}=\frac{1}{w}$ is a line whose direction vector has argument equal to $\frac{\pi}{2}+\arg(u)$. Note that $\frac{1}{w}$ is an orthogonal line with respect to vector $u$ (see figure 1.11). $\frac{1}{w}$ is a line that does not contain the origin, since $\frac{1}{u\cdot w}$ does not contain the origin. ■

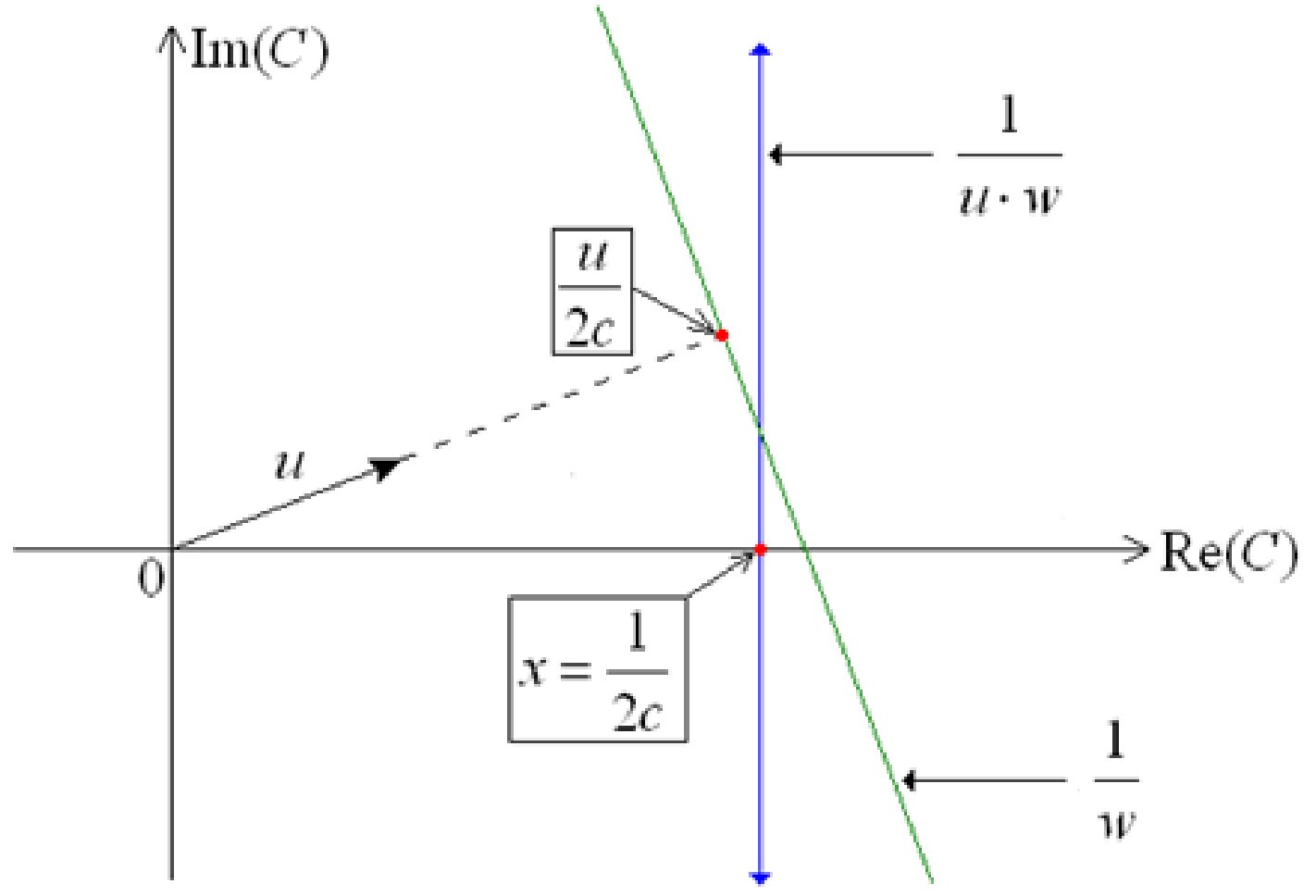

**Fig. 1.11**. The vertical line $1/(u\cdot w)$, when multiplied by complex number $u$ with unitary norm, produces a line $1/w$ that is orthogonal to $u$, and does not contain 0.

**Demonstration of property 1.1**. Since property 1.1A is valid, and the corresponding Möbius transformation $\omega \to 1/\omega$ that converts circumferences that pass through 0 into lines that do not pass through 0 is a bijective mapping, we have that a line $z = 1/\omega$ that does not pass through the origin is transformed, by means of mapping $1/\omega \to \omega$ (or equivalently, in terms of $z$, by means of mapping $z \to 1/z$), into a circumference $1/z = \omega$ that passes through the origin without containing it. ■

| |
|---|
| **Notation:** We will call a circumference that passes through the origin, a *z-circumference.* |

Note that if a z-circumference comes from converting a line $z$ that does not pass through the origin by a Möbius transformation $z \to 1/z$, then the z-circumference $1/z$ does not contain the origin.

Another property of the Möbius transformation that will be very useful to us arises if in mapping (1.9) we take $a = 0$, $b \neq 0$, $c = 1$, $d = 0$.

**Property 1.2**. If $z$ is a point on a line $L$ in the complex plane $\mathbb{C}$ that does not pass through the origin, then $b/z$ is a point on a z-circumference $zC$ in the complex plane $\mathbb{C}$. Conversely, if $z$ is a point on a z-circumference $zC$ in the complex plane $\mathbb{C}$, then $b/z$ is a point on a line $L$ in the complex plane $\mathbb{C}$ that does not pass through the origin.

The effect of the constant complex value $b$ on factor $1/z$ is a change of scale and a rotation with respect to the origin, but in essence property 1.2 states the same as properties 1.1 and 1.1A, just in a slightly more general way.

Under the context of property 1.2, we would then say, abusing of language, that $L \cdot zC = b$ (the line "multiplied" by the z-circumference is a constant).

From here, it is evident that $L = \frac{b}{zC}$ and $zC = \frac{b}{L} = \frac{b}{p+tv}$ ($p \in \mathbb{C}$ is the fixed point of line $L$, $v \in \mathbb{C}$ is the direction vector of $L$, and $t \in \mathbb{R}$ is a parameter that describes each point of $L$ and $zC$).

To conclude this chapter of preliminary concepts, we enunciate two very basic properties about pairs of points on a line in the complex plane that are at the same distance from a third point on that line. These properties will also be useful when developing the LC methodology.

**Property 1.3**. If a line $L$ in the complex plane $\mathbb{C}$ contains point $P = (a + b)/2$ and contains point $a$, then it also contains point $b$.

**Demonstration of property 1.3**. Let $V$ be the vector going from $P$ to $a$; then $P + V = a$, or equivalently, $V = a - P$. Therefore, vector $V$ can be considered as a direction vector for line $L$, just like vector $-V = P - a$. Now, $P - V = P + (P - a) = 2P - a = (a + b) - a = b$; this means that line $L$ contains point $b$. ■

**Property 1.3a**. If a line $L$ in complex plane $\mathbb{C}$ contains points $a$ and $b$, then it also contains point $P = (a + b)/2$.

**Demonstration of property 1.3a**. Let $W$ be the vector going from $a$ to $b$; then $a + W = b$, or equivalently, $W = b - a$. $W$ is a direction vector for line $L$, just like $W/2$. Now, $a + W/2 = a + (b - a)/2 = (a + b)/2 = P$; this means that line $L$ contains point $P$. ■

# Chapter Two: LC Method for Quadratic Polynomial Equations

## Theoretical Aspects

Our objective in this chapter is to explore an interesting geometric relationship that exists between the coefficients $C_1, C_2 \in \mathbb{C}$ of the univariate polynomial

$$z^2 + C_1 z + C_2 \tag{2.1}$$

and the roots $R_1, R_2 \in \mathbb{C} \setminus \{0\}$ of the corresponding equation

$$z^2 + C_1 z + C_2 = 0, \tag{2.2}$$

making use of lines and z-circumferences in the complex plane (see chapter 1: fundamentals of the LC method).

In this chapter we will see how, from this geometric relationship, it is possible to build a computational method that numerically approximates the roots of equation (2.2) from their coefficients; this will constitute the basis of the LC method. In later chapters we will see how to extend the basis of this numerical method towards univariate polynomials of degree greater than 2, taking advantage of the parallel processing capabilities in current computing devices.

As we saw in chapter 1, a direct algebraic relationship between the roots of equation (2.2), and the coefficients of polynomial (2.1), is given by Vieta's relations:

$$C_1 = -(R_1 + R_2) \tag{2.3-a}$$

$$C_2 = R_1 R_2; \tag{2.3-b}$$

Sometimes it is possible to use relationships (2.3-a), (2.3-b) directly to mentally find the roots of (2.2) without the help of an electronic computing device, especially when $R_1, R_2$ are small integers, and what is sought is to factor expression (2.1); in most practical situations, however, the fastest and most direct option is to use a calculator or computer to carry out the calculations indicated by the famous general formula for solving quadratic equations:

$$r = \left(-b \pm \sqrt{b^2 - 4ac}\,\right)/(2a), \tag{2.4}$$

where for the case of equation (2.2), the values $a = 1$, $b = C_1$, and $c = C_2$ are directly plugged into formula (2.4). It is a well-known fact that finding the roots of univariate polynomial equations by means of general analytic formulas like (2.4), is only possible for univariate polynomials of degree less than 5; for univariate polynomials of degree 5 or greater, generally it is only possible to approximate numerically the roots of their corresponding polynomial equations. To learn more about this fascinating scientific discovery, the reader is invited to consult the history of Gerolamo Cardano, Niccolò Fontana Tartaglia, Lodovico de Ferrari, and the Abel-Ruffini Theorem. An excellent and stimulating publication where you can consult the history of these famous mathematicians is [2.1]; a proof of the Abel-Ruffini theorem that makes use of Galois theory can be found in [2.2].

The LC method that we will describe throughout these chapters of course is not intended to be an analytical method for finding roots of polynomial equations of any degree; instead, it seeks to be a numerical method useful in obtaining good initial approximations to the roots of univariate polynomial equations of any degree, taking advantage of parallel processing capabilities in current computing devices; after the application of the LC method, the obtained initial approximations can be refined using the Newton-Raphson method or some other numerical method with local quadratic or super-quadratic convergence; this issue of refining initial approximations of polynomial roots, however, will not be addressed in this work, as our main objective is to describe the LC method per se. In this chapter we will look at the basic elements and concepts of the LC method; then, in the following chapters, when dealing with polynomial equations of degree greater than 2, we will address additional basic concepts that will complete the description of the LC method. For the remainder of this chapter, we will assume that $R_1 \neq R_2$; later we will deal with the case of repeated roots for univariate polynomials of degree two.

We begin our description of the LC method by referring to figure 2.1; if line $\ell_1$ in complex plane $\mathbb{C}$ contains roots $R_1, R_2$, then it also contains the point

$$P_1 = (R_1 + R_2)/2 = -C_1/2. \qquad (2.5)$$

This is a direct consequence of result 1.3a, seen in chapter 1. Note that point $P_1$ is precisely the midpoint of rectilinear segment $\overline{R_1 R_2}$. Moreover, if $0 \notin \ell_1$, the z-circumference $C_2/\ell_1$ contains points $R_2, R_1$, and $C_2/P_1$ (property 1.2 of Möbius transformations seen in chapter 1 guarantees that $C_2/\ell_1$ is a circumference).

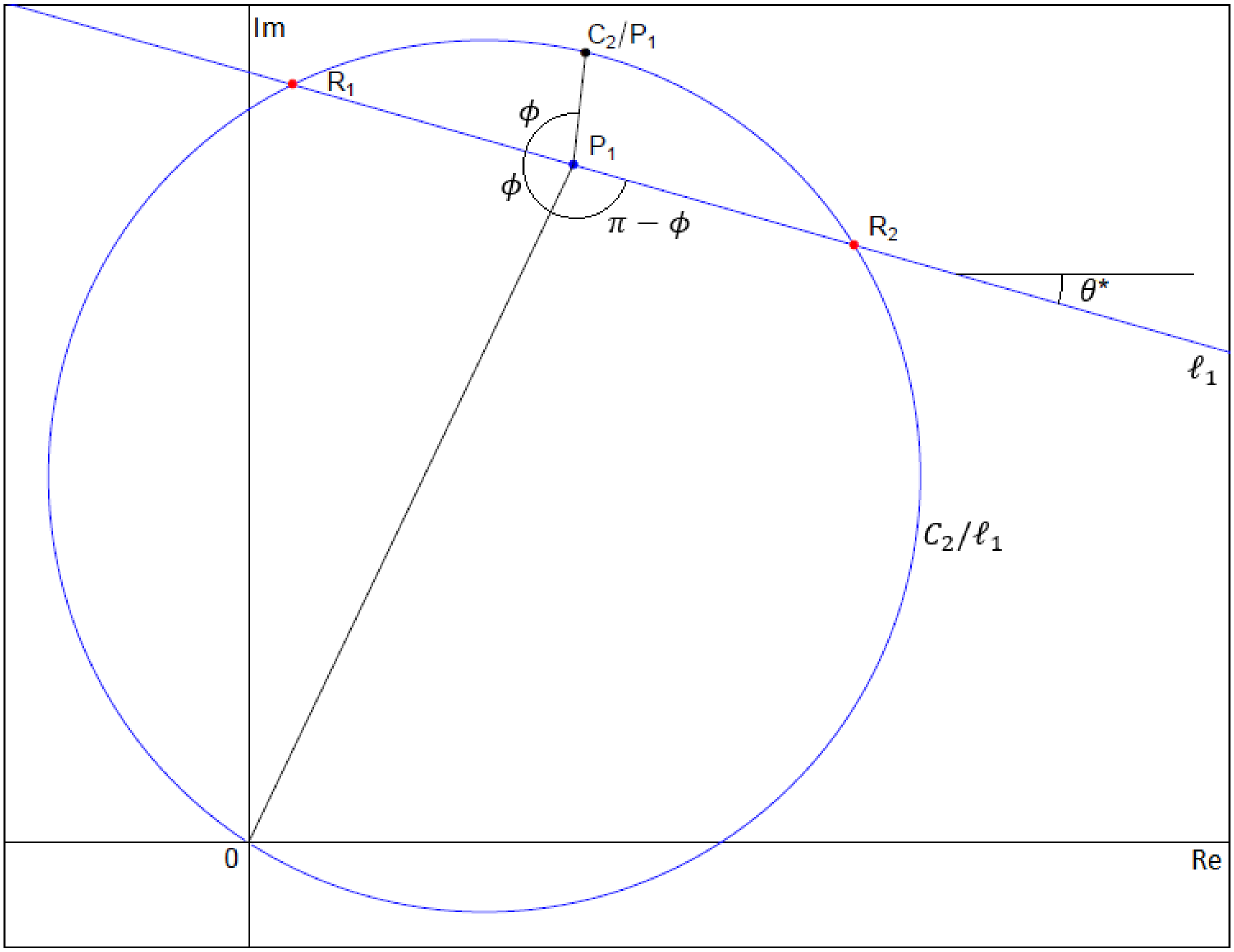

**Fig. 2.1**. If line $\ell_1 \subset \mathbb{C}\backslash\{0\}$ contains the roots $R_1, R_2$ of polynomial equation (2.2), then z-circumference $C_2/\ell_1$ also contains $R_1, R_2$.

Point $P_1$ defined by expression (2.5) is called *the fixed point* of line $\ell_1$; we can imagine that we rotate $\ell_1$, always keeping point $P_1$ motionless, until $\ell_1$ coincides with roots $R_1, R_2$ of equation (2.2); initially, we do not know the inclination angle $\theta^*$, with respect to the positive real axis in $\mathbb{C}$, which $\ell_1$ must have in order to coincide with the roots, since we do not know $R_1, R_2$; we only know coefficients $C_1, C_2$ of polynomial equation (2.2). In what follows, it will be useful to denote the line containing $R_1, R_2$ as $\ell_1(\theta^*)$.

Now we will see how to numerically approximate $\theta^*$ from the coefficients $C_1, C_2$ of polynomial (2.1). For this purpose, we will use the concept of *semi-line derived from another line*, seen in chapter 1, applying it directly to line $\ell_1(\theta^*)$, to obtain derived semi-line $\ell_d$ (see figure 2.2). Any point in semi-line $\ell_d$ is obtained by multiplying a unique pair of points at $\ell_1(\theta^*)$, each of which are at the same distance from fixed point $P_1$; in this way,

$\ell_d = \{p \in \mathbb{C} | p = (P_1 - tV)(P_1 + tV); t \in \mathbb{R} \text{ variable}, V \in \mathbb{C} \text{ fixed}\}$; $V$ is the direction vector of $\ell_1(\theta^*)$.

Note that if $\ell_1(\theta)$ contains roots $R_1$ and $R_2$, then its derived semi-line $\ell_d$ contains two key points: $P_1^2$ and $C_2$; see equation (2.3-b). These two key points are obtained directly from coefficients $C_1, C_2$ of polynomial (2.1). Under these circumstances, the inclination angle of $\ell_d$ with respect to the positive real axis in $\mathbb{C}$ is $2\theta^* + \pi$, according to the results obtained in chapter 1 for derived semi-lines.

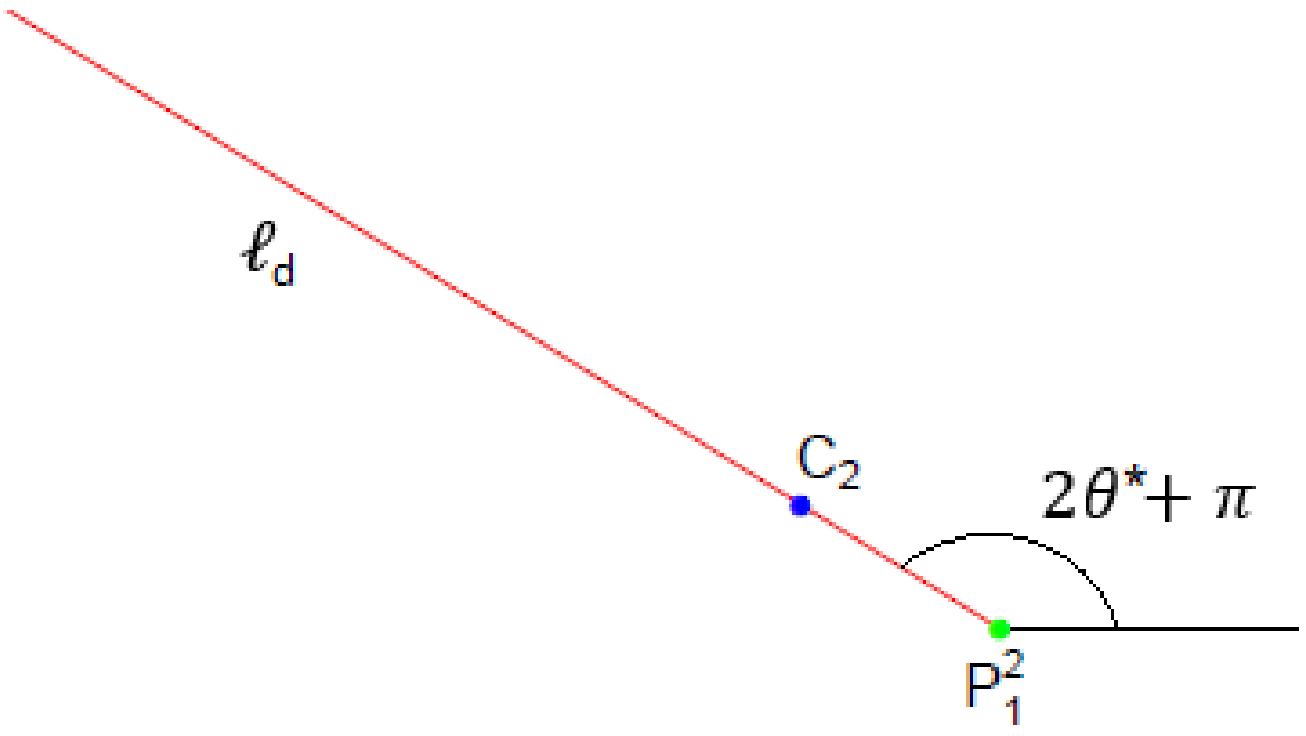

**Fig. 2.2**. Semi-line $\ell_d$ derived from $\ell_1(\theta^*)$, which geometrically expresses a direct relationship between coefficients $C_1, C_2$ of equation (2.2).

From the geometric-algebraic relationship illustrated in figure 2.2, we can obtain directly $\theta^*$; if we start from a direction vector $W$ for $\ell_d$ given by

$$W = C_2 - P_1^2, \tag{2.6}$$

then we have

$$2\theta^* + \pi = \arg(W)$$

$$2\theta^* = \arg(-W)$$

$$\theta^* = \arg(-W)/2. \tag{2.7}$$

Therefore, computationally we can obtain $\theta^*$ by means of the following expression:

$$\theta^* = \text{atan2}[\text{Im}(C_1^2/4 - C_2), \text{Re}(C_1^2/4 - C_2)]/2, \tag{2.8}$$

where $\text{atan2}(y, x)$ is an improved version of the conventional arctangent function that takes into account the signs of its arguments $y, x$ and is implemented in several programming languages; for more information, see [2.3]. Once we have obtained $\theta^*$, the next step is to locate $R_1$ and $R_2$ within line $\ell_1(\theta^*)$; for this, we simply compute the intersections between line $\ell_1(\theta^*)$ and z-circumference $C_2/\ell_1(\theta^*)$; see figure 2.1. This description constitutes, for the moment, the basic LC method for univariate polynomials of degree 2.

In order to construct z-circumference $C_2/\ell_1(\theta^*)$, we only need to know two points in $\ell_1(\theta^*)$; we know of course fixed point $P_1$, and once we know $\theta^*$, we can obtain a second point $P_2$ in $\ell_1(\theta^*)$ by means of the expression $P_2 = P_1 + t(\cos\theta^* + i\sin\theta^*)$, with $t \in \mathbb{R}\backslash\{0\}$; if we use $t = 1$, we will simplify the calculations a bit. After this, we already know three points in z-circumference $C_2/\ell_1(\theta^*)$: $C_2/P_1$, $C_2/P_2$ and 0 (the latter by the asymptotic properties of Möbius transformation $\ell_1 \to C_2/\ell_1$, as we saw in properties 1.1, 1.1A and 1.2 from chapter 1); with this information, the center and the radius of $C_2/\ell_1(\theta^*)$ are easily obtained (the center can be computed with the R function `center_zcircle` from annex 1 section 1; the radius is computed simply as the modulus of the center), and so we have enough information to find the intersections between line $\ell_1(\theta^*)$ and z-circumference $C_2/\ell_1(\theta^*)$. The algorithmic details for finding such intersections are described in a general way in annex 1 section 2, in which the script for the R function `intersect_semiLine_circle` is also listed; here, we can see that, ultimately, obtaining the intersections between line $\ell_1(\theta^*)$ and z-circumference $C_2/\ell_1(\theta^*)$ implies the use of the general formula (2.4) for solving quadratic equations; as we saw at the beginning of this chapter, formula (2.4) is traditionally used to directly find the roots of a univariate polynomial of degree 2, without complicating matters by using geometric constructions of lines and circumferences (drawn on the complex plane) related to the polynomial coefficients. Later we will see, however, that geometric constructions that take this idea as a basis, computed simultaneously on a machine with parallel processing capabilities for several values of the inclination angle $\theta$ of line $\ell_1$, are useful for producing *proximity maps* (real functions of a real variable where the independent variable is the inclination angle $\theta$ of $\ell_1$, and the dependent variable is a measure of proximity to any root of the associated polynomial equation), whose discretized representations in a computing device, will allow us to obtain reasonable initial approximations to the roots of univariate polynomials of degree $n$, with $n \geq 3$. In principle, obtaining each element of a discrete proximity map by using parallel computing techniques will improve the speed of the process. In chapter 3 we will develop this idea in greater detail.

Looking once again at figure 2.1, we observe that there is a very interesting additional geometric property between line $\ell_1(\theta^*)$ and z-circumference $C_2/\ell_1(\theta^*)$:

**Proposition 2.1**: If line $\ell_1 \subset \mathbb{C}\backslash\{0\}$ contains the roots $R_1, R_2$ of polynomial equation (2.2), and $R_1 \neq R_2$, then $\angle 0P_1R_1 = \angle R_1P_1(C_2/P_1)$, where $R_1$, without loss of generality, represents either of the two roots of equation (2.2).

**Demonstration**: let's consider the following line segments:

$A \equiv \overline{0P_1}$, $B \equiv \overline{P_1R_1}$, $C \equiv \overline{0R_1}$;

$A' \equiv \overline{R_1P_1}$, $B' \equiv \overline{P_1(C_2/P_1)}$, $C' \equiv \overline{R_1(C_2/P_1)}$.

In this way, $\angle 0P_1R_1$ is the angle formed by segments $A$ and $B$, and $\triangle\ 0P_1R_1$ is the triangle formed by segments $A$, $B$ and $C$; likewise, $\angle R_1P_1(C_2/P_1)$ is the angle formed by segments $A'$ and $B'$, and $\triangle\ R_1P_1(C_2/P_1)$ is the triangle formed by segments $A'$, $B'$ and $C'$; as an aid to the visualization of these segments and triangles, refer again to figure 2.1. We will now prove that $|A|/|A'| = |B|/|B'| = |C|/|C'|$, where $|A|$, for example, denotes the length of line segment $A$.

First, we have

$$\frac{|A|}{|A'|}=\frac{|\overline{0P_1}|}{|\overline{R_1P_1}|}=\frac{|(R_1+R_2)/2|}{|(R_1+R_2)/2-R_1|}=\frac{|(R_1+R_2)/2|}{|(R_2-R_1)/2|}=\frac{|R_1+R_2|}{|R_2-R_1|}=\frac{|R_1+R_2|}{|R_1-R_2|}.$$

Then,

$$\frac{|B|}{|B'|}=\frac{|\overline{P_1R_1}|}{|\overline{P_1(C_2/P_1)}|}=\frac{|R_1-(R_1+R_2)/2|}{\left|\frac{R_1R_2}{(R_1+R_2)/2}-(R_1+R_2)/2\right|}=\frac{|(R_1-R_2)/2|}{\left|\frac{2R_1R_2}{(R_1+R_2)}-\frac{(R_1+R_2)^2/2}{(R_1+R_2)}\right|}$$

$$\frac{|B|}{|B'|}=\frac{|(R_1-R_2)/2|}{\left|\frac{4R_1R_2}{2(R_1+R_2)}-\frac{(R_1^2+2R_1R_2+R_2^2)}{2(R_1+R_2)}\right|}=\frac{|2(R_1+R_2)||(R_1-R_2)/2|}{|-R_1^2+2R_1R_2-R_2^2|}=\frac{|R_1+R_2||R_1-R_2|}{|R_1-R_2|^2}$$

$$\frac{|B|}{|B'|}=\frac{|R_1+R_2|}{|R_1-R_2|}.$$

From here, we see that $|A|/|A'|=|B|/|B'|$. Finally, we have

$$\frac{|C|}{|C'|}=\frac{|\overline{0R_1}|}{|\overline{R_1(C_2/P_1)}|}=\frac{|R_1|}{\left|\frac{R_1R_2}{(R_1+R_2)/2}-R_1\right|}=\frac{|R_1|}{\left|\frac{2R_1R_2}{R_1+R_2}-\frac{R_1^2+R_1R_2}{R_1+R_2}\right|}=\frac{|R_1|}{\left|\frac{R_1R_2-R_1^2}{R_1+R_2}\right|}$$

$$\frac{|C|}{|C'|}=\frac{|R_1||R_1+R_2|}{|R_1||R_2-R_1|}=\frac{|R_1+R_2|}{|R_2-R_1|}=\frac{|R_1+R_2|}{|R_1-R_2|}.$$

With this, we have shown that $|A|/|A'|=|B|/|B'|=|C|/|C'|$. This result tells us that the triangles $\triangle\, 0P_1R_1$ and $\triangle\, R_1P_1(C_2/P_1)$ are similar, since their corresponding sides are proportional; therefore, the angle formed by segments $A$ and $B$ is equal to the angle formed by corresponding segments $A'$ and $B'$; or, equivalently, $\angle 0P_1R_1=\angle R_1P_1(C_2/P_1)$. ■

Proposition 2.1 provides an alternative geometric method for finding $\theta^*$; it would just suffice to bisect known angle $\angle 0P_1(C_2/P_1)$, which can be constructed directly from the coefficients of polynomial (2.1); in this way, we would immediately obtain line $\ell_1(\theta^*)$ and consequently z-circumference $C_2/\ell_1(\theta^*)$; as we saw above, the intersections between $\ell_1(\theta^*)$ and $C_2/\ell_1(\theta^*)$ are the geometric places of $R_1$ and $R_2$. It is more practical, from a computational perspective, to find $\theta^*$ by using expression (2.8), although it is also feasible to use the result from proposition 2.1 for the same purpose. At the end of this chapter, we will look at two numerical examples that use expression (2.8) to approximate the roots of quadratic equations in one variable. We will leave to the reader as an exercise to verify that proposition 2.1 holds in each of these numerical examples.

The basic LC method described above works if $0 \notin \ell_1(\theta^*)$, and $R_1 \neq R_2$. If $0 \in \ell_1(\theta^*)$, we can follow this alternative geometric procedure: we start from the fact that we know the inclination angle of derived semi-line $\ell_d$, and therefore we can still obtain $\theta^*$ (the inclination angle of $\ell_1$) in the same way as we did in the general method described above, by means of formula (2.8). Then, we use the direction vector of $\ell_1(\theta^*)$, $v_1 = \cos\theta^* + i\sin\theta^*$, and find values $t_r \in \mathbb{R}$ such that

$$(P_1 + t_r v_1)(P_1 - t_r v_1) = C_2. \tag{2.9}$$

We solve for $t_r$ in equation (2.9):

$$P_1^2 - t_r^2 v_1^2 = C_2$$

$$-t_r^2 v_1^2 = C_2 - P_1^2$$

$$-v_1^2 t_r^2 = W, \text{ according to (2.6)}$$

$$t_r^2 = -W/v_1^2.$$

From equation (2.7), $\arg(-W) = 2\theta^*$, and since $v_1 = \cos\theta^* + i\sin\theta^*$, we have $\arg(v_1^2) = 2\theta^*$; i.e., $-W$ and $v_1^2$ have the same direction (orientation and sense), so $t_r^2$ is a positive real number. Since $|v_1| = 1$, we have $t_r^2 = |-W| = |W| = |C_2 - P_1^2|$; therefore, $t_r = \pm\sqrt{|C_2 - P_1^2|}$.

Finally, we obtain $R_1$ and $R_2$ as two points on $\ell_1$ equidistant from $P_1 = -C_1/2$:

$$R_1 = -C_1/2 - \sqrt{|C_2 - (C_1/2)^2|}(\cos\theta^* + i\sin\theta^*), \tag{2.10}$$

$$R_2 = -C_1/2 + \sqrt{|C_2 - (C_1/2)^2|}(\cos\theta^* + i\sin\theta^*), \tag{2.11}$$

with $\theta^*$ given by formula (2.8).

Note that formulas (2.10) and (2.11) are also valid when $0 \notin \ell_1(\theta^*)$, so they can be used regardless of whether or not it is possible to obtain the center and radius of z-circumference $C_2/\ell_1(\theta^*)$; they also trivially cover the case in which $R_1 = R_2$, because in such circumstances, $t_r = \sqrt{|C_2 - (C_1/2)^2|} = 0$. Note also that formulas (2.10) and (2.11) are valid even if $C_1 = 0$. In fact, **formulas (2.10) and (2.11) can be used to solve any quadratic equation of the form (2.2)**. In later chapters, we will see that to extend the LC method towards univariate polynomials of degree greater than 2, it will be necessary to jointly use both a z-circumference and a variant of expression (2.11), which will allow us to precisely find intersections between this z-circumference and a generalization of semi-line $\ell_d$. For now, in the remainder of this chapter we will continue to use our basic initial assumptions, which will allow us to approximate roots of polynomial equations of the form (2.2) for typical situations where we can construct geometric elements $\ell_1(\theta^*)$, $C_2/\ell_1(\theta^*)$, and $\ell_d(\theta^*)$, without using formulas (2.10) and (2.11); these formulas can optionally be incorporated into a more robust method if one wants to cover all possible scenarios that arise when trying to approximate the roots of equation (2.2).

We will now see how to apply the LC method to two typical cases in which equation (2.2) has two distinct roots $R_1$, $R_2$, contained in a line $\ell_1(\theta^*)$ that does not pass through the origin, so it is possible to compute the center and radius of z-circumference $C_2/\ell_1(\theta^*)$.

## Numerical Example 2.1

We will use the LC method to approximate the roots of polynomial equation

$$z^2 + (1+i)z + (2+2i) = 0. \tag{2.12}$$

Here, $C_1 = (1+i)$, $C_2 = (2+2i)$; so, according to our described method,

the fixed point of $\ell_1(\theta^*)$ is $P_1 = -0.5 - 0.5i$;

the fixed point of $\ell_d$ is $P_1^2 = 0.5i$;

a direction vector for $\ell_d$ is $W = C_2 - P_1^2 = 2 + 1.5i$,

and the direction angle of $\ell_1(\theta^*)$, according to (2.8), is

$\theta^* = \text{atan2}\,(-1.5, -2)/2 = -1.249046\, rad\; = -71°\, 33'\, 54.18''$.

Then, two points on $\ell_1(\theta^*)$ are

$P_1 = -0.5 - 0.5i, \qquad P_2 = P_1 + \cos\theta^* + i\sin\theta^* = -0.183772 - 1.448683i.$

From here, we obtain 2 points on z-circumference $C_2/\ell_1(\theta^*)$:

$C_2/P_1 = -4, \qquad C_2/P_2 = -1.531057 + 1.186342i.$

These two values serve as input for the algorithm described in annex 1 section 1 (function `center_zcircle`), which finds the center $c$ of $C_2/\ell_1(\theta^*)$. By running this algorithm, we find

$c = -2 - i$,

and from here, we immediately obtain the radius $r$ of $C_2/\ell_1(\theta^*)$

$r = |c| = \sqrt{5}$.

We now have all the information we need to find the intersections between line $\ell_1(\theta^*)$ and z-circumference $C_2/\ell_1(\theta^*)$; we give function `intersect_semiLine_circle` from annex 1 section 2, which finds intersections between lines and circumferences, the following four arguments:

1. fixed point $P_1$;

2. unit direction vector of $\ell_1(\theta^*)$, given by $V(\theta^*) = \cos\theta^* + i\sin\theta^* = 0.3162278 - 0.9486833i$;

3. center $c$ and

4. radius $r$ of circumference $C_2/\ell_1(\theta^*)$.

By doing this, we obtain

$$z_1 = -1 + i \qquad z_2 = -2i. \tag{2.13}$$

These are the approximations that the LC method obtains for the roots of equation (2.12). In fact, it can be manually verified that the values in (2.13) satisfy the polynomial equation (2.12).

For polynomials of degree 2 such as those in equation (2.2), an implementation of the LC method on a personal computer usually produces approximations to their true roots with an absolute error of magnitude similar to that of the machine epsilon; this means that the approximation errors in this case are mostly attributable to the internal floating-point representations of such approximations, as well as to possible errors induced by floating-point arithmetic implemented in personal computer processors. Later we will see, informally, that the approximation errors associated to the LC method for univariate polynomials of degree greater than 2 depend, to a greater extent, on the degree of discretization of the generated proximity maps, and to a much lesser extent (so much so that they can be regarded as insignificant in comparison with discretization errors), on the errors induced by floating-point internal representations and floating-point arithmetic. In this numerical example, the absolute errors for approximations $z_1$ and $z_2$ in (2.13) with respect to the true roots are, respectively, $4.965068 \times 10^{-16}$ and $0$; the machine epsilon for a personal computer with a 64-bit multicore processor is $2^{-52} \approx 2.220446 \times 10^{-16}$, so we see that in this case one of the approximation errors has an order of magnitude similar to that of the machine epsilon, while the other has a magnitude smaller than that of the machine epsilon.

Figure 2.3 graphically shows the approximations obtained in this example; these are the intersections between trajectories $\ell_1(\theta^*) = P_1 + tV(\theta^*)$ and $C_2/\ell_1(\theta^*)$; $t \in \mathbb{R}$. The reader can verify the validity of proposition 2.1 with the results obtained in this numerical example.

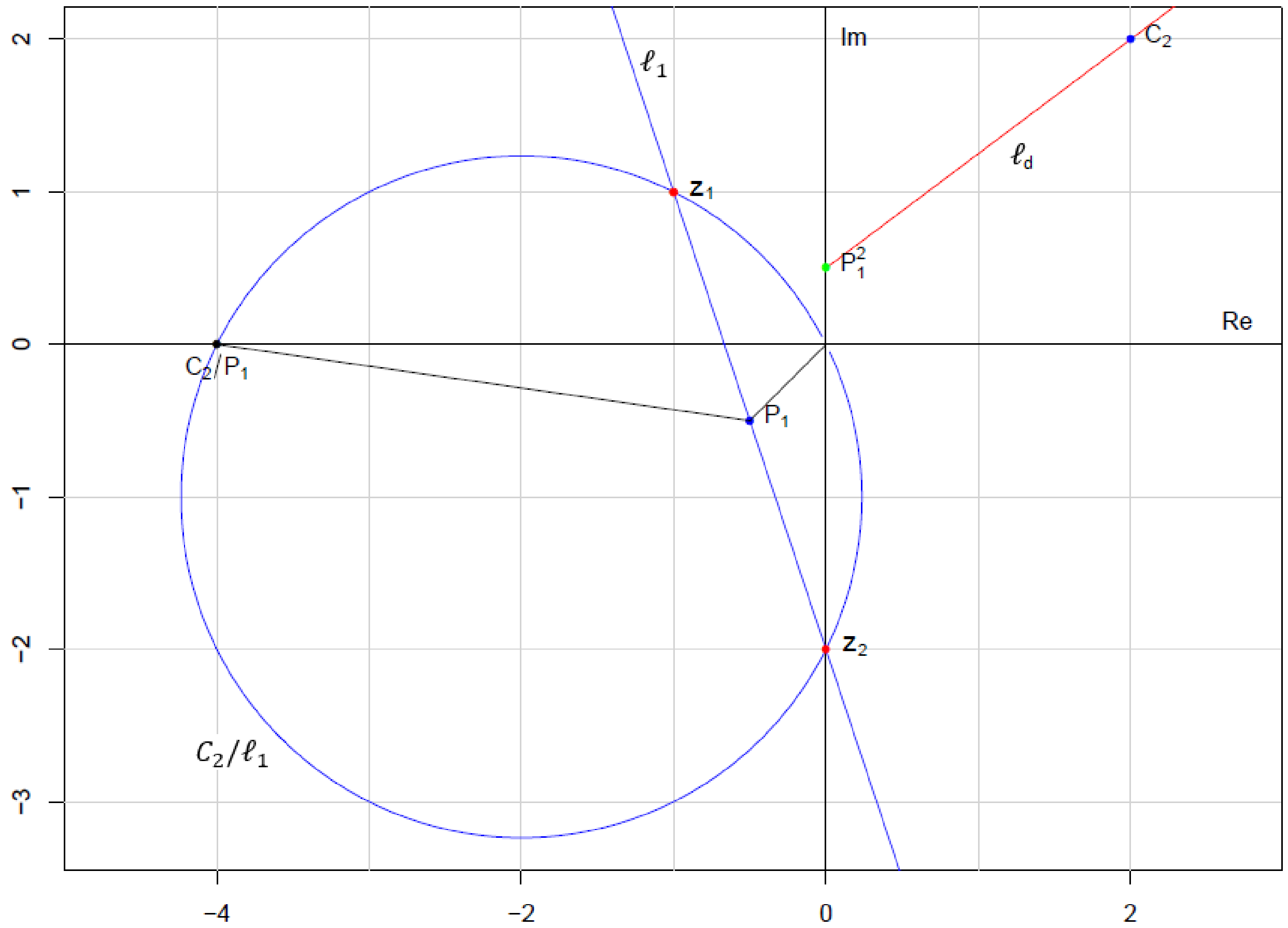

**Fig. 2.3**. Approximations $z_1, z_2$ to the roots $R_1, R_2$ of polynomial equation (2.12) obtained by calculating the intersections between line $\ell_1(\theta^*)$ and z-circumference $C_2/\ell_1(\theta^*)$ constructed with the LC method. Derived semi-line $\ell_d$ is also displayed.

## Numerical Example 2.2

In this example we generate two “random” complex numbers $R_1, R_2$ in the R programming language, by fixing R’s pseudorandom number generator seed, so we can reproduce the results. For more details, see annex 3 section 1, where you can find the R script used to produce the results of this numerical example. Additionally, since this script integrates functions `center_zcircle` and `intersect_semiLine_circle` from annex 1, sections 1 and 2, it is feasible to use it, with minor modifications, to generate the results described in numerical example 2.1.

The "randomly" generated complex numbers are:

$$R_1 = 0.5398031 + 0.4256795i, \qquad R_2 = -0.3932796 + 0.2364727i. \qquad (2.14)$$

The idea of this example is to construct, from generated roots (2.14), coefficients $C_1, C_2$ like those of polynomial (2.1), by using Vieta’s relations (2.3-a), (2.3-b). Then, we give as input to the LC method the constructed coefficients $C_1, C_2$, generating in this way two approximations $z_1$, $z_2$; note that the LC method does not have any knowledge of generated roots (2.14).

The coefficients calculated from generated values (2.14) are:

$$C_1 = -0.1465235 - 0.6621522i, \quad C_2 = -0.3129551 - 0.0397623i. \quad (2.15)$$

Below, we show some intermediate and final results from the LC method, after processing input coefficients (2.15):

- Fixed point $P_1^2$ and direction vector $W$ for derived semi-line $\ell_d$:

  $P_1^2 = -0.1042441 + 0.0485104i, \;\; W = C_2 - P_1^2 = -0.208711 - 0.0882728i.$

- Fixed point $P_1$, inclination angle $\theta^*$, direction vector $V(\theta^*)$ and a second point $P_2 = -C_1/2 + V(\theta^*)$ for line $\ell_1(\theta^*)$:

  $P_1 = 0.0732617 + 0.3310761i, \quad \theta^* = 0.2000633\, rad = 11°\, 27'\, 46.02'',$
  $V(\theta^*) = 0.980054 + 0.1987314i, \quad P_2 = 1.053316 + 0.529807i.$

- Fixed point $C_2/P_1$, center $c$, radius $r$, and a second point $C_2/P_2$ for z-circumference $C_2/\ell_1(\theta^*)$:

  $C_2/P_1 = -0.3139017 + 0.8758051i, \quad c = 0.0374698 + 0.5075859i,$
  $r = 0.508967, \quad C_2/P_2 = -0.2522763 + 0.0891428i.$

- Intersections between line $\ell_1(\theta^*)$ and z-circumference $C_2/\ell_1(\theta^*)$:

  $$z_1 = -0.3932796 + 0.2364727i \quad z_2 = 0.5398031 + 0.4256795i. \quad (2.16)$$

Values (2.16) are the approximations that our R language implementation of the LC method obtains when it receives as input the coefficients (2.15). By calculating the absolute errors of these approximations with respect to the corresponding generated roots (2.14), we obtain

$|R_1 - z_2| = 0, |R_2 - z_1| = 0.$

In this numerical example, we see that our implementation of the LC method achieved approximations to the roots each with an absolute error no greater than the machine epsilon. Figure 2.4 shows the geometric construction LC associated with input coefficients (2.15) in this numerical example. Again, the reader can verify the validity of proposition 2.1 by using the results obtained here.

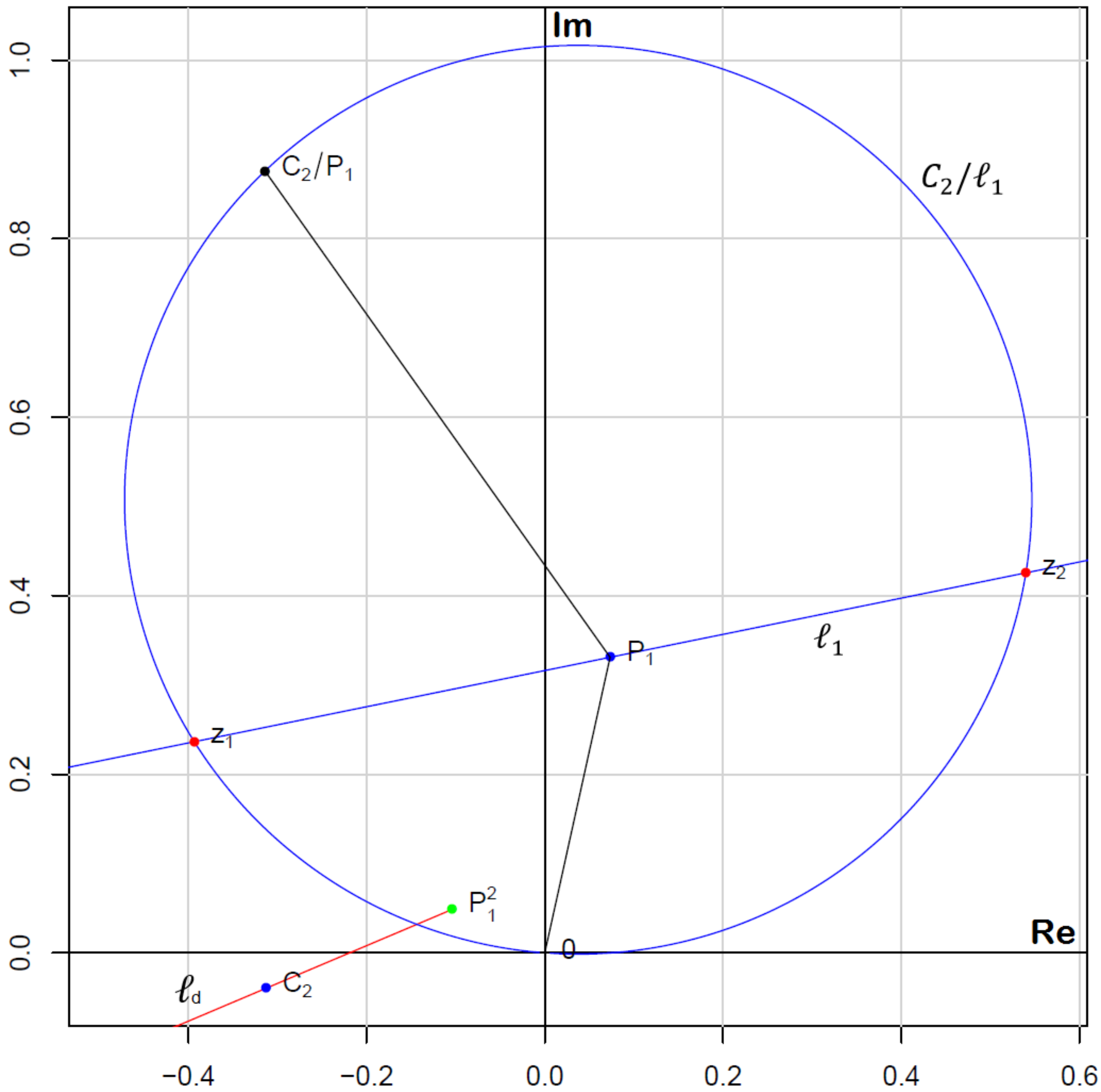

**Fig. 2.4**. Approximations $z_1, z_2$ to the roots of equation $z^2 + C_1 z + C_2 = 0$, with coefficients $C_1$, $C_2$ given by (2.15). $z_1, z_2$ were obtained by computing the intersections between line $\ell_1(\theta^*)$ and z-circumference $C_2/\ell_1(\theta^*)$ constructed by providing the LC algorithm with coefficients (2.15). Derived semi-line $\ell_d$ is also displayed.

# Chapter Three: LC Method for Cubic Polynomial Equations

## Theoretical Aspects

In this chapter we will extend the LC method to univariate polynomials of degree 3, trying to preserve the same notation of chapter 2.

The polynomial equation of degree 3, whose roots $R_1$, $R_2$, $R_3$ we seek to approximate is of the form

$$z^3 + C_1 z^2 + C_2 z + C_3 = 0 \tag{3.1}$$

In this case, Vieta's relations between the roots and the coefficients of equation (3.1) establish that

$$C_1 = -(R_1 + R_2 + R_3) \tag{3.2-a}$$

$$C_2 = R_1 R_2 + R_1 R_3 + R_2 R_3 \tag{3.2-b}$$

$$C_3 = -R_1 R_2 R_3 \tag{3.2-c}$$

Unlike the case of quadratic equations, if $R_1$, $R_2$, $R_3$ are small integers, this time it is no longer easy to mentally find numerical values for $R_1$, $R_2$, $R_3$ from Vieta's relations shown in equations (3.2-a), (3.2-b), (3.2-c). There are techniques based on the remainder theorem, the factor theorem, and the synthetic division to evaluate, by trial and error, root candidates in case that these are suspected to be integer or rational numbers (see [3.1]). On the other hand, there is an analytical formula for finding the roots of cubic equation (3.1), accredited to Gerolamo Cardano (in [3.2] you can consult the derivation of this formula). This formula is much more elaborate with respect to the general formula (2.4) for solving quadratic equations.

In a similar way to how we proceeded in chapter 2, for the case of cubic equations we will also seek to take advantage of existing geometric relationships between the coefficients of equation (3.1) and lines and circumferences in the complex plane $\mathbb{C}$, with the purpose of obtaining initial approximations to the solutions of the non-linear system of three equations (3.2-a), (3.2-b), (3.2-c) in three unknowns $R_1$, $R_2$, $R_3$. For the more general case, we will assume that the roots $R_1, R_2, R_3 \in \mathbb{C}\backslash\{0\}$ are distinct from each other; i.e., we will assume that we have no roots with multiplicity greater than 1: this will be our framework on which we will develop our LC methodology from now on.

Later in this chapter, we will briefly analyze, through numerical examples, what happens with our general LC methodology in the presence of repeated roots in univariate polynomials of degree 3; we will not delve into this analysis, since our main objective is to build a methodology that in general works well for univariate polynomials with roots of multiplicity 1, randomly distributed in the complex plane $\mathbb{C}$.

To begin with, let's consider a geometric construction, consisting of lines and a circumference in the complex plane $\mathbb{C}$, which contain fixed points directly related to coefficients $C_1$, $C_2$, $C_3$ in equation (3.1); without loss of generality, one of the lines in this geometric construction contains the root $R_1$ of equation (3.1); this geometric construction is illustrated, in a schematic way, in figure 3.1.

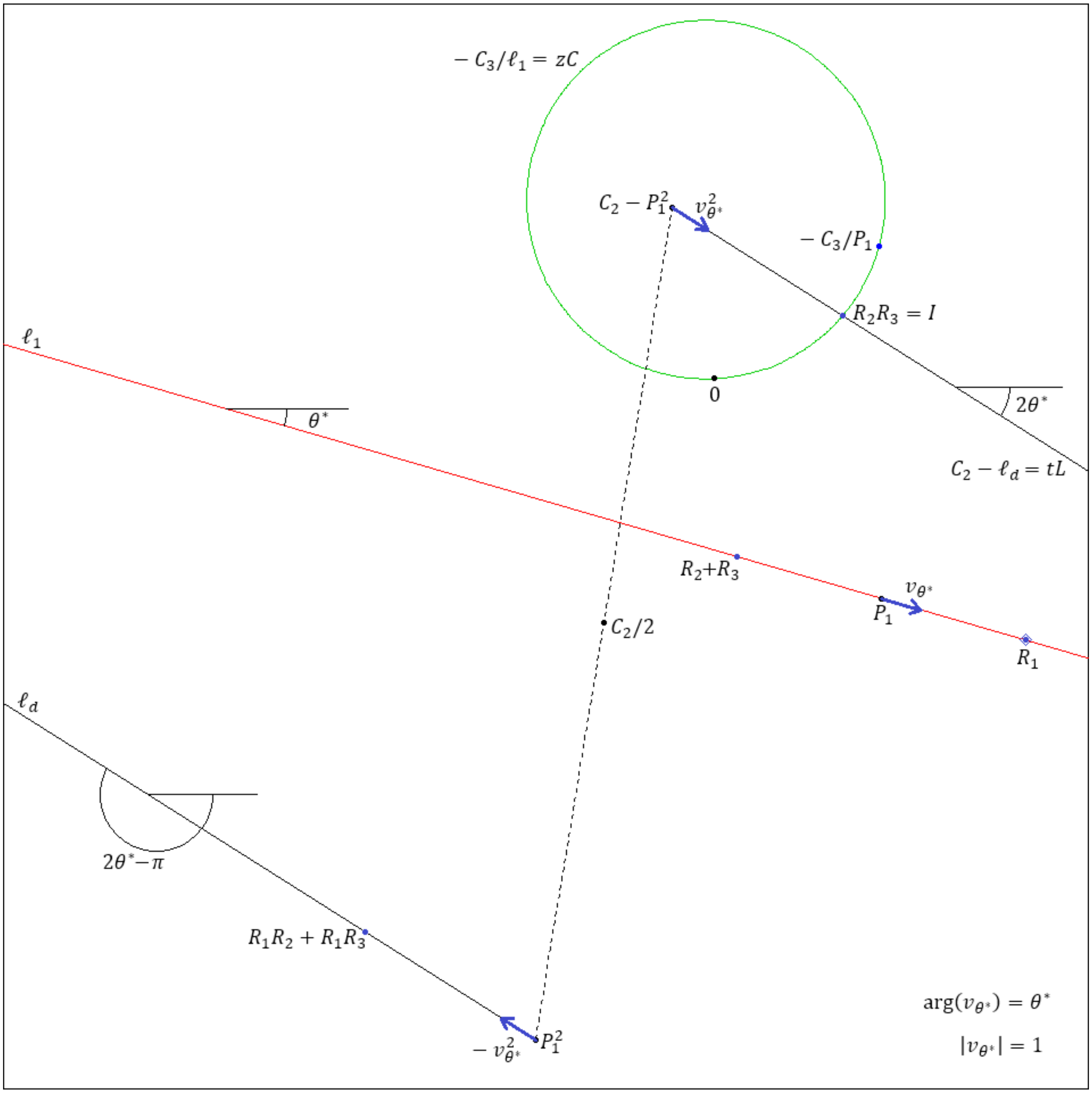

**Fig. 3.1**. If line $\ell_1 \subset \mathbb{C}\backslash\{0\}$ with fixed point $P_1 = -C_1/2$ contains the root $R_1$ of polynomial equation (3.1), then z-circumference $zC = -C_3/\ell_1$ contains the product $R_2R_3$; this product is located at the intersection $I$ between $zC$ and *terminal semi-line* $tL = C_2 - \ell_d$, where $\ell_d$ is the semi-line derived from $\ell_1$.

In the geometric construction of figure 3.1, we assume that $\ell_1$ contains both the root $R_1$ of equation (3.1), and fixed point $P_1 = -C_1/2$; z-circumference $zC$ is obtained by means of the Möbius transformation $zC = -C_3/\ell_1$, so both $R_2R_3$ and fixed point $-C_3/P_1$ are in $zC$. $v_{\theta^*}$ is a unit vector pointing from $P_1$ to $R_1$; its argument is $\theta^*$. We will say that $v_{\theta^*}$ is the direction vector of $\ell_1$, and that $P_1$ is the fixed point of $\ell_1$; observe that $R_2 + R_3$ is also in $\ell_1$ (according to property 1.3 from chapter 1).

Derived semi-line $\ell_d$, shown near the bottom-left corner of figure 3.1, and generated (as in chapter 2) from parametric expression $\ell_d$: $(P_1 - tv_{\theta^*})(P_1 + tv_{\theta^*})$, with $t \in \mathbb{R}$, contains points $P_1^2$ and $R_1R_2 + R_1R_3$; its direction vector is $-v_{\theta^*}^2$, with argument $2\theta^* - \pi$ (or equivalently, with argument $2\theta^* + \pi$). *The reflection of $\ell_d$ with respect to $C_2/2$*, i.e., semi-line $C_2 - \ell_d$, contains points $C_2 - P_1^2$ and $R_2R_3$; its direction vector is $v_{\theta^*}^2$, with argument $2\theta^*$. We will call $C_2 - \ell_d$ *the terminal semi-line $tL$*. Note that $tL$ contains point $I = R_2R_3$, as does z-circumference $zC$.

In practice, we will concentrate on finding the geometric intersection between $tL$ and $zC$.

We will denote the geometric construction of figure 3.1 as $LzC(C_1, C_2, C_3, \theta^*)$, to indicate that it depends on the coefficients in equation (3.1) as well as on the argument of $\ell_1$'s direction vector, $v_{\theta^*}$. Unfortunately, we cannot construct $LzC(C_1, C_2, C_3, \theta^*)$ directly using only the coefficients from equation (3.1), since with this information we would only know one point from $\ell_1$ ($P_1$), one from $tL$ ($C_2 - P_1^2$), and two from $zC$ (0 and $-C_3/P_1$); that is to say, if we use only the coefficients from equation (3.1), then we do not have enough information to obtain in a single step the numerical value of angle $\theta^*$ that leads us to $R_1$, just as we did for the case of univariate polynomials of degree 2. Note that, in general, for the case of polynomials of degree 3 in a complex variable, we must determine three different angles $\theta^*$, which can be distributed anywhere in the interval $[-\pi, \pi)$. Also note that, in a construction $LzC(C_1, C_2, C_3, \theta^*)$, $\ell_1$ and $zC$ may or may not intersect each other; in any case, this is irrelevant, since our interest now focuses on the intersection between $tL$ and $zC$.

From the previous paragraph, we see that we require additional information if we want to approximate the numerical value for any of the three angles $\theta^*$ in $LzC$ which leads us to any of the roots of equation (3.1). We will not attempt to discover geometric relationships among various elements from the 3 possible constructions $LzC(C_1, C_2, C_3, \theta^*)$ taken together, in an attempt to extend proposition 2.1 from chapter 2 towards polynomials of degree 3; instead, we will take a value arbitrarily close to angle $\theta^*$ from figure 3.1, let's say $\theta = \theta^* + \Delta\theta$, and construct a structure $LzC(C_1, C_2, C_3, \theta)$, using of course the same coefficients $C_1, C_2, C_3$ used in figure 3.1; such structure is shown schematically in figure 3.2:

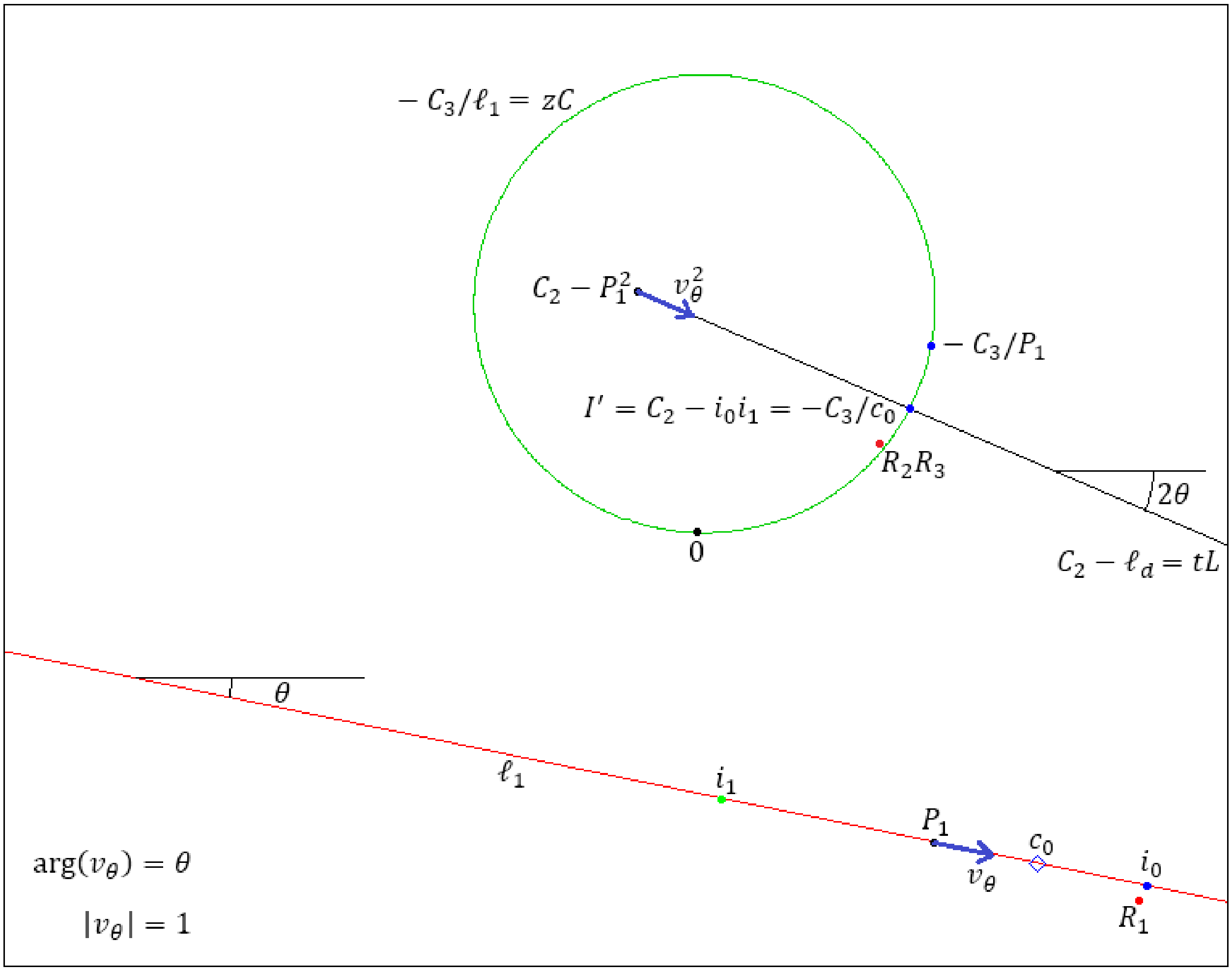

**Fig. 3.2.** Structure $LzC(C_1, C_2, C_3, \theta)$, in which $\theta$ is "close" to angle $\theta^*$ from figure 3.1, so that $\ell_1$ is "close" to $R_1$. The projection of $I'$ (the intersection between $tL$ and $zC$, which tends to $R_2R_3$ as $\theta$ tends to $\theta^*$) onto $\ell_1$ produces a pair of points, $i_0$ and $c_0$, which help us quantify in a more precise way the degree of proximity or "closeness" between $\ell_1$ and $R_1$, without knowing a priori where $R_1$ is located.

In the structure $LzC(C_1, C_2, C_3, \theta)$ from figure 3.2, $v_\theta$ is $\ell_1$'s direction vector. Point $I'$ is at the intersection between $tL$ and $zC$, and can be obtained in two different ways:

1. Via $tL$: from $C_2$, subtract product $i_0 i_1$; $i_0$, $i_1$ are both in $\ell_1$, $i_0 \neq i_1$ and $|i_0 - P_1| = |i_1 - P_1|$.

2. Via $zC$: divide $-C_3$ by $c_0$ ($c_0$ is also in $\ell_1$).

Now, the distance between $i_0$ and $c_0$ gives us information on how far, angularly, is $\ell_1$ from *one* of the roots of equation (3.1). We define the quantity

$$e(C_1, C_2, C_3, \theta) \coloneqq sgn\left(\frac{i_0 - c_0}{v_\theta}\right)\frac{|i_0 - c_0|}{|i_0 - P_1| + |c_0 - P_1|}, \tag{3.3}$$

where $i_0 = P_1 + \sqrt{|I' - (C_2 - P_1^2)|}v_\theta$, $c_0 = -C_3/I'$, and $sgn(x) = \begin{cases} -1 & \text{if } x < 0 \\ 0 & \text{if } x = 0. \\ 1 & \text{if } x > 0 \end{cases}$

The quantity $e(C_1, C_2, C_3, \theta)$ defined in expression (3.3) is a *weighted error* associated to structure $LzC(C_1, C_2, C_3, \theta)$. For simplicity's sake, we will denote this quantity as $e(\theta)$.

Note that $-1 < e(\theta) \leq 1$. Indeed, if $c_0$ were on the side of $i_1 = P_1 - \sqrt{|I' - (C_2 - P_1^2)|}v_\theta$ with respect to fixed point $P_1$, or $c_0 = P_1$, then $e(\theta) = 1$, which indicates that $\ell_1$ is very far from a root of equation (3.1); on the other hand, if $c_0$ were between $P_1$ and $i_0$, then $0 < e(\theta) < 1$; if $c_0$ were on the side of $i_0$ with respect fixed point $P_1$, but $|c_0 - P_1| > |i_0 - P_1|$, then $-1 < e(\theta) < 0$. Finally, if $c_0 = i_0$, then $e(\theta) = 0$, and in this case we may have found one of the angles $\theta^*$ that leads us to one of the roots of equation (3.1); in fact, in this latter situation, both $c_0$ and $i_0$ should coincide with $R_1$.

To see where the scalar parameter $\sqrt{|I' - (C_2 - P_1^2)|}$ (used in expression (3.3) to compute point $i_0$) comes from, just solve for parameter $t$ in the equation that is obtained by equating the parametric expression for terminal semi-line $tL$ with the intersection $I'$ between $tL$ and $zC$:

$$\begin{aligned} C_2 - (P_1 - t_{I'}v_\theta)(P_1 + t_{I'}v_\theta) &= I' \\ -(P_1 - t_{I'}v_\theta)(P_1 + t_{I'}v_\theta) &= I' - C_2 \\ -P_1^2 + t_{I'}^2 v_\theta^2 &= I' - C_2 \\ t_{I'}^2 v_\theta^2 &= I' - C_2 + P_1^2 \\ t_{I'}^2 &= [I' - (C_2 - P_1^2)]/v_\theta^2 \end{aligned}$$

Since $I' - (C_2 - P_1^2)$ and $v_\theta^2$ have the same direction and $v_\theta^2$ is a unit vector, we finally arrive at the expression $t_{I'} = \sqrt{|I' - (C_2 - P_1^2)|}$; this procedure for solving parameter $t$ is similar to the procedure followed in chapter 2 to obtain expression (2.11).

Also note that the function $e(\theta)$ is periodic, because $e(\theta) = e(\theta \pm 2\pi k), k \in \mathbb{Z}$; later, when we construct graphs for $e(\theta)$, we will restrict parameter $\theta$ (the inclination angle of $\ell_1$) to interval $[-\pi, \pi)$.

By integrating the concepts and elements we have so far, it is possible to design a strategy that allows us to find reasonable initial approximations to the roots of equation (3.1), taking advantage of the parallel processing capabilities of current computing devices:

**Strategy 3.1 to find initial approximations to the roots of a univariate polynomial of degree 3 with complex coefficients by means of parallel constructions $LzC$**

a) Obtain parallel values $e(\theta_k)$ associated with structures $LzC(C_1, C_2, C_3, \theta_k)$, where $\theta_k = -\pi + 2\pi k/N$, $k = 0,1,2, \ldots, N-1$; note that points $\theta_k$ constitute a regular partition of interval $[-\pi, \pi)$.

b) Join, by using rectilinear segments, adjacent points $(\theta_k, e(\theta_k)), (\theta_{k+1}, e(\theta_{k+1})) \in \mathbb{R}^2$ $\forall k$; in this way, we'll have constructed, by linear interpolation, a **discrete approximation** $\hat{e}(\theta)$ for function $e(\theta)$ at interval $-\pi \leq \theta < \pi$.

c) Find crossings of $\hat{e}(\theta)$ with horizontal axis $y = 0$; with this, we will obtain estimates $\hat{\theta}_i^*$ of direction angles $\theta_i^*$ ($i \in \{1,2,3\}$) that, for each $i$, cause line $\ell_1$ to contain one of the roots $R_i$ of equation (3.1).

d) Generate structures $LzC(\hat{\theta}_i^*)$, with which it is possible to obtain approximations to $R_i$, by means of expressions $\hat{R}_i = \left[i_0\left(\hat{\theta}_i^*\right) + c_0(\hat{\theta}_i^*)\right]/2$, where $i_0\left(\hat{\theta}_i^*\right)$ and $c_0(\hat{\theta}_i^*)$ denote points $i_0$ and $c_0$ on lines $\ell_1$ associated with structures $LzC(\hat{\theta}_i^*)$. Note that, for each $i$, $\hat{R}_i$ is also located on the line $\ell_1$ associated with structure $LzC(\hat{\theta}_i^*)$, by property 1.3a from chapter 1.

**Note 3.1**: In item a) of strategy 3.1, theoretically it is possible to distribute $N$ geometric constructions $LzC(\theta_k)$, in order to obtain values $e(\theta_k)$, among $M$ processors of a machine with parallel processing power; typical values of $N$ could be 1,000, 10,000, or even 100,000. In principle, it seems reasonable to assume, if for the moment we do not take into account Amdahl's law (see a description of this law in [3.3]), that $M$ should be equal to $N$ if we want to optimize execution time; a large-scale computer, or a network of interconnected computers, could trivially handle this situation, independently processing a structure $LzC(\theta_k)$ in each of its processing cores. Items c) and d) of strategy 3.1 can also be implemented by taking advantage of the parallel processing capabilities of current computing devices, although for item d), of course, a much smaller number of parallel processors would be required: only as many as estimates $\hat{\theta}_i^*$ are obtained in item c).

In practice, most people have access to personal computers for which $M \ll N$, although there is the possibility of access through the cloud from a personal computer to High Performance Computing (HPC) cluster services, such as those offered by Amazon Web Services (AWS), to increase the number of available processing cores. Another possibility that allows to increase the processing power of a personal computer is to equip it with a video card (two popular choices are NVIDIA GeForce and AMD Radeon), which contains physical arrays with hundreds or thousands of processing cores, which of course can be used for parallel computing tasks, in addition to their main function of speeding up the processing of video signals in a computer. The R statistical software, which also doubles as a programming language, has some advanced parallel processing capabilities, using either a personal computer with 4, 6, 8 or even 12 processing cores in its main CPU, or a personal computer with integrated video card, or a computer with access to HPC clusters. Readers interested in taking advantage of advanced parallel processing capabilities of the R language can try out one or more of the resources described in [3.4]. Of course, there are other programming languages (such as C, C++, Java or Fortran) which are much more efficient than R for programming numerical methods that take advantage of the parallel processing capabilities of modern hardware.

For our illustrative purposes, and to build numerical examples in this chapter and in later chapters that can be reproduced on a typical personal computer, it will be enough to use functions that take advantage of the *vectorization of operations*, as well as the `apply` function, both available in the basic version of R. In the vectorization of operations, an operation (such as addition) acts upon complete arrays of numerical elements (such as vectors or matrices), processing at once all the elements of the arrays involved in that operation. The vectorization of operations is an extension of the functionality of operators acting on individual scalars towards arrays of numerical elements; it can be seen as a rudimentary form of parallel processing. The `apply` function, on the other hand, simultaneously applies a function (predefined in R or defined by the user) to each of the elements in a row or column of a matrix, thus avoiding the explicit use of loops (control structures such as `for`, `while` or `repeat`). For more details about the concept of vectorization of operations, as well as the `apply` function in R, see [3.5].

Taking into account these basic parallel processing tools, together with notes 3.2 and 3.3 below, we designed an R script that implements strategy 3.1 described above; this script is listed in annex 3 section 2, and helps us to construct and reproduce the results from the numerical examples in this chapter; functions described in annexes 1 and 2 are used here: function `LzC3`, listed in annex 2 section 2, plays a key role in the implementation of items a) and d) of strategy 3.1, while function `get_angle_roots`, listed in annex 1 section 4, is useful in the implementation of items b) and c) of strategy 3.1. Both functions are called within function `approxP3`, which is listed in annex 2 section 5.

**Note 3.2**: A very important fact that we cannot overlook is the following: there are different scenarios when we look for intersections between a semi-line $tL$ and a z-circumference $zC$; these elements can intersect each other at 0, 1, or 2 points, as shown in figure 3.3.

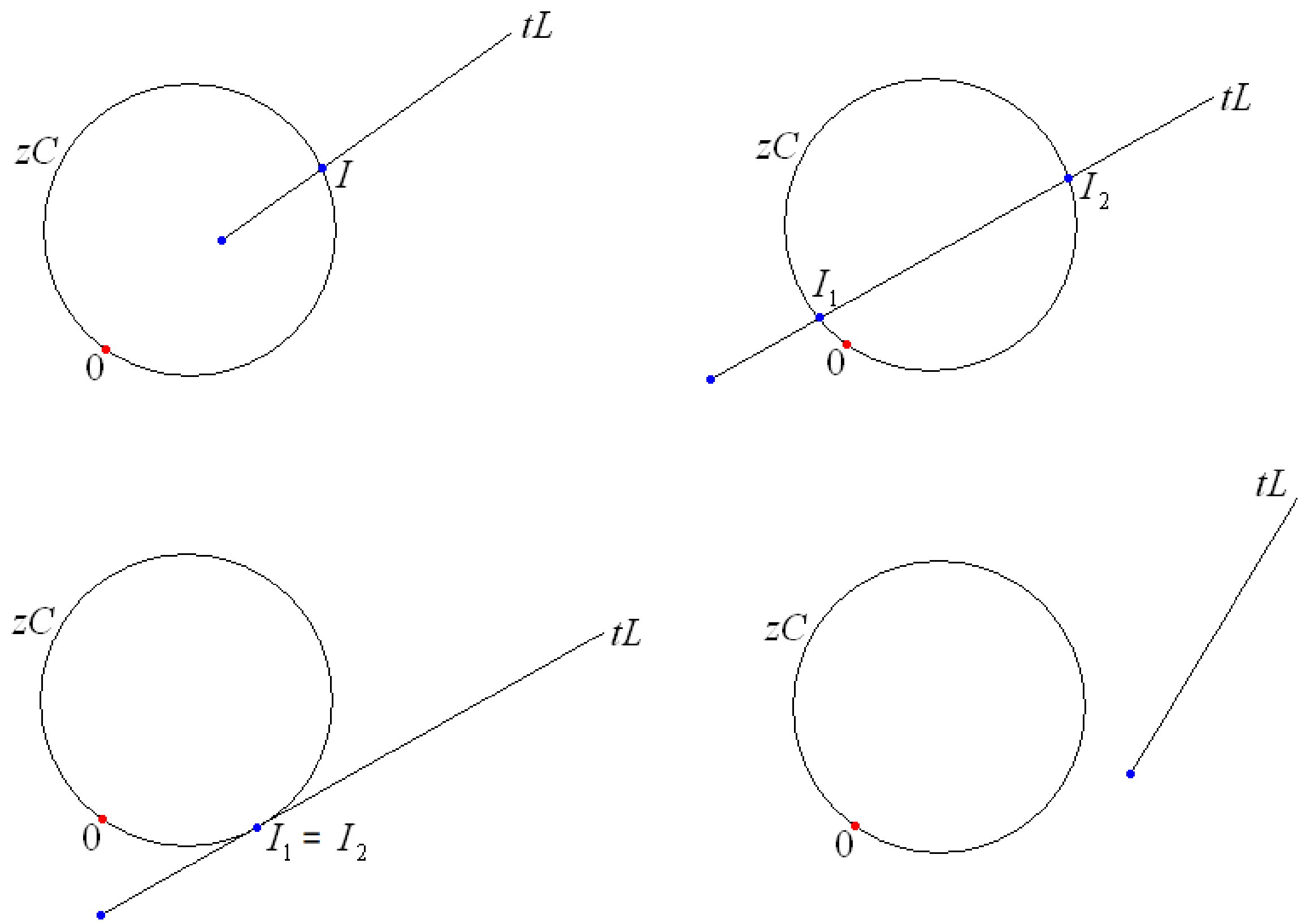

**Fig. 3.3.** Some possible scenarios for the intersections between a semi-line $tL$ and a z-circumference $zC$. The case where there are two intersections $I_1$, $I_2$ (upper right circle) can degenerate into a single tangential intersection when $I_1 = I_2$ (lower left circle). There may also be a single non-tangential intersection $I$ (upper left circle), or there may be no intersections at all (lower right circle).

This implies that there is a need to broaden our definition (3.3) of weighted error:

$$e_A(\theta) := \begin{cases} sgn\left(\frac{i_{0I_1} - c_{0I_1}}{v_\theta}\right) \frac{|i_{0I_1} - c_{0I_1}|}{|i_{0I_1} - P_1| + |c_{0I_1} - P_1|} & \text{if } \exists\, I_1 \\ \text{undefined} & \text{if } \nexists\, I_1 \end{cases} \tag{3.4a}$$

$$e_B(\theta) := \begin{cases} sgn\left(\frac{i_{0I_2} - c_{0I_2}}{v_\theta}\right) \frac{|i_{0I_2} - c_{0I_2}|}{|i_{0I_2} - P_1| + |c_{0I_2} - P_1|} & \text{if } \exists\, I_2 \\ \text{undefined} & \text{if } \nexists\, I_2 \end{cases} \tag{3.4b}$$

$\theta \in [-\pi, \pi)$.

In this extended definition (3.4a), (3.4b), $I_m$, $m = 1{,}2$, are the intersections (if any) between $tL$ and $zC$; $i_{0I_m}, c_{0I_m}$ are points on $\ell_1(\theta)$ associated with intersection $I_m$, analogous to points $i_0, c_0$ from definition (3.3). As we already know, $v_\theta$ and $P_1$ are, respectively, unit direction vector and fixed point of $\ell_1(\theta)$. Thus, we now have two weighted error functions $e_A(\theta)$ and $e_B(\theta)$, each of which is a real function of a real variable; as we vary parameter $\theta$ within structure $LzC$, we must decide who is $e_A$ and who is $e_B$, according to weighted errors obtained in a neighborhood of $\theta$; we must make sure that functions $e_A$ and $e_B$ are continuous at the regions of interval $[-\pi, \pi)$ where they are defined, if we want them to be useful when finding initial approximations to the roots of equation (3.1); fortunately, if in the algorithm used to find the intersections between $tL$ and $zC$, we always designate as $I_1$ the intersection associated to the solution with a negative sign in the quadratic equation involved, and as $I_2$ the intersection associated to the solution with a positive sign (see annex 1 section 2, function `intersect_semiLine_circle`), then $\hat{e}_A$ and $\hat{e}_B$ will be smooth and continuous enough to be used successfully in practice. As we will see in the following numerical examples, this criterion to designate $I_1$ and $I_2$ turns out to be useful in approximating, in a simple way, the *theta roots* $\theta^*$ of *proximity maps* $e_A(\theta)$ and $e_B(\theta)$; these approximations, in turn, help us estimate the true roots of equation (3.1).

**Note 3.3**: The introduction of extended definition (3.4a), (3.4b), implies that we must also refine item d) of strategy 3.1 to find initial approximations to the roots of a univariate polynomial of degree 3 with complex coefficients. To be more specific, when we have already obtained an initial approximation $\hat{\theta}_i^*$, $i \in \{1,2,3\}$ to an angular root, it is perfectly possible that our structure $LzC(\hat{\theta}_i^*)$ produces two estimates for $R_i$, let's say $\hat{R}_{iA}$ and $\hat{R}_{iB}$, obtained by the following expressions:

$$\hat{R}_{iA} = \left[i_{0I_1}\left(\hat{\theta}_i^*\right) + c_{0I_1}\left(\hat{\theta}_i^*\right)\right]/2, \qquad \hat{R}_{iB} = \left[i_{0I_2}\left(\hat{\theta}_i^*\right) + c_{0I_2}\left(\hat{\theta}_i^*\right)\right]/2.$$

In such circumstances, we keep the estimate with lowest associated absolute weighted error, i.e.,

$$\hat{R}_i = \begin{cases} \hat{R}_{iA} & \text{if } \left|e_A\left(\hat{\theta}_i^*\right)\right| < \left|e_B\left(\hat{\theta}_i^*\right)\right| \\ \hat{R}_{iB} & \text{if } \left|e_A\left(\hat{\theta}_i^*\right)\right| \geq \left|e_B\left(\hat{\theta}_i^*\right)\right| \end{cases} \qquad (3.5)$$

**Notation:** In the following numerical examples, the same subscript $i$ will be used to refer to roots $R_i$, $\theta_i^*$, and to root estimates $\hat{R}_i$, $\hat{\theta}_i^*$; the order in which values are assigned (by our R scripts) to root subscripts $i$ may differ from the order in which values are assigned to estimate subscripts $i$, so for example $\hat{R}_2$ does not necessarily represent the estimate of $R_2$. When we compare an individual estimate against its corresponding root, we will use the subscript values originally assigned to both elements.

## Numerical Example 3.1

In a similar way to how we proceeded in example 2.2 from chapter 2, in this numerical example we'll first generate three random complex numbers $R_1$, $R_2$, $R_3$, by fixing the pseudorandom number generator seed in our script, to be able to reproduce our results; later, with these randomly generated roots, we'll obtain coefficients $C_1$, $C_2$, $C_3$ by means of Vieta's relations (3.2-a), (3.2-b) and (3.2-c); such coefficients will serve as input to our procedure for constructing discrete proximity maps (described above in strategy 3.1), which of course does not know the "exact" values $R_1$, $R_2$, $R_3$; finally, the computed estimates $\hat{R}_1$, $\hat{R}_2$, $\hat{R}_3$ will be compared against the randomly generated reference roots, in an attempt to measure and interpret estimation errors. To analyze the development of computations for this numerical example in more detail, please refer to annex 3 section 2, which lists the R script that implements strategy 3.1 described above; this script can be used to reproduce the results of all numerical examples in this chapter and can also be used to generate and explore many other proximity maps associated with cubic polynomials in a complex variable.

Random complex numbers $R_1$, $R_2$, $R_3$, generated by fixing the random number generator seed, by means of the instruction `set.seed(2020)` at the beginning of the script in annex 3 section 2, together with their corresponding "theta roots" $\theta_i^*$, are listed next:

$$
\begin{aligned}
R_1 &= \phantom{-}0.2938057 - 0.2115485i & \theta_1^* &= \phantom{-}0.7283992 \\
R_2 &= \phantom{-}0.2370036 - 0.0462177i & \theta_2^* &= \phantom{-}0.9936450 \\
R_3 &= -0.7278056 - 0.8652312i & \theta_3^* &= -2.6919247
\end{aligned}
\tag{3.6}
$$

Remember that $\theta_i^*$ is the direction of line $\ell_1$, which contains fixed point $P_1 = -C_1/2$ and root $R_i$. The coefficients $C_1$, $C_2$, $C_3$ obtained by substituting reference roots $R_1$, $R_2$, $R_3$ (3.6) into Vieta's relations formulas (3.2-a), (3.2-b) and (3.2-c) are:

$$
\begin{aligned}
C_1 &= \phantom{-}0.19699632 + 1.1229974i, \\
C_2 &= -0.54949766 - 0.3353859i, \\
C_3 &= \phantom{-}0.09869309 + 0.0054156i.
\end{aligned}
\tag{3.7}
$$

The resulting proximity map generated from coefficients $C_1$, $C_2$, $C_3$ in (3.7) is shown in figure 3.4. For the construction of this discrete map, we used $N = 1{,}000$ elements $LzC(\theta_k)$ $(k = 0,1,2,\ldots,N-1)$, with $\theta_k = -\pi + 2\pi k/N$. The "smooth" crossings of this discrete map with horizontal axis $y = 0$ yield the estimates shown in table 3.1.

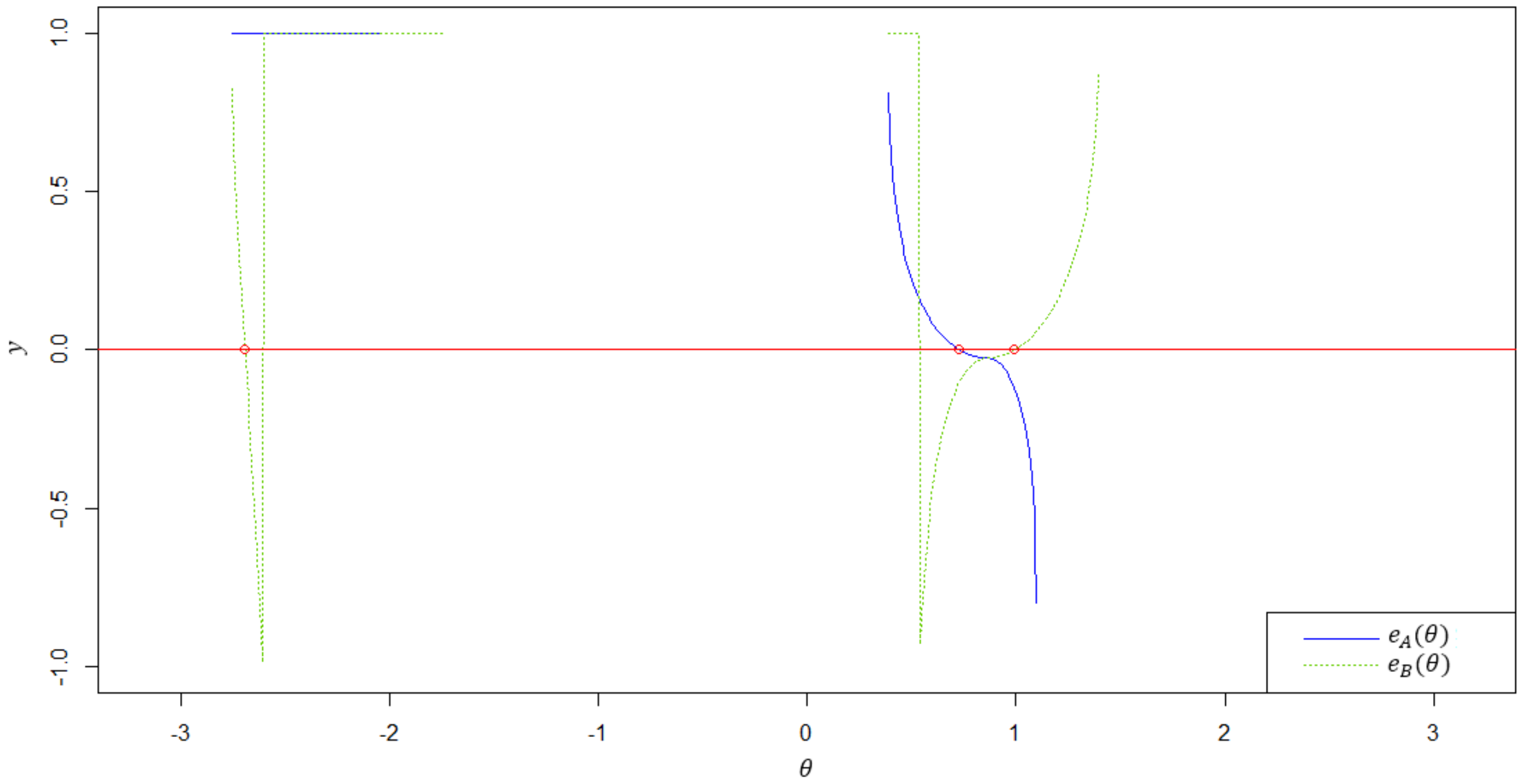

**Fig. 3.4.** Discrete angular proximity map for polynomial $p(z) = z^3 + C_1 z^2 + C_2 z + C_3$, with coefficients $C_1$, $C_2$, $C_3$ given by (3.7). This map was generated from $N = 1{,}000$ elements $LzC(\theta_k)$ associated to points $\theta_k$ in a regular partition of interval $[-\pi, \pi)$.

**Table 3.1**. Initial estimates for the roots of $p(z) = z^3 + C_1 z^2 + C_2 z + C_3$, with coefficients $C_1$, $C_2$, $C_3$ given by (3.7). These estimates were obtained from the map in figure 3.4.

| $i$ | $\hat{R}_i$ | $\hat{\theta}_i^*$ | $\lvert\Delta e_i\rvert$ |
|---|---|---|---|
| 1 | $0.2370095 - 0.0462304i$ | $0.9936258$ | $0.002449639$ |
| 2 | $0.2938045 - 0.2115416i$ | $0.7284104$ | $0.002559714$ |
| 3 | $-0.7278063 - 0.8652316i$ | $-2.6919246$ | $0.063539354$ |

In table 3.1, column $|\Delta e_i|$ shows the absolute vertical differences $\left|e(\theta_{k_i+1}) - e(\theta_{k_i})\right|$ between consecutive discrete elements $e(\theta_{k_i})$, $e(\theta_{k_i+1})$, for some $k_i \in \{0,1,2,\ldots,N-1\}$, such that point $(\hat{\theta}_i^*, 0) \in \mathbb{R}^2$ is contained in the line segment whose endpoints are $\left(\theta_{k_i}, e(\theta_{k_i})\right)$ and $\left(\theta_{k_i+1}, e(\theta_{k_i+1})\right)$; here, $e$ can be one of the functions $e_A$ or $e_B$ defined in (3.4a), (3.4b). For more details on how these values $|\Delta e_i|$ are computed, see annex 1 section 4, function `get_angle_roots`. Column $\hat{\theta}_i^*$ in table 3.1 shows the abscissae of "smooth" crossings of functions $\hat{e}_A(\theta)$, $\hat{e}_B(\theta)$ (obtained by linear interpolation) with horizontal axis $y = 0$. $\hat{\theta}_i^*$ are estimates for the "theta roots" of proximity map $e(\theta)$; that is to say, $\hat{\theta}_i^*$ are estimates for any of the values $\theta_i^*$ that lead, by means of a structure $LzC(\theta_i^*)$, to a root $R_i$ of equation (3.1); $i \in \{1,2,3\}$. Finally, column $\hat{R}_i$ in table 3.1 shows the estimates for roots $R_i$, obtained by implementing expression (3.5) in the `approxP3` function from annex 2 section 5.

It is worth mentioning here that in figure 3.4, reference values $\theta_1^*$, $\theta_2^*$, $\theta_3^*$ in (3.6) are represented by means of small circles drawn on the horizontal axis $y = 0$, with the purpose of graphically showing the proximity between theta roots $\theta_i^*$ and the corresponding "smooth" crossings $\hat{\theta}_i^*$.

An additional comment; in a certain way, the quantities $|\Delta e_i|$ listed in table 3.1 help us verify that the crossings of functions $e_A$, $e_B$ with horizontal axis $y = 0$ are "smooth" and not abrupt; we would expect that a value $|\Delta e_i|$ corresponding to a "smooth" crossing is not too large compared to the distance $\Delta\theta$ between two consecutive values $\theta_k$ and $\theta_{k+1}$ (which in this case is $\Delta\theta = 2\pi/1000 = 0.006283185$); the ratio of differences $\Delta e_i/\Delta\theta$ can be interpreted here as an approximation to the instantaneous rate of change of $e$ with respect to $\theta$, when $e$ crosses the horizontal axis $y = 0$ at point $\theta = \theta_i^*$. As we can see from table 3.1, values $|\Delta e_1|$ and $|\Delta e_2|$ are both less than $\Delta\theta = 0.006283$, while $|\Delta e_3|$ is about 10 times greater than $\Delta\theta$; this tells us that the variation of $e_A$ and $e_B$ around values $\theta_1^*$ and $\theta_2^*$ is small, while the instantaneous rate of change of $e_B$ around $\theta_3^*$ is large, although not too much to consider it an abrupt change, as we can see from figure 3.4. In the script from annex 2 section 5, within function `approxP3`, we arbitrarily defined the following criterion: all values $|\Delta e| \leq 1.0$ associated to crossings of $\hat{e}(\theta)$ with horizontal axis $y = 0$, where the horizontal distance between two consecutive discrete values $\theta_k$ and $\theta_{k+1}$ is $\Delta\theta = 2\pi/N$, correspond numerically to "smooth" crossings (this can be found in the part where function `get_angle_roots` is called within function `approxP3`). We believe that this criterion works well, and helps in most cases to distinguish between smooth crossings and abrupt crossings of $e_A$ or $e_B$; in figure 3.4, for example, two abrupt crossings of function $e_B$ can be seen: one slightly to the right of $\theta_3^*$, and the other, slightly to the left of $\theta_1^*$; in both cases, clearly $|\Delta e| > 1.0$ (In fact, in this situation, $|\Delta e|$ is close to 2, as can it be seen from figure 3.4). Due to this criterion, these two abrupt crossings were not included in table 3.1. The reader can modify this criterion, if the need arises when generating and studying some other proximity map; there are infinite possibilities, and the shapes of proximity maps can be very simple, or they can be quite capricious and complex. Our goal, however, is not to study different categories of proximity maps associated with polynomials of a complex variable, but to extend the LC method to univariate polynomials of any degree with a random distribution of their roots. Nevertheless, we believe that a more detailed study of the various categories of proximity maps would be a good research topic, worth of being explored in much more detail.

A basic question that arises is this: why do abrupt crossings of function $e(\theta)$ occur in proximity maps such as the one in figure 3.4? Informally, the reason is that terminal semi-line $tL$, when extremely close to touching point $0 + 0i$ in complex plane $\mathbb{C}$, intersects with z-circumference $zC$ at a point $I'$ with an extremely small modulus, which generates a value $c_0$ with an extremely large modulus, since $c_0 = -C_3/I'$; consequently, $c_0$ is far away from both fixed point $P_1$ and $i_0 = P_1 + \sqrt{|I' - (C_2 - P_1^2)|}v_\theta$ (but always at line $\ell_1$), while the distance between $i_0$ and $P_1$ is relatively small; therefore, function $e(\theta)$, according to its definition (3.3), assumes values very close to extreme values $\pm 1$ in its codomain; the change of sign, from one extreme value of $e(\theta)$ to another occurs when the intersection $I'$ between $tL$ and $zC$ passes from one side of the point $0 + 0i$ to the other on the $zC$ "bridge", as the value $\theta$ varies; this causes $c_0$ to make a "quantum jump" from one side of $P_1$ to the other within line $\ell_1$. Figure 3.5 schematically illustrates this informal explanation.

In figure 3.4, the two abrupt crossings of $e_B$ are separated by an angular distance of $\pi$ radians. It is to be expected that in most proximity maps for univariate polynomials of degree 3, we will observe this separation between abrupt crossings attributed to $tL$ being too close to the origin, since lines $\ell_1(\theta)$ and $\ell_1(\theta + \pi)$ occupy the same geometric place, and therefore both produce the same terminal semi-line $tL$ and the same z-circumference $zC$. For this same reason, as it can be seen from figure 3.4, these two abrupt crossings have opposite directions to each other; one jumps from $-1$ to $+1$, while the other jumps from $+1$ to $-1$, as $\theta$ increases its value.

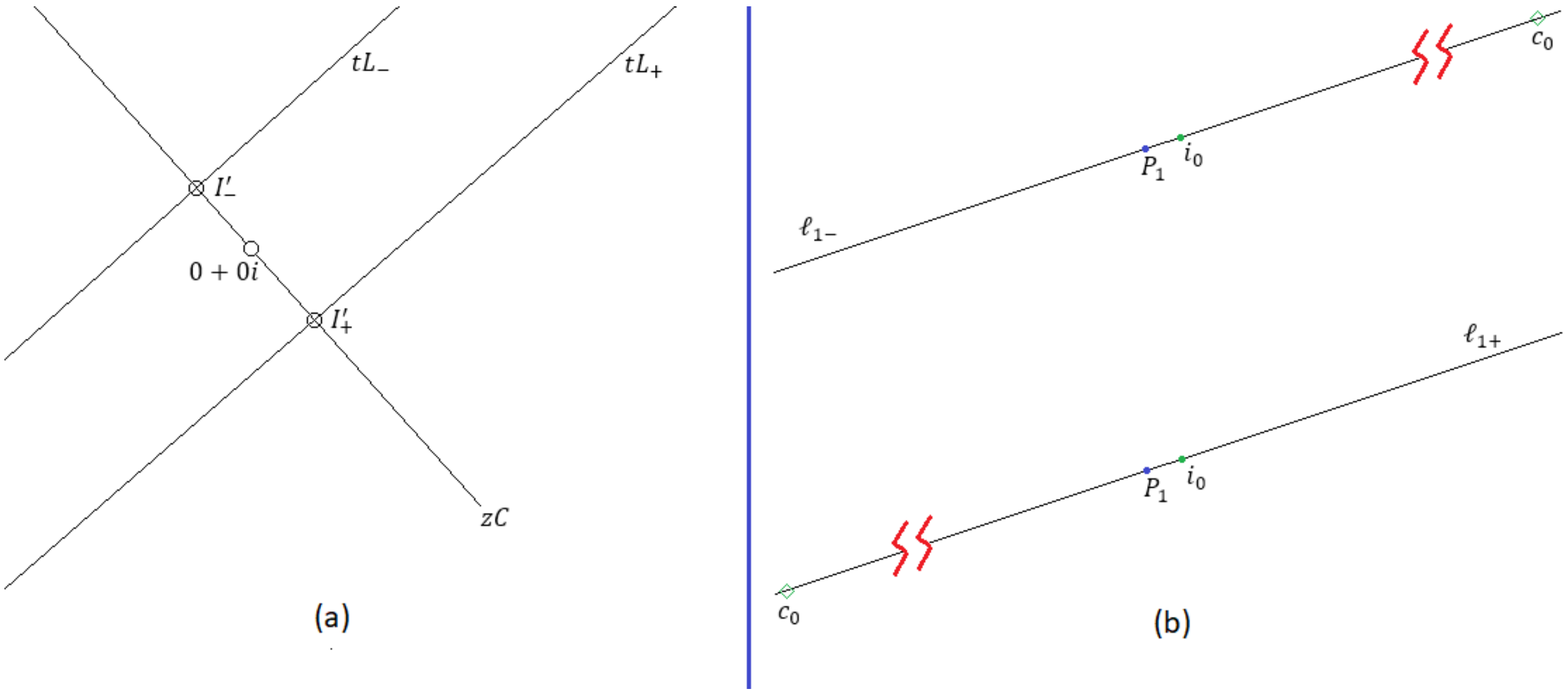

**Fig. 3.5.** Schematic (a) shows, on a "microscopic" scale (extreme zoom in), the intersection $I'$ between a terminal semi-line $tL$ and a z-circumference $zC$, which occurs extremely close to the origin $0 + 0i$; the passage from one side of the $zC$ "bridge" to the other, by infinitesimally changing value $2\theta$ (the direction of $tL$) is denoted by means of symbols $+$ and – in the subscripts for the names of the elements shown in both schematics (a) and (b); in schematic (a) the variation of $zC$ is negligible. Schematic (b), on the other hand, shows on a "macroscopic" scale (extreme zoom out) the states + and – of line $\ell_1$ which correspond to the elements $tL$ and $zC$ from schematic (a). As can be seen, the passage of $I'$ from one side of the $zC$ "bridge" to the other at the microscopic level produces a "quantum jump" of point $c_0$ at line $\ell_1$. The angular difference between lines $\ell_{1+}$ and $\ell_{1-}$ is infinitesimal; $P_1$ is of course the same point both at $\ell_{1+}$ and at $\ell_{1-}$ (one of these two lines is translated to another location, to avoid overlapping).

Let's go back to the estimates for the roots in this numerical example. How close are estimates $\hat{R}_i$, $\hat{\theta}_i^*$ from the "true" roots? Let's compute the absolute relative errors associated with the estimates from table 3.1, with respect to the corresponding reference values in (3.6):

$$\left|\frac{R_1-\hat{R}_2}{R_1}\right| = 1.920883\times10^{-5} \qquad \left|\frac{\theta_1^*-\hat{\theta}_2^*}{\theta_1^*}\right| = 1.541473\times10^{-5}$$
$$\left|\frac{R_2-\hat{R}_1}{R_2}\right| = 5.772065\times10^{-5} \qquad \left|\frac{\theta_2^*-\hat{\theta}_1^*}{\theta_2^*}\right| = 1.937968\times10^{-5} \qquad (3.8)$$
$$\left|\frac{R_3-\hat{R}_3}{R_3}\right| = 6.467441\times10^{-7} \qquad \left|\frac{\theta_3^*-\hat{\theta}_3^*}{\theta_3^*}\right| = 3.623760\times10^{-8}$$

The values in (3.8) tell us, for example, that the absolute percent deviation of $\hat{\theta}_1^*$ from corresponding theta root $\theta_2^*$ is about 0.002%; for calculations with a desired accuracy of 4 or 5 decimal places, the initial estimates in table 3.1 could be close enough to the desired solution, so it would no longer be necessary to use some other method (such as Newton-Raphson) to increase estimate accuracy; of course, the need to use a second stage of refinement for initial estimates depends on the problem in question, although for general purposes, it is almost always necessary to refine initial estimates, since what is sought is to obtain approximations as accurate as possible. One way to increase the accuracy of our initial estimates obtained by means of discrete proximity maps, is to use more elements $LzC$ in their construction; that is to say, to increase the resolution of proximity maps. In the following numerical example, we will explore this alternative.

## Numerical Example 3.2

This time we will build a discrete angular proximity map for estimating the roots of equation

$$z^3 + (1+i)z^2 + (2+2i)z + (3+3i) = 0. \qquad (3.9)$$

In order to obtain reference roots $R_1$, $R_2$, $R_3$ that help us compute approximation errors associated with estimates $\hat{R}_1$, $\hat{R}_2$, $\hat{R}_3$ to be generated in this example, we will use, within the script listed in annex 3 section 2, an R function called `polyroot`, which finds all the roots of a polynomial in a real or complex variable, by means of the Jenkins-Traub algorithm, which is, to our knowledge, one of the most robust algorithms for solving polynomial equations in a single variable. For more details on the `polyroot` function, see [3.6]. The comments in annex 3 section 2 explain how to apply the `polyroot` function to produce reference roots $R_1$, $R_2$, $R_3$ for equation (3.9), which are shown in expressions (3.10), together with their corresponding "theta roots" (inclination angles $\theta_i^*$ for line $\ell_1$ with fixed point $P_1 = -(1+i)/2$, which cause $\ell_1$ to contain any of the roots $R_i$):

$$R_1 = -0.264046 + 1.426772i \qquad \theta_1^* = 1.4489422rad$$
$$R_2 = -1.334806 - 0.429244i \qquad \theta_2^* = 3.0570371rad \qquad (3.10)$$
$$R_3 = 0.598852 - 1.997529i \qquad \theta_3^* = -0.9377592rad$$

The discrete proximity map for equation (3.9), obtained with the script in annex 3 section 2, and using $N = 1{,}000$ elements $LzC(\theta_k)$, $k = 0,1,2,\ldots,N-1$, with $\theta_k = -\pi + 2\pi k/N$, is shown in figure 3.6. The numerical values for the "smooth" crossings $\hat{\theta}_i^*$ of this discrete map $\hat{e}(\theta)$ with horizontal axis $y = 0$ are listed in table 3.2, together with their corresponding estimates $\hat{R}_i$.

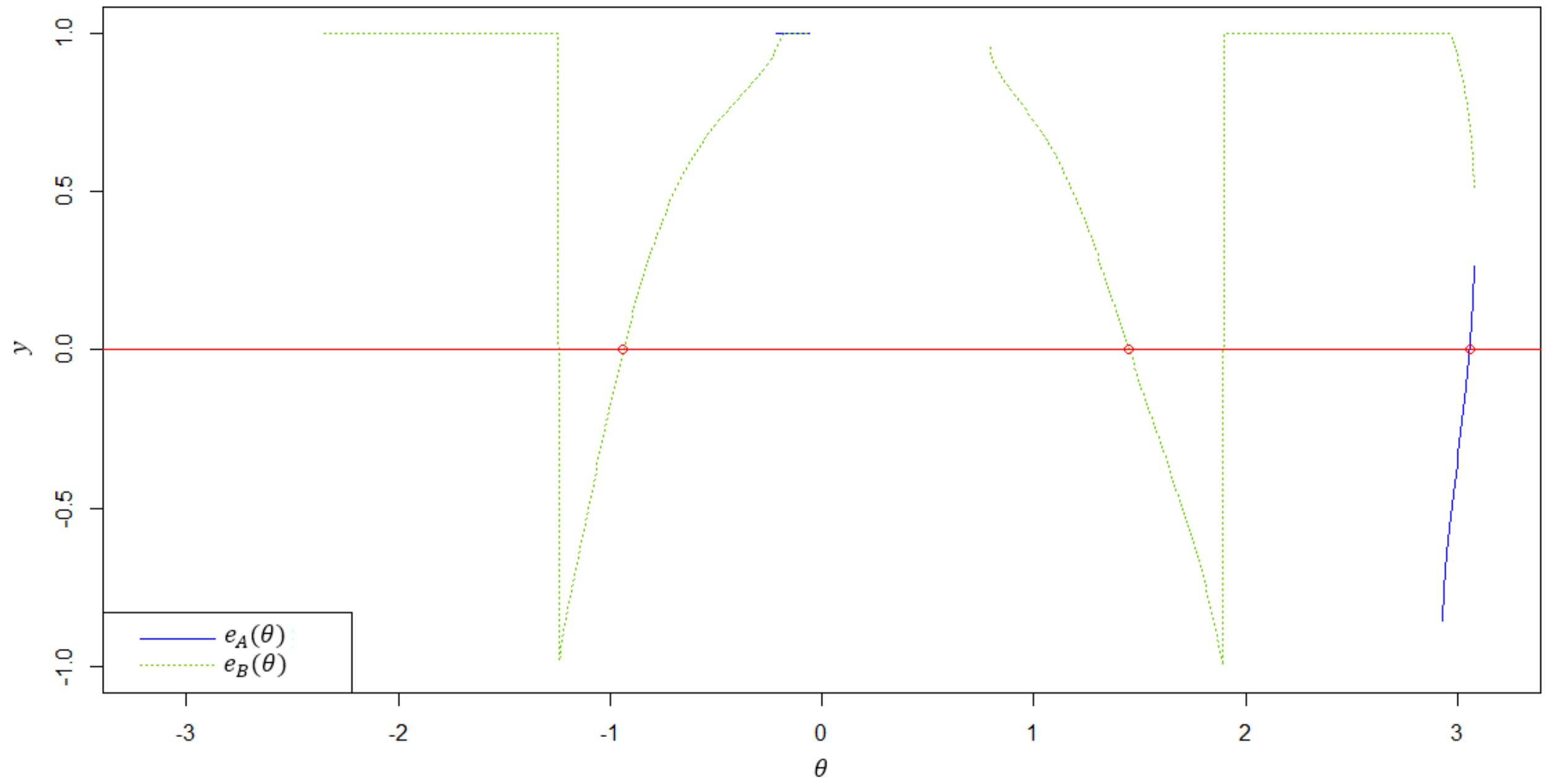

**Fig. 3.6.** Discrete angular proximity map for polynomial $p(z) = z^3 + (1+i)z^2 + (2+2i)z + (3+3i)$, generated from $N = 1{,}000$ elements $LzC(\theta_k)$ associated with points $\theta_k$ in a regular partition of interval $[-\pi, \pi)$.

**Table 3.2**. Initial estimates for the roots of $p(z) = z^3 + (1+i)z^2 + (2+2i)z + (3+3i)$, obtained from the map in figure 3.6.

| $i$ | $\hat{R}_i$ | $\hat{\theta}_i^*$ | $\lvert\Delta e_i\rvert$ |
|---|---|---|---|
| 1 | $-0.264048 + 1.426775i$ | 1.4489434 | 0.01274373 |
| 2 | $0.598848 - 1.997503i$ | −0.9377529 | 0.01699775 |
| 3 | $-1.334716 - 0.429171i$ | 3.0569412 | 0.05113541 |

In figure 3.6, reference values $\theta_1^*$, $\theta_2^*$, $\theta_3^*$ in (3.10) are also represented by means of small circles on the horizontal axis $y = 0$, for the purpose of graphically showing the proximity between reference theta roots $\theta_i^*$ and smooth crossings $\hat{\theta}_i^*$ of map $\hat{e}(\theta)$ with horizontal axis $y = 0$.

How close to the “real” roots are the estimates listed in table 3.2? To answer this question, we compute the absolute relative errors for the estimates in table 3.2 with respect to the corresponding reference values (3.10):

$$\left|\frac{R_1-\hat{R}_1}{R_1}\right| = 1.977020 \times 10^{-6} \qquad \left|\frac{\theta_1^*-\hat{\theta}_1^*}{\theta_1^*}\right| = 7.882217 \times 10^{-7}$$
$$\left|\frac{R_2-\hat{R}_3}{R_2}\right| = 8.275962 \times 10^{-5} \qquad \left|\frac{\theta_2^*-\hat{\theta}_3^*}{\theta_2^*}\right| = 3.136389 \times 10^{-5} \qquad (3.11)$$
$$\left|\frac{R_3-\hat{R}_2}{R_3}\right| = 1.224262 \times 10^{-5} \qquad \left|\frac{\theta_3^*-\hat{\theta}_2^*}{\theta_3^*}\right| = 6.687756 \times 10^{-6}$$

Comparing relative errors (3.8) and (3.11), we see that *the average value* of approximation errors for roots $R_i$ obtained in this example 3.2 is of similar magnitude to that of example 3.1; the same can be said about the average values of approximation errors for theta roots $\theta_i^*$ obtained in both examples. Note that in both instances, we used maps with $N = 1{,}000$ elements $LzC(\theta_k)$ associated with points $\theta_k$ in a regular partition of interval $[-\pi, \pi)$. Next, we will see what happens with the accuracy of the approximations in this example 3.2 if we use a greater number of elements $LzC$.

### Increasing discrete proximity map's resolution

The *resolution* of a discrete angular proximity map is inversely related to the distance between two consecutive discrete elements $\theta_k$ and $\theta_{k+1}$ in the regular partition of interval $[-\pi, \pi)$; for the case of the map shown in figure 3.6, its resolution is 1,000 elements/$2\pi$ radians ≅ 159 elements/rad. How do relative errors for estimates associated with proximity maps behave if we increase map resolution by an order of magnitude? Let us see what happens if we now use $N = 10{,}000$ elements $LzC(\theta_k)$ to construct the proximity map for equation (3.9); in such a situation, the map resolution is now around 1,592 elements/rad. Using again the script from annex 3 section 2 with the same reference values $R_1$, $R_2$, $R_3$ in (3.10), but now with argument `N=10000` in the call to function `approxP3`, we obtain the proximity map from figure 3.7 (also representing reference theta roots $\theta_1^*$, $\theta_2^*$, $\theta_3^*$ in (3.10) by means of small circles on the horizontal axis $y = 0$), and the numerical values for corresponding estimates, listed in table 3.3:

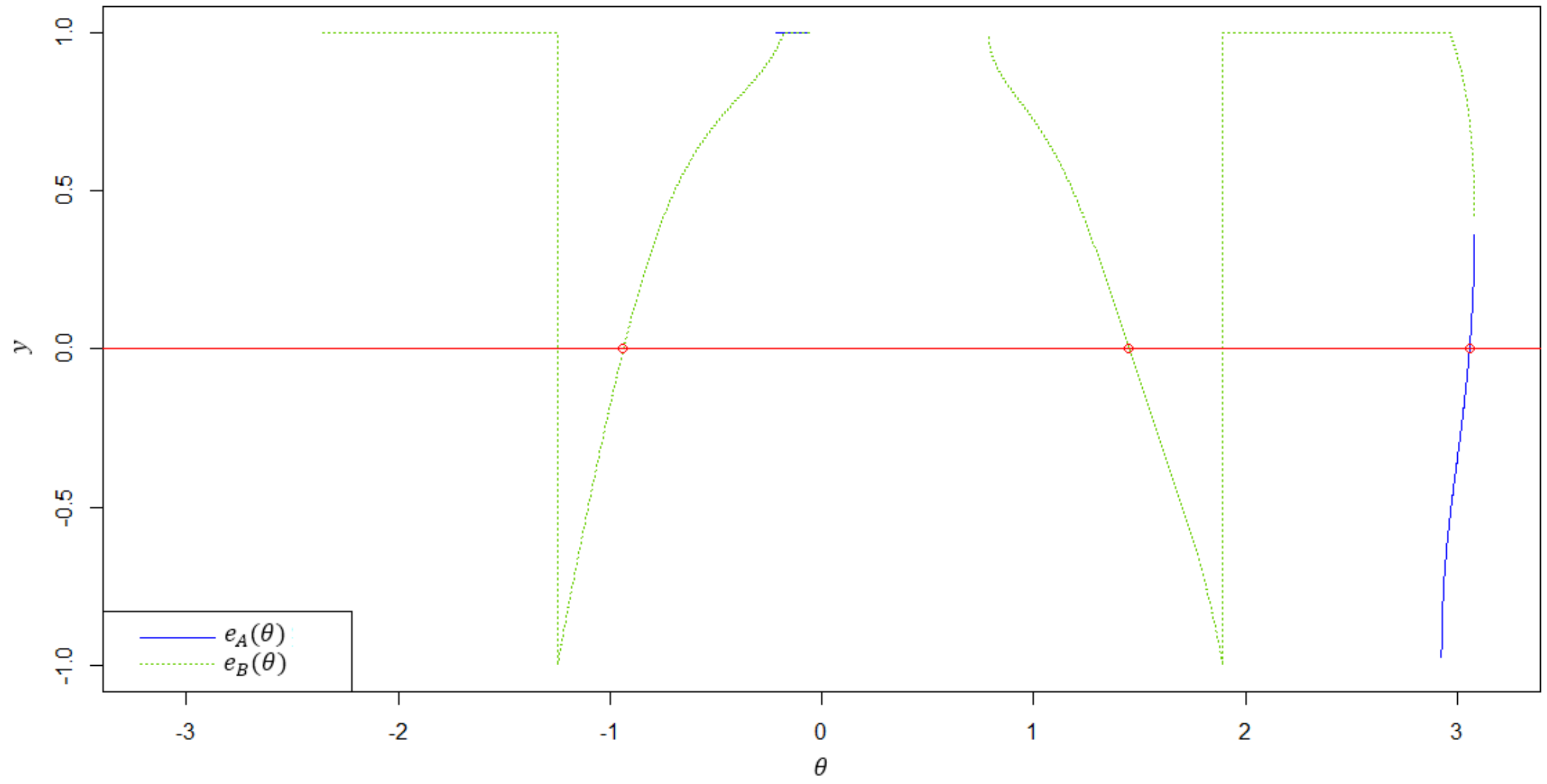

**Fig. 3.7.** Discrete angular proximity map for polynomial $p(z) = z^3 + (1 + i)z^2 + (2 + 2i)z + (3 + 3i)$, generated from $N = 10{,}000$ elements $LzC(\theta_k)$ associated with points $\theta_k$ in a regular partition of interval $[-\pi, \pi)$.

**Table 3.3**. Initial estimates for the roots of $p(z) = z^3 + (1 + i)z^2 + (2 + 2i)z + (3 + 3i)$, obtained from the map in figure 3.7.

| $i$ | $\hat{R}_i$ | $\hat{\theta}_i^*$ | $\lvert\Delta e_i\rvert$ |
|---|---|---|---|
| 1 | $-0.264046 + 1.426772i$ | $1.4489422$ | $0.001274091$ |
| 2 | $0.598852 - 1.997528i$ | $-0.9377591$ | $0.001695251$ |
| 3 | $-1.334805 - 0.429243i$ | $3.0570361$ | $0.005134830$ |

The absolute relative errors for estimates listed in table 3.3, with respect to reference values (3.10), are:

$$\left|\frac{R_1-\hat{R}_1}{R_1}\right| = 4.969709 \times 10^{-9} \qquad \left|\frac{\theta_1^*-\hat{\theta}_1^*}{\theta_1^*}\right| = 1.981388 \times 10^{-9}$$
$$\left|\frac{R_2-\hat{R}_3}{R_2}\right| = 8.186932 \times 10^{-7} \qquad \left|\frac{\theta_2^*-\hat{\theta}_3^*}{\theta_2^*}\right| = 3.101581 \times 10^{-7} \qquad (3.12)$$
$$\left|\frac{R_3-\hat{R}_2}{R_3}\right| = 1.646799 \times 10^{-7} \qquad \left|\frac{\theta_3^*-\hat{\theta}_2^*}{\theta_3^*}\right| = 8.995750 \times 10^{-8}$$

As expected, the accuracy of the estimates associated with the discrete proximity map increases if we increase the map's resolution; in this case, an increase of 10× in map resolution, from 159.2 elements/rad to 1,592 elements/rad, produces a decrease in absolute relative error of 398× for estimates $\hat{R}_1$ and $\hat{\theta}_1^*$, a decrease in absolute relative error of 74× for estimates $\hat{R}_2$ and $\hat{\theta}_2^*$, and a decrease in absolute relative error of 101× for estimates $\hat{R}_3$ and $\hat{\theta}_3^*$; on average, the absolute relative error for root estimates decreases 191× with a 10× increase in map resolution. This suggests a decrease of two orders of magnitude in relative errors of root estimates associated to a discrete proximity map with respect to an increase of an order of magnitude in its resolution.

What differences can we observe between the maps in figures 3.6 and 3.7? At first glance, the two look almost the same, except for the amplitude of a "gap" between the discrete versions of functions $e_A(\theta)$ and $e_B(\theta)$ located on the right side of both graphs; the amplitude of this gap is reduced by increasing map resolution. But why the existence of this gap?

When analyzing the numerical values for the discrete proximity map built with $N = 1{,}000$ elements, we observe that the gap occurs at $\theta_{990} = -\pi + 2\pi(990)/1000 = 3.078760801$, where the difference between $e_A$ and $e_B$ is $e_B(\theta_{990}) - e_A(\theta_{990}) = 0.5049939 - 0.26285594 = 0.2421379$ (this is the gap's amplitude in figure 3.6); for the next discrete angular value in the proximity map, i.e., for $\theta_{991} = -\pi + 2\pi(991)/1000 = 3.085043986$, there are no available values for $e_A$ and $e_B$, since in the corresponding geometric construction $LzC(\theta_{991})$, terminal semi-line $tL(\theta_{991})$ does not intersect with z-circumference $zC(\theta_{991})$, although both elements are very close to each other; see figures 3.8 and 3.9 to observe the behavior of respective structures $LzC(\theta_{990})$ and $LzC(\theta_{991})$; such behavior suggests that elements $tL$ and $zC$ associated with equation (3.9) must intersect tangentially at a single point for some angle $\theta$ between $\theta_{990}$ and $\theta_{991}$, much like the scenario shown in the lower left corner of figure 3.3.

With the help of the R script listed in annex 4 section 1, it is possible to interactively analyze structures $LzC$ associated with cubic polynomials in a complex variable; with some slight modifications to this script, it is possible to obtain graphs such as those shown in figures 3.8 and 3.9 for equation (3.9).

From the analysis in the above paragraphs, we can see that the gap on the right side of figures 3.6 and 3.7 is an artifact caused by the discretization of the proximity map, and indicates the presence of a tangential intersection between $tL$ and $zC$; this could be potentially problematic if a theta root occurs inside a gap, at the tangential intersection between $tL$ and $zC$, or very close to that intersection, beyond the resolution of the discrete proximity map; in this example this kind of problem does not arise, even though $\theta_2^* = 3.0570371$ is located relatively close to the gap. Even with the lower resolution of 159 elements/rad, it was sufficient in this case to obtain adequate initial estimates to all theta roots.

In the following numerical example, we will see a situation where we have the problem of a theta root inside a gap between the discrete versions of $e_A$ and $e_B$; we will see, however, that this problem can be addressed if we consider plotting a local region of the proximity map to a higher resolution at a neighborhood around such gap.

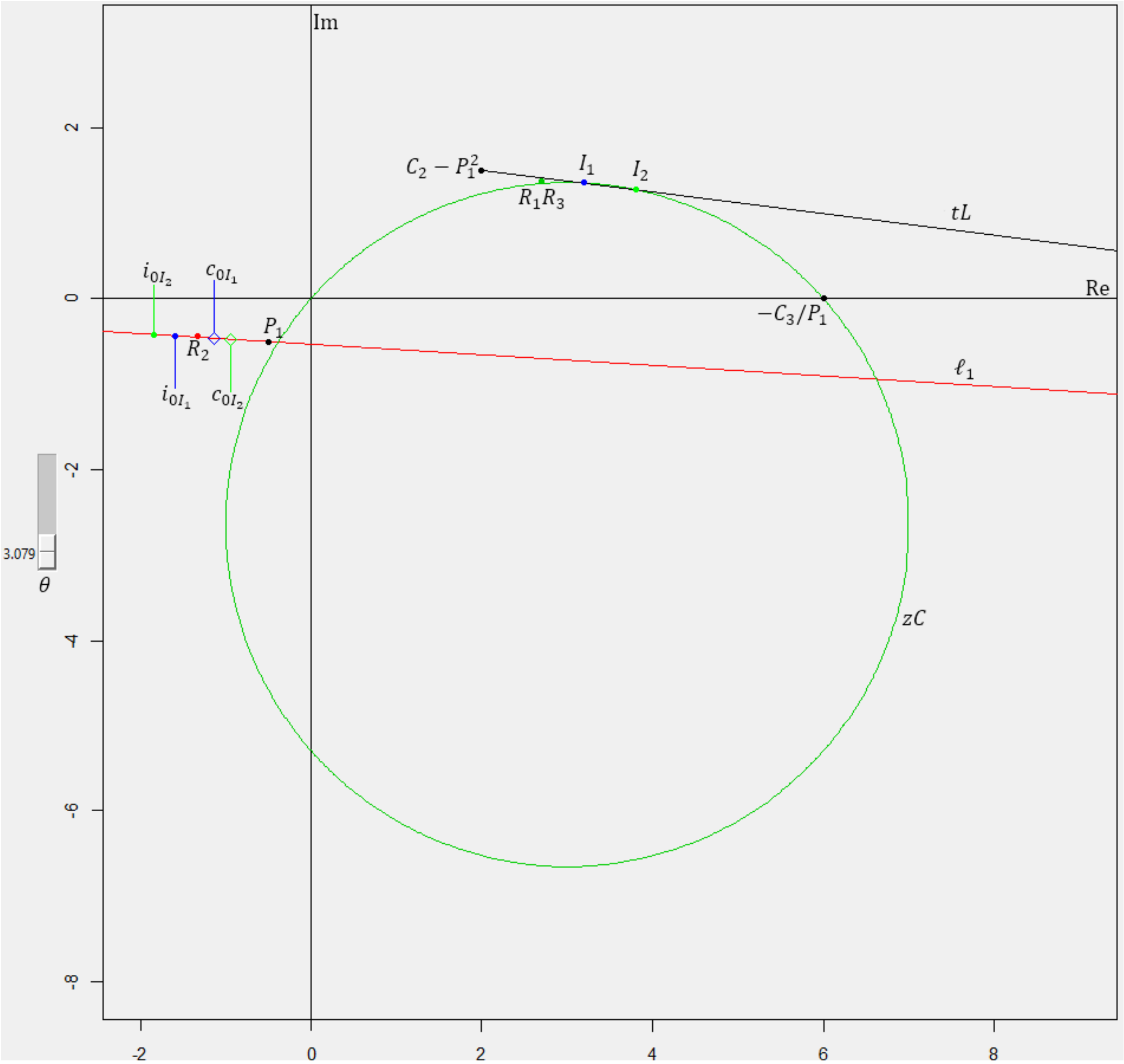

**Fig. 3.8.** Structure $LzC(\theta_{990} = 3.079)$ associated with polynomial $p(z) = z^3 + (1+i)z^2 + (2+2i)z + (3+3i)$. The intersections $I_1$, $I_2$ between $tL$ and $zC$ are "close" to becoming a single tangential intersection with a small increment of variable $\theta$ ($\ell_1$'s inclination angle). Although $R_2$ seems to be very close to $\ell_1$, $R_2 \notin \ell_1$.

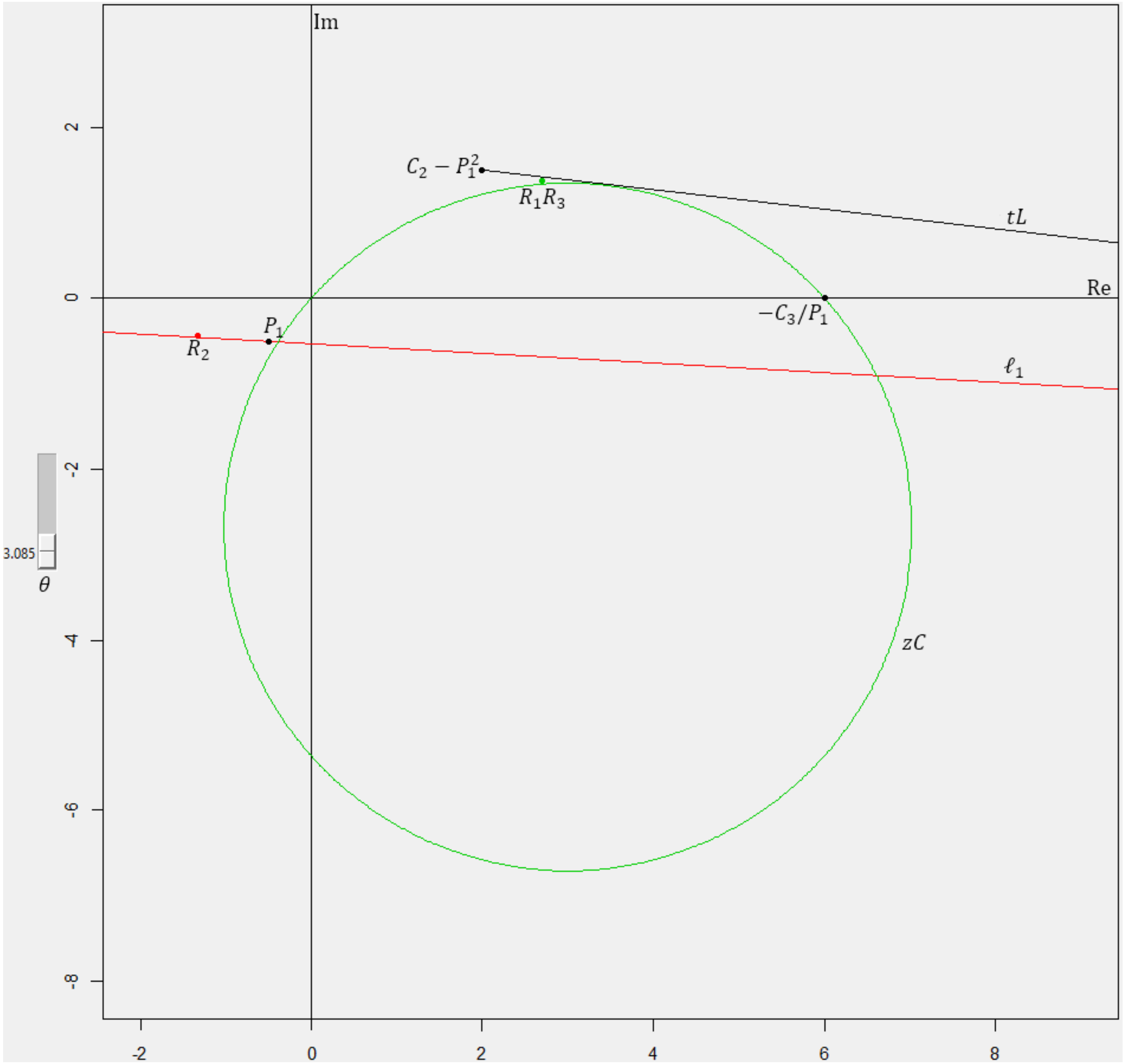

**Fig. 3.9.** Structure $LzC(\theta_{991} = 3.085)$ associated with polynomial $p(z) = z^3 + (1+i)z^2 + (2+2i)z + (3+3i)$. $tL$ and $zC$ do not touch each other, although they are about to doing so tangentially at a single point with a small decrement of variable $\theta$ ($\ell_1$'s inclination angle). Although $R_2$ seems to be very close to $\ell_1$, $R_2 \notin \ell_1$.

## Numerical Example 3.3

In this example we will use the script in annex 3 section 2 to generate a new discrete proximity map, this time by fixing the random number generator seed with the instruction `set.seed(2)` at the beginning of this script.

The randomly generated reference roots $R_i$, along with inclination angles $\theta_i^*$ for lines $\ell_1$ (with fixed point $P_1 = -C_1/2$) which contain roots $R_i$, are:

$$R_1 = -0.6302355 + 0.4047481i \qquad \theta_1^* = 3.0328690 rad$$
$$R_2 = 0.1466527 - 0.6638962i \qquad \theta_2^* = -1.6273890 rad \qquad (3.13)$$
$$R_3 = 0.8876787 + 0.8869499i \qquad \theta_3^* = 0.6961912 rad$$

When plugging generated roots $R_1, R_2, R_3$ into Vieta's relations (3.2-a), (3.2-b) and (3.2-c), we obtain these coefficients:

$$C_1 = -0.4040959 - 0.6278018i,$$
$$C_2 = -0.0231298 - 0.1811857i, \qquad (3.14)$$
$$C_3 = 0.2672721 - 0.5804607i.$$

The proximity map associated with coefficients (3.14) and generated from $N = 1{,}000$ elements $LzC(\theta_k)$, $k = 0,1,2,\dots,N-1$, with $\theta_k = -\pi + 2\pi k/N$, is shown in figure 3.10. Corresponding estimates $\hat{R}_i$, $\hat{\theta}_i^*$ are listed in table 3.4.

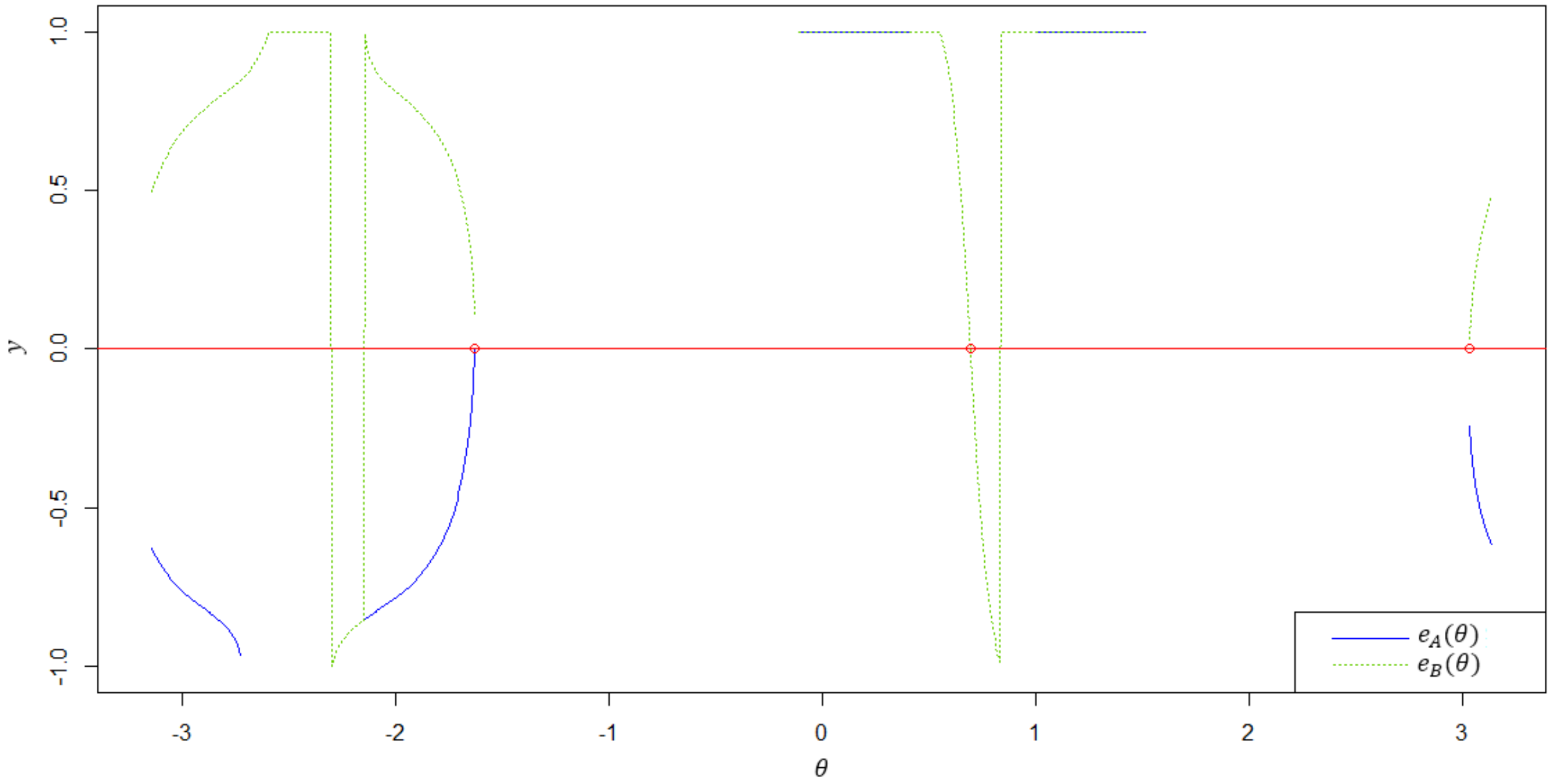

**Fig. 3.10.** Discrete angular proximity map for polynomial $p(z) = z^3 + C_1z^2 + C_2z + C_3$, with coefficients $C_1, C_2, C_3$ given by (3.14). This map was generated from $N = 1{,}000$ elements $LzC(\theta_k)$ associated to points $\theta_k$ in a regular partition of interval $[-\pi, \pi)$.

**Table 3.4**. Initial estimates for the roots of $p(z) = z^3 + C_1 z^2 + C_2 z + C_3$, with coefficients $C_1, C_2, C_3$ given by (3.14). These estimates were obtained from the map in figure 3.10.

| $i$ | $\hat{R}_i$ | $\hat{\theta}_i^*$ | $\lvert\Delta e_i\rvert$ |
|---|---|---|---|
| 1 | $0.8876652 + 0.8869355i$ | $0.6961885$ | $0.07887739$ |
| 2 | $0.1465865 - 0.6642112i$ | $-1.6274382$ | $0.11654792$ |
| 3 | $-0.6257991 + 0.4026645i$ | $3.0347785$ | $0.27385119$ |

As in previous numerical examples from this chapter, the theta roots $\theta_i^*$ in (3.13) are represented in figure 3.10 by means of small circles on the horizontal axis $y = 0$, in order to visually compare computed estimates with reference values; from here, we can see that theta roots $\theta_1^*$ and $\theta_2^*$ occur, each of them, very close to a gap (between functions $\hat{e}_A$ and $\hat{e}_B$) associated with a tangential intersection between $tL$ and $zC$; at the vicinity of $\theta_2^* = -1.6273890$, function $\hat{e}_A(\theta)$ barely crosses horizontal axis $y = 0$ to produce estimate $\hat{\theta}_2^* = -1.6274382$, before fading into the gap it forms with $\hat{e}_B(\theta)$; on the other hand, none of the functions $\hat{e}_A(\theta)$ and $\hat{e}_B(\theta)$ cross horizontal axis $y = 0$ at the vicinity of the gap they form near $\theta_1^* = 3.0328690$; thus, the discrete proximity map from figure 3.10, with a resolution of $1000/2\pi$ elements/rad, this time is not able to produce a crossing that approximates $\theta_1^*$. The third row in table 3.4, however, shows an estimate for $\theta_1^*$, namely $\hat{\theta}_3^* = 3.0347785$, thanks to an *additional routine* that was programmed within function `approxP3` listed in annex 2 section 5; this routine operates as follows:

***Additional routine* to rescue a root estimate not detected by a global discrete proximity map associated with a univariate cubic polynomial**

If the number $s$ of "smooth" crossings $\hat{\theta}_i^*$ produced by discrete proximity map at interval $[-\pi, \pi)$ is less than 3 (the degree of associated univariate cubic polynomial):

- Compute $d = \min_{\hat{e}_A \hat{e}_B < 0} |\hat{e}_A(\theta_m) - \hat{e}_B(\theta_m)|$; values $\theta_m$ are points in the regular partition of interval $[-\pi, \pi)$ that support the global discrete proximity map
- Determine, from values $\theta_m$ involved in the computation of $d$, which of them is the value $\theta_{m^*}$ such that $|\hat{e}_A(\theta_{m^*}) - \hat{e}_B(\theta_{m^*})| = d$
- Let $\hat{\theta}_{s+1}^* = \theta_{m^*}$
- Let $|\Delta e_{s+1}| = d$
- Compute $\hat{R}_{s+1}$ by applying criterion (3.5) on structure $LzC(\hat{\theta}_{s+1}^*)$

Note here that $|\Delta e_{s+1}|$ differs slightly from the definition we had originally given to quantity $|\Delta e_i|$ in numerical example 3.1.

In this case we were fortunate that this additional routine helped us find an initial estimate $\hat{\theta}_3^*$ for theta root $\theta_1^*$ (at best this additional routine recovers only one approximation); but how close is $\hat{\theta}_3^*$ from the true value $\theta_1^*$? Is $\hat{\theta}_3^*$ a lower quality estimate compared to the other two estimates $\hat{\theta}_1^*$ and $\hat{\theta}_2^*$? Let us compute the absolute relative errors of the estimates listed in table 3.4, with respect to the reference roots in (3.13):

$$\begin{aligned} \left|\frac{R_1-\hat{R}_3}{R_1}\right| &= 6.543734 \times 10^{-3} & \left|\frac{\theta_1^*-\hat{\theta}_3^*}{\theta_1^*}\right| &= 6.296100 \times 10^{-4} \\ \left|\frac{R_2-\hat{R}_2}{R_2}\right| &= 4.735346 \times 10^{-4} & \left|\frac{\theta_2^*-\hat{\theta}_2^*}{\theta_2^*}\right| &= 3.024715 \times 10^{-5} \\ \left|\frac{R_3-\hat{R}_1}{R_3}\right| &= 1.575987 \times 10^{-5} & \left|\frac{\theta_3^*-\hat{\theta}_1^*}{\theta_3^*}\right| &= 3.939760 \times 10^{-6} \end{aligned} \qquad (3.15)$$

From relative errors (3.15), we see that the absolute relative percent deviations for $\hat{\theta}_1^*$, $\hat{\theta}_2^*$, $\hat{\theta}_3^*$ are, respectively, 0.0004%, 0.003% and 0.06%; the relative error associated with $\hat{\theta}_3^*$ is larger than the relative errors associated with $\hat{\theta}_1^*$ and $\hat{\theta}_2^*$, but not so much as to label $\hat{\theta}_3^*$ a bad initial estimate; in fact, $\hat{\theta}_3^*$ coincides with $\theta_1^*$ at the first two decimal places, while the real and imaginary parts of $\hat{R}_3$ match the corresponding parts of $R_1$ if all quantities involved are rounded off to two decimal places; this could still be viewed as an acceptable initial estimate. On the other hand, the magnitudes of relative errors for estimates $\hat{R}_1$, $\hat{\theta}_1^*$, and $\hat{R}_2$, $\hat{\theta}_2^*$ in this numerical example 3.3 are similar to those observed in numerical examples 3.1 and 3.2 for proximity maps with a resolution of $1{,}000/2\pi$ elements/rad.

As for the abrupt crossings shown in figure 3.10, we see that there are now 3 abrupt crossings of $\hat{e}_B(\theta)$ with the horizontal axis $y = 0$. Two of these crossings are caused, as explained above at the schematics in figure 3.5, by the passage of terminal semi-line $tL$ through the origin $0 + 0i$; these two abrupt crossings have a horizontal separation of $\pi$ radians within the map. In fact, these two crossings occur, one between $\theta_{133} = -2.305929008$ and $\theta_{134} = -2.299645822$, and the other between $\theta_{633} = 0.835663646$ and $\theta_{634} = 0.841946831$. On the other hand, the third abrupt crossing of $\hat{e}_B(\theta)$ with horizontal axis $y = 0$, which occurs between $\theta_{158} = -2.148849375$ and $\theta_{159} = -2.142566190$, is caused by the passage of line $\ell_1$ through the origin $0 + 0i$; this event causes the center of $zC$ to experience an abrupt change of position (in these circumstances, the radius of $zC$ tends to be extremely large) and consequently, there is an abrupt change in the configuration of the intersections between $tL$ and $zC$, which finally triggers, in this case, an abrupt crossing of $\hat{e}_B(\theta)$ with horizontal axis $y = 0$. In figure 3.10, we can see that, by moving from $\theta_{158}$ to $\theta_{159}$, function $\hat{e}_B(\theta)$ suddenly jumps upwards, and at the same time, it "bifurcates", in a smooth transition, towards function $\hat{e}_A(\theta)$. As we said earlier in numerical example 3.2, annex 4 section 1 includes an R script that allows to interactively observe the evolution of a structure $LzC(C_1, C_2, C_3, \theta)$ associated with a polynomial of degree 3 in a complex variable; the user can vary, using a slider control, the angle $\theta$, and observe how the elements of the structure $LzC$ evolve. The reader is invited to use this script, along with reference roots $R_1$, $R_2$, $R_3$ in (3.13), to observe carefully what happens to the elements of the structure $LzC$ for values $\theta$ close to each of the three abrupt crossings found in figure 3.10.

Before concluding this numerical example, let us analyze in more detail the proximity map associated with coefficients (3.14); instead of constructing the "global" map at interval $[-\pi, \pi)$, this time we will focus on region [3.0,3.1), which contains theta root $\theta_1^* = 3.0328690$ not initially detected by the global map of figure 3.10. By making a few slight modifications to the script in annex 3 section 2, we obtain the map shown in figure 3.11; as for the modifications, it is only necessary to change the last two arguments within the call to function `approxP3`, as well as the two horizontal limits `xlim` within the call to function `plot`, to specify the boundaries of the region where the map is to be constructed and drawn.

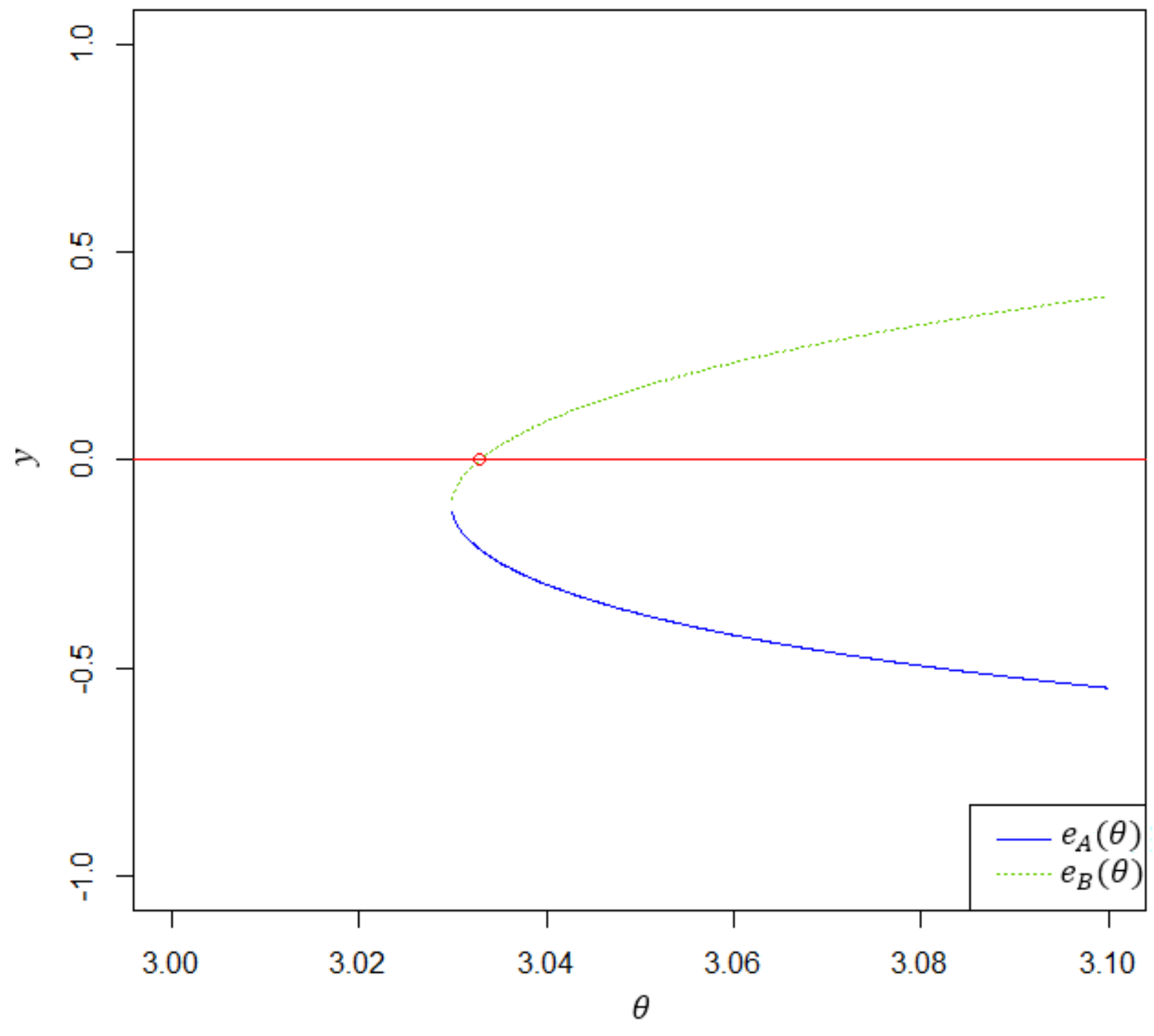

**Fig. 3.11.** Region of the proximity map for polynomial $p(z) = z^3 + C_1 z^2 + C_2 z + C_3$, with coefficients $C_1$, $C_2$, $C_3$ given by (3.14). This region was generated from $N = 1{,}000$ elements $LzC(\theta_k)$ associated with points $\theta_k$ in a regular partition of interval [3.0,3.1).

The region shown in figure 3.11 contains a smooth crossing of function $\hat{e}_B(\theta)$ with horizontal axis $y = 0$, to which corresponds the following numerical estimates shown in (3.16):

$$\hat{R} = -0.6302350 + 0.4047479i \qquad \hat{\theta}^* = 3.032869, \qquad |\Delta e| = 0.001870483, \quad (3.16)$$

Comparing estimates (3.16) with corresponding reference values (3.13) by means of absolute relative errors, we obtain:

$$\left|\frac{R_1 - \hat{R}}{R_1}\right| = 7.58328 \times 10^{-7} \qquad \left|\frac{\theta_1^* - \hat{\theta}^*}{\theta_1^*}\right| = 6.052239 \times 10^{-8} \qquad (3.17)$$

Clearly, relative errors (3.17) are much smaller (by 4 orders of magnitude) compared to the relative errors in (3.15) obtained thanks to the additional routine. Note from figure 3.11, that there is still a gap between $\hat{e}_A$ and $\hat{e}_B$, but it cannot longer contain a crossing with horizontal axis $y = 0$ (i.e., $\hat{e}_A$ and $\hat{e}_B$ have not opposite signs at this gap), and is much smaller compared to the gaps observed in figure 3.10; all this of course thanks to the resolution used, which in this case is 10,000 elements/rad ($63 \times$ more resolution with respect to the global map in figure 3.10). In conclusion, this demonstrates the usefulness of exploring specific regions of a global map in greater detail, especially regions centered around gaps where $\hat{e}_A\hat{e}_B < 0$, indicating probable undetected crossings of map $e(\theta)$ with horizontal axis $y = 0$.

## Numerical Example 3.4

In this example we will analyze the case of a proximity map associated to a cubic polynomial in one complex variable with a root of multiplicity 2. Let us now consider the polynomial equation

$$z^3 + (-5 - 8i)z^2 + (3 + 22i)z + (-19 - 62i) = 0. \qquad (3.18)$$

As a reference, we have that the roots $R_1$, $R_2$, $R_3$ of equation (3.18), together with their corresponding direction angles $\theta_1^*$, $\theta_2^*$, $\theta_3^*$ for lines $\ell_1$ (with fixed point $P_1 = 2.5 + 4i$) containing each of these roots, are:

$$\begin{aligned} R_1 &= 2 + 5i & \theta_1^* &= 2.034444 rad \\ R_2 &= 2 + 5i & \theta_2^* &= 2.034444 rad \\ R_3 &= 1 - 2i & \theta_3^* &= -1.815775 rad \end{aligned} \qquad (3.19)$$

We see that, in this case, $R_1 = R_2$; that is to say, equation (3.18) has a root of multiplicity 2. The proximity map associated with equation (3.18), constructed with $N = 1{,}000$ elements $LzC(\theta_k)$ associated to points $\theta_k$ in a regular partition of interval $[-\pi, \pi)$, by means of the script in annex 3 section 2, is shown in figure 3.12 (where reference values $\theta_i^*$ in (3.19) are also plotted by means of small circles on horizontal axis $y = 0$); the corresponding estimates obtained from the smooth crossings of functions $\hat{e}_A(\theta)$, $\hat{e}_B(\theta)$ with horizontal axis $y = 0$ within this map are listed in table 3.5.

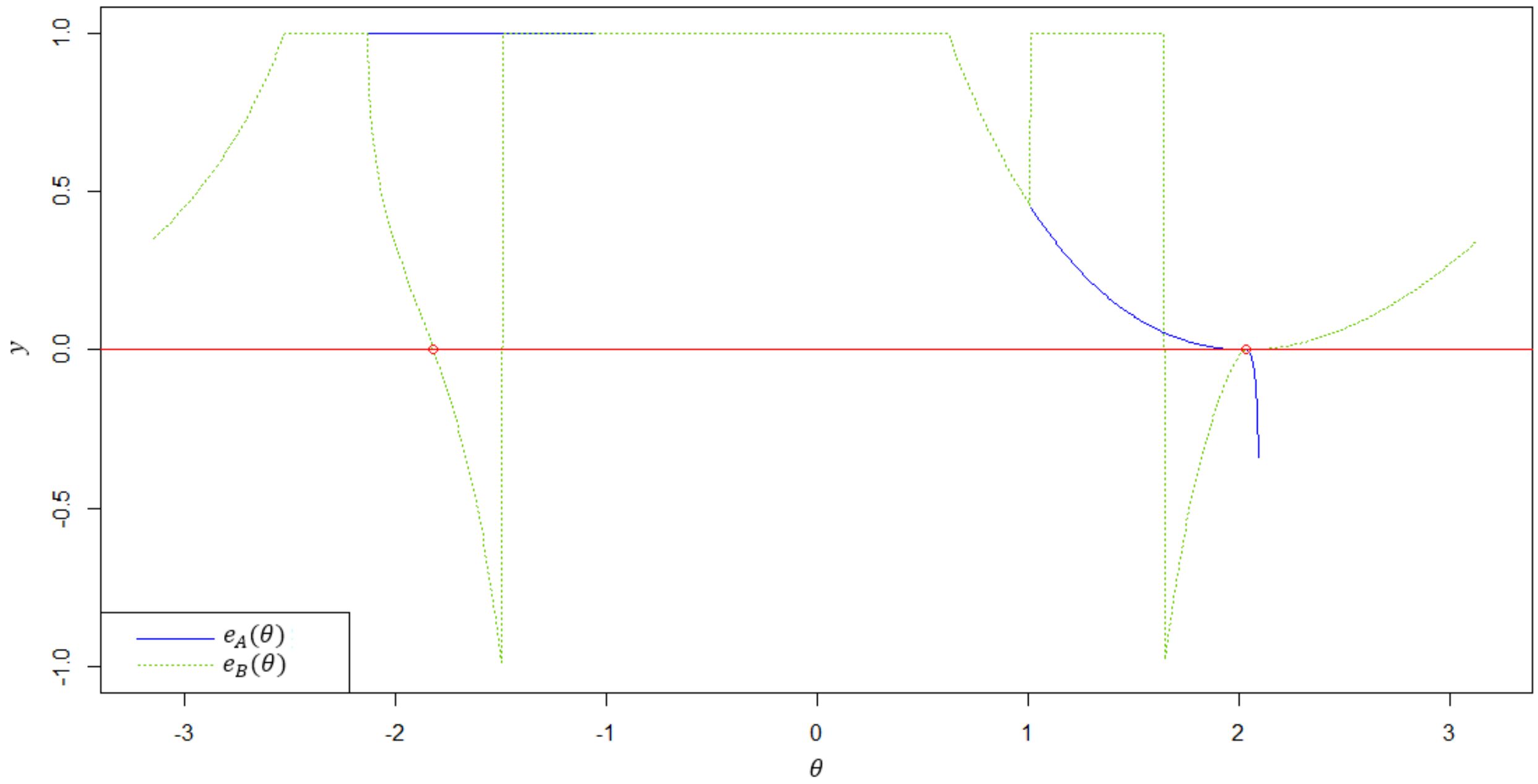

**Fig. 3.12.** Discrete angular proximity map for $p(z) = z^3 + (-5 - 8i)z^2 + (3 + 22i)z + (-19 - 62i)$. This map was generated from $N = 1{,}000$ elements $LzC(\theta_k)$ associated to points $\theta_k$ in a regular partition of interval $[-\pi, \pi)$.

**Table 3.5**. Initial estimates for the roots of $p(z) = z^3 + (-5 - 8i)z^2 + (3 + 22i)z + (-19 - 62i)$, obtained from the map in figure 3.12.

| $i$ | $\hat{R}_i$ | $\hat{\theta}_i^*$ | $\lvert\Delta e_i\rvert$ |
|---|---|---|---|
| 1 | $2.004145 + 5.001474i$ | $2.030536$ | $4.855191 \times 10^{-5}$ |
| 2 | $1.998621 + 4.999504i$ | $2.035745$ | $5.410481 \times 10^{-4}$ |
| 3 | $0.999998 - 2.000000i$ | $-1.815775$ | $1.184694 \times 10^{-2}$ |

As we can see from figure 3.12 and table 3.5, in this case the discrete angular proximity map was able to reasonably estimate all the roots of equation (3.18); in particular, both functions $\hat{e}_A$, $\hat{e}_B$ smoothly cross horizontal axis $y = 0$ at values close to $\theta_1^* = \theta_2^* = 2.034444$; this makes us suspect that, in the structure $LzC$ associated with the polynomial in equation (3.18), just at $\theta_1^*$, a tangential intersection between $tL$ and $zC$ occurs. With the help of the R script in annex 4 section 1, it is possible to arrive at a graph similar to that shown in figure 3.13, which reinforces our hypothesis of a tangential intersection between $tL(\theta_1^*)$ and $zC(\theta_1^*)$.

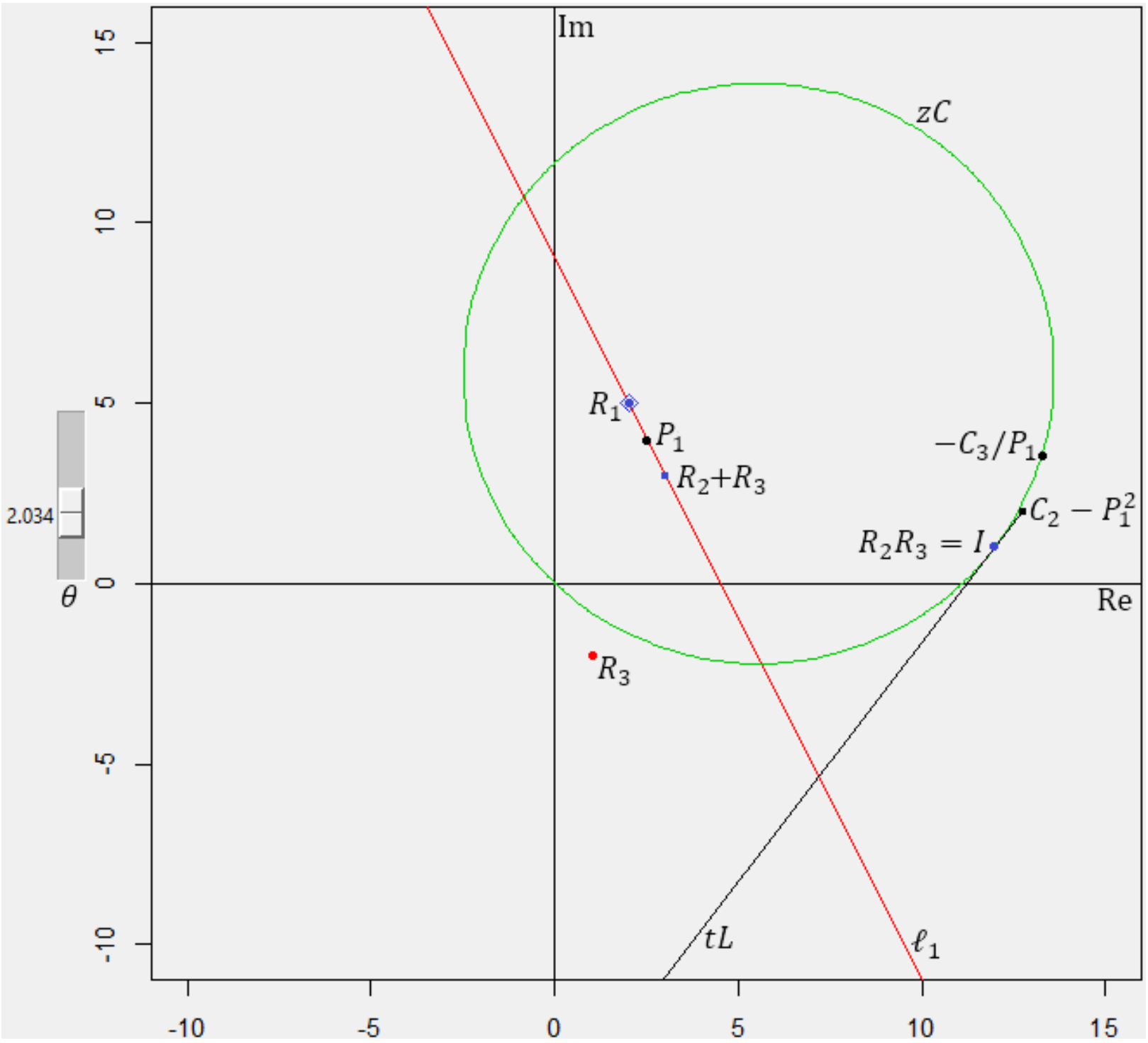

**Fig. 3.13.** Structure $LzC(\theta = 2.034)$ associated with $p(z) = z^3 + (-5-8i)z^2 + (3+22i)z + (-19-62i)$, which contains multiplicity-two root $R_1 = R_2$ in $\ell_1$. $tL$ and $zC$ intersect tangentially at a single point $I = R_2R_3$; it is possible to verify numerically that the vector that goes from the center of $zC$ and points towards $I$, in this case is orthogonal to $v_\theta^2$ ($tL$'s direction vector).

In the case of polynomial equation (3.18), although we have a tangential intersection between terminal semi-line $tL$ and z-circumference $zC$, there is no "gap" in the discrete angular proximity map, unlike discrete maps in examples 3.2 and 3.3. What we observe, instead, is a very small rate of change of functions $\hat{e}_A$ and $\hat{e}_B$ with respect to $\theta$ in regions close to $\theta_1^* = \theta_2^* = 2.034444$; this is confirmed by observing that values $|\Delta e_1|$ and $|\Delta e_2|$ in table 3.5 are smaller, by two or three orders of magnitude, compared to the other value $|\Delta e_3|$ in table 3.5, which seems to be a more "typical" value if we compare it, for example, against values $|\Delta e_i|$ in table 3.2 from example 3.2. This reduction in the rate of change $\Delta e/\Delta\theta$ is caused by the presence of the multiplicity-two root, and has as a cost an increase in the error of initial estimates, as we can see in expressions (3.20) for the absolute relative errors between the initial estimates of table 3.5 with respect to reference values (3.19):

$$\left|\frac{R_1-\hat{R}_1}{R_1}\right| = 8.169451 \times 10^{-4} \qquad \left|\frac{\theta_1^*-\hat{\theta}_1^*}{\theta_1^*}\right| = 1.920697 \times 10^{-3}$$
$$\left|\frac{R_2-\hat{R}_2}{R_2}\right| = 2.721549 \times 10^{-4} \qquad \left|\frac{\theta_2^*-\hat{\theta}_2^*}{\theta_2^*}\right| = 6.397269 \times 10^{-4} \qquad (3.20)$$
$$\left|\frac{R_3-\hat{R}_3}{R_3}\right| = 8.392409 \times 10^{-7} \qquad \left|\frac{\theta_3^*-\hat{\theta}_3^*}{\theta_3^*}\right| = 1.668354 \times 10^{-7}$$

From expressions (3.20), we see that the relative errors associated with the initial estimates for the multiplicity-2 root are larger, by three or four orders of magnitude, compared to the relative errors associated with the initial estimate for the multiplicity-1 root. Even so, the initial approximations listed in table 3.5 can be considered as quite reasonable, since the real and imaginary parts of $\hat{R}_1$ and $\hat{R}_2$, rounded to two decimal digits, coincide with reference roots $R_1$ and $R_2$ in (3.19).

We leave as an exercise for the reader to analyze, with the help of the script in annex 3 section 2, and also with the help of interactive script in annex 4 section 1, the case of a univariate polynomial of degree 3 where two of its roots, although different from each other, are contained in the same line $\ell_1$ with fixed point $-C_1/2$, so the two intersections between $tL$ and $zC$ correspond to those two roots; a particular example where this situation occurs is the polynomial

$$p(z) = z^3 + (-2 - 4i)z^2 + (0 - 19i)z + (-69 + 33i),$$

whose roots are $R_1 = 2 + 3i$, $R_2 = 3 + 4i$, and $R_3 = -3 - 3i$. In this case, the proximity map will show a theta root of multiplicity 2 ($\hat{e}_A$ and $\hat{e}_B$ both crossing horizontal axis $y = 0$ at the same value $\hat{\theta}_i^*$, but with rates of change $\Delta e/\Delta\theta$ having opposite signs), although, geometrically, this theta root corresponds to two different roots; the reader will be able to verify that there is no penalty on the accuracy of computed initial estimates, since corresponding rates of change $\Delta e/\Delta\theta$ are not too small in any case.

## Numerical Example 3.5

In this last example, we will analyze the case of a proximity map associated to a cubic univariate polynomial which has a single root of multiplicity 3. This time, the polynomial equation is

$$z^3 + (-6 - 9i)z^2 + (-15 + 36i)z + (46 - 9i) = 0. \qquad (3.21)$$

The multiplicity-3 root for equation (3.21), and the direction angle for line $\ell_1$ (with fixed point $P_1 = 3 + 4.5i$) containing this root, is:

$$R_1 = R_2 = R_3 = 2 + 3i \qquad \theta_1^* = \theta_2^* = \theta_3^* = -2.158799rad \qquad (3.22)$$

The global proximity map, constructed with $N = 1{,}000$ structures $LzC(\theta_k)$ associated to points $\theta_k = -\pi + 2\pi k/N$ $(k = 0,1,2,\dots,N-1)$, generated of course by means of the script in annex 3 section 2, is shown in figure 3.14:

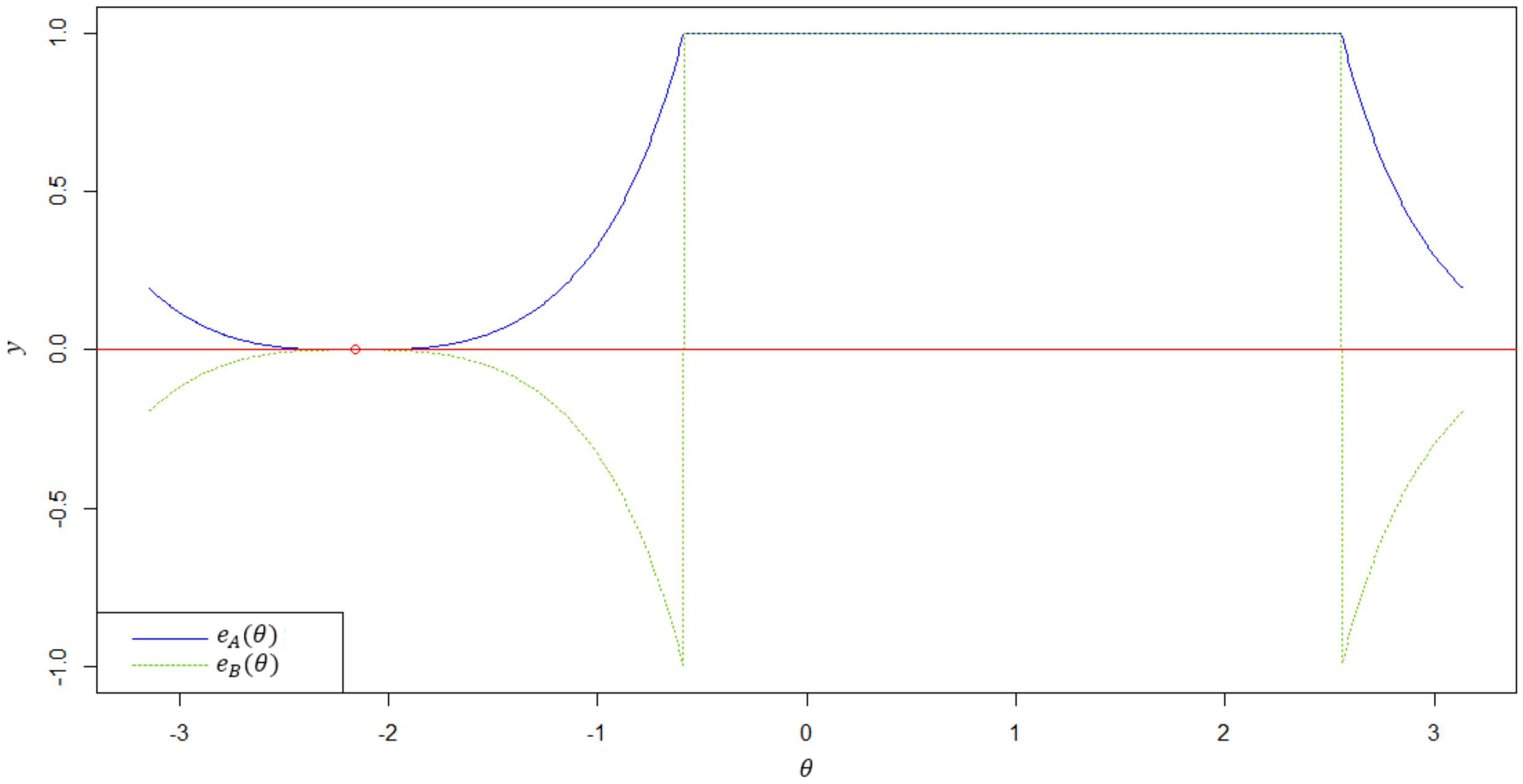

**Fig. 3.14.** Discrete angular proximity map for $p(z) = z^3 + (-6 - 9i)z^2 + (-15 + 36i)z + (46 - 9i)$. This map was generated from $N = 1{,}000$ elements $LzC(\theta_k)$ associated to points $\theta_k$ in a regular partition of interval $[-\pi, \pi)$.

This time, the estimated functions $\hat{e}_A$ and $\hat{e}_B$ in the proximity map of figure 3.14 do not cross in a smooth way the horizontal axis $y = 0$, although both are extremely close to this axis, but without touching it, in the vicinity of point $\theta_1^* = -2.158799$, which is also shown in figure 3.14 by means of a small circle on the horizontal axis $y = 0$. From here we would think that, in theory, function $e_A$ should touch at a single point $\theta_1^*$ the horizontal axis $y = 0$, and in a neighborhood around $\theta_1^*$, $e_A$ should be concave upward. In a completely analogous way, we would also assume that, theoretically, function $e_B$ should touch the horizontal axis $y = 0$ at a single point $\theta_1^*$, and in a neighborhood around $\theta_1^*$, $e_B$ should be concave downward.

We see then that the fact of discretizing an angular proximity map implies a practical impossibility when detecting intersections of estimated function $\hat{e}$ ($\hat{e}$ can be either $\hat{e}_A$ or $\hat{e}_B$) with horizontal axis $y = 0$, if there are no sign changes in $\hat{e}$. Fortunately, the additional routine described within numerical example 3.3 helps us approximate the point, within the discrete partition of $[-\pi, \pi)$, where the minimum vertical distance between $\hat{e}_A$ and $\hat{e}_B$ occurs, imposing the restriction that these functions are on opposing sides with respect to the horizontal axis $y = 0$. Applying this additional routine to the map in figure 3.14, we obtain the estimate shown in table 3.6.

**Table 3.6**. Initial estimation of the triple root of $p(z) = z^3 + (-6 - 9i)z^2 + (-15 + 36i)z + (46 - 9i)$, obtained from the additional routine (defined in example 3.3) applied to the map in figure 3.14.

| $i$ | $\hat{R}_i$ | $\hat{\theta}_i^*$ | $\lvert\Delta e_i\rvert$ |
|---|---|---|---|
| 1 | $1.99456 + 3.000358i$ | $-2.161416$ | $6.897105 \times 10^{-9}$ |

From table 3.6, we can see that the initial estimates $\hat{R}_1$, $\hat{\theta}_1^*$ are reasonably close to reference values (3.22). In this case, value $|\Delta e_1|$ in table 3.6 is the minimum vertical distance between estimated functions $\hat{e}_A$, $\hat{e}_B$ within the map in figure 3.14, where we can also observe that both estimated functions $\hat{e}_A$, $\hat{e}_B$ experience very low rates of change $\Delta e/\Delta\theta$ in the vicinity of triple root $\theta_1^*$. The reader is invited to experiment with the script in annex 3 section 2 in order to obtain closer views centered on the region around $\theta_1^*$, using of course the coefficients of equation (3.21), and changing the limits of vertical axis on the graph (bringing them closer to the origin) in order to better appreciate the behavior of estimated functions $\hat{e}_A$, $\hat{e}_B$.

The absolute relative errors associated with the estimate listed in table 3.6, with respect to reference values (3.22), are:

$$\left|\frac{R_1-\hat{R}_1}{R_1}\right| = 1.51196 \times 10^{-3} \qquad \left|\frac{\theta_1^*-\hat{\theta}_1^*}{\theta_1^*}\right| = 1.212163 \times 10^{-3} \qquad (3.23)$$

Relative errors (3.23) are larger than "typical" relative errors (observed in previous numerical examples in this chapter) associated with multiplicity-1 root estimates; this is partly due to the presence of multiplicity-3 theta root $\theta_1^*$ (linked to very low rates of change $\Delta e/\Delta\theta$ near $\theta_1^*$), and partly to the characteristics of the additional routine defined in numerical example 3.3. Even so, we consider that estimates $\hat{R}_1$, $\hat{\theta}_1^*$ are good initial approximations.

To conclude this example, we will see how the structure $LzC$ behaves at the vicinity of $\theta_1^*$. With the help of the script in annex 4 section 1, we can generate a graph similar to the one shown in figure 3.15.

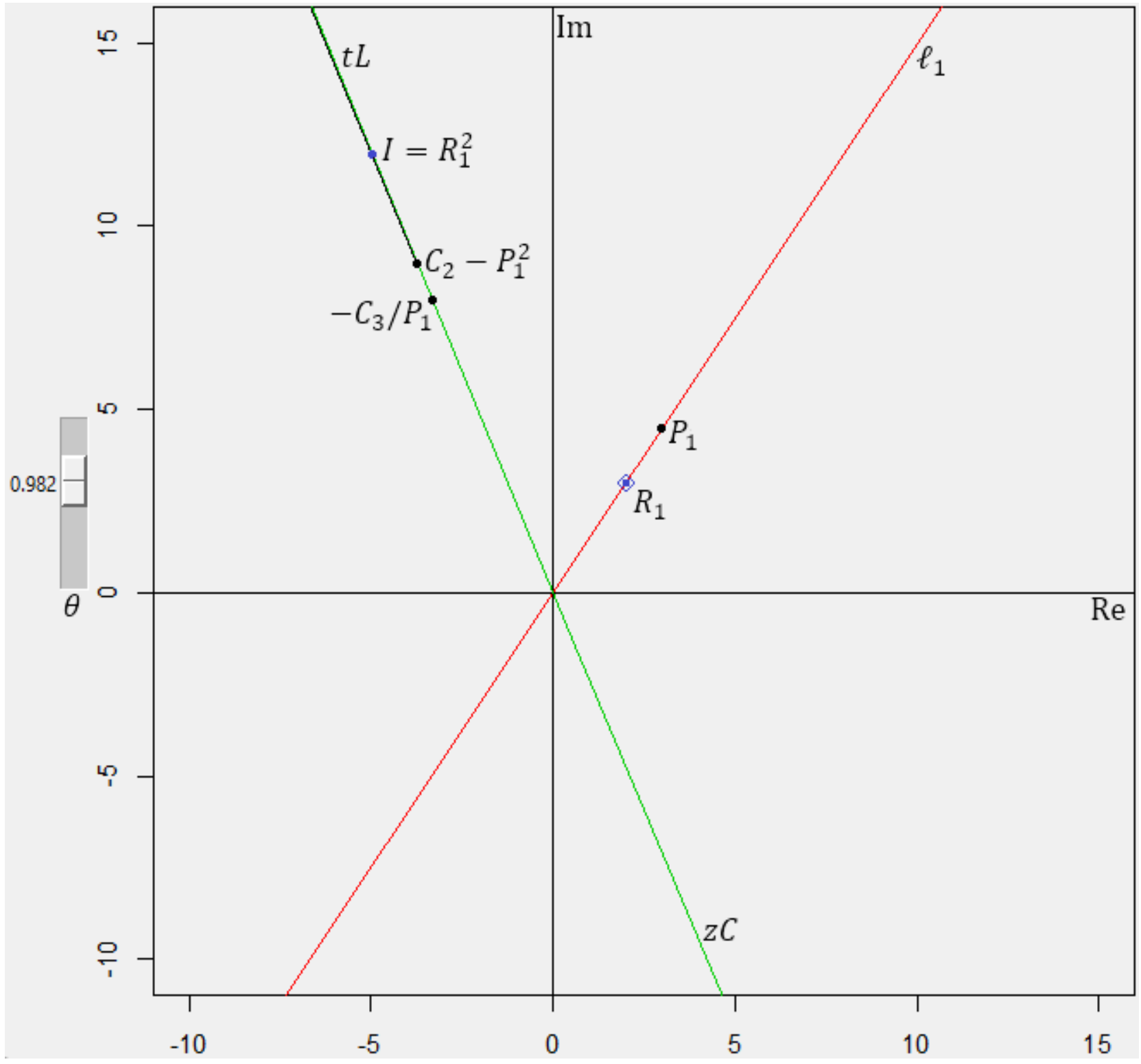

**Fig. 3.15.** Structure $LzC(\theta = 0.982)$ associated with $p(z) = z^3 + (-6 - 9i)z^2 + (-15 + 36i)z + (46 - 9i)$, which is extremely close to containing triple root $R_1$ in $\ell_1$. Because $\ell_1$ is extremely close to point $0 + 0i$, $zC$ has an extremely large radius, which makes it look more like a line rather than a circumference. $tL$ and $zC$ practically overlap along a common line, which contains point $I = R_1^2$, and fixed points of both $tL$ and $zC$.

In figure 3.15, inclination angle $\theta = 0.982$ for line $\ell_1$ contained in the structure $LzC$, is such that by subtracting $\pi$ radians from it, we obtain a value close to reference angle $\theta_1^* \cong -2.16$ in (3.22); of course, this reference angle also corresponds to the direction of line $\ell_1$ shown in figure 3.15. We see that z-circumference $zC$ tends to become a line, because $\ell_1$ is very close to the origin $0 + 0i$. We can also see that terminal semi-line $tL$ tends to be contained in $zC$. We see that the structure $LzC$ shown in figure 3.15 practically has root $R_1$ on $\ell_1$, since $c_0$ (the diamond-shaped point $\diamond$) and $i_0$ (the point • inside the diamond $\diamond$) are practically at the same place.

# Chapter Four: LC Method for Quartic Polynomial Equations

## Theoretical Aspects

In this chapter we will extend the LC method towards polynomials of degree four in a complex variable, preserving as much as possible the notation from previous chapters.

The polynomial equation of degree 4, whose roots $R_1$, $R_2$, $R_3$, $R_4$ we seek to approximate is of the form

$$z^4 + C_1 z^3 + C_2 z^2 + C_3 z + C_4 = 0. \qquad (4.1)$$

In this case, Vieta's relations between the roots and coefficients of equation (4.1) establish that

$$C_1 = -(R_1 + R_2 + R_3 + R_4), \qquad (4.2\text{-a})$$

$$C_2 = R_1R_2 + R_1R_3 + R_1R_4 + R_2R_3 + R_2R_4 + R_3R_4, \qquad (4.2\text{-b})$$

$$C_3 = -(R_1R_2R_3 + R_1R_2R_4 + R_1R_3R_4 + R_2R_3R_4), \qquad (4.2\text{-c})$$

$$C_4 = R_1R_2R_3R_4. \qquad (4.2\text{-d})$$

Lodovico Ferrari was the first to develop an analytical solution to find the roots of equation (4.1); the reader interested in knowing the details of this solution can consult [4.1]. An impractical, but direct, formula for obtaining the roots of equation (4.1) from its coefficients, without using intermediate variables or associated polynomials of lower degree, can be seen in [4.2]. In fact, there are many ways to find, by means of algebraic techniques, the roots of quartic equation (4.1); most of these techniques involve finding the roots of a degree 3 polynomial associated with equation (4.1), called *resolvent cubic*. In annex 5 section 1, we describe an algebraic technique to solve equation (4.1) by means of a resolvent cubic and associated quadratic equations.

Similarly to how we proceeded in chapter 3, for the case of polynomial equations of degree 4 in a complex variable of the form (4.1), we will seek to take advantage of existing geometric relationships between the coefficients of equation (4.1) and lines and circumferences contained in the complex plane $\mathbb{C}$, with the purpose of obtaining initial approximations to the solutions of the non-linear system of four equations (4.2-a), (4.2-b), (4.2-c), (4.2-d) in four unknowns $R_1$, $R_2$, $R_3$, $R_4$ by means of discrete proximity maps. We will assume that the roots of equation (4.1) are different from each other, and all are different from zero; this time we will not analyze examples with repeated roots; we will focus on cases in which the roots are randomly distributed in the complex plane $\mathbb{C}$.

We start with the construction of a structure $LzC(C_1, C_2, C_3, C_4, \theta^*)$ associated with equation (4.1), which we will assume (without loss of generality) that contains the root $R_1$ at some to-be-determined point inside line $\ell_1 \subset LzC(C_1, C_2, C_3, C_4, \theta^*)$ with fixed point $P_1 = -C_1/2$ and unit direction vector $v_{\theta^*} = e^{i\theta^*} = \cos\theta^* + i\sin\theta^*$, which starts at $P_1$ and points towards $R_1$. In order to geometrically determine the location of $R_1$, we will seek to obtain, by means of geometric elements directly constructible from coefficients $C_1, C_2, C_3$, a terminal semi-line $tL$ containing the product $R_2R_3R_4$, which of course is also contained in the z-circumference $zC: C_4/\ell_1$. Starting with line $\ell_1: P_1 + tv_{\theta^*}$, $t \in \mathbb{R}$, we obtain the following derived semi-lines, relying on Vieta's relations (4.2-a), (4.2-b), (4.2-c):

**(4.3-a)** $\boldsymbol{\ell_d}$: $P_1^2 + t^2(-v_{\theta^*}^2)$,
$R_1R_2 + R_1R_3 + R_1R_4 \in \ell_d$;

**(4.3-b)** $\boldsymbol{C_2 - \ell_d}$: $C_2 - P_1^2 + t^2 v_{\theta^*}^2$,
$R_2R_3 + R_2R_4 + R_3R_4 \in C_2 - \ell_d$;

**(4.3-c)** $R_1(\boldsymbol{C_2 - \ell_d})$: $R_1C_2 - R_1P_1^2 + t^2 R_1 v_{\theta^*}^2$,
$R_1R_2R_3 + R_1R_2R_4 + R_1R_3R_4 \in R_1(C_2 - \ell_d)$;

**(4.3-d)** $-\boldsymbol{C_3} - R_1(\boldsymbol{C_2 - \ell_d})$: $-C_3 - R_1C_2 + R_1P_1^2 + t^2(-R_1 v_{\theta^*}^2)$,
$R_2R_3R_4 \in -C_3 - R_1(C_2 - \ell_d)$.

In principle, expression (4.3-d), i.e., $-C_3 - R_1(C_2 - \ell_d)$, is the one that would help us to build the terminal semi-line $tL$. In expressions (4.3-a) to (4.3-d), the terms which do not contain parameter $t^2$, for example $R_1C_2 - R_1P_1^2$ in expression (4.3-c), define the *anchor point* of corresponding derived semi-line, while the factor associated with parameter $t^2$, for example $R_1v_{\theta^*}^2$, also in expression (4.3-c), represents a direction vector for the corresponding derived semi-line. Figure 4.1 schematically shows, within the complex plane $\mathbb{C}$, some geometric elements in structure $LzC(C_1, C_2, C_3, C_4, \theta^*)$ generated by expressions (4.3-a) to (4.3-d).

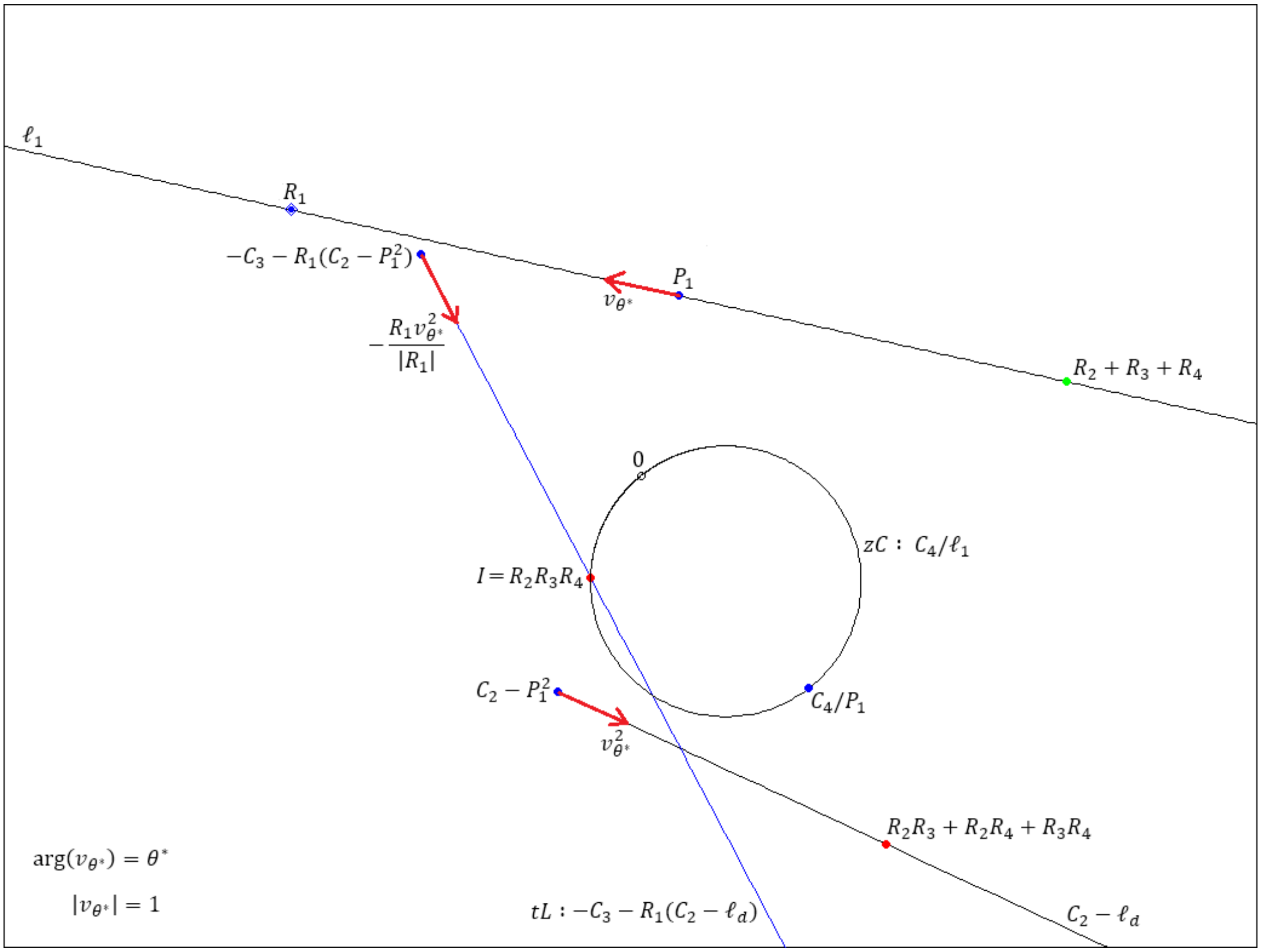

**Fig. 4.1**. If line $\ell_1 \subset \mathbb{C}\backslash\{0\}$ with fixed point $P_1 = -C_1/2$ contains the root $R_1$ of polynomial equation (4.1), then z-circumference $zC: C_4/\ell_1$ contains the product $R_2R_3R_4$; the latter is located at the intersection $I$ between $zC$ and terminal semi-line $tL: -C_3 - R_1(C_2 - \ell_d)$, where $\ell_d$ is the semi-line derived from $\ell_1$. All direction vectors shown in this schematic graph are unit vectors.

So far, everything looks great, except for one thing: there is a serious difficulty with expressions (4.3-c) and (4.3-d), since we are directly plugging the not-yet-determined value $R_1$ into these expressions to construct the terminal semi-line $tL(\theta^*)$. How can we obtain the anchor point and direction vector of $tL(\theta^*)$, without directly involving the unknown value $R_1$? A possible workaround to this problem is to use, within expressions (4.3-c) and (4.3-d), instead of $R_1$, the known factor $\ell_1$, which by assumption contains $R_1$. In this way, we obtain the following expressions:

| | | |
|---|---|---|
| **(4.3-c')** | $\boldsymbol{\ell_1(C_2 - \ell_d)}$: | $C_2(P_1 + tv_{\theta^*}) - (P_1 - tv_{\theta^*})(P_1 + tv_{\theta^*})^2$,<br>$R_1R_2R_3 + R_1R_2R_4 + R_1R_3R_4 \in \ell_1(C_2 - \ell_d)$; |
| **(4.3-d')** | $\boldsymbol{-C_3 - \ell_1(C_2 - \ell_d)}$: | $-C_3 - C_2(P_1 + tv_{\theta^*}) + (P_1 - tv_{\theta^*})(P_1 + tv_{\theta^*})^2$,<br>$R_2R_3R_4 \in -C_3 - \ell_1(C_2 - \ell_d)$. |

Now we see that expressions (4.3-c') and (4.3-d') no longer explicitly contain the unknown root $R_1$, but are constructed from known geometric elements within structure $LzC(C_1, C_2, C_3, C_4, \theta^*)$; the cost of this adaptation is that now $\ell_1(C_2 - \ell_d)$ and $-C_3 - \ell_1(C_2 - \ell_d)$ are no longer semi-lines, but parametric curves: curve (4.3-d') is in fact a reflection of curve (4.3-c') with respect to point $-C_3/2$. Figure 4.2 shows, within the complex plane $\mathbb{C}$, the same z-circumference $zC$ from figure 4.1, together with the corresponding *terminal curve* $t\mathfrak{C}$ defined by expression (4.3-d').

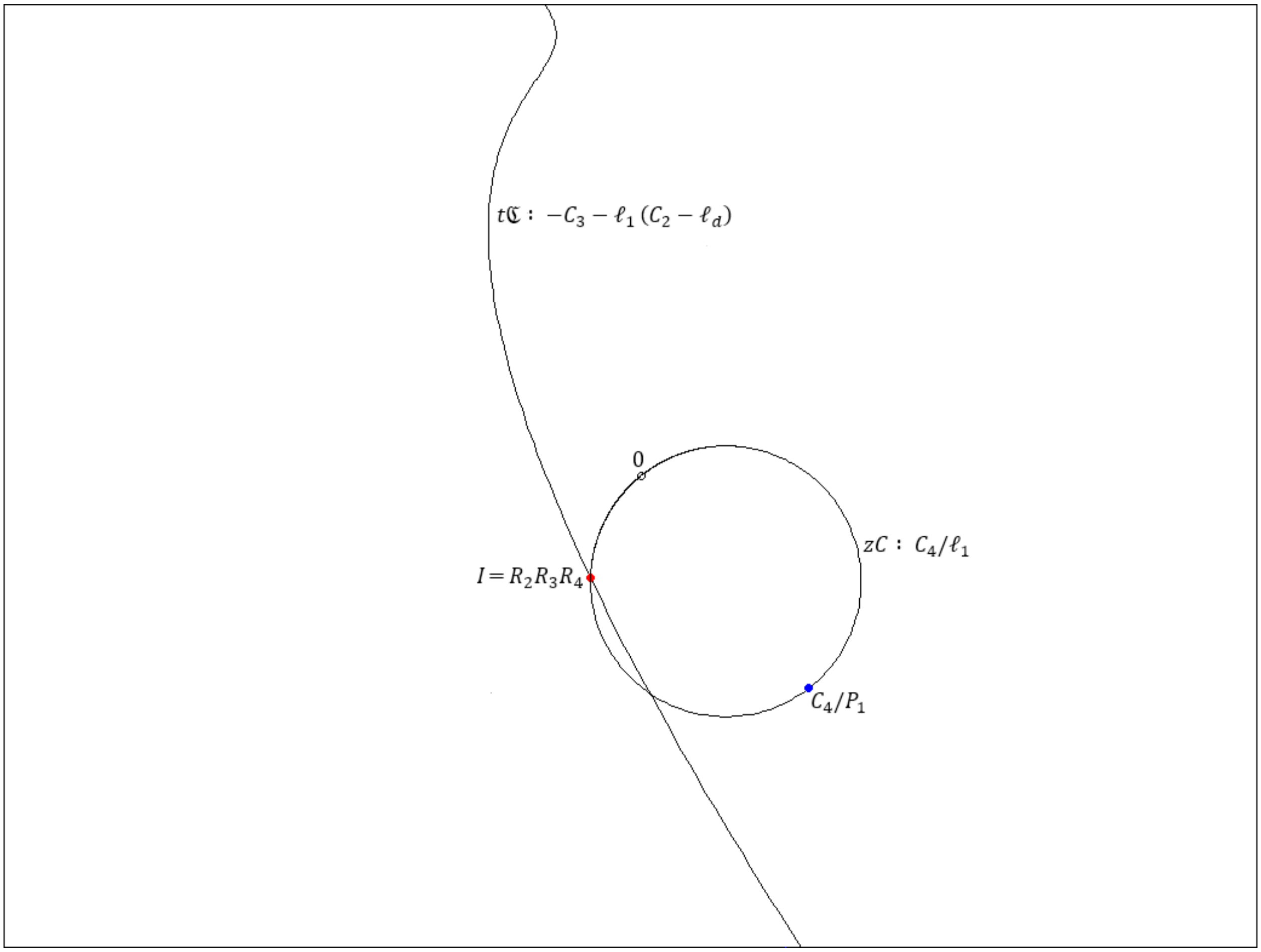

**Fig. 4.2**. z-circumference $zC(\theta^*)$: $C_4/\ell_1(\theta^*)$ from figure 4.1 intersects with corresponding terminal curve $t\mathfrak{C}(\theta^*)$, defined by parametric expression (4.3-d'), at point $I = R_2R_3R_4$.

To find the intersection $I = R_2R_3R_4$ between $t\mathfrak{C}(\theta^*)$ and $zC(\theta^*)$, we can no longer use the function which computes intersections between a semi-line and a circumference (function `intersect_semiLine_circle` in annex 1 section 2); in this case, it will be necessary to use an alternative computational resource. To do this, let's consider the function

$$d_\theta^2(t) = \left[\mathrm{Re}\left(t\mathfrak{C}(\theta,t) - zC(\theta,t)\right)\right]^2 + \left[\mathrm{Im}\left(t\mathfrak{C}(\theta,t) - zC(\theta,t)\right)\right]^2, \qquad (4.4)$$

where $t\mathfrak{C}(\theta,t) = -C_3 - C_2(P_1 + tv_\theta) + (P_1 - tv_\theta)(P_1 + tv_\theta)^2,$

$zC(\theta,t) = C_4/(P_1 + tv_\theta);$

$\theta, t \in \mathbb{R}, v_\theta = e^{i\theta} = \cos\theta + i\sin\theta, P_1 = -C_1/2.$

**We will call function $d_\theta^2(t)$ in (4.4)** ***the dynamic squared distance function between*** $\boldsymbol{t\mathfrak{C}(\theta)}$ ***and*** $\boldsymbol{zC(\theta)}$ (note that this function is defined for any real value $\theta$, not just for $\theta^*$). Figure 4.3 shows an important region of the graph of $d_{\theta^*}^2(t)$ associated with the terminal curve $t\mathfrak{C}(\theta^*)$ and z-circumference $zC(\theta^*)$ from figure 4.2; this region contains the point $t^*$ in the domain of $d_{\theta^*}^2(t)$ where its **global minimum** $d_{\theta^*}^2(t^*)$ occurs.

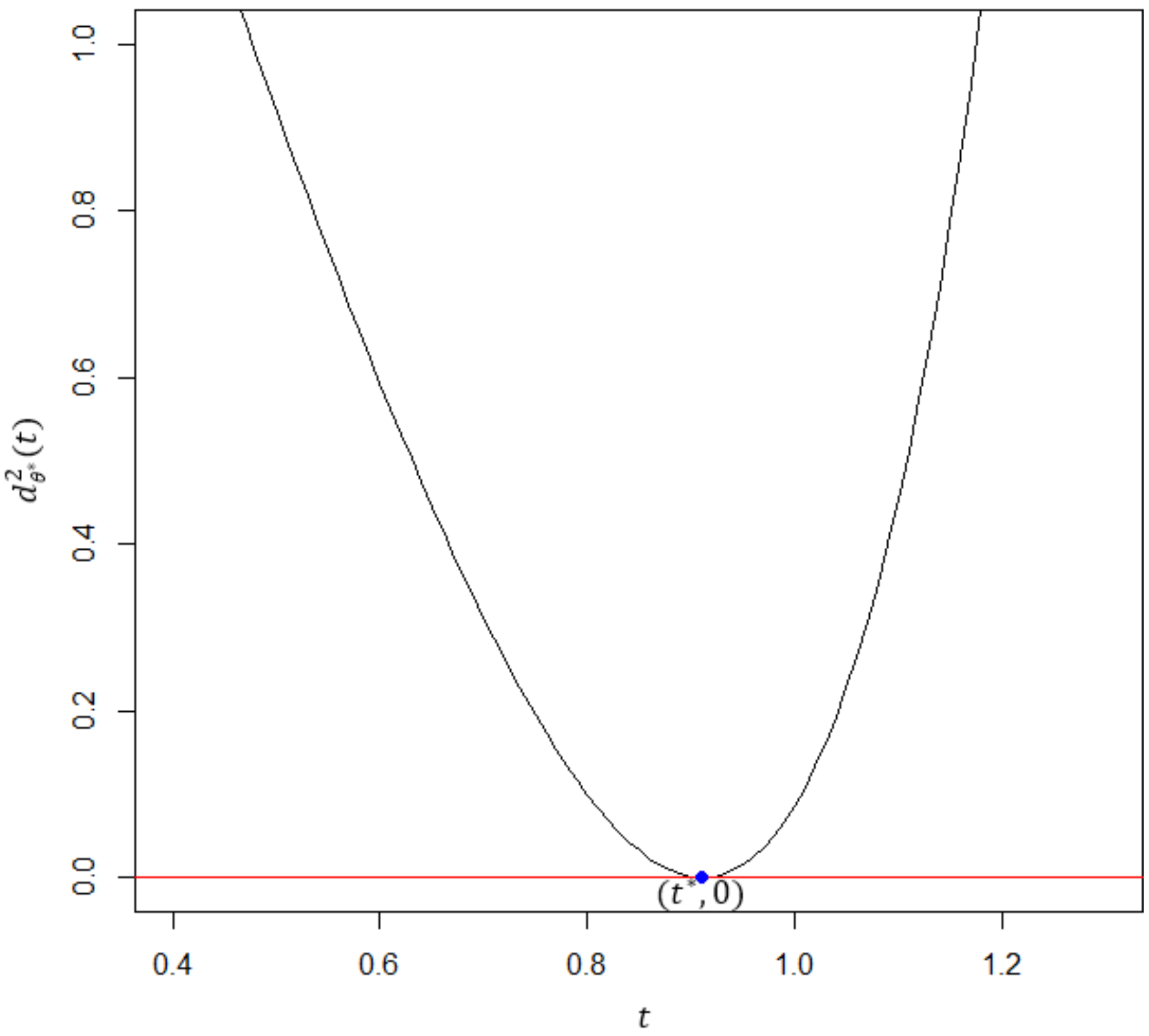

**Fig. 4.3**. Dynamic squared distance function $d_{\theta^*}^2(t)$ associated with z-circumference $zC(\theta^*)$ and terminal curve $t\mathfrak{C}(\theta^*)$ from figure 4.2. $d_{\theta^*}^2(t)$ attains its global minimum value at $t^*$; additionally, $d_{\theta^*}^2(t^*) = 0$.

Note: The parametric value $t^*$ where the global minimum of the dynamic squared distance $d_\theta^2(t)$ defined in (4.4) occurs, coincides with the parametric value where the global minimum of the Euclidean distance $\sqrt{d_\theta^2(t)}$ between $t\mathfrak{C}(\theta,t)$ and $zC(\theta,t)$ occurs, since $\sqrt{x}$ is a monotonically increasing function $\forall x \in \mathbb{R}^+ \cup \{0\}$; therefore, it is valid to use the squared distance (4.4) instead of the conventional Euclidean distance as an objective function within optimization algorithms.

In terms of parameter $t$, the intersection $I = R_2R_3R_4$ between terminal curve $t\mathfrak{C}(\theta^*, t)$ and z-circumference $zC(\theta^*, t)$, shown in figure 4.2, occurs at $t^* = \arg\min_{t\in\mathbb{R}} d^2_{\theta^*}(t)$. In this case, since $I = t\mathfrak{C}(\theta^*, t^*) = zC(\theta^*, t^*)$, we have $d^2_{\theta^*}(t^*) = 0$. This also implies that value $t^*$, applied to line $\ell_1(\theta^*, t)$, produces root $R_1$. Therefore,

$$R_1 = P_1 + \arg\min_{t\in\mathbb{R}} d^2_{\theta^*}(t)\, v_{\theta^*}. \tag{4.5}$$

Under the assumption that we know the inclination angle $\theta^*$ necessary for $\ell_1$ to contain $R_1$, but that we do not know $R_1$ itself, equation (4.5) suggests that it is convenient to first find the *minimizing argument* $t^*$ of function $d^2_{\theta^*}(t)$; only then can we build the terminal semi-line $tL(\theta^*)$ shown in figure 4.1 by means of value $R_1$ given by equation (4.5). Annex 2 section 4 presents an initial proposal on how to estimate $t^*$ by means of optimization functions native to the basic version of the R language; such proposal, as we will see later, works reasonably well for our illustrative purposes. In principle, it would seem redundant to obtain $tL(\theta^*)$ when we have already obtained $R_1$ by means of formula (4.5); however, since in practice we generally do not know $\theta^*$, in such circumstances expression (4.5) becomes an estimator for $R_1$, and $tL$ is used to measure the degree of error of estimate $\hat{R}_1$ by first finding the intersection between $tL$ and $zC$, and then by applying the weighted error function $e(\theta)$ defined in chapter 3.

In summary, the parametric trajectories characterized by expressions (4.3-a) to (4.3-d'), together with $\ell_1(\theta^*)$ and $zC(\theta^*)$, form an integral part of the structure $LzC(C_1, C_2, C_3, C_4, \theta^*)$ associated with equation (4.1); the computation of the minimizing argument $t^*(\theta^*)$ for the dynamic squared distance function $d^2_{\theta^*}(t)$ given by (4.4), allows us to explicitly know all the elements of $LzC(C_1, C_2, C_3, C_4, \theta^*)$. As we saw in chapter 3, in practice it is necessary to consider also neighboring structures around $LzC(C_1, C_2, C_3, C_4, \theta^*)$ to estimate any of the four values $\theta^*$, which are distributed anywhere in the angular interval $[-\pi, \pi)$; this is what we will do next.

**Construction of structures $LzC$ that do not contain a root of equation (4.1)**. In a similar way to how we proceeded in chapter 3, we will now translate our building procedure for the structure $LzC(C_1, C_2, C_3, C_4, \theta^*)$ to an arbitrary structure $LzC(C_1, C_2, C_3, C_4, \theta)$ associated with equation (4.1), where $\theta$ is a value arbitrarily close to $\theta^*$; i.e., $\theta = \theta^* + \Delta\theta$. For this, we only need to substitute, into expressions (4.3-a) to (4.3-d'), the real parameter $\theta$ instead of the angle $\theta^*$. In this way, we obtain an expression for the terminal curve $t\mathfrak{C}(\theta, t)$, based on expression (4.3-d'):

$$t\mathfrak{C}(\theta, t) = -C_3 - C_2(P_1 + tv_\theta) + (P_1 - tv_\theta)(P_1 + tv_\theta)^2 \tag{4.6}$$

In expression (4.6), $v_\theta = e^{i\theta} = \cos\theta + i\sin\theta$ is the unit direction vector for line $\ell_1$ associated with the structure $LzC(C_1, C_2, C_3, C_4, \theta)$. In figure 4.4, we can see terminal curve $t\mathfrak{C}(\theta, t)$, defined by expression (4.6), for a value $\theta$ close to value $\theta^*$ from figure 4.1, along with associated elements $\ell_1(\theta, t) = P_1 + tv_\theta$ and $zC(\theta, t) = C_4/\ell_1(\theta, t)$.

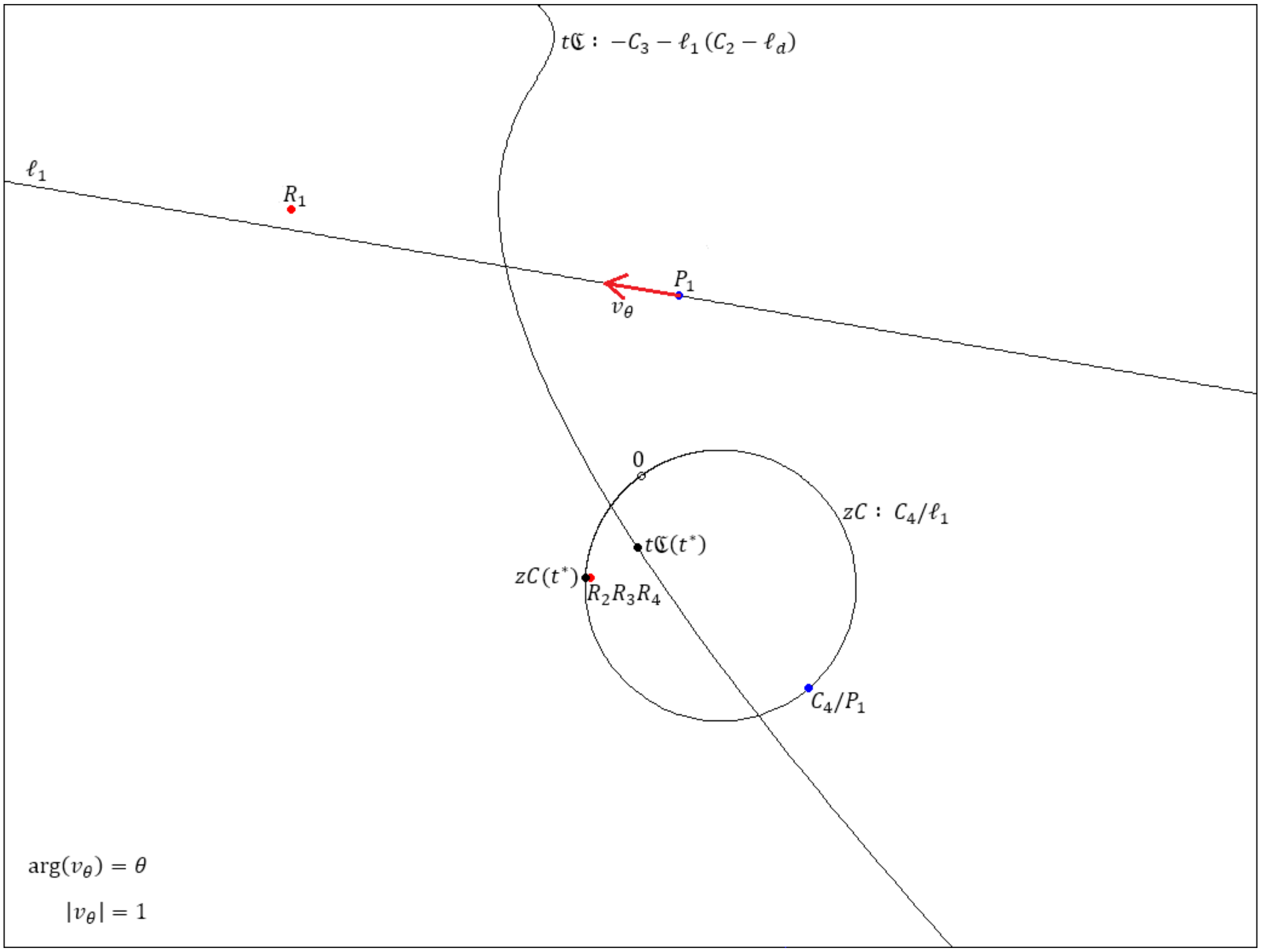

**Fig. 4.4**. The minimum distance between $zC: C_4/\ell_1(\theta, t)$ and $t\mathfrak{C}: -C_3 - \ell_1(\theta, t)[C_2 - \ell_d(\theta, t)]$ is given by the length of the line segment whose ends are the points $zC(\theta, t^*)$ and $t\mathfrak{C}(\theta, t^*)$, both of them shown in the graph. Inclination angle $\theta$ for line $\ell_1(\theta, t)$ in this graph is close to angle $\theta^*$ from figure 4.1.

In figure 4.4, $t^*(\theta) = \arg\min_{t\in\mathbb{R}} d_\theta^2(t)$, where $d_\theta^2(t)$ is the dynamic squared distance function between trajectories $t\mathfrak{C}(\theta, t)$ and $zC(\theta, t)$; see definition (4.4). Note that this time we are no longer interested in the "static" intersection between $t\mathfrak{C}(\theta, t)$ and $zC(\theta, t)$, but in the minimum distance between two points, one on $t\mathfrak{C}(\theta, t)$ and the other on $zC(\theta, t)$, which move along these curves in synchronicity with "time" parameter $t$. As we saw for the case of the structure $LzC(C_1, C_2, C_3, C_4, \theta^*)$, the global minimum of the dynamic squared distance occurs when those two points on curves $t\mathfrak{C}(\theta^*, t)$ and $zC(\theta^*, t)$ coincide spatially, so the "static" intersection between $t\mathfrak{C}(\theta, t)$ and $zC(\theta, t)$ only makes sense when the structure $LzC$ contains one of the roots of equation (4.1). What happens if the structure $LzC$ contains two or more distinct roots of equation (4.1) on $\ell_1$? Function $d^2_{\theta^*}(t)$ will probably have more than one global minimum value, touching horizontal axis $y = 0$ on more than one occasion; this situation, however, will not be addressed in this exposition.

In figure 4.5, we can see part of the graph for the function $d^2_\theta(t)$ associated with the geometric elements $t\mathfrak{C}(\theta, t)$ and $zC(\theta, t)$ from figure 4.4.

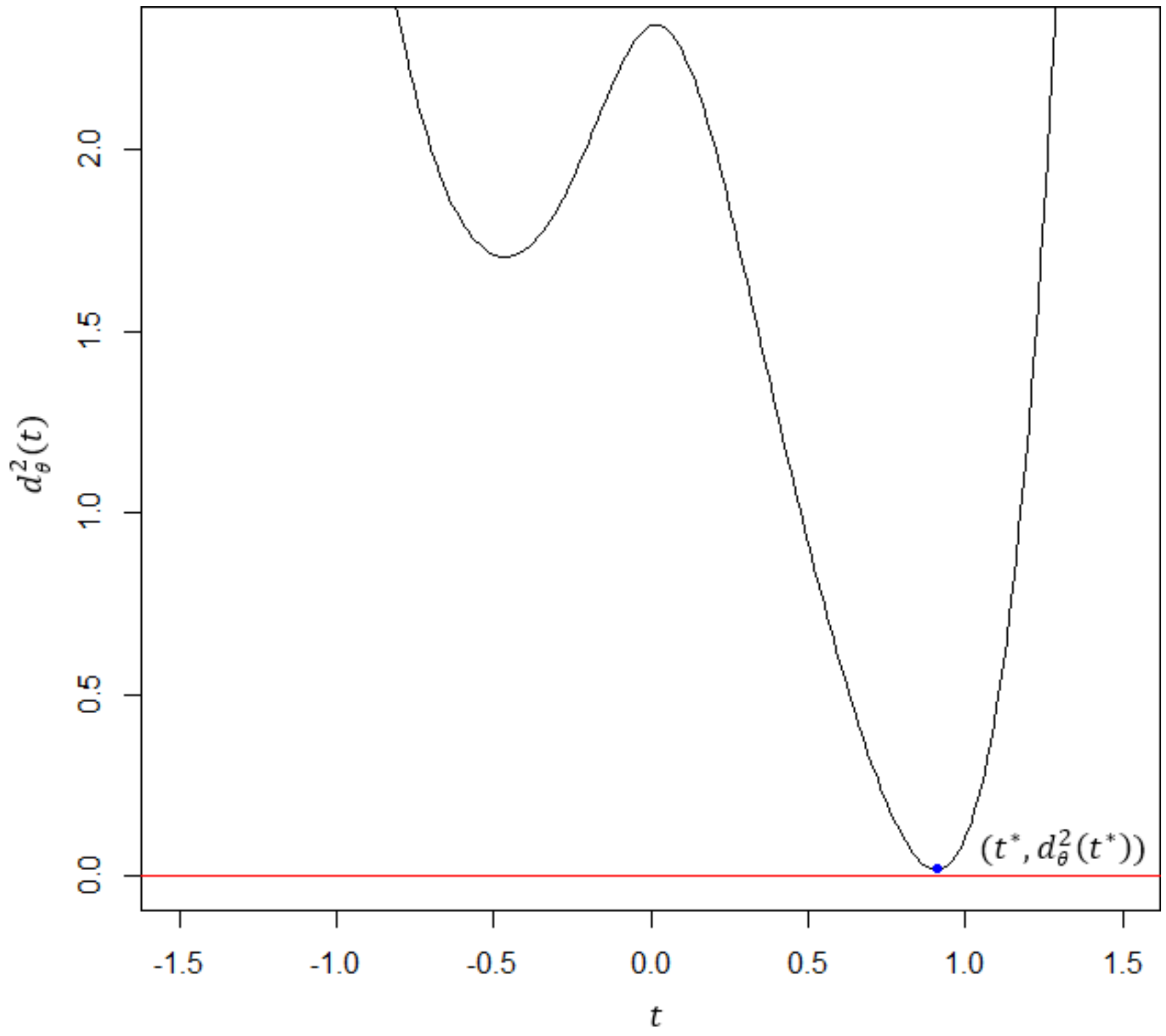

**Fig. 4.5**. Dynamic squared distance function $d^2_\theta(t)$ associated with z-circumference $zC\colon C_4/\ell_1(\theta, t)$ and terminal curve $t\mathfrak{C}\colon -C_3 - \ell_1(\theta, t)[C_2 - \ell_d(\theta, t)]$ from figure 4.4. $d^2_\theta(t)$ attains its global minimum value at $t^*$. Note that $d^2_\theta(t^*) > 0$.

From figure 4.5 we can see that, in general, the function $d^2_\theta(t)$ can contain local minima; we emphasize, however, that our interest focuses on obtaining the argument $t^*(\theta)$ where the global minimum of $d^2_\theta(t)$ occurs, which must additionally meet the requirement $t^*(\theta) > 0$, since we expect to find, on $\ell_1$, one of the roots of equation (4.1) located on the same side to which the direction vector $v_\theta$ points to with respect to fixed point $P_1$; the function `min_D2` in annex 2 section 4 takes into account all these considerations in order to reliably obtain an estimate of $t^*$, although there is certainly much to be done here to improve the optimization process, especially in terms of temporal efficiency.

Under these circumstances, the value $t^*(\theta)$ associated with structure $LzC(C_1, C_2, C_3, C_4, \theta)$ helps us to find an estimate **of one of the roots** of equation (4.1), by means of expression (4.7):

$$R_x(\theta) = P_1 + t^*(\theta)v_\theta. \tag{4.7}$$

Approximation $R_x(\theta)$ in (4.7) now allows us to construct the terminal semi-line $tL(\theta)$, analogous to how it was done in definition (4.3-d) for $tL(\theta^*)$:

$$tL(\theta): -C_3 - R_x(\theta)[C_2 - \ell_d(\theta)] = -C_3 - R_x(\theta)C_2 + R_x(\theta)P_1^2 + t^2[-R_x(\theta)v_\theta^2] \tag{4.8}$$

We can now compute intersections between terminal semi-line $tL(\theta)$ and z-circumference $zC(\theta)$ in the usual way, by using the function which computes intersections between a semi-line and a circumference (function `intersect_semiLine_circle` in annex 1 section 2). In figure 4.6 we can visualize, in a schematic way, some elements related to the structure $LzC(C_1, C_2, C_3, C_4, \theta)$, analogous to those shown in figure 4.1.

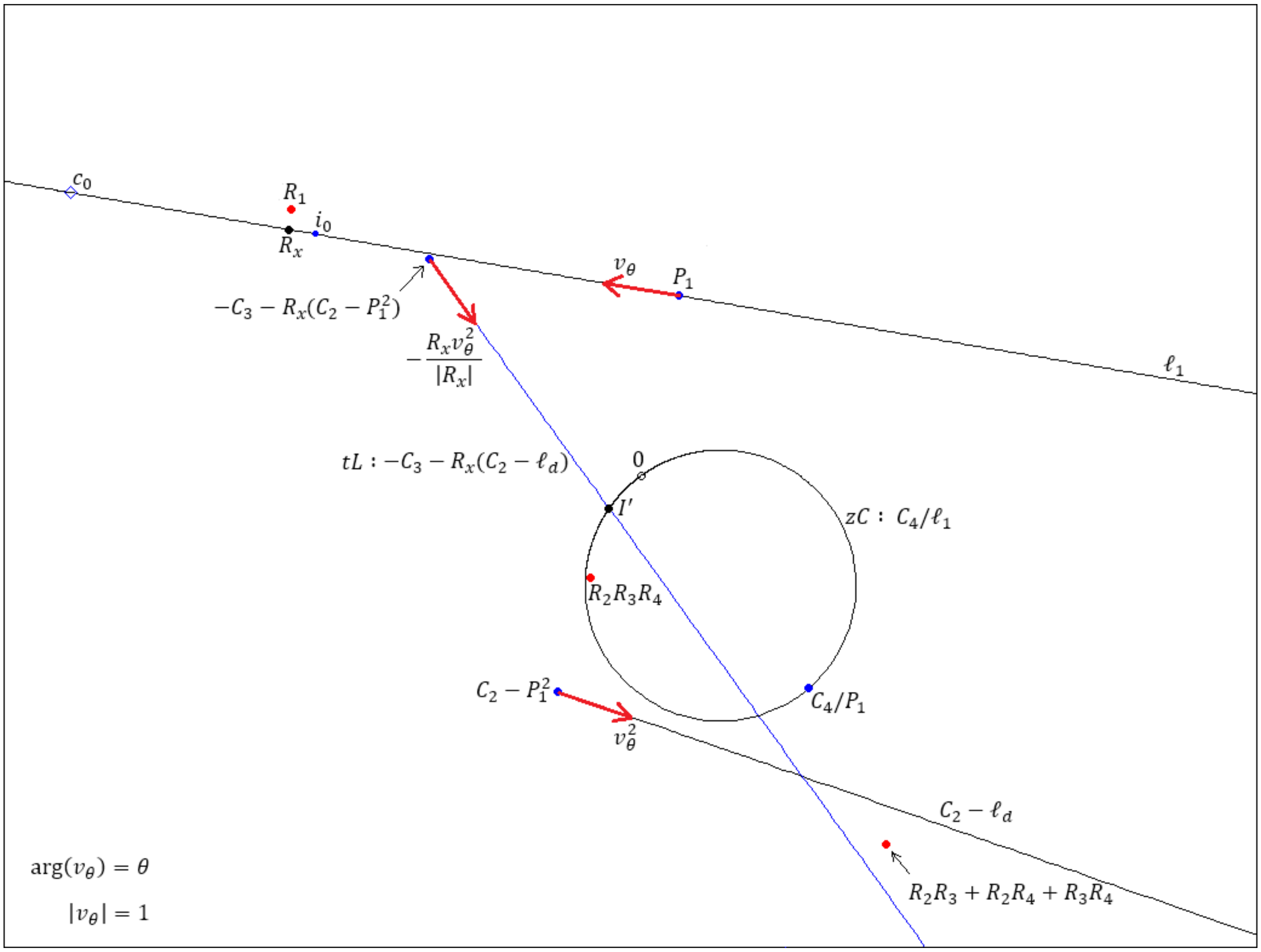

**Fig. 4.6**. Some elements of structure $LzC(C_1, C_2, C_3, C_4, \theta)$, in which $\theta$ is "close" to value $\theta^*$ shown in figure 4.1, so that $\ell_1$ is "close" to $R_1$. The projection of $I'$ (the intersection between $zC$ and $tL$, which tends to approximate value $R_2R_3R_4$ as $\theta$ tends to approximate $\theta^*$) onto $\ell_1$ produces a couple of points, $i_0$ and $c_0$, which help us quantify the degree of "closeness" between $\ell_1$ and $R_1$, without knowing where $R_1$ is located. All direction vectors shown in this schematic graph are unit vectors.

Based on the diagram from figure 4.6, we will now find algebraic expressions for $i_0(\theta)$ and $c_0(\theta)$, which are the projections of $I'$ onto line $\ell_1(\theta)$, similar to those defined in expression (3.3) from chapter 3.

The projection of $I'$ onto line $\ell_1(\theta)$ via the z-circumference $zC(\theta)$ is simply

$$c_0(\theta) = C_4/I'. \qquad (4.9)$$

On the other hand, the projection of $I'$ onto line $\ell_1(\theta)$ via the terminal semi-line $tL(\theta)$ is given by $i_0(\theta) = P_1 + t_{I'}(\theta)v_\theta$, where $t_{I'}(\theta)$ is the positive value that satisfies the equation that is obtained by equating the parametric expression for $tL(\theta)$ in (4.8) with $I'$:

$$-C_3 - R_xC_2 + R_xP_1^2 + t_{I'}^2(\theta)(-R_xv_\theta^2) = I'$$

$$-C_3 - R_x(C_2 - P_1^2) + t_{I'}^2(\theta)(-R_xv_\theta^2) = I'$$

If we set $aP(\theta) = -C_3 - R_x(\theta)(C_2 - P_1^2)$ and $v_{tL}(\theta) = -R_x(\theta)v_\theta^2$, then we have

$$aP(\theta) + t_{I'}^2(\theta)v_{tL}(\theta) = I'$$

$$\Rightarrow t_{I'}(\theta) = \sqrt{\frac{I'-aP(\theta)}{v_{tL}(\theta)}}.$$

Therefore, the projection of $I'$ onto line $\ell_1(\theta)$ via the terminal semi-line $tL(\theta)$ is

$$i_0(\theta) = P_1 + \sqrt{\frac{I'-aP(\theta)}{v_{tL}(\theta)}}\,v_\theta, \qquad (4.10)$$

where

$aP(\theta) = -C_3 - R_x(\theta)(C_2 - P_1^2)$ is the *anchor point* for terminal semi-line $tL(\theta)$, and

$v_{tL}(\theta) = -R_x(\theta)v_\theta^2$ is a non-normalized direction vector for terminal semi-line $tL(\theta)$.

We will show that the points $c_0(\theta), i_0(\theta) \in \ell_1(\theta)$ defined by expressions (4.9) and (4.10), under the relevant inverse transformations, each of them corresponds to the point $I'$ shown in the diagram from figure 4.6. First, from expression (4.9), we immediately see that $I' = C_4/c_0(\theta)$; i.e., $c_0(\theta) \in \ell_1(\theta)$ corresponds to $I' \in zC(\theta)$, by the properties of the Möbius transformation seen in chapter 1. Second, to prove that $i_0(\theta) \in \ell_1(\theta)$ corresponds to $I' \in tL(\theta)$, we must consider the reflection of $i_0(\theta)$ on $\ell_1(\theta)$ with respect to fixed point $P_1$; i.e., we must consider

$$i_1(\theta) = P_1 - \sqrt{[I' - aP(\theta)]/v_{tL}(\theta)}\,v_\theta;$$

hence, we have

$$i_0(\theta) \in \ell_1(\theta)$$

$$\Rightarrow \qquad i_0(\theta)i_1(\theta) \in \ell_d(\theta)$$

$$\Rightarrow \qquad C_2 - i_0(\theta)i_1(\theta) \in C_2 - \ell_d(\theta)$$

$$\Rightarrow \qquad R_x(\theta)[C_2 - i_0(\theta)i_1(\theta)] \in R_x(\theta)[C_2 - \ell_d(\theta)]$$

$$\Rightarrow -C_3 - R_x(\theta)[C_2 - i_0(\theta)i_1(\theta)] \in -C_3 - R_x(\theta)[C_2 - \ell_d(\theta)]$$

In this last expression, $-C_3 - R_x(\theta)[C_2 - \ell_d(\theta)] = tL(\theta)$, according to (4.8); if we use here the expressions for $i_0(\theta)$ and $i_1(\theta)$ related to (4.10), we have

$$-C_3 - R_x(\theta)\left[C_2 - \left(P_1 + \sqrt{\frac{I' - aP(\theta)}{v_{tL}(\theta)}}v_\theta\right)\left(P_1 - \sqrt{\frac{I' - aP(\theta)}{v_{tL}(\theta)}}v_\theta\right)\right] \in tL(\theta)$$

$$\Rightarrow -C_3 - R_x(\theta)\left[C_2 - P_1^2 + \frac{I' - aP(\theta)}{v_{tL}(\theta)}v_\theta^2\right] \in tL(\theta)$$

$$\Rightarrow -C_3 - R_x(\theta)\left[C_2 - P_1^2 + \frac{I' + C_3 + R_x(\theta)(C_2 - P_1^2)}{-R_x(\theta)v_\theta^2}v_\theta^2\right] \in tL(\theta)$$

$$\Rightarrow \cancel{-C_3 - R_x(\theta)(C_2 - P_1^2)} + I' + \cancel{C_3 + R_x(\theta)(C_2 - P_1^2)} \in tL(\theta)$$

That is to say, $i_0(\theta) \in \ell_1(\theta)$ corresponds to $I' \in tL(\theta)$.

Now, once we know the elements $\ell_1(\theta)$, $tL(\theta)$ and $zC(\theta)$ from structure $LzC(C_1, C_2, C_3, C_4, \theta)$, as well as the possible intersections $I_1$, $I_2$ between $tL(\theta)$ and $zC(\theta)$ and their associated projections $i_{0I_1}$, $c_{0I_1}$, $i_{0I_2}$, $c_{0I_2}$ onto line $\ell_1(\theta)$, we can define weighted error functions $e_A(\theta)$, $e_B(\theta)$ in exactly the same way as it was done in expressions (3.4a), (3.4b) from chapter 3:

**Weighted error measurements associated with structure $\boldsymbol{LzC(C_1, C_2, C_3, C_4, \theta)}$, $\boldsymbol{\theta \in [-\pi, \pi)}$**

$$e_A(\theta) := \begin{cases} sgn\left(\frac{i_{0I_1} - c_{0I_1}}{v_\theta}\right) \frac{|i_{0I_1} - c_{0I_1}|}{|i_{0I_1} - P_1| + |c_{0I_1} - P_1|} & \text{if } \exists\, I_1 \\ \text{undefined} & \text{if } \nexists\, I_1 \end{cases} \tag{4.11a}$$

$$e_B(\theta) := \begin{cases} sgn\left(\frac{i_{0I_2} - c_{0I_2}}{v_\theta}\right) \frac{|i_{0I_2} - c_{0I_2}|}{|i_{0I_2} - P_1| + |c_{0I_2} - P_1|} & \text{if } \exists\, I_2 \\ \text{undefined} & \text{if } \nexists\, I_2 \end{cases} \tag{4.11b}$$

where:

$$sgn(x) = \begin{cases} -1 & \text{if } x < 0 \\ 0 & \text{if } x = 0 \\ 1 & \text{if } x > 0 \end{cases}, \qquad i_{0I_m} = P_1 + \sqrt{\frac{[I_m - aP(\theta)]}{v_{tL}(\theta)}}\, v_\theta\,, \qquad c_{0I_m} = C_4/I_m\,, \qquad m = 1,2.$$

$I_m$ is one of two possible intersections between terminal semi-line $tL(\theta) = aP(\theta) + t^2 v_{tL}(\theta)$ and z-circumference $zC(\theta) = C_4/(P_1 + tv_\theta)$, with $aP(\theta) = -C_3 - R_x(\theta)(C_2 - P_1^2)$ the anchor point of $tL(\theta)$, and $v_{tL}(\theta) = -R_x(\theta)v_\theta^2$ a direction vector for $tL(\theta)$. $v_\theta$ is the unit direction vector of line $\ell_1(\theta)$ with fixed point $P_1 = -C_1/2$.

$R_x(\theta) = P_1 + t^*(\theta)v_\theta$ estimates one of the roots of equation (4.1), with $t^*(\theta) = \underset{t \in \mathbb{R}}{\arg\min}\, d_\theta^2(t)$.

The dynamic squared distance function $d_\theta^2(t)$ is defined by (4.4).

We already have all the necessary concepts to extend the strategy 3.1 of chapter 3 towards polynomials of degree 4 in a complex variable; in this extended strategy, we would construct discrete proximity maps, in principle taking advantage of the parallel processing capabilities of modern computers. As we mentioned in chapter 3, our implementations in R (which are for purely illustrative purposes and in no way intend to implement general-usage procedures) are designed to run multiple functions sequentially, and each of these implemented functions uses vectorized operations, which can be considered as a form of implicit parallelism in R, adapted according to the capabilities of the available hardware. Next, we list the strategy 4.1 to construct discrete proximity maps associated to polynomials of degree 4 in a complex variable:

**Strategy 4.1 to find initial approximations to the roots of a univariate polynomial of degree 4 with complex coefficients by means of parallel constructions $LzC(C_1, C_2, C_3, C_4, \theta)$**

a) Find (in parallel) minimizing arguments $t^*(\theta_k) = \arg\min_{t \in \mathbb{R}} d^2_{\theta_k}(t)$ associated with structures $LzC(C_1, C_2, C_3, C_4, \theta_k)$ for values $\theta_k = -\pi + 2\pi k/N$, $k = 0,1,2, \dots, N-1$; note that points $\theta_k$ constitute a regular partition of interval $[-\pi, \pi)$.

b) Evaluate (in parallel) weighted error functional values $e_A(\theta_k)$, $e_B(\theta_k)$ given by (4.11a) and (4.11b), which require as input, the values $t^*(\theta_k)$ previously computed in a).

c) Join, by using rectilinear segments, adjacent points $(\theta_k, e_A(\theta_k)), (\theta_{k+1}, e_A(\theta_{k+1})) \in \mathbb{R}^2 \ \forall k$; in this way, we will have constructed, by linear interpolation, a **discrete approximation** $\hat{e}_A(\theta)$ of function $e_A(\theta)$ at interval $-\pi \leq \theta < \pi$. Similarly, obtain a discrete approximation $\hat{e}_B(\theta)$ of function $e_B(\theta)$ at interval $-\pi \leq \theta < \pi$.

d) Find crossings of $\hat{e}_A(\theta)$ and $\hat{e}_B(\theta)$ with horizontal axis $y = 0$; with this, we will obtain *candidates* for estimates $\hat{\theta}_i^*$ of direction angles $\theta_i^*$ ($i \in \{1,2,3,4\}$) that, for each $i$, cause line $\ell_1(\theta_i^*)$ to contain one of the roots $R_i$ of equation (4.1).

e) By means of parallel constructions $LzC(\hat{\theta}_i^*)$, compute estimates $\hat{R}_i$ of roots $R_i$ of equation (4.1). For this, simply obtain points $R_x(\hat{\theta}_i^*) = P_1 + t^*(\hat{\theta}_i^*) v_{\hat{\theta}_i^*}$, which serve as possible initial approximations to roots $R_i$. For each $i$, $t^*(\hat{\theta}_i^*)$ is the minimizing argument of the dynamic squared distance function $d^2_{\hat{\theta}_i^*}(t)$ defined in (4.4).

Of course, obtaining all the necessary elements in this strategy for the construction of a discrete proximity map $e(\theta)$ implies first obtaining the minimizing argument $t^*(\theta_k)$, within each structure $LzC(C_1, C_2, C_3, C_4, \theta_k)$, for the dynamic squared distance between terminal curve $t\mathfrak{C}(\theta_k)$ and z-circumference $zC(\theta_k)$. This adds a cost in computational complexity; in addition to finding intersections between lines and circumferences, it is now also necessary to minimize a function whose complexity depends on the degree of the polynomial for which we try to obtain a discrete proximity map $e(\theta)$; for this optimization task, we must define an implementation strategy, so that our estimate of the minimizing argument $t^*(\theta)$ is reliable enough not to distort the true form of the weighted error functions $e(\theta)$; at the same time, our implemented optimization algorithm must be efficient in terms of temporal complexity.

In our strategy to compute $t^*(\theta)$, we need an algorithm that preferably does not get stuck at local minima; here we can use the *simulated annealing* (SANN) metaheuristic; this algorithm produces an initial estimate $\hat{t}_0^*(\theta)$ that contains a certain amount of stochastic noise $\varepsilon$; i.e., $\hat{t}_0^*(\theta) = t^*(\theta) + \varepsilon$. Therefore, it is necessary to add a second optimization phase, which takes as input SANN's output $\hat{t}_0^*(\theta)$; here, we can use the BFGS method (a quasi-Newton method), in order to refine $\hat{t}_0^*(\theta)$, reducing $\varepsilon$ as much as possible. For further details, see annex 2 section 4 (function `min_D2`). As we said previously, this is an initial proposal for estimating $t^*(\theta)$ that, nevertheless, produces reasonable approximations, as we will see next. If more accurate approximations (computed by more efficient optimization algorithms) are sought, the concepts and ideas used in this work serve as a starting point for further investigation.

## Numerical Example 4.1

In this example, we will use the same strategy used in the numerical examples from chapter 3; first, we will randomly generate a set of roots $R_1$, $R_2$, $R_3$, $R_4$, which will serve as reference when measuring estimation errors associated with the generated proximity maps. For this, we will use the script listed in annex 3 section 3, initializing the random number generator seed by means of the instruction `set.seed(1987)`. Then, from the randomly generated roots, we will construct coefficients $C_1$, $C_2$, $C_3$, $C_4$, applying Vieta's relations (4.2-a), (4.2-b), (4.2-c) and (4.2-d); these coefficients will serve as inputs to function `approxLC`, which generates a discrete proximity map $e(\theta)$ (also called LC map) and the corresponding approximations to the roots for the equation of the form (4.1); of course, `approxLC` (which is called within the script from annex 3 section 3) never uses direct information about the true roots $R_1$, $R_2$, $R_3$, $R_4$. For more details on the inner workings of function `approxLC`, which implements strategy 4.1, see annex 2 section 6.

The roots $R_i$ generated in this example, together with their corresponding theta roots (inclination angles $\theta_i^*$ for lines $\ell_1$ with fixed point $-C_1/2$ containing roots $R_i$), are:

$$R_1 = -0.8024716 + 0.6100969i \qquad \theta_1^* = 2.923580$$

$$R_2 = 0.2151479 - 0.2343356i \qquad \theta_2^* = -1.375446 \qquad (4.12)$$

$$R_3 = 0.9956074 + 0.9848524i \qquad \theta_3^* = 0.561732$$

$$R_4 = -0.2341845 - 0.5345379i \qquad \theta_4^* = -1.897644$$

The coefficients constructed by plugging generated roots $R_1$, $R_2$, $R_3$, $R_4$ into Vieta's relations (4.2-a), (4.2-b), (4.2-c) and (4.2-d) are:

$$C_1 = -0.1740992 - 0.8260758i$$

$$C_2 = -0.3528101 - 0.4218854i \qquad (4.13)$$

$$C_3 = 0.0520028 - 0.7879940i$$

$$C_4 = 0.2348718 + 0.1162912i$$

By constructing a global LC proximity map with $N = 2{,}500$ elements $LzC(C_1, C_2, C_3, C_4, \theta_k)$ associated with points $\theta_k$ in a regular partition of interval $[-\pi, \pi)$, we obtain the graph in figure 4.7, which also displays reference values $\theta_1^*$, $\theta_2^*$, $\theta_3^*$, $\theta_4^*$ as small circles on the horizontal axis $y = 0$.

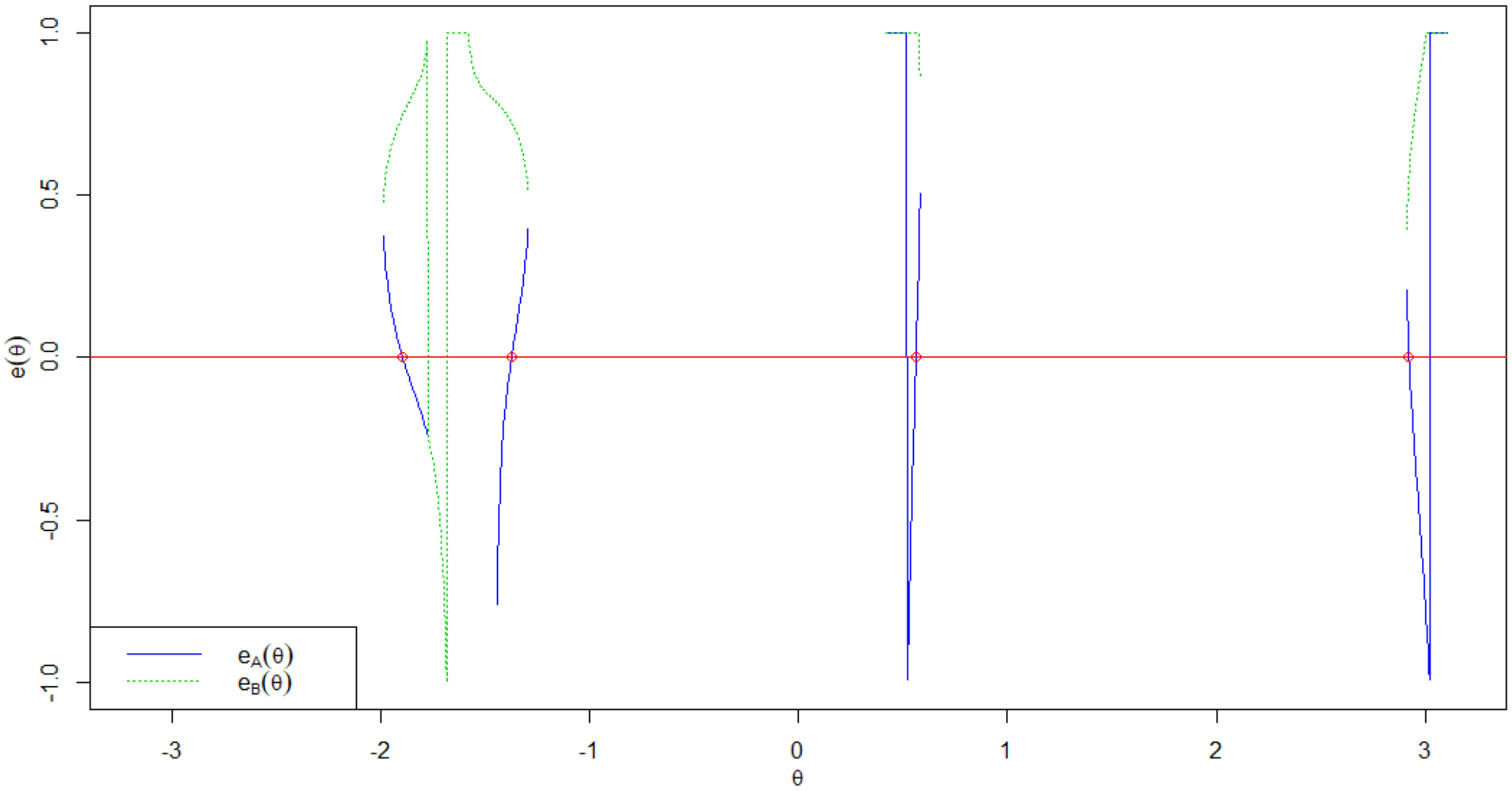

**Fig. 4.7.** Discrete angular proximity LC map for polynomial $p(z) = z^4 + C_1 z^3 + C_2 z^2 + C_3 z + C_4$ with coefficients $C_1$, $C_2$, $C_3$, $C_4$ given by (4.13). This map was generated from $N = 2{,}500$ elements $LzC(\theta_k)$ associated with points $\theta_k$ in a regular partition of interval $[-\pi, \pi)$.

As we can see in figure 4.7, in this case there are four "smooth" crossings $\hat{\theta}_i^*$ of map $e(\theta)$ with horizontal axis $y = 0$; we will see in numeric example 4.2, however, that there may be more than four smooth crossings in a map $e(\theta)$ associated with a univariate degree 4 polynomial; this implies that a smooth crossing in an LC proximity map associated with a univariate degree 4 polynomial is a necessary, but not sufficient, condition for obtaining an approximation $\hat{\theta}_i^*$ to a true theta root $\theta_i^*$. In light of this empirical evidence, we see that it is necessary to establish a numerical criterion that helps us identify smooth crossings $\hat{\theta}_i^*$ associated with true theta roots $\theta_i^*$ within a proximity map $e(\theta)$ such as the one shown in figure 4.13 from numerical example 4.2. One possibility is to preferentially select smooth crossings $\hat{\theta}_i^*$ with small associated global minima $d^2_{\hat{\theta}_i^*}\left(t^*(\hat{\theta}_i^*)\right)$; the smaller the global minimum of the dynamic squared distance associated with smooth crossing $\hat{\theta}_i^*$, the better. Table 4.1 shows the numerical details (analogous to those defined in chapter 3) of the smooth crossings $\hat{\theta}_i^*$ identified in the proximity map $e(\theta)$ from figure 4.7, together with their corresponding global minima $d^2_{\hat{\theta}_i^*}\left(t^*(\hat{\theta}_i^*)\right)$.

**Table 4.1**. Initial estimates for the roots of $p(z) = z^4 + C_1 z^3 + C_2 z^2 + C_3 z + C_4$ with coefficients $C_1$, $C_2$, $C_3$, $C_4$ given by (4.13). These estimates were obtained from the LC map in figure 4.7.

| $i$ | $\hat{R}_i$ | $\hat{\theta}_i^*$ | $\lvert\Delta e_i\rvert$ | $d^2_{\hat{\theta}_i^*}\left(t^*(\hat{\theta}_i^*)\right)$ |
|---|---|---|---|---|
| 1 | $-0.2341833 - 0.5345378i$ | $-1.8976427$ | $0.005575651$ | $7.111838 \times 10^{-12}$ |
| 2 | $0.2151517 - 0.2343367i$ | $-1.3754407$ | $0.010994010$ | $1.642397 \times 10^{-10}$ |
| 3 | $0.9956091 + 0.9848480i$ | $0.5617277$ | $0.055988753$ | $3.048712 \times 10^{-10}$ |
| 4 | $-0.8024751 + 0.6100794i$ | $2.9235994$ | $0.033677561$ | $2.995278 \times 10^{-9}$ |

From table 4.1, we can see that the rows with smooth crossings $\hat{\theta}_i^*$ associated with the LC map in figure 4.7 are arranged in ascending order (from the smallest to the largest values), in accordance with variable $d^2_{\hat{\theta}_i^*}\left(t^*(\hat{\theta}_i^*)\right)$; this ordering will be useful when identifying valid estimates in cases where we observe more than $n$ smooth crossings $\hat{\theta}_i^*$, where $n \geq 4$ is the degree of the univariate polynomial in question. For now, if we compare the initial approximations $\hat{R}_i$, $\hat{\theta}_i^*$ in table 4.1 with reference values $R_i$, $\theta_i^*$ in (4.12), we can see that the approximations listed in table 4.1 indeed look quite similar to the true roots; we see that at least the first five digits in each of the real and imaginary parts for approximations $\hat{R}_i$, as well as of values $\hat{\theta}_i^*$ in table 4.1, coincide with the corresponding digits in reference values (4.12).

What is the degree of accuracy for the approximations in table 4.1? Let us compute the absolute relative errors for these approximations with respect to corresponding reference values (4.12):

$$\begin{aligned}
\left|\frac{R_1-\hat{R}_4}{R_1}\right| &= 1.763275 \times 10^{-5} & \left|\frac{\theta_1^*-\hat{\theta}_4^*}{\theta_1^*}\right| &= 6.671271 \times 10^{-6} \\
\left|\frac{R_2-\hat{R}_2}{R_2}\right| &= 1.238908 \times 10^{-5} & \left|\frac{\theta_2^*-\hat{\theta}_2^*}{\theta_2^*}\right| &= 3.842554 \times 10^{-6} \\
\left|\frac{R_3-\hat{R}_3}{R_3}\right| &= 3.359573 \times 10^{-6} & \left|\frac{\theta_3^*-\hat{\theta}_3^*}{\theta_3^*}\right| &= 7.650780 \times 10^{-6} \\
\left|\frac{R_4-\hat{R}_1}{R_4}\right| &= 1.975550 \times 10^{-6} & \left|\frac{\theta_4^*-\hat{\theta}_1^*}{\theta_4^*}\right| &= 5.716482 \times 10^{-7}
\end{aligned} \qquad (4.14)$$

Absolute relative errors (4.14) are similar, in terms of order of magnitude, to absolute relative errors (3.8) and (3.11) observed in the first two examples of chapter 3 (cubic polynomials) where we have non-repeated roots; in those examples, however, we used only $N = 1{,}000$ elements $LzC(\theta_k)$, while in this example we are using $N = 2{,}500$ elements $LzC(\theta_k)$. It is worth mentioning here, that the degree of accuracy for initial approximations $\hat{R}_i$, $\hat{\theta}_i^*$ obtained from a proximity map associated with a univariate polynomial of degree 4 or higher, could be affected not only by the map's degree of discretization, but also by the fact that such map is constructed from numerical approximations to minimizing arguments $t^*(\theta_k)$, according to strategy 4.1. These approximations to minimizing arguments $t^*$ are computed by means of function `min_D2`, which uses, in a first optimization phase, the simulated annealing (SANN) algorithm (whose stochastic behavior is regulated by a series of parameters) to obtain an initial approximation $\hat{t}_0^*$; and in a second optimization phase, the BFGS algorithm, to refine $\hat{t}_0^*$. For further details, see annex 2 section 4.

Note: If initial approximations $\hat{R}_i$, $\hat{\theta}_i^*$ from table 4.1 are ranked according to the absolute relative error $\left|\frac{R-\hat{R}}{R}\right|$, we obtain a different ordering than that shown in said table; the same happens if the criterion for ranking approximations is based on the absolute relative error $\left|\frac{\theta^*-\hat{\theta}^*}{\theta^*}\right|$; furthermore, the two ordering criteria $\left|\frac{R-\hat{R}}{R}\right|$, $\left|\frac{\theta^*-\hat{\theta}^*}{\theta^*}\right|$ produce orderings that are different from each other.

**Additional results related to the dynamic squared distance function.** In this numerical example (and in the following numerical examples hereafter), in addition to the discrete LC proximity map formed by weighted errors $e(\theta_k)$ associated with intersections between terminal semi-lines $tL(\theta_k)$ and z-circumferences $zC(\theta_k)$, we will also generate some other discrete graphs related to the dynamic squared distance function $d_\theta^2(t)$, which after all turns out to be fundamental in the construction of structures $LzC(C_1, C_2, C_3, C_4, \theta_k)$. A first graph is the trajectory, on the plane of complex numbers $\mathbb{C}$, traced by the "best" approximations $R_x(\theta_k)$ defined by (4.7), together with the trajectory of anchor points $aP(\theta_k)$ for terminal semi-lines $tL(\theta_k)$ defined by expression $aP(\theta_k) = -C_3 - R_x(\theta_k)(C_2 - P_1^2)$. The tracings of these trajectories are shown in figure 4.8; we see here that both $aP(\theta_k)$ and $R_x(\theta_k)$ are disjoint trajectories, and at the same time they are similar to each other, although of course with different spatial orientation and different size between them.

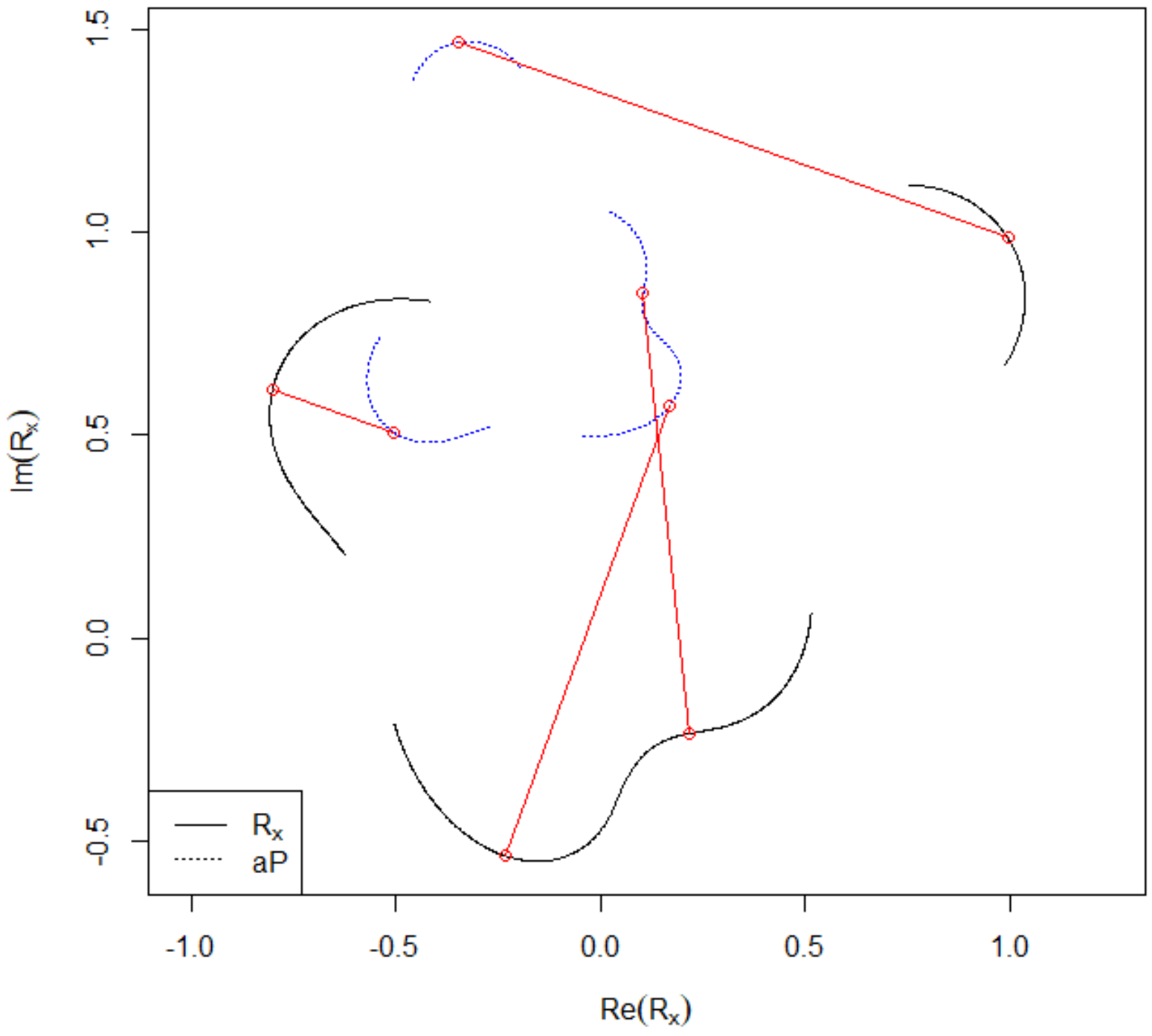

**Fig. 4.8.** Trajectory of best approximations $R_x(\theta_k)$ (solid line), and trajectory of anchor points $aP(\theta_k)$ for terminal semi-lines $tL(\theta_k)$ (dotted line). The trajectory $R_x(\theta_k)$ shown here contains the true roots of polynomial $p(z) = z^4 + C_1z^3 + C_2z^2 + C_3z + C_4$ with coefficients $C_1$, $C_2$, $C_3$, $C_4$ given by (4.13). The true roots are represented in this figure by means of small circles, which are joined, by means of line segments, with their corresponding anchor points.

Additional discrete graphs related to the dynamic squared distance function $d^2_{\theta_k}(t)$ are shown in figures 4.9 and 4.10, in which we can see, respectively, global minima $d^2_{\theta_k}\big(t^*(\theta_k)\big) = \min_t d^2_{\theta_k}(t)$ vs. $\theta_k$, and corresponding minimizing arguments $t^*(\theta_k) = \arg\min_t d^2_{\theta_k}(t)$ vs. $\theta_k$. It is worth mentioning that **the graphs of functions $\boldsymbol{d^2_\theta\big(t^*(\theta)\big)}$ and $\boldsymbol{t^*(\theta)}$ can be visualized as piecewise continuous and differentiable real functions, which depend on a real variable $\boldsymbol{\theta}$**. In reality, under strategy 4.1, the numerical values for the graphs shown in figures 4.9 and 4.10 were obtained by function `min_D2` (annex 2 section 4), within function `approxLC` (annex 2 section 6), before computing the LC proximity map shown in figure 4.7. From what can be seen in figures 4.9 and 4.10, these graphs alone also provide clues that help locate the roots of polynomial equations of the form (4.1). From this, we can say that **the graphs of functions $\boldsymbol{d^2_\theta\big(t^*(\theta)\big)}$ and $\boldsymbol{t^*(\theta)}$ can also be considered as a distinct class of discrete proximity maps.**

If we analyze the discrete arrays $d^2_{\theta_k}\big(t^*(\theta_k)\big)$ and $t^*(\theta_k)$ generated by the script from annex 3 section 3 (arrays `MinF` and `MinT`, respectively), we can see that half of the elements in each array contain numerical values, while remaining elements contain `NA` values (not available, or null values). This is partly so because the function $d^2_\theta(t)$ defined in (4.4) has the property that $d^2_\theta(t) = d^2_{\theta+\pi}(-t)$; i.e., function $d^2_\theta(t)$ is the mirror reflection of $d^2_{\theta+\pi}(t)$ with respect to vertical axis $t = 0$, since both functions are associated to the same line $\ell_1(\theta) = \ell_1(\theta + \pi)$ with fixed point $P_1 = -C_1/2$ in complex plane $\mathbb{C}$. This symmetry property for function $d^2_\theta(t)$, together with the imposed constraint $t^*(\theta) > 0$, cause that only half of the support values $\theta$ in the graphs from figures 4.9 and 4.10 have associated functional values. The same is true for the graphs shown in figures 4.11 and 4.12, which correspond, respectively, to the first-order discrete derivatives of $d^2_\theta\big(t^*(\theta)\big)$ and $t^*(\theta)$ with respect to $\theta$.

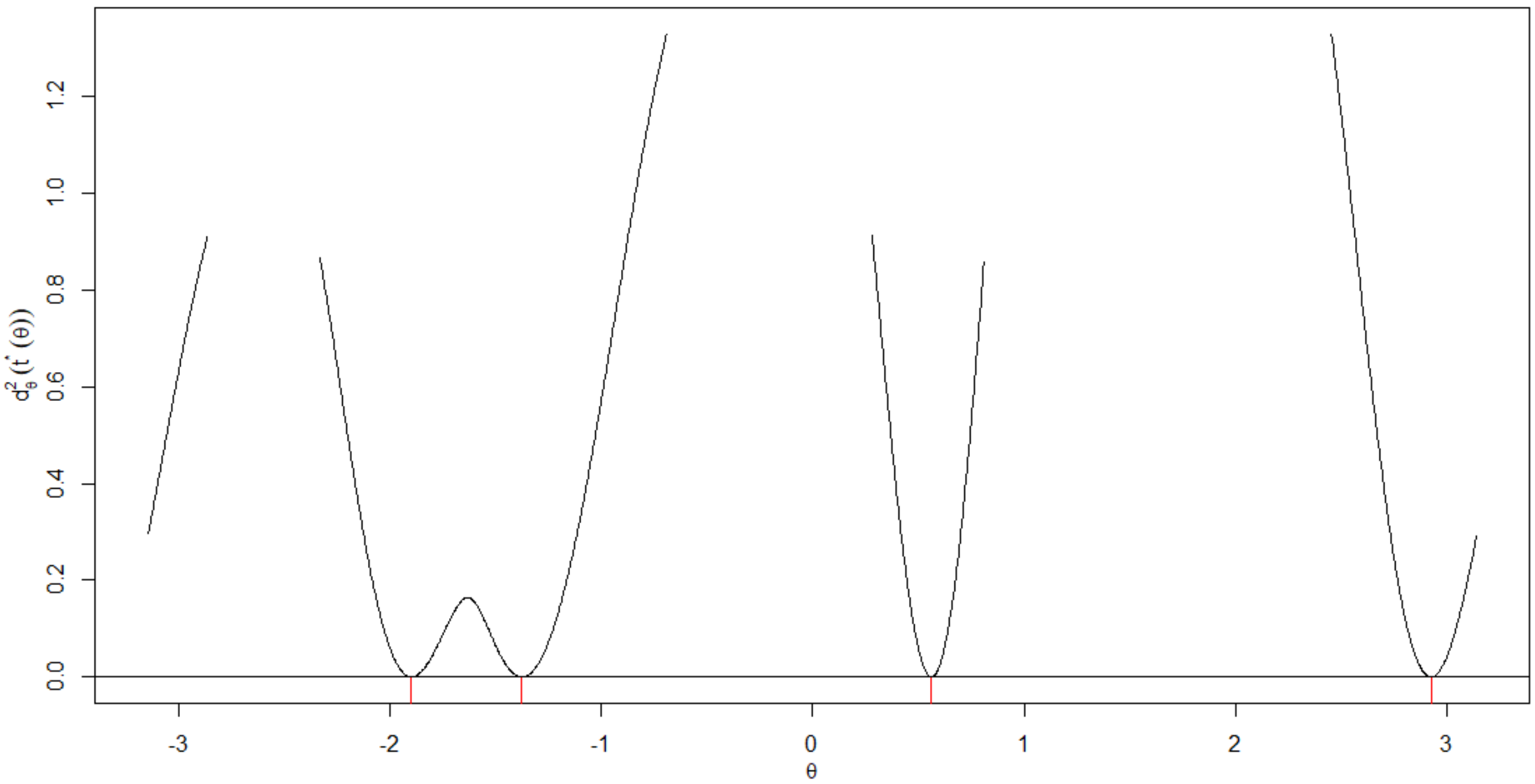

**Fig. 4.9.** Global minima of the dynamic squared distance functions $d^2_{\theta_k}(t^*(\theta_k)) = \min_t d^2_{\theta_k}(t)$ associated to polynomial $p(z) = z^4 + C_1z^3 + C_2z^2 + C_3z + C_4$ with coefficients $C_1$, $C_2$, $C_3$, $C_4$ given by (4.13), vs. values $\theta_k$. This graph can be thought of as a real function $d^2_\theta(t^*(\theta))$ which depends on variable $\theta$. At the bottom of the graph there are vertical line segments which indicate the location of true theta roots $\theta_i^*$.

The behavior of the discrete graph $d^2_{\theta_k}(t^*(\theta_k))$ vs. $\theta_k$ in figure 4.9, suggests that theta roots $\theta_i^*$ are the minimizing arguments of continuous and differentiable function $d^2_\theta(t^*(\theta))$, which is to be expected, according to the process of construction (described at the beginning of this chapter) of a structure $LzC(C_1, C_2, C_3, C_4, \theta^*)$ and of neighboring structures $LzC$; the behavior of the graph in figure 4.9 also suggests that function $d^2_\theta(t^*(\theta))$ is concave upward in the vicinities around true theta roots $\theta_i^*$, which implies that the derivative of $d^2_\theta(t^*(\theta))$ with respect to $\theta$ is an increasing function in the vicinities around values $\theta_i^*$, and exactly zero at values $\theta_i^*$: this hypothesis is reinforced by analyzing the behavior shown in figure 4.11, which displays a discrete approximation of $\frac{d}{d\theta} d^2_\theta(t^*(\theta))$ (computed by means of function `approxDMin` listed in annex 2 section 7). In summary, the joint behavior of the graphs in figures 4.9 and 4.11, suggests the following properties for the function $d^2_\theta(t^*(\theta))$:

**Properties of mapping $\boldsymbol{d_\theta^2(t^*(\theta))}$:**

a) Theta roots $\theta_i^*$ occur at points where $\frac{d}{d\theta} d_\theta^2(t^*(\theta)) = 0$, and $d_\theta^2(t^*(\theta)) = 0$.

b) At neighborhoods around values $\theta_i^*$, $\frac{d}{d\theta} d_\theta^2(t^*(\theta))$ is an increasing function, which implies that $\frac{d^2}{d\theta^2} d_\theta^2(t^*(\theta)) > 0$ at $\theta_i^*$.

c) Since $\frac{d^2}{d\theta^2} d_\theta^2(t^*(\theta)) > 0$ and $\frac{d}{d\theta} d_\theta^2(t^*(\theta)) = 0$ at $\theta_i^*$, this implies that $d_\theta^2(t^*(\theta))$ attains a minimum value at $\theta_i^*$, by the second derivative criterion.

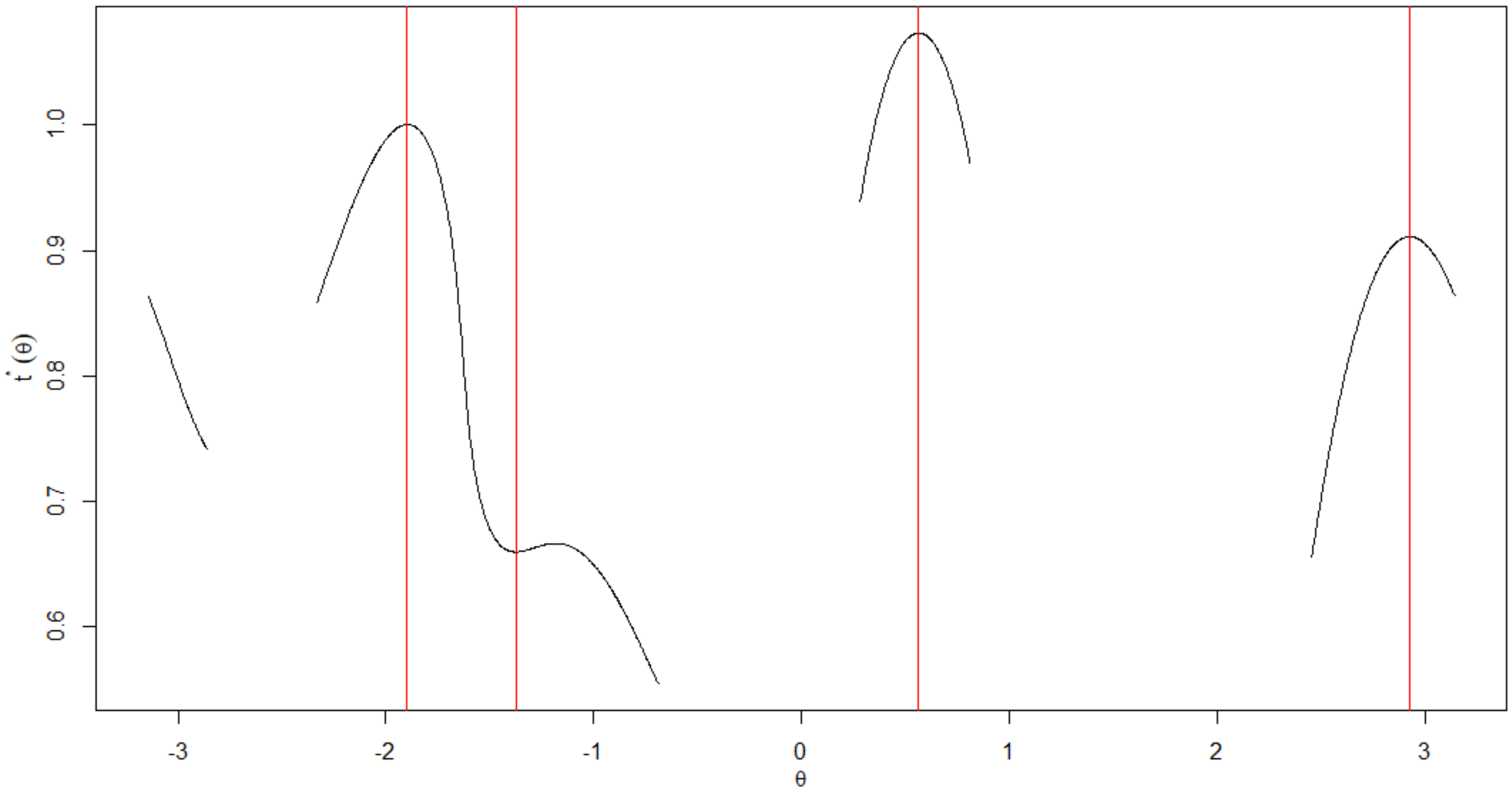

**Fig. 4.10.** Minimizing arguments of dynamic squared distance functions $t^*(\theta_k) = \arg\min_t d_{\theta_k}^2(t)$ associated to polynomial $p(z) = z^4 + C_1 z^3 + C_2 z^2 + C_3 z + C_4$ with coefficients $C_1, C_2, C_3, C_4$ given by (4.13), vs. values $\theta_k$. This graph can be thought of as a real function $t^*(\theta)$ which depends on variable $\theta$. This graph includes vertical line segments which indicate the location of true theta roots $\theta_i^*$.

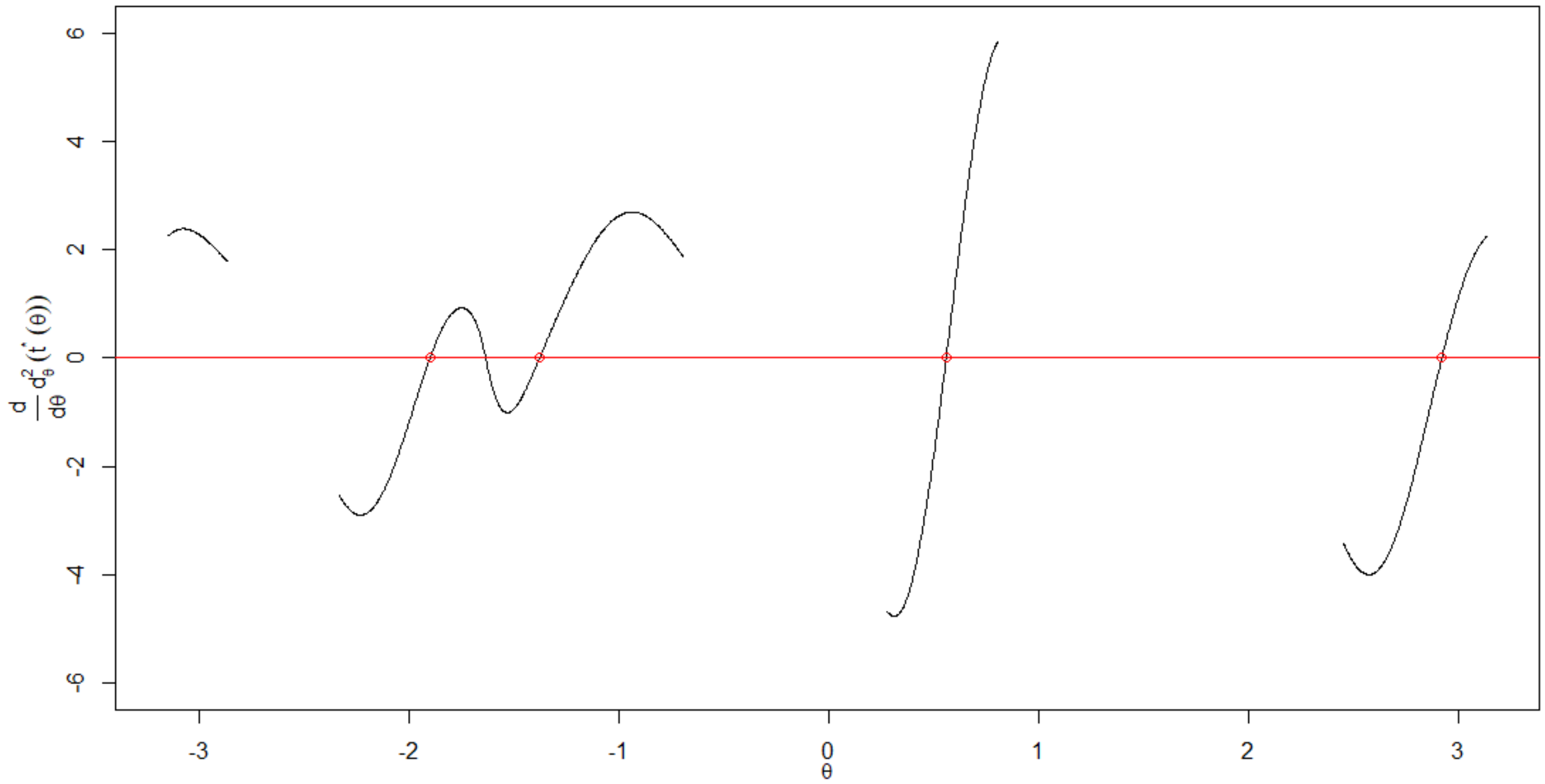

**Fig. 4.11**. Discrete approximation to the derivative of function $d_\theta^2(t^*(\theta))$ from figure 4.9 with respect to $\theta$; this function is associated to polynomial $p(z) = z^4 + C_1z^3 + C_2z^2 + C_3z + C_4$ with coefficients $C_1$, $C_2$, $C_3$, $C_4$ given by (4.13). Small circles on horizontal axis $y = 0$ indicate the location of true theta roots $\theta_i^*$.

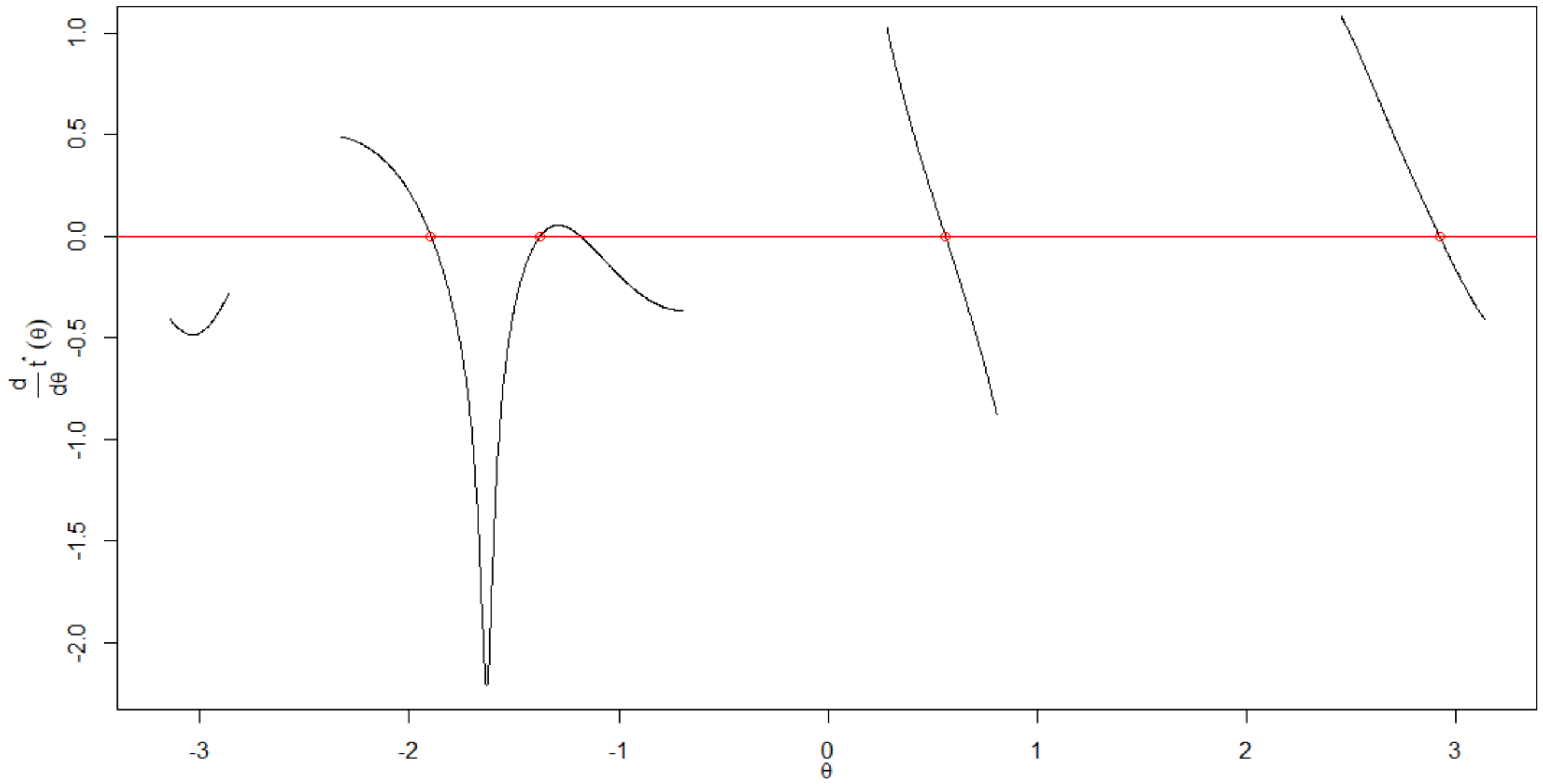

**Fig. 4.12.** Discrete approximation to the derivative of function $t^*(\theta)$ from figure 4.10 with respect to $\theta$; this function is associated to polynomial $p(z) = z^4 + C_1z^3 + C_2z^2 + C_3z + C_4$ with coefficients $C_1$, $C_2$, $C_3$, $C_4$ given by (4.13). Small circles on horizontal axis $y = 0$ indicate the location of true theta roots $\theta_i^*$.

Similarly, the joint behavior of graphs in figures 4.10 and 4.12, which show, respectively, discrete versions of functions $t^*(\theta)$ and $\frac{d}{d\theta}t^*(\theta)$ (the latter also computed by means of function `approxDMin` listed in annex 2 section 7), suggests the following properties for the function $t^*(\theta)$:

**Properties of mapping $t^*(\theta)$:**

a) Theta roots $\theta_i^*$ occur at critical points of function $t^*(\theta)$; i.e., at points where $\frac{d}{d\theta}t^*(\theta) = 0$.

b) Critical points of function $t^*(\theta)$ where theta roots $\theta_i^*$ are located can correspond to:

- local maxima (when $\frac{d^2}{d\theta^2}t^*(\theta) < 0$ at $\theta_i^*$, or equivalently, when $\frac{d}{d\theta}t^*(\theta)$ is a decreasing function at a vicinity around $\theta_i^*$),
- local minima (when $\frac{d^2}{d\theta^2}t^*(\theta) > 0$ at $\theta_i^*$, or equivalently, when $\frac{d}{d\theta}t^*(\theta)$ is an increasing function at a vicinity around $\theta_i^*$), or
- saddle points (when $\frac{d^2}{d\theta^2}t^*(\theta) = 0$ at $\theta_i^*$).

In this example, the joint behavior of the graphs in figures 4.10 and 4.12 suggest that three of the theta roots $\theta_i^*$ correspond to local maxima of $t^*(\theta)$, while one of the theta roots, specifically $\theta_2^* = -1.375446$, corresponds to a local minimum of function $t^*(\theta)$. In example 4.2, we will see a case where the joint behavior of the discrete versions of functions $t^*(\theta)$ and $\frac{d}{d\theta}t^*(\theta)$ suggests that one of the roots corresponds to a saddle point of $t^*(\theta)$.

<u>Note</u>: One thing in common that we can observe between the discrete versions of trajectories $R_x(\theta)$ and $aP(\theta)$ (figure 4.8), together with the discrete versions of real-valued functions of a real variable $\theta$ $d_\theta^2(t^*(\theta))$ (figure 4.9), $t^*(\theta)$ (figure 4.10), $\frac{d}{d\theta}d_\theta^2(t^*(\theta))$ (figure 4.11), and $\frac{d}{d\theta}t^*(\theta)$ (figure 4.12), is that all these graphs show the same number of "continuous sections", in this case 3 (considering that functions from figures 4.9, 4.10, 4.11 and 4.12 are periodic with period $2\pi$); we also see that each continuous section from any of these trajectories or functions has its counterpart (or image) at any other trajectory or function in this set of graphs. Each continuous section has certain particularities, which are reflected in its corresponding counterparts. For example, in any of figures 4.8 to 4.12, we can find a specific continuous section containing two roots.

In conclusion, we see that the graphs in figures 4.9, 4.10, 4.11 and 4.12, together with the graph in figure 4.7, represent examples of different classes of discrete proximity maps, which help us identify the location of theta roots $\theta_i^*$. It is worth mentioning that proximity maps $e(\theta)$, $d_\theta^2\big(t^*(\theta)\big)$, $\frac{d}{d\theta}d_\theta^2\big(t^*(\theta)\big)$, $t^*(\theta)$, $\frac{d}{d\theta}t^*(\theta)$, are all periodic functions with fundamental period $2\pi$; it is for this reason that we construct their graphs on the fundamental interval $\theta \in [-\pi, \pi)$. From this numerical example we further observe that the functional values of map $d_\theta^2\big(t^*(\theta)\big)$ in principle are in the range of non-negative real numbers $[0, +\infty)$, while the functional values of map $t^*(\theta)$ are in the range of positive real numbers $(0, +\infty)$; as we have already seen, the functional values of map $e(\theta)$ are in interval $(-1,1]$. On the other hand, the functional values of maps $\frac{d}{d\theta}d_\theta^2\big(t^*(\theta)\big)$ and $\frac{d}{d\theta}t^*(\theta)$, can be positive or negative numbers, although they are not limited to interval $(-1,1]$, as is the case with map $e(\theta)$. From the latter, we see that **the discrete maps $\boldsymbol{\frac{d}{d\theta}d_\theta^2\big(t^*(\theta)\big)}$ and $\boldsymbol{\frac{d}{d\theta}t^*(\theta)}$ can be processed in the same way than weighted error maps $\boldsymbol{e(\theta)}$: by identifying smooth crossings of $\boldsymbol{\frac{d}{d\theta}d_\theta^2\big(t^*(\theta)\big)}$ or $\boldsymbol{\frac{d}{d\theta}t^*(\theta)}$ with horizontal axis $\boldsymbol{y = 0}$, we obtain theta roots candidates $\boldsymbol{\hat{\theta}_i^*}$, which can be ranked by their associated values $\boldsymbol{d_{\hat{\theta}_i^*}^2\left(t^*(\hat{\theta}_i^*)\right)}$**, as it was done with approximations $\hat{\theta}_i^*$ in table 4.1. In this particular case, from figures 4.11 and 4.12, we see that there are 5 smooth crossings of the discrete version of $\frac{d}{d\theta}d_\theta^2\big(t^*(\theta)\big)$ with horizontal axis $y = 0$, and there are also 5 smooth crossings of the discrete version of $\frac{d}{d\theta}t^*(\theta)$ with horizontal axis $y = 0$. Tables 4.2 and 4.3 show, respectively, the numerical characteristics of smooth crossings $\hat{\theta}_i^*$ of the discrete versions of $\frac{d}{d\theta}d_\theta^2\big(t^*(\theta)\big)$ and $\frac{d}{d\theta}t^*(\theta)$ with horizontal axis $y = 0$; these smooth crossings $\hat{\theta}_i^*$ are ranked in ascending order, from smallest to largest value, according to their associated values $d_{\hat{\theta}_i^*}^2\left(t^*(\hat{\theta}_i^*)\right)$.

**Table 4.2**. Initial estimates for the roots of $p(z) = z^4 + C_1z^3 + C_2z^2 + C_3z + C_4$ with coefficients $C_1$, $C_2$, $C_3$, $C_4$ given by (4.13). These estimates were obtained from the discrete map $\frac{d}{d\theta}d_\theta^2\big(t^*(\theta)\big)$ in figure 4.11.

| $i$ | $\hat{R}_i$ | $\hat{\theta}_i^*$ | $\left\lvert\Delta\frac{d}{d\theta}d_{\hat{\theta}_i^*}^2(t^*)\right\rvert$ | $d_{\hat{\theta}_i^*}^2\left(t^*(\hat{\theta}_i^*)\right)$ |
|---|---|---|---|---|
| 1 | $0.9956068 + 0.9848517i$ | $0.5617318$ | $0.07985537$ | $1.240413 \times 10^{-11}$ |
| 2 | $-0.8024712 + 0.6100954i$ | $2.9235813$ | $0.03959989$ | $2.192571 \times 10^{-11}$ |
| 3 | $0.2151482 - 0.2343377i$ | $-1.3754462$ | $0.02314992$ | $4.986238 \times 10^{-11}$ |
| 4 | $-0.2341807 - 0.5345383i$ | $-1.8976400$ | $0.02691153$ | $7.985967 \times 10^{-11}$ |
| 5 | $0.0371744 - 0.4051295i$ | $-1.6316806$ | $0.04667642$ | $1.643973 \times 10^{-1}$ |

**Table 4.3**. Initial estimates for the roots of $p(z) = z^4 + C_1 z^3 + C_2 z^2 + C_3 z + C_4$ with coefficients $C_1$, $C_2$, $C_3$, $C_4$ given by (4.13). These estimates were obtained from the discrete map $\frac{d}{d\theta} t^*(\theta)$ in figure 4.12.

| $i$ | $\hat{R}_i$ | $\hat{\theta}_i^*$ | $\left\lvert \Delta \frac{d}{d\theta} t^*(\hat{\theta}_i^*) \right\rvert$ | $d^2_{\hat{\theta}_i^*}\left(t^*(\hat{\theta}_i^*)\right)$ |
|---|---|---|---|---|
| 1 | $-0.2341886 - 0.5345356i$ | $-1.8976484$ | 0.006510733 | $1.148044 \times 10^{-10}$ |
| 2 | $0.9956081 + 0.9848495i$ | $0.5617294$ | 0.008003726 | $1.176466 \times 10^{-10}$ |
| 3 | $0.2151742 - 0.2343327i$ | $-1.3754060$ | 0.003277041 | $7.394599 \times 10^{-9}$ |
| 4 | $-0.8024620 + 0.6101369i$ | $2.9235347$ | 0.005600884 | $1.605620 \times 10^{-8}$ |
| 5 | $0.3385304 - 0.2043986i$ | $-1.1840143$ | 0.002358787 | $1.634235 \times 10^{-1}$ |

We see that values $d^2_{\hat{\theta}_i^*}\left(t^*(\hat{\theta}_i^*)\right)$ in the first four rows of tables 4.2 and 4.3 are much smaller in comparison with value $d^2_{\hat{\theta}_5^*}\left(t^*(\hat{\theta}_5^*)\right)$ from the last row in each of these tables; if we also compare the approximations contained in the first four rows of tables 4.2 and 4.3 with reference values (4.12), we see that these approximations correspond to the true roots. **Since values $d^2_{\hat{\theta}_i^*}\left(t^*(\hat{\theta}_i^*)\right)$ help us identify candidates for approximations to true roots of a polynomial, we can use them as a quality indicator for approximation $\hat{\theta}_i^*$: the smaller the value of $d^2_{\hat{\theta}_i^*}\left(t^*(\hat{\theta}_i^*)\right)$, the better the quality of $\hat{\theta}_i^*$.**

Next, we compute the absolute relative errors for the first four approximations in table 4.2 with respect to reference values (4.12):

$$
\begin{aligned}
&\left|\frac{R_1-\hat{R}_2}{R_1}\right| = 1.508607 \times 10^{-6} &\quad& \left|\frac{\theta_1^*-\hat{\theta}_2^*}{\theta_1^*}\right| = 4.966485 \times 10^{-7} \\
&\left|\frac{R_2-\hat{R}_3}{R_2}\right| = 6.826323 \times 10^{-6} && \left|\frac{\theta_2^*-\hat{\theta}_3^*}{\theta_2^*}\right| = 1.634514 \times 10^{-7} \\
&\left|\frac{R_3-\hat{R}_1}{R_3}\right| = 6.776553 \times 10^{-7} && \left|\frac{\theta_3^*-\hat{\theta}_1^*}{\theta_3^*}\right| = 3.697053 \times 10^{-7} \\
&\left|\frac{R_4-\hat{R}_4}{R_4}\right| = 6.620055 \times 10^{-6} && \left|\frac{\theta_4^*-\hat{\theta}_4^*}{\theta_4^*}\right| = 1.994664 \times 10^{-6}
\end{aligned}
\tag{4.15}
$$

Similarly, we compute the absolute relative errors for the first four approximations in table 4.3 with respect to reference values (4.12):

$$
\begin{aligned}
&\left|\frac{R_1-\hat{R}_4}{R_1}\right| = 4.082394 \times 10^{-5} &\quad& \left|\frac{\theta_1^*-\hat{\theta}_4^*}{\theta_1^*}\right| = 1.544730 \times 10^{-5} \\
&\left|\frac{R_2-\hat{R}_3}{R_2}\right| = 8.312983 \times 10^{-5} && \left|\frac{\theta_2^*-\hat{\theta}_3^*}{\theta_2^*}\right| = 2.902837 \times 10^{-5} \\
&\left|\frac{R_3-\hat{R}_2}{R_3}\right| = 2.086966 \times 10^{-6} && \left|\frac{\theta_3^*-\hat{\theta}_2^*}{\theta_3^*}\right| = 4.598819 \times 10^{-6} \\
&\left|\frac{R_4-\hat{R}_1}{R_4}\right| = 7.937340 \times 10^{-6} && \left|\frac{\theta_4^*-\hat{\theta}_1^*}{\theta_4^*}\right| = 2.404189 \times 10^{-6}
\end{aligned}
\tag{4.16}
$$

From values $d^2_{\hat{\theta}_i^*}\left(t^*(\hat{\theta}_i^*)\right)$ in the first four rows of tables 4.1, 4.2, and 4.3, and from the absolute relative errors in (4.14), (4.15), and (4.16), we see that, ultimately, each of the three discrete proximity maps $e(\theta)$, $\frac{d}{d\theta}d^2_\theta\big(t^*(\theta)\big)$, and $\frac{d}{d\theta}t^*(\theta)$, produce initial approximations of acceptable quality. If we compute arithmetic means and standard deviations for these error measures, we obtain the results shown in table 4.4.

**Table 4.4**. Arithmetic means and standard deviations for error measures associated with valid initial approximations (first four rows in each of tables 4.1, 4.2 and 4.3) to the roots of $p(z) = z^4 + C_1z^3 + C_2z^2 + C_3z + C_4$ with coefficients $C_1$, $C_2$, $C_3$, $C_4$ given by (4.13), obtained by means of discrete proximity maps $e(\theta)$, $\frac{d}{d\theta}d^2_\theta\big(t^*(\theta)\big)$, and $\frac{d}{d\theta}t^*(\theta)$ from figures 4.7, 4.11 and 4.12.

| Map | Error measure | | | | | |
|---|---|---|---|---|---|---|
| | $d^2_{\hat{\theta}_i^*}\left(t^*(\hat{\theta}_i^*)\right)$ | | $\left\lvert\frac{R-\hat{R}}{R}\right\rvert$ | | $\left\lvert\frac{\theta^*-\hat{\theta}^*}{\theta^*}\right\rvert$ | |
| | Mean | Standard deviation | Mean | Standard deviation | Mean | Standard deviation |
| $e(\theta)$ | $8.678753 \times 10^{-10}$ | $1.423474 \times 10^{-9}$ | $8.839238 \times 10^{-6}$ | $7.462444 \times 10^{-6}$ | $4.684063 \times 10^{-6}$ | $3.181745 \times 10^{-6}$ |
| $\frac{d}{d\theta}d^2_\theta\big(t^*(\theta)\big)$ | $4.101297 \times 10^{-11}$ | $3.038732 \times 10^{-11}$ | $3.908160 \times 10^{-6}$ | $3.269254 \times 10^{-6}$ | $7.561173 \times 10^{-7}$ | $8.370362 \times 10^{-7}$ |
| $\frac{d}{d\theta}t^*(\theta)$ | $5.920812 \times 10^{-9}$ | $7.578139 \times 10^{-9}$ | $3.349452 \times 10^{-5}$ | $3.722449 \times 10^{-5}$ | $1.286967 \times 10^{-5}$ | $1.218853 \times 10^{-5}$ |

The results from table 4.4 tell us that the initial approximations obtained with discrete map $\frac{d}{d\theta}d^2_\theta\big(t^*(\theta)\big)$ are slightly more consistent (smaller errors on average with smaller variability) than the initial approximations obtained with discrete map $e(\theta)$, which in turn are slightly more consistent than the initial approximations obtained with discrete map $\frac{d}{d\theta}t^*(\theta)$.

We see then, in this case, that any of the discrete maps $e(\theta)$, $\frac{d}{d\theta}d^2_\theta\big(t^*(\theta)\big)$, and $\frac{d}{d\theta}t^*(\theta)$, is in itself a useful tool that helps us identify reasonable initial approximations to the roots $R_i$, $\theta_i^*$ of the equation of the form (4.1) with coefficients (4.13); even with the map $\frac{d}{d\theta}t^*(\theta)$, we see that the real and imaginary parts of their associated initial estimates, shown in the first four rows of table 4.3, match at least the first four digits of the corresponding components of the true roots in (4.12).

## Numerical Example 4.2

In this example we proceed in exactly the same way as we did in example 4.1, using the same operation parameters for the script in annex 3 section 3, but this time initializing the random number generator seed by means of the instruction `set.seed(2020)` at the beginning of said script.

The roots $R_i$ generated this time, together with their corresponding theta roots $\theta_i^*$ (which depend on $P_1$, the semi-sum of all roots $R_i$), are:

$$
\begin{aligned}
R_1 &= \phantom{-}0.2938057 - 0.2115485i & \theta_1^* &= \phantom{-}0.5393991 \\
R_2 &= \phantom{-}0.2370036 - 0.0462177i & \theta_2^* &= \phantom{-}0.7221141 \\
R_3 &= -0.7278056 - 0.8652312i & \theta_3^* &= -2.4906911 \\
R_4 &= -0.7416948 - 0.2137641i & \theta_4^* &= \phantom{-}2.1105350
\end{aligned}
\tag{4.17}
$$

The coefficients constructed from Vieta's relations (4.2-a), (4.2-b), (4.2-c) and (4.2-d) are:

$$
\begin{aligned}
C_1 &= \phantom{-}0.93869109 + 1.33676158i \\
C_2 &= -0.64344310 + 0.53964620i \\
C_3 &= -0.23717297 - 0.36080124i \\
C_4 &= \phantom{-}0.07204248 + 0.02511377i
\end{aligned}
\tag{4.18}
$$

The global proximity LC map constructed from $N = 2{,}500$ elements $LzC(C_1, C_2, C_3, C_4, \theta_k)$ associated with points $\theta_k$ in a regular partition of interval $[-\pi, \pi)$, together with reference values $\theta_i^*$ in (4.17), is shown in figure 4.13.

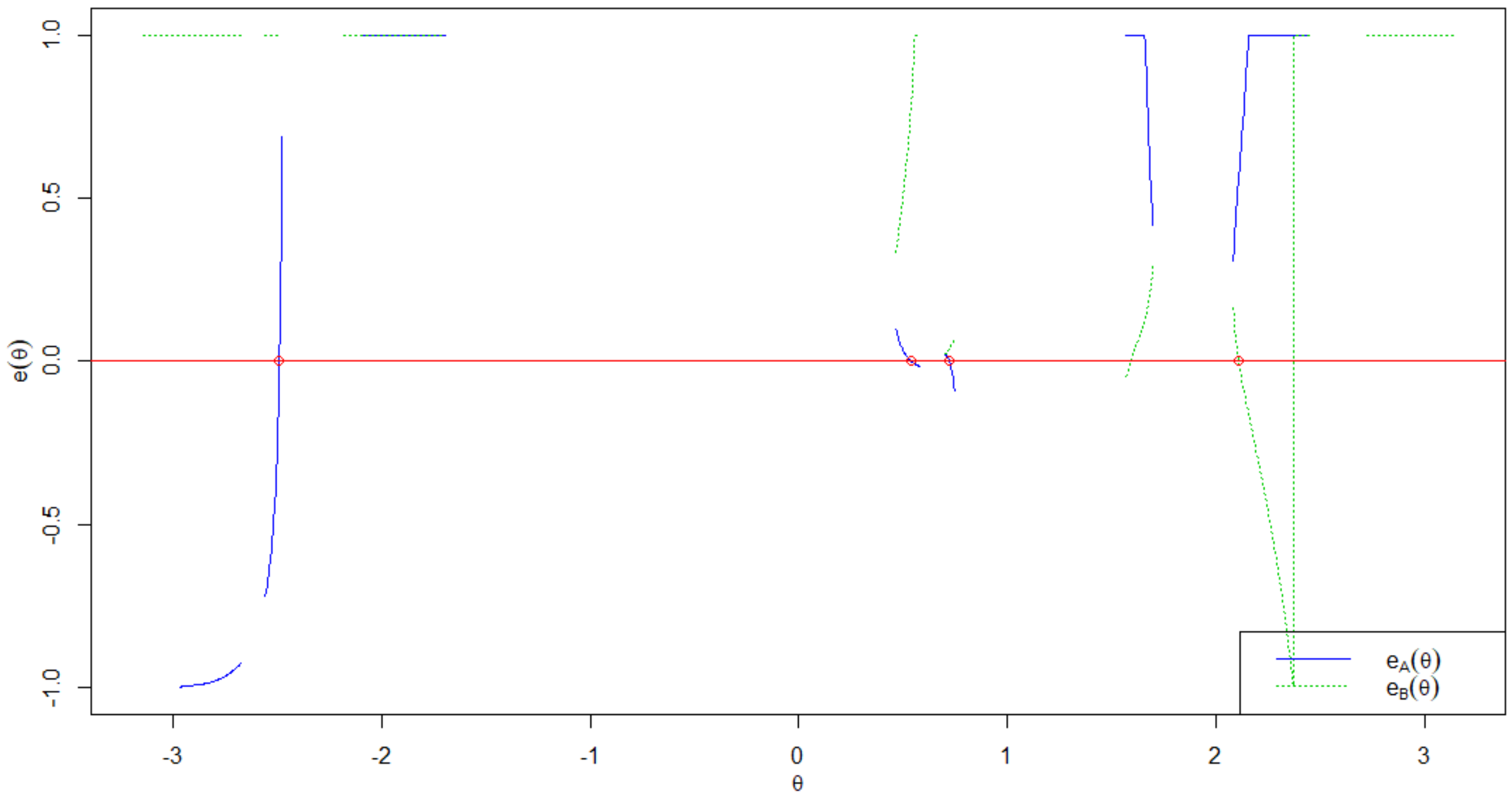

**Fig. 4.13.** Discrete angular proximity LC map for polynomial $p(z) = z^4 + C_1 z^3 + C_2 z^2 + C_3 z + C_4$ with coefficients $C_1$, $C_2$, $C_3$, $C_4$ given by (4.18). This map was generated from $N = 2{,}500$ elements $LzC(\theta_k)$ associated with points $\theta_k$ in a regular partition of interval $[-\pi, \pi)$. True theta roots $\theta_i^*$ are shown by means of small circles on horizontal axis $y = 0$.

As we can see in figure 4.13, in this case there are more than four smooth crossings $\hat{\theta}_i^*$ of functions $\hat{e}_A(\theta)$ and $\hat{e}_B(\theta)$ with horizontal axis $y = 0$; table 4.5 shows numerical details for these smooth crossings $\hat{\theta}_i^*$, ranked according to the criterion $d^2_{\hat{\theta}_i^*}\left(t^*(\hat{\theta}_i^*)\right)$ defined in example 4.1 (small values in column $d^2_{\hat{\theta}_i^*}\left(t^*(\hat{\theta}_i^*)\right)$ are better, and therefore appear in higher positions in the table).

**Table 4.5**. Initial estimates for the roots of $p(z) = z^4 + C_1 z^3 + C_2 z^2 + C_3 z + C_4$ with coefficients $C_1$, $C_2$, $C_3$, $C_4$ given by (4.18). These estimates were obtained from the discrete LC map in figure 4.13.

| $i$ | $\hat{R}_i$ | $\hat{\theta}_i^*$ | $\lvert\Delta e_i\rvert$ | $d^2_{\hat{\theta}_i^*}\left(t^*(\hat{\theta}_i^*)\right)$ |
|---|---|---|---|---|
| 1 | $-0.7278045 - 0.8652300i$ | $-2.4906919$ | $0.078488211$ | $2.184817 \times 10^{-12}$ |
| 2 | $-0.7416957 - 0.2137667i$ | $2.1105391$ | $0.010185580$ | $5.775807 \times 10^{-12}$ |
| 3 | $0.2937995 - 0.2115369i$ | $0.5394138$ | $0.001390601$ | $6.339265 \times 10^{-11}$ |
| 4 | $0.2370268 - 0.0462459i$ | $0.7220753$ | $0.004903348$ | $1.101618 \times 10^{-9}$ |
| 5 | $-0.4747049 - 0.4106948i$ | $1.5915914$ | $0.004768822$ | $3.046815 \times 10^{-2}$ |

From table 4.5, we see that the first four rows, ordered from smallest to largest value in column $d^2_{\hat{\theta}_i^*}\left(t^*(\hat{\theta}_i^*)\right)$, contain, in effect, the approximations to the true roots shown in reference values (4.17). To know the degree of accuracy for the approximations $\hat{R}_i$, $\hat{\theta}_i^*$ in the first four rows of table 4.5, we compute their absolute relative errors with respect to reference values (4.17):

$$\begin{aligned}
&\left|\frac{R_1-\hat{R}_3}{R_1}\right| = 3.626616 \times 10^{-5} &\qquad& \left|\frac{\theta_1^*-\hat{\theta}_3^*}{\theta_1^*}\right| = 2.733656 \times 10^{-5} \\
&\left|\frac{R_2-\hat{R}_4}{R_2}\right| = 1.511324 \times 10^{-4} && \left|\frac{\theta_2^*-\hat{\theta}_4^*}{\theta_2^*}\right| = 5.365845 \times 10^{-5} \\
&\left|\frac{R_3-\hat{R}_1}{R_3}\right| = 1.477850 \times 10^{-6} && \left|\frac{\theta_3^*-\hat{\theta}_1^*}{\theta_3^*}\right| = 3.146420 \times 10^{-7} \\
&\left|\frac{R_4-\hat{R}_2}{R_4}\right| = 3.587106 \times 10^{-6} && \left|\frac{\theta_4^*-\hat{\theta}_2^*}{\theta_4^*}\right| = 1.934729 \times 10^{-6}
\end{aligned} \qquad (4.19)$$

We see that the orders of magnitude of absolute relative errors (4.19) are similar to the orders of magnitude of absolute relative errors (4.14) computed in numerical example 4.1.

Let us now look at the discrete graphs related to the optimal values (global minima) of the dynamic squared distances $d^2_\theta$ associated with the structures $LzC$ in this example: figure 4.14 shows the tracings of trajectories $aP(\theta)$ and $R_x(\theta)$; figure 4.15 shows function $d^2_\theta\left(t^*(\theta)\right)$; figure 4.16 shows function $t^*(\theta)$, figure 4.17 shows function $\frac{d}{d\theta} d^2_\theta\left(t^*(\theta)\right)$, and figure 4.18 shows function $\frac{d}{d\theta} t^*(\theta)$. From this set of discrete graphs, we see that any of them contains five continuous sections (considering the periodicity of functions $d^2_\theta\left(t^*(\theta)\right)$ and $t^*(\theta)$, as well as of their derivatives with respect to $\theta$), although in figure 4.14, for both trajectories $aP(\theta)$ and $R_x(\theta)$, the ends of two of these five sections are joined at the point indicated by a diamond ◈.

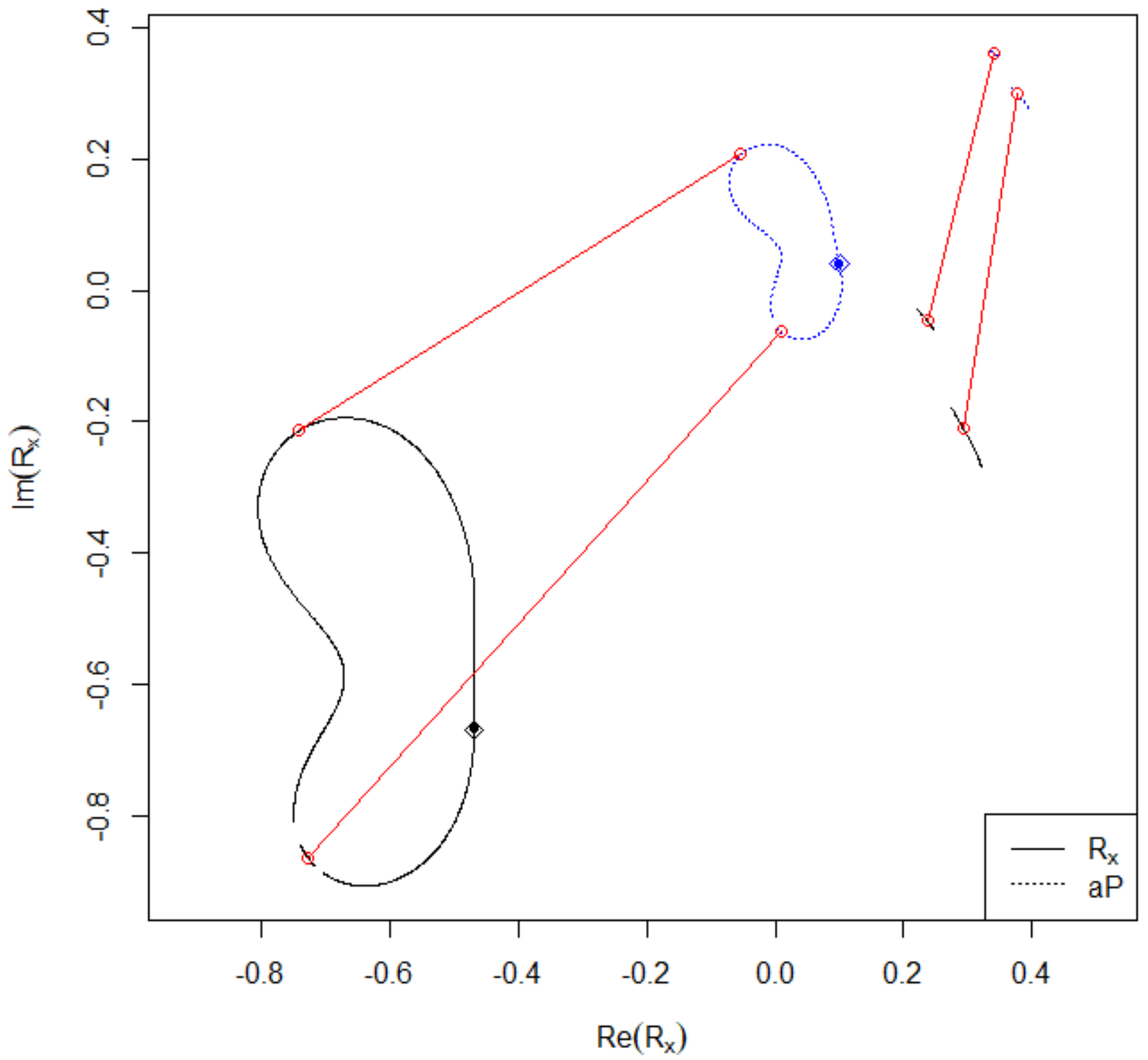

**Fig. 4.14.** Trajectory of best approximations $R_x(\theta_k)$ (solid line), and trajectory of anchor points $aP(\theta_k)$ for terminal semi-lines $tL(\theta_k)$ (dotted line). The trajectory $R_x(\theta_k)$ shown here contains the true roots of polynomial $p(z) = z^4 + C_1 z^3 + C_2 z^2 + C_3 z + C_4$ with coefficients $C_1$, $C_2$, $C_3$, $C_4$ given by (4.18). The true roots are represented in this graph by means of small circles, which are joined, by means of line segments, with their corresponding anchor points.

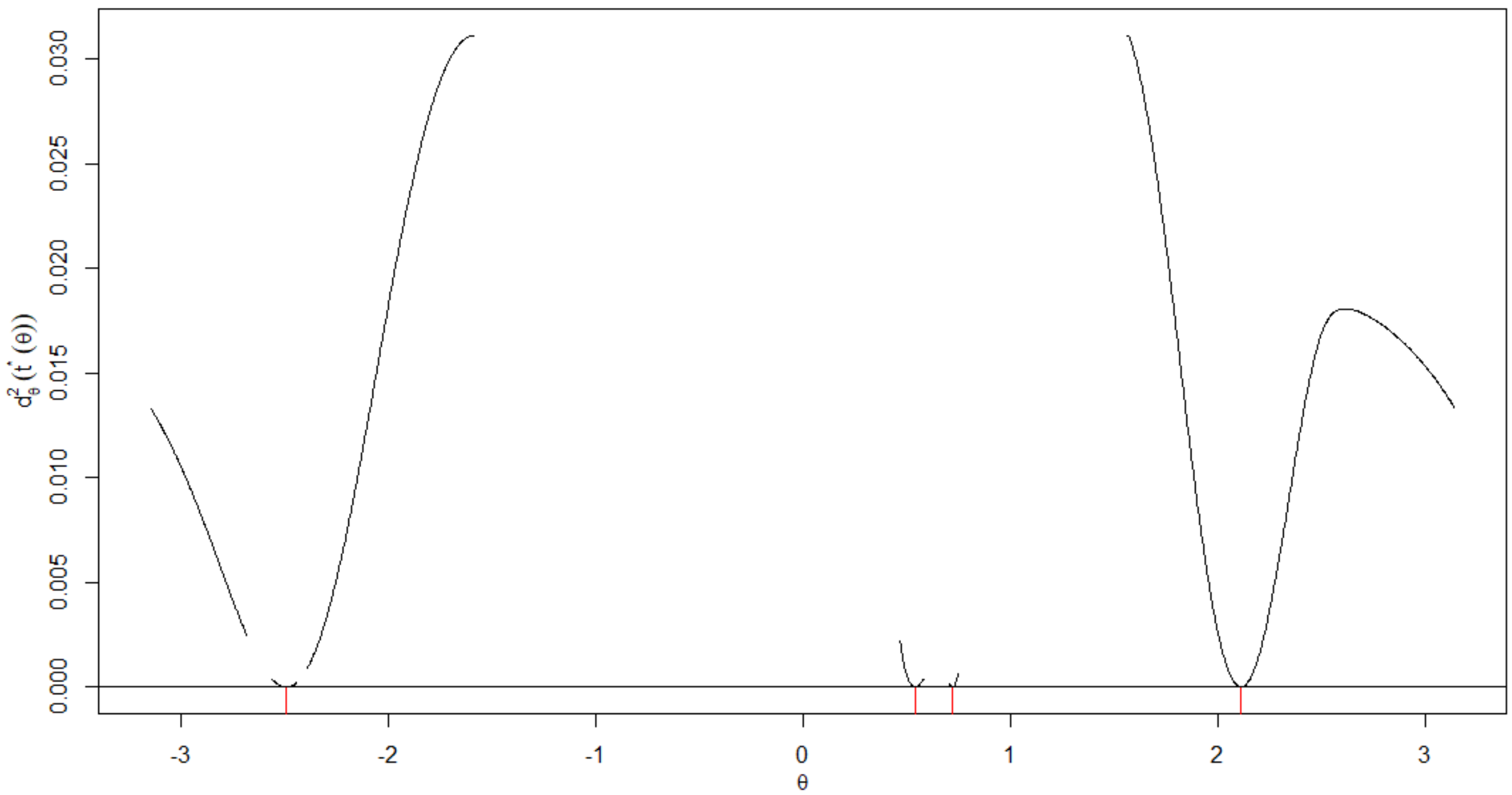

**Fig. 4.15.** Global minima of the dynamic squared distance functions $d^2_{\theta_k}(t^*(\theta_k)) = \min_t d^2_{\theta_k}(t)$ associated to polynomial $p(z) = z^4 + C_1 z^3 + C_2 z^2 + C_3 z + C_4$ with coefficients $C_1$, $C_2$, $C_3$, $C_4$ given by (4.18), vs. values $\theta_k$. At the bottom of the graph there are vertical line segments which indicate the location of true theta roots $\theta_i^*$.

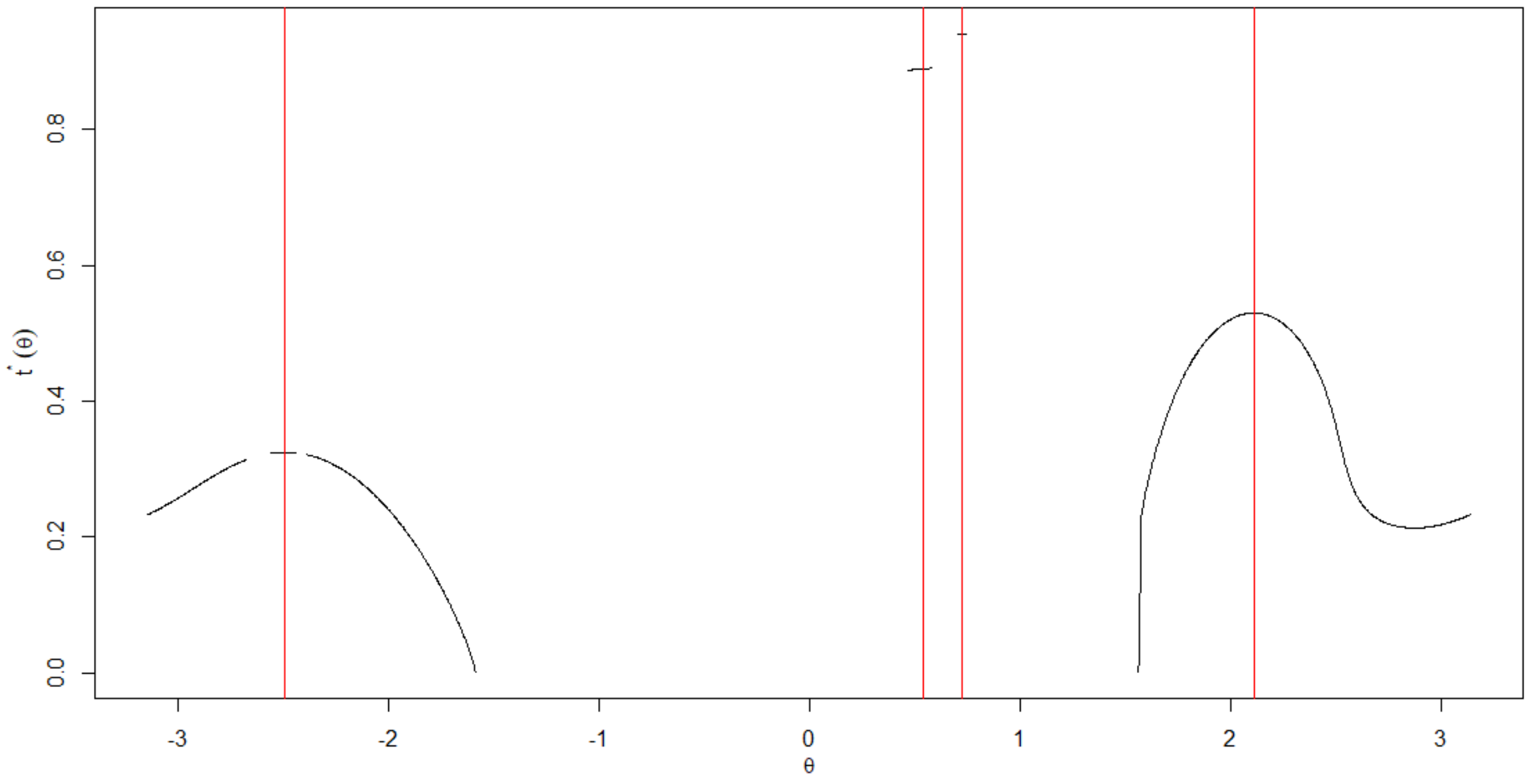

**Fig. 4.16.** Minimizing arguments of dynamic squared distance functions $t^*(\theta_k) = \arg\min_t d^2_{\theta_k}(t)$ associated to polynomial $p(z) = z^4 + C_1 z^3 + C_2 z^2 + C_3 z + C_4$ with coefficients $C_1$, $C_2$, $C_3$, $C_4$ given by (4.18), vs. values $\theta_k$. This graph includes vertical line segments which indicate the location of true theta roots $\theta_i^*$.

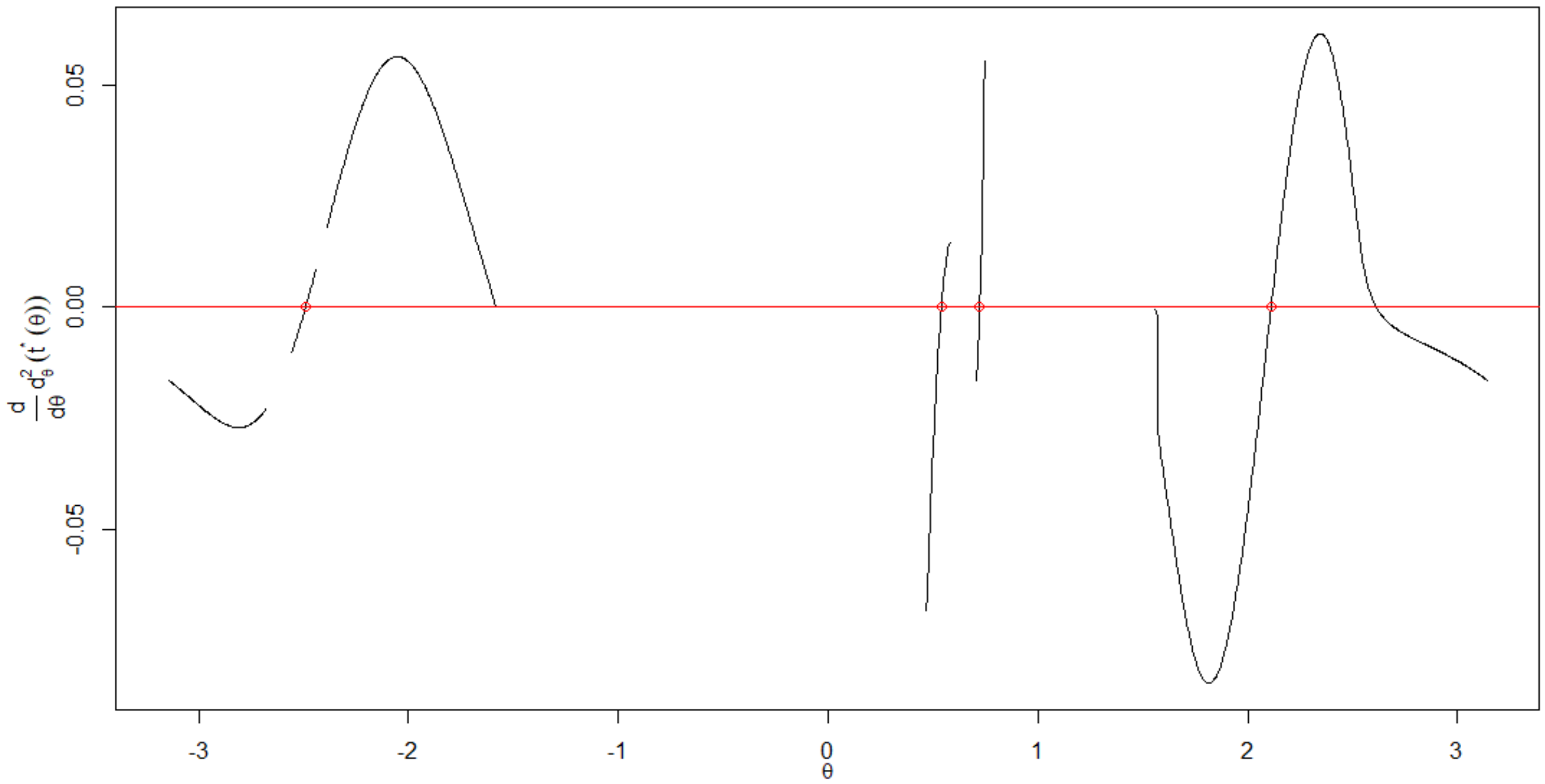

**Fig. 4.17**. Discrete approximation to the derivative of function $d^2_\theta(t^*(\theta))$ from figure 4.15 with respect to $\theta$; this function is associated to polynomial $p(z) = z^4 + C_1 z^3 + C_2 z^2 + C_3 z + C_4$ with coefficients $C_1$, $C_2$, $C_3$, $C_4$ given by (4.18). Small circles on horizontal axis $y = 0$ indicate the location of true theta roots $\theta_i^*$.

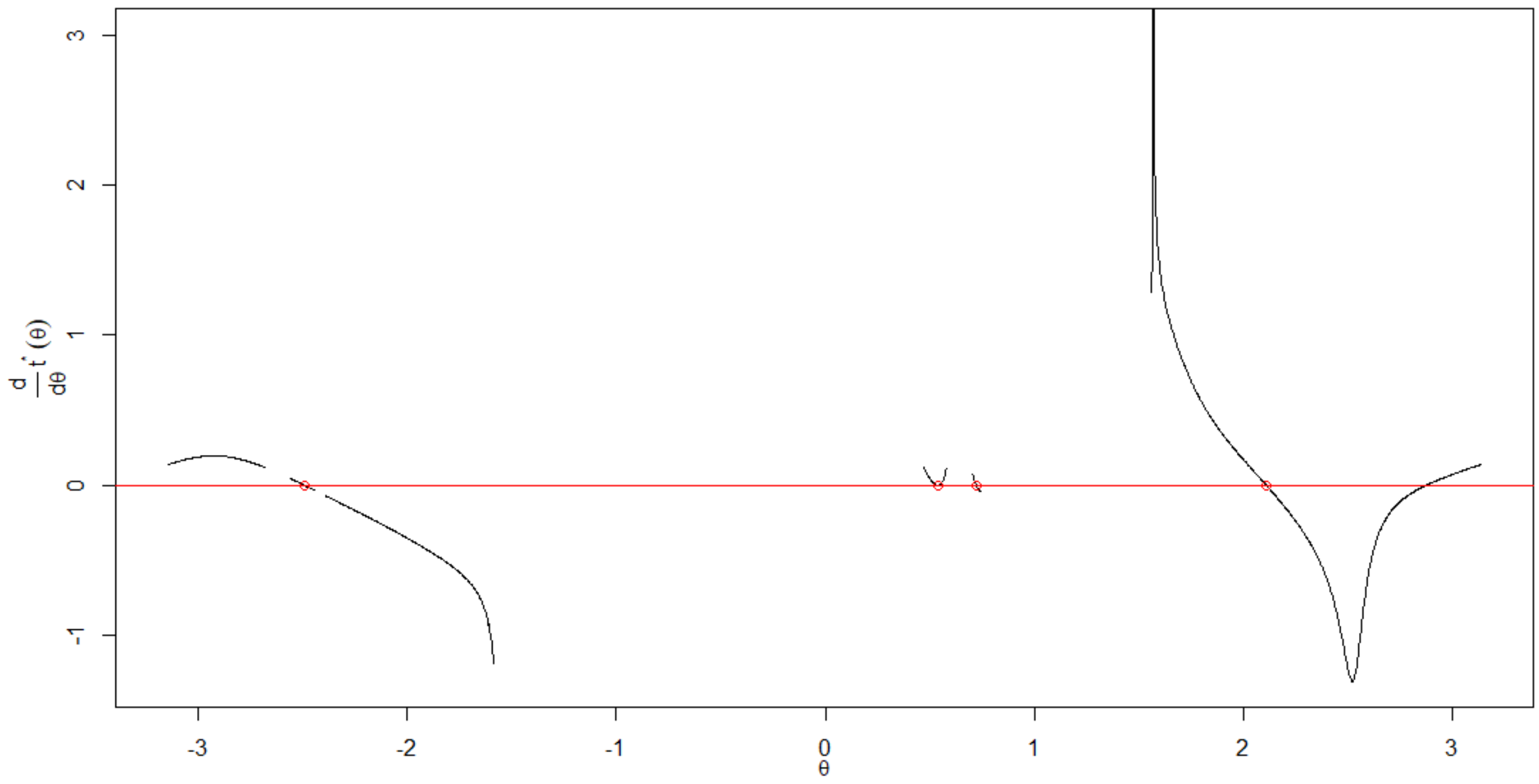

**Fig. 4.18.** Discrete approximation to the derivative of function $t^*(\theta)$ from figure 4.16 with respect to $\theta$; this function is associated to polynomial $p(z) = z^4 + C_1 z^3 + C_2 z^2 + C_3 z + C_4$ with coefficients $C_1$, $C_2$, $C_3$, $C_4$ given by (4.18). Small circles on horizontal axis $y = 0$ indicate the location of true theta roots $\theta_i^*$.

By analyzing the numerical arrays `aP`, `Rx`, `MinF`, `MinT`, `DD_MinF`, `DD_MinT` (generated by the script of annex 3 section 3), which respectively contain the discrete versions of $aP(\theta), R_x(\theta)$, $d_\theta^2(t^*(\theta))$, $t^*(\theta)$, $\frac{d}{d\theta} d_\theta^2(t^*(\theta))$ and $\frac{d}{d\theta} t^*(\theta)$, together with their corresponding graphs in figures 4.14, 4.15, 4.16, 4.17 and 4.18, we can see that, for this particular example, the observed behavior seems to suggest the following:

- The ends of the two continuous sections that appear joined in the trajectories $aP(\theta)$ and $R_x(\theta)$ from figure 4.14, exhibit a separation of $\pi$ radians in the remaining graphs (figures 4.15, 4.16, 4.17 and 4.18).
- In the graph of function $d_\theta^2(t^*(\theta))$ (figure 4.15), these two ends with separation of $\pi$ radians, both correspond to the maximum functional value of $d_\theta^2(t^*(\theta))$.
- In the graph of $t^*(\theta)$ (figure 4.16), these two ends with separation of $\pi$ radians, both correspond to the minimum functional value of $t^*(\theta)$, which is close to zero, but positive.
- In the graph of $\frac{d}{d\theta} d_\theta^2(t^*(\theta))$ (figure 4.17), these two ends with separation of $\pi$ radians, both correspond to functional values of $\frac{d}{d\theta} d_\theta^2(t^*(\theta))$ that tend to zero: one through positive values, and other through negative values.
- In the graph of $\frac{d}{d\theta} t^*(\theta)$ (figure 4.18), these two ends with separation of $\pi$ radians, both correspond to functional values of $\frac{d}{d\theta} t^*(\theta)$ of equal magnitude, but opposite signs.

This observed behavior for the continuous sections in figures 4.14, 4.15, 4.16, 4.17 and 4.18 (as well as the fact that half of the support values $\theta$, of total length $\pi$, are associated with valid functional values) is a consequence both of the symmetry property of function $d_\theta^2(t)$ (described in example 4.1), and of the constraint $t^*(\theta) > 0$.

In this numerical example, the joint behaviors observed in the graphs from figures 4.15 and 4.17, as well as in the graphs from figures 4.16 and 4.18, reinforce the hypotheses about the properties posed in example 4.1 for mappings $d_\theta^2(t^*(\theta))$ and $t^*(\theta)$. From figures 4.16 and 4.18, it seems that theta root $\theta_1^* = 0.5393991$ corresponds to a saddle point of $t^*(\theta)$, while theta root $\theta_2^* = 0.7221141$ apparently corresponds to a local maximum of $t^*(\theta)$; unfortunately, the corresponding continuous sections cannot be clearly seen in figures 4.16 and 4.18, given their relatively reduced length at global level; later in this numerical example we will zoom in on these two continuous sections, in order to corroborate their qualitative behavior.

From the discrete map $\frac{d}{d\theta} d_\theta^2(t^*(\theta))$ shown in figure 4.17, we obtain 5 smooth crossings $\hat{\theta}_i^*$ with horizontal axis $y = 0$; the numerical characteristics of these smooth crossings $\hat{\theta}_i^*$ are shown in table 4.6, where rows are sorted in ascending order according to values in column $d_{\hat{\theta}_i^*}^2\left(t^*(\hat{\theta}_i^*)\right)$.

**Table 4.6**. Initial estimates for the roots of $p(z) = z^4 + C_1 z^3 + C_2 z^2 + C_3 z + C_4$ with coefficients $C_1$, $C_2$, $C_3$, $C_4$ given by (4.18). These estimates were obtained from the map $\frac{d}{d\theta} d_\theta^2(t^*(\theta))$ in figure 4.17.

| $i$ | $\hat{R}_i$ | $\hat{\theta}_i^*$ | $\left\lvert\Delta \frac{d}{d\theta} d_{\hat{\theta}_i^*}^2(t^*)\right\rvert$ | $d_{\hat{\theta}_i^*}^2\left(t^*(\hat{\theta}_i^*)\right)$ |
|---|---|---|---|---|
| 1 | $-0.7278047 - 0.8652297i$ | $-2.4906932$ | $0.0004152196$ | $2.443717 \times 10^{-12}$ |
| 2 | $-0.7416941 - 0.2137659i$ | $2.1105357$ | $0.0010630391$ | $2.703311 \times 10^{-12}$ |
| 3 | $0.2938010 - 0.2115394i$ | $0.5394105$ | $0.0014360912$ | $3.823739 \times 10^{-11}$ |
| 4 | $0.2370238 - 0.0462426i$ | $0.7220801$ | $0.0035551513$ | $8.481752 \times 10^{-10}$ |
| 5 | $-0.6882267 - 0.5402932i$ | $2.6121323$ | $0.0002175127$ | $1.807938 \times 10^{-2}$ |

Similarly, table 4.7 shows the numerical characteristics of the smooth crossings $\hat{\theta}_i^*$ detected in the discrete map $\frac{d}{d\theta} t^*(\theta)$ from figure 4.18.

**Table 4.7**. Initial estimates for the roots of $p(z) = z^4 + C_1 z^3 + C_2 z^2 + C_3 z + C_4$ with coefficients $C_1$, $C_2$, $C_3$, $C_4$ given by (4.18). These estimates were obtained from the map $\frac{d}{d\theta} t^*(\theta)$ in figure 4.18.

| $i$ | $\hat{R}_i$ | $\hat{\theta}_i^*$ | $\left\lvert\Delta \frac{d}{d\theta} t^*(\hat{\theta}_i^*)\right\rvert$ | $d_{\hat{\theta}_i^*}^2\left(t^*(\hat{\theta}_i^*)\right)$ |
|---|---|---|---|---|
| 1 | $0.2370005 - 0.0462161i$ | $0.7221175$ | $0.0089308208$ | $1.006958 \times 10^{-11}$ |
| 2 | $-0.7416902 - 0.2137636i$ | $2.1105271$ | $0.0039440589$ | $1.590274 \times 10^{-11}$ |
| 3 | $-0.7278143 - 0.8652171i$ | $-2.4907418$ | $0.0016334961$ | $2.150274 \times 10^{-10}$ |
| 4 | $0.2936607 - 0.2113052i$ | $0.5397174$ | $0.0009564820$ | $2.942764 \times 10^{-8}$ |
| 5 | $0.2946947 - 0.2130361i$ | $0.5374507$ | $0.0003602247$ | $1.117984 \times 10^{-6}$ |
| 6 | $-0.6750676 - 0.6117350i$ | $2.8729002$ | $0.0017016144$ | $1.664505 \times 10^{-2}$ |

From tables 4.6 and 4.7, we see that the estimates in the first four rows are indeed the ones that best approximate reference values (4.17); in table 4.7 we see that rows 4 and 5 contain values $\hat{\theta}_4^*$ and $\hat{\theta}_5^*$ similar to each other, which are close to theta root $\theta_1^* = 0.5393991$, so that $\hat{\theta}_5^* < \theta_1^* < \hat{\theta}_4^*$; in figure 4.18, $\hat{\theta}_4^*$ and $\hat{\theta}_5^*$ are so close to each other, that it is not possible to distinguish them from the true theta root $\theta_1^*$, which seems to coincide with the only point where the curve $\frac{d}{d\theta}t^*(\theta)$ touches, without crossing, the horizontal axis $y = 0$, precisely at a neighborhood around $\theta_1^*$.

We now compute the absolute relative errors for the initial estimates in the first four rows of table 4.6 with respect to reference values (4.17):

$$\begin{aligned} \left|\frac{R_1-\hat{R}_3}{R_1}\right| &= 2.816578 \times 10^{-5} & \left|\frac{\theta_1^*-\hat{\theta}_3^*}{\theta_1^*}\right| &= 2.121973 \times 10^{-5} \\ \left|\frac{R_2-\hat{R}_4}{R_2}\right| &= 1.326079 \times 10^{-4} & \left|\frac{\theta_2^*-\hat{\theta}_4^*}{\theta_2^*}\right| &= 4.706974 \times 10^{-5} \\ \left|\frac{R_3-\hat{R}_1}{R_3}\right| &= 1.562962 \times 10^{-6} & \left|\frac{\theta_3^*-\hat{\theta}_1^*}{\theta_3^*}\right| &= 8.410823 \times 10^{-7} \\ \left|\frac{R_4-\hat{R}_2}{R_4}\right| &= 2.454063 \times 10^{-6} & \left|\frac{\theta_4^*-\hat{\theta}_2^*}{\theta_4^*}\right| &= 3.219408 \times 10^{-7} \end{aligned} \qquad (4.20)$$

Similarly, we compute the absolute relative errors for the initial estimates in the first four rows of table 4.7:

$$\begin{aligned} \left|\frac{R_1-\hat{R}_4}{R_1}\right| &= 7.821215 \times 10^{-4} & \left|\frac{\theta_1^*-\hat{\theta}_4^*}{\theta_1^*}\right| &= 5.902113 \times 10^{-4} \\ \left|\frac{R_2-\hat{R}_1}{R_2}\right| &= 1.444473 \times 10^{-5} & \left|\frac{\theta_2^*-\hat{\theta}_1^*}{\theta_2^*}\right| &= 4.790197 \times 10^{-6} \\ \left|\frac{R_3-\hat{R}_3}{R_3}\right| &= 1.466145 \times 10^{-5} & \left|\frac{\theta_3^*-\hat{\theta}_3^*}{\theta_3^*}\right| &= 2.038348 \times 10^{-5} \\ \left|\frac{R_4-\hat{R}_2}{R_4}\right| &= 5.952142 \times 10^{-6} & \left|\frac{\theta_4^*-\hat{\theta}_2^*}{\theta_4^*}\right| &= 3.756020 \times 10^{-6} \end{aligned} \qquad (4.21)$$

If we compute arithmetic means and standard deviations for the first four values of column $d^2_{\hat{\theta}_i^*}\left(t^*(\hat{\theta}_i^*)\right)$ in tables 4.5, 4.6 and 4.7, as well as for the absolute relative errors in expressions (4.19), (4.20) and (4.21), we obtain the results shown in table 4.8; from here we see that the initial approximations associated with the map $\frac{d}{d\theta}d^2_\theta\left(t^*(\theta)\right)$ are slightly more consistent (smaller errors on average with smaller variability) than the approximations from map $e(\theta)$, and in turn the approximations from map $e(\theta)$ are slightly more consistent than the approximations from map $\frac{d}{d\theta}t^*(\theta)$. In the worst case, if we compare estimate $\hat{R}_4$ in fourth row of table 4.7 with corresponding reference value $R_1$ in (4.17), we see that the real and imaginary part of estimate $\hat{R}_4$ generated by map $\frac{d}{d\theta}t^*(\theta)$, each matches the first four digits of the corresponding components of the true root $R_1$; thus, estimate $\hat{R}_4$ from map $\frac{d}{d\theta}t^*(\theta)$ can still be thought of as a reasonable initial approximation. In conclusion, we consider that the global discrete maps $e(\theta)$, $\frac{d}{d\theta}d_\theta^2\big(t^*(\theta)\big)$, and $\frac{d}{d\theta}t^*(\theta)$ constructed in this numerical example, produce initial estimates of good quality.

**Table 4.8**. Arithmetic means and standard deviations for error measures associated with valid initial approximations (first four rows in each of tables 4.5, 4.6 and 4.7) to the roots of $p(z) = z^4 + C_1z^3 + C_2z^2 + C_3z + C_4$ with coefficients $C_1$, $C_2$, $C_3$, $C_4$ given by (4.18), obtained by means of discrete proximity maps $e(\theta)$, $\frac{d}{d\theta}d_\theta^2\big(t^*(\theta)\big)$, and $\frac{d}{d\theta}t^*(\theta)$ from figures 4.13, 4.17 and 4.18.

| Map | Error measure | | | | | |
|---|---|---|---|---|---|---|
| | $d^2_{\hat{\theta}_i^*}\left(t^*(\hat{\theta}_i^*)\right)$ | | $\left\lvert\frac{R-\hat{R}}{R}\right\rvert$ | | $\left\lvert\frac{\theta^*-\hat{\theta}^*}{\theta^*}\right\rvert$ | |
| | **Mean** | **Standard deviation** | **Mean** | **Standard deviation** | **Mean** | **Standard deviation** |
| $e(\theta)$ | $2.932428 \times 10^{-10}$ | $5.396459 \times 10^{-10}$ | $4.811588 \times 10^{-5}$ | $7.049998 \times 10^{-5}$ | $2.081110 \times 10^{-5}$ | $2.515255 \times 10^{-5}$ |
| $\frac{d}{d\theta}d_\theta^2\big(t^*(\theta)\big)$ | $2.228899 \times 10^{-10}$ | $4.171958 \times 10^{-10}$ | $4.119768 \times 10^{-5}$ | $6.217619 \times 10^{-5}$ | $1.736312 \times 10^{-5}$ | $2.206609 \times 10^{-5}$ |
| $\frac{d}{d\theta}t^*(\theta)$ | $7.417159 \times 10^{-9}$ | $1.467396 \times 10^{-8}$ | $2.042950 \times 10^{-4}$ | $3.852391 \times 10^{-4}$ | $1.547853 \times 10^{-4}$ | $2.903837 \times 10^{-4}$ |

**Close-ups of maps $t^*(\theta)$ and $\frac{d}{d\theta}t^*(\theta)$ around theta roots $\theta_1^*$ and $\theta_2^*$.** Let us now take a closer look at maps $t^*(\theta)$ and $\frac{d}{d\theta}t^*(\theta)$ in order to better appreciate their qualitative behavior in the vicinities of theta roots $\theta_1^* = 0.5393991$ and $\theta_2^* = 0.7221141$; by modifying a few operation parameters in the script from annex 3 section 3 (see "specific conditions for reproducing results from the examples in chapter 4" in that subsection for more details), we generate graphs for these maps in the region $[0.46, 0.75)$, which are shown in figures 4.19 and 4.20.

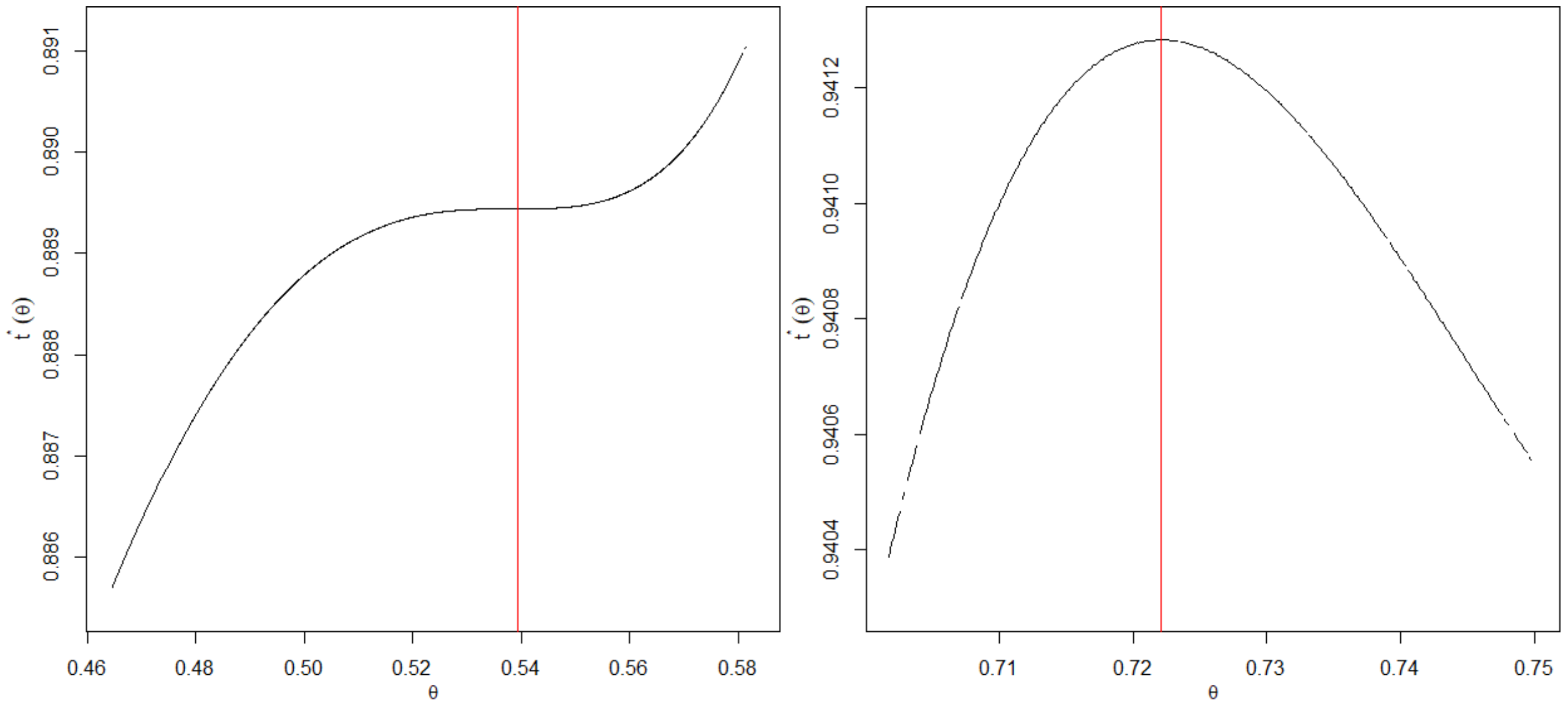

**Fig. 4.19.** Two close-ups of map $t^*(\theta)$ associated to the polynomial $p(z) = z^4 + C_1z^3 + C_2z^2 + C_3z + C_4$ with coefficients $C_1$, $C_2$, $C_3$, $C_4$ given by (4.18), in the vicinity of $\theta_1^* = 0.5393991$ and $\theta_2^* = 0.7221141$. These two close-ups include vertical line segments that indicate the location of the true theta roots $\theta_1^*$ and $\theta_2^*$. These two close-ups were generated from a single regional map $t^*(\theta)$ constructed with $N = 2{,}500$ elements $LzC(\theta_k)$ associated with points $\theta_k$ in a regular partition of interval $[0.46, 0.75)$.

The joint behavior of the graphs in figures 4.19 and 4.20 suggests that:

a) $\theta_1^*$ corresponds to a saddle point of function $t^*(\theta)$, since apparently, if we look beyond the noise present in figure 4.20, and we focus on the trend of graph $\frac{d}{d\theta}t^*(\theta)$, we see that $\frac{d}{d\theta}t^*(\theta_1^*) = 0$ and $\frac{d^2}{d\theta^2}t^*(\theta_1^*) = 0$;

b) $\theta_2^*$ corresponds to a local maximum of function $t^*(\theta)$, since it does appear that the trend of graph $\frac{d}{d\theta}t^*(\theta)$ in figure 4.20 indicates that $\frac{d}{d\theta}t^*(\theta_2^*) = 0$ and $\frac{d^2}{d\theta^2}t^*(\theta_2^*) < 0$.

The graph in figure 4.20, as we can see, contains highly oscillatory stochastic noise, which makes its analysis difficult; the presence of oscillatory noise in this discrete map $\frac{d}{d\theta}t^*(\theta)$ is due to the combined effect of two factors: 1) the quotient of first differences $\frac{\Delta t^*(\theta)}{\Delta\theta}$, which is what was used to numerically approximate the derivative of $t^*(\theta)$ with respect to $\theta$ (see the `approxDMin` function in annex 2 section 7), amplifies the errors (from the optimization process) contained in the discrete points $\hat{t}^*(\theta_k)$; 2) the high sampling rate (i.e., a very small value $\Delta\theta$), used in this close-up of map $\frac{d}{d\theta}t^*(\theta)$ to region $[0.46, 0.75)$.

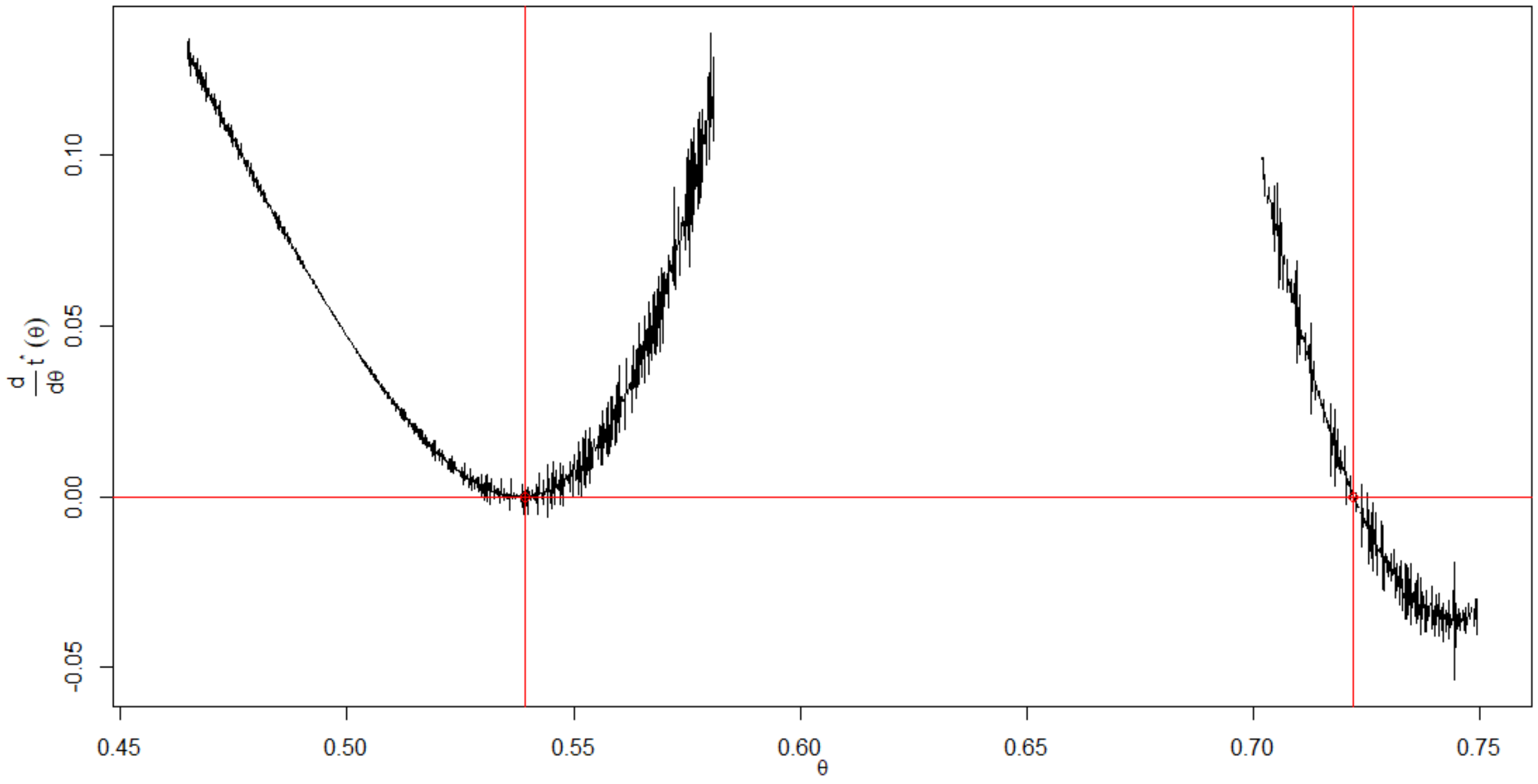

**Fig. 4.20.** Close-up of map $\frac{d}{d\theta}t^*(\theta)$ associated to the polynomial $p(z) = z^4 + C_1z^3 + C_2z^2 + C_3z + C_4$ with coefficients $C_1$, $C_2$, $C_3$, $C_4$ given by (4.18), in the vicinity of $\theta_1^* = 0.5393991$ and $\theta_2^* = 0.7221141$. This graph includes vertical line segments that indicate the location of the true theta roots $\theta_1^*$ and $\theta_2^*$. This close-up corresponds to a regional map $\frac{d}{d\theta}t^*(\theta)$ constructed with $N = 2{,}500$ elements $LzC(\theta_k)$ associated with points $\theta_k$ in a regular partition of interval $[0.46, 0.75)$.

From the empirical results obtained in this numerical example, we see that, **the higher the resolution of a discrete map $\frac{d}{d\theta}t^*(\theta)$ constructed with quotients of first differences $\frac{\Delta t^*(\theta)}{\Delta\theta}$, the greater the amount of oscillatory noise in this same map**.

Because of the oscillatory noise contained in the regional discrete map $\frac{d}{d\theta}t^*(\theta)$ shown in figure 4.20, a large number of spurious crossings with the horizontal axis $y = 0$ are generated; table 4.9 shows the numerical values of some of the crossings detected in the regional map from figure 4.20, ranked according to their associated values $d^2_{\hat{\theta}_i^*}\left(t^*\left(\hat{\theta}_i^*\right)\right)$. The first two rows of table 4.9 contain the best approximations to $\theta_1^* = 0.5393991$ and $\theta_2^* = 0.7221141$: such approximations are correct in their first 4 digits; the remaining rows in table 4.9 are "replicas" of the first two rows and are associated with spurious crossings. If we compare, by means of values $d^2_{\hat{\theta}_i^*}\left(t^*(\hat{\theta}_i^*)\right)$, the first two rows of table 4.9 with the corresponding rows 1 and 4 of table 4.7 associated with the global map $\frac{d}{d\theta}t^*(\theta)$ from figure 4.18, we see that the approximation to $\theta_2^*$ in table 4.7 (global map) is better, while the approximation to $\theta_1^*$ in table 4.9 (regional map) is better; in this case, as a matter of fact, we were fortunate that the global map $\frac{d}{d\theta}t^*(\theta)$ in figure 4.18 contained an approximation to $\theta_1^*$ at all, ironically thanks to the approximation errors (from the optimization process that seeks to estimate minimizing arguments $t^*$) present in the functional values of said map. In conclusion, **a high-resolution close-up of a discrete map $\frac{d}{d\theta}t^*(\theta)$ constructed with quotients of first differences $\frac{\Delta t^*(\theta)}{\Delta\theta}$, does not necessarily produce approximations of better quality, due to the possible presence of oscillatory noise**; fortunately, the close-up made in this case gives us the opportunity to corroborate the qualitative properties of mapping $t^*(\theta)$ listed in numerical example 4.1.

**Table 4.9**. Initial estimates for the roots of $p(z) = z^4 + C_1z^3 + C_2z^2 + C_3z + C_4$ with coefficients $C_1$, $C_2$, $C_3$, $C_4$ given by (4.18), obtained from the discrete map $\frac{d}{d\theta}t^*(\theta)$ in figure 4.20. A total of 71 estimates were produced; only a few of them are shown in this table.

| $i$ | $\hat{R}_i$ | $\hat{\theta}_i^*$ | $\left\lvert\Delta\frac{d}{d\theta}t^*(\hat{\theta}_i^*)\right\rvert$ | $d^2_{\hat{\theta}_i^*}\left(t^*(\hat{\theta}_i^*)\right)$ |
|---|---|---|---|---|
| 1 | $0.2370439 - 0.0462654i$ | 0.7220478 | 0.007162218 | $3.215141 \times 10^{-9}$ |
| 2 | $0.2937342 - 0.2114291i$ | 0.5395555 | 0.000285155 | $7.113933 \times 10^{-9}$ |
| 3 | $0.2938975 - 0.2117007i$ | 0.5391993 | 0.001172575 | $1.162821 \times 10^{-8}$ |
| ⋮ | ⋮ | ⋮ | ⋮ | ⋮ |
| 6 | $0.2367609 - 0.0459445i$ | 0.7225023 | 0.002161431 | $1.111888 \times 10^{-7}$ |
| ⋮ | ⋮ | ⋮ | ⋮ | ⋮ |
| 70 | $0.2890168 - 0.2035944i$ | 0.5498375 | 0.006647609 | $2.962778 \times 10^{-5}$ |
| 71 | $0.2890128 - 0.2035875i$ | 0.5498465 | 0.006645524 | $2.967724 \times 10^{-5}$ |

From here, a question arises: will the discrete map $\frac{d}{d\theta}d^2_\theta\big(t^*(\theta)\big)$ also contain oscillatory noise, if we build it on the same region $[0.46, 0.75)$ and under the same conditions as the maps from figures 4.19 and 4.20? Let us take a look at figures 4.21 and 4.22, which respectively show the discrete graphs of $d^2_\theta\big(t^*(\theta)\big)$ and $\frac{d}{d\theta}d^2_\theta\big(t^*(\theta)\big)$, each of them constructed from $N = 2{,}500$ elements $LzC(\theta_k)$ associated with points $\theta_k$ in a regular partition of interval $[0.46, 0.75)$. From these figures we can see that the answer to the question is negative; the discrete map $\frac{d}{d\theta}d^2_\theta\big(t^*(\theta)\big)$, at least in this example, is more robust to the presence of oscillatory noise compared to the discrete map $\frac{d}{d\theta}t^*(\theta)$.

The graphs in figures 4.21 and 4.22 help us to see with a little more detail what happens in figures 4.15 and 4.17, and their behavior reinforces the hypotheses posed in numerical example 4.1 on the properties of mapping $d_\theta^2\big(t^*(\theta)\big)$. Figures 4.21 and 4.22 suggest that: 1) theta roots $\theta_1^*$ and $\theta_2^*$ are minimizing arguments of function $d_\theta^2\big(t^*(\theta)\big)$; 2) $d_\theta^2\big(t^*(\theta_1^*)\big) = d_\theta^2\big(t^*(\theta_2^*)\big) = 0$; 3) function $d_\theta^2\big(t^*(\theta)\big)$ is concave upward at the neighborhoods around $\theta_1^*$ and $\theta_2^*$, which implies 4) $\frac{d}{d\theta} d_\theta^2\big(t^*(\theta_1^*)\big) = \frac{d}{d\theta} d_\theta^2\big(t^*(\theta_2^*)\big) = 0$, and 5) $\frac{d^2}{d\theta^2} d_\theta^2\big(t^*(\theta_1^*)\big) > 0, \frac{d^2}{d\theta^2} d_\theta^2\big(t^*(\theta_2^*)\big) > 0$.

The discrete regional map $\frac{d}{d\theta} d_\theta^2\big(t^*(\theta)\big)$ in figure 4.22 contains two smooth crossings $\hat{\theta}_i^*$ with horizontal axis $y = 0$. Are these two crossings better approximations compared to their global counterparts in table 4.6? Table 4.10 shows the numerical characteristics of the smooth crossings $\hat{\theta}_i^*$ associated with regional map $\frac{d}{d\theta} d_\theta^2\big(t^*(\theta)\big)$ in figure 4.22. When comparing these with their global counterparts in table 4.6, we see that the approximations in table 4.10 are more accurate, since their values $d_{\hat{\theta}_i^*}^2\left(t^*\big(\hat{\theta}_i^*\big)\right)$ are comparatively smaller. We therefore conclude that, in this case, **a high-resolution close-up of a discrete map $\frac{d}{d\theta} d_\theta^2\big(t^*(\theta)\big)$ produces more accurate approximations relative to its global counterpart**.

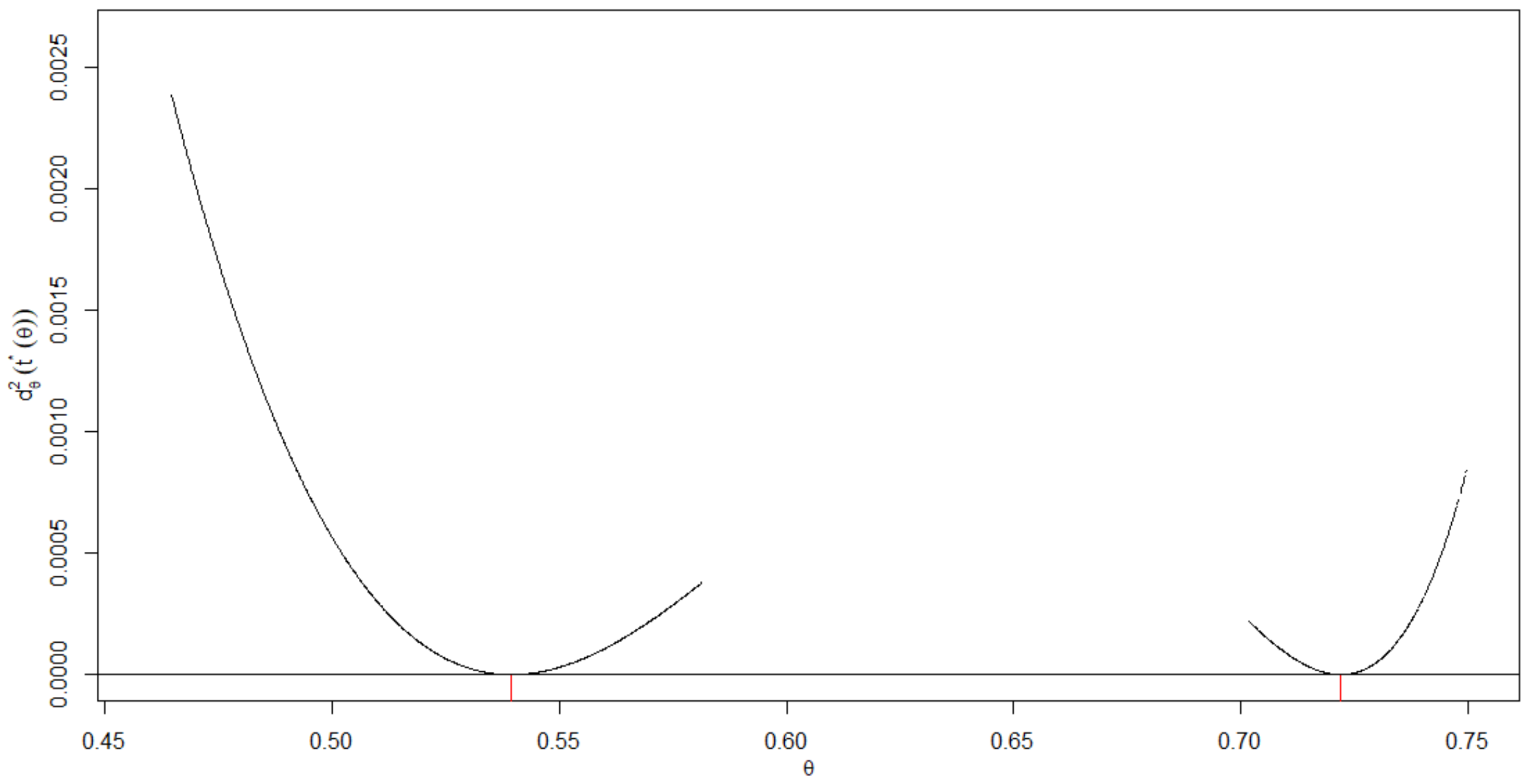

**Fig. 4.21.** Close-up of map $d_\theta^2\big(t^*(\theta)\big)$ associated to the polynomial $p(z) = z^4 + C_1 z^3 + C_2 z^2 + C_3 z + C_4$ with coefficients $C_1$, $C_2$, $C_3$, $C_4$ given by (4.18), in the vicinity of $\theta_1^* = 0.5393991$ and $\theta_2^* = 0.7221141$. This graph includes vertical line segments that indicate the location of the true theta roots $\theta_1^*$ and $\theta_2^*$. In this close-up, we used $N = 2{,}500$ elements $LzC(\theta_k)$ associated with points $\theta_k$ in a regular partition of interval $[0.46, 0.75)$.

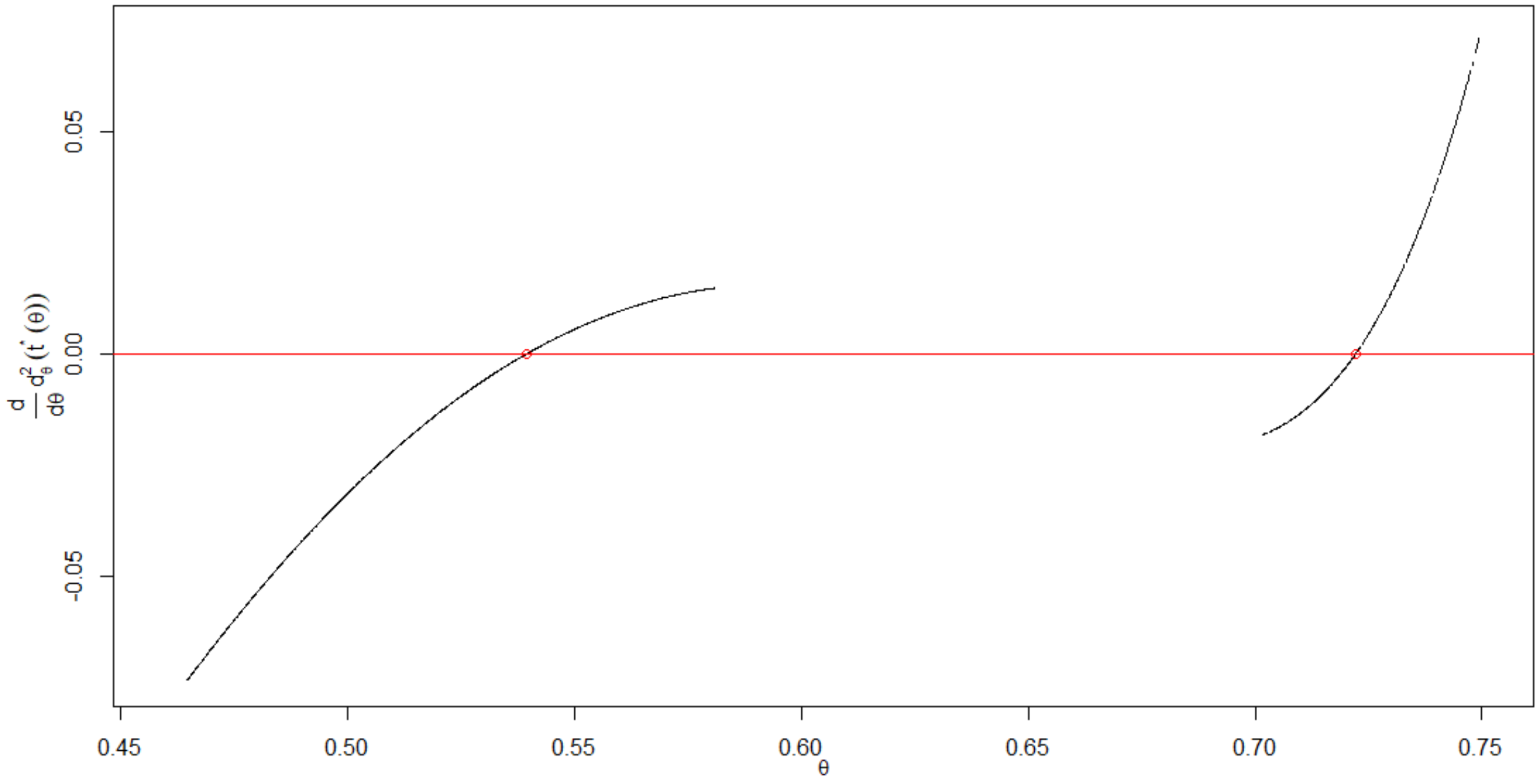

**Fig. 4.22.** Close-up of map $\frac{d}{d\theta}d_\theta^2\big(t^*(\theta)\big)$ associated to the polynomial $p(z) = z^4 + C_1z^3 + C_2z^2 + C_3z + C_4$ with coefficients $C_1$, $C_2$, $C_3$, $C_4$ given by (4.18), in the vicinity of $\theta_1^* = 0.5393991$ and $\theta_2^* = 0.7221141$. This graph includes small circles on the horizontal axis $y = 0$, that indicate the location of the true theta roots $\theta_1^*$ and $\theta_2^*$. For this regional map, we used $N = 2{,}500$ elements $LzC(\theta_k)$ associated with points $\theta_k$ in a regular partition of interval $[0.46, 0.75)$.

**Table 4.10**. Initial estimates for the roots of $p(z) = z^4 + C_1z^3 + C_2z^2 + C_3z + C_4$ with coefficients $C_1$, $C_2$, $C_3$, $C_4$ given by (4.18), obtained from the regional map $\frac{d}{d\theta}d_\theta^2\big(t^*(\theta)\big)$ in figure 4.22.

| $i$ | $\hat{R}_i$ | $\hat{\theta}_i^*$ | $\left\lvert\Delta\frac{d}{d\theta}d_{\hat{\theta}_i^*}^2(t^*)\right\rvert$ | $d_{\hat{\theta}_i^*}^2\left(t^*(\hat{\theta}_i^*)\right)$ |
|---|---|---|---|---|
| 1 | $0.2938060 - 0.2115483i$ | 0.5393991 | $6.756669 \times 10^{-5}$ | $4.030057 \times 10^{-14}$ |
| 2 | $0.2370027 - 0.0462186i$ | 0.7221140 | $1.704788 \times 10^{-4}$ | $1.405094 \times 10^{-12}$ |

Note: In figures 4.19 and 4.21, it can be seen that there are small gaps in the graphs of $t^*(\theta)$ and $d_\theta^2\big(t^*(\theta)\big)$; this is due to errors in the optimization process when finding the minimizing argument $t^*(\theta)$; specifically, what happens is that, for some angle $\theta_j$, sometimes there are two values $t_1(\theta_j) < 0$ and $t_2(\theta_j) > 0$ such that $d_{\theta_j}^2\left(t_1(\theta_j)\right)$ and $d_{\theta_j}^2\left(t_2(\theta_j)\right)$ are both local minima of function $d_{\theta_j}^2(t)$, very similar to each other, but where $d_{\theta_j}^2\left(t_2(\theta_j)\right)$ is actually the global minimum of $d_{\theta_j}^2(t)$; the simulated annealing algorithm in function `min_D2` of annex 2 section 4, however, erroneously selects $t_1(\theta_j)$ as the minimizing argument, so a `NA` value is assigned both to $d_{\theta_j}^2\left(\hat{t}^*(\theta_j)\right)$ and $\hat{t}^*(\theta_j)$, in accordance with restriction $\hat{t}^* > 0$, thus creating a gap in maps $t^*(\theta)$

and $d_\theta^2(t^*(\theta))$ at $\theta = \theta_j$. Figure 4.23 illustrates this situation for the function $d_{0.733064}^2(t)$ associated to the polynomial $p(z) = z^4 + C_1z^3 + C_2z^2 + C_3z + C_4$ with coefficients $C_1$, $C_2$, $C_3$, $C_4$ given by (4.18), where this type of optimization error occurs.

To conclude this example, let us see if a close-up of discrete map $e(\theta)$ at region $[0.46, 0.75)$ produces initial estimates of better quality than those obtained with its global counterpart. Figure 4.24 shows this close-up.

The numerical characteristics of the smooth crossings associated with map $e(\theta)$ (or LC map) in figure 4.24 are shown in table 4.11.

**Table 4.11**. Initial estimates for the roots of $p(z) = z^4 + C_1z^3 + C_2z^2 + C_3z + C_4$ with coefficients $C_1$, $C_2$, $C_3$, $C_4$ given by (4.18), obtained from the regional LC map in figure 4.24.

| $i$ | $\hat{R}_i$ | $\hat{\theta}_i^*$ | $\lvert\Delta e_i\rvert$ | $d_{\hat{\theta}_i^*}^2\left(t^*(\hat{\theta}_i^*)\right)$ |
|---|---|---|---|---|
| 1 | $0.2938062 - 0.2115485i$ | 0.5393988 | $6.408172 \times 10^{-5}$ | $1.072895 \times 10^{-13}$ |
| 2 | $0.2370039 - 0.0462200i$ | 0.7221121 | $2.216795 \times 10^{-4}$ | $4.192534 \times 10^{-12}$ |

When comparing corresponding rows in tables 4.5 and 4.11, we see that the values $d_{\hat{\theta}_i^*}^2\left(t^*(\hat{\theta}_i^*)\right)$ in table 4.11 are smaller than those in table 4.5; this tells us, at least in this case, that **the higher the resolution of an LC map, the more accurate its corresponding initial estimates to the roots**. Finally, if we compare the corresponding quantities $d_{\hat{\theta}_i^*}^2\left(t^*(\hat{\theta}_i^*)\right)$ between tables 4.10 and 4.11, we see that the initial estimates associated with the map $\frac{d}{d\theta}d_\theta^2(t^*(\theta))$ from figure 4.22 are slightly more accurate than the initial estimates associated with the LC map from figure 4.24.

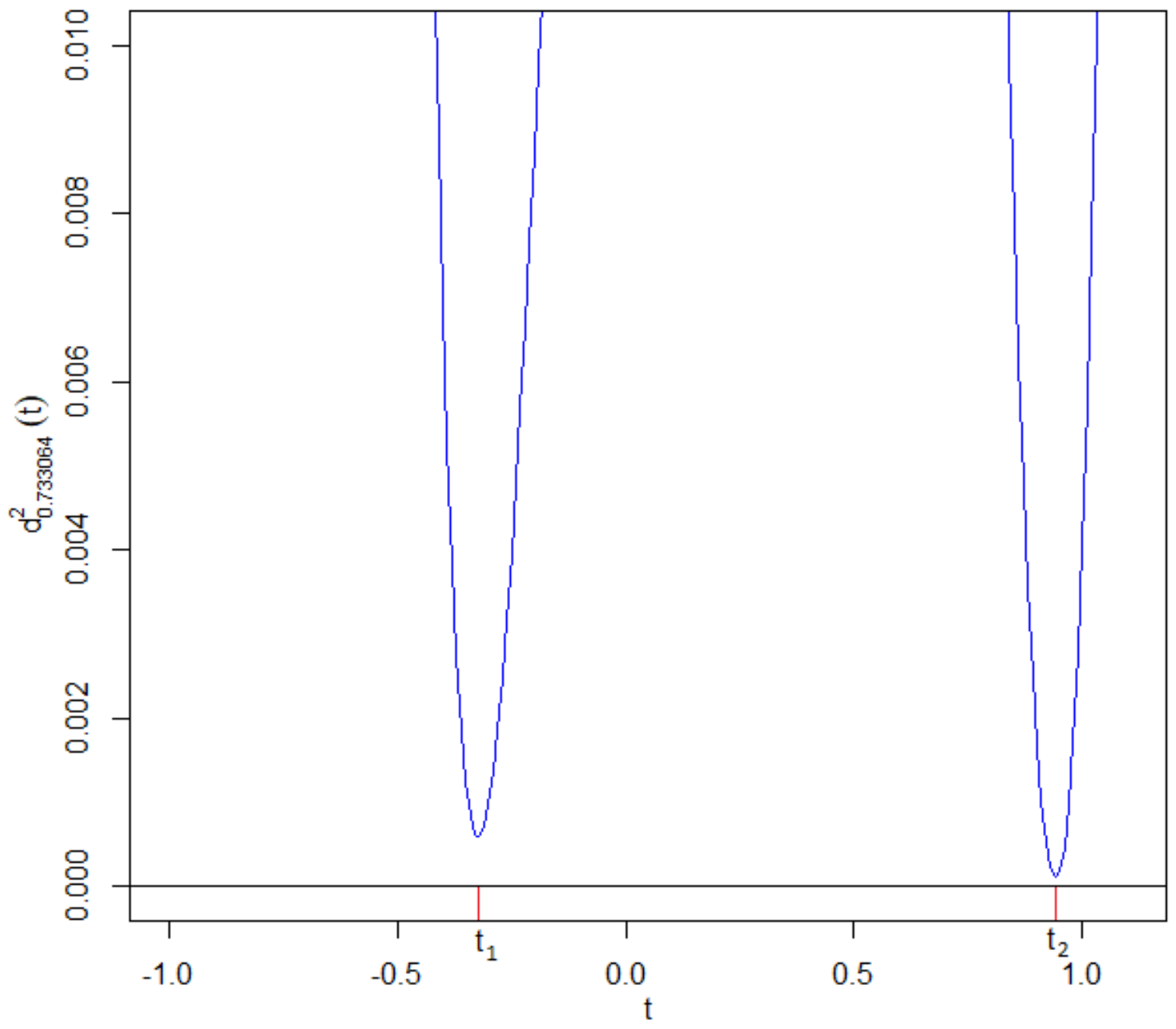

**Fig. 4.23.** Close-up of function $d^2_{0.733064}(t)$ associated to the polynomial $p(z) = z^4 + C_1z^3 + C_2z^2 + C_3z + C_4$ with coefficients $C_1$, $C_2$, $C_3$, $C_4$ given by (4.18), at a region containing its minimizing argument $t^* = t_2 > 0$; the simulated annealing algorithm used in `min_D2` misidentified the minimizing argument $t^*$ as $t_1 < 0$, so both the functional value of the global minimum $d^2_{0.733064}(t^*)$ and minimizing argument $t^*$ were assigned as `NA`, in accordance to constraint $t^* > 0$. In this case, the simulated annealing algorithm got stuck at a local minimum.

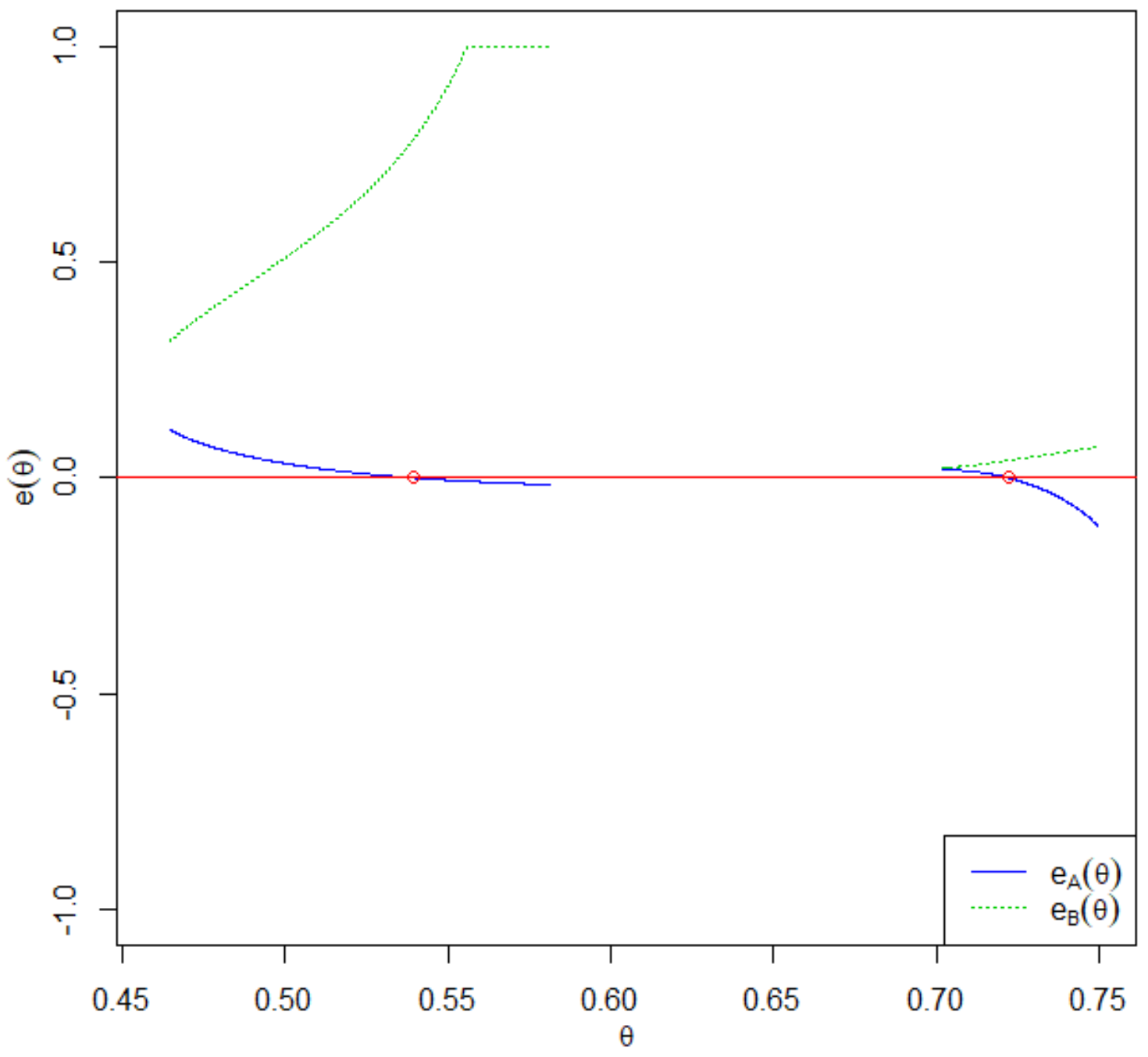

**Fig. 4.24.** Close-up of discrete map $e(\theta)$ (LC map) associated to polynomial $p(z) = z^4 + C_1z^3 + C_2z^2 + C_3z + C_4$ with coefficients $C_1$, $C_2$, $C_3$, $C_4$ given by (4.18). This close-up was generated from $N = 2{,}500$ elements $LzC(\theta_k)$ associated with points $\theta_k$ in a regular partition of interval $[0.46, 0.75)$. Small circles on horizontal axis $y = 0$ indicate the location of true theta roots $\theta_1^*$ and $\theta_2^*$.

## Numerical Example 4.3

This time we will use the LC method to find approximations to the roots of equation

$$z^4 + (1 + i)z^3 + (2 + 2i)z^2 + (3 + 3i)z + (4 + 4i) = 0 \qquad (4.22)$$

The reference roots in this case will be generated in a similar way as we did in example 3.2 of chapter 3, by using the R function `polyroot`. Additionally, we will initialize the random number generator seed by using the instruction `set.seed(1987)` at the beginning of the script listed in annex 3 section 3, in order to reproduce the numerical outputs from the stochastic optimization processes used in the search of minimizing arguments $t^*$ for the functions $d_\theta^2$ associated with the structures $LzC(\theta)$ generated in this numerical example. For further details, see subsection "specific conditions for reproducing results from the examples in chapter 4" within annex 3 section 3.

In this case, the reference roots $R_i$ and $\theta_i^*$ are:

$$\begin{aligned} R_1 &= \phantom{-}0.217902 + 1.406896i & \theta_1^* &= \phantom{-}1.2107319 \\ R_2 &= -1.231898 + 0.586985i & \theta_2^* &= \phantom{-}2.1633969 \\ R_3 &= -0.846674 - 1.167477i & \theta_3^* &= -2.0498271 \\ R_4 &= \phantom{-}0.860670 - 1.826404i & \theta_4^* &= -0.7726466 \end{aligned} \qquad (4.23)$$

First, let us take a look at the discrete global proximity maps (figures 4.25, 4.26, 4.27, 4.28, 4.29, 4.30), all of them generated from $N = 5{,}000$ elements $LzC(\theta_k)$ associated with points $\theta_k$ in a regular partition of interval $[-\pi, \pi)$.

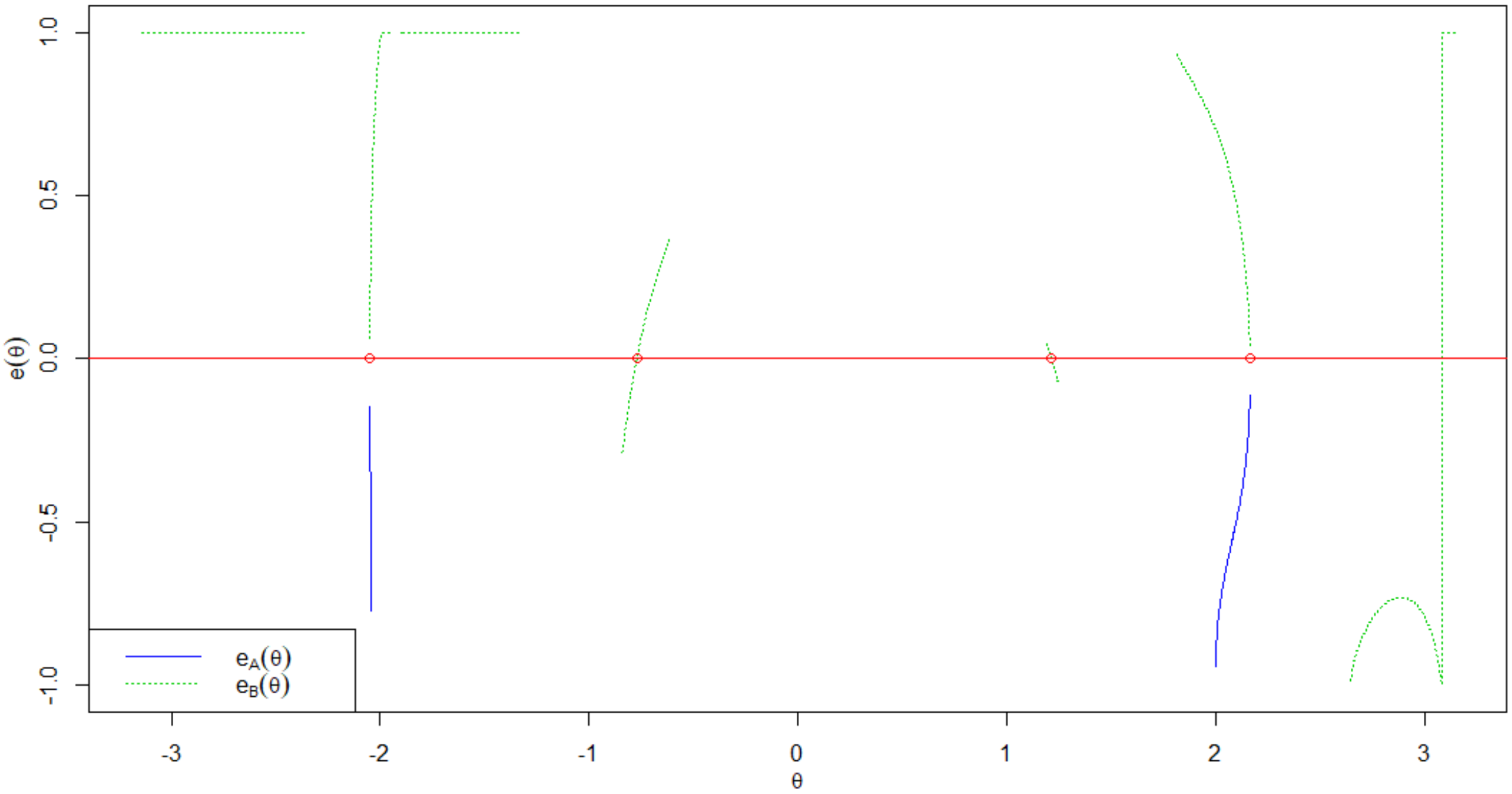

**Fig. 4.25.** Global map $e(\theta)$ associated with $p(z) = z^4 + (1 + i)z^3 + (2 + 2i)z^2 + (3 + 3i)z + (4 + 4i)$. This LC map was generated from $N = 5{,}000$ elements $LzC(\theta_k)$ associated with points $\theta_k$ in a regular partition of interval $[-\pi, \pi)$. Small circles on horizontal axis $y = 0$ mark the location of true theta roots $\theta_i^*$.

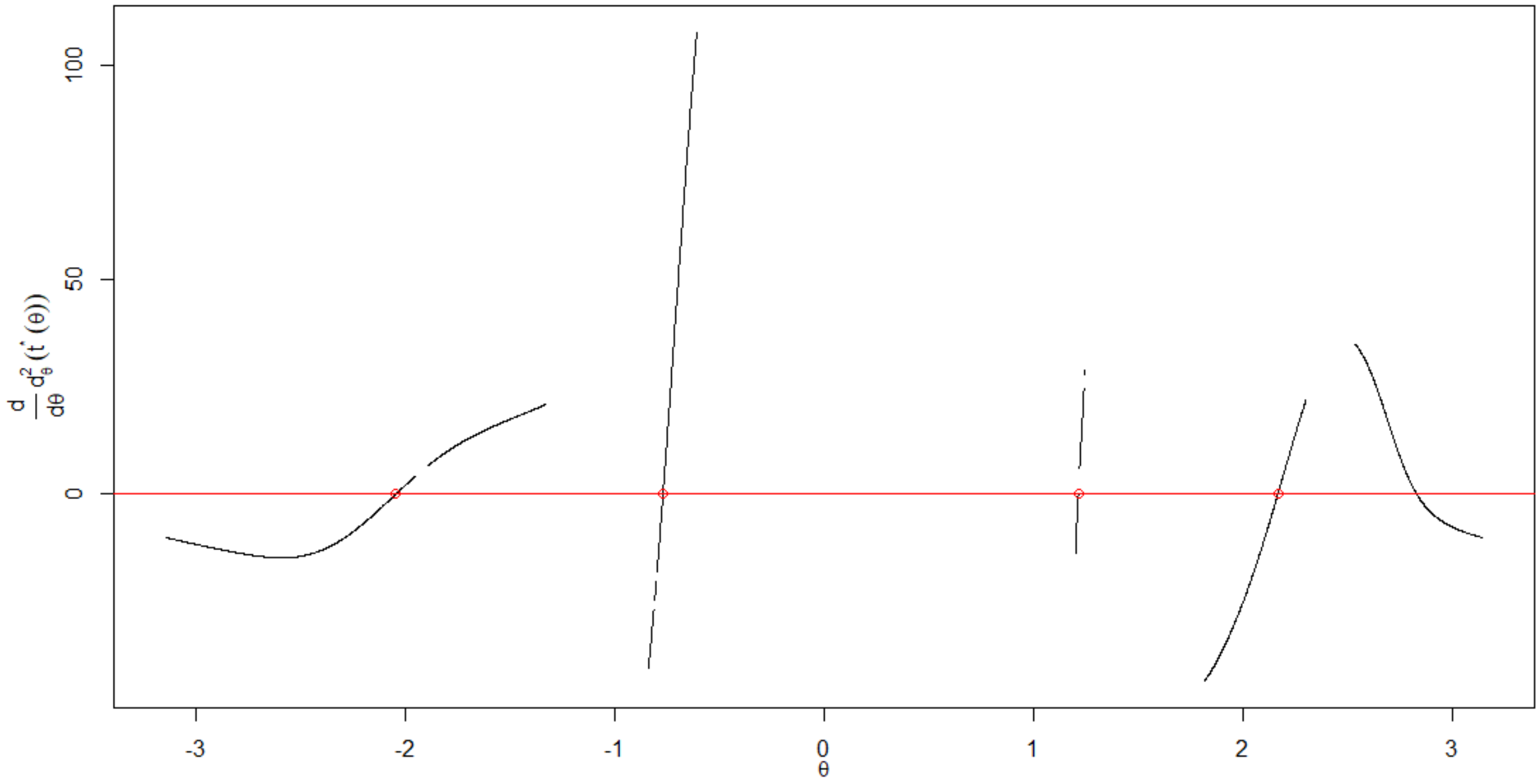

**Fig. 4.26.** Global map $\frac{d}{d\theta}d_\theta^2\big(t^*(\theta)\big)$ associated with $p(z) = z^4 + (1+i)z^3 + (2+2i)z^2 + (3+3i)z + (4+4i)$. This map was generated from $N = 5{,}000$ elements $LzC(\theta_k)$ associated with points $\theta_k$ in a regular partition of interval $[-\pi, \pi)$. Small circles on horizontal axis $y = 0$ mark the location of true theta roots $\theta_i^*$.

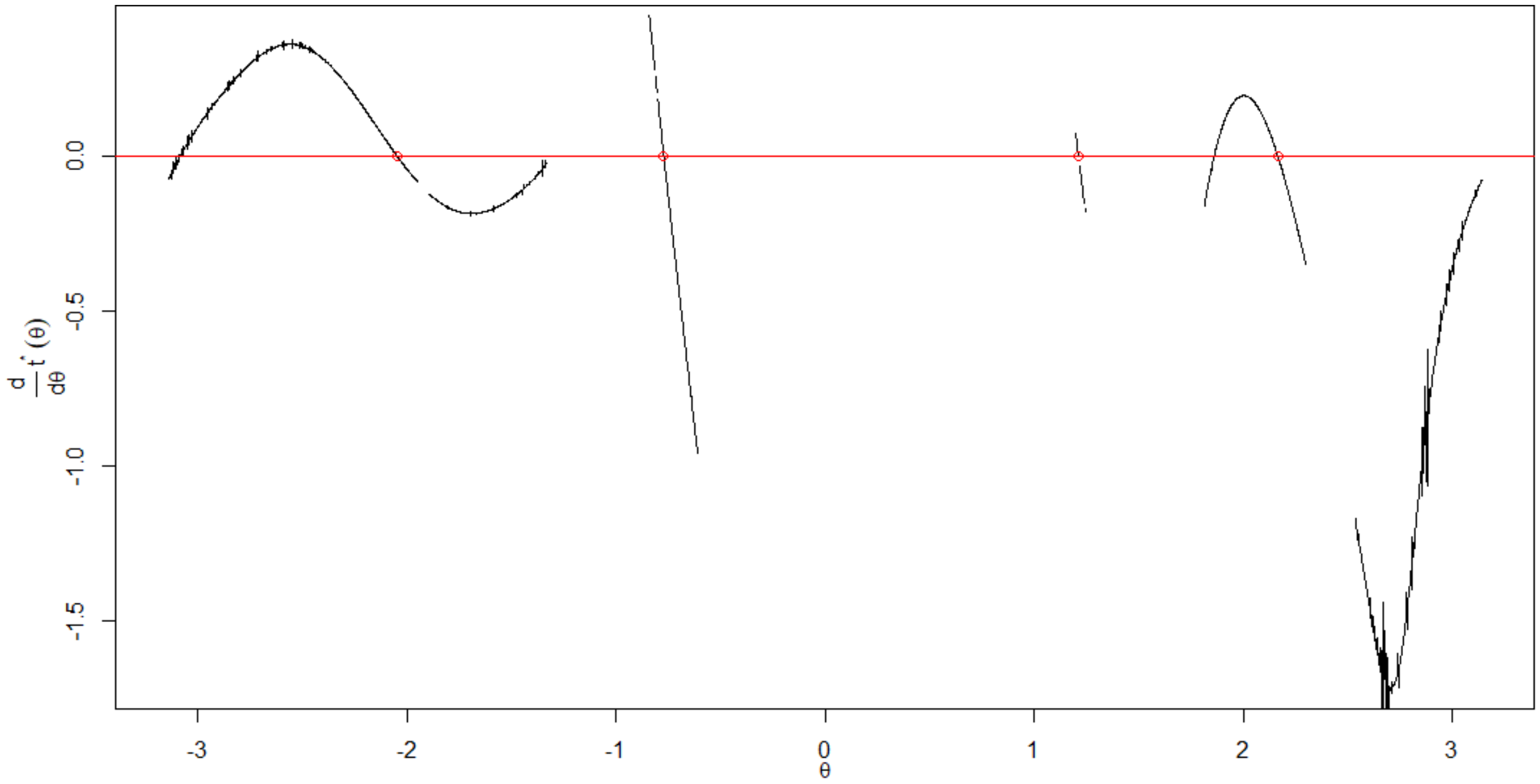

**Fig. 4.27.** Global map $\frac{d}{d\theta}t^*(\theta)$ associated with $p(z) = z^4 + (1+i)z^3 + (2+2i)z^2 + (3+3i)z + (4+4i)$. This map was generated from $N = 5{,}000$ elements $LzC(\theta_k)$ associated with points $\theta_k$ in a regular partition of interval $[-\pi, \pi)$. Small circles on horizontal axis $y = 0$ mark the location of true theta roots $\theta_i^*$.

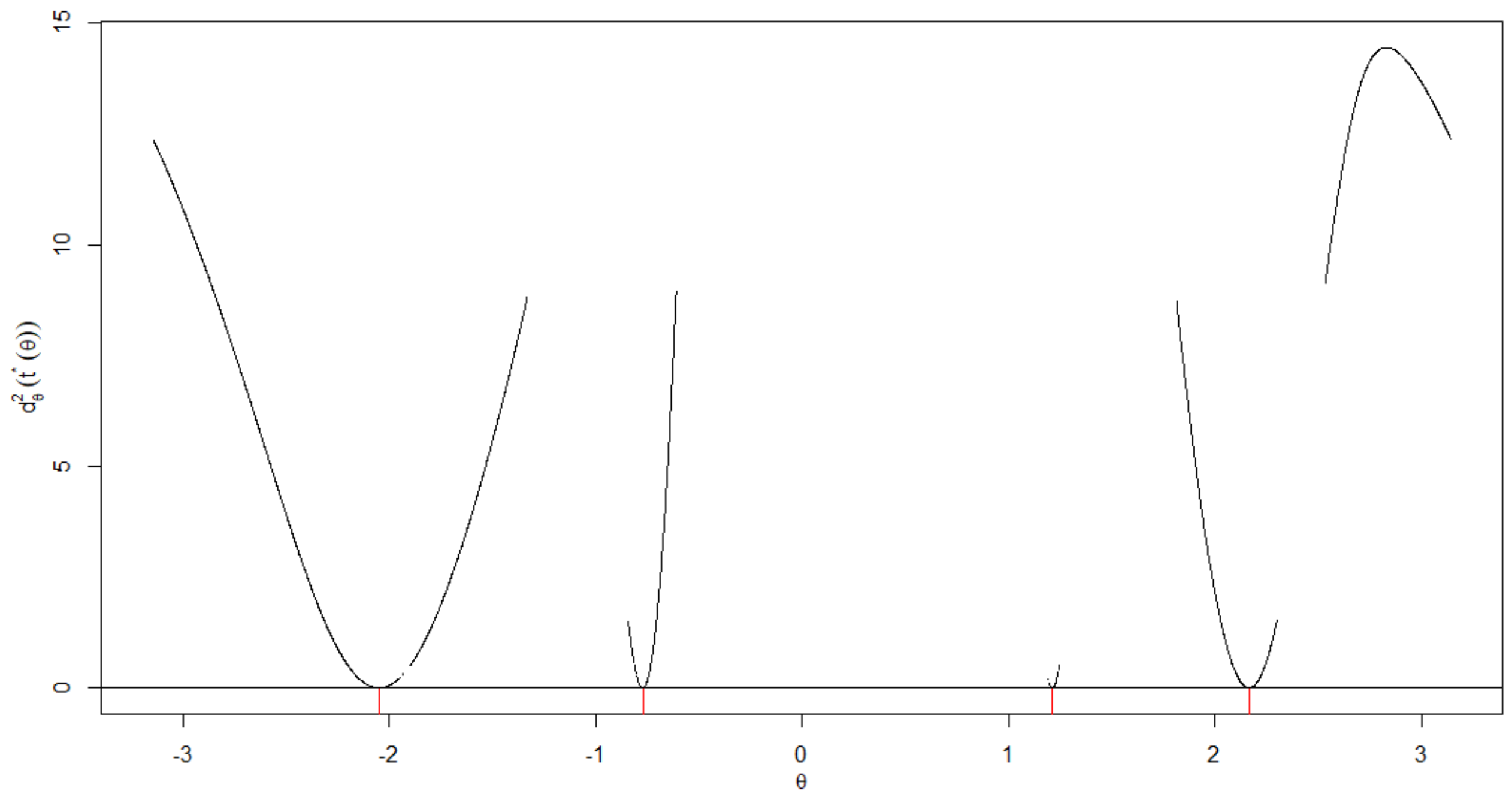

**Fig. 4.28.** Global minima $d_\theta^2(t^*(\theta))$ vs. $\theta$ for $p(z) = z^4 + (1+i)z^3 + (2+2i)z^2 + (3+3i)z + (4+4i)$. This map was generated from $N = 5{,}000$ elements $LzC(\theta_k)$ associated with points $\theta_k$ in a regular partition of interval $[-\pi, \pi)$. This graph includes vertical line segments that indicate the location of true theta roots $\theta_i^*$.

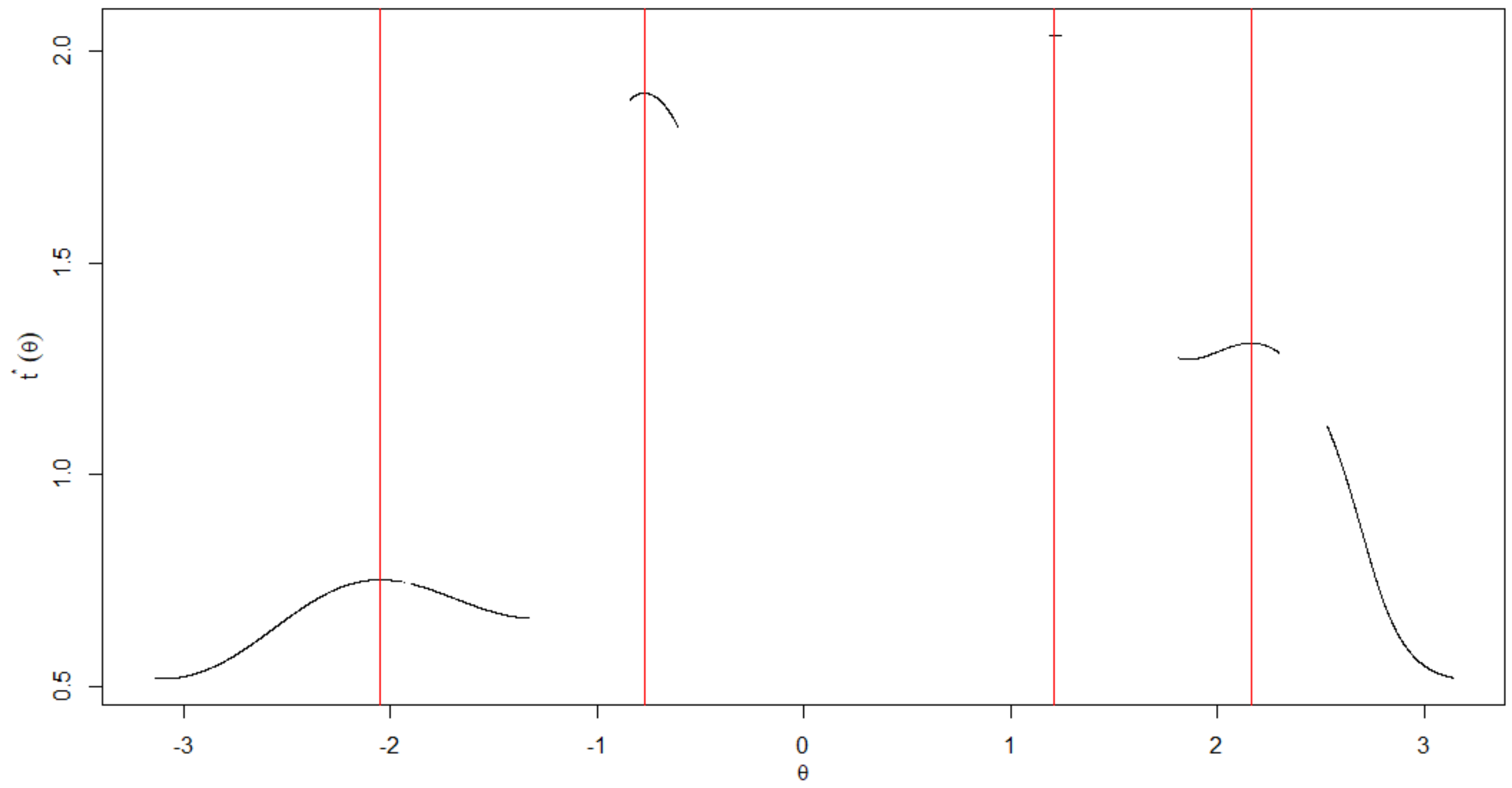

**Fig. 4.29.** Minimizing arguments $t^*(\theta)$ vs. $\theta$ for $p(z) = z^4 + (1+i)z^3 + (2+2i)z^2 + (3+3i)z + (4+4i)$. This map was generated from $N = 5{,}000$ elements $LzC(\theta_k)$ associated with points $\theta_k$ in a regular partition of interval $[-\pi, \pi)$. This graph includes vertical line segments that indicate the location of true theta roots $\theta_i^*$.

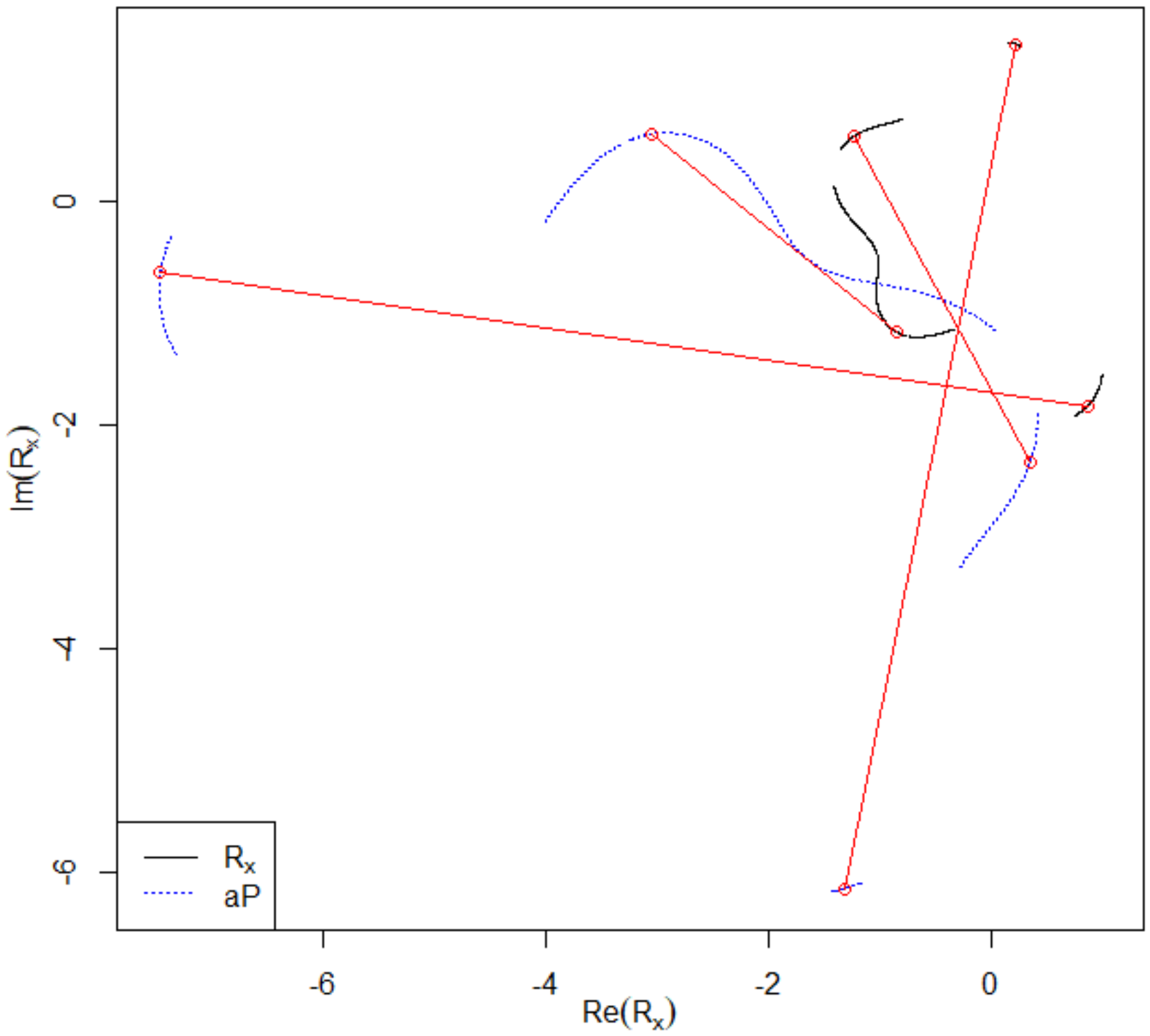

**Fig. 4.30.** Trajectories $R_x(\theta)$ and $aP(\theta)$ associated with $p(z) = z^4 + (1+i)z^3 + (2+2i)z^2 + (3+3i)z + (4+4i)$. These trajectories were generated from $N = 5{,}000$ elements $LzC(\theta_k)$ associated with points $\theta_k$ in a regular partition of interval $[-\pi, \pi)$. True roots $R_i$ are represented in this graph by means of small circles, which are joined, by means of line segments, with their corresponding anchor points.

Now let us see, in tables 4.12, 4.13 and 4.14, the numerical characteristics of smooth crossings associated with the respective maps from figures 4.25, 4.26, and 4.27. As usual, the rows in these tables are ranked according to values $d^2_{\hat{\theta}_i^*}\left(t^*(\hat{\theta}_i^*)\right)$.

**Table 4.12**. Initial estimates of the roots of $p(z) = z^4 + (1+i)z^3 + (2+2i)z^2 + (3+3i)z + (4+4i)$, obtained from the map $e(\theta)$ in figure 4.25.

| $i$ | $\hat{R}_i$ | $\hat{\theta}_i^*$ | $\lvert\Delta e_i\rvert$ | $d^2_{\hat{\theta}_i^*}\left(t^*(\hat{\theta}_i^*)\right)$ |
|---|---|---|---|---|
| 1 | $0.217902 + 1.406896i$ | 1.210732 | 0.002695238 | $1.558697 \times 10^{-11}$ |
| 2 | $0.860671 - 1.826403i$ | $-0.772646$ | 0.004180607 | $1.450015 \times 10^{-10}$ |

**Table 4.13**. Initial estimates of the roots of $p(z) = z^4 + (1+i)z^3 + (2+2i)z^2 + (3+3i)z + (4+4i)$, obtained from the map $\frac{d}{d\theta} d^2_\theta(t^*(\theta))$ in figure 4.26.

| $i$ | $\hat{R}_i$ | $\hat{\theta}_i^*$ | $\left\lvert\Delta \frac{d}{d\theta} d^2_{\hat{\theta}_i^*}(t^*)\right\rvert$ | $d^2_{\hat{\theta}_i^*}\left(t^*(\hat{\theta}_i^*)\right)$ |
|---|---|---|---|---|
| 1 | $-1.231898 + 0.586985i$ | 2.1633968 | 0.21172420 | $7.610396 \times 10^{-13}$ |
| 2 | $-0.846674 - 1.167477i$ | $-2.0498271$ | 0.05732612 | $1.095553 \times 10^{-12}$ |
| 3 | $0.217901 + 1.406897i$ | 1.2107323 | 1.20777904 | $8.134912 \times 10^{-11}$ |
| 4 | $0.860669 - 1.826405i$ | $-0.7726471$ | 0.82638505 | $1.270799 \times 10^{-10}$ |
| 5 | $-1.137827 - 0.292914i$ | 2.8276553 | 0.10207474 | $1.445777 \times 10^{1}$ |

**Table 4.14**. Initial estimates of the roots of $p(z) = z^4 + (1+i)z^3 + (2+2i)z^2 + (3+3i)z + (4+4i)$, obtained from the map $\frac{d}{d\theta}t^*(\theta)$ in figure 4.27.

| $i$ | $\hat{R}_i$ | $\hat{\theta}_i^*$ | $\left\|\Delta \frac{d}{d\theta} t^*(\hat{\theta}_i^*)\right\|$ | $d^2_{\hat{\theta}_i^*}\left(t^*(\hat{\theta}_i^*)\right)$ |
|---|---|---|---|---|
| 1 | $0.8606750 - 1.826399i$ | $-0.7726428$ | 0.007801266 | $4.742278 \times 10^{-9}$ |
| 2 | $-0.8466385 - 1.167495i$ | $-2.0497740$ | 0.001343209 | $6.453961 \times 10^{-8}$ |
| 3 | $-1.2319771 + 0.586931i$ | 2.1634697 | 0.002914981 | $4.472522 \times 10^{-7}$ |
| 4 | $0.2179867 + 1.406864i$ | 1.2106872 | 0.006993710 | $9.580237 \times 10^{-7}$ |
| 5 | $-0.8623155 + 0.718781i$ | 1.8597531 | 0.003541973 | $6.736168 \times 10^{0}$ |
| 6 | $-1.0181834 - 0.535593i$ | $-3.0730125$ | 0.017711527 | $1.165101 \times 10^{1}$ |
| 7 | $-1.0181863 - 0.535517i$ | $-3.0731592$ | 0.031076200 | $1.165262 \times 10^{1}$ |
| 8 | $-1.0183058 - 0.532976i$ | $-3.0780553$ | 0.011078846 | $1.170599 \times 10^{1}$ |
| 9 | $-1.0182990 - 0.532910i$ | $-3.0781810$ | 0.017608043 | $1.170735 \times 10^{1}$ |
| 10 | $-1.0184919 - 0.529151i$ | $-3.0854296$ | 0.001584458 | $1.178586 \times 10^{1}$ |
| 11 | $-1.0186510 - 0.526764i$ | $-3.0900350$ | 0.010437256 | $1.183543 \times 10^{1}$ |
| 12 | $-1.0186549 - 0.526278i$ | $-3.0909697$ | 0.025580442 | $1.184546 \times 10^{1}$ |

Table 4.15 shows a summary of the approximation errors for valid initial estimates in tables 4.12, 4.13 and 4.14, with respect to corresponding reference values (4.23).

**Table 4.15**. Arithmetic means and standard deviations for error measures associated with valid initial approximations (first four rows in each of tables 4.13 and 4.14, and the two rows of table 4.12) to the roots of $p(z) = z^4 + (1+i)z^3 + (2+2i)z^2 + (3+3i)z + (4+4i)$, obtained by means of discrete proximity maps $e(\theta)$, $\frac{d}{d\theta}d^2_\theta\left(t^*(\theta)\right)$, and $\frac{d}{d\theta}t^*(\theta)$ from figures 4.25, 4.26 and 4.27.

| **Map** | **Error measure** | | | | | |
|---|---|---|---|---|---|---|
| | $d^2_{\hat{\theta}_i^*}\left(t^*(\hat{\theta}_i^*)\right)$ | | $\left\|\frac{R-\hat{R}}{R}\right\|$ | | $\left\|\frac{\theta^*-\hat{\theta}^*}{\theta^*}\right\|$ | |
| | **Mean** | **Standard deviation** | **Mean** | **Standard deviation** | **Mean** | **Standard deviation** |
| $e(\theta)$ | $8.029422 \times 10^{-11}$ | $9.150988 \times 10^{-11}$ | $4.414834 \times 10^{-7}$ | $2.593329 \times 10^{-7}$ | $4.606924 \times 10^{-7}$ | $5.174777 \times 10^{-7}$ |
| $\frac{d}{d\theta}d^2_\theta(t^*(\theta))$ | $5.257140 \times 10^{-11}$ | $6.248666 \times 10^{-11}$ | $3.450509 \times 10^{-7}$ | $2.799033 \times 10^{-7}$ | $2.715014 \times 10^{-7}$ | $3.173839 \times 10^{-7}$ |
| $\frac{d}{d\theta}t^*(\theta)$ | $3.686395 \times 10^{-7}$ | $4.391096 \times 10^{-7}$ | $4.132022 \times 10^{-5}$ | $3.132802 \times 10^{-5}$ | $2.536428 \times 10^{-5}$ | $1.440080 \times 10^{-5}$ |

As we can see from table 4.15, the estimates in the first four rows of table 4.13, corresponding to map $\frac{d}{d\theta}d^2_\theta\left(t^*(\theta)\right)$ in figure 4.26, are the most consistent, closely followed by the two estimates in table 4.12 corresponding to map $e(\theta)$ in figure 4.25; in third place, we have the estimates in the first four rows of table 4.14, corresponding to map $\frac{d}{d\theta}t^*(\theta)$ in figure 4.27. In this case, we can say that all valid initial estimates obtained for the roots of equation (4.22) by means of discrete maps $e(\theta)$, $\frac{d}{d\theta}d^2_\theta\left(t^*(\theta)\right)$ and even $\frac{d}{d\theta}t^*(\theta)$, are of very good quality. With map $e(\theta)$, however, it was not possible to estimate the 4 roots of equation (4.22), as it was with maps $\frac{d}{d\theta}d_\theta^2\big(t^*(\theta)\big)$ and $\frac{d}{d\theta}t^*(\theta)$. This demonstrates the usefulness of any of these three types of proximity maps as complementary root-finding tools, when used jointly.

**Remarks on LC map:** As we can see from map $e(\theta)$ in figure 4.25, it only contains two smooth crossings with horizontal axis $y = 0$; the numerical characteristics of these two crossings, shown in table 4.12, approach reference values $R_1$, $\theta_1^*$ and $R_4$, $\theta_4^*$ in (4.23). In this case, the existence of gaps between $\hat{e}_A(\theta)$ and $\hat{e}_B(\theta)$ prevents obtaining the other two expected estimates; this kind of difficulty was mentioned in examples 3.2 and 3.3 of chapter 3: it is difficult to detect theta roots $\theta_i^*$ that lie between two consecutive sample angles $\theta_k = -\pi + 2\pi k/N$ and $\theta_{k+1} = \theta_k + 2\pi/N$ such that, for one of them, say for $\theta_k$, $\hat{e}_A(\theta_k)$ and $\hat{e}_B(\theta_k)$ are of opposite signs, while for the other sample angle $\theta_{k+1}$, $\hat{e}_A(\theta_{k+1})$ and $\hat{e}_B(\theta_{k+1})$ are not defined; as we saw in chapter 3, this happens because $\theta_i^*$ is close to an angle $\theta_T$ such that terminal semi-line $tL(\theta_T)$ and z-circumference $zC(\theta_T)$ intersect tangentially at a single point; in such circumstances, of course, it is to be expected that both $\theta_i^*$ and $\theta_T$ be located between $\theta_k$ (where $tL(\theta_k)$ and $zC(\theta_k)$ intersect at two points very close to each other) and $\theta_{k+1}$ (where $tL(\theta_{k+1})$ and $zC(\theta_{k+1})$ do not intersect). In this example, the gaps between $\hat{e}_A(\theta)$ and $\hat{e}_B(\theta)$ occur near $\theta_2^* = 2.1633969$ and $\theta_3^* = -2.0498271$; for the case of $\theta_2^*$, its associated terminal semi-line $tL(\theta_2^*)$ intersects z-circumference $zC(\theta_2^*)$ at two points very close to each other; this suggests that a small variation $\delta\theta$ in $\theta_2^*$ will cause $tL(\theta_2^* + \delta\theta)$ and $zC(\theta_2^* + \delta\theta)$ to intersect tangentially at a single point. See figure 4.31, where we can confirm this situation described for theta root $\theta_2^*$, which was not detected by map $e(\theta)$ in figure 4.25, due to the angular resolution used.

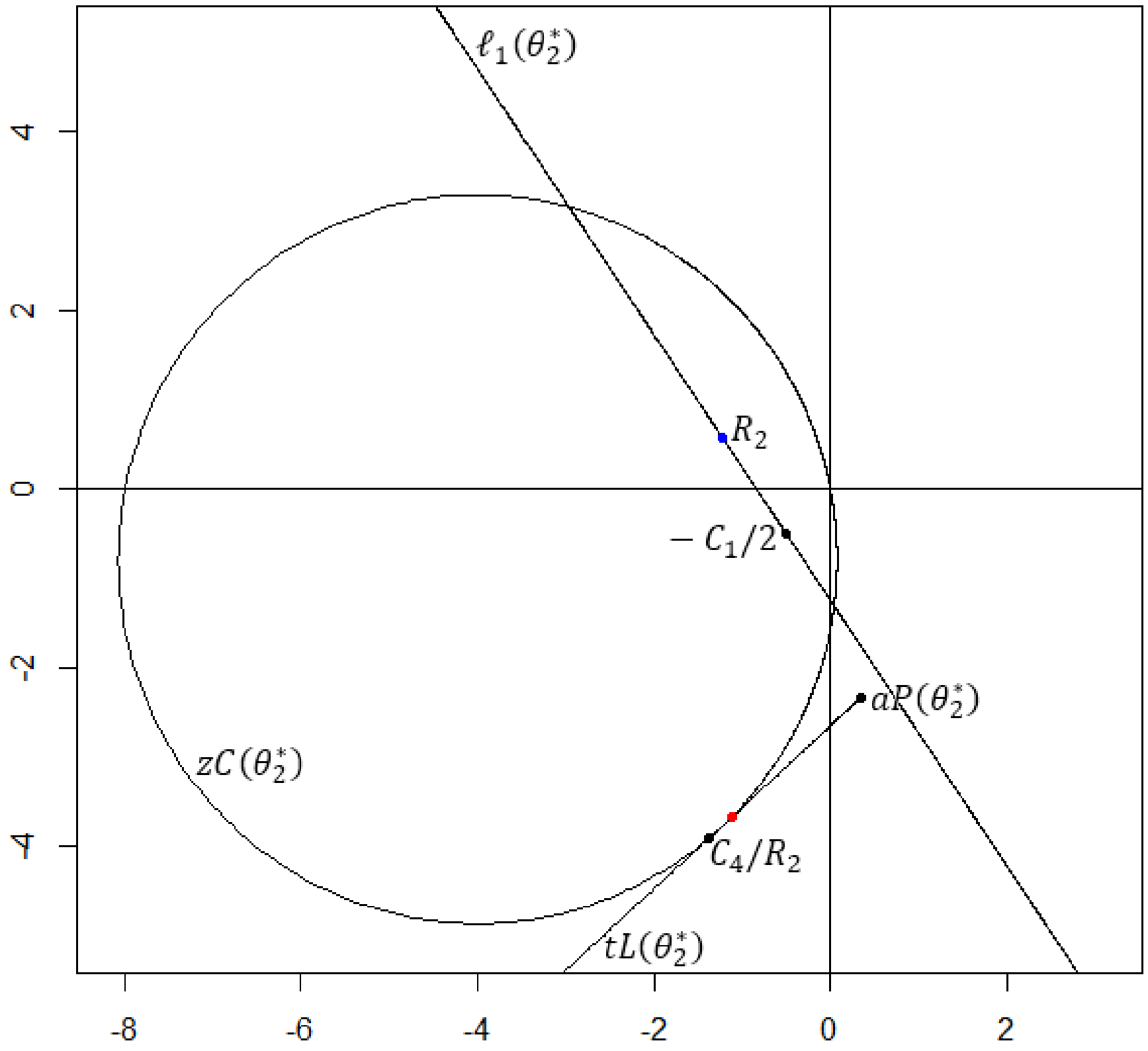

**Fig. 4.31.** Some elements of the structure $LzC(\theta_2^*)$ associated with equation (4.22). We see that terminal semi-line $tL(\theta_2^*)$ almost touches z-circumference $zC(\theta_2^*)$ tangentially at a single point, although it is actually touching it at two points (shown in this graph) very close to each other; this suggests that a small change in $\theta_2^*$ will cause $tL$ and $zC$ to touch tangentially at a single point. This graph also shows line $\ell_1(\theta_2^*)$, which contains fixed point $-C_1/2$, as well as root $R_2$.

**Remarks on map $\frac{d}{d\theta}d_\theta^2(t^*(\theta))$:** From the discrete map $\frac{d}{d\theta}d_\theta^2(t^*(\theta))$ in figure 4.26, we see that there are 5 smooth crossings $\hat{\theta}_i^*$ with the horizontal axis $y = 0$; table 4.13 shows the numerical characteristics of these crossings; from here, we see that the estimates in the first 4 rows correspond respectively to reference roots $R_2$, $R_3$, $R_1$, $R_4$ in (4.23). Another detail we can see from table 4.13 is that quantities $\left|\Delta\frac{d}{d\theta}d_{\hat{\theta}_i^*}^2(t^*)\right|$, analogous to quantities $|\Delta e_i|$ in tables with estimates from LC maps, and to quantities $\left|\Delta\frac{d}{d\theta}t^*(\hat{\theta}_i^*)\right|$ in tables with estimates from maps $\frac{d}{d\theta}t^*(\theta)$, are some orders of magnitude greater than their counterparts in tables 4.2 and 4.6, even though the resolution of the global map $\frac{d}{d\theta}d_\theta^2(t^*(\theta))$ associated with table 4.13 is two times the resolution of global maps $\frac{d}{d\theta}d_\theta^2(t^*(\theta))$ associated with tables 4.2 and 4.6; this discrepancy is due to the fact that the roots of equation (4.22) are not limited to the region $\{x + iy : -1 \le x \le 1, -1 \le y \le 1\}$, as was the case for examples 4.1 and 4.2. This suggests that the map $\frac{d}{d\theta}d_\theta^2(t^*(\theta))$ is sensitive to the norms of the coefficients of the univariate polynomial from which we want to approximate its roots (as we can see by comparing the vertical ranges of functions $\frac{d}{d\theta}d_\theta^2(t^*(\theta))$ from figures 4.11, 4.17 and 4.26), so it would be prudent to be careful when selecting the argument `Tol` which defines the upper limit for the variation of two consecutive discrete functional values that determine a smooth crossing; this argument `Tol` serves as an input to the function `approxDMin`, which implements the construction of maps $\frac{d}{d\theta} d_\theta^2\big(t^*(\theta)\big)$ and $\frac{d}{d\theta} t^*(\theta)$ (see annex 2 section 7).

**Remarks on map $\frac{d}{d\theta} t^*(\theta)$:** From map $\frac{d}{d\theta} t^*(\theta)$ in figure 4.27, apparently 6 smooth crossings with horizontal axis $y = 0$ are observed, but there are actually more than 6, given the presence of oscillatory noise, which appears here due to the high sampling rate (5,000 elements $LzC(\theta_k)$ associated with points $\theta_k$ in a regular partition of global interval $[-\pi, \pi)$) combined with the amplification effect of optimization errors induced by the finite difference operators used to approximate differential operator $\frac{d}{d\theta}$. The numerical characteristics of these crossings are shown in table 4.14; we see here that the first 4 rows correspond to approximations of respective reference roots $R_4$, $R_3$, $R_2$, $R_1$ in (4.23). Row 5 in table 4.14 corresponds to a smooth crossing that does not approach any root of equation (4.22), and rows 6 to 12 correspond to "replicas", due to the oscillatory noise present on the map in figure 4.27, of another smooth crossing that does not approach any root of equation (4.22).

**Remarks on the joint behavior of proximity maps generated in this example:** By analyzing jointly the graphs from figures 4.26, 4.27, 4.28, 4.29 and 4.30, we see once again that these provide empirical evidence in favor of the hypotheses posed in example 4.1 on the relationships of functions $d_\theta^2\big(t^*(\theta)\big)$, $t^*(\theta)$, $\frac{d}{d\theta} d_\theta^2\big(t^*(\theta)\big)$ and $\frac{d}{d\theta} t^*(\theta)$ with the theta roots of their associated univariate polynomial; we also see that, in this set of graphs, we can identify 5 continuous sections with support of total length $\pi$. Small gaps are observed in some of these continuous sections, due to errors when assigning the location of global minima within the stochastic optimization process, in the presence of local minima with similar functional values, as described in example 4.2 (see figure 4.23). In this case, the four theta roots $\theta_i^*$ of equation (4.22) correspond to local maxima of mapping $t^*(\theta)$ (see figures 4.27 and 4.29).

**Close-up to theta root $\boldsymbol{\theta_2^*}$ on LC map:** We will see if by means of the LC map in this numeric example it is possible to detect theta root $\theta_2^*$ by making a high-resolution close-up in the vicinity of the gap observed between $\hat{e}_A(\theta)$ and $\hat{e}_B(\theta)$ on the right side of the LC map from figure 4.25; for this, we will consider $N = 2{,}500$ elements $LzC(\theta_k)$ associated with points $\theta_k$ in a regular partition of interval $[2.0, 2.3)$. The graph of this close-up can be seen in figure 4.32. From here, we see that there is a smooth crossing of function $\hat{e}_B(\theta)$ with horizontal axis $y = 0$; the numerical characteristics of this crossing are shown in table 4.16.

**Table 4.16**. Estimate of theta root $\theta_2^* = 2.1633969$ associated with $p(z) = z^4 + (1+i)z^3 + (2+2i)z^2 + (3+3i)z + (4+4i)$, obtained from regional map $e(\theta)$ in figure 4.32.

| $i$ | $\hat{R}_i$ | $\hat{\theta}_i^*$ | $\lvert\Delta e_i\rvert$ | $d^2_{\hat{\theta}_i^*}\left(t^*(\hat{\theta}_i^*)\right)$ |
|---|---|---|---|---|
| 1 | $-1.231896 + 0.586986i$ | 2.163395 | 0.005618004 | $2.422679 \times 10^{-10}$ |

If we compare the value $d^2_{\hat{\theta}_i^*}\left(t^*(\hat{\theta}_i^*)\right)$ in table 4.16 with any of the values $d^2_{\hat{\theta}_i^*}\left(t^*(\hat{\theta}_i^*)\right)$ in table 4.12 associated with the global map $e(\theta)$ in figure 4.25, we see that the quality of the estimate associated with the regional map $e(\theta)$ in figure 4.32 is no better than the quality of either of the two estimates associated with the global map $e(\theta)$ in figure 4.25; however, this close-up of map $e(\theta)$ to region $[2.0, 2.3)$ was able to detect, with good estimation quality, a theta root that the global map $e(\theta)$ in figure 4.25 had not detected. It is left as an exercise for the reader to verify that the close-up of map $\frac{d}{d\theta} d^2_\theta\big(t^*(\theta)\big)$ to region $[2.0, 2.3)$ by using $N = 2{,}500$ elements $LzC(\theta_k)$ associated with points $\theta_k$ in a regular partition of interval $[2.0, 2.3)$, produces an approximation to $\theta_2^*$ of slightly better quality than that of the corresponding estimate associated with the global map $\frac{d}{d\theta} d^2_\theta\big(t^*(\theta)\big)$ from figure 4.26; it is also possible to verify that the close-up of map $\frac{d}{d\theta} t^*(\theta)$ to region $[2.0, 2.3)$ by using $N = 2{,}500$ elements $LzC(\theta_k)$ associated with points $\theta_k$ in a regular partition of interval $[2.0, 2.3)$, produces an approximation to $\theta_2^*$ of better quality than that of the corresponding estimate obtained with the global map $\frac{d}{d\theta} t^*(\theta)$ from figure 4.27, although such regional approximation is also accompanied by several “replicas” caused by oscillatory noise.

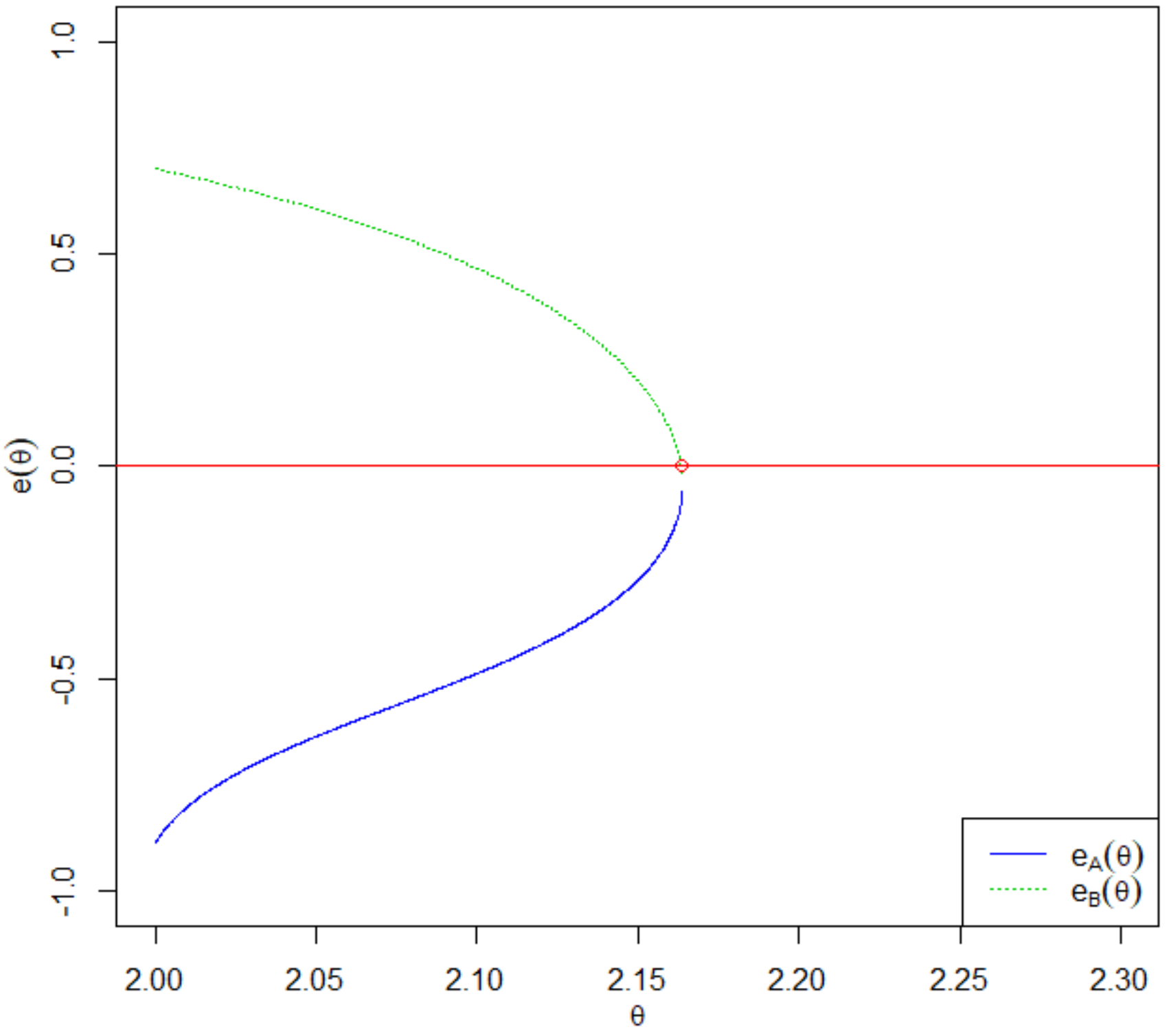

**Fig. 4.32.** Close-up of map $e(\theta)$ associated with $p(z) = z^4 + (1+i)z^3 + (2+2i)z^2 + (3+3i)z + (4+4i)$; this close-up was generated from $N = 2{,}500$ elements $LzC(\theta_k)$ associated with points $\theta_k$ in a regular partition of interval $[2.0, 2.3)$. The small circle on horizontal axis $y = 0$ indicates the location of true theta root $\theta_2^*$.

## Numerical Example 4.4

In this last example we will randomly generate four roots in the same way as we did in the first two examples of this chapter, by using the script of annex 3 section 3, but this time initializing the random number generator seed at the beginning of this script with the instruction `set.seed(2038)`; the global proximity maps will be generated by using $N = 5{,}000$ elements $LzC(\theta_k)$ associated with points $\theta_k$ in a regular partition of interval $[-\pi, \pi)$.

The reference roots $R_i$, $\theta_i^*$ generated in this example are:

$$
\begin{aligned}
R_1 &= 0.1408556 + 0.2165732i & \theta_1^* &= -2.982106 \\
R_2 &= 0.7478233 + 0.4507773i & \theta_2^* &= 2.996166 \\
R_3 &= 0.7685930 + 0.9011391i & \theta_3^* &= 2.258735 \\
R_4 &= 0.7271858 - 0.7971104i & \theta_4^* &= -1.945402
\end{aligned}
\tag{4.24}
$$

The coefficients constructed from Vieta's relations (4.2-a), (4.2-b), (4.2-c) and (4.2-d) are:

$$
\begin{aligned}
C_1 &= -2.384458 - 0.7713792i \\
C_2 &= \ \ 2.544769 + 1.3587520i \\
C_3 &= -1.094656 - 1.2282765i \\
C_4 &= \ \ 0.000232 + 0.2882812i
\end{aligned}
\tag{4.25}
$$

The generated proximity maps are shown in figures 4.33, 4.34, 4.35, 4.36, 4.37, and 4.38.

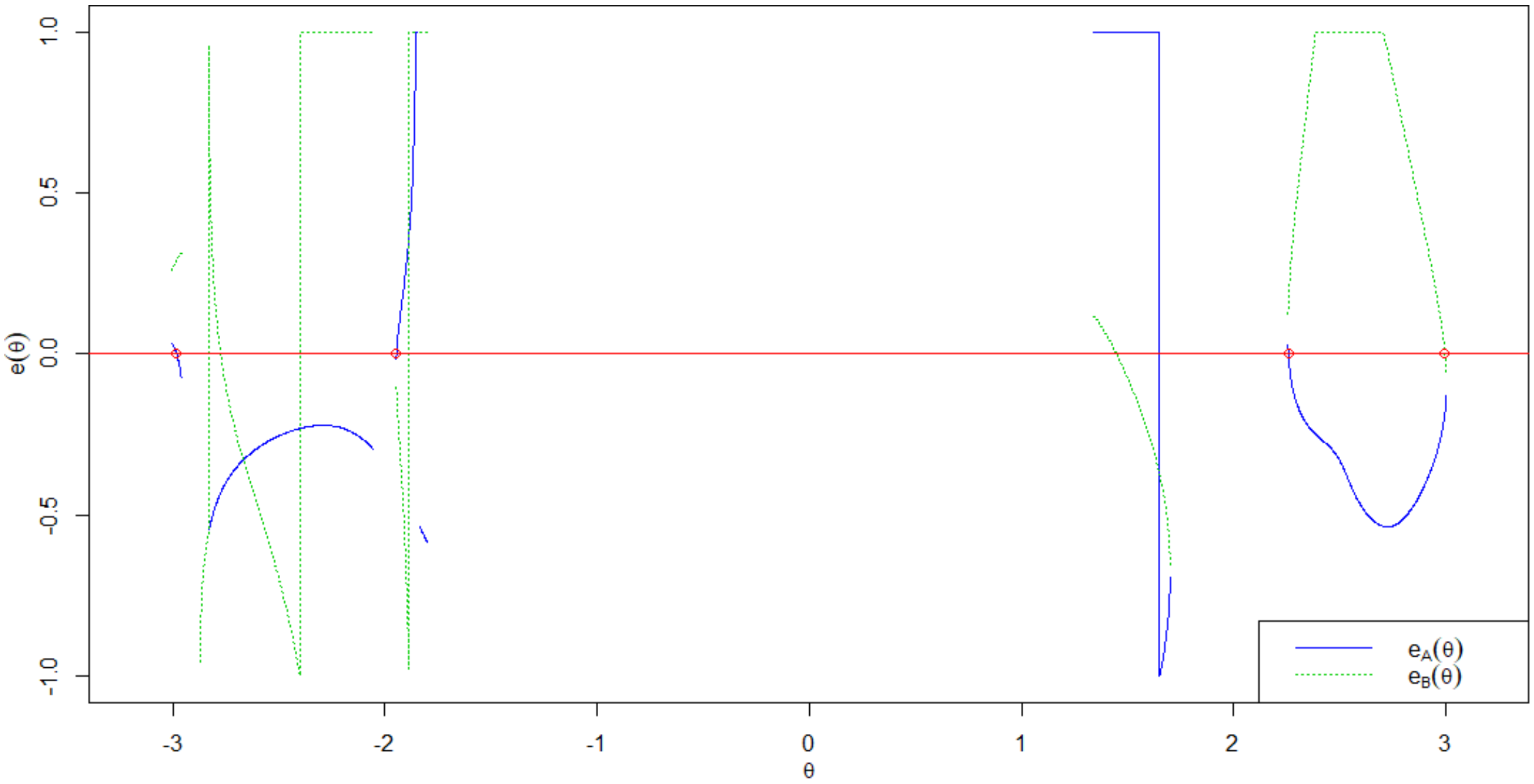

**Fig. 4.33.** Global map $e(\theta)$ associated to $p(z) = z^4 + C_1z^3 + C_2z^2 + C_3z + C_4$ with coefficients $C_1$, $C_2$, $C_3$, $C_4$ given by (4.25). This LC map was generated from $N = 5{,}000$ elements $LzC(\theta_k)$ associated with points $\theta_k$ in a regular partition of interval $[-\pi, \pi)$. Small circles on horizontal axis $y = 0$ mark the location of true theta roots $\theta_i^*$.

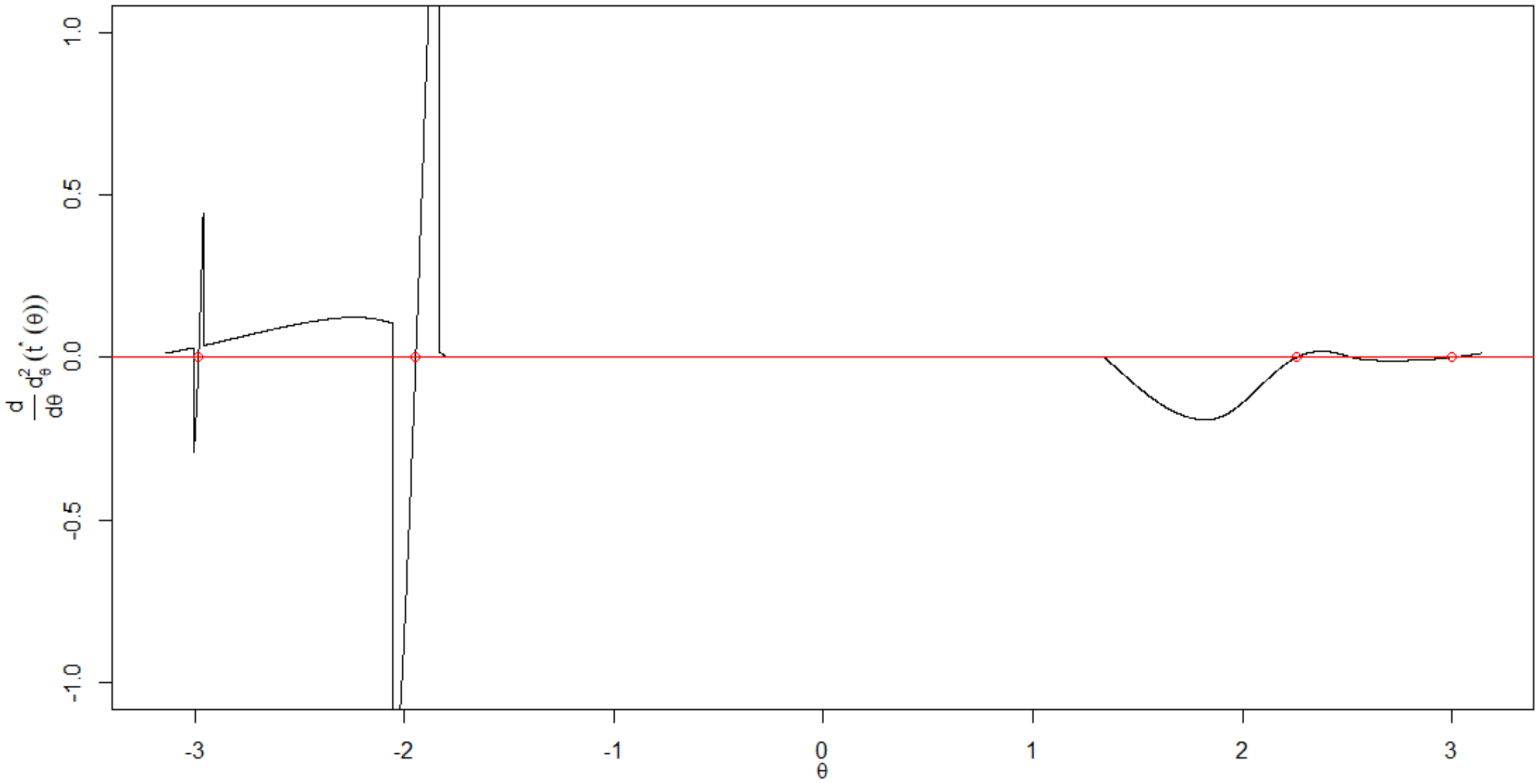

**Fig. 4.34.** Global map $\frac{d}{d\theta} d_\theta^2\left(t^*(\theta)\right)$ associated to $p(z) = z^4 + C_1 z^3 + C_2 z^2 + C_3 z + C_4$ with coefficients $C_1$, $C_2$, $C_3$, $C_4$ given by (4.25). This map was generated from $N = 5{,}000$ elements $LzC(\theta_k)$ associated with points $\theta_k$ in a regular partition of interval $[-\pi, \pi)$. Small circles on horizontal axis $y = 0$ mark the location of true theta roots $\theta_i^*$.

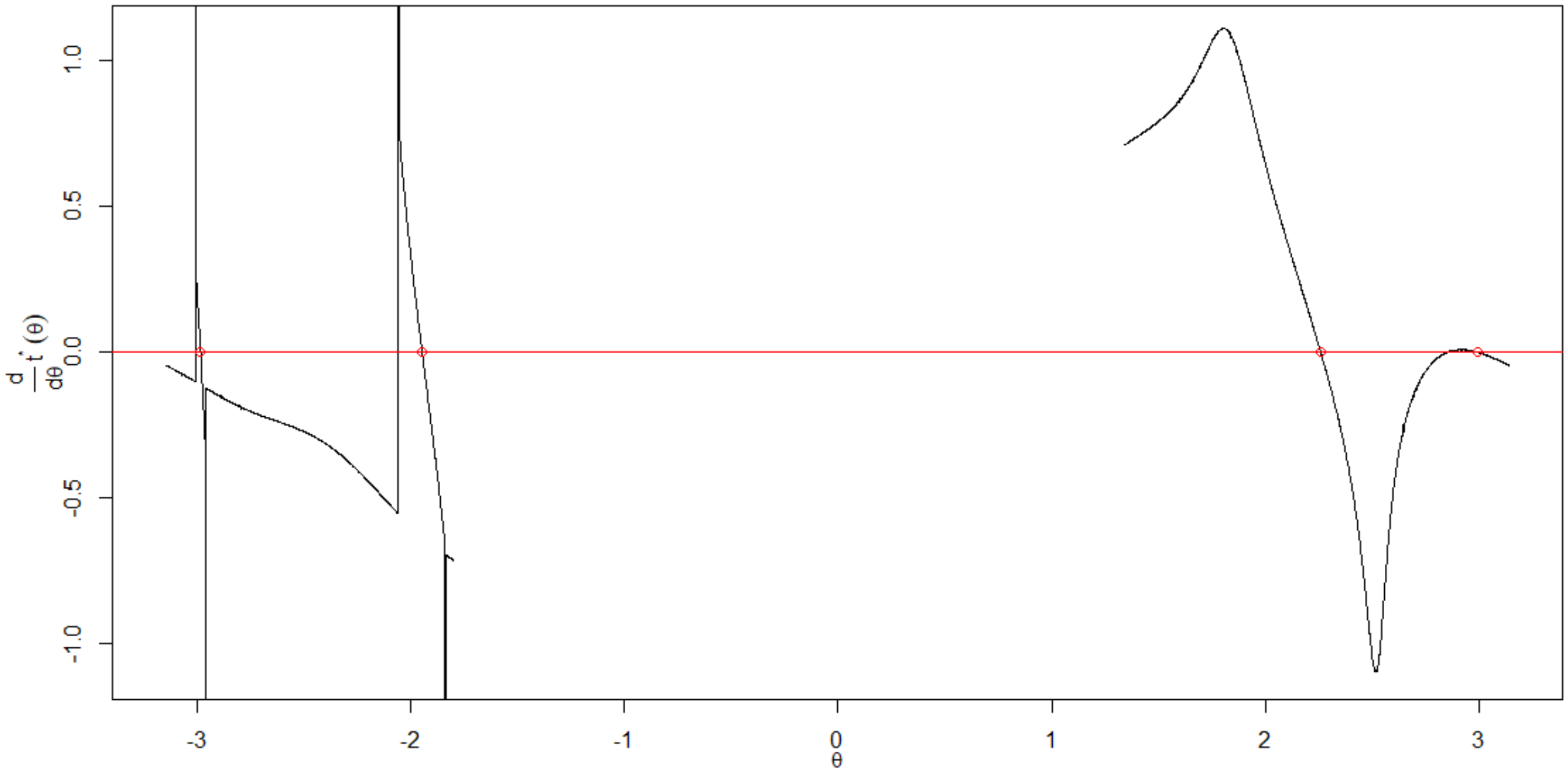

**Fig. 4.35.** Global map $\frac{d}{d\theta} t^*(\theta)$ associated to $p(z) = z^4 + C_1 z^3 + C_2 z^2 + C_3 z + C_4$ with coefficients $C_1$, $C_2$, $C_3$, $C_4$ given by (4.25). This map was generated from $N = 5{,}000$ elements $LzC(\theta_k)$ associated with points $\theta_k$ in a regular partition of interval $[-\pi, \pi)$. Small circles on horizontal axis $y = 0$ mark the location of true theta roots $\theta_i^*$.

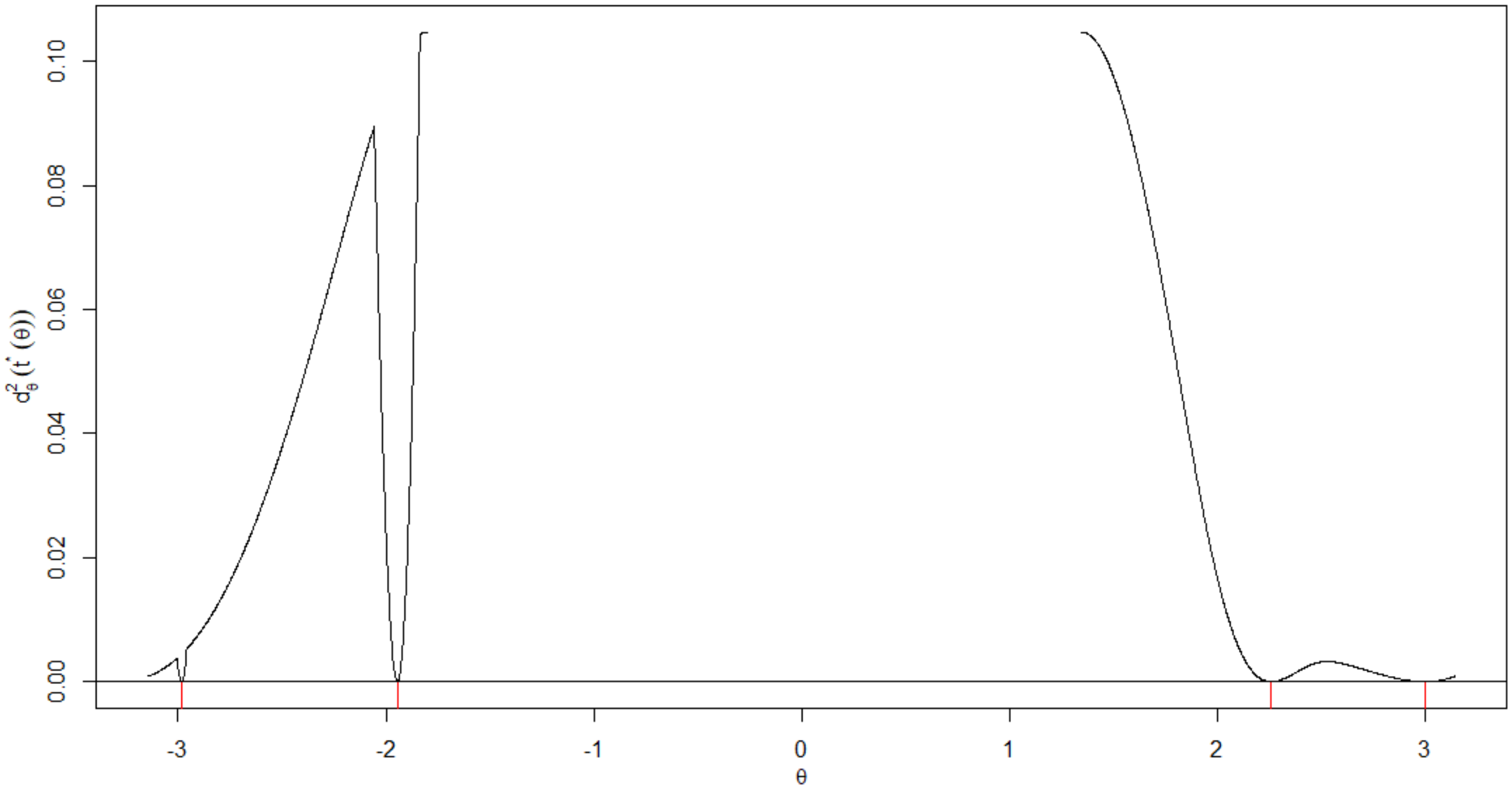

**Fig. 4.36.** Global minima $d_\theta^2(t^*(\theta))$ vs. $\theta$ associated to $p(z) = z^4 + C_1z^3 + C_2z^2 + C_3z + C_4$ with coefficients $C_1$, $C_2$, $C_3$, $C_4$ given by (4.25). This map was generated from $N = 5{,}000$ elements $LzC(\theta_k)$ associated with points $\theta_k$ in a regular partition of interval $[-\pi, \pi)$. This graph includes vertical line segments that indicate the location of true theta roots $\theta_i^*$.

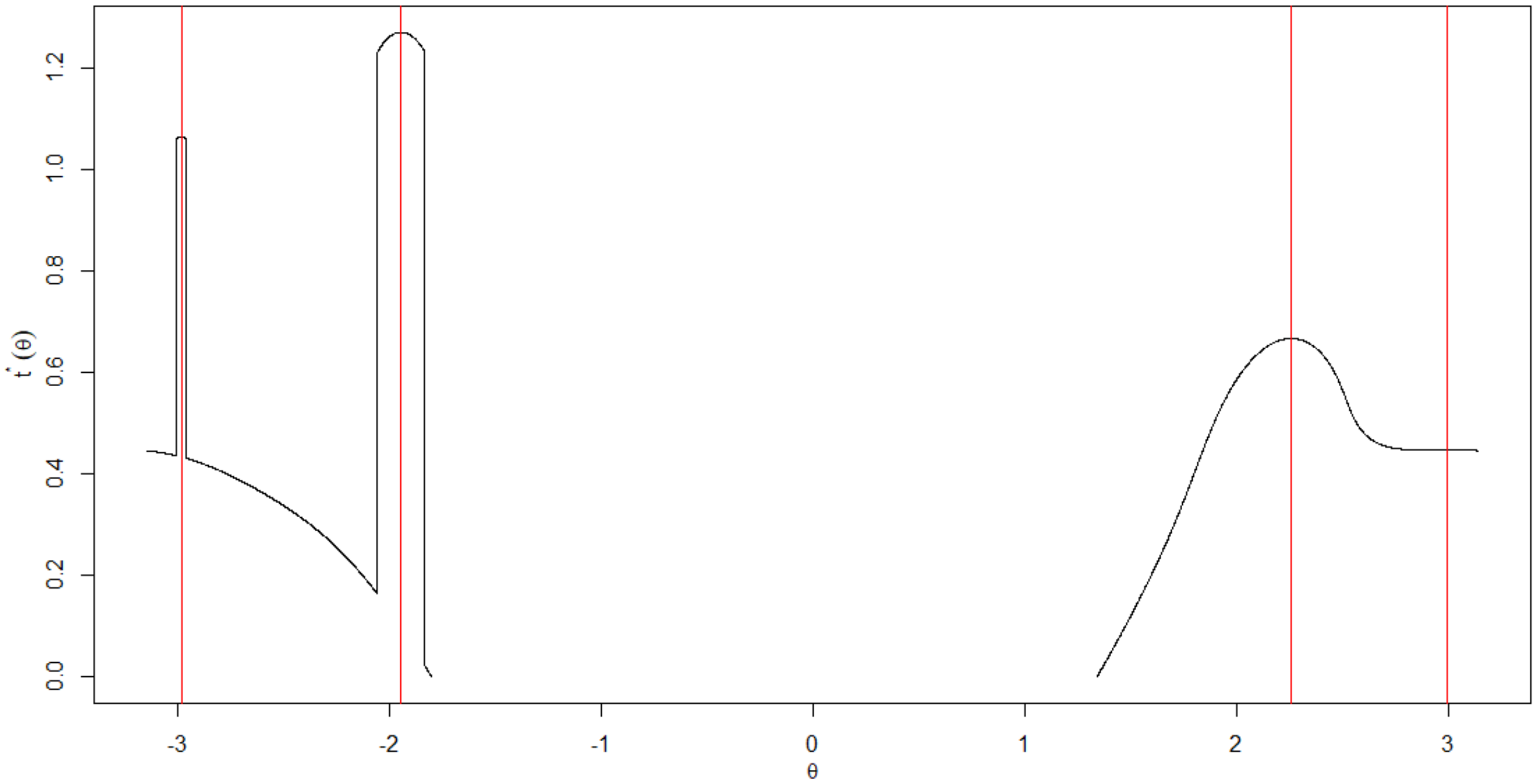

**Fig. 4.37.** Minimizing arguments $t^*(\theta)$ vs. $\theta$ associated to $p(z) = z^4 + C_1z^3 + C_2z^2 + C_3z + C_4$ with coefficients $C_1$, $C_2$, $C_3$, $C_4$ given by (4.25). This map was generated from $N = 5{,}000$ elements $LzC(\theta_k)$ associated with points $\theta_k$ in a regular partition of interval $[-\pi, \pi)$. This graph includes vertical line segments that indicate the location of true theta roots $\theta_i^*$.

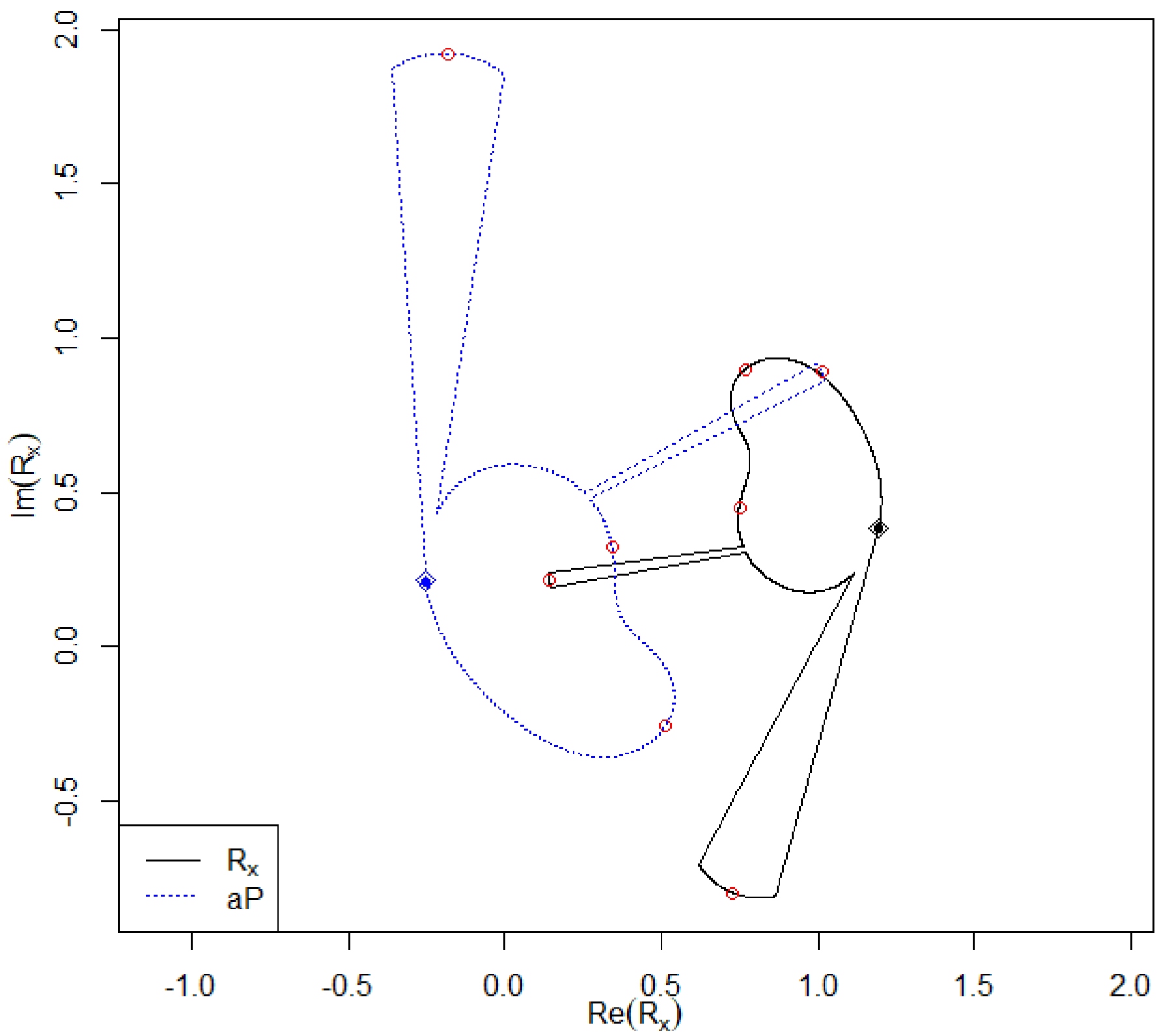

**Fig. 4.38.** Trajectories $R_x(\theta)$ and $aP(\theta)$ associated to $p(z) = z^4 + C_1z^3 + C_2z^2 + C_3z + C_4$ with coefficients $C_1$, $C_2$, $C_3$, $C_4$ given by (4.25). These trajectories were generated from $N = 5{,}000$ elements $LzC(\theta_k)$ associated with points $\theta_k$ in a regular partition of interval $[-\pi, \pi)$. The true roots $R_i$ and their corresponding anchor points are represented in this graph by means of small circles.

**Remarks on proximity maps:** From the maps in figures 4.33, 4.34, and 4.35, we see that theta roots $\theta_1^*$, $\theta_2^*$, $\theta_3^*$, $\theta_4^*$ occur in vicinities of points where the discrete versions of functions $e(\theta)$, $\frac{d}{d\theta}d_\theta^2\big(t^*(\theta)\big)$, and $\frac{d}{d\theta}t^*(\theta)$ cross the horizontal axis $y = 0$; this reinforces the hypotheses about the behavior of mappings $d_\theta^2\big(t^*(\theta)\big)$ and $t^*(\theta)$, as well as of their derivatives with respect to $\theta$, posed in example 4.1. The graph from figure 4.36, suggests that all theta roots correspond to the global minima of function $d_\theta^2\big(t^*(\theta)\big)$, while figure 4.37 suggests that these theta roots all correspond to local maxima of function $t^*(\theta)$, since, as we can see from the map in figure 4.35, function $\frac{d}{d\theta}t^*(\theta)$ is decreasing in the neighborhoods around all theta roots $\theta_i^*$; it is left as an exercise for the reader to make a high resolution close-up of mapping $t^*(\theta)$ around $\theta_2^* = 2.996166$ in order to corroborate graphically that $\theta_2^*$ corresponds to a local maximum of $t^*(\theta)$.

Figures 4.36 and 4.37 also show "peaks" around $\theta_1^* = -2.982106$ and $\theta_4^* = -1.945402$, which are transient behaviors that deviate from the usual trend of functions $d_\theta^2\big(t^*(\theta)\big)$ and $t^*(\theta)$, precisely at the neighborhoods around these roots; it is worth mentioning that the presence of "peaks with roots" is rather common in proximity maps; the amplitude of support $\theta$ where such peaks occur can be arbitrarily small, which would require a drastic increase in the resolution of proximity maps, in order to properly detect transient peaks. In a possible future improvement to the LC method, it would be necessary to consider the possibility of implementing proximity maps with dynamic resolution, in order to locally increase the resolution of these maps in the presence of transient peaks.

In figure 4.38, we see that trajectories $R_x(\theta)$ and $aP(\theta)$ consist each of a single continuous section, which closes at the point marked with a diamond ◈; these trajectories also reflect the transient peaks observed in functions $d_\theta^2\big(t^*(\theta)\big)$ and $t^*(\theta)$, and in this case the transient peaks of $R_x(\theta)$ contain roots $R_1$ and $R_4$. Maps $\frac{d}{d\theta}d_\theta^2\big(t^*(\theta)\big)$, $\frac{d}{d\theta}t^*(\theta)$, $d_\theta^2\big(t^*(\theta)\big)$ and $t^*(\theta)$ in figures 4.34, 4.35, 4.36 and 4.37 also reflect the behavior observed in the trajectories $R_x(\theta)$ and $aP(\theta)$ from figure 4.38:

- Each map is composed of a single continuous section with support $\theta$ of length $\pi$.
- At map $\frac{d}{d\theta}d_\theta^2\big(t^*(\theta)\big)$ in figure 4.34, the points that correspond to the closure ◈ of trajectories $R_x(\theta)$ and $aP(\theta)$, are separated horizontally by $\pi$ radians, and one of them approaches zero through positive values, while the other approaches zero through negative values.
- At map $\frac{d}{d\theta}t^*(\theta)$ in figure 4.35, the points that correspond to the closure ◈ of trajectories $R_x(\theta)$ and $aP(\theta)$, are separated horizontally by $\pi$ radians, and have functional values with the same magnitude, although of signs opposite to each other.
- At function $d_\theta^2\big(t^*(\theta)\big)$ in figure 4.36, the points that correspond to the closure ◈ of trajectories $R_x(\theta)$ and $aP(\theta)$, are separated horizontally by $\pi$ radians, and both correspond to the global maximum functional value.
- At function $t^*(\theta)$ in figure 4.37, the points that correspond to the closure ◈ of trajectories $R_x(\theta)$ and $aP(\theta)$, are separated horizontally by $\pi$ radians, and both correspond to the global minimum functional value.

From what is stated here, it is evident that the two values $\theta$ involved in the closure ◈ of trajectories $R_x(\theta)$ and $aP(\theta)$ are the same in any of these proximity maps; moreover, these two points are the ends of the continuous section shown in figures 4.34, 4.35, 4.36 and 4.37.

**Remarks on the root estimates generated by the proximity maps:** The numerical characteristics for the smooth crossings in the proximity maps from figures 4.33, 4.34 and 4.35 are shown in tables 4.17, 4.18 and 4.19, in which rows are ranked according to variable $d_{\hat{\theta}_i^*}^2\left(t^*\left(\hat{\theta}_i^*\right)\right)$.

**Table 4.17**. Initial estimates of the roots of $p(z) = z^4 + C_1z^3 + C_2z^2 + C_3z + C_4$ with coefficients $C_1$, $C_2$, $C_3$, $C_4$ given by (4.25), obtained from map $e(\theta)$ (LC map) in figure 4.33.

| $i$ | $\hat{R}_i$ | $\hat{\theta}_i^*$ | $\lvert\Delta e_i\rvert$ | $d^2_{\hat{\theta}_i^*}\left(t^*(\hat{\theta}_i^*)\right)$ |
|---|---|---|---|---|
| 1 | $0.7478237 + 0.4507816i$ | 2.996157 | 0.009423658 | $3.196961 \times 10^{-12}$ |
| 2 | $0.7685745 + 0.9011212i$ | 2.258773 | 0.014607805 | $2.387861 \times 10^{-10}$ |
| 3 | $0.1408590 + 0.2165791i$ | −2.982111 | 0.002640887 | $3.468903 \times 10^{-10}$ |
| 4 | $0.7272073 - 0.7971188i$ | −1.945384 | 0.018159596 | $2.824243 \times 10^{-9}$ |
| 5 | $0.8172218 + 0.2415929i$ | −2.774737 | 0.006226965 | $1.456718 \times 10^{-2}$ |
| 6 | $1.2019145 + 0.4646419i$ | 1.448728 | 0.001686796 | $1.014673 \times 10^{-1}$ |

**Table 4.18**. Initial estimates of the roots of $p(z) = z^4 + C_1z^3 + C_2z^2 + C_3z + C_4$ with coefficients $C_1$, $C_2$, $C_3$, $C_4$ given by (4.25), obtained from map $\frac{d}{d\theta}d^2_\theta\left(t^*(\theta)\right)$ in figure 4.34.

| $i$ | $\hat{R}_i$ | $\hat{\theta}_i^*$ | $\left\lvert\Delta \frac{d}{d\theta} d^2_{\hat{\theta}_i^*}(t^*)\right\rvert$ | $d^2_{\hat{\theta}_i^*}\left(t^*(\hat{\theta}_i^*)\right)$ |
|---|---|---|---|---|
| 1 | $0.7478231 + 0.4507777i$ | 2.996166 | $8.899748 \times 10^{-5}$ | $3.135157 \times 10^{-14}$ |
| 2 | $0.7685938 + 0.9011371i$ | 2.258736 | $4.054657 \times 10^{-4}$ | $1.622182 \times 10^{-12}$ |
| 3 | $0.7271856 - 0.7971088i$ | −1.945403 | $2.150170 \times 10^{-2}$ | $1.352073 \times 10^{-11}$ |
| 4 | $0.1408598 + 0.2165793i$ | −2.982111 | $2.142514 \times 10^{-2}$ | $4.109254 \times 10^{-10}$ |
| 5 | $0.7524220 + 0.6968413i$ | 2.525868 | $2.152893 \times 10^{-4}$ | $3.265132 \times 10^{-3}$ |
| 6 | $0.7598892 + 0.3267639i$ | −3.006132 | $1.004921 \times 10^{-1}$ | $3.742073 \times 10^{-3}$ |
| 7 | $1.1155414 + 0.2407919i$ | −2.057571 | $7.655708 \times 10^{-1}$ | $8.943054 \times 10^{-2}$ |

**Table 4.19**. Initial estimates of the roots of $p(z) = z^4 + C_1z^3 + C_2z^2 + C_3z + C_4$ with coefficients $C_1$, $C_2$, $C_3$, $C_4$ given by (4.25), obtained from map $\frac{d}{d\theta}t^*(\theta)$ in figure 4.35.

| $i$ | $\hat{R}_i$ | $\hat{\theta}_i^*$ | $\left\lvert\Delta \frac{d}{d\theta} t^*(\hat{\theta}_i^*)\right\rvert$ | $d^2_{\hat{\theta}_i^*}\left(t^*(\hat{\theta}_i^*)\right)$ |
|---|---|---|---|---|
| 1 | $0.7685687 + 0.9011190i$ | 2.258782 | 0.0014934152 | $3.575552 \times 10^{-10}$ |
| 2 | $0.7478754 + 0.4511336i$ | 2.995365 | 0.0003226535 | $2.276700 \times 10^{-8}$ |
| 3 | $0.7477676 + 0.4503992i$ | 2.997017 | 0.0006190431 | $2.569399 \times 10^{-8}$ |
| 4 | $0.7272582 - 0.7971377i$ | −1.945341 | 0.0070001910 | $3.170704 \times 10^{-8}$ |
| 5 | $0.1408713 + 0.2165084i$ | −2.982043 | 0.0144145755 | $3.343996 \times 10^{-8}$ |
| 6 | $0.7478907 + 0.4512393i$ | 2.995127 | 0.0008132729 | $3.827011 \times 10^{-8}$ |
| 7 | $0.7630109 + 0.5147078i$ | 2.849596 | 0.0017218014 | $6.240191 \times 10^{-4}$ |

In table 4.17, the first 4 rows show initial estimates derived from the map $e(\theta)$ that correspond respectively to reference roots $R_2$, $R_3$, $R_1$, $R_4$ in (4.24); rows 5 and 6 are smooth crossings of $\hat{e}_B(\theta)$ with horizontal axis $y = 0$ that do not approach any root of polynomial $p(z) = z^4 + C_1z^3 + C_2z^2 + C_3z + C_4$ with coefficients $C_1$, $C_2$, $C_3$, $C_4$ given by (4.25).

In table 4.18, the first 4 rows show initial estimates derived from the map $\frac{d}{d\theta}d_\theta^2\big(t^*(\theta)\big)$ that correspond respectively to reference roots $R_2$, $R_3$, $R_4$, $R_1$ in (4.24); rows 5, 6 and 7 are crossings of $\frac{d}{d\theta}d_\theta^2\big(t^*(\theta)\big)$ with horizontal axis $y = 0$ that do not approach any root of polynomial $p(z) = z^4 + C_1z^3 + C_2z^2 + C_3z + C_4$ with coefficients $C_1$, $C_2$, $C_3$, $C_4$ given by (4.25).

With respect to reference roots in (4.24), table 4.19 shows in its first row an estimate of $R_3$ generated by map $\frac{d}{d\theta}t^*(\theta)$, while its rows 2, 3 and 6 show estimates of $R_2$ (of course the best approximation is in the second row, while the other two rows are "replicas" attributed to oscillatory noise in map $\frac{d}{d\theta}t^*(\theta)$); its row 4 shows an estimate of $R_4$, its row 5 shows an estimate of $R_1$, and its row 7 shows a smooth crossing of $\frac{d}{d\theta}t^*(\theta)$ with horizontal axis $y = 0$ that does not correspond to any root of polynomial $p(z) = z^4 + C_1z^3 + C_2z^2 + C_3z + C_4$ with coefficients $C_1$, $C_2$, $C_3$, $C_4$ given by (4.25).

In this example, the amount of oscillatory noise in the map $\frac{d}{d\theta}t^*(\theta)$ from figure 4.35 is not as large as in its counterpart from example 4.3, also generated from $N = 5{,}000$ elements $LzC(\theta_k)$ associated with points $\theta_k$ in a regular partition of interval $[-\pi, \pi)$ (see figure 4.27). This suggests a more stable numerical behavior for the maps in this example.

Table 4.20 shows a summary of the approximation errors for the valid initial estimates in tables 4.17, 4.18 and 4.19, with respect to corresponding reference values (4.24).

**Table 4.20**. Arithmetic means and standard deviations for error measures associated with valid initial approximations (first four rows in each of tables 4.17 and 4.18, and rows 1, 2, 4 and 5 of table 4.19) to the roots of $p(z) = z^4 + C_1 z^3 + C_2 z^2 + C_3 z + C_4$ with coefficients $C_1$, $C_2$, $C_3$, $C_4$ given by (4.25), obtained by means of discrete proximity maps $e(\theta)$, $\frac{d}{d\theta} d_\theta^2\big(t^*(\theta)\big)$, and $\frac{d}{d\theta} t^*(\theta)$ from figures 4.33, 4.34 and 4.35.

| Map | Error measure | | | | | |
|---|---|---|---|---|---|---|
| | $d^2_{\widehat{\theta}_i^*}\left(t^*(\widehat{\theta}_i^*)\right)$ | | $\left\lvert\frac{R-\widehat{R}}{R}\right\rvert$ | | $\left\lvert\frac{\theta^*-\widehat{\theta}^*}{\theta^*}\right\rvert$ | |
| | **Mean** | **Standard deviation** | **Mean** | **Standard deviation** | **Mean** | **Standard deviation** |
| $e(\theta)$ | $8.532791 \times 10^{-10}$ | $1.321788 \times 10^{-9}$ | $1.858379 \times 10^{-5}$ | $9.405103 \times 10^{-6}$ | $7.792315 \times 10^{-6}$ | $6.992159 \times 10^{-6}$ |
| $\frac{d}{d\theta} d_\theta^2(t^*(\theta))$ | $1.065249 \times 10^{-10}$ | $2.030229 \times 10^{-10}$ | $8.100859 \times 10^{-6}$ | $1.370959 \times 10^{-5}$ | $6.633863 \times 10^{-7}$ | $6.813677 \times 10^{-7}$ |
| $\frac{d}{d\theta} t^*(\theta)$ | $2.206789 \times 10^{-8}$ | $1.521035 \times 10^{-8}$ | $1.922672 \times 10^{-4}$ | $1.777707 \times 10^{-4}$ | $8.518092 \times 10^{-5}$ | $1.216986 \times 10^{-4}$ |

From table 4.20, we see once again that the map $\frac{d}{d\theta} d_\theta^2\big(t^*(\theta)\big)$ is the one which produces, on average, the most accurate estimates; in close second place, there is the map $e(\theta)$, and in third place is the map $\frac{d}{d\theta} t^*(\theta)$. We consider all initial estimates obtained by the three types of proximity maps in this example to be reasonably good; in the worst case, the initial estimate in row 5 of table 4.19 from the map $\frac{d}{d\theta} t^*(\theta)$, coincides with the first 4 digits of each of the components from reference root $R_1$, $\theta_1^*$ in (4.24).

# Chapter Five: LC Method for Polynomials of Degree $n \geq 4$

## Theoretical Aspects

In this chapter we will generalize the concepts seen in chapter 4; as we shall see, all the elements of the LC method used in the approximation to roots for univariate polynomials of degree 4 can be generalized directly to any univariate polynomial of degree $n \geq 5$, without the need to define additional new concepts, as was done in the previous chapters, when we went from polynomials of degree 2 to polynomials of degree 3, and then from polynomials of degree 3 to polynomials of degree 4. Therefore, all the elements seen in chapter 4 for univariate polynomials of degree 4 constitute the basis of the *complete* LC method; it only remains for us to generalize that basis in order to cover univariate polynomials of any degree $n \geq 4$.

We begin by establishing, as in previous chapters, that we seek to approximate the roots $R_1, R_2, \ldots, R_n \in \mathbb{C}, n \geq 4$, of equation

$$z^n + C_1 z^{n-1} + C_2 z^{n-2} + \cdots + C_{n-2} z^2 + C_{n-1} z + C_n = 0 \tag{5.1}$$

As we have already seen in chapter 1, the coefficient $C_k$ in equation (5.1), with $k = 1,2, \ldots, n$, is equal, according to Vieta's relations, to the sum of all possible $k$-products of the roots of equation (5.1), multiplied by $(-1)^k$; i.e.,

$$C_k = (-1)^k \sum_{\substack{i_1 < i_2 < \cdots < i_k \\ i_1, i_2, \ldots, i_k \in \{1,2,\ldots,n\}}} R_{i_1} \cdot R_{i_2} \cdot \cdots \cdot R_{i_k} \tag{5.2}$$

The number of $k$-products in the summation on the right side of equation (5.2) is $\binom{n}{k} = \frac{n!}{k!(n-k)!}$; that is to say, it is the number of combinations that can be formed with $n$ different roots, taking $k$ roots at the same time; in this formula for counting combinations,

$k! = \begin{cases} k \cdot (k-1)! & \text{if } k > 0 \\ 1 & \text{if } k = 0 \end{cases}$ is the factorial of $k$.

Of course, in the exposition of this chapter we will also suppose that the roots $R_1, R_2, \ldots, R_n \in \mathbb{C}$ are different from each other, and all are different from 0.

As in previous chapters, we seek to build structures $LzC(C_1, C_2, \dots, C_n, \theta)$ such that

1. $\ell_1(\theta): -C_1/2 + te^{i\theta} \to R_1$, when $\theta \to \theta^*$ (line $\ell_1(\theta)$ with fixed point $P_1 = -C_1/2$ tends to $R_1$, as $\ell_1$'s inclination angle $\theta$ tends to $\theta^*$).

2. $(-1)^n C_n/\ell_1(\theta) \to R_2 \cdot R_3 \cdot \dots \cdot R_n$, $tL(\theta) \to R_2 \cdot R_3 \cdot \dots \cdot R_n$ and $t\mathfrak{C}(\theta) \to R_2 R_3 \cdot \dots \cdot R_n$, when $\theta \to \theta^*$ (z-circumference $zC(\theta) = (-1)^n C_n/\ell_1(\theta)$, terminal semi-line $tL(\theta)$, and terminal curve $t\mathfrak{C}(\theta)$ tend to $R_2 \cdot R_3 \cdot \dots \cdot R_n$ as $\theta$ tends to $\theta^*$).

3. $\min_{t\in\mathbb{R}} d_\theta^2(t) \to 0$, when $\theta \to \theta^*$ (the global minimum of the dynamic squared distance between $t\mathfrak{C}(\theta)$ and $zC(\theta)$ tends to 0 as $\theta$ tends to $\theta^*$).

4. $e(\theta) \to 0$, when $\theta \to \theta^*$ (weighted error $e(\theta)$ derived from the intersection between $tL(\theta)$ and $zC(\theta)$ tends to 0 as $\theta$ tends to $\theta^*$).

Note: Without loss of generality, we use $R_1$, $\theta^*$ in statements 1 to 4 above to refer to any of the roots of equation (5.1).

Of course, if in statement 1 the containment condition of root $R_1$ in line $\ell_1(\theta)$ is fulfilled, automatically in statement 2 the containment conditions of $(n-1)$-product $R_2 \cdot R_3 \cdot \dots \cdot R_n$ in trajectories $zC(\theta)$, $tL(\theta)$ and $t\mathfrak{C}(\theta)$ will also be fulfilled; likewise, in statement 3 equation $\min_{t\in\mathbb{R}} d_\theta^2(t) = 0$ will be satisfied, and in statement 4 equation $e(\theta) = 0$ will also be satisfied.

If the equation $\min_{t\in\mathbb{R}} d_\theta^2(t) = 0$ in statement 3 is satisfied, this is a strong indicator that the associated structure $LzC(\theta)$ contains $R_1$ in $\ell_1(\theta)$, and $R_2 \cdot R_3 \cdot \dots \cdot R_n$ in $zC(\theta)$, $tL(\theta)$ and $t\mathfrak{C}(\theta)$. On the other hand, for univariate polynomials of degree $n \geq 4$, the fulfillment of condition $e(\theta) = 0$ in statement 4 does not necessarily guarantee that the associated structure $LzC(\theta)$ contains $R_1$ (see examples 4.2 and 4.4 in chapter 4).

**General algebraic expressions for the components of geometric structure $\boldsymbol{LzC(C_1, C_2, \dots, C_n, \theta)}$**

Next, we will find general algebraic expressions for the main geometric components of $LzC(C_1, C_2, \dots, C_n, \theta)$. We will start with the terminal curve $t\mathfrak{C}(\theta, t)$. As we saw in expression (4.6) of chapter 4, for $n = 4$ we have

$$t\mathfrak{C}(\theta, t) = -C_3 - (P_1 + tv_\theta)[C_2 - (P_1 + tv_\theta)(P_1 - tv_\theta)] \tag{5.3}$$

Now, for $n=5$, we can see that, when $\theta \to \theta^*$,

$\ell_1(\theta,t) = P_1 + tv_\theta \to R_1$

$\ell_d(\theta,t) = (P_1+tv_\theta)(P_1-tv_\theta) \to R_1R_2 + R_1R_3 + R_1R_4 + R_1R_5$

$C_2 - \ell_d(\theta,t) \to R_2R_3 + R_2R_4 + R_2R_5 + R_3R_4 + R_3R_5 + R_4R_5$

$(P_1+tv_\theta)[C_2 - \ell_d(\theta,t)] \to R_1R_2R_3 + R_1R_2R_4 + R_1R_2R_5 + R_1R_3R_4 + R_1R_3R_5 + R_1R_4R_5$

$-C_3 - (P_1+tv_\theta)[C_2 - \ell_d(\theta,t)] \to R_2R_3R_4 + R_2R_3R_5 + R_2R_4R_5 + R_3R_4R_5$

$-(P_1+tv_\theta)C_3 - (P_1+tv_\theta)^2[C_2 - \ell_d(\theta,t)] \to R_1R_2R_3R_4 + R_1R_2R_3R_5 + R_1R_2R_4R_5 + R_1R_3R_4R_5$

$C_4 + (P_1+tv_\theta)C_3 + (P_1+tv_\theta)^2[C_2 - \ell_d(\theta,t)] \to R_2R_3R_4R_5$

Therefore, for $n=5$,

$$t\mathfrak{C}(\theta,t) = C_4 + (P_1+tv_\theta)C_3 + (P_1+tv_\theta)^2[C_2 - \ell_d(\theta,t)] \tag{5.4}$$

It is easy to see that, for $n=6$, following the same reasoning,

$$t\mathfrak{C}(\theta,t) = -C_5 - (P_1+tv_\theta)C_4 - (P_1+tv_\theta)^2C_3 - (P_1+tv_\theta)^3[C_2 - \ell_d(\theta,t)] \tag{5.5}$$

If we analyze the structure of equations (5.3), (5.4) and (5.5), we conclude inductively that a general pattern emerges for the algebraic expression of terminal curve $t\mathfrak{C}(\theta,t)$:

**General algebraic expression for terminal curve $t\mathfrak{C}(\theta,t)$ associated with equation (5.1)**

$$t\mathfrak{C}(\theta,t) = (-1)^{n+1}\left[\sum_{i=1}^{n-3} C_{n-i}\,(P_1+tv_\theta)^{i-1} + (P_1+tv_\theta)^{n-3}(C_2 - \ell_d(\theta,t))\right]. \tag{5.6}$$

$v_\theta$ is the unit direction vector of line $\ell_1(\theta,t)$, $P_1 = -C_1/2$ is the fixed point of $\ell_1(\theta,t)$, and $\ell_d(\theta,t) = (P_1+tv_\theta)(P_1-tv_\theta)$.

The general algebraic expression for the z-circumference is obtained immediately:

**General algebraic expression for the z-circumference $zC(\theta,t)$ associated with equation (5.1)**

$$zC(\theta,t) = (-1)^n C_n/(P_1+tv_\theta) \tag{5.7}$$

Now, we can define in a general way the dynamic squared distance $d_\theta^2(t)$ between $t\mathfrak{C}(\theta,t)$ and $zC(\theta,t)$:

**General algebraic expression for the dynamic squared distance $\boldsymbol{d_\theta^2(t)}$ between terminal curve $\boldsymbol{t\mathfrak{C}(\theta,t)}$ and z-circumference $\boldsymbol{zC(\theta,t)}$ associated with equation (5.1)**

$$d_\theta^2(t) = \left[\mathrm{Re}\big(t\mathfrak{C}(\theta,t) - zC(\theta,t)\big)\right]^2 + \left[\mathrm{Im}\big(t\mathfrak{C}(\theta,t) - zC(\theta,t)\big)\right]^2, \tag{5.8}$$

where $t\mathfrak{C}(\theta,t) = (-1)^{n+1}\left[\sum_{i=1}^{n-3} C_{n-i}\,(P_1 + tv_\theta)^{i-1} + (P_1 + tv_\theta)^{n-3}\big(C_2 - \ell_d(\theta,t)\big)\right]$,

$zC(\theta,t) = (-1)^n C_n/(P_1 + tv_\theta)$; $\qquad \theta, t \in \mathbb{R}$.

The "best" approximation to $R_i$, given the inclination angle $\theta$ of line $\ell_1$, is given by:

**General algebraic expression for the "best" approximation to one of the roots $\boldsymbol{R_i}$ $(\boldsymbol{i = 1, 2, \dots, n})$, given the inclination angle $\boldsymbol{\theta}$ of line $\boldsymbol{\ell_1}$ associated with equation (5.1)**

$$R_x(\theta) = P_1 + t^*(\theta)v_\theta, \tag{5.9}$$

where $t^*(\theta) = \underset{t\in\mathbb{R}}{\arg\min}\, d_\theta^2(t)$ is the minimizing argument of function (5.8).

Once we have obtained $R_x(\theta)$, now we can derive general expressions for the anchor point and direction vector of the terminal semi-line $tL(\theta,t)$ associated with equation (5.1).

From chapter 4, we saw that $tL(\theta,t)$, for the case $n = 4$, is given by:

$$tL(\theta,t)\colon -C_3 - R_x(\theta)[C_2 - \ell_d(\theta,t)] \tag{5.10}$$

$$tL(\theta,t)\colon -C_3 - R_x(\theta)C_2 + R_x(\theta)(P_1 + tv_\theta)(P_1 - tv_\theta)$$

$$tL(\theta,t)\colon [-C_3 - R_x(\theta)(C_2 - P_1^2)] + t^2[-R_x(\theta)v_\theta^2] \tag{5.11}$$

From (5.11), we see that an anchor point $aP(\theta)$ and a direction vector $v_{tL}(\theta)$ of $tL(\theta,t)$, for $n = 4$, is

$$aP(\theta) = -C_3 - R_x(\theta)(C_2 - P_1^2), \qquad v_{tL}(\theta) = -R_x(\theta)v_\theta^2 \tag{5.12}$$

Let us first notice that $R_x(\theta)$ acts as a fixed value, as do the coefficients $P_1^2$, $C_2$, $C_3$ and vector $v_\theta^2$, while $t^2$ is a non-negative real variable parameter; these conditions guarantee $tL(\theta,t)$ has the form of a semi-line. Note also that if we rewrite the right-hand side of expression (5.3) as

$-C_3 - (P_1 + tv_\theta)[C_2 - \ell_d(\theta,t)]$,

and then replace in this last expression the factor $(P_1 + tv_\theta)$ with the "fixed" value $R_x(\theta)$, we arrive at expression (5.10) and therefore to expression (5.11); this suggests that, from the general expression (5.6) for the terminal curve $t\mathfrak{C}(\theta,t)$, we can obtain a general expression for terminal semi-line $tL(\theta,t)$, by replacing factors $(P_1 + tv_\theta)^{i-1}$ and $(P_1 + tv_\theta)^{n-3}$ with $R_x(\theta)^{i-1}$ and $R_x(\theta)^{n-3}$. By doing these replacements on the right-hand side of expression (5.6), we have

$(-1)^{n+1}\left[\sum_{i=1}^{n-3} C_{n-i}\, R_x(\theta)^{i-1} + R_x(\theta)^{n-3}(C_2 - \ell_d(\theta,t))\right] =$

$(-1)^{n+1}\left[\sum_{i=1}^{n-3} C_{n-i}\, R_x(\theta)^{i-1} + R_x(\theta)^{n-3}(C_2 - (P_1 + tv_\theta)(P_1 - tv_\theta))\right] =$

$(-1)^{n+1}\left[\sum_{i=1}^{n-3} C_{n-i}\, R_x(\theta)^{i-1} + R_x(\theta)^{n-3}(C_2 - P_1^2 + t^2v_\theta^2)\right] =$

$(-1)^{n+1}\left[\sum_{i=1}^{n-3} C_{n-i}\, R_x(\theta)^{i-1} + R_x(\theta)^{n-3}(C_2 - P_1^2) + t^2R_x(\theta)^{n-3}v_\theta^2\right] =$

$(-1)^{n+1}\left[\sum_{i=1}^{n-3} C_{n-i}\, R_x(\theta)^{i-1} + R_x(\theta)^{n-3}(C_2 - P_1^2)\right] + (-1)^{n+1}[t^2R_x(\theta)^{n-3}v_\theta^2] =$

$(-1)^{n+1}\left[\sum_{i=1}^{n-3} C_{n-i}\, R_x(\theta)^{i-1} + R_x(\theta)^{n-3}(C_2 - P_1^2)\right] + t^2[(-1)^{n+1}R_x(\theta)^{n-3}v_\theta^2]$.

In this way, we have arrived at the desired general expression. We enunciate it formally below:

**General algebraic expression for terminal semi-line $\boldsymbol{tL(\theta,t)}$ associated with equation (5.1)**

$tL(\theta,t)$: $(-1)^{n+1}\left[\sum_{i=1}^{n-3} C_{n-i}\, R_x(\theta)^{i-1} + R_x(\theta)^{n-3}(C_2 - P_1^2)\right] + t^2\left[(-1)^{n+1}R_x(\theta)^{n-3}v_\theta^2\right]$, (5.13)

where

$aP(\theta) = (-1)^{n+1}\left[\sum_{i=1}^{n-3} C_{n-i}\, R_x(\theta)^{i-1} + R_x(\theta)^{n-3}(C_2 - P_1^2)\right]$ is the anchor point of $tL(\theta,t)$, and $v_{tL}(\theta) = (-1)^{n+1}R_x(\theta)^{n-3}v_\theta^2$ is a direction vector for $tL(\theta,t)$.

$v_\theta$ is the unit direction vector of line $\ell_1(\theta,t)$, and $P_1 = -C_1/2$ is the fixed point of $\ell_1(\theta,t)$.

$R_x(\theta)$ is the best estimate of $R_i$, defined in (5.9), given the inclination angle $\theta$ of line $\ell_1$.

With the general expressions of terminal semi-line $tL(\theta)$ and z-circumference $zC(\theta)$, we can obtain numerical intersections $I_1(\theta)$, $I_2(\theta)$ between $tL(\theta)$ and $zC(\theta)$ by means of function `intersect_semiLine_circle` listed in annex 1 section 2. Once we obtain these numerical intersections, we can compute projections of $I_1(\theta)$, $I_2(\theta)$ onto $\ell_1(\theta)$ via the z-circumference $zC(\theta)$ as

$c_{0I_1}(\theta) = (-1)^nC_n/I_1(\theta)$, $\qquad c_{0I_2}(\theta) = (-1)^nC_n/I_2(\theta)$. (5.14)

Also, we can compute projections of $I_1(\theta)$, $I_2(\theta)$ onto $\ell_1(\theta)$ via the terminal semi-line $tL(\theta)$ as

$$i_{0I_1}(\theta) = P_1 + \sqrt{\frac{[I_1(\theta)-aP(\theta)]}{v_{tL}(\theta)}}\,v_\theta\ , \qquad i_{0I_2}(\theta) = P_1 + \sqrt{\frac{[I_2(\theta)-aP(\theta)]}{v_{tL}(\theta)}}\,v_\theta\ . \qquad (5.15)$$

In order to glimpse the path that leads to the expressions in (5.15), see the derivation of expression (4.10) in chapter 4.

With these elements, we can now generalize the definition of weighted error to univariate polynomials of degree $n \geq 4$:

**Weighted error measurements associated with structure $\boldsymbol{LzC(C_1, C_2, C_3, \dots, C_n, \theta)}$, $\boldsymbol{\theta \in [-\pi, \pi)}$**

$$e_A(\theta) := \begin{cases} sgn\left(\frac{i_{0I_1}-c_{0I_1}}{v_\theta}\right)\frac{|i_{0I_1}-c_{0I_1}|}{|i_{0I_1}-P_1|+|c_{0I_1}-P_1|} & \text{if } \exists\ I_1(\theta) \\ \text{undefined} & \text{if } \nexists\ I_1(\theta) \end{cases} \qquad (5.16a)$$

$$e_B(\theta) := \begin{cases} sgn\left(\frac{i_{0I_2}-c_{0I_2}}{v_\theta}\right)\frac{|i_{0I_2}-c_{0I_2}|}{|i_{0I_2}-P_1|+|c_{0I_2}-P_1|} & \text{if } \exists\ I_2(\theta) \\ \text{undefined} & \text{if } \nexists\ I_2(\theta) \end{cases} \qquad (5.16b)$$

where:

$$i_{0I_m} = P_1 + \sqrt{\frac{[I_m(\theta)-aP(\theta)]}{v_{tL}(\theta)}}\,v_\theta, \qquad c_{0I_m} = \frac{(-1)^n C_n}{I_m(\theta)}, \qquad m = 1,2.$$

$I_m(\theta)$ is one of two possible intersections between terminal semi-line $tL(\theta) = aP(\theta) + t^2 v_{tL}(\theta)$ and z-circumference $zC(\theta) = (-1)^n C_n/(P_1 + tv_\theta)$, with $aP(\theta)$ and $v_{tL}(\theta)$ given by (5.13).

$v_\theta$ is the unit direction vector of line $\ell_1(\theta)$ with fixed point $P_1 = -C_1/2$.

Finally, with all these elements, we can define a general strategy that helps us to process in parallel $N$ structures $LzC(C_1, C_2, C_3, \dots, C_n, \theta_k)$ associated with points $\theta_k$ in a regular partition of interval $[-\pi, \pi)$, in order to build a map $e(\theta)$ (LC map) that leads to initial approximations to the roots of polynomial equation (5.1):

**Strategy 5.1 to find initial approximations to the roots of polynomial equation (5.1) of degree $\boldsymbol{n \geq 4}$ by means of a discrete map $\boldsymbol{e(\theta)}$:**

a) Find (in parallel) minimizing arguments $t^*(\theta_k) = \arg\min_{t\in\mathbb{R}} d^2_{\theta_k}(t)$ associated with structures $LzC(C_1, C_2, C_3, \dots, C_n, \theta_k)$ for values $\theta_k = -\pi + 2\pi k/N$, $k = 0,1,2, \dots, N-1$; note that points $\theta_k$ constitute a regular partition of interval $[-\pi, \pi)$.

b) Evaluate (in parallel) weighted error functional values $e_A(\theta_k)$, $e_B(\theta_k)$, given by (5.16a) and (5.16b), which require as input, the values $t^*(\theta_k)$ previously computed in a).

c) Join, by using rectilinear segments, adjacent points $\big(\theta_k, e_A(\theta_k)\big), \big(\theta_{k+1}, e_A(\theta_{k+1})\big) \in \mathbb{R}^2\ \forall k$; in this way, we will have constructed, by linear interpolation, a **discrete approximation** $\hat{e}_A(\theta)$ of function $e_A(\theta)$ at interval $-\pi \leq \theta < \pi$. Similarly, obtain a discrete approximation $\hat{e}_B(\theta)$ of function $e_B(\theta)$ at interval $-\pi \leq \theta < \pi$.

d) Find crossings of $\hat{e}_A(\theta)$ and $\hat{e}_B(\theta)$ with horizontal axis $y = 0$; with this, we will obtain *candidates* for estimates $\hat{\theta}_i^*$ of direction angles $\theta_i^*$ that, for each $i$, cause line $\ell_1(\theta_i^*)$ to contain one of the roots $R_i$ of equation (5.1).

e) By means of parallel constructions $LzC\big(\hat{\theta}_i^*\big)$, it is also possible to obtain estimates $\hat{R}_1$, $\hat{R}_2$, ..., $\hat{R}_n$ to the roots $R_1$, $R_2$, ..., $R_n$ of equation (5.1). For this, simply obtain points $R_x\big(\hat{\theta}_i^*\big) = P_1 + t^*\big(\hat{\theta}_i^*\big)v_{\hat{\theta}_i^*}$ as possible initial approximations to the roots. For each $i$, $t^*\big(\hat{\theta}_i^*\big)$ is the minimizing argument of the dynamic squared distance function $d^2_{\hat{\theta}_i^*}(t)$ defined in (5.8).

Note that strategy 5.1 is completely analogous to strategy 4.1 defined in chapter 4, and as in that case, strategy 5.1 serves as a general guide for implementing a program that uses advanced parallel processing techniques, in order to efficiently compute the elements involved in the construction of discrete proximity maps $e(\theta)$ associated to univariate polynomials of degree $n \geq 4$. For our purposes, the numerical examples in this chapter will be constructed from the same basic script used in chapter 4 (see annex 3 section 3), which in fact was designed from the start taking into account the generalizations seen in this chapter.

As we saw in the numerical examples from chapter 4, from the implementation of strategy 4.1 to construct the discrete proximity map $e(\theta)$ based on the static intersections between terminal semi-lines $tL(\theta, t)$ and z-circumferences $zC(\theta, t)$, two new proximity maps arise, both based on the dynamic squared distances $d_\theta^2(t)$ between terminal curves $t\mathfrak{C}(\theta, t)$ and z-circumferences $zC(\theta, t)$; one of these maps is obtained by a discrete approximation to function $\frac{d}{d\theta} d_\theta^2\big(t^*(\theta)\big)$, while the other is obtained by a discrete approximation to function $\frac{d}{d\theta} t^*(\theta)$. As with the map $e(\theta)$, theta roots candidates $\hat{\theta}_i^*$ are found by means of crossings of the map $\frac{d}{d\theta} d_\theta^2\big(t^*(\theta)\big)$ or of the map $\frac{d}{d\theta} t^*(\theta)$ with horizontal axis $y = 0$. Next, we define the strategies to construct each of these additional proximity maps $\frac{d}{d\theta} d_\theta^2\big(t^*(\theta)\big)$ and $\frac{d}{d\theta} t^*(\theta)$.

**Strategy 5.2 to find initial approximations to the roots of polynomial equation (5.1) of degree $\boldsymbol{n \geq 4}$ by means of a discrete map $\boldsymbol{\frac{d}{d\theta} d_\theta^2\big(t^*(\theta)\big)}$:**

a) Find (in parallel) minimizing arguments $t^*(\theta_k) = \arg\min_{t \in \mathbb{R}} d_{\theta_k}^2(t)$ associated with structures $LzC(C_1, C_2, C_3, \dots, C_n, \theta_k)$ for values $\theta_k = -\pi + 2\pi k/N$, $k = 0,1,2, \dots, N-1$; note that points $\theta_k$ constitute a regular partition of interval $[-\pi, \pi)$.

b) Approximate (in parallel) the function $\frac{d}{d\theta} d_\theta^2\big(t^*(\theta)\big)$ by means of discrete difference quotients $\Delta d_{\theta_j}^2/\Delta\theta = \left[d_{\theta_j}^2\left(t^*(\theta_j)\right) - d_{\theta_{j-1}}^2\left(t^*(\theta_{j-1})\right)\right]/\Delta\theta$, with $j = 1,2, \dots, N-1$; $\Delta\theta = \theta_j - \theta_{j-1} = 2\pi/N$. Since $\frac{d}{d\theta} d_\theta^2\big(t^*(\theta)\big)$ is periodic with fundamental period $2\pi$ at interval $[-\pi, \pi)$, we also consider the difference quotient $\Delta d_{\theta_0}^2/\Delta\theta = \left[d_{\theta_0}^2\big(t^*(\theta_0)\big) - d_{\theta_{N-1}}^2\big(t^*(\theta_{N-1})\big)\right]/(2\pi/N)$; with this, we have the set of values $\left\{\Delta d_{\theta_0}^2/\Delta\theta, \Delta d_{\theta_1}^2/\Delta\theta, \Delta d_{\theta_2}^2/\Delta\theta, \dots, \Delta d_{\theta_{N-1}}^2/\Delta\theta\right\}$, whose support is the discrete set $\{\theta_k - \Delta\theta/2\}$; thus, $\Delta d_{\theta_k}^2/\Delta\theta$ is the approximation of $\frac{d}{d\theta} d_\theta^2\big(t^*(\theta)\big)$, evaluated at point $\theta = \theta_k - \Delta\theta/2$.

c) Join, by using rectilinear segments, adjacent points $\left(\theta_k - \Delta\theta/2, \Delta d_{\theta_k}^2/\Delta\theta\right)$, $\left(\theta_{k+1} - \Delta\theta/2, \Delta d_{\theta_{k+1}}^2/\Delta\theta\right) \in \mathbb{R}^2$ $\forall k$; in this way, we will have constructed, by linear interpolation, a **discrete approximation** $\frac{\Delta}{\Delta\theta} d_\theta^2\big(t^*(\theta)\big)$ of function $\frac{d}{d\theta} d_\theta^2\big(t^*(\theta)\big)$ at interval $-\pi \leq \theta < \pi$.

d) Find crossings of $\frac{\Delta}{\Delta\theta} d_\theta^2\big(t^*(\theta)\big)$ with horizontal axis $y = 0$; with this, we will obtain *candidates* for estimates $\hat{\theta}_i^*$ of direction angles $\theta_i^*$ that, for each $i$, cause line $\ell_1(\theta_i^*)$ to contain one of the roots $R_i$ of equation (5.1).

e) By means of parallel constructions $LzC\big(\hat{\theta}_i^*\big)$, it is also possible to obtain estimates $\hat{R}_1$, $\hat{R}_2$, ..., $\hat{R}_n$ to the roots $R_1$, $R_2$, ..., $R_n$ of equation (5.1). For this, simply obtain points $R_x\big(\hat{\theta}_i^*\big) = P_1 + t^*\big(\hat{\theta}_i^*\big) v_{\hat{\theta}_i^*}$ as possible initial approximations to the roots. For each $i$, $t^*\big(\hat{\theta}_i^*\big)$ is the minimizing argument of the dynamic squared distance function $d_{\hat{\theta}_i^*}^2(t)$ defined in (5.8).

**Strategy 5.3 to find initial approximations to the roots of polynomial equation (5.1) of degree $\boldsymbol{n \geq 4}$ by means of a discrete map $\boldsymbol{\frac{d}{d\theta} t^*(\theta)}$:**

a) Find (in parallel) minimizing arguments $t^*(\theta_k) = \arg\min_{t\in\mathbb{R}} d^2_{\theta_k}(t)$ associated with structures $LzC(C_1, C_2, C_3, \ldots, C_n, \theta_k)$ for values $\theta_k = -\pi + 2\pi k/N$, $k = 0,1,2, \ldots, N-1$; note that points $\theta_k$ constitute a regular partition of interval $[-\pi, \pi)$.

b) Approximate (in parallel) the function $\frac{d}{d\theta} t^*(\theta)$ by means of discrete difference quotients $\Delta t_j^*/\Delta\theta = \left[t^*(\theta_j) - t^*(\theta_{j-1})\right]/\Delta\theta$, with $j = 1,2, \ldots, N-1$; $\Delta\theta = \theta_j - \theta_{j-1} = 2\pi/N$. Since $\frac{d}{d\theta} t^*(\theta)$ is periodic with fundamental period $2\pi$ at interval $[-\pi, \pi)$, we also consider the difference quotient $\Delta t_0^*/\Delta\theta = [t^*(\theta_0) - t^*(\theta_{N-1})]/(2\pi/N)$; with this, we have the set of values $\{\Delta t_0^*/\Delta\theta, \Delta t_1^*/\Delta\theta, \Delta t_2^*/\Delta\theta, \ldots, \Delta t_{N-1}^*/\Delta\theta\}$, whose support is the discrete set $\{\theta_k - \Delta\theta/2\}$; thus, $\Delta t_k^*/\Delta\theta$ is the approximation of $\frac{d}{d\theta} t^*(\theta)$, evaluated at point $\theta = \theta_k - \Delta\theta/2$.

c) Join, by using rectilinear segments, adjacent points $(\theta_k - \Delta\theta/2, \Delta t_k^*/\Delta\theta)$, $(\theta_{k+1} - \Delta\theta/2, \Delta t_{k+1}^*/\Delta\theta) \in \mathbb{R}^2$ $\forall k$; in this way, we will have constructed, by linear interpolation, a **discrete approximation** $\frac{\Delta}{\Delta\theta} t^*(\theta)$ of function $\frac{d}{d\theta} t^*(\theta)$ at interval $-\pi \leq \theta < \pi$.

d) Find crossings of $\frac{\Delta}{\Delta\theta} t^*(\theta)$ with horizontal axis $y = 0$; with this, we will obtain *candidates* for estimates $\hat{\theta}_i^*$ of direction angles $\theta_i^*$ that, for each $i$, cause line $\ell_1(\theta_i^*)$ to contain one of the roots $R_i$ of equation (5.1).

e) By means of parallel constructions $LzC(\hat{\theta}_i^*)$, it is also possible to obtain estimates $\hat{R}_1$, $\hat{R}_2$, ..., $\hat{R}_n$ to the roots $R_1$, $R_2$, ..., $R_n$ of equation (5.1). For this, simply obtain points $R_x(\hat{\theta}_i^*) = P_1 + t^*(\hat{\theta}_i^*) v_{\hat{\theta}_i^*}$ as possible initial approximations to the roots. For each $i$, $t^*(\hat{\theta}_i^*)$ is the minimizing argument of the dynamic squared distance function $d^2_{\hat{\theta}_i^*}(t)$ defined in (5.8).

Any of these three strategies can be slightly modified in order to obtain a regional map at interval $[a, b)$; for this, we make the following general modifications:

- Step a) must be modified in such a way that $\theta_k = a + (b - a)k/N$, $k = 0,1,2, \dots, N - 1$.
- In step b) of strategies 5.2 and 5.3, we set $\Delta\theta = (b - a)/N$.
- In step b) of strategy 5.2 we do not compute $\Delta d_{\theta_0}^2/\Delta\theta$; similarly, in step b) of strategy 5.3 we do not compute $\Delta t_0^*/\Delta\theta$. In this way, from the support $\{\theta_k - \Delta\theta/2\}$ of functional values $\Delta d_{\theta_k}^2/\Delta\theta$ (strategy 5.2) or $\Delta t_k^*/\Delta\theta$ (strategy 5.3), we exclude $\theta_0 - \Delta\theta/2$.

Next, we will consider four numerical examples, using univariate polynomials of degree $n \geq 5$, in order to test the operation of these three strategies for the construction of proximity maps, which in theory allow obtaining estimates of the roots of equation (5.1). In the first example we will analyze the results for a univariate polynomial of degree $n = 7$ with randomly generated roots inside the usual region $\{x + iy : -1 \leq x \leq 1, -1 \leq y \leq 1\}$, while in the second example we will work with a univariate polynomial of degree $n = 10$ with coefficients of the form $C_k = k + ik$, in order to verify that the LC method also works for roots in regions of the complex plane $\mathbb{C}$ different from the usual region of the first example. In a third example we will analyze the case of a univariate polynomial of degree $n = 5$ with real roots, carrying out a translation of the polynomial variable prior to the application of the LC method, in order to avoid having all the roots on the real axis, which would surely cause problems in our root approximation strategy; thus, in this third example we will show that the LC method is also useful in the resolution of univariate polynomials with real roots, and is not limited exclusively to univariate polynomials with complex roots. Finally, in the fourth example we will work with a univariate polynomial of degree $n = 15$ with real coefficients of the form $C_k = k + 0i$; in this example we will see that it is not necessary to carry out a translation of the polynomial variable prior to the application of the LC method, since it turns out that this univariate polynomial only has one real root; the rest of the roots are complex values arranged in conjugate pairs.

## Numerical Example 5.1

This example follows the same dynamics as several of the examples from previous chapters: first we randomly generate roots $R_i$ which serve as reference values; then, we generate coefficients $C_i$ from Vieta's relations (5.2) between $R_i$ and $C_i$; finally, by using the constructed coefficients $C_i$ as inputs for the LC method, which is implemented by means of the script in annex 3 section 3, we generate proximity maps and initial estimates of roots $R_i$. As we already mentioned above, the functions used in the script from annex 3 section 3 take into account the generalizations seen in this chapter for structures $LzC(C_1, C_2, C_3, \dots, C_n, \theta)$; such generalizations work for univariate polynomials of degree greater than or equal to 4, and in fact they were already used when we constructed the examples in chapter 4.

So, we begin the execution of the script from annex 3 section 3 by inserting at the beginning of it the instruction `set.seed(2020)` for initializing the random number generator seed. This time, we seek to approximate the roots of an equation of degree $n = 7$ of the form

$$z^7 + C_1 z^6 + C_2 z^5 + C_3 z^4 + C_4 z^3 + C_5 z^2 + C_6 z + C_7 = 0. \tag{5.17}$$

The reference roots generated "randomly" by the script in annex 3 section 3 are:

$$\begin{aligned}
R_1 &= \phantom{-}0.2938057 - 0.2115485i & \theta_1^* &= \phantom{-}0.2485777\\
R_2 &= \phantom{-}0.2370036 - 0.0462177i & \theta_2^* &= \phantom{-}0.5005108\\
R_3 &= -0.7278056 - 0.8652312i & \theta_3^* &= -2.1994546\\
R_4 &= -0.7416948 - 0.2137641i & \theta_4^* &= \phantom{-}2.7111078\\
R_5 &= -0.9948346 + 0.2404119i & \theta_5^* &= \phantom{-}2.3534779\\
R_6 &= \phantom{-}0.5288280 + 0.4876715i & \theta_6^* &= \phantom{-}0.7652635\\
R_7 &= \phantom{-}0.6523314 - 0.1545418i & \theta_7^* &= \phantom{-}0.2172874
\end{aligned} \tag{5.18}$$

The coefficients of equation (5.17) generated from Vieta's relations (5.2) applied to reference roots $R_i$ in (5.18) are:

$$\begin{aligned}
C_1 &= \phantom{-}0.75236627 + 0.7632200i\\
C_2 &= -0.88646602 - 0.0588508i\\
C_3 &= -0.36898968 - 0.8967600i\\
C_4 &= \phantom{-}0.61107101 + 0.4171449i\\
C_5 &= -0.11086361 + 0.3804181i\\
C_6 &= -0.12915236 - 0.2105427i\\
C_7 &= \phantom{-}0.03085174 + 0.0215915i
\end{aligned} \tag{5.19}$$

With coefficients (5.19) used as inputs for the LC method, we generate discrete proximity maps built with $N = 2{,}500$ elements $LzC(C_1, C_2, C_3, \dots, C_7, \theta_k)$ associated to points $\theta_k$ in a regular partition of interval $[-\pi, \pi)$, following strategies 5.1, 5.2 and 5.3, and using the same operating conditions used in all examples from chapter 4 (for more details, see subsection "specific conditions for reproducing results from the examples in chapter 5" within annex 3 section 3). Figures 5.1, 5.2, and 5.3 show, respectively, weighted error map $e(\theta)$, map $\frac{d}{d\theta} d_\theta^2\big(t^*(\theta)\big)$, and map $\frac{d}{d\theta} t^*(\theta)$.

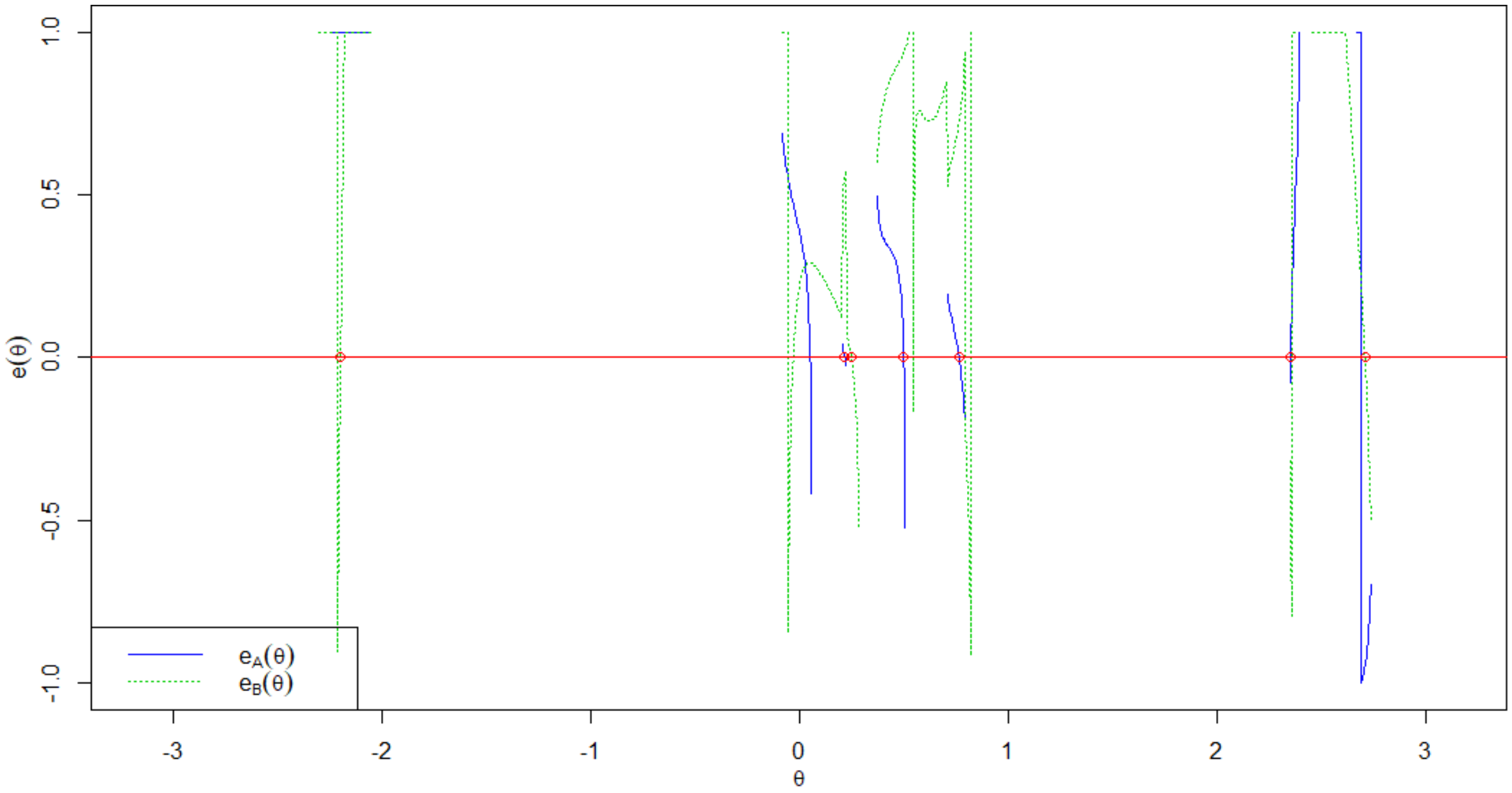

**Fig. 5.1.** Discrete angular proximity map $e(\theta)$ for the polynomial of degree $n = 7$ in equation (5.17) with coefficients given by (5.19). This map was generated from $N = 2{,}500$ elements $LzC(\theta_k)$ associated with points $\theta_k$ in a regular partition of interval $[-\pi, \pi)$. True theta roots $\theta_i^*$ are shown in this graph by means of small circles on horizontal axis $y = 0$.

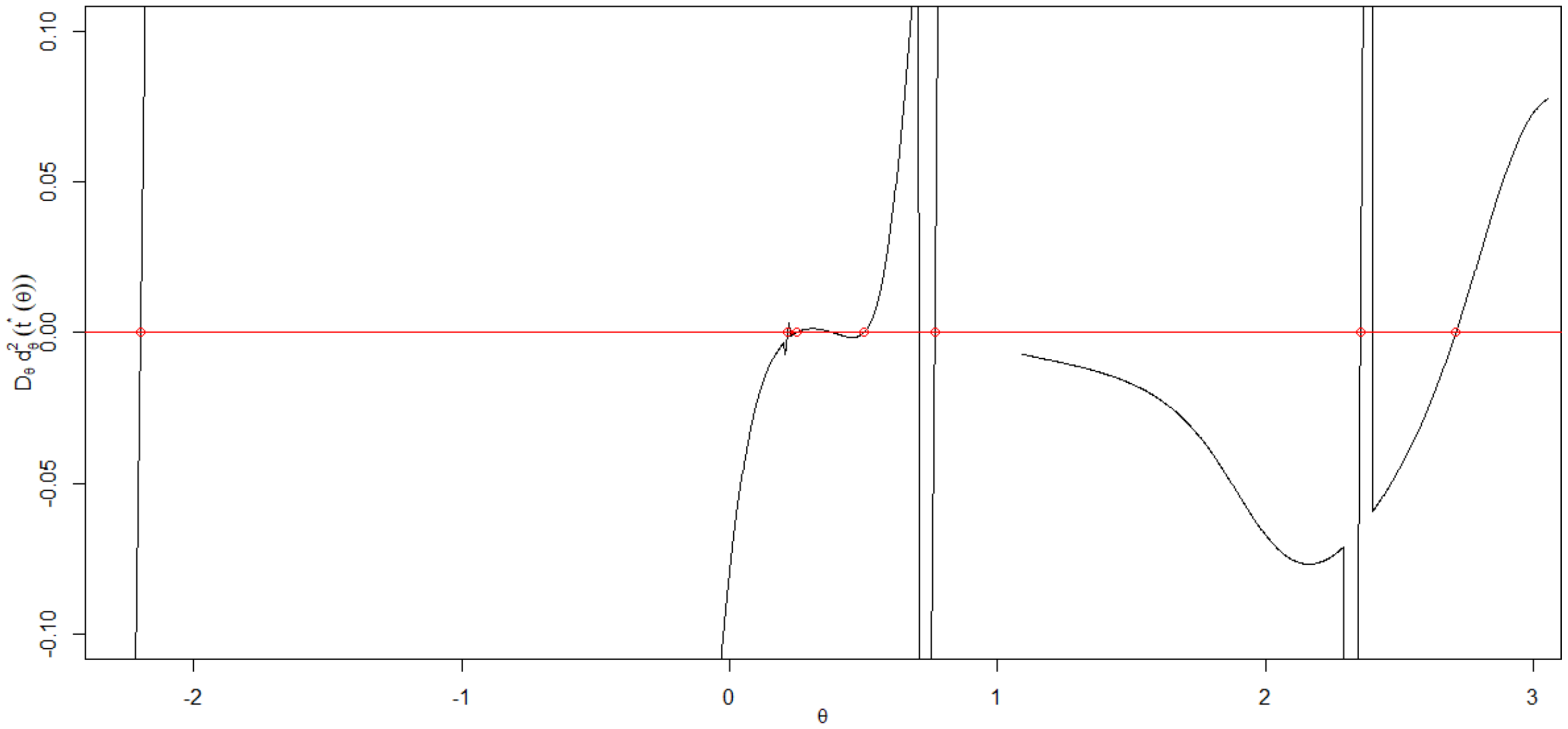

**Fig. 5.2.** Global map $\frac{d}{d\theta} d_\theta^2(t^*(\theta))$ for the polynomial of degree $n = 7$ in equation (5.17) with coefficients given by (5.19). This map was generated from $N = 2{,}500$ elements $LzC(\theta_k)$ associated with points $\theta_k$ in a regular partition of interval $[-\pi, \pi)$. True theta roots $\theta_i^*$ are shown in this graph by means of small circles on horizontal axis $y = 0$.

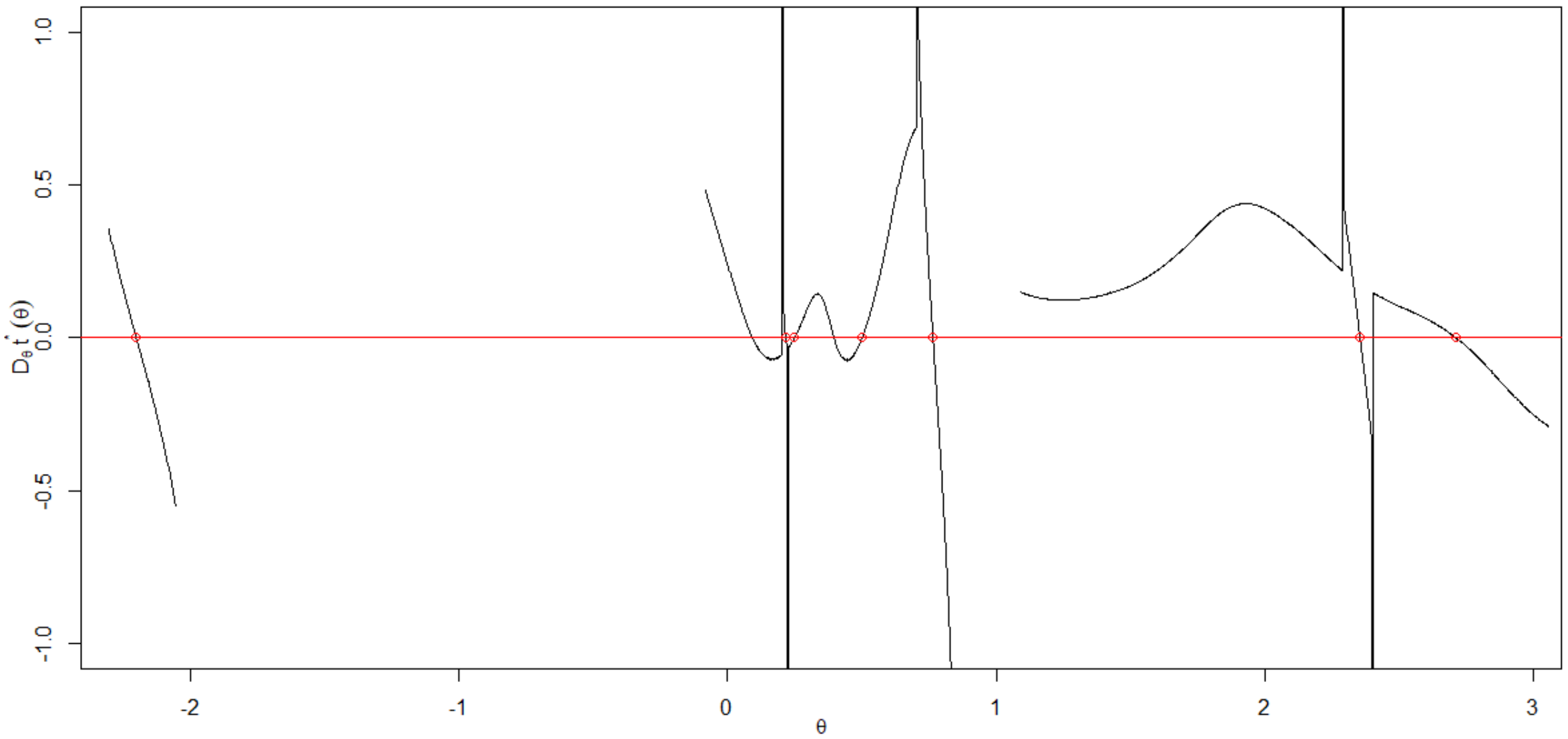

**Fig. 5.3.** Global map $\frac{d}{d\theta}t^*(\theta)$ for the polynomial of degree $n = 7$ in equation (5.17) with coefficients given by (5.19). This map was generated from $N = 2{,}500$ elements $LzC(\theta_k)$ associated with points $\theta_k$ in a regular partition of interval $[-\pi, \pi)$. True theta roots $\theta_i^*$ are shown in this graph by means of small circles on horizontal axis $y = 0$.

Tables 5.1, 5.2 and 5.3 show the numerical characteristics of the smooth crossings with the horizontal axis $y = 0$ detected in the respective proximity maps from figures 5.1, 5.2 and 5.3; such crossings of course determine the initial approximations to the roots of equation (5.17) with coefficients given by (5.19).

**Table 5.1**. Initial estimates for the roots of the polynomial of degree $n = 7$ in equation (5.17) with coefficients given by (5.19). These estimates were obtained from the map $e(\theta)$ (LC map) in figure 5.1.

| $i$ | $\hat{R}_i$ | $\hat{\theta}_i^*$ | $\lvert\Delta e_i\rvert$ | $d^2_{\hat{\theta}_i^*}\left(t^*(\hat{\theta}_i^*)\right)$ |
|---|---|---|---|---|
| 1 | $0.2938081 - 0.2115515i$ | 0.24857258 | 0.010403795 | $7.256933 \times 10^{-13}$ |
| 2 | $-0.7416927 - 0.2137618i$ | 2.71110045 | 0.033026375 | $8.006905 \times 10^{-12}$ |
| 3 | $0.6523252 - 0.1545383i$ | 0.21729197 | 0.009058721 | $1.583809 \times 10^{-11}$ |
| 4 | $-0.7278020 - 0.8652327i$ | −2.19944821 | 0.149406510 | $1.294250 \times 10^{-10}$ |
| 5 | $0.2370235 - 0.0462519i$ | 0.50045431 | 0.080021302 | $1.407863 \times 10^{-10}$ |
| 6 | $0.5288356 + 0.4876579i$ | 0.76525145 | 0.011271610 | $6.935222 \times 10^{-10}$ |
| 7 | $-0.9950221 + 0.2402194i$ | 2.35378411 | 0.141074334 | $4.915882 \times 10^{-7}$ |
| 8 | $0.2233646 - 0.0172098i$ | 0.54612909 | 0.571142076 | $1.509104 \times 10^{-4}$ |
| 9 | $0.3194888 - 0.3459742i$ | 0.05118024 | 0.077215562 | $2.877938 \times 10^{-3}$ |
| 10 | $0.3027097 - 0.4021197i$ | −0.03020135 | 0.035081695 | $8.773489 \times 10^{-3}$ |

**Table 5.2**. Initial estimates for the roots of the polynomial of degree $n = 7$ in equation (5.17) with coefficients given by (5.19). These estimates were obtained from the discrete map $\frac{d}{d\theta} d_{\theta}^{2}\big(t^{*}(\theta)\big)$ in figure 5.2.

| $i$ | $\hat{R}_i$ | $\hat{\theta}_i^*$ | $\left\lvert\Delta \frac{d}{d\theta} d_{\hat{\theta}_i^*}^2(t^*)\right\rvert$ | $d_{\hat{\theta}_i^*}^2\left(t^*(\hat{\theta}_i^*)\right)$ |
|---|---|---|---|---|
| 1 | $-0.7416937 - 0.2137640i$ | 2.7111064 | $6.800327 \times 10^{-4}$ | $9.533200 \times 10^{-13}$ |
| 2 | $0.2938045 - 0.2115340i$ | 0.2485984 | $1.175669 \times 10^{-4}$ | $1.037770 \times 10^{-11}$ |
| 3 | $0.6523279 - 0.1545365i$ | 0.2172931 | $1.763044 \times 10^{-3}$ | $1.279626 \times 10^{-11}$ |
| 4 | $0.2370154 - 0.0462388i$ | 0.5004762 | $2.189385 \times 10^{-4}$ | $5.277580 \times 10^{-11}$ |
| 5 | $-0.7278055 - 0.8652280i$ | −2.1994575 | $1.496352 \times 10^{-2}$ | $8.496743 \times 10^{-11}$ |
| 6 | $-0.9948317 + 0.2404093i$ | 2.3534776 | $2.632469 \times 10^{-2}$ | $1.063003 \times 10^{-10}$ |
| 7 | $0.5288290 + 0.4876648i$ | 0.7652591 | $2.214106 \times 10^{-2}$ | $1.310068 \times 10^{-10}$ |
| 8 | $0.2980537 - 0.2273455i$ | 0.2249271 | $3.419639 \times 10^{-3}$ | $1.532369 \times 10^{-5}$ |
| 9 | $0.2759660 - 0.1175893i$ | 0.3846779 | $7.145937 \times 10^{-5}$ | $1.180011 \times 10^{-4}$ |
| 10 | $0.2113400 + 0.1200644i$ | 0.7067425 | $1.944396 \times 10^{-1}$ | $9.889894 \times 10^{-3}$ |
| 11 | $-0.6549955 - 0.1266221i$ | 2.4007968 | $2.337855 \times 10^{-1}$ | $1.035961 \times 10^{-2}$ |

**Table 5.3**. Initial estimates for the roots of the polynomial of degree $n = 7$ in equation (5.17) with coefficients given by (5.19). These estimates were obtained from the discrete map $\frac{d}{d\theta} t^{*}(\theta)$ in figure 5.3.

| $i$ | $\hat{R}_i$ | $\hat{\theta}_i^*$ | $\left\lvert\Delta \frac{d}{d\theta} t^*(\hat{\theta}_i^*)\right\rvert$ | $d_{\hat{\theta}_i^*}^2\left(t^*(\hat{\theta}_i^*)\right)$ |
|---|---|---|---|---|
| 1 | $-0.7416936 - 0.2137636i$ | 2.71110541 | 0.001722451 | $1.447895 \times 10^{-12}$ |
| 2 | $0.5288217 + 0.4876725i$ | 0.76526754 | 0.033093921 | $1.185879 \times 10^{-10}$ |
| 3 | $0.2938388 - 0.2116691i$ | 0.24839675 | 0.003806301 | $7.703245 \times 10^{-10}$ |
| 4 | $0.2370517 - 0.0463055i$ | 0.50036765 | 0.004738137 | $9.024158 \times 10^{-10}$ |
| 5 | $-0.7278163 - 0.8652201i$ | −2.19947989 | 0.008146723 | $1.961897 \times 10^{-9}$ |
| 6 | $0.6523558 - 0.1546738i$ | 0.21716007 | 0.036193568 | $5.606980 \times 10^{-9}$ |
| 7 | $-0.9947727 + 0.2404679i$ | 2.35338284 | 0.018490947 | $4.744234 \times 10^{-8}$ |
| 8 | $0.2727815 - 0.1091656i$ | 0.39746986 | 0.008013320 | $1.156450 \times 10^{-4}$ |
| 9 | $0.3192382 - 0.3162669i$ | 0.09368678 | 0.002271653 | $1.402532 \times 10^{-3}$ |

From tables 5.1, 5.2 and 5.3, we see that the first 7 rows in each of these tables correspond with valid approximations to reference roots (5.18), thanks to the fact that rows are arranged in ascending order according to their corresponding quality indicators $d_{\hat{\theta}_i^*}^2\left(t^*(\hat{\theta}_i^*)\right)$.

Next, we compute the absolute relative errors for the first seven approximations shown in table 5.1, corresponding to map $e(\theta)$, with respect to reference values (5.18):

$$
\begin{aligned}
&\left|\frac{R_1-\hat{R}_1}{R_1}\right| = 1.060988\times10^{-5} && \left|\frac{\theta_1^*-\hat{\theta}_1^*}{\theta_1^*}\right| = 2.041465\times10^{-5}\\
&\left|\frac{R_2-\hat{R}_5}{R_2}\right| = 1.635820\times10^{-4} && \left|\frac{\theta_2^*-\hat{\theta}_5^*}{\theta_2^*}\right| = 1.128774\times10^{-4}\\
&\left|\frac{R_3-\hat{R}_4}{R_3}\right| = 3.488718\times10^{-6} && \left|\frac{\theta_3^*-\hat{\theta}_4^*}{\theta_3^*}\right| = 2.905686\times10^{-6}\\
&\left|\frac{R_4-\hat{R}_2}{R_4}\right| = 4.011305\times10^{-6} && \left|\frac{\theta_4^*-\hat{\theta}_2^*}{\theta_4^*}\right| = 2.719205\times10^{-6}\\
&\left|\frac{R_5-\hat{R}_7}{R_5}\right| = 2.625238\times10^{-4} && \left|\frac{\theta_5^*-\hat{\theta}_7^*}{\theta_5^*}\right| = 1.301179\times10^{-4}\\
&\left|\frac{R_6-\hat{R}_6}{R_6}\right| = 2.169096\times10^{-5} && \left|\frac{\theta_6^*-\hat{\theta}_6^*}{\theta_6^*}\right| = 1.571895\times10^{-5}\\
&\left|\frac{R_7-\hat{R}_3}{R_7}\right| = 1.064499\times10^{-5} && \left|\frac{\theta_7^*-\hat{\theta}_3^*}{\theta_7^*}\right| = 2.089517\times10^{-5}
\end{aligned}
\tag{5.20}
$$

Similarly, the absolute relative errors associated with the first seven rows of table 5.2, corresponding to map $\frac{d}{d\theta}d_\theta^2\big(t^*(\theta)\big)$, are:

$$
\begin{aligned}
&\left|\frac{R_1-\hat{R}_2}{R_1}\right| = 4.012572\times10^{-5} && \left|\frac{\theta_1^*-\hat{\theta}_2^*}{\theta_1^*}\right| = 8.333080\times10^{-5}\\
&\left|\frac{R_2-\hat{R}_4}{R_2}\right| = 1.001421\times10^{-4} && \left|\frac{\theta_2^*-\hat{\theta}_4^*}{\theta_2^*}\right| = 6.912211\times10^{-5}\\
&\left|\frac{R_3-\hat{R}_5}{R_3}\right| = 2.826752\times10^{-6} && \left|\frac{\theta_3^*-\hat{\theta}_5^*}{\theta_3^*}\right| = 1.321283\times10^{-6}\\
&\left|\frac{R_4-\hat{R}_1}{R_4}\right| = 1.384114\times10^{-6} && \left|\frac{\theta_4^*-\hat{\theta}_1^*}{\theta_4^*}\right| = 5.404571\times10^{-7}\\
&\left|\frac{R_5-\hat{R}_6}{R_5}\right| = 3.860511\times10^{-6} && \left|\frac{\theta_5^*-\hat{\theta}_6^*}{\theta_5^*}\right| = 1.114570\times10^{-7}\\
&\left|\frac{R_6-\hat{R}_7}{R_6}\right| = 9.427304\times10^{-6} && \left|\frac{\theta_6^*-\hat{\theta}_7^*}{\theta_6^*}\right| = 5.747118\times10^{-6}\\
&\left|\frac{R_7-\hat{R}_3}{R_7}\right| = 9.568190\times10^{-6} && \left|\frac{\theta_7^*-\hat{\theta}_3^*}{\theta_7^*}\right| = 2.618673\times10^{-5}
\end{aligned}
\tag{5.21}
$$

Finally, the absolute relative errors associated with the first seven rows of table 5.3, corresponding to map $\frac{d}{d\theta}t^*(\theta)$, are:

$$
\begin{aligned}
&\left|\frac{R_1-\hat{R}_3}{R_1}\right| = 3.454741\times 10^{-4} && \left|\frac{\theta_1^*-\hat{\theta}_3^*}{\theta_1^*}\right| = 7.277809\times 10^{-4}\\
&\left|\frac{R_2-\hat{R}_4}{R_2}\right| = 4.143581\times 10^{-4} && \left|\frac{\theta_2^*-\hat{\theta}_4^*}{\theta_2^*}\right| = 2.860196\times 10^{-4}\\
&\left|\frac{R_3-\hat{R}_5}{R_3}\right| = 1.358327\times 10^{-5} && \left|\frac{\theta_3^*-\hat{\theta}_5^*}{\theta_3^*}\right| = 1.149445\times 10^{-5}\\
&\left|\frac{R_4-\hat{R}_1}{R_4}\right| = 1.705773\times 10^{-6} && \left|\frac{\theta_4^*-\hat{\theta}_1^*}{\theta_4^*}\right| = 8.884381\times 10^{-7} \qquad (5.22)\\
&\left|\frac{R_5-\hat{R}_7}{R_5}\right| = 8.155748\times 10^{-5} && \left|\frac{\theta_5^*-\hat{\theta}_7^*}{\theta_5^*}\right| = 4.038258\times 10^{-5}\\
&\left|\frac{R_6-\hat{R}_2}{R_6}\right| = 8.969165\times 10^{-6} && \left|\frac{\theta_6^*-\hat{\theta}_2^*}{\theta_6^*}\right| = 5.311865\times 10^{-6}\\
&\left|\frac{R_7-\hat{R}_6}{R_7}\right| = 2.002243\times 10^{-4} && \left|\frac{\theta_7^*-\hat{\theta}_6^*}{\theta_7^*}\right| = 5.861431\times 10^{-4}
\end{aligned}
$$

In table 5.4 we can see some statistics for absolute relative errors (5.20), (5.21) and (5.22), as well as for quality indicators $d^2_{\hat{\theta}_i^*}\left(t^*\left(\hat{\theta}_i^*\right)\right)$ associated with the first seven rows in each of tables 5.1, 5.2, and 5.3. From here, we can observe that, in this example, once again the map $\frac{d}{d\theta}d^2_\theta\left(t^*(\theta)\right)$ is the one who generates the most accurate and consistent approximations; in second place, is the map $e(\theta)$, although its quality indicator $d^2_{\hat{\theta}_i^*}\left(t^*\left(\hat{\theta}_i^*\right)\right)$ has a slightly lower consistency compared with its counterpart from map $\frac{d}{d\theta}t^*(\theta)$; this is due to the value $d^2_{\hat{\theta}_i^*}\left(t^*\left(\hat{\theta}_i^*\right)\right)$ in row 7 of table 5.1, which is unusually high compared to the rest of values $d^2_{\hat{\theta}_i^*}\left(t^*\left(\hat{\theta}_i^*\right)\right)$ associated with valid approximations in table 5.1. In the end, however, we can say that the three maps considered in this example produce reasonable initial approximations to the roots of equation (5.17) with coefficients given by (5.19), if we compare the results obtained here with the results obtained in the numerical examples from previous chapters.

**Table 5.4**. Arithmetic means and standard deviations for error measures associated with valid initial approximations (first seven rows in each of tables 5.1, 5.2 and 5.3) to the roots of the polynomial of degree $n = 7$ in equation (5.17) with coefficients given by (5.19), obtained by means of discrete proximity maps $e(\theta)$, $\frac{d}{d\theta} d_\theta^2\big(t^*(\theta)\big)$, and $\frac{d}{d\theta} t^*(\theta)$ from figures 5.1, 5.2 and 5.3.

| Map | Error measure | | | | | |
|---|---|---|---|---|---|---|
| | $d_{\hat{\theta}_i^*}^2\left(t^*(\hat{\theta}_i^*)\right)$ | | $\left\lvert\frac{R-\hat{R}}{R}\right\rvert$ | | $\left\lvert\frac{\theta^*-\hat{\theta}^*}{\theta^*}\right\rvert$ | |
| | **Mean** | **Standard deviation** | **Mean** | **Standard deviation** | **Mean** | **Standard deviation** |
| $e(\boldsymbol{\theta})$ | $7.036807 \times 10^{-8}$ | $1.857408 \times 10^{-7}$ | $6.807882 \times 10^{-5}$ | $1.032468 \times 10^{-4}$ | $4.366414 \times 10^{-5}$ | $5.391717 \times 10^{-5}$ |
| $\frac{d}{d\boldsymbol{\theta}} d_{\boldsymbol{\theta}}^2(t^*(\boldsymbol{\theta}))$ | $5.702537 \times 10^{-11}$ | $5.160165 \times 10^{-11}$ | $2.390495 \times 10^{-5}$ | $3.615617 \times 10^{-5}$ | $2.662285 \times 10^{-5}$ | $3.531119 \times 10^{-5}$ |
| $\frac{d}{d\boldsymbol{\theta}} t^*(\boldsymbol{\theta})$ | $8.114856 \times 10^{-9}$ | $1.744768 \times 10^{-8}$ | $1.522675 \times 10^{-4}$ | $1.711513 \times 10^{-4}$ | $2.368601 \times 10^{-4}$ | $3.066238 \times 10^{-4}$ |

In figures 5.4, 5.5 and 5.6 we can see graphs related to the global minima of dynamic squared distances between terminal curves and z-circumferences associated to equation (5.17) with coefficients (5.19). From here, we can verify that once again we have evidence in favor of the hypotheses on the properties for mappings $d_\theta^2\big(t^*(\theta)\big)$ and $t^*(\theta)$, posed in example 4.1 from chapter 4.

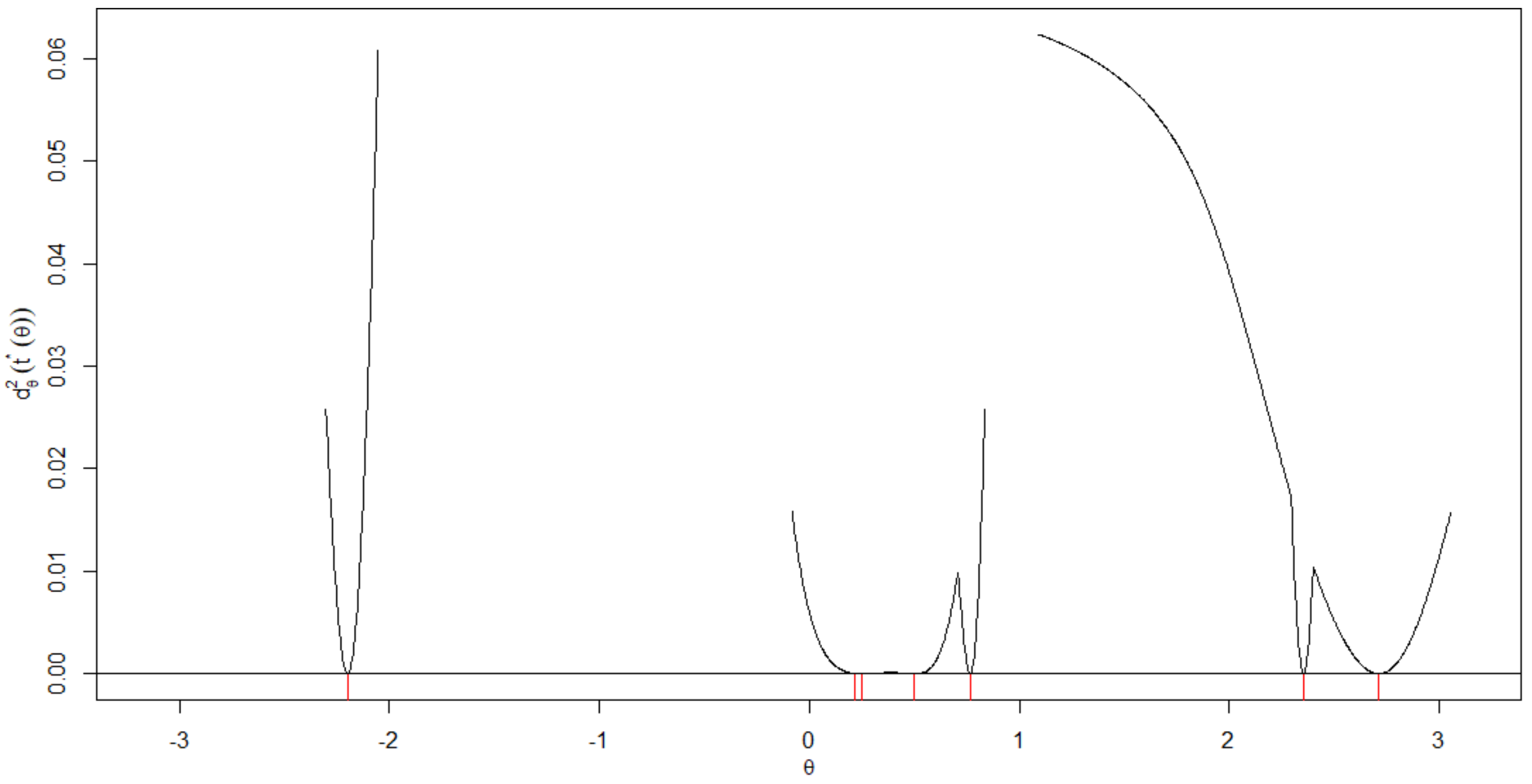

**Fig. 5.4.** Global minima $d_\theta^2\big(t^*(\theta)\big)$ vs. $\theta$ for the polynomial of degree $n = 7$ in equation (5.17) with coefficients given by (5.19). This map was generated from $N = 2{,}500$ elements $LzC(\theta_k)$ associated with points $\theta_k$ in a regular partition of interval $[-\pi, \pi)$. This graph includes vertical line segments that indicate the location of true theta roots $\theta_i^*$.

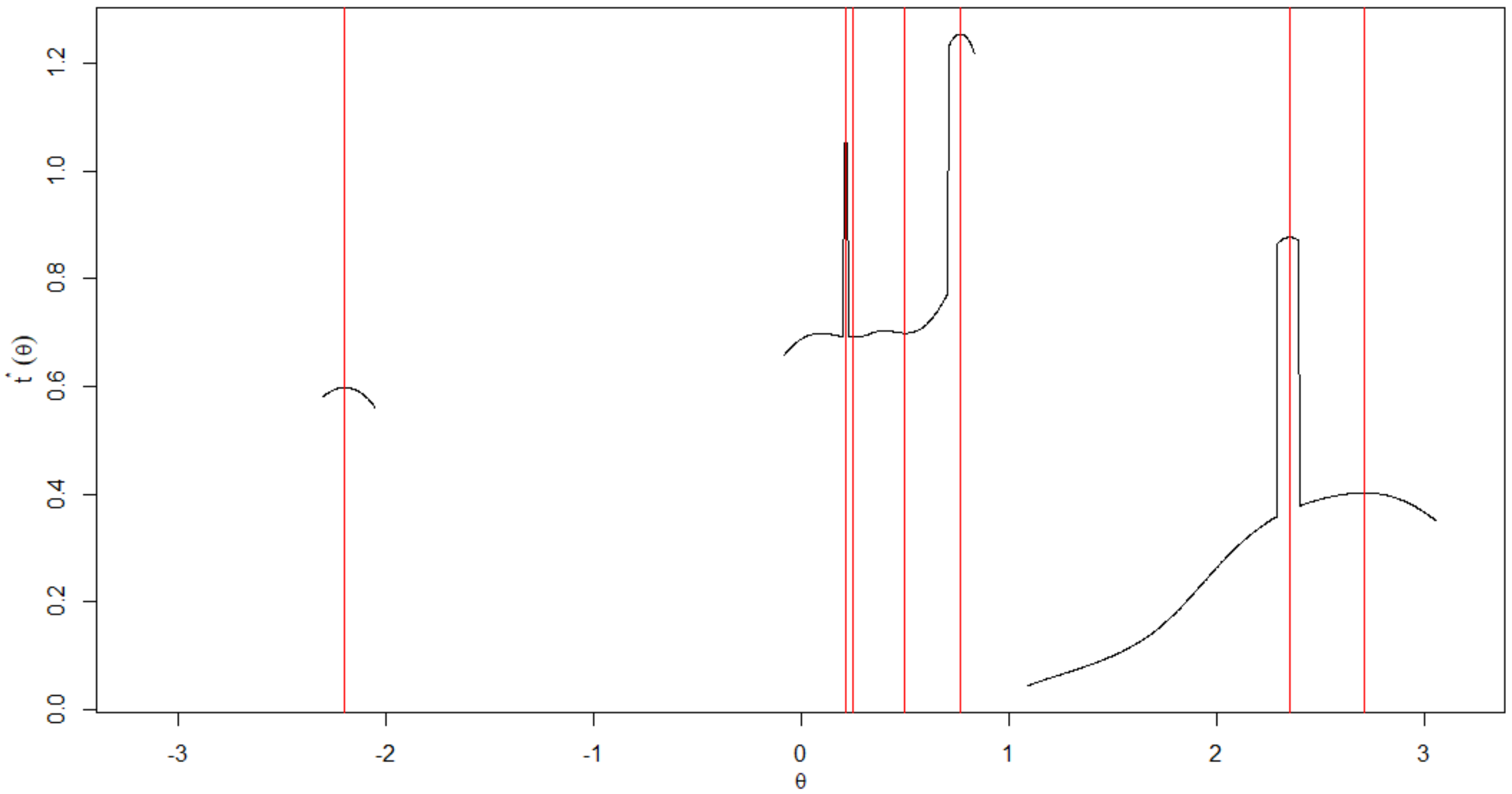

**Fig. 5.5.** Minimizing arguments $t^*(\theta)$ vs. $\theta$ for the polynomial of degree $n = 7$ in equation (5.17) with coefficients given by (5.19). This map was generated from $N = 2{,}500$ elements $LzC(\theta_k)$ associated with points $\theta_k$ in a regular partition of interval $[-\pi, \pi)$. This graph includes vertical line segments that indicate the location of true theta roots $\theta_i^*$.

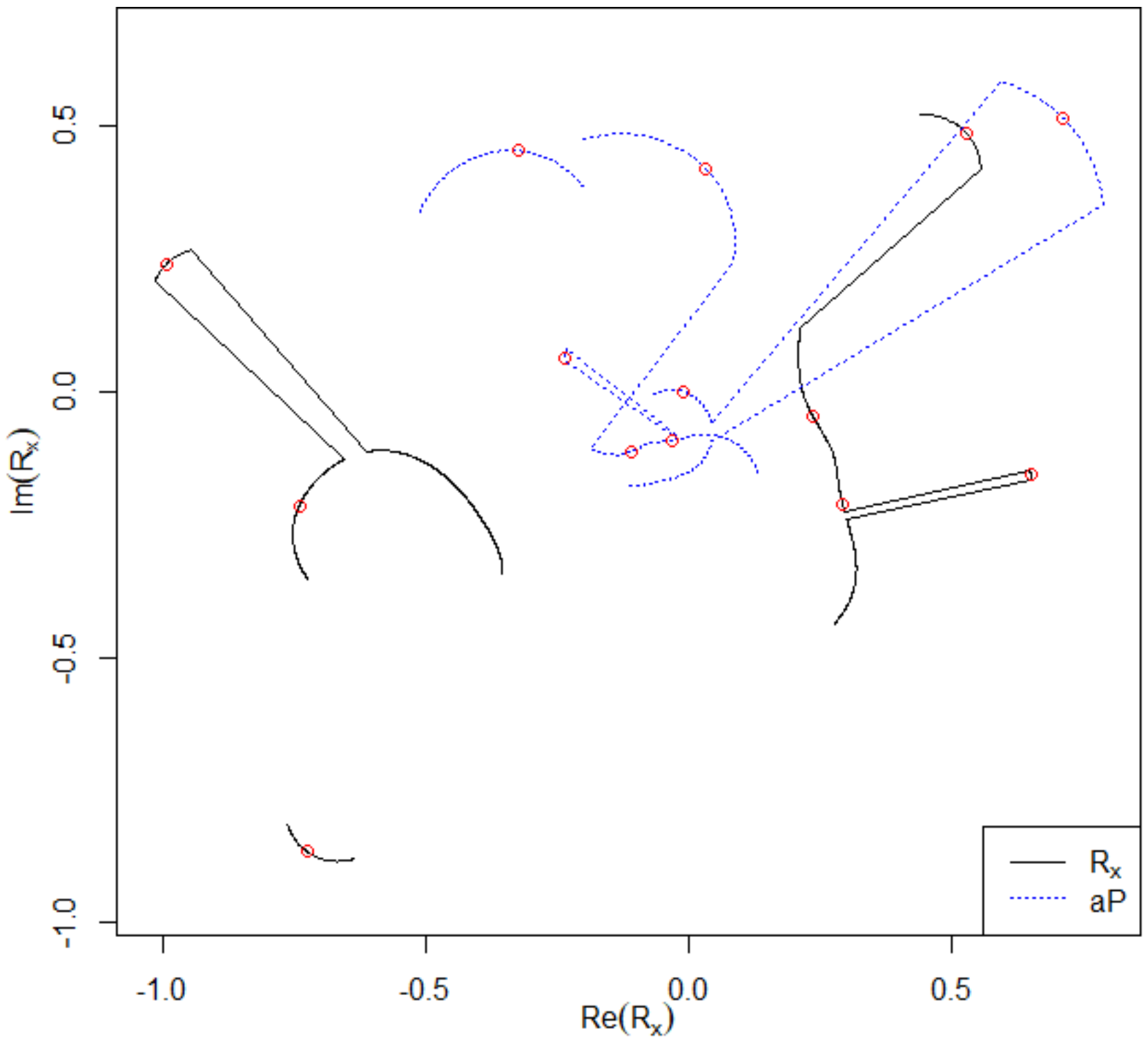

**Fig. 5.6.** Trajectories $R_x(\theta)$ and $aP(\theta)$ associated with the polynomial of degree $n = 7$ in equation (5.17) with coefficients given by (5.19). These trajectories were generated from $N = 2{,}500$ elements $LzC(\theta_k)$ associated with points $\theta_k$ in a regular partition of interval $[-\pi, \pi)$. The true roots $R_i(\theta_i^*)$, and their corresponding anchor points $aP(\theta_i^*)$, are represented in this figure by means of small circles on the trajectories.

Figures 5.2 and 5.4 suggest that theta roots $\theta_i^*$ occur at the global minima of function $d_\theta^2\big(t^*(\theta)\big)$, while figures 5.3 and 5.5 suggest that five of the theta roots correspond to local maxima of function $t^*(\theta)$, and the other two theta roots correspond to local minima of $t^*(\theta)$.

Figures 5.2, 5.3, 5.4 and 5.5 each show 3 continuous sections with support $\theta$ of total length $\pi$; in any of these figures, one of these sections contains 4 theta roots, while another contains 2 theta roots, and the third one contains one theta root. Additionally, each of these 3 sections projects its own characteristics onto the trajectories $aP(\theta)$ and $R_x(\theta)$ shown in figure 5.6.

Note from figure 5.6 that trajectory $aP(\theta)$ is no longer similar to trajectory $R_x(\theta)$, as was the case in all numerical examples seen for univariate polynomials of degree $n = 4$ (see figures 4.8, 4.14, 4.30 and 4.38 in chapter 4). This is because, for the case of $n = 4$,

$$aP(\theta) = -C_3 - R_x(\theta)(C_2 - P_1^2);$$

i.e., in this case, $aP(\theta)$ is clearly an expression that depends linearly on $R_x(\theta)$. On the other hand, according to the expressions in (5.13), for $n = 7$,

$$aP(\theta) = C_6 + C_5R_x(\theta) + C_4R_x(\theta)^2 + C_3R_x(\theta)^3 + R_x(\theta)^4(C_2 - P_1^2); \quad (5.23)$$

we see from expression (5.23) that $aP(\theta)$ depends on factors of $R_x(\theta)$ raised to different powers, additional to power 1, so $aP(\theta)$, in this case, depends on $R_x(\theta)$, but not linearly.

Before concluding this example, let us look at the behavior of a structure $LzC(\theta)$ associated to equation (5.17) of degree $n = 7$ with coefficients given by (5.19), angularly close to reference root $R_6$ in (5.18). The graphs in figures 5.7 and 5.8, which allow us to analyze the situation of this structure $LzC(\theta)$, were generated with the help of the script listed in annex 4 section 2; this script can also be used to interactively analyze, in a general way (although with some necessary modifications to its code), structures $LzC(\theta)$ associated with univariate polynomials of degree $n \geq 4$.

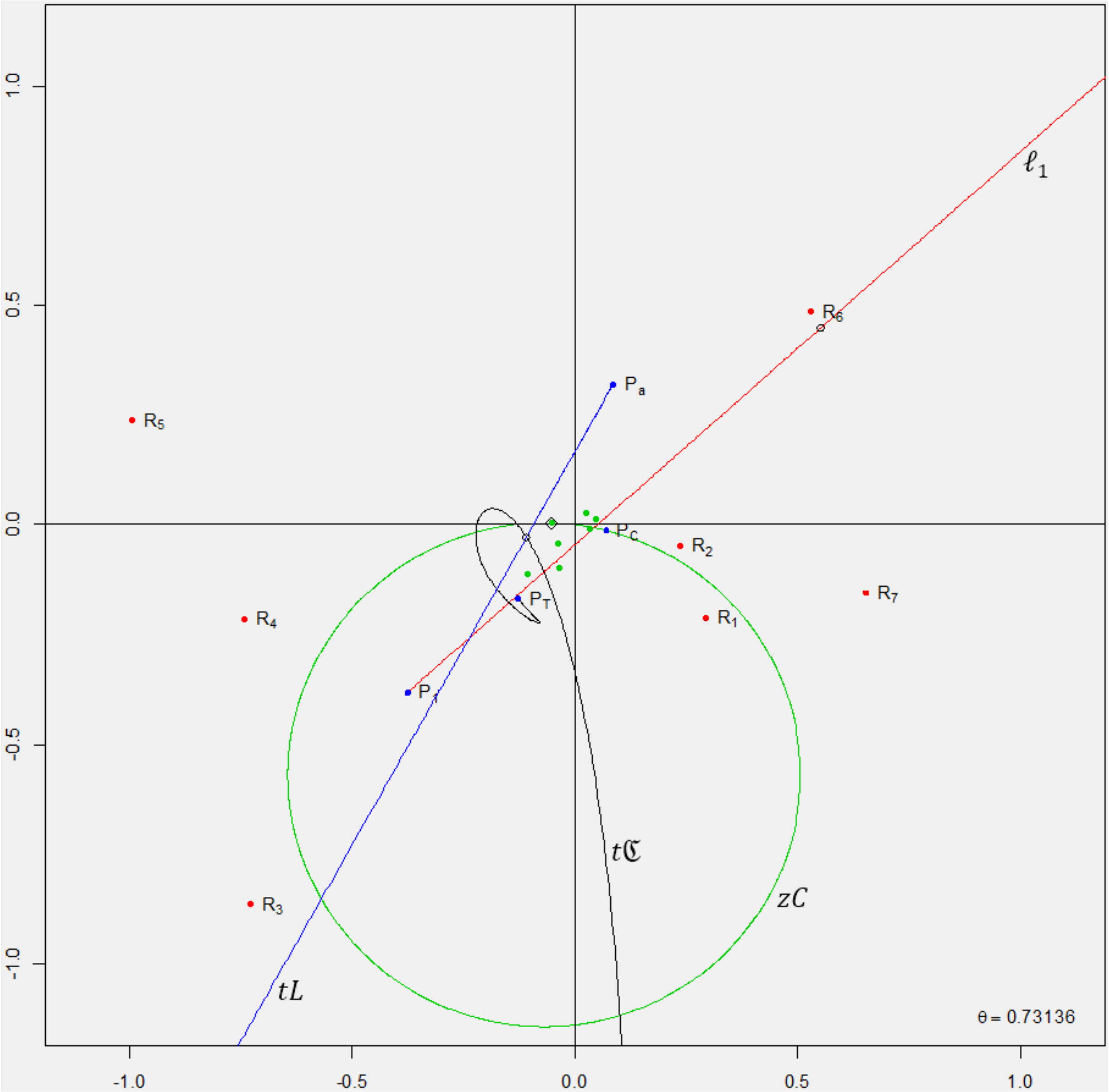

**Fig. 5.7.** Some elements of the structure $LzC(\theta)$ associated to equation (5.17) of degree $n = 7$ with coefficients given by (5.19), for $\theta = 0.73136$. Small circles on line $\ell_1$, terminal curve $t\mathfrak{C}$, terminal semi-line $tL$, as well as the small diamond ◊ on z-circumference $zC$, correspond to the minimizing argument $t^*$ of associated function $d_\theta^2$ (see figure 5.8). On the other hand, points $P_a$ in $tL$, $P_T$ in $t\mathfrak{C}$, $P_1$ in $\ell_1$ and $P_c$ in $zC$ correspond to $t = 0$ in these parametric trajectories. The position of the anchor point of $tL$, $P_a$, depends on $\theta$ (see expression (5.23)), but not the positions of points $P_1 = -C_1/2$, $P_c = -C_7/P_1$, $P_T = C_6 + C_5P_1 + C_4P_1^2 + C_3P_1^3 + P_1^4(C_2 - P_1^2)$, which are true fixed points.

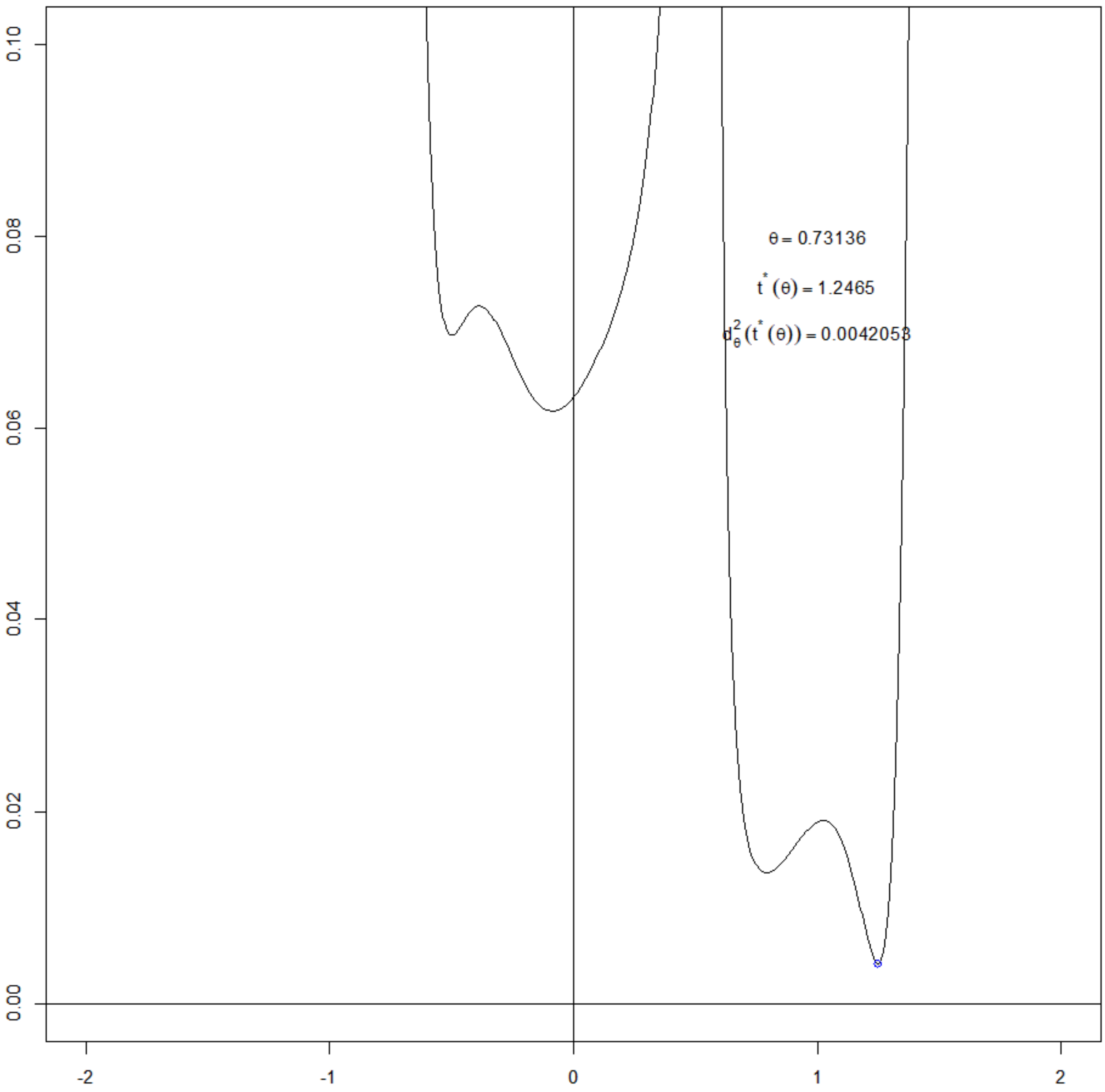

**Fig. 5.8.** Close-up of the graph for the dynamic squared distance function $d_{\theta}^{2}(t)$ associated to equation (5.17) of degree $n=7$ with coefficients given by (5.19), where $\theta=0.73136$. The global minimum of this function occurs at minimizing argument $t^{*}=1.2465$. At $t^{*}$, the minimum functional value is $d_{\theta}^{2}(t^{*})=4.2053\times10^{-3}$.

From figure 5.7, it appears that trajectories $t\mathfrak{C}$ and $tL$ intersect each other when their common parameter $t$ coincides with the minimizing argument $t^{*}$ of the function $d_{\theta}^{2}(t)$ shown in figure 5.8. This can be verified algebraically; by evaluating expression (5.6) for the terminal curve $t\mathfrak{C}$ at $t=t^{*}$, we have

$$t\mathfrak{C}(\theta,t^*) = (-1)^{n+1}\left[\sum_{i=1}^{n-3} C_{n-i}\,(P_1+t^*v_\theta)^{i-1} + (P_1+t^*v_\theta)^{n-3}(C_2-\ell_d(\theta,t^*))\right] \quad (5.24)$$

Let's remember that $R_x(\theta) = P_1 + t^*(\theta)v_\theta$ (see expression (5.9)); substituting this into (5.24), we have

$$t\mathfrak{C}(\theta,t^*) = (-1)^{n+1}\left[\sum_{i=1}^{n-3} C_{n-i}\,R_x(\theta)^{i-1} + R_x(\theta)^{n-3}(C_2-\ell_d(\theta,t^*))\right]$$

$$t\mathfrak{C}(\theta,t^*) = (-1)^{n+1}\left[\sum_{i=1}^{n-3} C_{n-i}\,R_x(\theta)^{i-1} + R_x(\theta)^{n-3}(C_2-(P_1^2-t^{*2}v_\theta^2))\right] \quad (5.25)$$

On the other hand, if we evaluate expression (5.13) for the terminal semi-line at $t = t^*$, we have

$$tL(\theta,t^*) = (-1)^{n+1}\left[\sum_{i=1}^{n-3} C_{n-i}\,R_x(\theta)^{i-1} + R_x(\theta)^{n-3}(C_2-P_1^2)\right] + t^{*2}[(-1)^{n+1}R_x(\theta)^{n-3}v_\theta^2]$$

$$tL(\theta,t^*) = (-1)^{n+1}\left[\sum_{i=1}^{n-3} C_{n-i}\,R_x(\theta)^{i-1} + R_x(\theta)^{n-3}(C_2-P_1^2) + t^{*2}R_x(\theta)^{n-3}v_\theta^2\right]$$

$$tL(\theta,t^*) = (-1)^{n+1}\left[\sum_{i=1}^{n-3} C_{n-i}\,R_x(\theta)^{i-1} + R_x(\theta)^{n-3}((C_2-P_1^2)+t^{*2}v_\theta^2)\right]$$

$$tL(\theta,t^*) = (-1)^{n+1}\left[\sum_{i=1}^{n-3} C_{n-i}\,R_x(\theta)^{i-1} + R_x(\theta)^{n-3}(C_2-(P_1^2-t^{*2}v_\theta^2))\right] \quad (5.26)$$

From (5.25) and (5.26) we see that, in effect, $t\mathfrak{C}(\theta,t^*) = tL(\theta,t^*)$.

From figure 5.7, we can also see that, when line $\ell_1$ with fixed point $P_1$ approaches root $R_6$, points $tL(\theta,t^*)$, $t\mathfrak{C}(\theta,t^*)$ and $zC(\theta,t^*)$ approach product $R_1R_2R_3R_4R_5R_7$, which appears in figure 5.7 as a point inside the diamond ◊ that marks the location of $zC(\theta,t^*)$; this tells us that point $zC(\theta,t^*)$ is at moment $\theta = 0.73136$ much closer to product $R_1R_2R_3R_4R_5R_7$ than points $tL(\theta,t^*)$ and $t\mathfrak{C}(\theta,t^*)$. In fact, expressions (5.25) and (5.26) tell us that trajectories $tL$ and $t\mathfrak{C}$ have an identical approaching behavior towards product $R_1R_2R_3R_4R_5R_7$, with respect to value $t = t^*$, at least at any moment $\theta$ close to $\theta_6^* = 0.7652635$. Note also that one of the static intersections between $tL(\theta = 0.73136)$ and $zC(\theta = 0.73136)$, which determines weighted error value $e_A(\theta = 0.73136)$, is also close to the product $R_1R_2R_3R_4R_5R_7$, although its approaching behavior is different from that of points $tL(\theta,t^*)$, $t\mathfrak{C}(\theta,t^*)$ and $zC(\theta,t^*)$.

From figure 5.8, we can see that a dynamic squared distance function $d_\theta^2(t)$ associated with a univariate polynomial of degree $n = 7$ can contain several local minima and maxima, which means that finding the minimizing argument $t^*$ required for the construction of the complete structure $LzC(\theta)$ becomes more difficult as the degree of the polynomial grows.

## Numerical Example 5.2

This time we will approximate, by using the LC method, the roots of equation

$$z^{10} + (1+i)z^9 + (2+2i)z^8 + (3+3i)z^7 + (4+4i)z^6 + (5+5i)z^5 + (6+6i)z^4 + (7+7i)z^3 + (8+8i)z^2 + (9+9i)z + (10+10i) = 0 \qquad (5.27)$$

The procedure we will follow here is similar to that of example 4.3 in chapter 4; first, we will use the R instruction `polyroot` in order to obtain reference roots $R_i$, and then we will reconstruct the coefficients $C_k$ indicated in (5.1) from reference roots $R_i$, by means of Vieta's relations (5.2); we will use these reconstructed coefficients $C_k$ as inputs for the LC method in order to generate discrete proximity maps and corresponding initial estimates $\hat{R}_i$ for the roots of equation (5.27). Finally, initial estimates $\hat{R}_i$ will be compared with reference roots $R_i$, in order to obtain absolute relative errors.

Note: In theory, the reconstructed coefficients $C_k$ in this example should be exactly of the form $C_k = k + ik$, but in reality they differ from their exact values by extremely small quantities in their real and imaginary parts, of the order of $10^{-14}$ or less; this is due to errors induced by floating-point operations carried out when applying Vieta's formulas (5.2) on reference roots $R_i$, which in turn, are approximations (found by the R function `polyroot`) to the true roots of equation (5.27). In this numerical example, and, to tell the truth, in all the examples that we have considered in this work, we proceed under two implicit assumptions: 1) the floating-point-arithmetic-induced errors contained in the reconstructed coefficients $C_k$ are extremely small, and 2) initial estimates $\hat{R}_i$ produced by the LC method are robust to the extremely small errors contained in reconstructed coefficients $C_k$. The validity of the second assumption may be challenged by some types of polynomials, particularly if they have multiple roots; this issue is related to the case of Wilkinson's polynomial $p(x) = \prod_{i=1}^{20}(x-i)$, a classic example of a single-variable polynomial with real coefficients in which a small perturbation in one of its coefficients has a large effect on its roots. See [5.1] for more details.

Returning to our example, we will change some operation parameters of the script listed in annex 3 section 3, with the purpose of reducing (although not completely eliminating) the problem of gaps in maps $\frac{d}{d\theta} d_\theta^2\big(t^*(\theta)\big)$ and $\frac{d}{d\theta} t^*(\theta)$, due to errors in the simulated annealing phase within the optimization process that approximates values $d_\theta^2\big(t^*(\theta)\big)$ and $t^*(\theta)$; this problem was previously mentioned in the final part of example 4.2 in chapter 4 (see figure 4.23). If we used the same operation parameters of example 5.1 in our attempt to find approximations to the roots of equation (5.27), we would see that the gap problem in maps $\frac{d}{d\theta} d_\theta^2\big(t^*(\theta)\big)$ and $\frac{d}{d\theta} t^*(\theta)$ is so severe that it prevents obtaining all expected reasonable approximations; only 6 or 7 reasonable approximations would be obtained, instead of the 10 expected. What happens in this case is that the temperature with which the simulated annealing algorithm starts out is not high enough, which increases the probability that the algorithm gets stuck in a local minimum, instead of correctly approximating the required global minimum of function $d_\theta^2(t)$. We will therefore raise the initial temperature of the simulated annealing algorithm, hoping that the probability of finding the global minimum of $d_\theta^2(t)$ increases.

Let us talk about the specific details on how we will produce the results in this example. Unlike example 5.1 in this chapter, where we invoked the function `approxLC` in the script of annex 3 section 3 with the argument `TEMP=10`, this time we will do it with the argument `TEMP=500`; with this, we considerably increase the initial temperature for the simulated annealing algorithm, giving it more opportunity to explore the search space in an exhaustive way, with which we seek to increase the probability of getting close enough to the global minimum, without getting stuck in local minima. In the two calls to the function `approxDMin` used to generate the maps $\frac{d}{d\theta}d_\theta^2\big(t^*(\theta)\big)$ and $\frac{d}{d\theta}t^*(\theta)$, we will also change the argument `TEMP`, from `12` (used in example 5.1 of this chapter) to `500`. At the beginning of the script from annex 3 section 3, in order to initialize the random number generator seed, we will use the instruction `set.seed(2022)`; with this, the script will always produce the same results from the "stochastic" optimization processes.

Another adjustment that we will make in the script from annex 3 section 3 in order to detect "smooth" crossings within the map $\frac{d}{d\theta}d_\theta^2\big(t^*(\theta)\big)$, consists in changing argument `Tol`, in the corresponding call to function `approxDMin` (listed in annex 2 section 7), from `2.0` (used in example 5.1 of this chapter) to `1000.0`. It is necessary to do this adjustment because of the sensitivity of functional values $\frac{d}{d\theta}d_\theta^2\big(t^*(\theta)\big)$ to the norms of the coefficients in equation (5.27); this issue was briefly mentioned in example 4.3, chapter 4. Later in this example we will comment more about it.

For more details on how to reproduce the results of this example, see the specific conditions for the examples in chapter 5, located in annex 3 section 3, after the script listing.

The reference roots $R_i$ generated by the R function `polyroot`, together with their corresponding theta roots $\theta_i^*$, are:

$$\begin{aligned}
R_1 &= \phantom{-}0.9472406 + 0.7629888i & \theta_1^* &= \phantom{-}0.7175188\\
R_2 &= -0.8833566 + 0.8776421i & \theta_2^* &= \phantom{-}1.8422002\\
R_3 &= -0.6774270 - 1.0722056i & \theta_3^* &= -1.8714710\\
R_4 &= \phantom{-}0.9645605 - 0.8405533i & \theta_4^* &= -0.2284693\\
R_5 &= \phantom{-}0.4239585 + 1.1940415i & \theta_5^* &= \phantom{-}1.0714788\\
R_6 &= -1.2142847 + 0.2371991i & \theta_6^* &= \phantom{-}2.3404089\\
R_7 &= \phantom{-}0.0810748 - 1.3396423i & \theta_7^* &= -0.9654256\\
R_8 &= \phantom{-}0.7678941 - 1.5650217i & \theta_8^* &= -0.6986555\\
R_9 &= -0.2665917 + 1.2337601i & \theta_9^* &= \phantom{-}1.4369754\\
R_{10} &= -1.1430685 - 0.4882087i & \theta_{10}^* &= \phantom{-}3.1232587
\end{aligned} \tag{5.28}$$

By reconstructing the coefficients $C_1$, $C_2$, ..., $C_{10}$ via Vieta's relations (5.2) applied to reference roots $R_i$ in (5.28), we obtain numerical approximations to the coefficients of equation (5.27); these approximate coefficients $C_k$ are the ones that will be used as inputs for the LC method, which will generate discrete proximity maps built from $N = 5{,}000$ elements $LzC(\theta_k)$ associated to points $\theta_k$ in a regular partition of interval $[-\pi, \pi)$, following the strategies 5.1, 5.2 and 5.3 listed above in this chapter. It is worth mentioning that the absolute relative error between the approximate coefficient $C_{10}$ and the true coefficient $10 + 10i$ is equal to $4.387278 \times 10^{-15}$; this is the largest absolute relative error for the 10 coefficients reconstructed in this example, and is one order of magnitude greater than the machine epsilon for a personal computer with a 64-bit processor, which has a value of $2^{-52} \approx 2.220446 \times 10^{-16}$. We will continue with the usual methodology in this example, hoping that the minuscule errors (induced by the limited precision inherent in our personal computer) contained in the approximate coefficients $C_k$ that feed the LC method do not significantly influence the approximations $\hat{R}_i$ we intend to obtain. It is left as an exercise to the reader to verify that if we provide directly to the LC method the exact coefficients $C_k = k + ik$, $k = 1,2,\ldots,10$, instead of using the coefficients approximated via Vieta's relations (5.2), the results obtained will be practically the same as those shown in this example.

Figures 5.9, 5.10, and 5.11 show, respectively, the weighted error map $e(\theta)$ (LC map), map $\frac{d}{d\theta} d_\theta^2(t^*(\theta))$, and map $\frac{d}{d\theta} t^*(\theta)$.

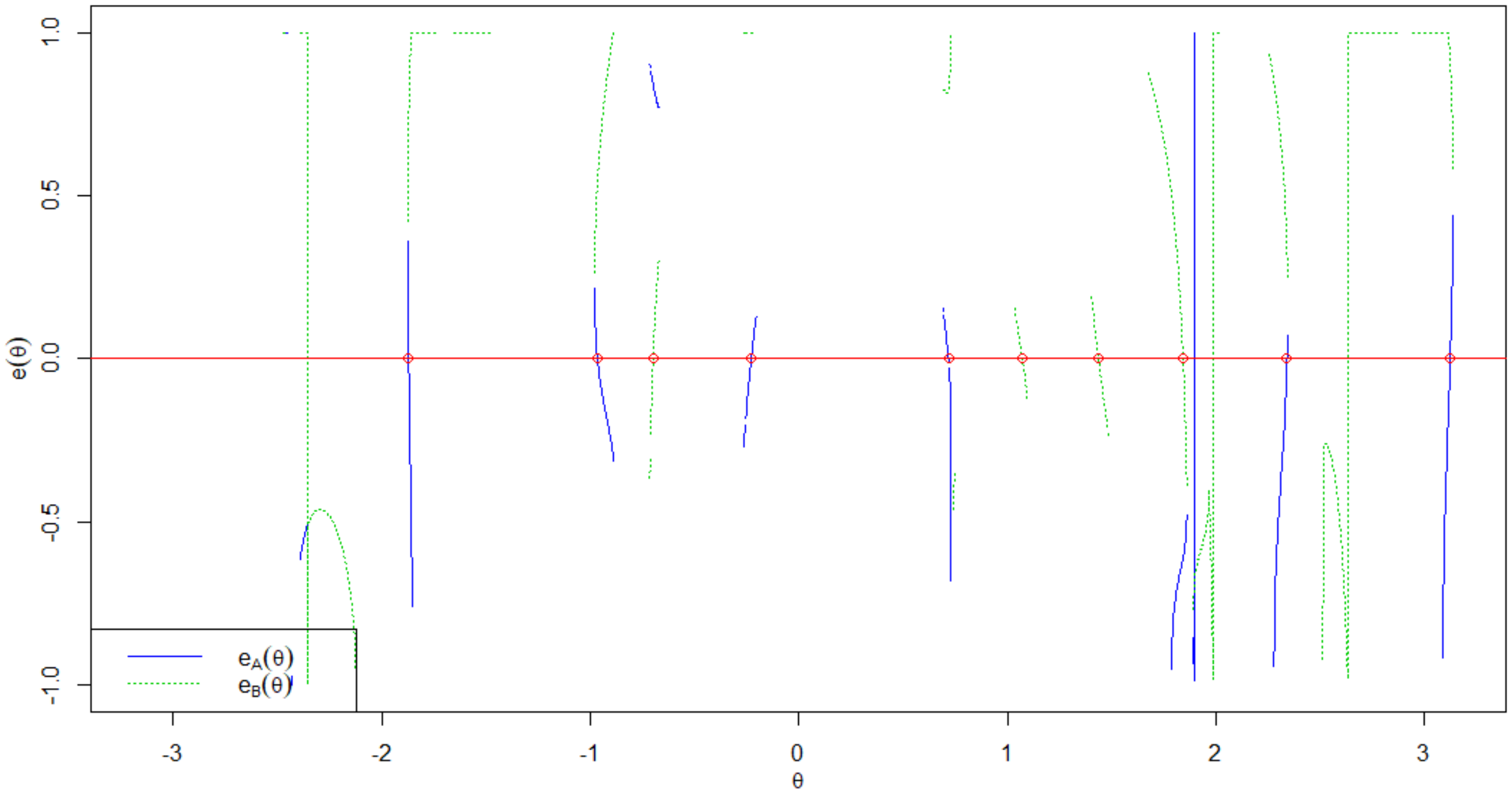

**Fig. 5.9.** Discrete angular proximity map $e(\theta)$ (LC map) for the polynomial equation (5.27) of degree $n = 10$. This map was generated from $N = 5{,}000$ elements $LzC(\theta_k)$ associated with points $\theta_k$ in a regular partition of interval $[-\pi, \pi)$. The reference theta roots $\theta_i^*$ are shown in this graph by means of small circles on horizontal axis $y = 0$.

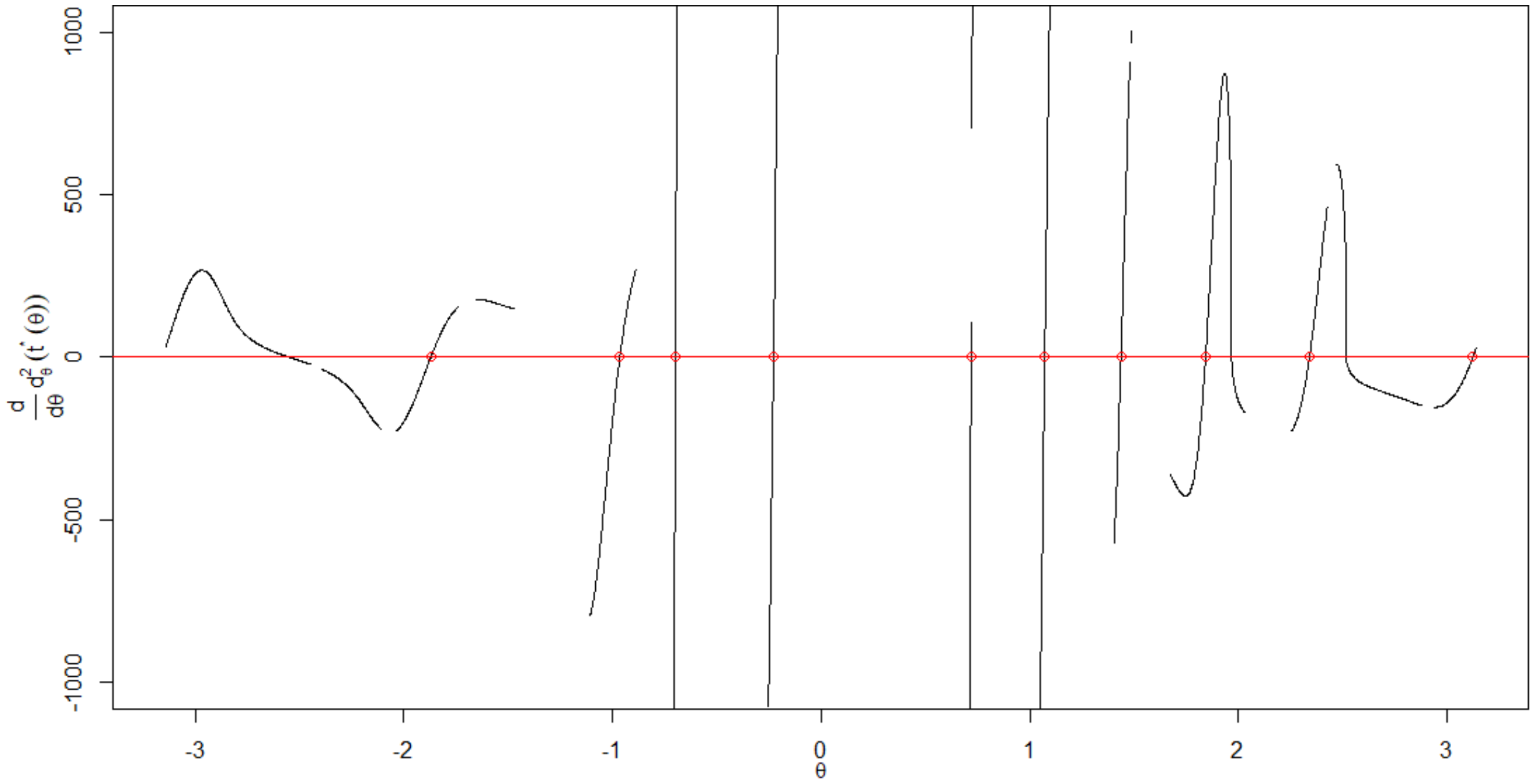

**Fig. 5.10.** Global map $\frac{d}{d\theta}d_\theta^2\big(t^*(\theta)\big)$ associated with the polynomial equation (5.27) of degree $n = 10$. This map was generated by means of $N = 5{,}000$ elements $LzC(\theta_k)$ associated with points $\theta_k$ in a regular partition of interval $[-\pi, \pi)$. True theta roots $\theta_i^*$ are shown in this graph by means of small circles on horizontal axis $y = 0$.

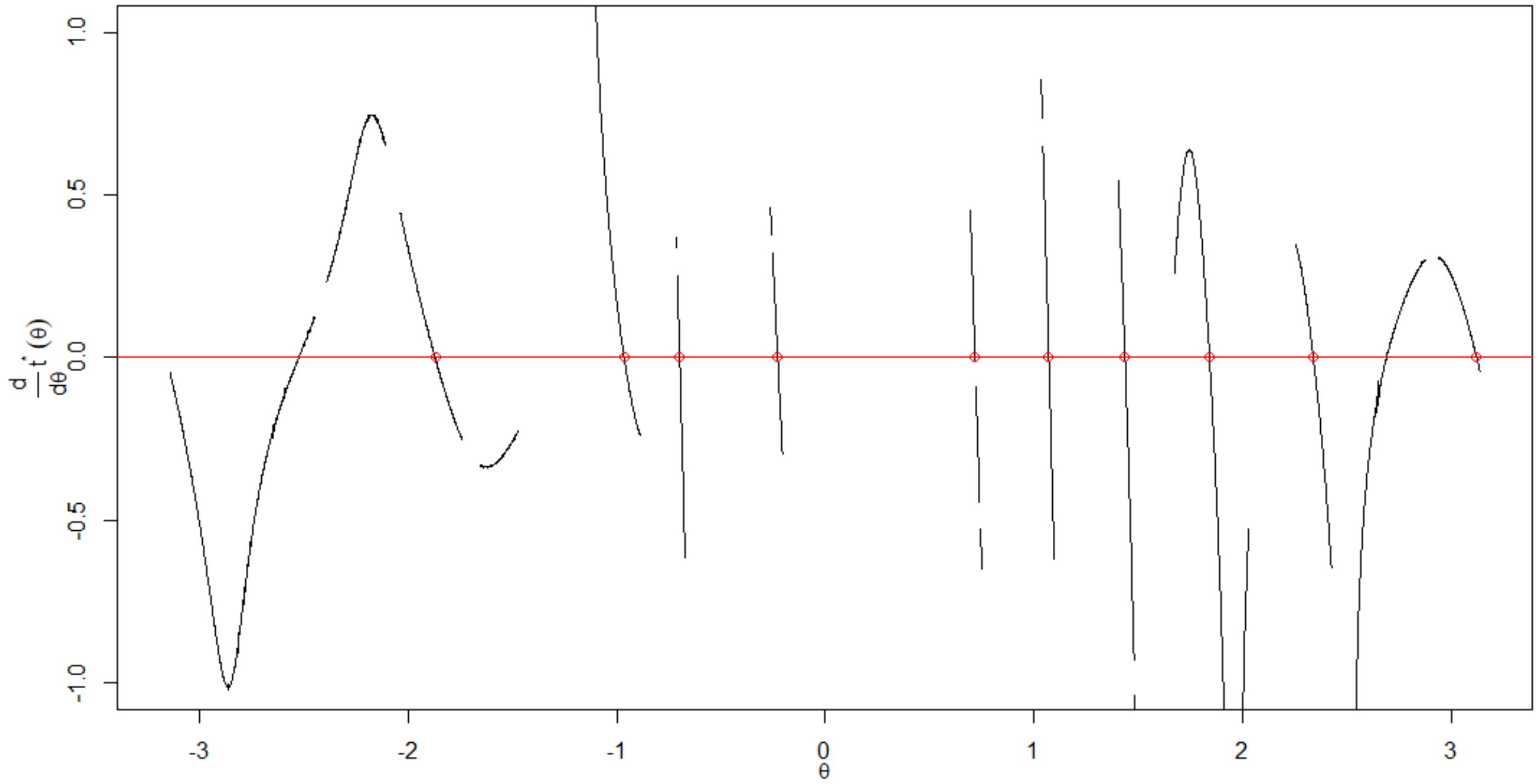

**Fig. 5.11.** Global map $\frac{d}{d\theta}t^*(\theta)$ associated with the polynomial equation (5.27) of degree $n = 10$. This map was generated by means of $N = 5{,}000$ elements $LzC(\theta_k)$ associated with points $\theta_k$ in a regular partition of interval $[-\pi, \pi)$. True theta roots $\theta_i^*$ are shown in this graph by means of small circles on horizontal axis $y = 0$.

Tables 5.5, 5.6 and 5.7 show the numerical characteristics of the smooth crossings with horizontal axis $y = 0$ detected in maps $e(\theta)$, $\frac{d}{d\theta} d_\theta^2(t^*(\theta))$, and $\frac{d}{d\theta} t^*(\theta)$ from respective figures 5.9, 5.10 and 5.11.

**Table 5.5**. Initial estimates for the roots of the polynomial of degree $n = 10$ in equation (5.27). These estimates were obtained from the map $e(\theta)$ (LC map) in figure 5.9.

| $i$ | $\hat{R}_i$ | $\hat{\theta}_i^*$ | $\lvert\Delta e_i\rvert$ | $d_{\hat{\theta}_i^*}^2\left(t^*(\hat{\theta}_i^*)\right)$ |
|---|---|---|---|---|
| 1 | $-0.8833565 + 0.8776421i$ | 1.8422001 | 0.018610925 | $4.741100 \times 10^{-11}$ |
| 2 | $-0.2665930 + 1.2337603i$ | 1.4369762 | 0.006661283 | $6.360430 \times 10^{-9}$ |
| 3 | $-1.1430660 - 0.4882054i$ | 3.1232534 | 0.025959192 | $3.597147 \times 10^{-8}$ |
| 4 | $0.4239610 + 1.1940399i$ | 1.0714772 | 0.005398199 | $5.944249 \times 10^{-8}$ |
| 5 | $0.0810789 - 1.3396393i$ | −0.9654207 | 0.008671058 | $6.085763 \times 10^{-8}$ |
| 6 | $0.9645595 - 0.8405508i$ | −0.2284679 | 0.007703705 | $7.552519 \times 10^{-8}$ |
| 7 | $0.7678934 - 1.5650196i$ | −0.6986548 | 0.020191203 | $1.216847 \times 10^{-7}$ |
| 8 | $-0.6774197 - 1.0722059i$ | −1.8714591 | 0.057621364 | $1.242850 \times 10^{-7}$ |
| 9 | $0.9472419 + 0.7629835i$ | 0.7175162 | 0.008829351 | $6.829469 \times 10^{-7}$ |
| 10 | $-1.2142567 + 0.2372232i$ | 2.3403730 | 0.035194044 | $2.926880 \times 10^{-6}$ |

**Table 5.6**. Initial estimates for the roots of the polynomial of degree $n = 10$ in equation (5.27). These estimates were obtained from the discrete map $\frac{d}{d\theta} d_\theta^2(t^*(\theta))$ in figure 5.10.

| $i$ | $\hat{R}_i$ | $\hat{\theta}_i^*$ | $\left\lvert\Delta \frac{d}{d\theta} d_{\hat{\theta}_i^*}^2(t^*)\right\rvert$ | $d_{\hat{\theta}_i^*}^2\left(t^*(\hat{\theta}_i^*)\right)$ |
|---|---|---|---|---|
| 1 | $-0.6774265 - 1.0722051i$ | −1.8714705 | 2.0887064 | $1.219591 \times 10^{-9}$ |
| 2 | $-0.2665908 + 1.2337598i$ | 1.4369749 | 26.9913992 | $3.219025 \times 10^{-9}$ |
| 3 | $0.4239577 + 1.1940419i$ | 1.0714792 | 62.6568360 | $5.068989 \times 10^{-9}$ |
| 4 | $-1.1430669 - 0.4882084i$ | 3.1232582 | 2.2009814 | $5.561905 \times 10^{-9}$ |
| 5 | $-1.2142824 + 0.2371983i$ | 2.3404079 | 5.6624053 | $1.276349 \times 10^{-8}$ |
| 6 | $0.0810754 - 1.3396401i$ | −0.9654239 | 6.2312970 | $1.293548 \times 10^{-8}$ |
| 7 | $-0.8833521 + 0.8776408i$ | 1.8421974 | 12.7907929 | $5.427298 \times 10^{-8}$ |
| 8 | $0.9645573 - 0.8405558i$ | −0.2284714 | 62.9592047 | $1.829769 \times 10^{-7}$ |
| 9 | $0.7678940 - 1.5650160i$ | −0.6986529 | 175.4249892 | $8.348369 \times 10^{-7}$ |
| 10 | $0.9472320 + 0.7629948i$ | 0.7175241 | 210.0621834 | $2.475134 \times 10^{-6}$ |
| 11 | $-1.0371018 - 0.1117528i$ | 2.5156912 | 6.8403781 | $6.623983 \times 10^{+1}$ |
| 12 | $-0.9183139 + 0.4996647i$ | 1.9671096 | 22.7829559 | $6.957491 \times 10^{+1}$ |
| 13 | $-0.7934549 - 0.6968381i$ | −2.5507610 | 0.2702611 | $7.323267 \times 10^{+1}$ |

**Table 5.7**. Initial estimates for the roots of the polynomial of degree $n = 10$ in equation (5.27). These estimates were obtained from the discrete map $\frac{d}{d\theta}t^*(\theta)$ in figure 5.11.

| $i$ | $\hat{R}_i$ | $\hat{\theta}_i^*$ | $\left\lvert\Delta\frac{d}{d\theta}t^*(\hat{\theta}_i^*)\right\rvert$ | $d^2_{\hat{\theta}_i^*}\left(t^*(\hat{\theta}_i^*)\right)$ |
|---|---|---|---|---|
| 1 | $-1.1430661 - 0.4882082i$ | 3.1232579 | 0.003265717 | $1.248898 \times 10^{-8}$ |
| 2 | $0.4239661 + 1.1940370i$ | 1.0714742 | 0.028744965 | $5.169979 \times 10^{-7}$ |
| 3 | $-0.8833319 + 0.8776487i$ | 1.8421823 | 0.011881402 | $1.636674 \times 10^{-6}$ |
| 4 | $0.7678853 - 1.5650263i$ | −0.6986610 | 0.034965612 | $2.474036 \times 10^{-6}$ |
| 5 | $0.9645634 - 0.8405377i$ | −0.2284588 | 0.014996786 | $2.794531 \times 10^{-6}$ |
| 6 | $0.9472584 + 0.7629649i$ | 0.7175033 | 0.024839967 | $2.006078 \times 10^{-5}$ |
| 7 | $0.0809470 - 1.3397288i$ | −0.9655767 | 0.003251080 | $5.679994 \times 10^{-5}$ |
| 8 | $-1.2141658 + 0.2373139i$ | 2.3402479 | 0.007153389 | $5.865969 \times 10^{-5}$ |
| 9 | $-0.2667591 + 1.2337799i$ | 1.4370718 | 0.024765968 | $9.933535 \times 10^{-5}$ |
| 10 | $-0.6777714 - 1.0720959i$ | −1.8720743 | 0.004791710 | $3.029048 \times 10^{-4}$ |
| 11 | $-1.0006775 - 0.2575788i$ | 2.6906757 | 0.003483489 | $5.258528 \times 10^{+1}$ |
| 12 | $-0.7879306 - 0.7039383i$ | −2.5253247 | 0.001328513 | $7.316346 \times 10^{+1}$ |

From tables 5.5, 5.6 and 5.7, whose rows are in ascending order with respect to column $d^2_{\hat{\theta}_i^*}\left(t^*(\hat{\theta}_i^*)\right)$, we can see that the first 10 rows within each table contain close approximations to each of the 10 reference roots in (5.28). From table 5.6 we can see that quantities $\left\lvert\Delta\frac{d}{d\theta}d^2_{\hat{\theta}_i^*}(t^*)\right\rvert$ are significantly larger in comparison with quantities $|\Delta e_i|$ and $\left\lvert\Delta\frac{d}{d\theta}t^*(\hat{\theta}_i^*)\right\rvert$ from tables 5.5 and 5.7; this indicates that, in this case, the map $\frac{d}{d\theta}d^2_\theta\left(t^*(\theta)\right)$ is much more sensitive to norm differences among polynomial coefficients in (5.27), compared to maps $e(\theta)$ and $\frac{d}{d\theta}t^*(\theta)$; it is for this reason that the argument `Tol` had to be set to `1000.0` in the corresponding call to function `approxDMin` within the script from annex 3 section 3. In theory, map $e(\theta)$ is more robust to these differences compared to maps $\frac{d}{d\theta}d^2_\theta\left(t^*(\theta)\right)$ and $\frac{d}{d\theta}t^*(\theta)$, since the functional values of $e(\theta)$ are always restricted to the interval $(-1,1]$. Certainly, there are also situations where it is necessary to make similar adjustments to argument `Tol` when calling function `approxDMin`, in order to correctly detect smooth crossings of map $\frac{d}{d\theta}t^*(\theta)$ with horizontal axis $y = 0$; see, for instance, table 5.11 in numerical example 5.3 of this chapter.

Carrying on with the description of results in this numerical example, let us now see the absolute relative errors for the first 10 rows of tables 5.5, 5.6, and 5.7, with respect to reference values (5.28).

For the estimates in table 5.5 (corresponding to map $e(\theta)$, also called LC map), the absolute relative errors, with respect to reference values (5.28), are:

$$\left|\frac{R_1-\hat{R}_9}{R_1}\right| = 4.513227 \times 10^{-6} \qquad \left|\frac{\theta_1^*-\hat{\theta}_9^*}{\theta_1^*}\right| = 3.554942 \times 10^{-6}$$
$$\left|\frac{R_2-\hat{R}_1}{R_2}\right| = 1.108561 \times 10^{-7} \qquad \left|\frac{\theta_2^*-\hat{\theta}_1^*}{\theta_2^*}\right| = 5.095804 \times 10^{-8}$$
$$\left|\frac{R_3-\hat{R}_8}{R_3}\right| = 5.774161 \times 10^{-6} \qquad \left|\frac{\theta_3^*-\hat{\theta}_8^*}{\theta_3^*}\right| = 6.311645 \times 10^{-6}$$
$$\left|\frac{R_4-\hat{R}_6}{R_4}\right| = 2.038412 \times 10^{-6} \qquad \left|\frac{\theta_4^*-\hat{\theta}_6^*}{\theta_4^*}\right| = 6.270144 \times 10^{-6}$$
$$\left|\frac{R_5-\hat{R}_4}{R_5}\right| = 2.350102 \times 10^{-6} \qquad \left|\frac{\theta_5^*-\hat{\theta}_4^*}{\theta_5^*}\right| = 1.434912 \times 10^{-6} \tag{5.29}$$
$$\left|\frac{R_6-\hat{R}_{10}}{R_6}\right| = 2.982952 \times 10^{-5} \qquad \left|\frac{\theta_6^*-\hat{\theta}_{10}^*}{\theta_6^*}\right| = 1.533496 \times 10^{-5}$$
$$\left|\frac{R_7-\hat{R}_5}{R_7}\right| = 3.764446 \times 10^{-6} \qquad \left|\frac{\theta_7^*-\hat{\theta}_5^*}{\theta_7^*}\right| = 5.122418 \times 10^{-6}$$
$$\left|\frac{R_8-\hat{R}_7}{R_8}\right| = 1.254335 \times 10^{-6} \qquad \left|\frac{\theta_8^*-\hat{\theta}_7^*}{\theta_8^*}\right| = 1.044933 \times 10^{-6}$$
$$\left|\frac{R_9-\hat{R}_2}{R_9}\right| = 1.068964 \times 10^{-6} \qquad \left|\frac{\theta_9^*-\hat{\theta}_2^*}{\theta_9^*}\right| = 5.367404 \times 10^{-7}$$
$$\left|\frac{R_{10}-\hat{R}_3}{R_{10}}\right| = 3.312876 \times 10^{-6} \qquad \left|\frac{\theta_{10}^*-\hat{\theta}_3^*}{\theta_{10}^*}\right| = 1.679262 \times 10^{-6}$$

For the estimates in table 5.6, corresponding to map $\frac{d}{d\theta}d_\theta^2\big(t^*(\theta)\big)$, the absolute relative errors, with respect to reference values (5.28), are:

$$\left|\frac{R_1-\hat{R}_{10}}{R_1}\right| = 8.592202 \times 10^{-6} \qquad \left|\frac{\theta_1^*-\hat{\theta}_{10}^*}{\theta_1^*}\right| = 7.367247 \times 10^{-6}$$
$$\left|\frac{R_2-\hat{R}_7}{R_2}\right| = 3.750734 \times 10^{-6} \qquad \left|\frac{\theta_2^*-\hat{\theta}_7^*}{\theta_2^*}\right| = 1.523598 \times 10^{-6}$$
$$\left|\frac{R_3-\hat{R}_1}{R_3}\right| = 5.719834 \times 10^{-7} \qquad \left|\frac{\theta_3^*-\hat{\theta}_1^*}{\theta_3^*}\right| = 2.579602 \times 10^{-7}$$
$$\left|\frac{R_4-\hat{R}_8}{R_4}\right| = 3.172836 \times 10^{-6} \qquad \left|\frac{\theta_4^*-\hat{\theta}_8^*}{\theta_4^*}\right| = 9.273666 \times 10^{-6}$$
$$\left|\frac{R_5-\hat{R}_3}{R_5}\right| = 6.862764 \times 10^{-7} \qquad \left|\frac{\theta_5^*-\hat{\theta}_3^*}{\theta_5^*}\right| = 4.191217 \times 10^{-7} \tag{5.30}$$
$$\left|\frac{R_6-\hat{R}_5}{R_6}\right| = 1.969725 \times 10^{-6} \qquad \left|\frac{\theta_6^*-\hat{\theta}_5^*}{\theta_6^*}\right| = 4.390896 \times 10^{-7}$$
$$\left|\frac{R_7-\hat{R}_6}{R_7}\right| = 1.735537 \times 10^{-6} \qquad \left|\frac{\theta_7^*-\hat{\theta}_6^*}{\theta_7^*}\right| = 1.750448 \times 10^{-6}$$
$$\left|\frac{R_8-\hat{R}_9}{R_8}\right| = 3.285497 \times 10^{-6} \qquad \left|\frac{\theta_8^*-\hat{\theta}_9^*}{\theta_8^*}\right| = 3.756293 \times 10^{-6}$$
$$\left|\frac{R_9-\hat{R}_2}{R_9}\right| = 7.604714 \times 10^{-7} \qquad \left|\frac{\theta_9^*-\hat{\theta}_2^*}{\theta_9^*}\right| = 3.726508 \times 10^{-7}$$
$$\left|\frac{R_{10}-\hat{R}_4}{R_{10}}\right| = 1.302674 \times 10^{-6} \qquad \left|\frac{\theta_{10}^*-\hat{\theta}_4^*}{\theta_{10}^*}\right| = 1.713556 \times 10^{-7}$$

For the estimates in table 5.7, corresponding to map $\frac{d}{d\theta}t^*(\theta)$, the absolute relative errors, with respect to reference values (5.28), are:

$$\begin{aligned}
&\left|\frac{R_1-\hat{R}_6}{R_1}\right| = 2.445956\times10^{-5} && \left|\frac{\theta_1^*-\hat{\theta}_6^*}{\theta_1^*}\right| = 2.152361\times10^{-5}\\
&\left|\frac{R_2-\hat{R}_3}{R_2}\right| = 2.059752\times10^{-5} && \left|\frac{\theta_2^*-\hat{\theta}_3^*}{\theta_2^*}\right| = 9.735897\times10^{-6}\\
&\left|\frac{R_3-\hat{R}_{10}}{R_3}\right| = 2.849896\times10^{-4} && \left|\frac{\theta_3^*-\hat{\theta}_{10}^*}{\theta_3^*}\right| = 3.223760\times10^{-4}\\
&\left|\frac{R_4-\hat{R}_5}{R_4}\right| = 1.239920\times10^{-5} && \left|\frac{\theta_4^*-\hat{\theta}_5^*}{\theta_4^*}\right| = 4.613651\times10^{-5}\\
&\left|\frac{R_5-\hat{R}_2}{R_5}\right| = 6.930781\times10^{-6} && \left|\frac{\theta_5^*-\hat{\theta}_2^*}{\theta_5^*}\right| = 4.244679\times10^{-6}\\
&\left|\frac{R_6-\hat{R}_8}{R_6}\right| = 1.335656\times10^{-4} && \left|\frac{\theta_6^*-\hat{\theta}_8^*}{\theta_6^*}\right| = 6.878652\times10^{-5}\\
&\left|\frac{R_7-\hat{R}_7}{R_7}\right| = 1.149814\times10^{-4} && \left|\frac{\theta_7^*-\hat{\theta}_7^*}{\theta_7^*}\right| = 1.565306\times10^{-4}\\
&\left|\frac{R_8-\hat{R}_4}{R_8}\right| = 5.655844\times10^{-6} && \left|\frac{\theta_8^*-\hat{\theta}_4^*}{\theta_8^*}\right| = 7.884945\times10^{-6}\\
&\left|\frac{R_9-\hat{R}_9}{R_9}\right| = 1.335745\times10^{-4} && \left|\frac{\theta_9^*-\hat{\theta}_9^*}{\theta_9^*}\right| = 6.706050\times10^{-5}\\
&\left|\frac{R_{10}-\hat{R}_1}{R_{10}}\right| = 1.952038\times10^{-6} && \left|\frac{\theta_{10}^*-\hat{\theta}_1^*}{\theta_{10}^*}\right| = 2.633427\times10^{-7}
\end{aligned} \tag{5.31}$$

Table 5.8 summarizes the errors of the estimates corresponding to the proximity maps obtained here; it contains arithmetic means and standard deviations for the approximation measures $d^2_{\hat{\theta}_i^*}\left(t^*(\hat{\theta}_i^*)\right)$ from tables 5.5, 5.6 and 5.7, as well as for the absolute relative errors (5.29), (5.30) and (5.31).

**Table 5.8**. Arithmetic means and standard deviations for error measures associated with valid initial estimates (first ten rows in each of tables 5.5, 5.6 and 5.7) of the roots of the polynomial of degree $n = 10$ in equation (5.27), obtained by means of proximity maps $e(\theta)$, $\frac{d}{d\theta}d^2_\theta\left(t^*(\theta)\right)$, and $\frac{d}{d\theta}t^*(\theta)$ from figures 5.9, 5.10 and 5.11.

| Map | Error measure | | | | | |
|---|---|---|---|---|---|---|
| | $d^2_{\hat{\theta}_i^*}\left(t^*(\hat{\theta}_i^*)\right)$ | | $\left\lvert\frac{R-\hat{R}}{R}\right\rvert$ | | $\left\lvert\frac{\theta^*-\hat{\theta}^*}{\theta^*}\right\rvert$ | |
| | Mean | Standard deviation | Mean | Standard deviation | Mean | Standard deviation |
| $e(\theta)$ | $4.094001\times10^{-7}$ | $9.068558\times10^{-7}$ | $5.401690\times10^{-6}$ | $8.751384\times10^{-6}$ | $4.134092\times10^{-6}$ | $4.577041\times10^{-6}$ |
| $\frac{d}{d\theta}d^2_\theta(t^*(\theta))$ | $3.587990\times10^{-7}$ | $7.868856\times10^{-7}$ | $2.582794\times10^{-6}$ | $2.403988\times10^{-6}$ | $2.533143\times10^{-6}$ | $3.269514\times10^{-6}$ |
| $\frac{d}{d\theta}t^*(\theta)$ | $5.451953\times10^{-5}$ | $9.367885\times10^{-5}$ | $7.391060\times10^{-5}$ | $9.234839\times10^{-5}$ | $7.045426\times10^{-5}$ | $1.004355\times10^{-4}$ |

This time, we see that the approximations $\hat{R}_i$, $\hat{\theta}_i^*$ from map $\frac{d}{d\theta}d_\theta^2\left(t^*(\theta)\right)$ are the most consistent according to the error measures considered, followed very closely by the approximations from map $e(\theta)$; the error measures for the approximations from map $\frac{d}{d\theta}t^*(\theta)$ are on average and on standard deviation one or two orders of magnitude larger than their counterparts from the other two maps; even so, we consider that the initial approximations $\hat{R}_i$, $\hat{\theta}_i^*$ corresponding to map $\frac{d}{d\theta}t^*(\theta)$ are of reasonable quality, because in the worst case, the least accurate approximation $\hat{R}_{10}$, $\hat{\theta}_{10}^*$ corresponding to map $\frac{d}{d\theta}t^*(\theta)$ in table 5.7 coincides with at least the first three digits (both for the real and imaginary parts) of reference root $R_3$, $\theta_3^*$ in (5.28), while in a better scenario the most accurate approximation $\hat{R}_1$, $\hat{\theta}_1^*$ corresponding to map $\frac{d}{d\theta}t^*(\theta)$ in table 5.7 coincides with at least the first six digits (both for the real and imaginary parts) of reference root $R_{10}$, $\theta_{10}^*$ in (5.28).

To continue with our analysis of the results in this numerical example let us see, in figures 5.12, 5.13, 5.14 and 5.15, the graphs related to the global minima of the dynamic squared distances between terminal curves and z-circumferences associated with equation (5.27), as a function of the inclination angle $\theta$ for line $\ell_1$ with fixed point $P_1 = -C_1/2$.

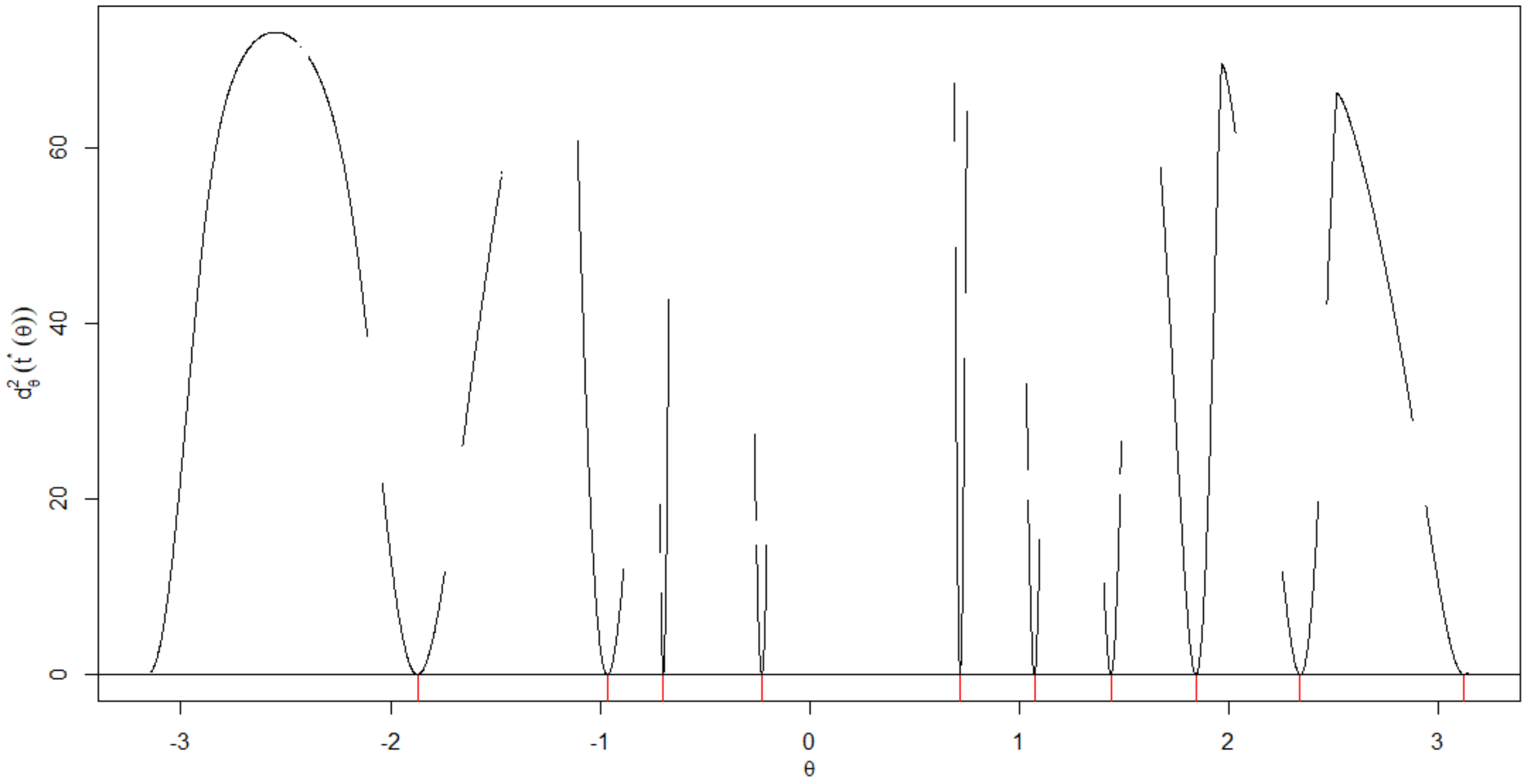

**Fig. 5.12.** Global minima $d_\theta^2\left(t^*(\theta)\right)$ vs. $\theta$ associated with equation (5.27) of degree $n = 10$. This map was generated from $N = 5{,}000$ elements $LzC(\theta_k)$ associated with points $\theta_k$ in a regular partition of interval $[-\pi, \pi)$. This graph includes vertical line segments that indicate the location of true theta roots $\theta_i^*$.

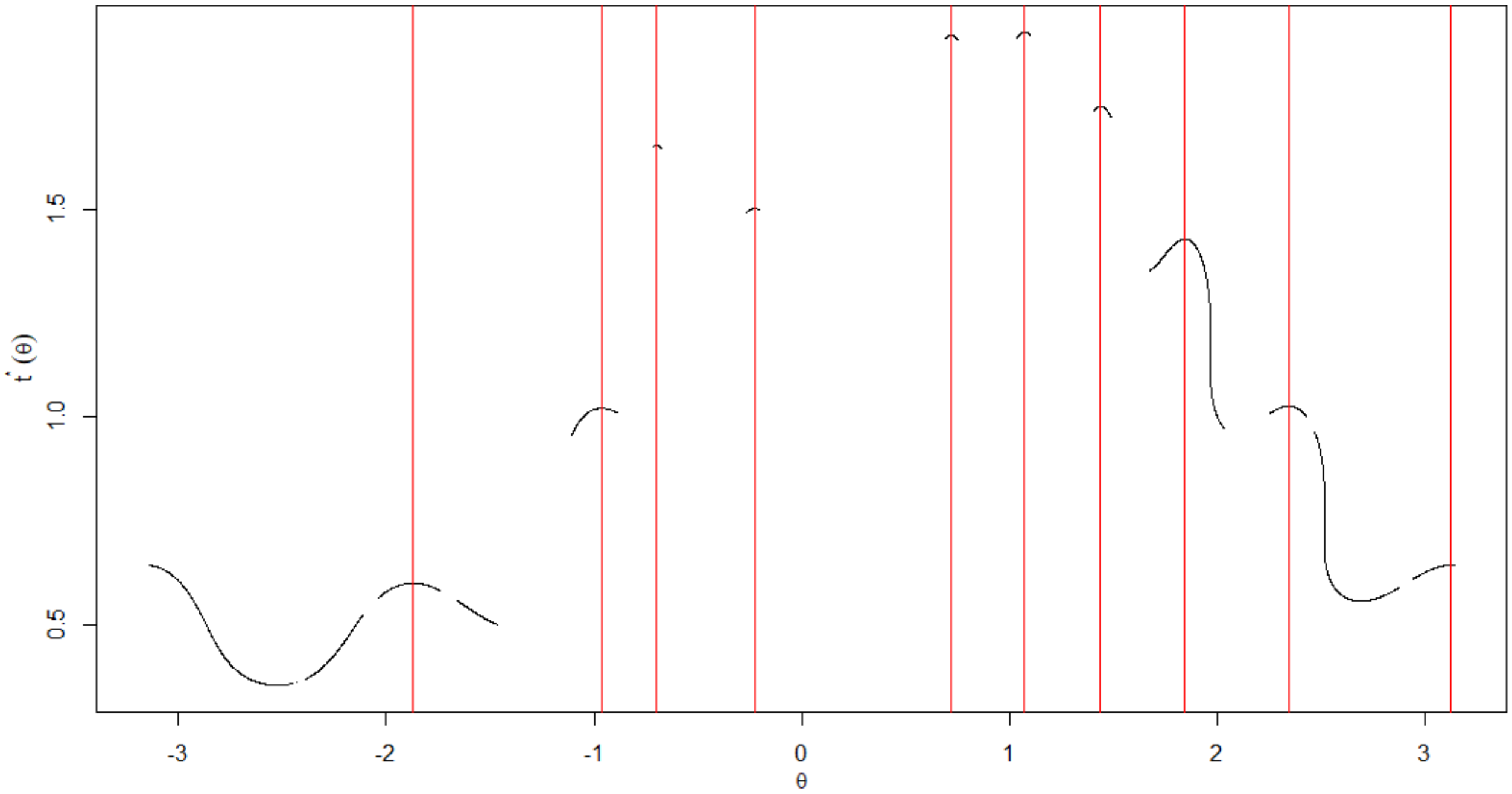

**Fig. 5.13.** Minimizing arguments $t^*(\theta)$ vs. $\theta$ associated with equation (5.27) of degree $n = 10$. This map was generated from $N = 5{,}000$ elements $LzC(\theta_k)$ associated with points $\theta_k$ in a regular partition of interval $[-\pi, \pi)$. This graph includes vertical line segments that indicate the location of true theta roots $\theta_i^*$.

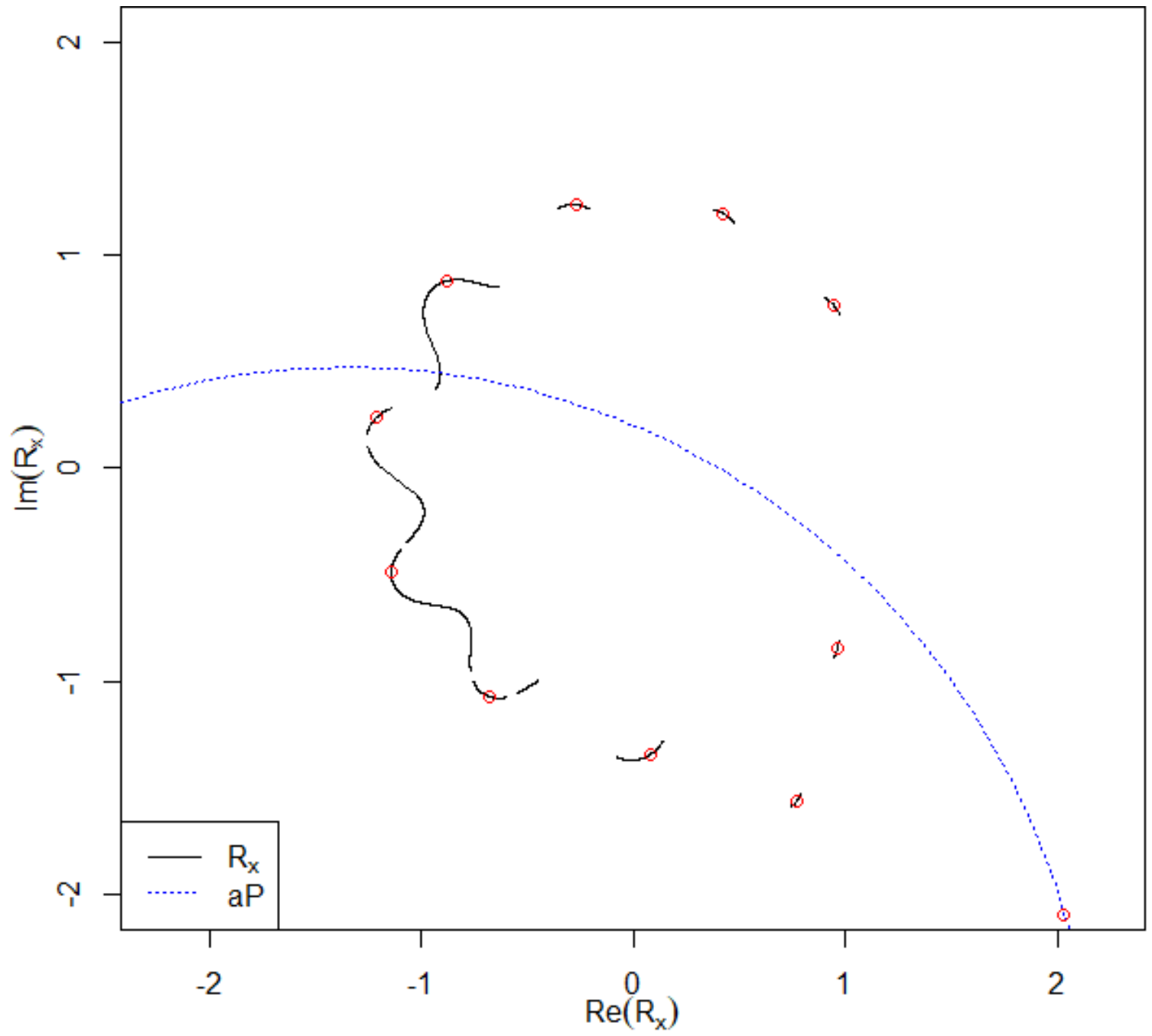

**Fig. 5.14.** Trajectory $R_x(\theta)$ and part of trajectory $aP(\theta)$, both associated with equation (5.27) of degree $n = 10$. These trajectories were generated from $N = 5{,}000$ elements $LzC(\theta_k)$ associated with points $\theta_k$ in a regular partition of interval $[-\pi, \pi)$. The true roots $R_i(\theta_i^*)$, and their corresponding anchor points $aP(\theta_i^*)$, are represented in this graph by means of small circles on the trajectories $R_x(\theta)$, $aP(\theta)$.

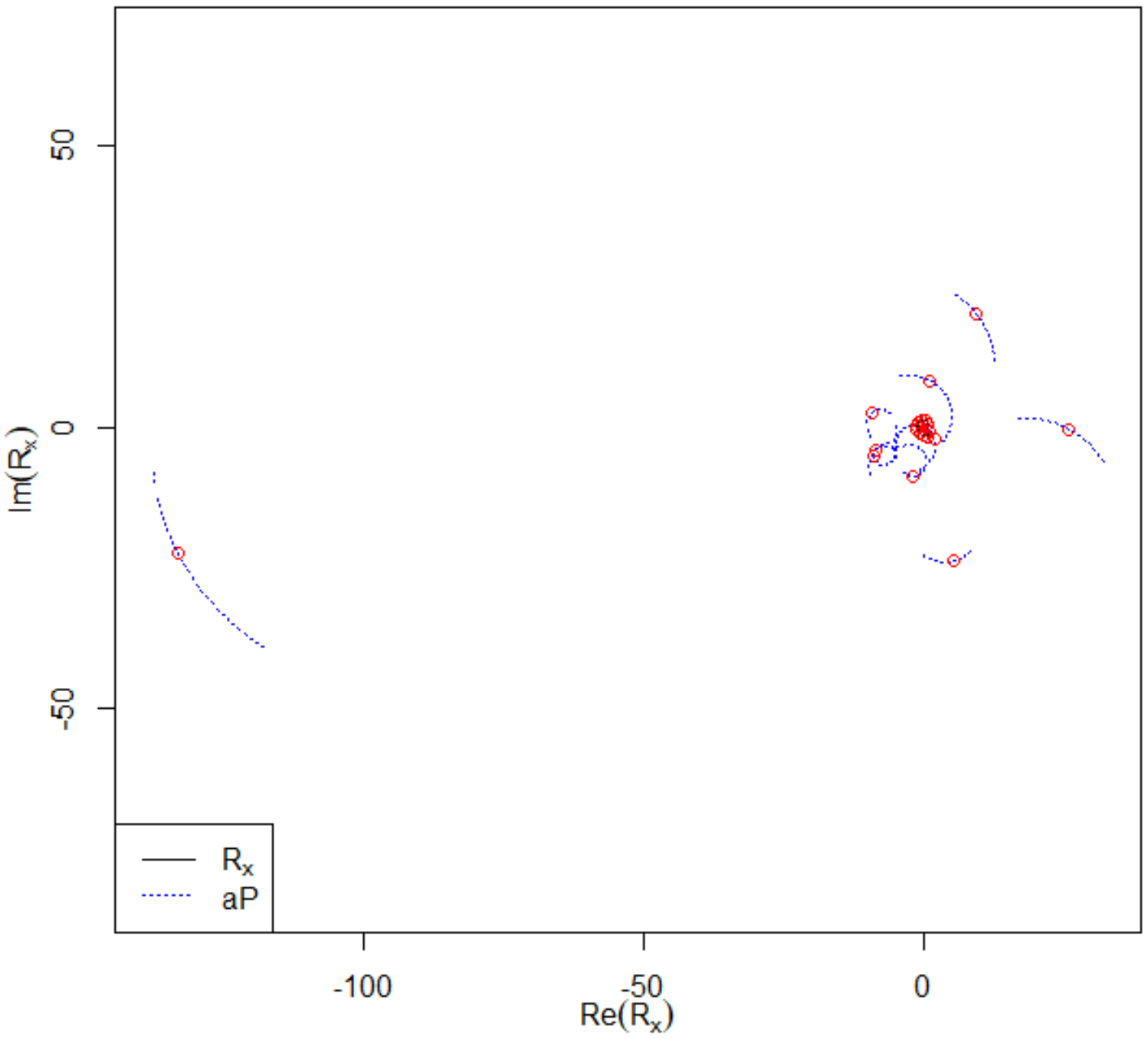

**Fig. 5.15.** Trajectories $R_x(\theta)$ and $aP(\theta)$ associated with equation (5.27) of degree $n = 10$. These trajectories are the same as those in figure 5.14 but displayed over a larger spatial extent in order to fully visualize trajectory $aP(\theta)$. Trajectory $R_x(\theta)$ is seen here as a small cluster of points centered at the origin.

The graphs in figures 5.10, 5.11, 5.12 and 5.13 provide evidence in favor of the hypotheses about the properties of the mappings $d_\theta^2(t^*(\theta))$ and $t^*(\theta)$ posed in numerical example 4.1 from chapter 4: figures 5.10 and 5.12 suggest that theta roots $\theta_i^*$ occur at the global minima of map $d_\theta^2(t^*(\theta))$, which by the way has a significantly wider vertical range relative to its counterpart from example 5.1 (compare figures 5.4 and 5.12); on the other hand, figures 5.11 and 5.13 suggest that, in this case, all theta roots $\theta_i^*$ correspond to local maxima of map $t^*(\theta)$.

In any of the figures 5.10, 5.11, 5.12, 5.13 and 5.14, we can see that, apparently, the resulting proximity maps are made up of at least 13 continuous sections; in reality, existing gaps attributed to errors in the optimization process make it difficult to determine with certainty the number of continuous sections in these maps. Figure 5.15, on the other hand, shows that the trajectory $aP(\theta)$ evolves over a much larger spatial extent with respect to the trajectory $R_x(\theta)$, and if we compare figures 5.14 and 5.15, we see that trajectories $aP(\theta)$ and $R_x(\theta)$ are not similar; as we saw in example 5.1 of this chapter, this similarity property between $aP(\theta)$ and $R_x(\theta)$ holds only for univariate polynomials of degree $n = 4$, since in that case $aP(\theta)$ depends linearly on $R_x(\theta)$, while for $n > 4$, $aP(\theta)$ is a polynomial expression in terms of powers of $R_x(\theta)$ greater than 1; see for instance expression (5.23) for $n = 7$.

To conclude this example, let us look at the graphs of a structure $LzC(\theta_k)$ and of its corresponding function $d^2_{\theta_k}(t)$, displayed with the help of the script listed in annex 4 section 2 (although it is necessary to modify some lines in this script, and add some others, in order to properly display the graphs). Figure 5.16 shows the structure $LzC(\theta = 1.4313)$ associated with equation (5.27), and figure 5.17 shows the corresponding dynamic squared distance function $d^2_{\theta=1.4313}(t)$ with some relevant details, such as the specific values of the global minimum $d^2_{\theta=1.4313}(t^*)$ and the corresponding minimizing argument $t^*$. From figure 5.16, we can observe the following:

- Structure $LzC(\theta = 1.4313)$ is close to reference root $R_9$ in (5.28), and simultaneously to the 9-product $R_1R_2R_3R_4R_5R_6R_7R_8R_{10}$, shown as a small point very near to z-circumference $zC$.

- Roots $R_i$ in this case are closer to the origin, while the 9-products are significantly farther from the origin.

- Trajectories $tL$ and $t\mathfrak{C}$ intersect (both statically and dynamically) at a point corresponding to the minimizing argument $t^*$ of associated function $d^2_{\theta=1.4313}(t)$; this intersection tends to the 9-product $R_1R_2R_3R_4R_5R_6R_7R_8R_{10}$. As we demonstrated in example 5.1 of this chapter, this intersection between $tL$ and $t\mathfrak{C}$ does not happen by chance.

- Terminal semi-line $tL$ intersects z-circumference $zC$ at a point that also tends to the 9-product $R_1R_2R_3R_4R_5R_6R_7R_8R_{10}$, although with its own approaching behavior, different from that of the intersection between $tL$ and $t\mathfrak{C}$, and also different from that of the diamond ◊ that represents point $zC(\theta = 1.4313, t^*)$.

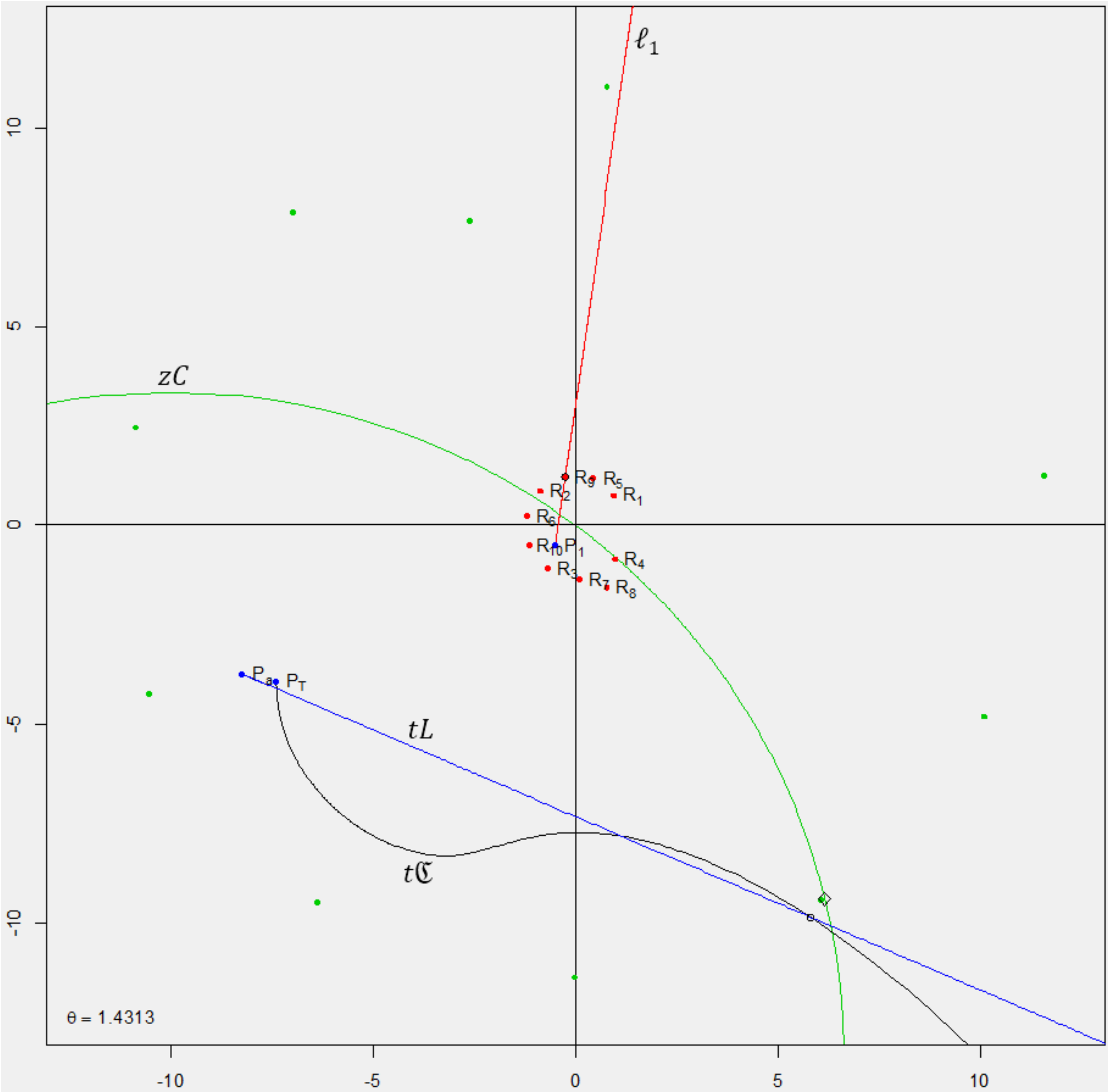

**Fig. 5.16.** Some elements of the structure $LzC(\theta)$ associated with equation (5.27) of degree $n = 10$, for $\theta = 1.4313$. The small circles on line $\ell_1$, terminal curve $t\mathfrak{C}$, terminal semi-line $tL$, as well as the small diamond ◊ on z-circumference $zC$, correspond to the minimizing argument $t^*$ of associated function $d_\theta^2$ (see figure 5.17). On the other hand, points $P_a$ in $tL$, $P_T$ in $t\mathfrak{C}$, $P_1$ in $\ell_1$ and $P_c$ in $zC$ (the latter not shown here) correspond to $t = 0$ in these parametric trajectories. The position of the anchor point of $tL$, $P_a$, depends on $\theta$, but not the positions of points $P_1$, $P_c$, and $P_T$, which are true fixed points, independent of $\theta$.

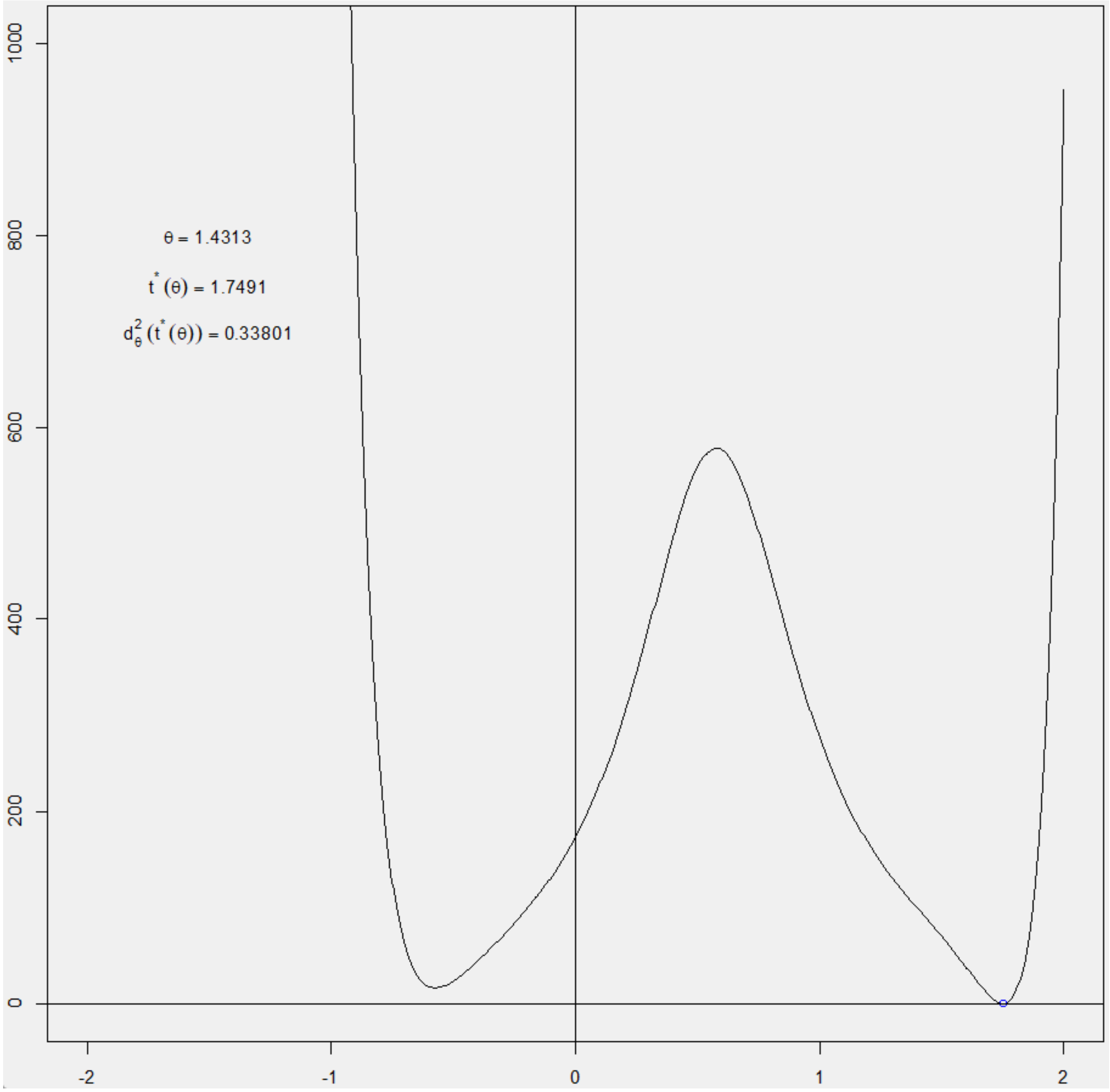

**Fig. 5.17.** Dynamic squared distance function $d_\theta^2(t)$ associated with equation (5.27) of degree $n = 10$, for $\theta = 1.4313$. The global minimum of this function occurs at minimizing argument $t^* = 1.7491$. At $t^*$, the minimum functional value is $d_\theta^2(t^*) = 0.33801$.

From figure 5.17, we see that the functional values of $d_\theta^2(t)$ are typically four orders of magnitude larger than their counterparts in example 5.1 of this chapter (figure 5.8); this, of course, has to do with the norms of both the coefficients (5.19) involved in example 5.1, and the coefficients in equation (5.27).

## Numerical Example 5.3

In this example we will try to approximate, by means of the LC method, the roots of an equation of degree $n = 5$ that contains only real roots. Specifically, we will work with the polynomial equation

$$x^5 - 15x^4 + 85x^3 - 225x^2 + 274x - 120 = 0, \qquad (5.32)$$

whose roots (which of course will be the reference values that will help us evaluate the approximations obtained) are $R_1 = 1$, $R_2 = 2$, $R_3 = 3$, $R_4 = 4$, and $R_5 = 5$; the polynomial on the left side of equation (5.32) is known as *Wilkinson's polynomial* of degree $n = 5$. Note that, under the LC method, all the theta roots $\theta_i^*$ of equation (5.32) are equal to $\pi$, since the line $\ell_1(\pi)$ with fixed point $P_1 = -C_1/2 = 7.5$ contains all the roots. In a certain way, in this example we will face a situation where we have a polynomial with repeated roots (at least from the point of view of proximity maps associated with univariate polynomials).

Initially, we will try to approximate the roots of equation (5.32) by using directly as input for the LC method the coefficients reconstructed from the reference roots and Vieta's relations (5.2), and whose theoretical values are $C_1 = -15 + 0i$, $C_2 = 85 + 0i$, $C_3 = -225 + 0i$, $C_4 = 274 + 0i$ and $C_5 = -120 + 0i$. These will therefore be the inputs for the script listed in annex 3 section 3; we will use the same operation parameters used in example 5.2 of this chapter. In this way, we obtain the proximity maps shown in figures 5.18, 5.19 and 5.20.

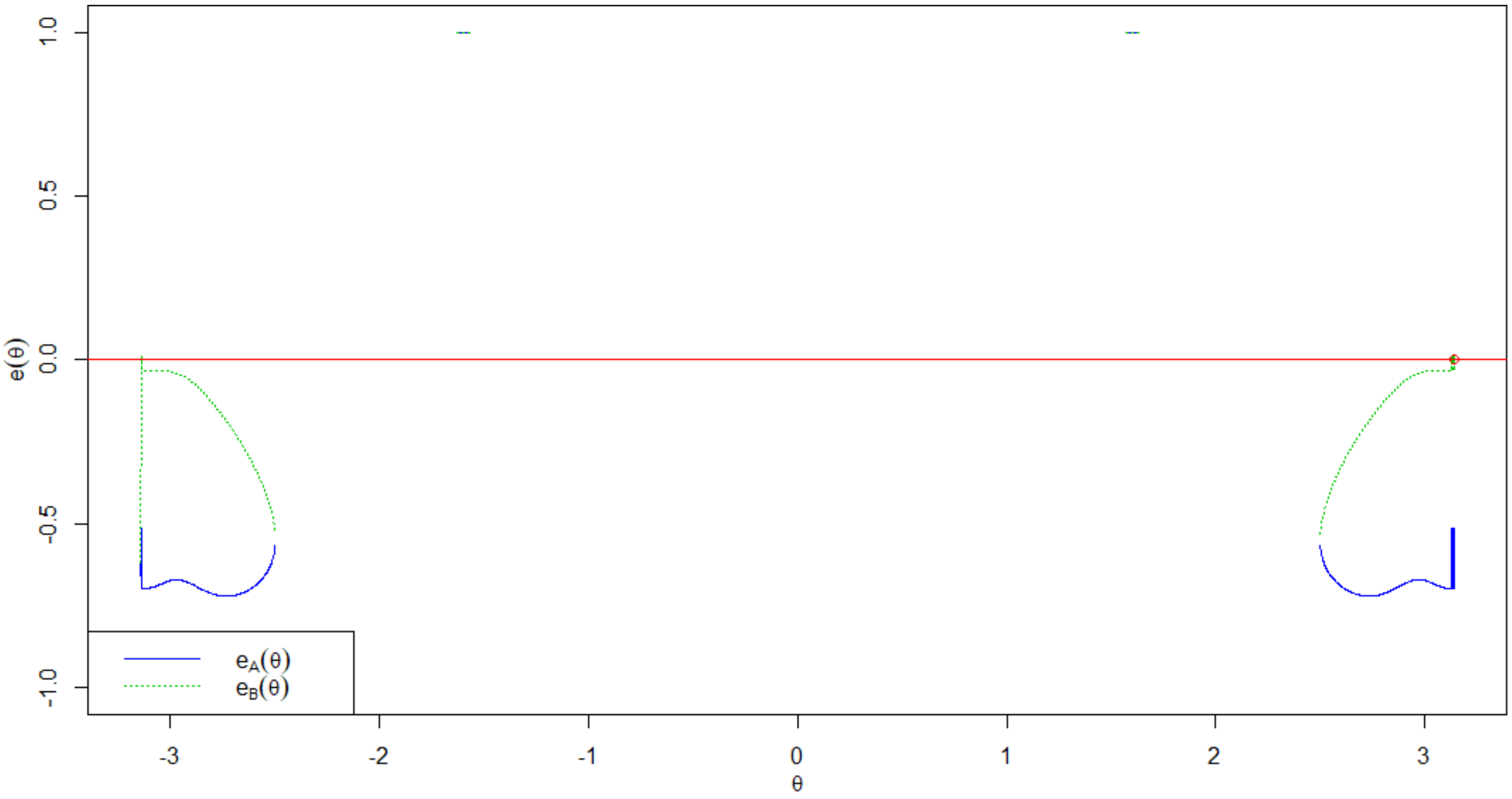

**Fig. 5.18.** Discrete proximity map $e(\theta)$ for the polynomial $x^5 - 15x^4 + 85x^3 - 225x^2 + 274x - 120$. This map was generated from $N = 5{,}000$ elements $LzC(\theta_k)$ associated with points $\theta_k$ in a regular partition of interval $[-\pi, \pi)$. Reference theta root $\theta_i^* = \pi$ for $i = 1,2,3,4,5$, is shown by means of a small circle on horizontal axis $y = 0$, just on the right edge of the proximity map.

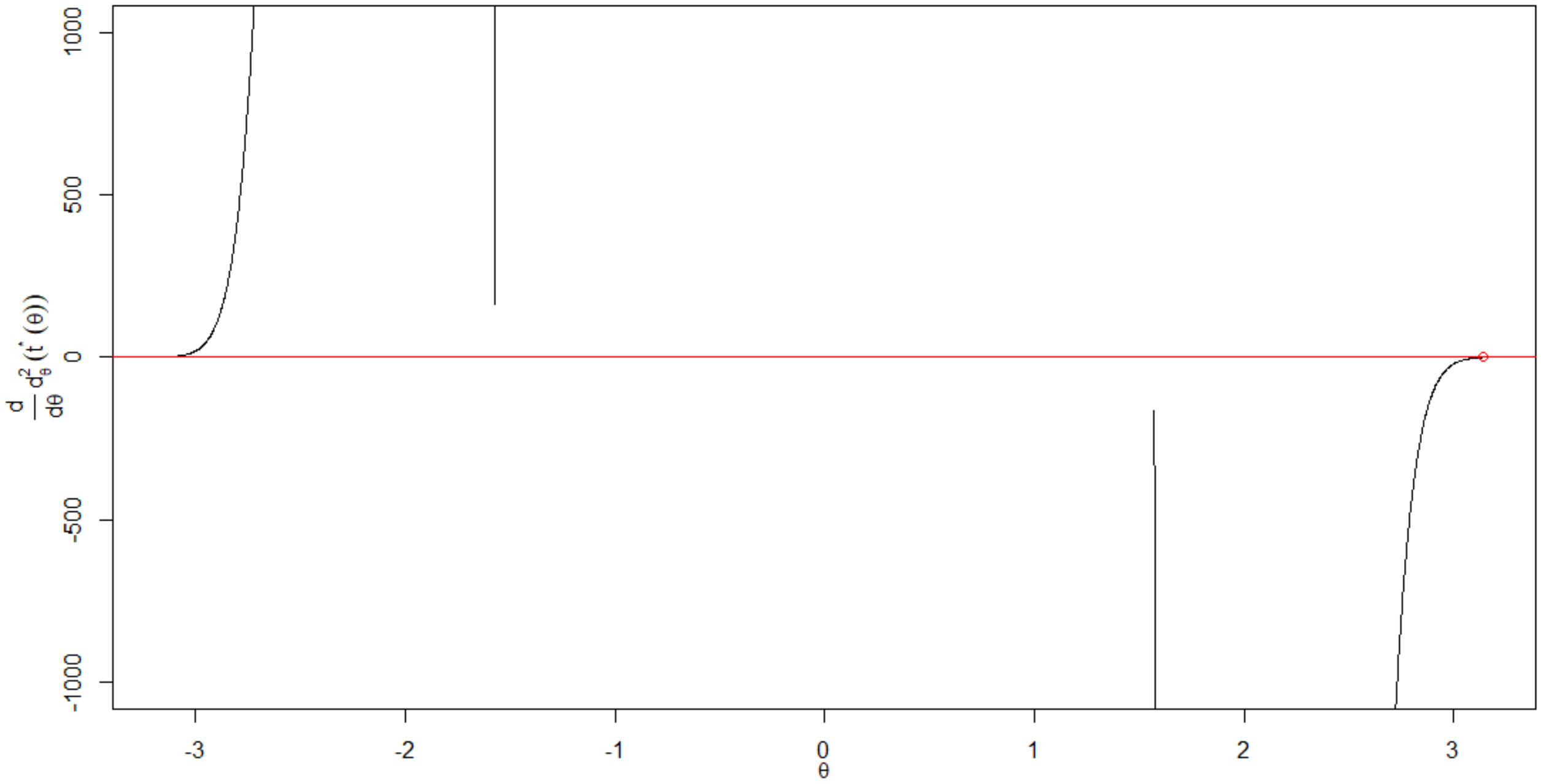

**Fig. 5.19.** Global map $\frac{d}{d\theta}d_\theta^2\big(t^*(\theta)\big)$ associated with the polynomial $x^5 - 15x^4 + 85x^3 - 225x^2 + 274x - 120$. This map was generated from $N = 5{,}000$ elements $LzC(\theta_k)$ associated with points $\theta_k$ in a regular partition of interval $[-\pi, \pi)$. Reference theta root $\theta_i^* = \pi$ for $i = 1{,}2{,}3{,}4{,}5$, is shown by means of a small circle on horizontal axis $y = 0$, just on the right edge of the proximity map.

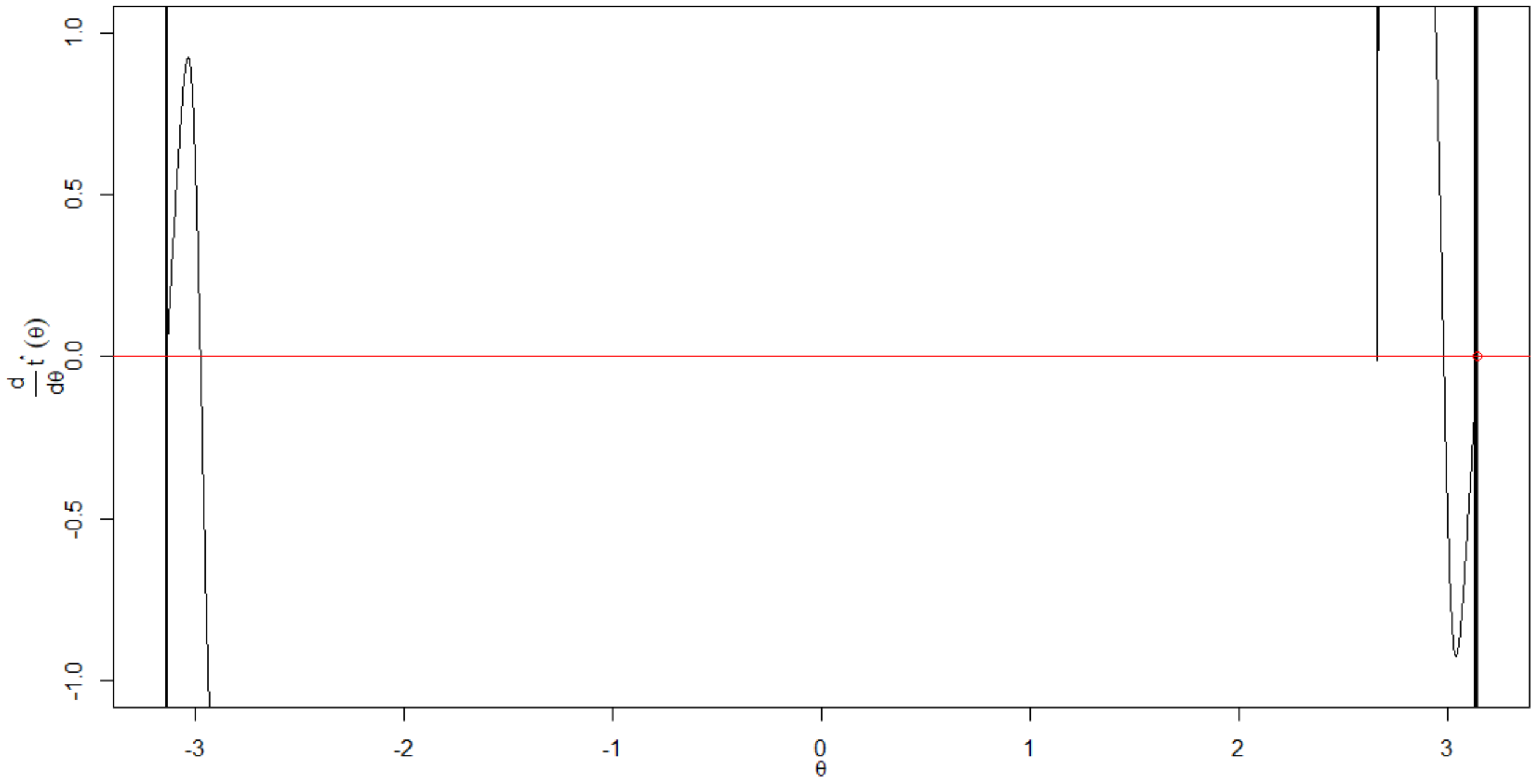

**Fig. 5.20.** Global map $\frac{d}{d\theta}t^*(\theta)$ associated with the polynomial $x^5 - 15x^4 + 85x^3 - 225x^2 + 274x - 120$. This map was generated from $N = 5{,}000$ elements $LzC(\theta_k)$ associated with points $\theta_k$ in a regular partition of interval $[-\pi, \pi)$. Reference theta root $\theta_i^* = \pi$ for $i = 1{,}2{,}3{,}4{,}5$, is shown by means of a small circle on horizontal axis $y = 0$, just on the right edge of the proximity map.

The proximity maps of figures 5.18, 5.19 and 5.20 produce the initial approximations whose numerical characteristics are shown in tables 5.9, 5.10 and 5.11. As usual, the rows in these tables are in ascending order according to variable $d^2_{\widehat{\theta}^*_i}\left(t^*(\widehat{\theta}^*_i)\right)$.

**Table 5.9**. Initial estimates of the roots for polynomial $x^5-15x^4+85x^3-225x^2+274x-120$. These estimates were obtained from map $e(\theta)$ (LC map) in figure 5.18.

| $i$ | $\widehat{R}_i$ | $\widehat{\theta}^*_i$ | $\lvert\Delta e_i\rvert$ | $d^2_{\widehat{\theta}^*_i}\left(t^*(\widehat{\theta}^*_i)\right)$ |
|---|---|---|---|---|
| 1 | $3.999981+0.007837i$ | 3.139353 | 0.04085724 | 0.0001382235 |
| 2 | $3.000087+0.012543i$ | 3.138805 | 0.04085900 | 0.0002797770 |
| 3 | $3.999955-0.012236i$ | −3.138097 | 0.04086043 | 0.0003369499 |
| 4 | $3.000185-0.018197i$ | −3.137549 | 0.04086308 | 0.0005891475 |
| 5 | $3.999897+0.018553i$ | 3.136292 | 0.04086849 | 0.0007750350 |
| 6 | $3.000255+0.021385i$ | 3.136840 | 0.04086497 | 0.0008138521 |
| 7 | $3.999330+0.047429i$ | 3.128045 | 0.04093706 | 0.0050846525 |
| 8 | $3.999275+0.049357i$ | 3.127494 | 0.04094667 | 0.0055085853 |

**Table 5.10**. Initial estimates of the roots for polynomial $x^5-15x^4+85x^3-225x^2+274x-120$. These estimates were obtained from discrete map $\frac{d}{d\theta}d^2_\theta\big(t^*(\theta)\big)$ in figure 5.19.

| $i$ | $\widehat{R}_i$ | $\widehat{\theta}^*_i$ | $\left\lvert\Delta\frac{d}{d\theta}d^2_{\widehat{\theta}^*_i}(t^*)\right\rvert$ | $d^2_{\widehat{\theta}^*_i}\left(t^*(\widehat{\theta}^*_i)\right)$ |
|---|---|---|---|---|
| 1 | $4.999998+0.000000i$ | −3.141593 | 0.06927515 | $8.171848\times10^{-11}$ |
| 2 | $3.999264+0.049728i$ | 3.127389 | 2.53939976 | $5.592121\times10^{-3}$ |
| 3 | $3.999189+0.052191i$ | 3.126686 | 1.37735679 | $6.162811\times10^{-3}$ |

**Table 5.11**. Initial estimates of the roots for polynomial $x^5-15x^4+85x^3-225x^2+274x-120$. These estimates were obtained from discrete map $\frac{d}{d\theta}t^*(\theta)$ in figure 5.20.

| $i$ | $\widehat{R}_i$ | $\widehat{\theta}^*_i$ | $\left\lvert\Delta\frac{d}{d\theta}t^*(\widehat{\theta}^*_i)\right\rvert$ | $d^2_{\widehat{\theta}^*_i}\left(t^*(\widehat{\theta}^*_i)\right)$ |
|---|---|---|---|---|
| 1 | $3.999883-0.019793i$ | −3.135938 | 795.62572466 | $8.821321\times10^{-4}$ |
| 2 | $3.999825+0.024191i$ | 3.134681 | 795.49771922 | $1.318330\times10^{-3}$ |
| 3 | $3.999100+0.054996i$ | 3.125885 | 793.57640229 | $6.847204\times10^{-3}$ |
| 4 | $3.956179+0.582776i$ | 2.978603 | 0.03454418 | $1.447769\times10^{0}$ |
| 5 | $3.956192-0.582856i$ | −2.978580 | 0.03513612 | $1.448423\times10^{0}$ |
| 6 | $4.794334+1.406538i$ | 2.662192 | 1.91318605 | $1.759392\times10^{2}$ |
| 7 | $4.793445+1.407033i$ | 2.662182 | 3.82812132 | $1.759589\times10^{2}$ |

From table 5.9, we see that with the map $e(\theta)$ from figure 5.18 it was only possible to obtain approximations of two roots, $R_4=4+0i$ and $R_3=3+0i$, which are shown in the first two rows; the remaining rows in table 5.9 are replicas of the first two. It is worth mentioning that all the approximations shown in table 5.9 were generated by structures $LzC$ close to $LzC(\theta=\pi)$.

Evidently, line $\ell_1(\pi) = 7.5 - t$, z-circumference $zC(\pi) = 120/(7.5 - t)$, and terminal curve $t\mathfrak{C}(\pi) = 274 - 225(7.5 - t) + (7.5 - t)^2(85 - 7.5^2 + t^2)$ (see expression (5.6)), are trajectories that move only on the real axis in complex plane $\mathbb{C}$; this prevents us from finding, in a usual way, intersections between z-circumference $zC(\pi)$ (of infinite radius) and terminal semi-lines of the form $tL(\pi) = [274 - 225 \cdot R_x(\pi) + R_x(\pi)^2(85 - 7.5^2)] + t^2 R_x(\pi)^2$ (see expression (5.13)), which are also trajectories that move only on the real axis in $\mathbb{C}$; moreover, $R_x(\pi)$ does not have a unique value, since all roots of equation (5.32) have the same theta root $\theta_i^* = \pi$. This does not prevent us, however, from obtaining the dynamic squared distance $d_\pi^2(t) = [t\mathfrak{C}(\pi) - zC(\pi)]^2$, whose graph is shown in figure 5.21.

From table 5.10, we see that with the map $\frac{d}{d\theta} d_\theta^2\big(t^*(\theta)\big)$ from figure 5.19, we obtain approximations to two roots, $R_5 = 5 + 0i$ and $R_4 = 4 + 0i$; these approximations are shown in the first two rows. The third row is a replica of the second.

From table 5.11, we see that with the map $\frac{d}{d\theta} t^*(\theta)$ from figure 5.20 it was only possible to approximate a single root, $R_4 = 4 + 0i$; this approximation is shown in the first row; the remaining rows are replicas of the first one, or correspond to other crossings of $\frac{d}{d\theta} t^*(\theta)$ not associated with true roots.

As for the accuracy of the best initial approximations to the root $R_4 = 4$ contained in tables 5.9, 5.10 and 5.11, we see that map $e(\theta)$ (LC map) from figure 5.18 is the one that generates the most accurate approximation; in second place is the map $\frac{d}{d\theta} t^*(\theta)$ from figure 5.20, and in third place is the map $\frac{d}{d\theta} d_\theta^2\big(t^*(\theta)\big)$ from figure 5.19. In all cases, these initial approximations can be considered reasonable, especially that of the LC map from figure 5.18. With the same LC map from figure 5.18, a good initial approximation to root $R_3 = 3$ was obtained by chance, while with map $\frac{d}{d\theta} d_\theta^2\big(t^*(\theta)\big)$ from figure 5.19, an excellent initial approximation to root $R_5 = 5$ was obtained, also by chance.

Since global maps $e(\theta)$, $\frac{d}{d\theta} d_\theta^2\big(t^*(\theta)\big)$, and $\frac{d}{d\theta} t^*(\theta)$, shown in figures 5.18, 5.19 and 5.20, are periodic functions, the valid initial approximations $\hat{\theta}_i^*$ shown in the top rows of tables 5.9, 5.10 and 5.11 are all close to $\theta = \pi$, sometimes with a positive sign and sometimes with a negative sign; this is because the true theta root $\theta_i^* = \pi$ occurs precisely at the edge of the proximity maps, and given the periodicity of these maps, it is valid to assume that $\theta_i^* = \pm\pi$.

Why was the LC method not able to approximate all the roots of equation (5.32)? The answer lies in the way the optimization method we use to find the optimal arguments $t^*$ of function $d_\theta^2(t)$ is designed. Let us look at figure 5.21, which shows an important region of function $d_{\theta=\pi}^2(t)$ associated with equation (5.32).

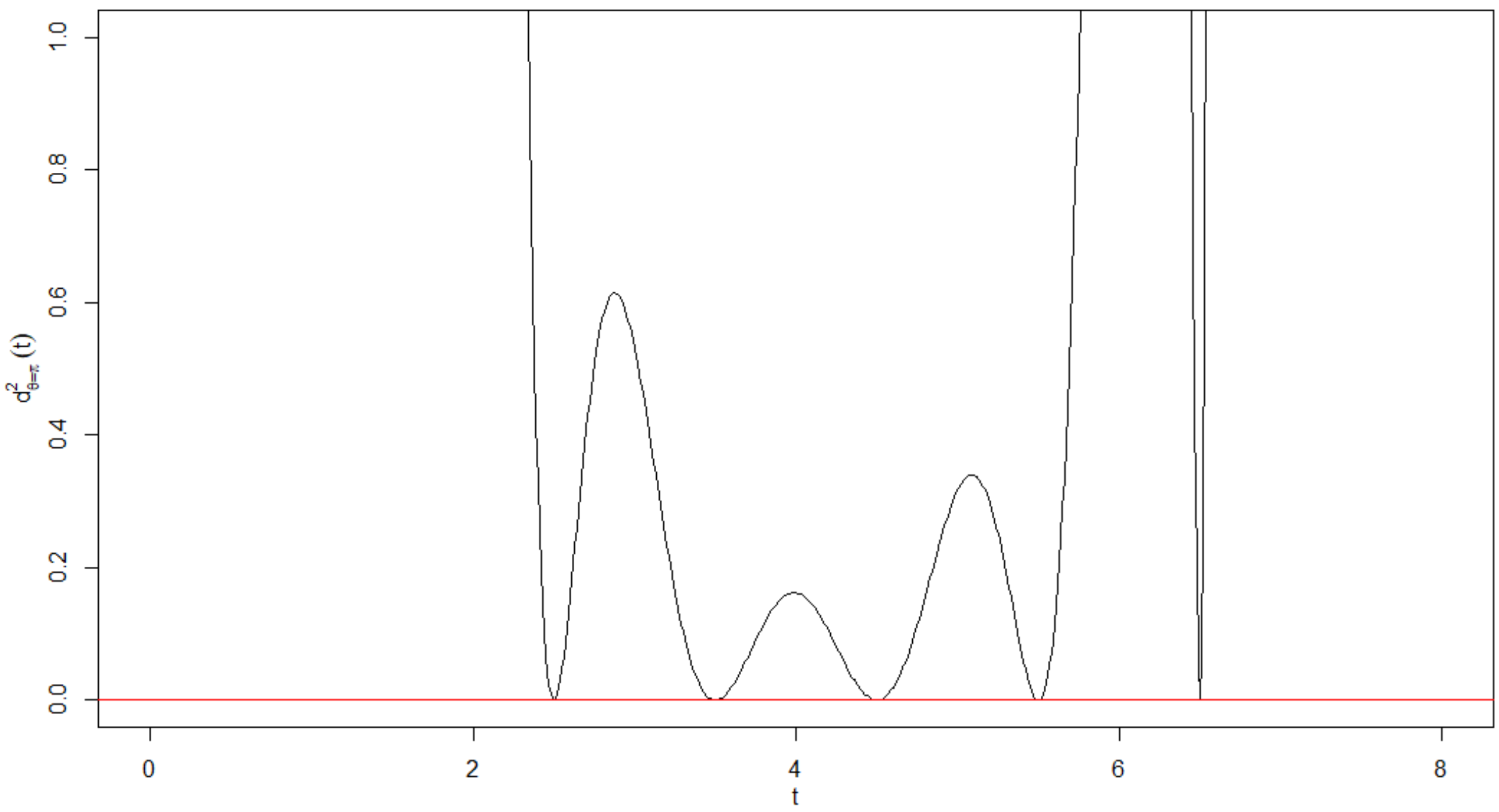

**Fig. 5.21.** Graph of function $d^2_{\theta=\pi}(t)$ associated with polynomial $x^5 - 15x^4 + 85x^3 - 225x^2 + 274x - 120$ at interval $t \in [0,8]$.

The function $d^2_{\theta=\pi}(t)$ associated with polynomial $x^5 - 15x^4 + 85x^3 - 225x^2 + 274x - 120$ contains, in fact, five minimizing arguments $t^*$, where, of course, $d^2_{\theta=\pi}(t^*) = 0$ at each of them. The optimization method we use to find $t^*$ (see function `min_D2` listed in annex 2 section 4) is not designed to find multiple minimizing arguments $t^*$ of $d^2_\theta(t)$, so it can only approximate one of the five optimal arguments $t^*$ of function $d^2_{\theta=\pi}(t)$ shown in figure 5.21. The stochastic nature of the simulated annealing (SANN) algorithm used in the first phase of our optimization method is the reason why it was possible to obtain initial approximations of more than one root (but not all roots) by means of the proximity maps shown in figures 5.18, 5.19 and 5.20.

**A generalizable observation about the graph in figure 5.21**; we can see that $d^2_{\theta=\pi}(t)$ behaves as a continuous and differentiable function, so it is possible to approximate its global minima by calculating a discrete approximation to its corresponding derivative $\frac{d}{dt}d^2_{\theta=\pi}(t)$, and then find the crossings of this derivative with horizontal axis $y = 0$, in an analogous way to how we have been proceeding when finding numerical approximations (crossings with $y = 0$) associated with maps $e(\theta)$, $\frac{d}{d\theta}d^2_\theta(t^*(\theta))$, and $\frac{d}{d\theta}t^*(\theta)$. All this suggests that, in principle, it is possible to design an alternative algorithm to find optimal arguments $t^*$ of any dynamic squared distance function $d^2_\theta(t)$, based on discrete approximations of the derivative $\frac{d}{dt}d^2_\theta(t)$, and considering the possibility of finding more than one minimizing argument $t^*$. This alternative is proposed as future research work, and promises to obtain a more efficient, faster, and more complete optimization algorithm, adapted to the problem of root approximation for polynomials of a single variable that we are considering. For now, if we obtain a discrete approximation to the derivative $\frac{d}{dt}d^2_{\theta=\pi}(t)$ of function $d^2_{\theta=\pi}(t)$ in figure 5.21 (figure 5.22 shows the graph of this derivative), it is possible to obtain, from this discrete derivative, approximations to the optimal arguments $t^*$ of $d^2_{\theta=\pi}(t)$; in fact, these optimal arguments occur at values $t^*$ equal to 2.5, 3.5, 4.5, 5.5 and 6.5. With these optimal arguments $t^*$, direction vector $v_\pi = -1 + 0i$, and fixed point $P_1 = 7.5$, we obtain the roots $R_i = P_1 + t^* v_\pi$, which coincide with reference roots of equation (5.32).

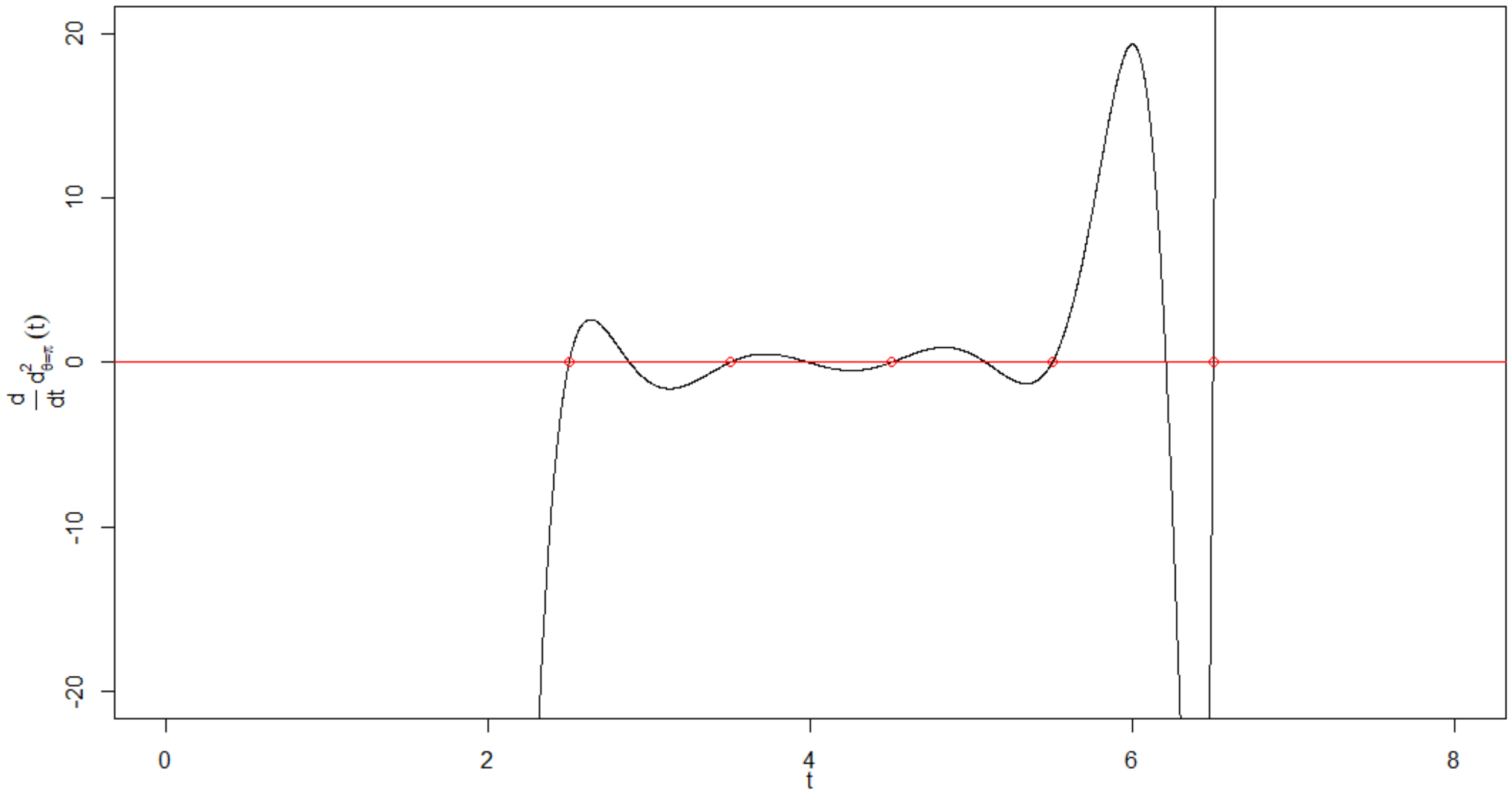

**Fig. 5.22.** Graph of function $\frac{d}{dt}d^2_{\theta=\pi}(t)$ associated with polynomial $x^5 - 15x^4 + 85x^3 - 225x^2 + 274x - 120$ at interval $t \in [0,8]$. The global minima of the corresponding function $d^2_{\theta=\pi}(t)$ from figure 5.21 are located at the "increasing" crossings of $\frac{d}{dt}d^2_{\theta=\pi}(t)$ with horizontal axis $y = 0$, which are identified in this figure by means of small circles.

In order to complete the picture of the results obtained until this moment, let us analyze some additional graphs related to the family of functions $d^2_\theta(t)$ associated with equation (5.32).

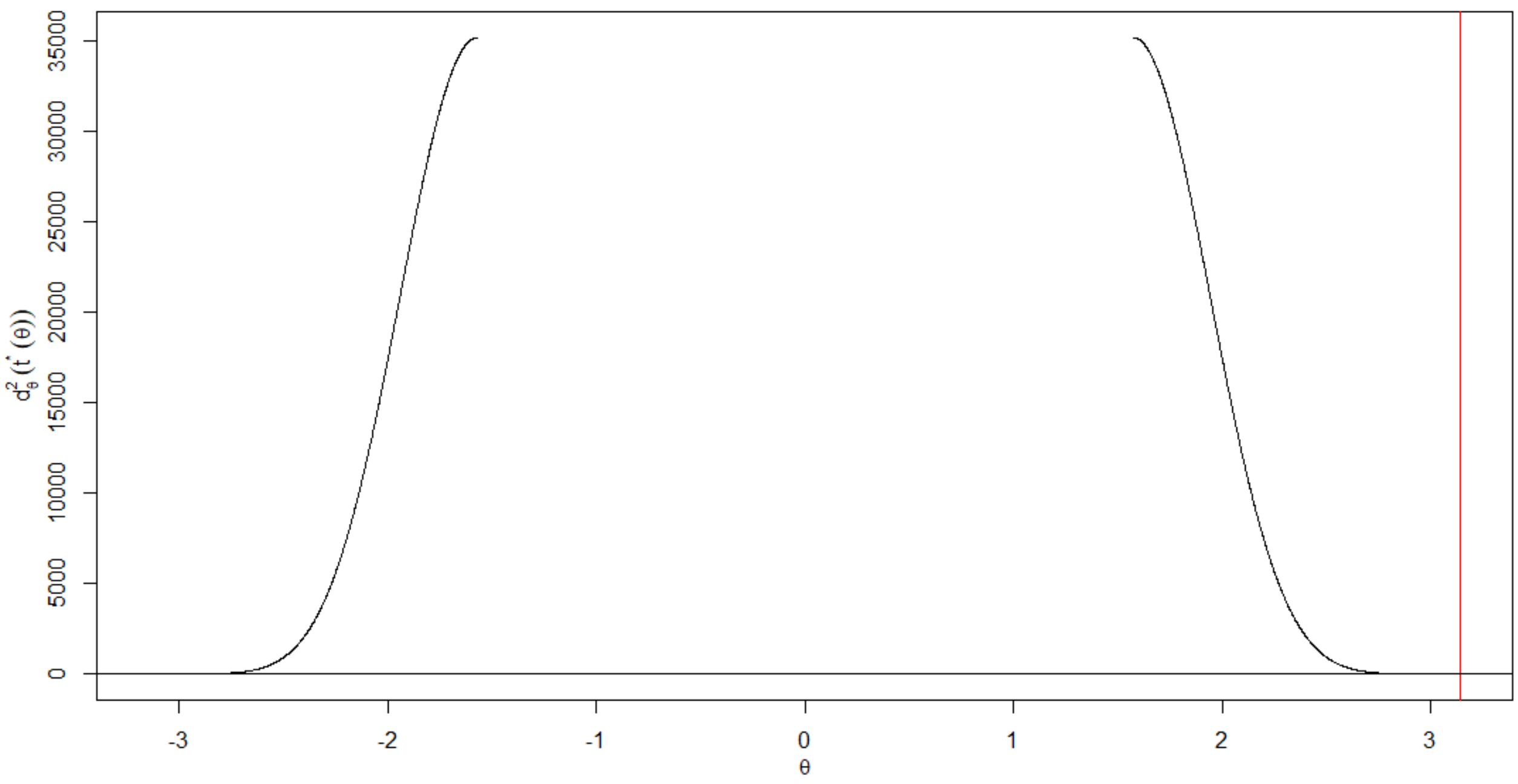

**Fig. 5.23.** Global minima $d_\theta^2\big(t^*(\theta)\big)$ vs. $\theta$ associated to equation $x^5 - 15x^4 + 85x^3 - 225x^2 + 274x - 120 = 0$. This graph includes a vertical line segment that marks the location of the true theta root $\theta_i^* = \pi$ for $i = 1,2,3,4,5$.

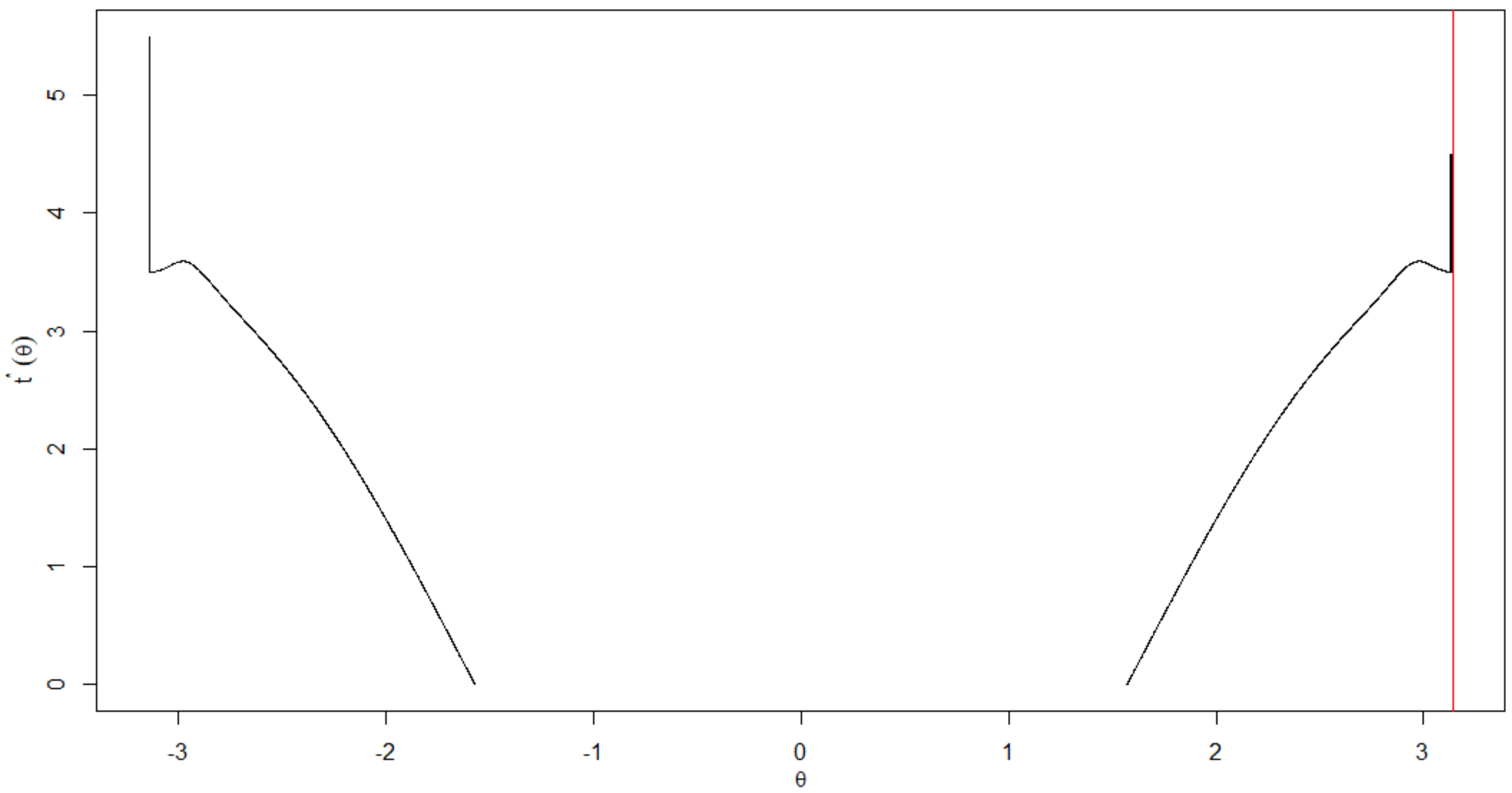

**Fig. 5.24.** Minimizing arguments $t^*(\theta)$ vs. $\theta$ associated to equation $x^5 - 15x^4 + 85x^3 - 225x^2 + 274x - 120 = 0$. This graph includes a vertical line segment that marks the location of the true theta root $\theta_i^* = \pi$ for $i = 1,2,3,4,5$.

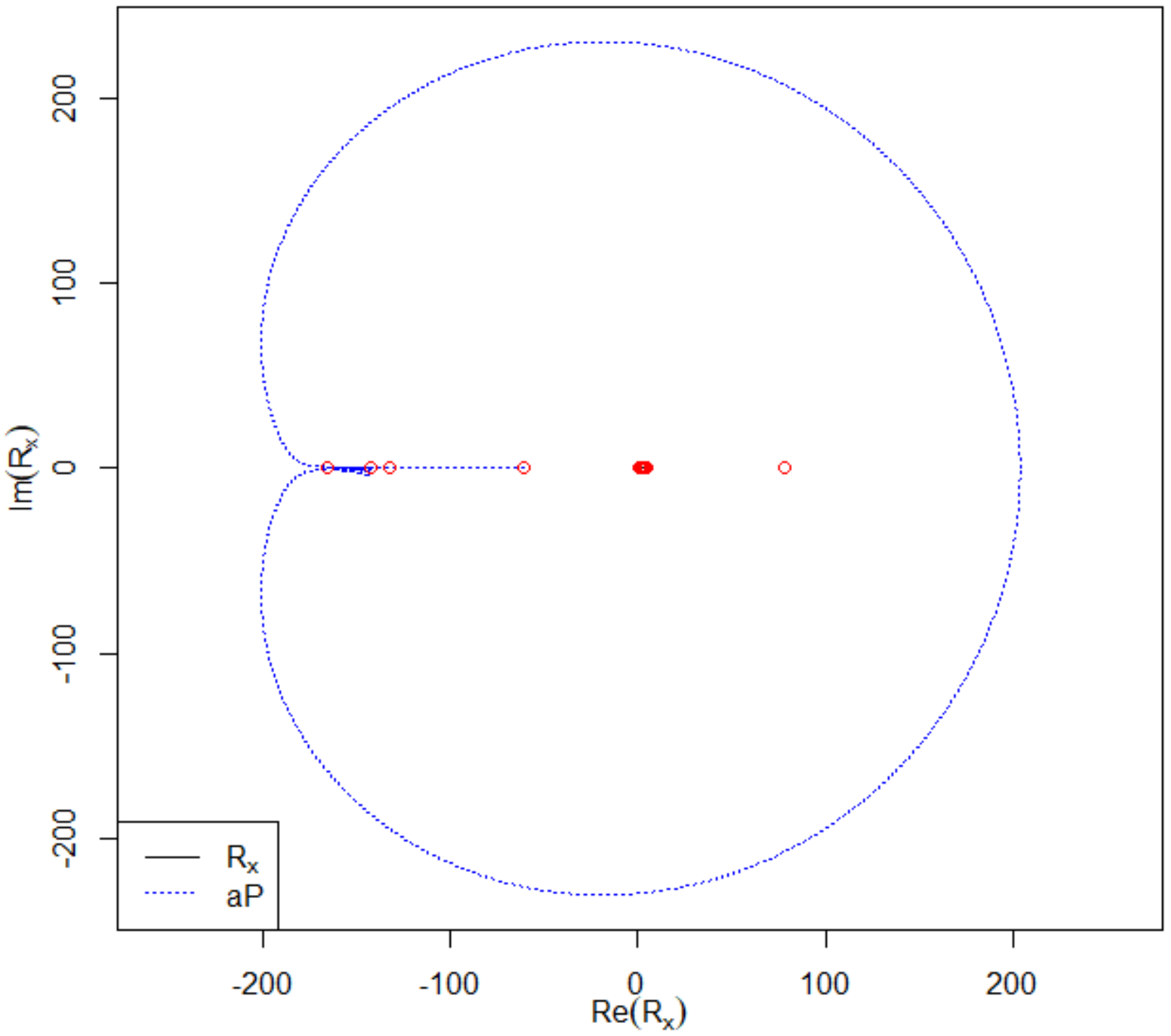

**Fig. 5.25.** Trajectory $R_x(\theta)$ (small circles clustered near the origin, in a "microscopic" view) and trajectory $aP(\theta)$ associated with equation $x^5 - 15x^4 + 85x^3 - 225x^2 + 274x - 120 = 0$. The small circles on $aP(\theta)$ are anchor points associated with the true roots.

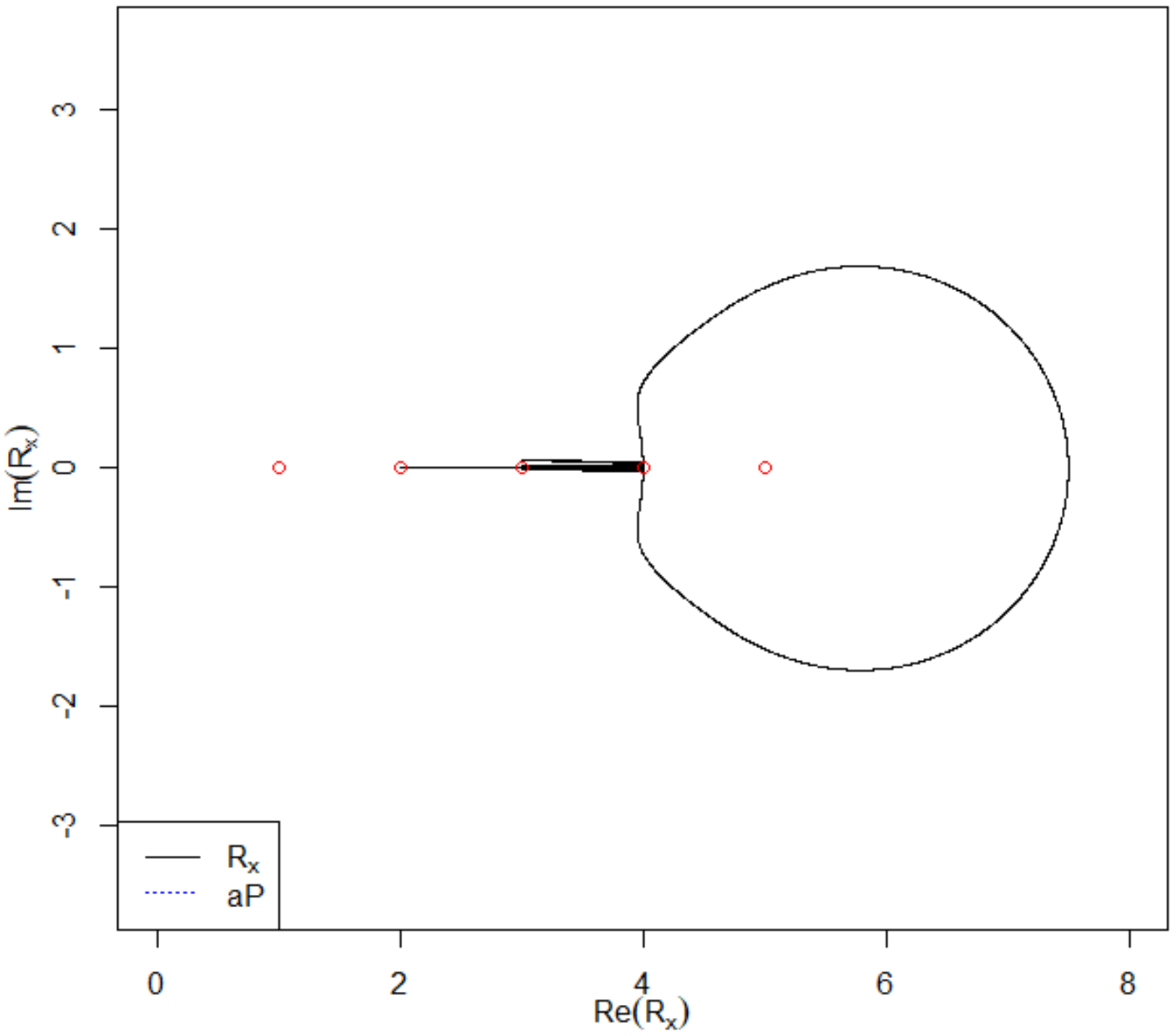

**Fig. 5.26.** Close-up of figure 5.25, showing in better detail the trajectory $R_x(\theta)$ associated with equation $x^5 - 15x^4 + 85x^3 - 225x^2 + 274x - 120 = 0$. Small circles are the true roots $R_i$.

Figures 5.23, 5.24, 5.25 and 5.26 show that each one of the maps $\frac{d}{d\theta}d_\theta^2\big(t^*(\theta)\big)$ and $\frac{d}{d\theta}t^*(\theta)$ from figures 5.19 and 5.20 are actually made up of a single continuous section (this can be verified visually if we extend the vertical range in the graphs from figures 5.19 and 5.20). Figures 5.23 and 5.25 show that the quantities $d_\theta^2\big(t^*(\theta)\big)$ and $|aP(\theta)|$ reach relatively large values in this case. On the other hand, figures 5.24 and 5.26, together with table 5.11, show that the discrete values of map $\frac{d}{d\theta}t^*(\theta)$ from figure 5.20 involved in the best approximation to root $R_4 = 4 + 0i$, present very large differences in their functional values; this means that at the edge $\theta = \pm\pi$ (where all the roots of equation (5.32) are located) we have difficulties in obtaining a well-defined proximity map $\frac{d}{d\theta}t^*(\theta)$, due to the problems associated with our optimization algorithm used in approximating optimal arguments $t^*(\theta)$. Figures 5.25 and 5.26 show that the trajectories $aP(\theta)$ and $R_x(\theta)$ are not similar to each other, since, in this case, $aP(\theta)$ does not depend linearly on $R_x(\theta)$ (see expressions in (5.13)).

## Variable change, and root-shifting in equation (5.32)

Carrying on with the analysis of this numerical example, we will now make a variable change in equation (5.32), in order to approximate all its roots via the LC method, but this time avoiding the problem of multiple theta roots. According to the result in annex 5 section 2, if we make the variable change $z = x + a$, where $x$ is the variable term in equation (5.32), and $a = 0 - 2i$ is a displacement (or shifting) of variable $x$, we will obtain the equation

$$z^5 + (-15 + 10i)z^4 + (45 - 120i)z^3 + (135 + 430i)z^2 + (-666 - 420i)z + (540 - 100i) = 0, \qquad (5.33)$$

whose reference roots $r_i$ and $\vartheta_i^*$ are

$$\begin{aligned}
r_1 &= 1 - 2i, & \vartheta_1^* &= 2.709185,\\
r_2 &= 2 - 2i, & \vartheta_2^* &= 2.642246,\\
r_3 &= 3 - 2i, & \vartheta_3^* &= 2.553590, \qquad (5.34)\\
r_4 &= 4 - 2i, & \vartheta_4^* &= 2.432966,\\
r_5 &= 5 - 2i, & \vartheta_5^* &= 2.265535.
\end{aligned}$$

This time, theta root $\vartheta_i^*$ is the inclination angle of line $\ell_1(\vartheta_i^*)$ with fixed point $P_1 = 7.5 - 5i$ containing root $r_i$. As we can see, now the fixed point $P_1$ is located below and to the right of roots $r_1$, $r_2$, $r_3$, $r_4$ $r_5$, which are all located on horizontal line $x - 2i$, so $\ell_1(\vartheta)$ rotates within angular interval $(\pi/2, \pi)$ as it "sweeps" the roots $r_i$. In this case, the root-shifting choice $a = -2i$ is arbitrary, and does not follow any strategy, except that of preventing the roots from being located on the real axis of the complex plane $\mathbb{C}$, as it was the case in the first part of this numerical example. The operating parameters for reproducing the results shown below are specified in annex 3 section 3, within subsection “specific conditions for reproducing results from the examples in chapter 5”.

Figure 5.27 shows the "semi-global" map (region from 0 to $\pi$) $e(\vartheta)$ associated with equation (5.33); additionally, figure 5.28 shows a close-up of this same map $e(\vartheta)$, narrowing down the vertical range in order to better appreciate the behavior of weighted errors in the vicinity of theta roots $\vartheta_i^*$. Figures 5.29 and 5.30 show, respectively, the "semi-global" maps $\frac{d}{d\vartheta}d_\vartheta^2\big(t^*(\vartheta)\big)$ and $\frac{d}{d\vartheta}t^*(\vartheta)$. Due to technical limitations, the text within these graphs (as well as within the remaining graphs in this numerical example) shows the angular variable $\theta$, but we actually refer to the angular variable $\vartheta$ associated with the variable change $z = x + a$. Tables 5.12, 5.13 and 5.14 show the numerical characteristics of the smooth crossings of these semi-global maps with vertical axis $y = 0$, which serve as initial estimates to the roots of equation (5.33).

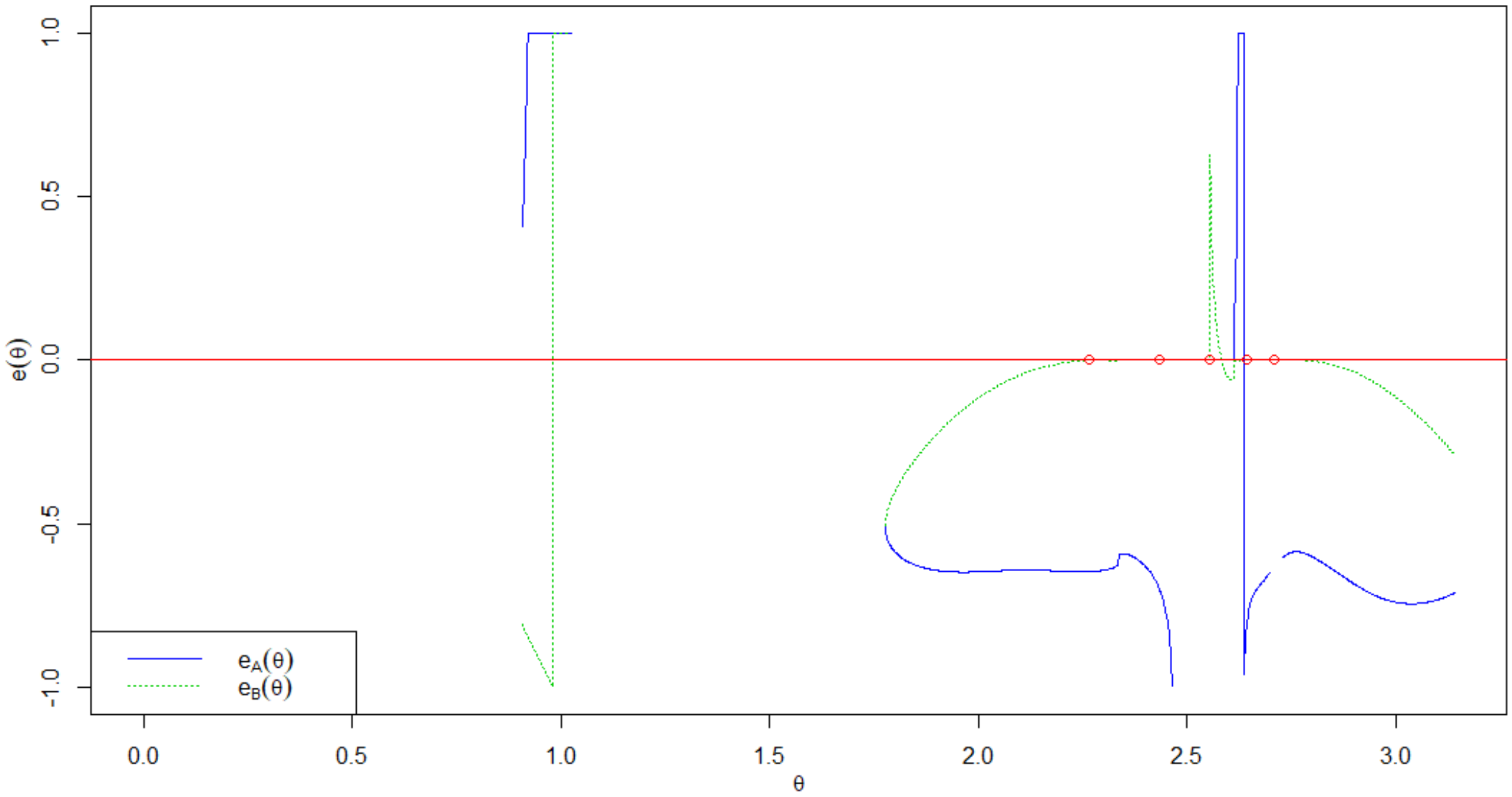

**Fig. 5.27.** Discrete angular proximity map $e(\vartheta)$ for the polynomial of degree $n = 5$ in equation (5.33). This map was generated from $N = 2{,}500$ elements $LzC(\vartheta_k)$ associated with points $\vartheta_k$ in a regular partition of interval $[0, \pi)$. Reference theta roots $\vartheta_i^*$ in (5.34) are shown here by means of small circles on horizontal axis $y = 0$.

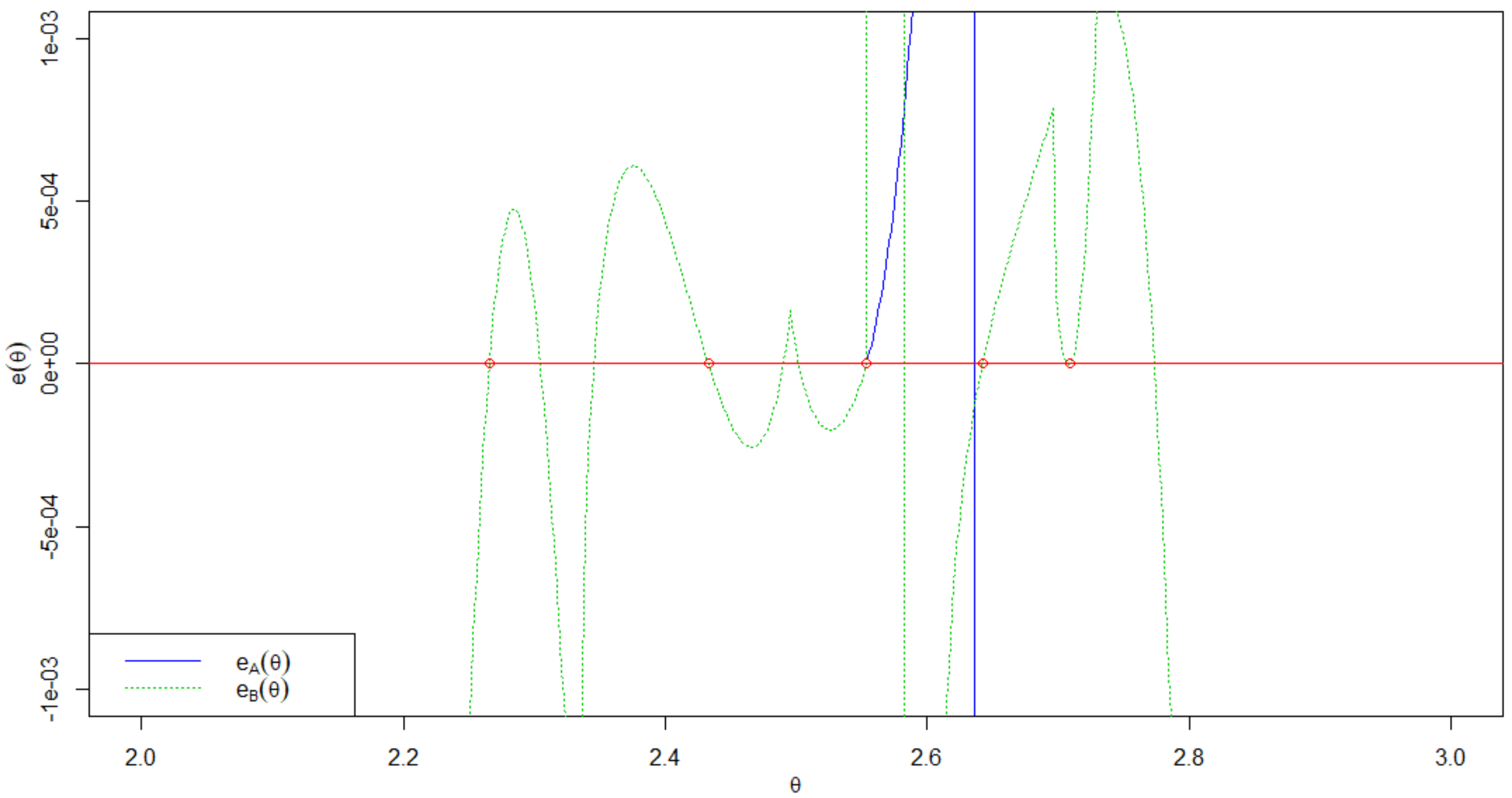

**Fig. 5.28.** Close-up of map $e(\vartheta)$ from figure 5.27 to region $[2, 3)$, narrowing down the vertical axis range to interval $[-0.001, 0.001]$, in order to better appreciate the map's behavior in the vicinity of theta roots $\vartheta_i^*$.

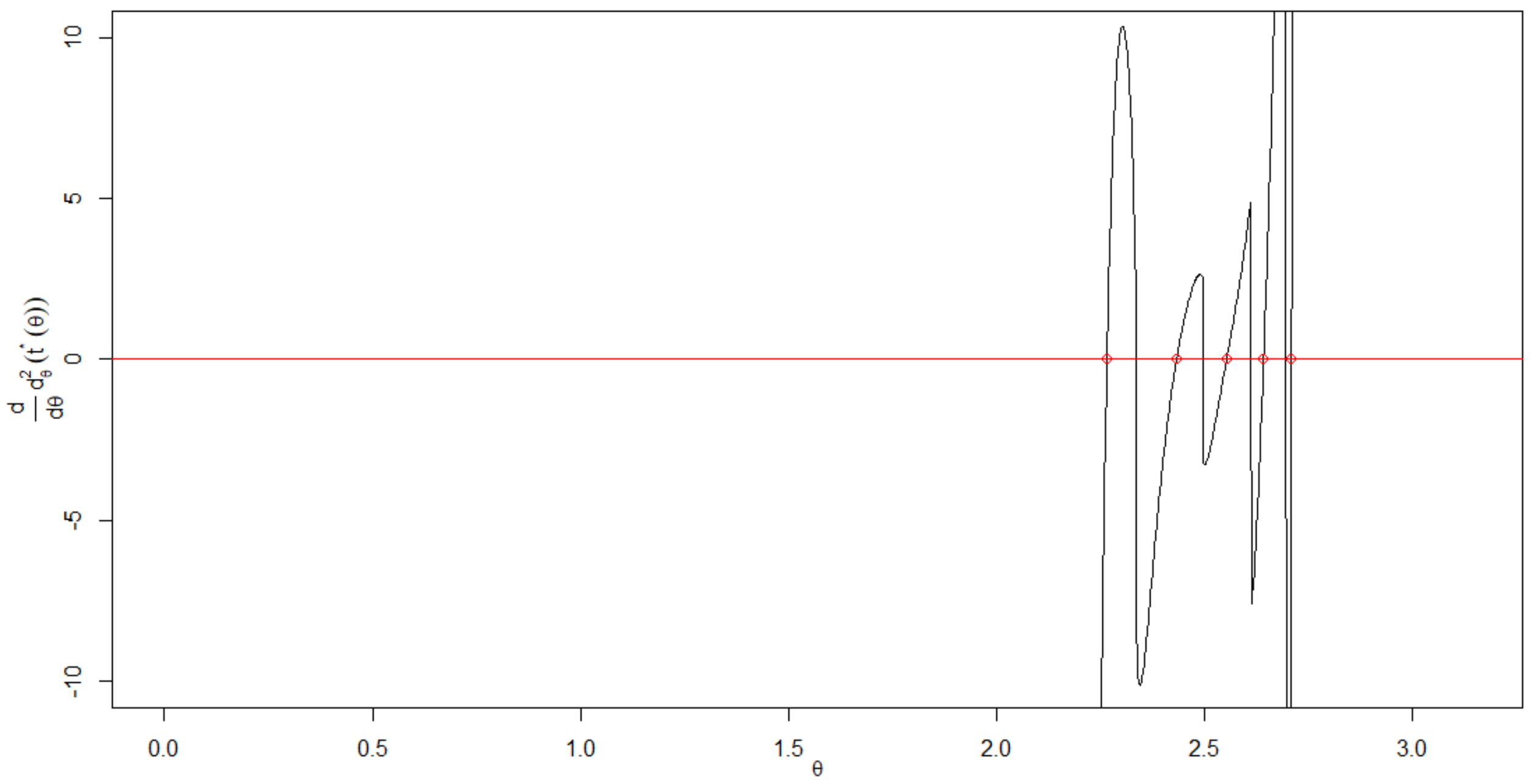

**Fig. 5.29.** Map $\frac{d}{d\vartheta} d_\vartheta^2(t^*(\vartheta))$ corresponding to the polynomial of degree $n = 5$ in equation (5.33). This map was generated from $N = 2{,}500$ elements $LzC(\vartheta_k)$ associated with points $\vartheta_k$ in a regular partition of interval $[0, \pi)$. Reference theta roots $\vartheta_i^*$ in (5.34) are shown in this graph by means of small circles on horizontal axis $y = 0$.

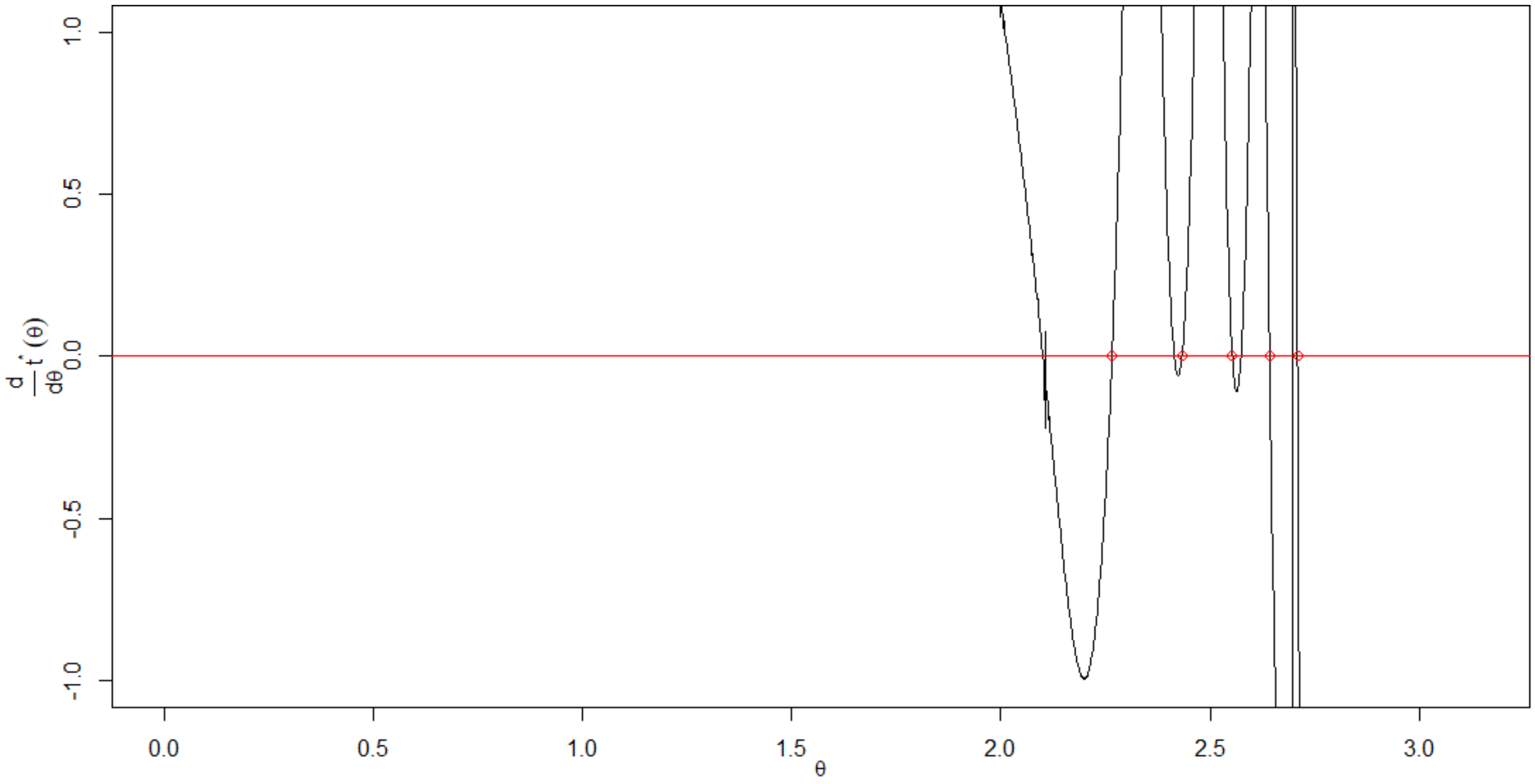

**Fig. 5.30.** Map $\frac{d}{d\vartheta}t^*(\vartheta)$ corresponding to the polynomial of degree $n = 5$ in equation (5.33). This map was generated from $N = 2{,}500$ elements $LzC(\vartheta_k)$ associated with points $\vartheta_k$ in a regular partition of interval $[0, \pi)$. Reference theta roots $\vartheta_i^*$ in (5.34) are shown in this graph by means of small circles on horizontal axis $y = 0$.

**Table 5.12**. Initial estimates of the roots for the polynomial of degree $n = 5$ in equation (5.33). These estimates were obtained from the map $e(\vartheta)$ (LC map) in figure 5.27.

| $i$ | $\hat{r}_i$ | $\hat{\vartheta}_i^*$ | $\lvert\Delta e_i\rvert$ | $d^2_{\hat{\vartheta}_i^*}\left(t^*(\hat{\vartheta}_i^*)\right)$ |
|---|---|---|---|---|
| 1 | $3.999997 - 2.000003i$ | 2.432967 | $1.597593 \times 10^{-5}$ | $4.066744 \times 10^{-11}$ |
| 2 | $1.999987 - 2.000026i$ | 2.642251 | $2.477056 \times 10^{-5}$ | $3.887469 \times 10^{-9}$ |
| 3 | $4.999983 - 2.000013i$ | 2.265540 | $6.639185 \times 10^{-5}$ | $8.960825 \times 10^{-9}$ |
| 4 | $3.000311 - 1.999534i$ | 2.553487 | $6.265772 \times 10^{-1}$ | $3.859252 \times 10^{-7}$ |
| 5 | $1.000346 - 1.999257i$ | 2.709070 | $3.589863 \times 10^{-6}$ | $7.726386 \times 10^{-5}$ |
| 6 | $1.006975 - 1.986251i$ | 2.707033 | $3.580384 \times 10^{-6}$ | $2.682199 \times 10^{-2}$ |
| 7 | $2.914696 - 2.132925i$ | 2.582796 | $6.904884 \times 10^{-3}$ | $3.141582 \times 10^{-2}$ |
| 8 | $3.238375 - 1.827287i$ | 2.501631 | $2.483913 \times 10^{-5}$ | $9.437074 \times 10^{-2}$ |
| 9 | $3.772182 - 2.157763i$ | 2.490177 | $3.049956 \times 10^{-5}$ | $9.446098 \times 10^{-2}$ |
| 10 | $4.863209 - 2.075243i$ | 2.304463 | $5.515150 \times 10^{-5}$ | $2.787635 \times 10^{-1}$ |
| 11 | $4.367306 - 1.797707i$ | 2.345209 | $6.955586 \times 10^{-5}$ | $4.268703 \times 10^{-1}$ |
| 12 | $2.054934 - 2.899006i$ | 2.773341 | $8.541956 \times 10^{-5}$ | $9.141682 \times 10^{+0}$ |

**Table 5.13**. Initial estimates of the roots for the polynomial of degree $n = 5$ in equation (5.33). These estimates were obtained from the discrete map $\frac{d}{d\vartheta} d_{\vartheta}^2\left(t^*(\vartheta)\right)$ in figure 5.29.

| $i$ | $\hat{r}_i$ | $\hat{\vartheta}_i^*$ | $\left\|\Delta \frac{d}{d\vartheta} d_{\hat{\vartheta}_i^*}^2(t^*)\right\|$ | $d_{\hat{\vartheta}_i^*}^2\left(t^*(\hat{\vartheta}_i^*)\right)$ |
|---|---|---|---|---|
| 1 | $3.000000 - 2.000000i$ | 2.553590 | 0.09048086 | $8.519054 \times 10^{-14}$ |
| 2 | $3.999989 - 2.000012i$ | 2.432970 | 0.09628729 | $4.574795 \times 10^{-10}$ |
| 3 | $2.000006 - 1.999991i$ | 2.642244 | 0.44620845 | $5.468520 \times 10^{-10}$ |
| 4 | $4.999973 - 2.000021i$ | 2.265543 | 0.75642471 | $2.337559 \times 10^{-8}$ |
| 5 | $1.000020 - 1.999958i$ | 2.709178 | 14.88034964 | $2.461564 \times 10^{-7}$ |
| 6 | $3.723976 - 2.158759i$ | 2.496530 | 3.92112872 | $1.109228 \times 10^{-1}$ |
| 7 | $2.803418 - 2.250919i$ | 2.612025 | 8.29179484 | $1.348901 \times 10^{-1}$ |
| 8 | $4.650651 - 2.027396i$ | 2.335027 | 3.50054256 | $5.164211 \times 10^{-1}$ |
| 9 | $1.918176 - 2.338445i$ | 2.696656 | 71.14956798 | $6.973010 \times 10^{-1}$ |

**Table 5.14**. Initial estimates of the roots for the polynomial of degree $n = 5$ in equation (5.33). These estimates were obtained from the discrete map $\frac{d}{d\vartheta} t^*(\vartheta)$ in figure 5.30.

| $i$ | $\hat{r}_i$ | $\hat{\vartheta}_i^*$ | $\left\|\Delta \frac{d}{d\vartheta} t^*(\hat{\vartheta}_i^*)\right\|$ | $d_{\hat{\vartheta}_i^*}^2\left(t^*(\hat{\vartheta}_i^*)\right)$ |
|---|---|---|---|---|
| 1 | $4.999990 - 2.000007i$ | 2.265538 | 0.04193443 | $3.093276 \times 10^{-9}$ |
| 2 | $1.999985 - 2.000029i$ | 2.642251 | 0.11186867 | $4.882548 \times 10^{-9}$ |
| 3 | $4.000059 - 1.999930i$ | 2.432947 | 0.01581659 | $1.511718 \times 10^{-8}$ |
| 4 | $2.999921 - 2.000118i$ | 2.553616 | 0.02748934 | $2.492832 \times 10^{-8}$ |
| 5 | $0.999990 - 2.000026i$ | 2.709189 | 0.30051235 | $8.962769 \times 10^{-8}$ |
| 6 | $4.055371 - 1.935467i$ | 2.414522 | 0.01725878 | $1.417395 \times 10^{-2}$ |
| 7 | $2.940950 - 2.093276i$ | 2.574004 | 0.02538819 | $1.516186 \times 10^{-2}$ |
| 8 | $5.451480 - 1.555441i$ | 2.107318 | 0.17330129 | $4.330904 \times 10^{+1}$ |
| 9 | $5.454554 - 1.553588i$ | 2.106423 | 0.29801650 | $4.421650 \times 10^{+1}$ |
| 10 | $5.474486 - 1.541583i$ | 2.100617 | 0.01609811 | $5.049160 \times 10^{+1}$ |

From tables 5.12, 5.13, and 5.14, which, as usual, are in ascending order according to the values in column $d_{\hat{\vartheta}_i^*}^2\left(t^*(\hat{\vartheta}_i^*)\right)$, we can easily see that, in each of them, the first five rows contain very good initial approximations to each of the reference roots in (5.34). If we want to associate these estimates with equation (5.32), we just need to apply the inverse transformation $x = z - a$, which implies adding the amount $2i$ to the approximations $\hat{r}_i$ shown in tables 5.12, 5.13 and 5.14, with which we would obtain initial estimates $\hat{R}_i$. In this case, to simplify the analysis and not create more tables, it is enough to verify that the initial approximations $\hat{r}_i$ in tables 5.12, 5.13 and 5.14 have an imaginary part close to $-2$.

Table 5.15 shows a summary of the approximation errors associated with the initial estimates contained in the first five rows of tables 5.12, 5.13 and 5.14. As we can see, if we rely on the criterion $d_{\hat{\vartheta}_i^*}^2\left(t^*(\hat{\vartheta}_i^*)\right)$, this time the map $\frac{d}{d\vartheta} t^*(\vartheta)$ is the one that generates on average the most accurate and consistent (least variable) initial estimates, followed very closely by map $\frac{d}{d\vartheta}d_\vartheta^2\big(t^*(\vartheta)\big)$, and in third place by map $e(\vartheta)$. As for the average magnitudes of the relative errors of the initial estimates with respect to the reference roots in (5.34), the map with the most accurate and consistent initial estimates is $\frac{d}{d\vartheta}d_\vartheta^2\big(t^*(\vartheta)\big)$, followed in second place by map $\frac{d}{d\vartheta}t^*(\vartheta)$, and in third place by map $e(\vartheta)$. As we mentioned before, with any of the three proximity maps $e(\vartheta)$, $\frac{d}{d\vartheta}d_\vartheta^2\big(t^*(\vartheta)\big)$, $\frac{d}{d\vartheta}t^*(\vartheta)$, we obtain good initial estimates of the roots of equation (5.33).

**Table 5.15**. Arithmetic means and standard deviations for error measures associated with valid initial estimates (first five rows in each of tables 5.12, 5.13 and 5.14) of the roots of the polynomial of degree $n = 5$ in equation (5.33), obtained by means of proximity maps $e(\vartheta)$, $\frac{d}{d\vartheta}d_\vartheta^2\big(t^*(\vartheta)\big)$, and $\frac{d}{d\vartheta}t^*(\vartheta)$ from figures 5.27, 5.29, and 5.30.

| Map | Error measure | | | | | |
|---|---|---|---|---|---|---|
| | $d^2_{\vartheta_i^*}\left(t^*(\hat{\vartheta}_i^*)\right)$ | | $\left\|\frac{R-\hat{R}}{R}\right\|$ | | $\left\|\frac{\vartheta^*-\hat{\vartheta}^*}{\vartheta^*}\right\|$ | |
| | **Mean** | **Standard deviation** | **Mean** | **Standard deviation** | **Mean** | **Standard deviation** |
| $e(\vartheta)$ | $1.553253 \times 10^{-5}$ | $3.450925 \times 10^{-5}$ | $1.074271 \times 10^{-4}$ | $1.587551 \times 10^{-4}$ | $1.747735 \times 10^{-5}$ | $2.185473 \times 10^{-5}$ |
| $\frac{d}{d\vartheta}d_\vartheta^2(t^*(\vartheta))$ | $5.410728 \times 10^{-8}$ | $1.078215 \times 10^{-7}$ | $6.915851 \times 10^{-6}$ | $8.011338 \times 10^{-6}$ | $1.668728 \times 10^{-6}$ | $1.519372 \times 10^{-6}$ |
| $\frac{d}{d\vartheta}t^*(\vartheta)$ | $2.752980 \times 10^{-8}$ | $3.580164 \times 10^{-8}$ | $1.728000 \times 10^{-5}$ | $1.397772 \times 10^{-5}$ | $4.660625 \times 10^{-6}$ | $4.252334 \times 10^{-6}$ |

To conclude this example, we will display the graphs related to the family of dynamic squared distances $d_\vartheta^2$ as a function of their optimal arguments $t^*(\vartheta)$. The graphs in figures 5.31, 5.32, 5.33 and 5.34, together with the graphs in figures 5.29 and 5.30, show that the maps associated with $d_\vartheta^2$ in this case are made up of a single continuous section, within region $[0,\pi)$. In accordance with what was initially stated in numerical example 4.1 of chapter 4, figures 5.30 and 5.32 suggest that theta roots $\vartheta_i^*$ correspond to two local minima, and to three local maxima of mapping $t^*(\vartheta)$; similarly, figures 5.29 and 5.31 suggest that theta roots $\vartheta_i^*$ occur at the five global minima of mapping $d_\vartheta^2\big(t^*(\vartheta)\big)$. Figures 5.33 and 5.34 show that trajectories $r_x(\vartheta)$ and $aP(\vartheta)$ are not similar in shape, but they are in behavior. Note that trajectory $aP(\vartheta)$ requires much more space to evolve, compared with trajectory $r_x(\vartheta)$.

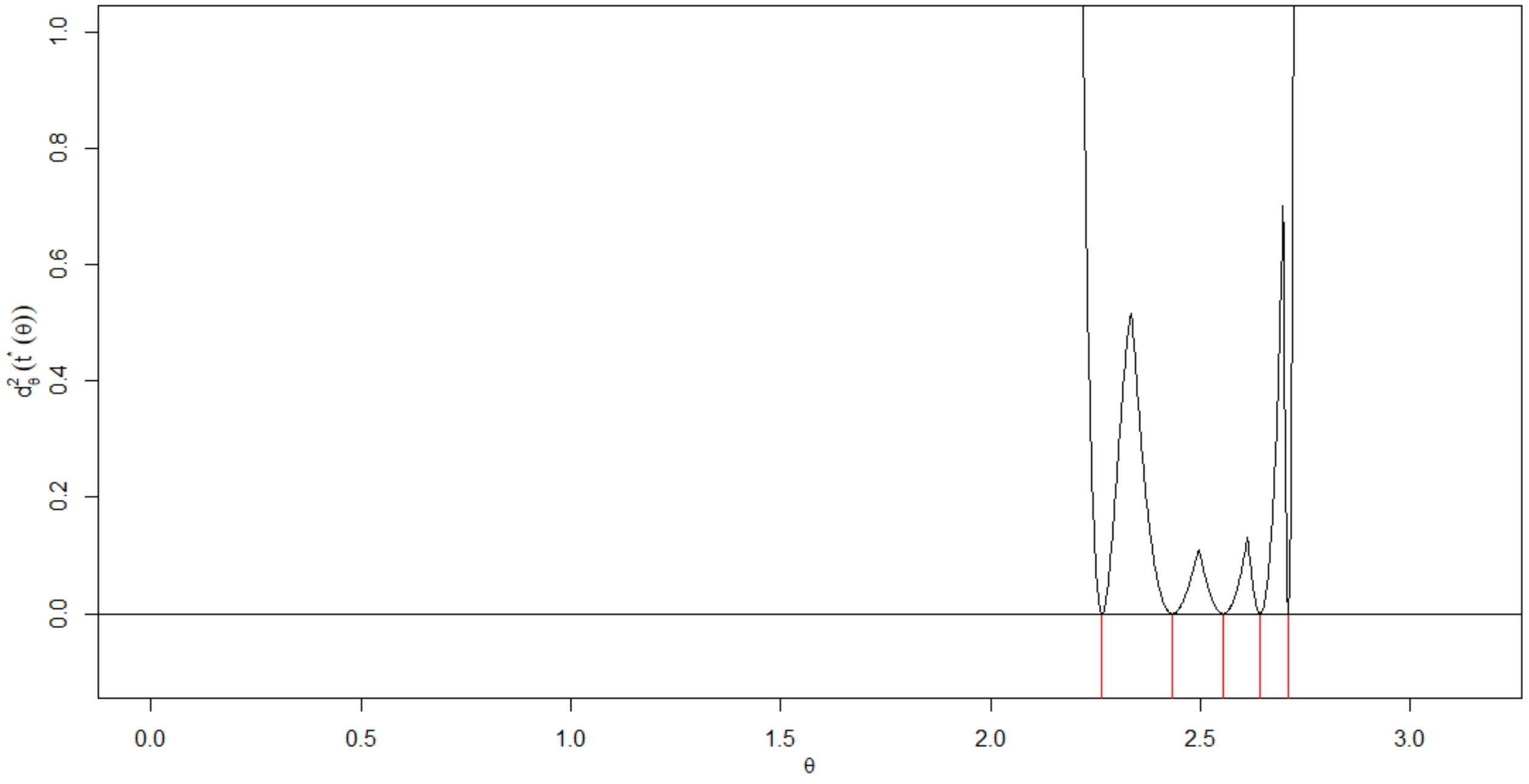

**Fig. 5.31.** Global minima $d_\vartheta^2\big(t^*(\vartheta)\big)$ vs. $\vartheta$ associated with the polynomial of degree $n = 5$ in equation (5.33). This map was generated from $N = 2{,}500$ elements $LzC(\vartheta_k)$ associated with points $\vartheta_k$ in a regular partition of interval $[0, \pi)$. This graph includes vertical line segments that indicate the location of true theta roots $\vartheta_i^*$.

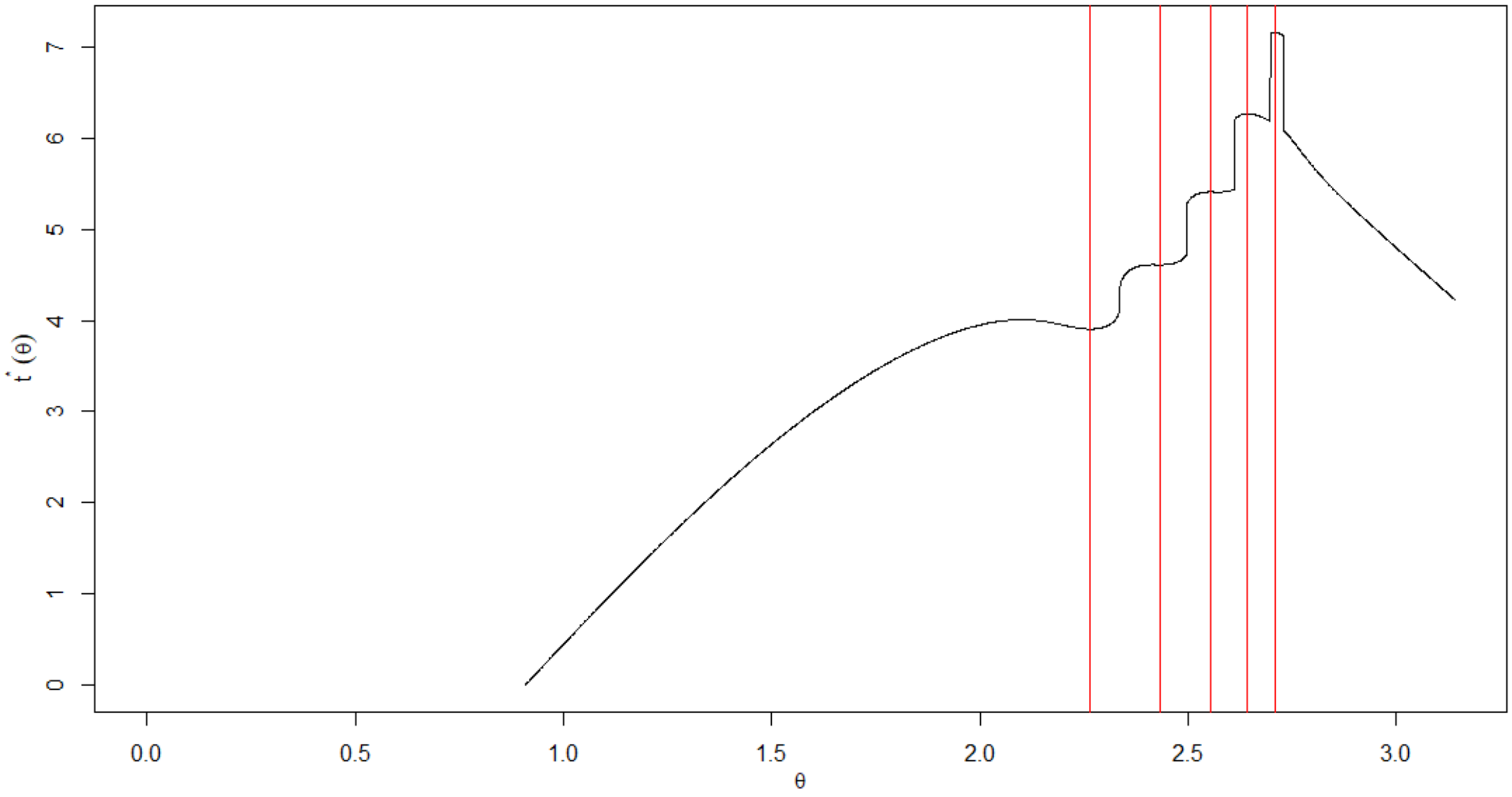

**Fig. 5.32.** Minimizing arguments $t^*(\vartheta)$ vs. $\vartheta$ associated with the polynomial of degree $n = 5$ in equation (5.33). This map was generated from $N = 2{,}500$ elements $LzC(\vartheta_k)$ associated with points $\vartheta_k$ in a regular partition of interval $[0, \pi)$. This graph includes vertical line segments that indicate the location of true theta roots $\vartheta_i^*$.

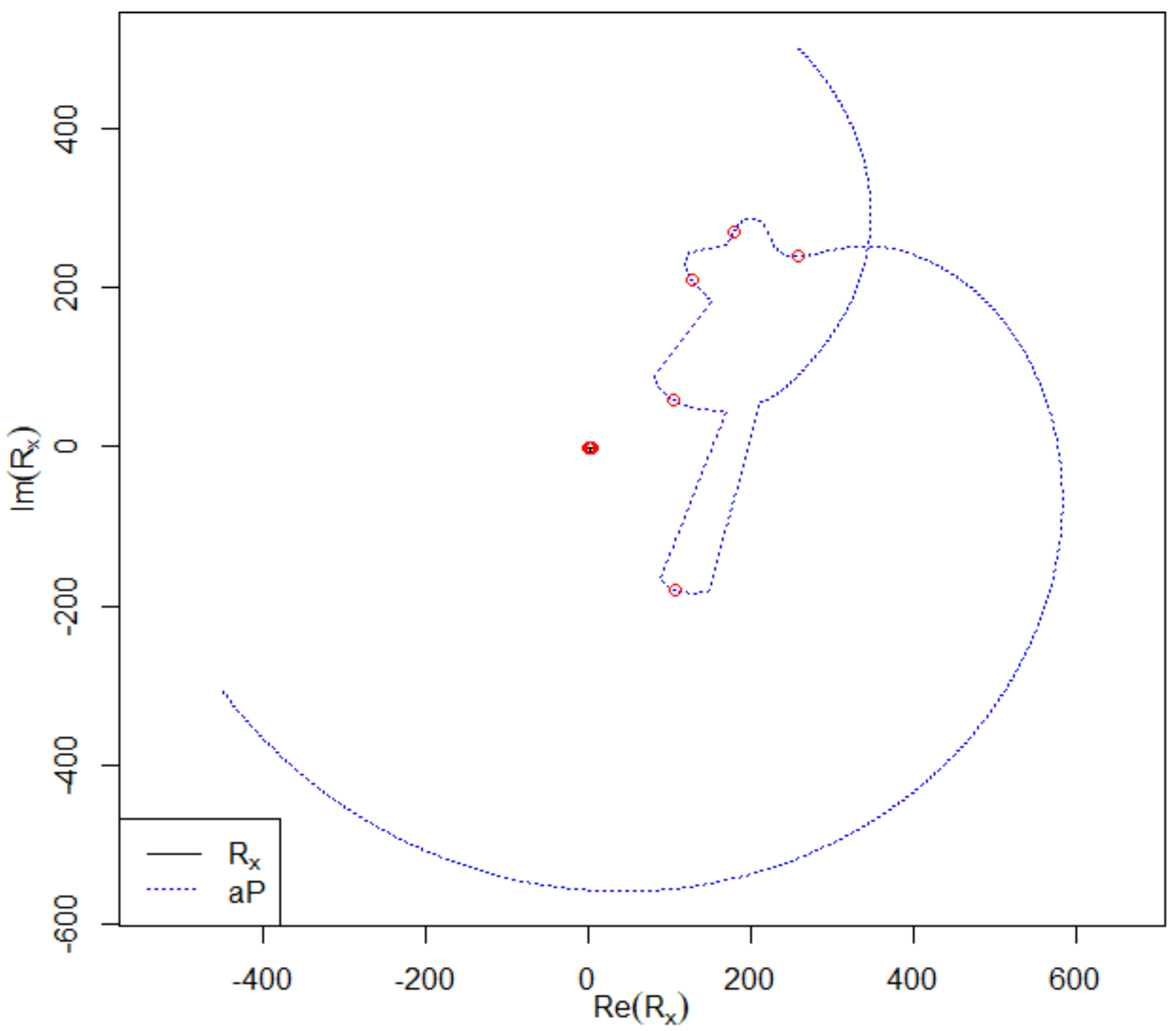

**Fig. 5.33.** Trajectory $r_x(\vartheta)$ (small circles clustered near the origin, in a "microscopic" view) and trajectory $aP(\vartheta)$ associated with the polynomial of degree $n = 5$ in equation (5.33). These trajectories were generated from $N = 2{,}500$ elements $LzC(\vartheta_k)$ associated with points $\vartheta_k$ in a regular partition of interval $[0, \pi)$. Small circles on $aP(\vartheta)$ are anchor points associated with the true roots $r_i$.

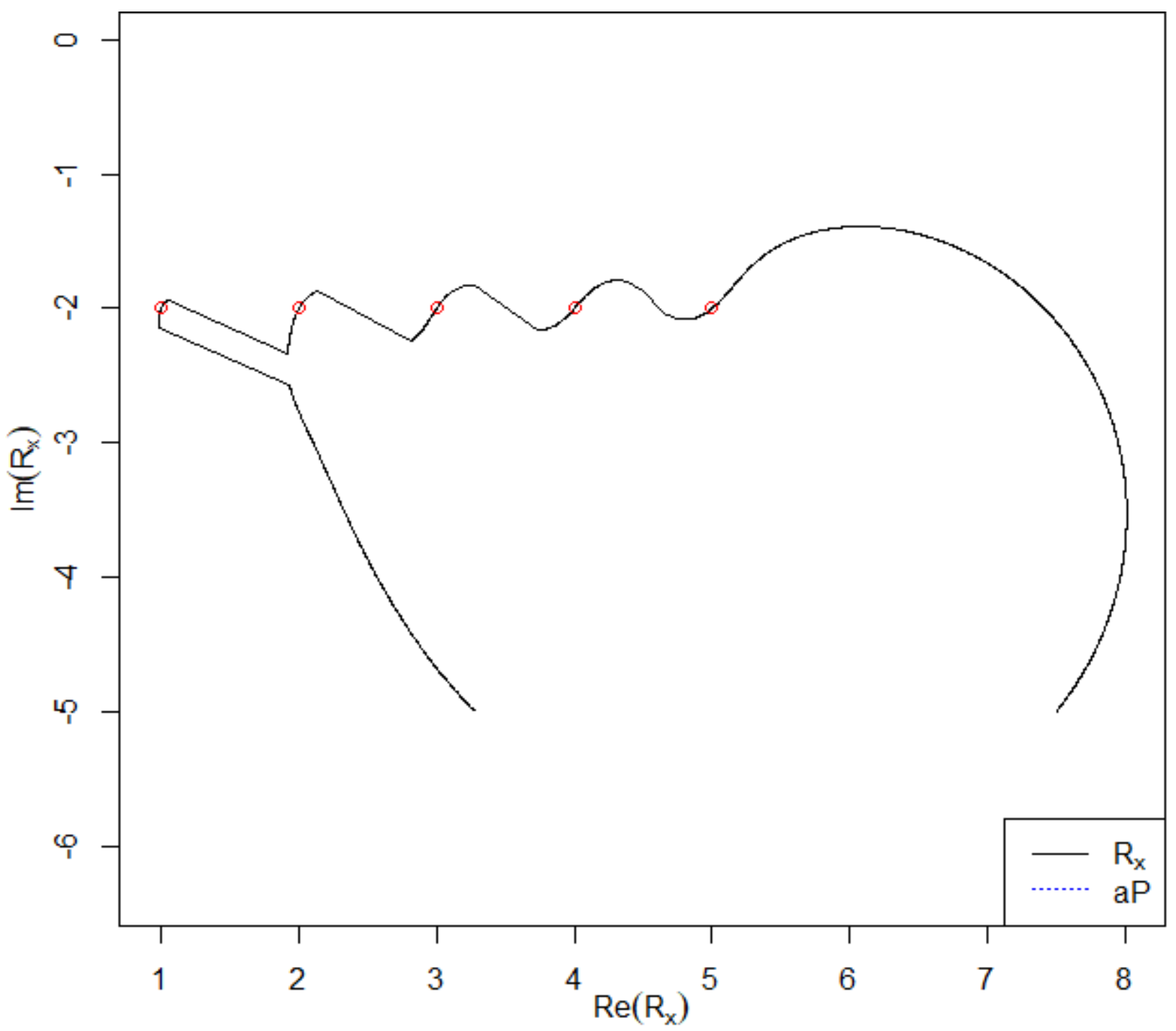

**Fig. 5.34.** Close-up of graph in figure 5.33 in order to better observe trajectory $r_x(\vartheta)$ associated with the polynomial of degree $n = 5$ in equation (5.33). Small circles on $r_x(\vartheta)$ are the true roots $r_i$.

## Numerical Example 5.4

In this final example we will obtain, by means of the LC method, approximations to the roots of the following equation of degree $n = 15$:

$$x^{15} + \sum_{m=1}^{15} m x^{15-m} = 0. \tag{5.35}$$

The reference roots $R_i$ of equation (5.35), which we obtain, as on some previous examples, by means of the R instruction `polyroot`, are:

$$\begin{aligned}
R_1 &= \phantom{-}1.0328432 + 0.5237970i & \theta_1^* &= \phantom{-}0.3292759 \\
R_2 &= -1.0711779 + 0.4580460i & \theta_2^* &= \phantom{-}2.4656747 \\
R_3 &= -0.8110108 - 0.8478387i & \theta_3^* &= -1.9223831 \\
R_4 &= \phantom{-}1.0328432 - 0.5237970i & \theta_4^* &= -0.3292759 \\
R_5 &= \phantom{-}0.4884613 + 1.1640926i & \theta_5^* &= \phantom{-}0.8668099 \\
R_6 &= -1.1624405 - 0.0000000i & \theta_6^* &= -3.1415927 \\
R_7 &= -1.0711779 - 0.4580460i & \theta_7^* &= -2.4656747 \\
R_8 &= \phantom{-}0.8223161 - 0.9272969i & \theta_8^* &= -0.6115759 \\
R_9 &= \phantom{-}0.0405263 + 1.2186607i & \theta_9^* &= \phantom{-}1.1533265 \\
R_{10} &= -0.8110108 + 0.8478387i & \theta_{10}^* &= \phantom{-}1.9223831 \\
R_{11} &= \phantom{-}0.0405263 - 1.2186607i & \theta_{11}^* &= -1.1533265 \\
R_{12} &= \phantom{-}0.8223161 + 0.9272969i & \theta_{12}^* &= \phantom{-}0.6115759 \\
R_{13} &= -0.4207379 + 1.1126076i & \theta_{13}^* &= \phantom{-}1.4996765 \\
R_{14} &= -0.4207379 - 1.1126076i & \theta_{14}^* &= -1.4996765 \\
R_{15} &= \phantom{-}0.4884613 - 1.1640926i & \theta_{15}^* &= -0.8668099
\end{aligned} \tag{5.36}$$

Theta root $\theta_i^*$ in (5.36), for $i = 1,2,\ldots,15$, as we already know, is the inclination angle of line $\ell_1$, which contains both the fixed point $P_1 = -C_1/2 = -1/2$ and root $R_i$. We see that all roots $R_i$ in (5.36), except for one, $R_6$, are complex numbers with non-zero imaginary part, which occur in conjugate pairs, i.e., we have seven pairs of roots, where each pair shares the same real part, and possesses imaginary parts of the same magnitude, but of opposite signs. One of such pairs in (5.36) is $R_{13}$ and $R_{14}$.

In this numerical example, we will test the LC method by using as input the original coefficients $C_m = m$ of equation (5.35), although as always, reconstructed from the reference roots (5.36) via Vieta's relations (5.2). The operating conditions under which the results will be obtained in this numerical example are again specified in the final part of annex 3 section 3, within subsection "specific conditions for reproducing results from the examples in chapter 5".

We will construct proximity maps by using $N_1 = 1{,}000$ and $N_2 = 50{,}000$ elements $LzC(\theta_k)$ associated to points $\theta_k$ in regular partitions of interval $[-\pi, \pi)$, in order to give us an idea of the effect that the resolution of proximity maps has on their corresponding root estimates. Figures 5.35 and 5.36 show the proximity maps $e(\theta)$ generated with these two resolutions, while figures 5.37 and 5.38 show corresponding maps $\frac{d}{d\theta} d_\theta^2(t^*(\theta))$, and figures 5.39 and 5.40 show corresponding maps $\frac{d}{d\theta} t^*(\theta)$. From figures 5.35 and 5.36, we see that with $N_1 = 1{,}000$ elements $LzC(\theta_k)$ we have a problem of gaps between $e_A(\theta)$ and $e_B(\theta)$, to the extent that some roots in the vicinity of these gaps cannot be detected; with $N_2 = 50{,}000$ elements $LzC(\theta_k)$, on the other hand, this gap problem subsides, and the map $e(\theta)$ is able to detect all complex roots in (5.36), but not the real root $R_6$.

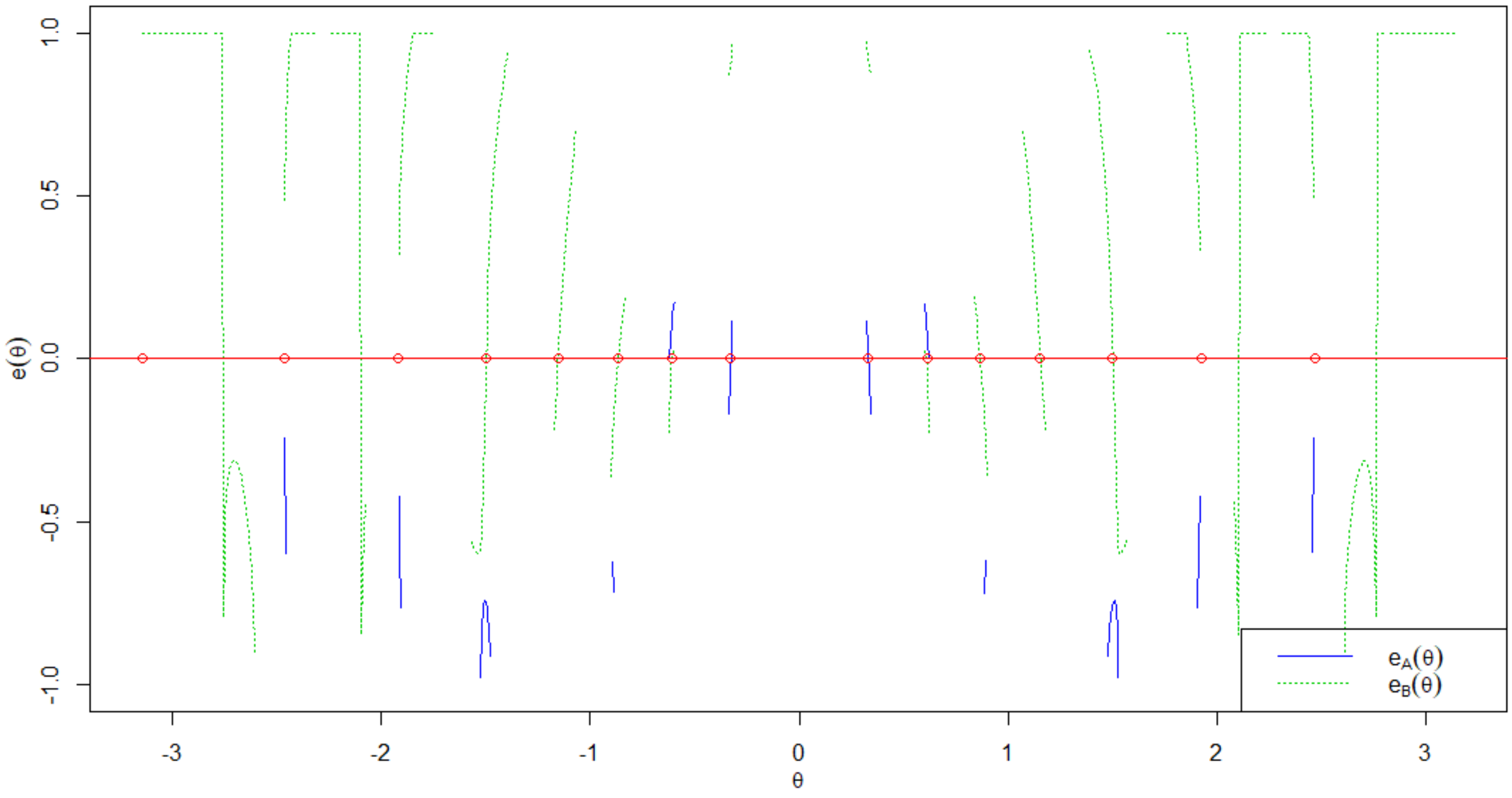

**Fig. 5.35.** Map $e(\theta)$ for polynomial $x^{15} + \sum_{m=1}^{15} mx^{15-m}$ in equation (5.35). This map was generated from $N_1 = 1{,}000$ elements $LzC(\theta_k)$ associated with points $\theta_k$ in a regular partition of interval $[-\pi, \pi)$. Reference theta roots $\theta_i^*$ in (5.36) are shown by means of small circles on horizontal axis $y = 0$.

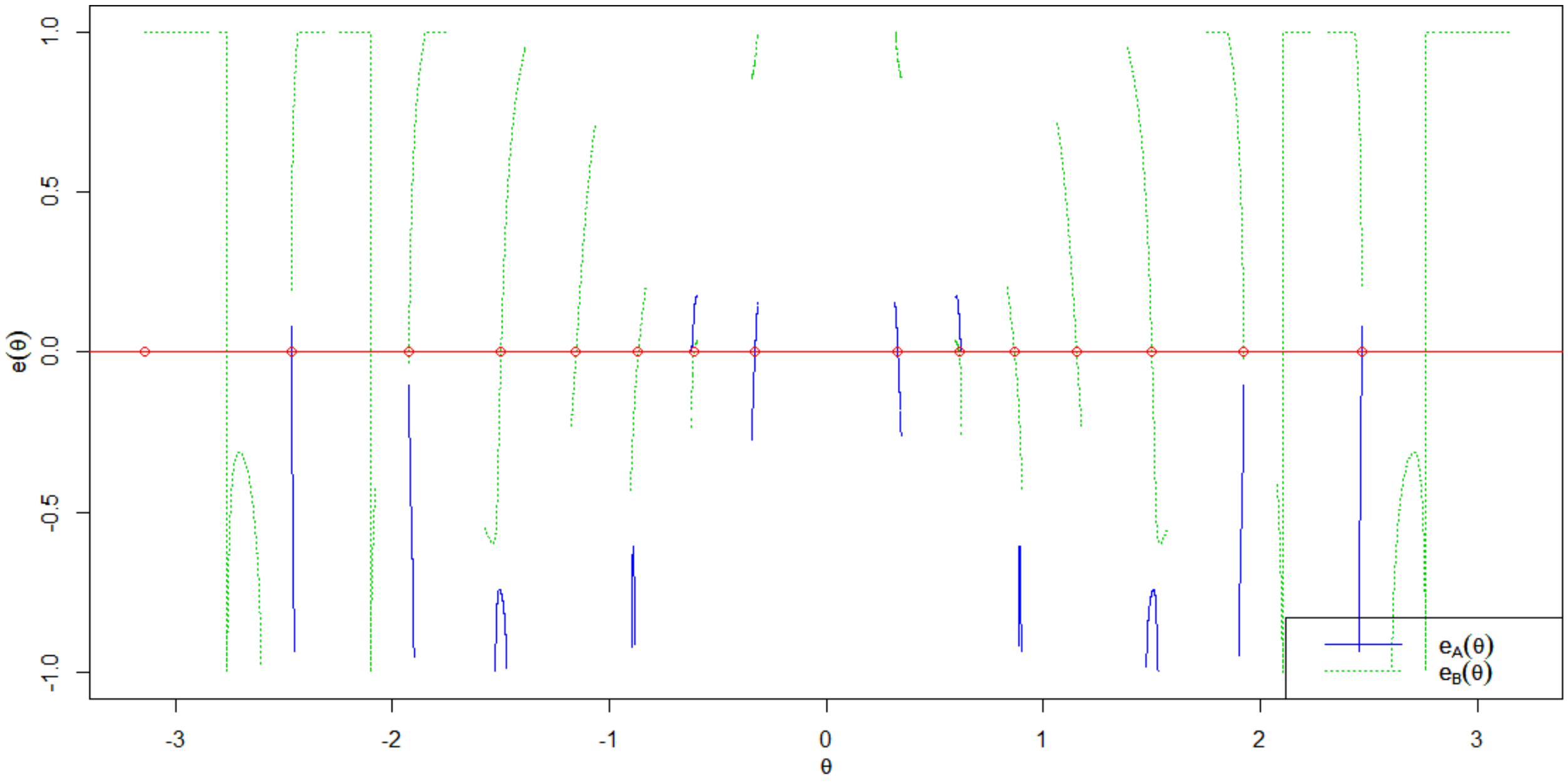

**Fig. 5.36.** Map $e(\theta)$ for polynomial $x^{15} + \sum_{m=1}^{15} mx^{15-m}$ in equation (5.35). This map was generated from $N_2 = 50{,}000$ elements $LzC(\theta_k)$ associated with points $\theta_k$ in a regular partition of interval $[-\pi, \pi)$. Reference theta roots $\theta_i^*$ in (5.36) are shown by means of small circles on horizontal axis $y = 0$.

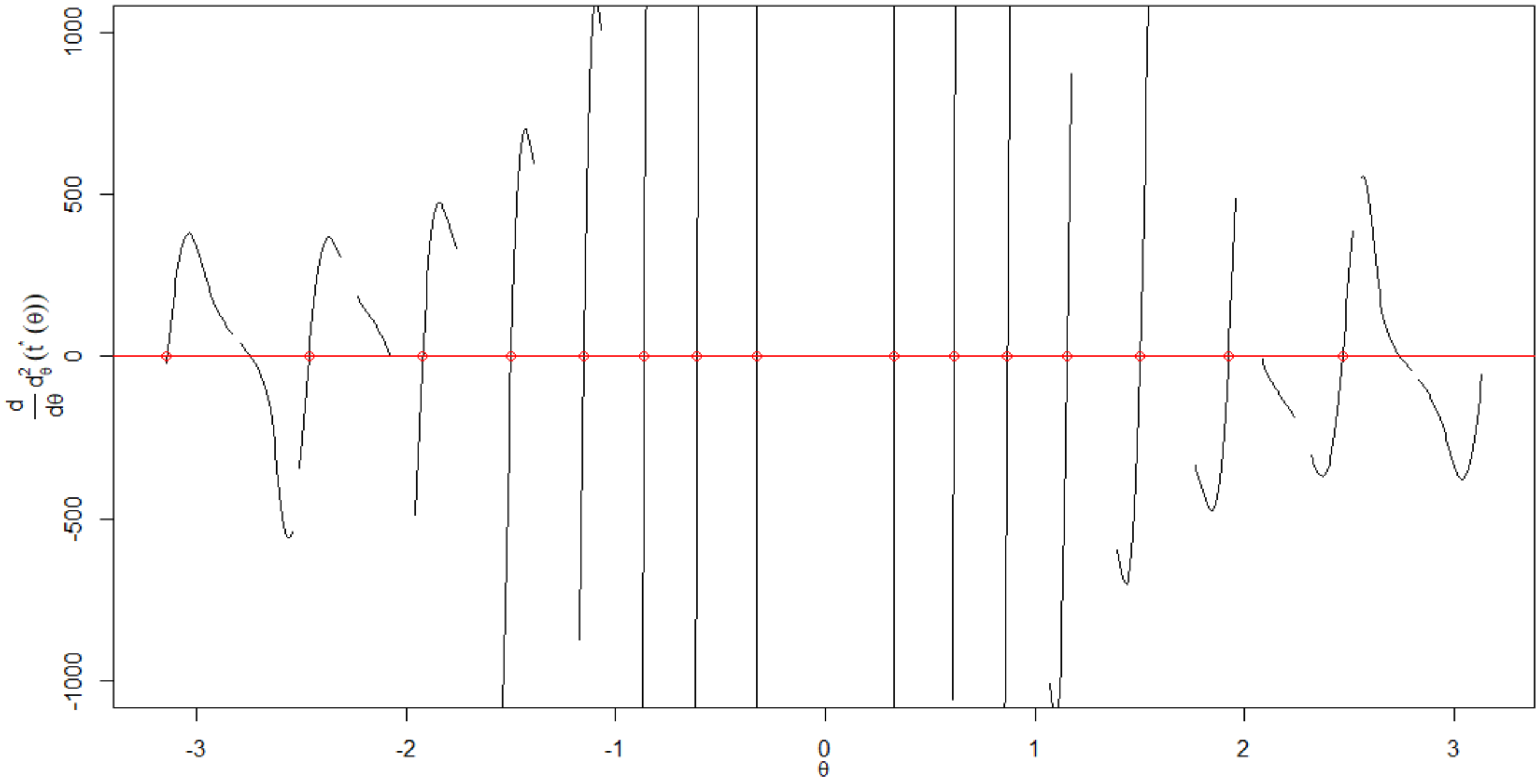

**Fig. 5.37.** Map $\frac{d}{d\theta}d_\theta^2(t^*(\theta))$ for polynomial $x^{15}+\sum_{m=1}^{15}mx^{15-m}$ in equation (5.35). This map was generated from $N_1 = 1{,}000$ elements $LzC(\theta_k)$ associated with points $\theta_k$ in a regular partition of interval $[-\pi,\pi)$. Reference theta roots $\theta_i^*$ in (5.36) are shown by means of small circles on horizontal axis $y = 0$.

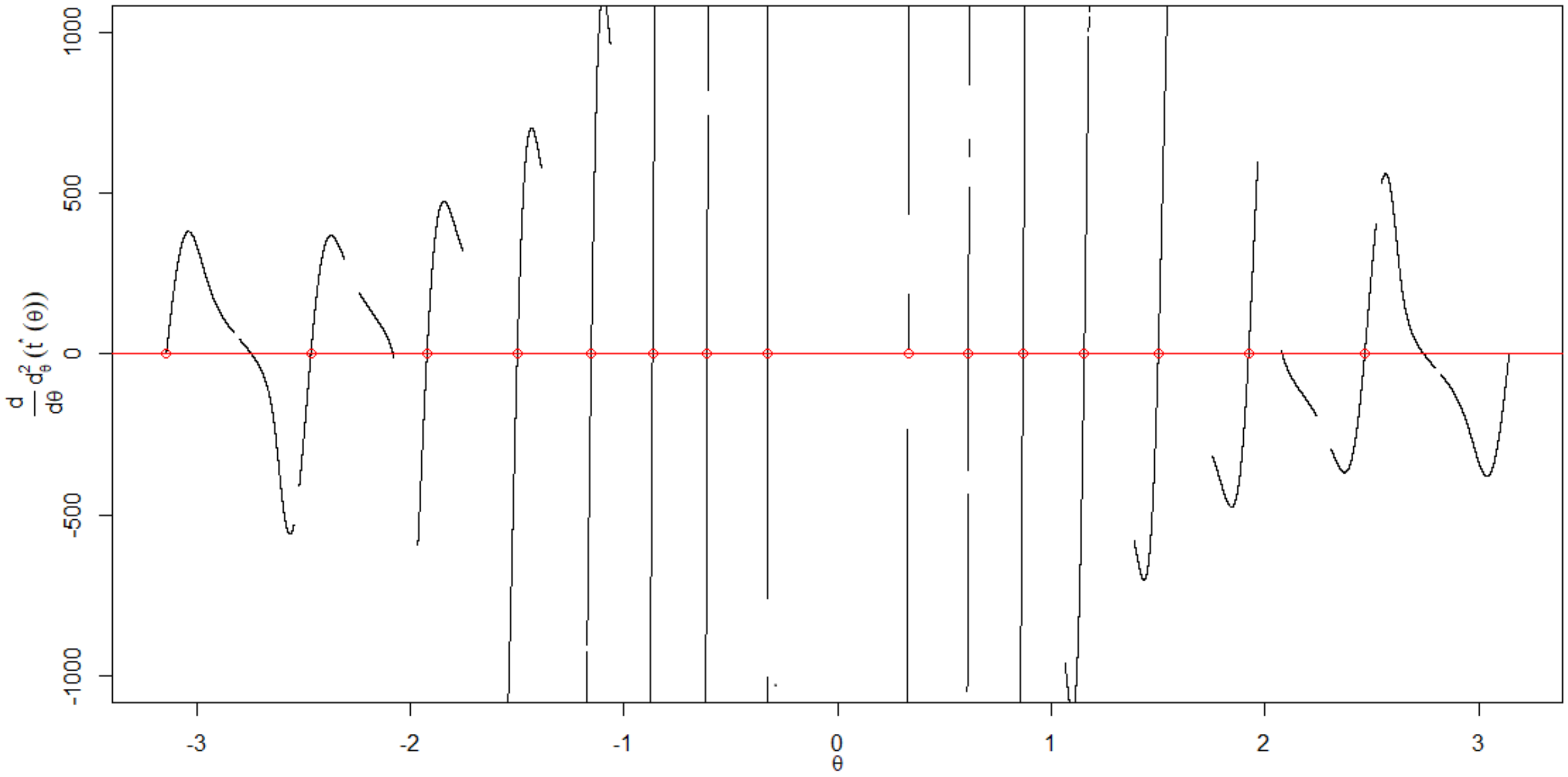

**Fig. 5.38.** Map $\frac{d}{d\theta}d_\theta^2(t^*(\theta))$ for polynomial $x^{15}+\sum_{m=1}^{15}mx^{15-m}$ in equation (5.35). This map was generated from $N_2 = 50{,}000$ elements $LzC(\theta_k)$ associated with points $\theta_k$ in a regular partition of interval $[-\pi,\pi)$. Reference theta roots $\theta_i^*$ in (5.36) are shown by means of small circles on horizontal axis $y = 0$.

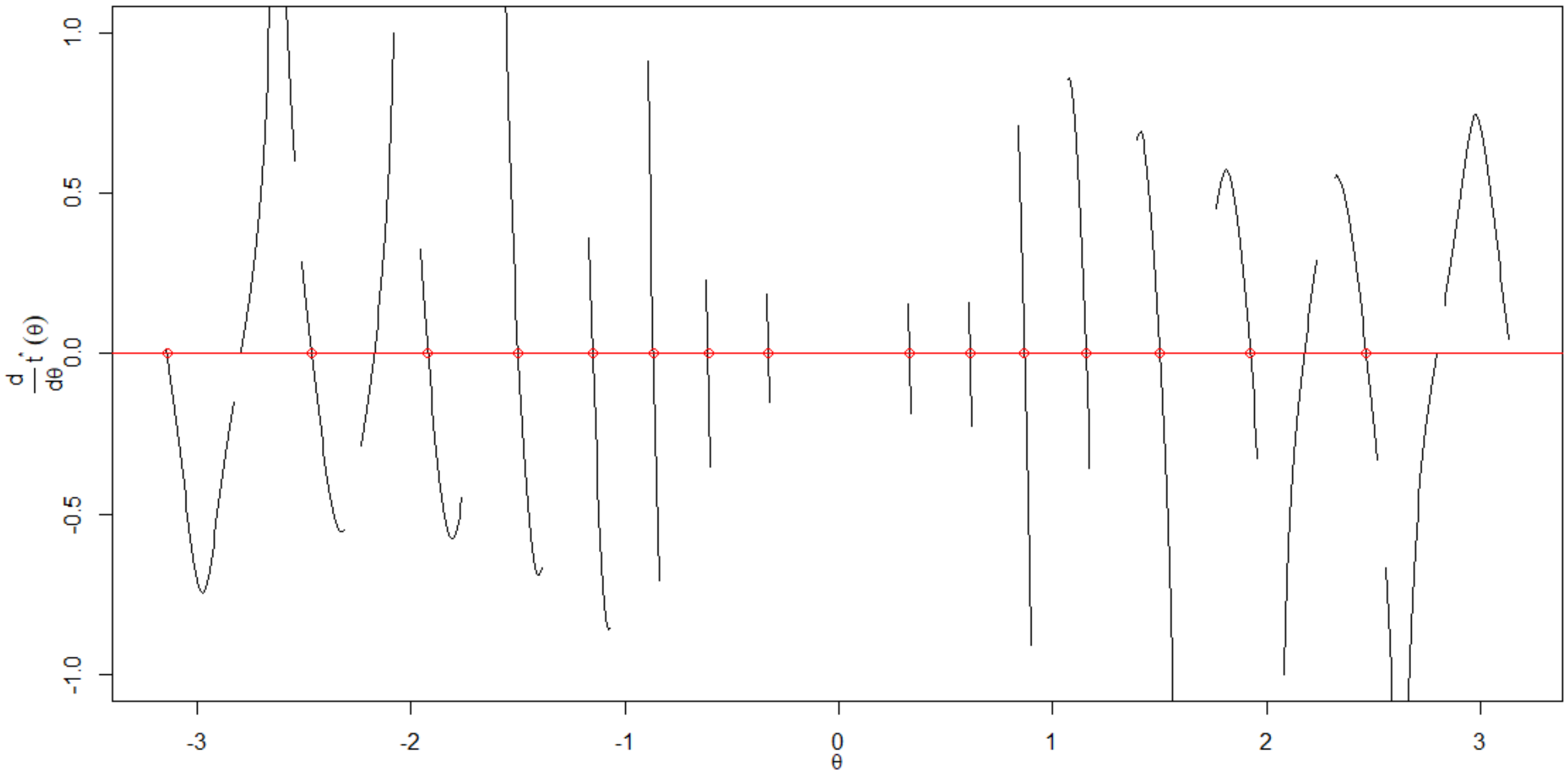

**Fig. 5.39.** Map $\frac{d}{d\theta}t^*(\theta)$ for polynomial $x^{15}+\sum_{m=1}^{15} m x^{15-m}$ in equation (5.35). This map was generated from $N_1 = 1{,}000$ elements $LzC(\theta_k)$ associated with points $\theta_k$ in a regular partition of interval $[-\pi,\pi)$. Reference theta roots $\theta_i^*$ in (5.36) are shown by means of small circles on horizontal axis $y = 0$.

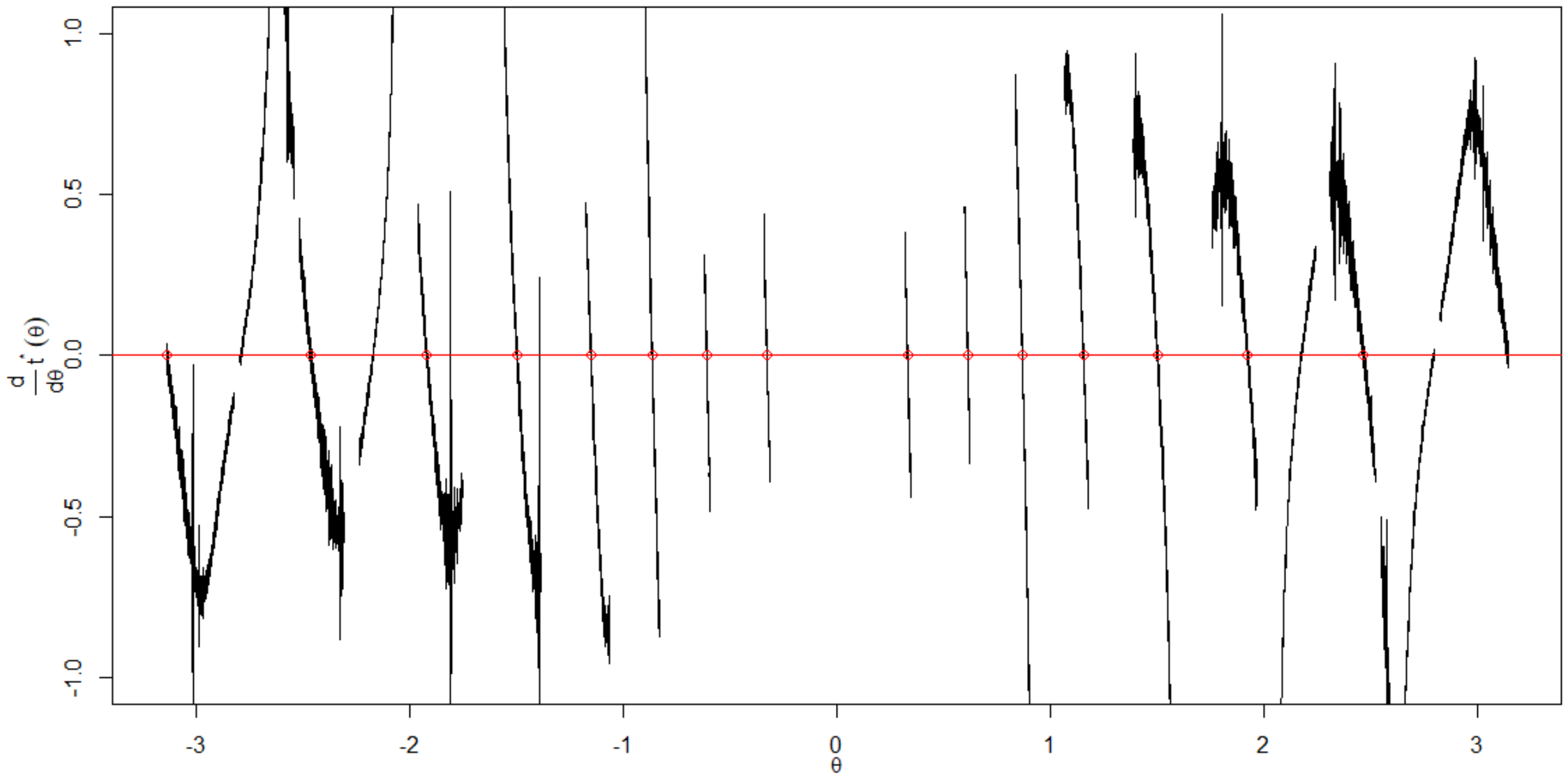

**Fig. 5.40.** Map $\frac{d}{d\theta}t^*(\theta)$ for polynomial $x^{15}+\sum_{m=1}^{15} m x^{15-m}$ in equation (5.35). This map was generated from $N_2 = 50{,}000$ elements $LzC(\theta_k)$ associated with points $\theta_k$ in a regular partition of interval $[-\pi,\pi)$. Reference theta roots $\theta_i^*$ in (5.36) are shown by means of small circles on horizontal axis $y = 0$.

From figures 5.37 and 5.38 we can see that, apparently, there is not much difference between the two maps $\frac{d}{d\theta}d_\theta^2(t^*(\theta))$ generated with the two resolutions $N_1 = 1{,}000$ and $N_2 = 50{,}000$, although for the map generated with $N_2 = 50{,}000$ elements $LzC(\theta_k)$ (figure 5.38), there seems to be a little more difficulty in detecting all the roots of equation (5.35), due to the presence of gaps caused by errors in the process for obtaining the minimizing arguments $t^*(\theta)$: the simulated annealing algorithm gets "stuck" at local minima for some values $\theta_k$, as we saw in a little more detail in numerical example 4.2 of chapter 4. On the other hand, for the map $\frac{d}{d\theta}d_\theta^2(t^*(\theta))$ generated from $N_1 = 1{,}000$ elements $LzC(\theta_k)$ (figure 5.37), there was better luck in identifying the 15 roots of equation (5.35). On both maps $\frac{d}{d\theta}d_\theta^2(t^*(\theta))$ from figures 5.37 and 5.38, the real root $R_6$ was correctly identified.

From figures 5.39 and 5.40, we see that with $N_1 = 1{,}000$ elements $LzC(\theta_k)$, we obtain a map $\frac{d}{d\theta}t^*(\theta)$ of acceptable quality, which correctly identifies all the roots of equation (5.35), including the real root $R_6$; however, with $N_2 = 50{,}000$ elements $LzC(\theta_k)$, the quality of map $\frac{d}{d\theta}t^*(\theta)$ declines markedly, due to the presence of high-frequency oscillatory noise caused, as we have seen before in example 4.2 of chapter 4, by the combination of two factors: 1) the amplification of small optimization errors caused by the BFGS algorithm in approximating $t^*(\theta)$, attributable to the discrete operator $\frac{\Delta}{\Delta\theta}$ used in the construction of map $\frac{d}{d\theta}t^*(\theta)$, and 2) the high resolution of said map. Consequently, although valid initial estimates generated by map $\frac{d}{d\theta}t^*(\theta)$ in figure 5.40 (built with $N_2 = 50{,}000$ elements $LzC(\theta_k)$) cover all 15 roots of equation (5.35), they are found interspersed with spurious estimates associated with high-frequency oscillatory noise; it is worth mentioning that, if we arrange in ascending order the 360 initial estimates generated by map $\frac{d}{d\theta}t^*(\theta)$ in figure 5.40 according to their corresponding values $d_{\hat{\theta}_i^*}^2\left(t^*\left(\hat{\theta}_i^*\right)\right)$, we will find that the valid initial estimates (i.e., the initial estimates that best approximate reference roots) are located in the following positions: 1, 2, 3, 4, 5, 6, 7, 11, 13, 22, 24, 31, 42, 57 and 178. In contrast, valid initial estimates obtained with all other maps in figures 5.35, 5.36, 5.37, 5.38, and 5.39, are always at the top positions and not interspersed with spurious estimates, when all the initial estimates are arranged in ascending order according to variable $d_{\hat{\theta}_i^*}^2\left(t^*\left(\hat{\theta}_i^*\right)\right)$.

Tables 5.16 and 5.17 summarize the approximation errors associated with valid initial estimates of the roots of equation (5.35), obtained by means of the proximity maps shown in figures 5.35, 5.36, 5.37, 5.38, 5.39, and 5.40. Of course, the approximation errors $\left|(R-\hat{R})/R\right|$, $\left|(\theta^*-\hat{\theta}^*)/\theta^*\right|$ for valid initial estimates $\hat{R}$, $\hat{\theta}^*$, are relative to the reference roots (5.36).

**Table 5.16**. Arithmetic means and standard deviations for error measures associated with valid initial estimates of the roots of the polynomial $x^{15}+\sum_{m=1}^{15} mx^{15-m}$ in equation (5.35), detected by means of proximity maps from figures 5.35, 5.37 and 5.39, built with $N_1 = 1{,}000$ elements $LzC(\theta_k)$.

| Map | Number of valid initial estimates | Error measure | | | | | |
|---|---|---|---|---|---|---|---|
| | | $d^2_{\hat{\theta}_i^*}\left(t^*(\hat{\theta}_i^*)\right)$ | | $\left\lvert\frac{R-\hat{R}}{R}\right\rvert$ | | $\left\lvert\frac{\theta^*-\hat{\theta}^*}{\theta^*}\right\rvert$ | |
| | | Mean | Standard deviation | Mean | Standard deviation | Mean | Standard deviation |
| $e(\boldsymbol{\theta})$ | 10 | $1.744153 \times 10^{-2}$ | $3.432715 \times 10^{-2}$ | $2.955913 \times 10^{-4}$ | $4.709806 \times 10^{-4}$ | $3.876965 \times 10^{-4}$ | $5.928253 \times 10^{-4}$ |
| $\frac{d}{d\boldsymbol{\theta}} d^2_{\boldsymbol{\theta}}(t^*(\boldsymbol{\theta}))$ | 15 | $2.560758 \times 10^{-4}$ | $4.995276 \times 10^{-4}$ | $5.067946 \times 10^{-5}$ | $3.282447 \times 10^{-5}$ | $5.860315 \times 10^{-5}$ | $6.359460 \times 10^{-5}$ |
| $\frac{d}{d\boldsymbol{\theta}} t^*(\boldsymbol{\theta})$ | 15 | $4.284718 \times 10^{-5}$ | $6.298015 \times 10^{-5}$ | $3.161279 \times 10^{-5}$ | $2.117546 \times 10^{-5}$ | $2.897753 \times 10^{-5}$ | $2.175426 \times 10^{-5}$ |

**Table 5.17**. Arithmetic means and standard deviations for error measures associated with valid initial estimates of the roots of the polynomial $x^{15}+\sum_{m=1}^{15} mx^{15-m}$ in equation (5.35), detected by means of proximity maps from figures 5.36, 5.38 and 5.40, built with $N_2 = 50{,}000$ elements $LzC(\theta_k)$.

| Map | Number of valid initial estimates | Error measure | | | | | |
|---|---|---|---|---|---|---|---|
| | | $d^2_{\hat{\theta}_i^*}\left(t^*(\hat{\theta}_i^*)\right)$ | | $\left\lvert\frac{R-\hat{R}}{R}\right\rvert$ | | $\left\lvert\frac{\theta^*-\hat{\theta}^*}{\theta^*}\right\rvert$ | |
| | | Mean | Standard deviation | Mean | Standard deviation | Mean | Standard deviation |
| $e(\boldsymbol{\theta})$ | 14 | $1.070257 \times 10^{-5}$ | $2.964008 \times 10^{-5}$ | $6.862822 \times 10^{-6}$ | $1.212532 \times 10^{-5}$ | $7.156495 \times 10^{-6}$ | $1.569346 \times 10^{-5}$ |
| $\frac{d}{d\boldsymbol{\theta}} d^2_{\boldsymbol{\theta}}(t^*(\boldsymbol{\theta}))$ | 13 | $2.406454 \times 10^{-7}$ | $3.134019 \times 10^{-7}$ | $3.072763 \times 10^{-6}$ | $1.735169 \times 10^{-6}$ | $3.733288 \times 10^{-8}$ | $2.944051 \times 10^{-8}$ |
| $\frac{d}{d\boldsymbol{\theta}} t^*(\boldsymbol{\theta})$ | 15 | $2.442306 \times 10^{-3}$ | $8.029688 \times 10^{-3}$ | $9.267667 \times 10^{-5}$ | $1.103970 \times 10^{-4}$ | $1.292787 \times 10^{-4}$ | $2.347825 \times 10^{-4}$ |

When comparing tables 5.16 and 5.17, we see that the increase in resolution, from $N_1 = 1{,}000$ to $N_2 = 50{,}000$ elements $LzC(\theta_k)$, improves the accuracy of the valid initial estimates obtained by maps $e(\theta)$ and $\frac{d}{d\theta} d^2_\theta(t^*(\theta))$, but also causes a slight reduction in the quality of valid initial estimates obtained by map $\frac{d}{d\theta} t^*(\theta)$. Additionally, with this comparison we confirm the difficulty observed on the map $\frac{d}{d\theta} d^2_\theta(t^*(\theta))$ from figure 5.38 to detect all roots of equation (5.35), due to the characteristics of our optimization process used in obtaining approximations to minimizing arguments $t^*$ for dynamic squared distance functions $d^2_\theta$.

To conclude with the analysis of results in this example, let us look at some graphs related to the global minima of functions $d^2_{\theta_k}$ associated with elements $LzC(\theta_k)$ derived from equation (5.35). Figures 5.41, 5.42, 5.43, 5.44 and 5.45 were obtained by using $N_2 = 50{,}000$ elements $LzC(\theta_k)$ associated with points $\theta_k$ in a regular partition of interval $[-\pi, \pi)$.

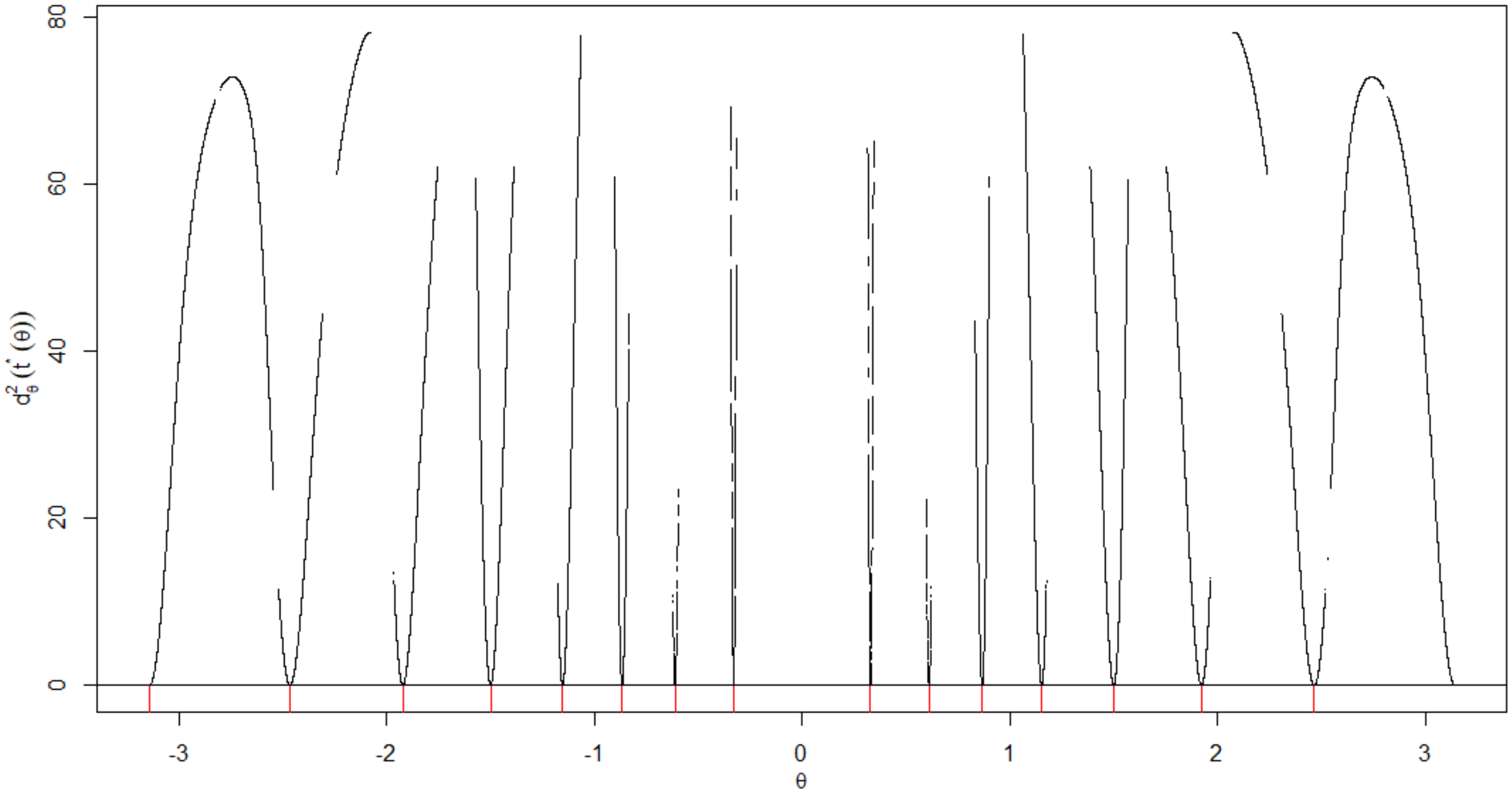

**Fig. 5.41.** Global minima $d_\theta^2(t^*(\theta))$ vs. $\theta$ associated with polynomial $x^{15} + \sum_{m=1}^{15} m x^{15-m}$ in equation (5.35). This map was generated from $N_2 = 50{,}000$ elements $LzC(\theta_k)$ associated with points $\theta_k$ in a regular partition of interval $[-\pi, \pi)$. This graph includes vertical line segments that indicate the location of true theta roots $\theta_i^*$.

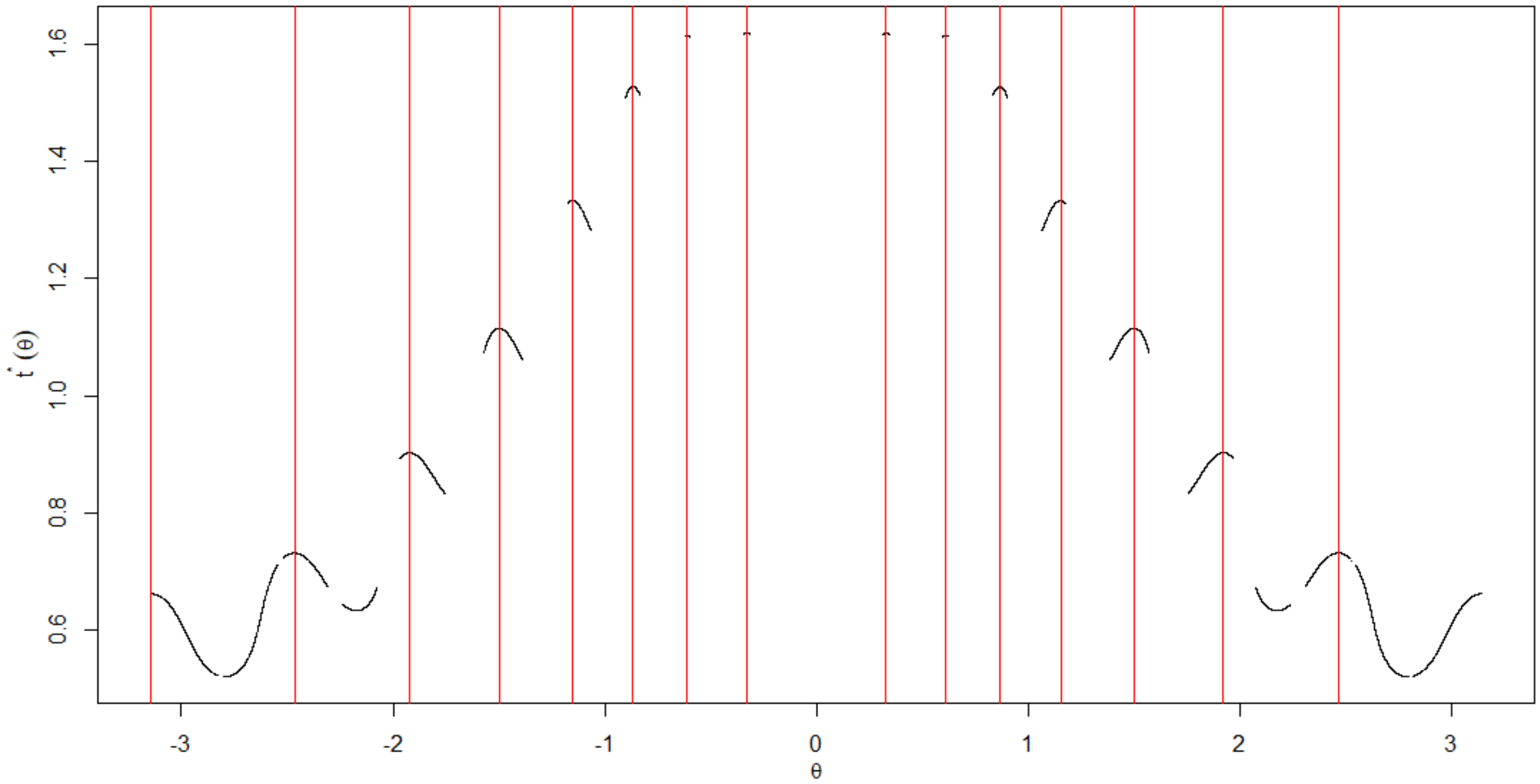

**Fig. 5.42.** Minimizing arguments $t^*(\theta)$ vs. $\theta$ associated with polynomial $x^{15} + \sum_{m=1}^{15} m x^{15-m}$ in equation (5.35). This map was generated from $N_2 = 50{,}000$ elements $LzC(\theta_k)$ associated with points $\theta_k$ in a regular partition of interval $[-\pi, \pi)$. This graph includes vertical line segments that indicate the location of true theta roots $\theta_i^*$.

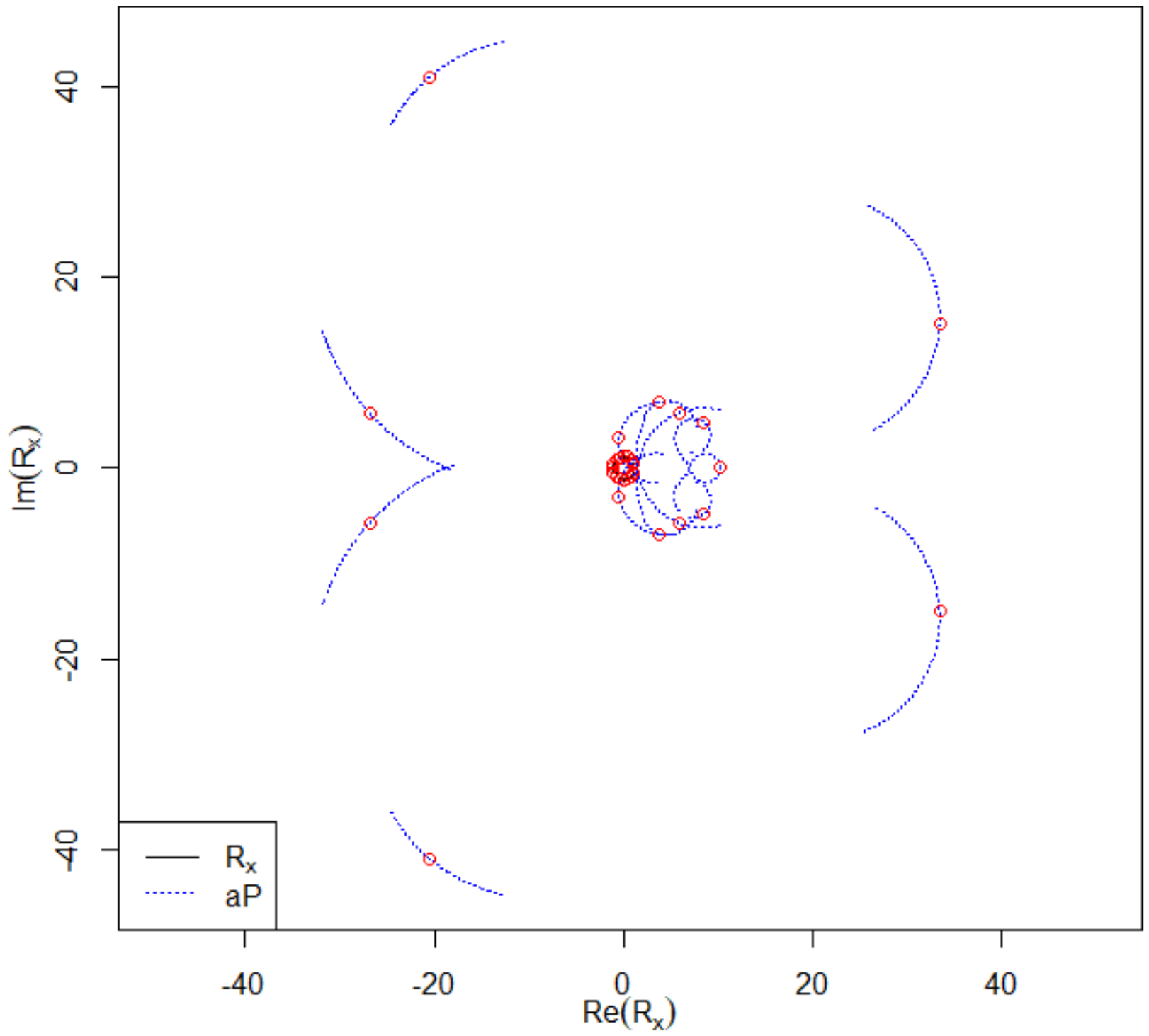

**Fig. 5.43.** Trajectory $R_x(\theta)$ (small circles clustered around the origin, in a "microscopic" view) and trajectory $aP(\theta)$ associated with polynomial $x^{15} + \sum_{m=1}^{15} mx^{15-m}$ in equation (5.35). These trajectories were generated from $N_2 = 50{,}000$ elements $LzC(\theta_k)$ associated with points $\theta_k$ in a regular partition of interval $[-\pi, \pi)$. The small circles on $aP(\theta)$ are anchor points associated with true roots $R_i$.

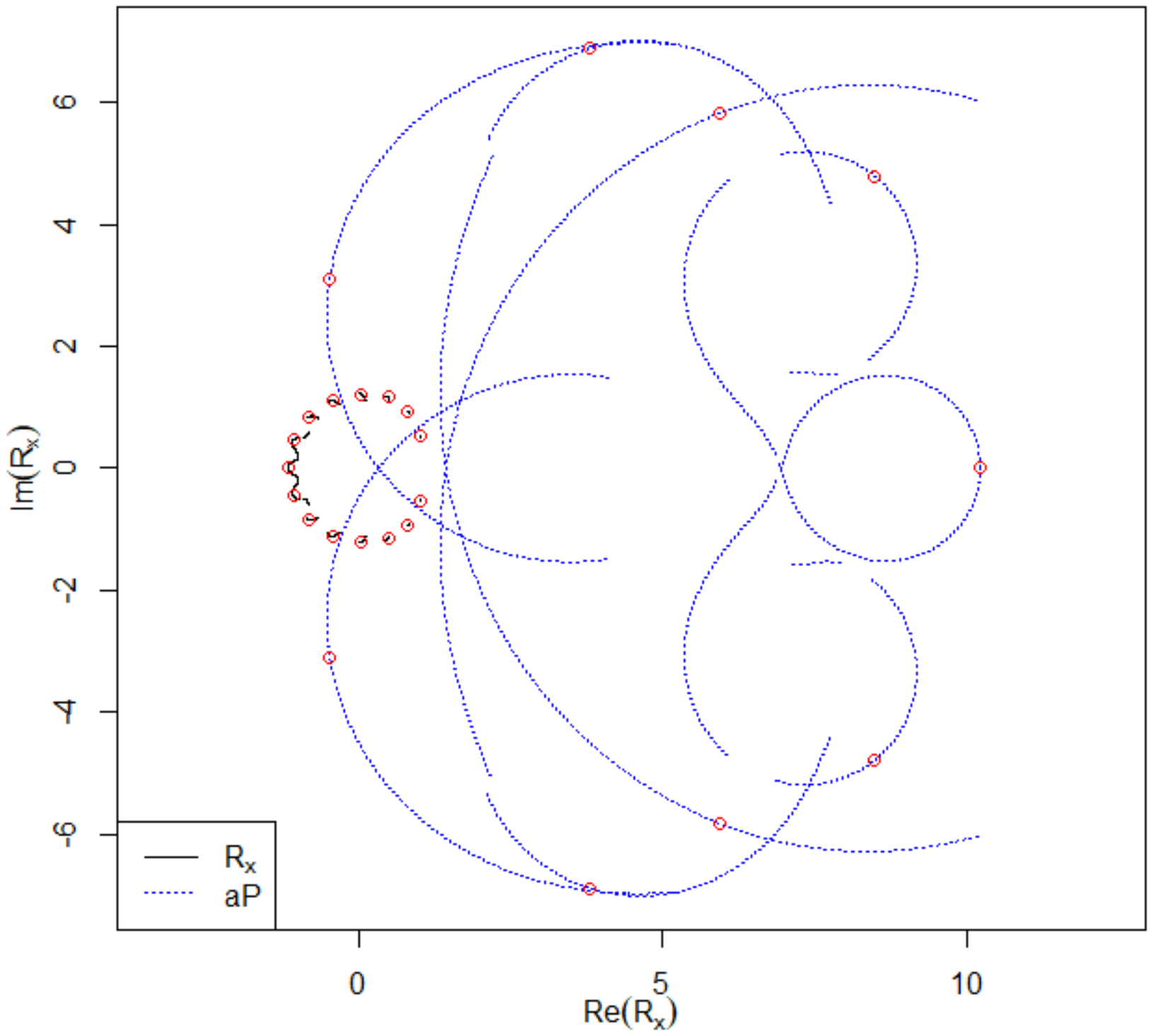

**Fig. 5.44.** Close-up of figure 5.43 in order to better observe the central structure of trajectories $R_x(\theta)$ and $aP(\theta)$ associated with polynomial $x^{15} + \sum_{m=1}^{15} mx^{15-m}$ in equation (5.35).

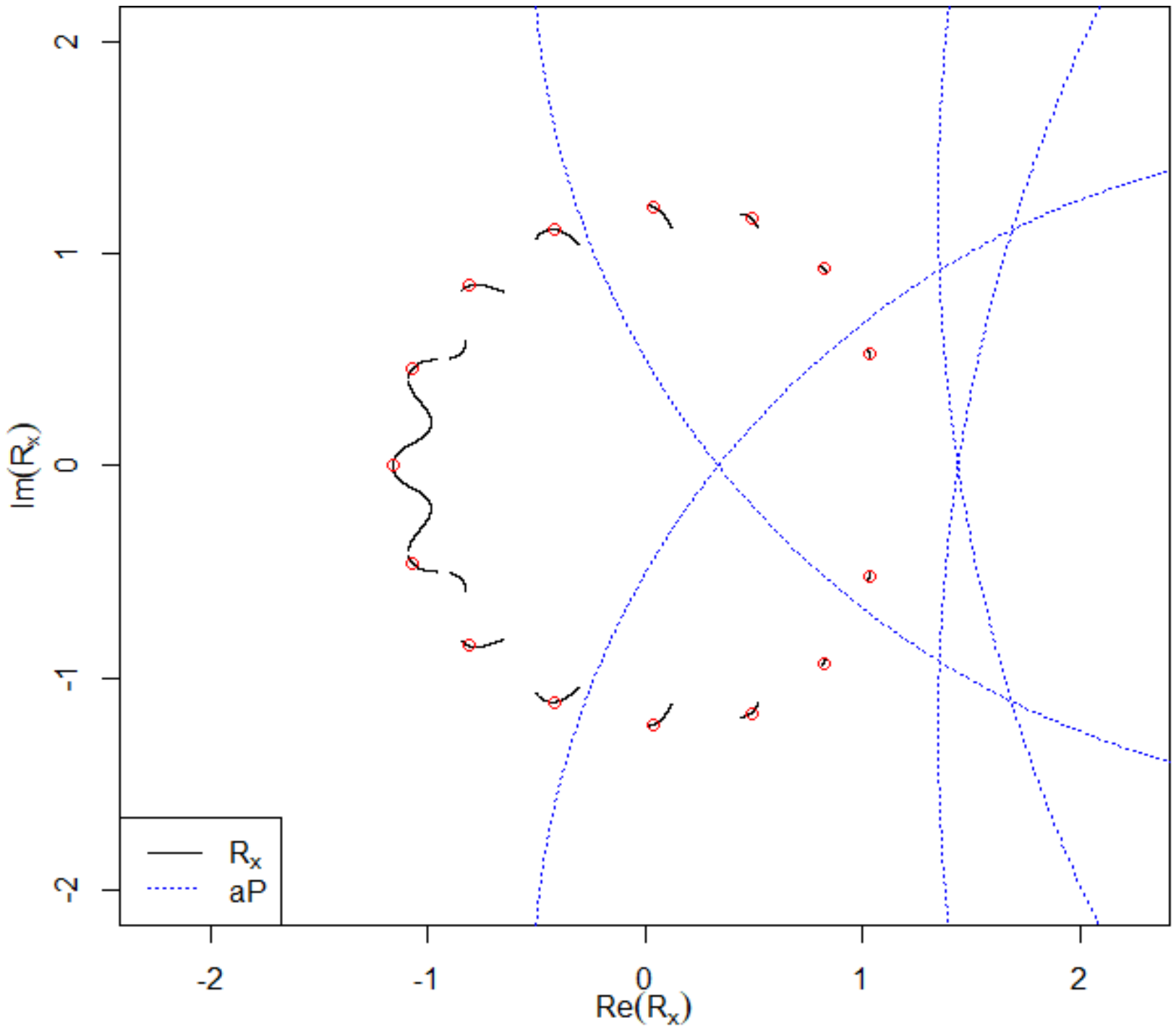

**Fig. 5.45.** Second close-up of figure 5.43 in order to better observe the trajectory $R_x(\theta)$ associated with polynomial $x^{15} + \sum_{m=1}^{15} mx^{15-m}$ in equation (5.35). Small circles on trajectory $R_x(\theta)$ are the true roots $R_i$.

In summary, from the graphs generated in this numerical example, we observe the following:

- In figures 5.38 and 5.41, gaps can be seen within maps $\frac{d}{d\theta} d_\theta^2(t^*(\theta))$ and $d_\theta^2(t^*(\theta))$; these gaps, as we have already mentioned, correspond to errors that occur when it comes to approximating the global minimum of function $d_{\theta_k}^2(t)$. Figures 5.38 and 5.41 also reinforce the hypothesis (posed in the numerical example 4.1 of chapter 4) that theta roots $\theta_i^*$ occur at the global minima of mapping $d_\theta^2(t^*(\theta))$.
- Figures 5.39 and 5.42 suggest that, in this case, all theta roots $\theta_i^*$ correspond to local maxima of mapping $t^*(\theta)$.
- When comparing figure 5.43 (and its close-ups in figures 5.44 and 5.45) with the rest of the graphs obtained in this example, we see that, apparently, the trajectories and most of the maps associated with polynomial $x^{15} + \sum_{m=1}^{15} mx^{15-m}$ are made up of 19 continuous sections. Symmetries can also be seen among continuous sections: trajectories $R_x(\theta)$ and $aP(\theta)$ are symmetric with respect to the real axis in the complex plane $\mathbb{C}$, maps $d_\theta^2(t^*(\theta))$, $t^*(\theta)$, $e(\theta)$ are symmetric with respect to vertical axis $\theta = 0$, and maps $\frac{d}{d\theta} t^*(\theta)$, $\frac{d}{d\theta} d_\theta^2(t^*(\theta))$ are symmetric with respect to the origin $(\theta = 0, y = 0)$. These symmetries are probably associated with the nature of the integer coefficients $C_m = m$ in equation (5.35). We will not investigate such associations here, as they are beyond our objectives.

## Conclusions and Future Work

Throughout this work we have empirically verified that the LC method, under suitable operating conditions, generates good initial estimates of roots of polynomials in a complex variable. The empirical evidence obtained here shows that the LC method works properly when all the roots are sufficiently separated from each other; there are difficulties, of course, as in any numerical method designed to solve polynomial equations, when there are groups of two or more roots very close to each other, or when there are multiple roots. In example 5.3, we saw that the LC method can also be used in the resolution of polynomials of a real variable, if we previously extend the polynomial variable to the complex plane $\mathbb{C}$, and on this variable we make a displacement, in order to move the roots out of the real axis, with which the LC method can operate without difficulties in the usual way.

For the case of quadratic polynomials, we saw that the LC method simply finds intersections between a line and a circumference in the complex plane $\mathbb{C}$, using of course the coefficients of the quadratic polynomial in question for the construction of these geometric elements; the intersections obtained in this way serve as the estimates of the roots. This simple geometric construction, however, lays the foundations for extending the LC method to cubic polynomials, where we introduced the concept of proximity map $e(\theta)$, generated by means of parallel processing techniques, and whose intersections with horizontal axis $y = 0$ provide the estimates of the roots. This extension to the case $n = 3$, in turn, lays the foundations for extending the LC method towards univariate polynomials of degree $n \geq 4$: this is achieved by minimizing (in parallel) dynamic squared distances $d_\theta^2(t)$ between terminal curves and circumferences in $\mathbb{C}$, with which we finally managed to construct the corresponding proximity maps. At first, we only had contemplated generating the weighted error proximity map $e(\theta)$; however, the theoretical-practical development carried out in this work showed us that it is also possible to generate the proximity maps $\frac{d}{d\theta}\, d_\theta^2(t^*)$ and $\frac{d}{d\theta}\, t^*$, obtained from the concept of dynamic squared distance.

The optimization algorithm that we implemented in this work (see annex 2 section 4) to find the minimizing arguments $t^*$ of the dynamic squared distance functions $d_\theta^2(t)$, turned out to be suitable for our purposes, which focused primarily on testing the concepts and functionality of the LC method; however, we also found a deficiency in this implementation, which specifically manifests itself by the presence of gaps in some computed proximity maps; from here, we see that it is necessary to improve this initial proposal if we want to have a more efficient and accurate optimization algorithm. Likewise, the routines and scripts in R that we used to produce the results in the numerical examples of this work, and whose implementations are based on the vectorization of operations, were designed mainly to illustrate, on a typical personal computer, the parallel processing potential of the LC method, without attempting to achieve a general-purpose implementation. An implementation of the LC method that substantially takes advantage of the parallel processing capabilities of a modern advanced hardware configuration will necessarily require data structures and software matched to the characteristics of the hardware where an efficient implementation of the LC method is sought.

As future work, it remains the task to strengthen the basic LC method described in this work, in order to use it as part of a complete, general-purpose system for solving polynomial equations. Below, we list some lines of work that should be carried out in order to achieve important progress in this task:

- Test whether the accuracy of initial root estimates generated by discrete proximity maps improves with the use of cubic splines instead of rectilinear segments, when interpolating discrete functional values $e(\theta_k)$, $\frac{d}{d\theta}d^2_{\theta_k}(t^*(\theta_k))$ or $\frac{d}{d\theta}t^*(\theta_k)$.

- Devise a strategy to change the scale and position of the variable from the input polynomial, so that the expected location of its roots on a specific region in the complex plane $\mathbb{C}$ favors the numerical stability of the LC method. In the end, an inverse transformation would be made on the root estimates obtained by the LC method, in order to restore the original scale and position of the polynomial variable.

- Design a sampling strategy in the global region $[-\pi,\pi)$, in order to obtain a distribution of angular values $\theta_k$ such that their density is greater only where theta roots are suspected to exist (as in transient peaks, often present in proximity maps); with this, we would obtain a dynamic and more efficient distribution of values $\theta_k$ within $[-\pi,\pi)$, relative to the uniform distribution of equidistant values $\theta_k$. Clearly, this line of work requires further investigation.

- Improve the optimization process that approximates the minimizing arguments $t^*$ of function $d^2_\theta(t)$; it is necessary to use a better alternative to the simulated annealing algorithm (`SANN`) in the first phase of our optimization process implemented in this work, in order to reduce the probability of getting stuck at local minima of function $d^2_\theta(t)$; there are several optimization alternatives that we could try, some of these based on evolutionary heuristics, such as genetic algorithms (GA), particle swarm optimization (PSO), differential evolution (DE), etc. In the first part of numerical example 5.3 in this chapter, we saw that function $d^2_{\theta^*}(t)$ can have more than one minimizing argument $t^*$ for an angular value $\theta^*$ where in fact $d^2_{\theta^*}\big(t^*(\theta^*)\big)=0$; additionally, a solution to this problem was also proposed, without attempting a computational implementation: if we "discretize" the function $\frac{d}{dt}d^2_\theta(t)$ for a given value $\theta$, we can estimate one or more possible values $t^*$, using a strategy similar to the one utilized throughout this work to find crossings of proximity maps $e(\theta)$, $\frac{d}{d\theta}d^2_\theta(t^*)$, or $\frac{d}{d\theta}t^*$ with horizontal axis $y=0$ (of course, it is also necessary to verify that the possible estimated values $t^*$ do not correspond to local minima). It would be worth investigating whether this proposed optimization-by-discretization strategy is a better alternative to the two-phase optimization strategy we used, where the first phase uses the simulated annealing algorithm `SANN` to obtain an initial estimate of $t^*$, and the second phase uses the BFGS algorithm to polish the initial estimate obtained in the first phase. In principle, the optimization-by-discretization proposal seems to be a promising alternative, since it attacks the problem of multiple values $t^*$, which is not considered in the two-phase optimization strategy, where only a minimizing value $t^*$ is sought.

- Reduce the amount of noise present in discrete proximity maps, especially in maps $\frac{d}{d\theta}t^*$; this could be achieved, for example, through the use of smoothing filters, frequently used in signal processing applications, or in time series trend analysis.

- Delimit regions of the complex plane $\mathbb{C}$ where roots of a univariate polynomial could be located, based on information from its coefficients; see for example [5.2]. This would further improve the efficiency of the optimization method that searches for optimal arguments $t^*$ of function $d_\theta^2(t)$, as this would reduce the size of the search space $t$ for function $d_\theta^2(t)$.

- Extend the basic LC method to properly address cases of univariate polynomial equations with multiple roots, or with roots very close to each other, whether they are theta roots $\theta_i^*$ in proximity maps $e(\theta)$, $\frac{d}{d\theta}d_\theta^2(t^*)$ and $\frac{d}{d\theta}t^*$, or the roots $R_i$ in $\mathbb{C}$ themselves. Here it would be very useful, once again, to improve, among other things, our optimization algorithm so that it can find (if any) multiple optimal arguments $t^*$ of function $d_\theta^2(t)$.

- Implement an extended and robust version of the LC method on hardware architectures that support massively parallel processing, such as those mentioned in note 3.1 of chapter 3, and preferably using more efficient programming languages, such as C, C++, or Fortran. In this work, the R language was used for illustrative purposes only; in no way did we seek to achieve an implementation that could compete with designs built with more efficient and faster programming languages; we only used vectorization of functions, in order to improve the performance of our routines written in R on a personal computer. Ideally, such an extended version of the LC method would use the proposed strategy of optimization by discretization of function $\frac{d}{dt}d_\theta^2(t)$ in order to find possible multiple values $t^*$ for a given angle $\theta$, via a scheme of vectorization of operations *within* each of the $N$ processors that handle each of the structures $LzC(\theta)$; in other words, a hardware architecture with two levels of parallel processing (a first level to build structures $LzC(\theta)$ and a second level to find the global minima of corresponding functions $d_\theta^2(t)$) would in principle provide an ideal platform to obtain proximity maps $e(\theta)$, $\frac{d}{d\theta}d_\theta^2(t^*)$ or $\frac{d}{d\theta}t^*$ in the shortest possible time.

- Carry out a detailed spatiotemporal complexity analysis for the LC method, in order to gain a deeper understanding of the relationship between the LC method and other numerical methods that address various aspects of the root finding problem for polynomial equations of a single variable. Informally, we could say that if the LC method is implemented on a machine with massively parallel processing capabilities, its spatial complexity will be relatively large, although manageable, given the capabilities of modern hardware, while its temporal complexity will be relatively small; this may be an incentive to use the LC method as a preprocessing module that identifies initial estimates of polynomial roots, within a robust, general-purpose system that finds roots of polynomial equations of one variable.

In this work we have seen that the basic LC method, under normal circumstances, does a good job of providing initial estimates that can then be refined by means of iterative algorithms with rapid local convergence, such as Newton-Raphson. An additional scenario is to consider the functional values of any of the proximity maps $e(\theta)$, $\frac{d}{d\theta}d_\theta^2(t^*)$ or $\frac{d}{d\theta}t^*$ as an alternative to polynomial evaluations $p(x)$ within existing iterative methods, taking $\theta$ as an independent variable instead of $x$; this would also be a possible line of future research. Since functional values of proximity maps give us an idea of how close we are to a possible root, they could be used within an iterative method to make decisions about how to carry out the next iteration in a more efficient way.

# PART II: Supplementary Material

# Annex 1: Basic R Functions to Implement the LC Method

## 1. Function to Find the Center of a Z-circumference

### Source Code

```
center_zcircle <- function(p1,p2) {

   # computes the coordinates h + ki in the complex plane that correspond
   # to the center of a z-circumference, given two points p1 and p2
   # from that z-circumference

   x1 <- Re(p1)
   y1 <- Im(p1)

   x2 <- Re(p2)
   y2 <- Im(p2)

   h <- (y2*((x1^2)+(y1^2))-y1*((x2^2)+(y2^2)))/(2*x1*y2-2*y1*x2)
   k <- (x1*((x2^2)+(y2^2))-x2*((x1^2)+(y1^2)))/(2*x1*y2-2*y1*x2)

   return( h + k*1i )
}
```

### Description

For a z-circumference $C_z(t) = b/(p + tv)$, where $b, p, v \in \mathbb{C}$ are fixed points, and $t \in \mathbb{R}$ is a variable parameter, we have that when $t \to \pm\infty$, $C_z(t) \to 0 + i0$, although $C_z(t)$ never takes the value $0 + i0$. For computational implementation purposes, if we know two non-zero points from $C_z(t)$, say, $p_1, p_2 \in \mathbb{C}$, then we can assume that $0 + i0$ is also contained in z-circumference $C_z(t)$; this means that we know 3 distinct points of a circumference, and, therefore, we can obtain the coordinates of its center $c = h + ik$ by solving a system of 3 equations in 3 unknowns: $h$, the abscissa of $c$; $k$, the ordinate of $c$; and $r$, the radius of $C_z(t)$. If

$p_1 = x_1 + iy_1,$

$p_2 = x_2 + iy_2,$

then the system of equations is posed as follows (see figure A1.1 to better understand the origin of these equations):

$$
\begin{aligned}
(x_1 - h)^2 + (y_1 - k)^2 &= r^2 \\
(x_2 - h)^2 + (y_2 - k)^2 &= r^2 \qquad \text{(A1.1)} \\
(0 - h)^2 + (0 - k)^2 &= r^2
\end{aligned}
$$

We carry out the operations indicated by equations (A1.1):

$$x_1^2 - 2x_1h + h^2 + y_1^2 - 2y_1k + k^2 = r^2$$
$$x_2^2 - 2x_2h + h^2 + y_2^2 - 2y_2k + k^2 = r^2 \qquad \text{(A1.2)}$$
$$h^2 + k^2 = r^2$$

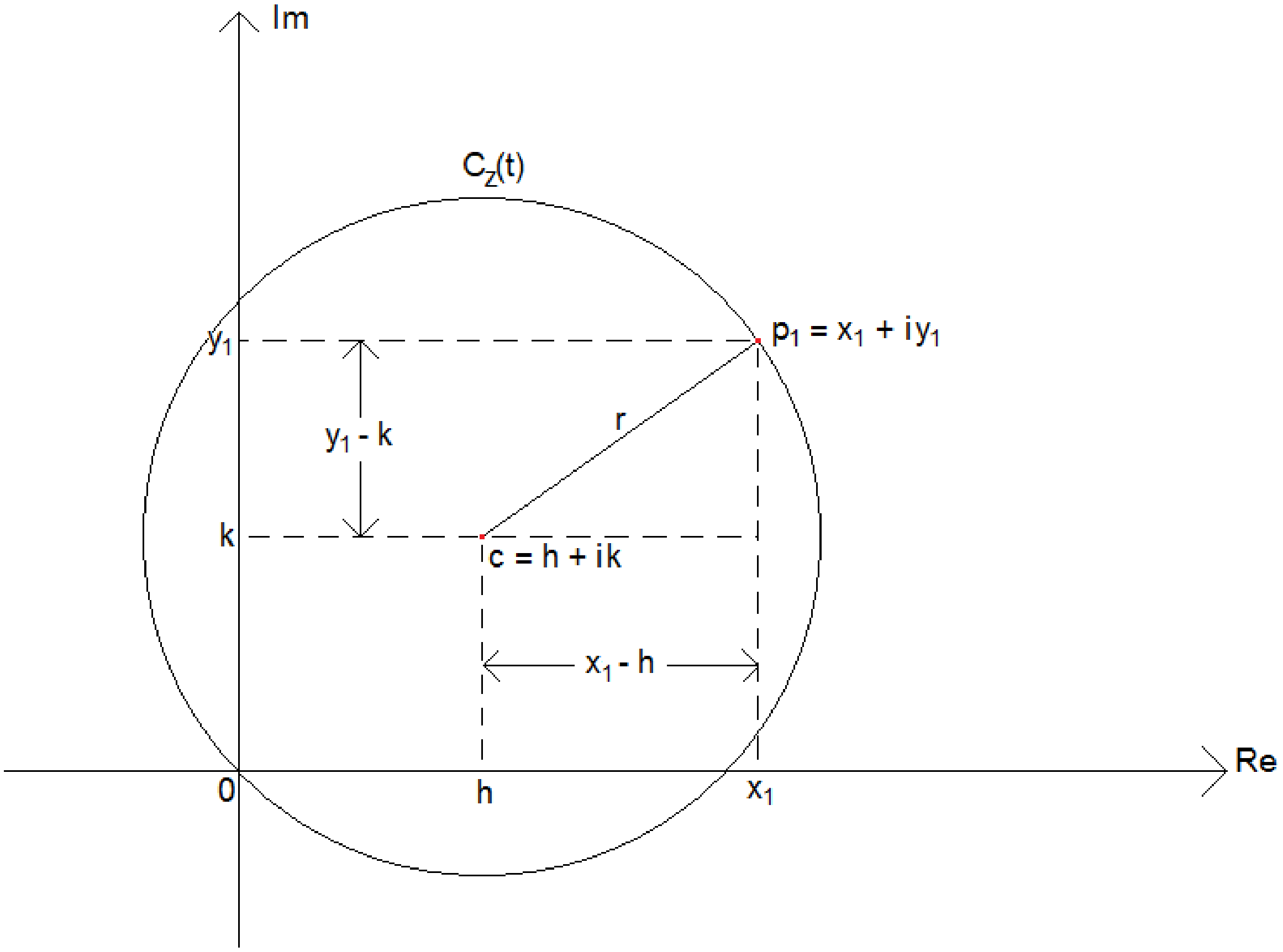

**Fig. A1.1**. A z-circumference $C_z(t)$ and a point $p_1$ contained in $C_z(t)$. The right triangle whose vertices are $c$, $p_1$, and $x_1 + ik$ is the key to stating each of the equations in (A1.1), considering the relationship between the length of its sides given by the Pythagorean theorem. Each of the equations in (A1.1) corresponds to the three known points of $C_z(t)$: $p_1, p_2$ and $0 + i0$.

If we subtract the last equation from the first two in equation system (A1.2), we reduce the problem to a system of two linear equations in 2 unknowns, $h, k$:

$$2x_1h + 2y_1k = x_1^2 + y_1^2 \qquad \text{(A1.3)}$$
$$2x_2h + 2y_2k = x_2^2 + y_2^2$$

Now, we proceed to use Gaussian elimination in order to solve equation system (A1.3).

**Step 1**: The first equation in (A1.3) is divided by $2x_1$, which give us the first equation in (A1.4); then, to the second equation in (A1.3) we add $-2x_2$ times the first equation in (A1.4); with this, we complete the equivalent system (A1.4):

$$h + (y_1/x_1)k = (x_1^2 + y_1^2)/2x_1 \qquad \text{(A1.4)}$$

$$\left(\frac{2x_1y_2 - 2x_2y_1}{x_1}\right)k = (x_2^2 + y_2^2) - x_2(x_1^2 + y_1^2)/x_1$$

**Step 2**: Divide the second equation in (A1.4) by $\left(\frac{2x_1y_2-2x_2y_1}{x_1}\right)$, and solve for $h$ in the first equation in (A1.4), with which we obtain the equivalent system (A1.5):

$$h = \frac{x_1^2+y_1^2-2y_1k}{2x_1} \qquad \text{(A1.5)}$$

$$k = \frac{x_1(x_2^2+y_2^2)-x_2(x_1^2+y_1^2)}{2x_1y_2-2x_2y_1}$$

**Step 3**: Plug the value $k$ indicated by the second equation in (A1.5) into the first equation in (A1.5), with which we finally obtain the expressions for $h, k$ in terms of known values:

$$h = \frac{-y_1(x_2^2+y_2^2)+y_2(x_1^2+y_1^2)}{2x_1y_2-2x_2y_1} \qquad \text{(A1.6)}$$

$$k = \frac{x_1(x_2^2+y_2^2)-x_2(x_1^2+y_1^2)}{2x_1y_2-2x_2y_1}$$

The equations in (A1.6) are the coordinates sought. These formulas are implemented directly into the source code of the function `center_zcircle`.

**Remarks on the source code:**

1) `p1`, `p2` can each be a vector of $N$ complex numbers, so the function `center_zcircle` processes "in parallel" $N$ z-circumferences (subject to hardware resources available to the R programming language, and under a scheme of vectorization of operations).

2) The points `p1`, `p2`, as well as the computed centers of z-circumferences, are treated as complex numbers; this is perfectly compatible with the notions of points and vectors in the Cartesian plane $\mathbb{R}^2$.

3) This function fails if `p1`, `p2` are two complex numbers such that `Re(p1)*Im(p2)=Re(p2)*Im(p1)`. This occurs, for example, if the z-circumference in question is of infinite radius (i.e., if the z-circumference degenerates into a line, which occurs if the generating line $p + tv$ for z-circumference $C_z(t)$ contains zero).

## 2. Function to Find Intersections Between Semi-lines and Circumferences

### Source Code

```
intersect_semiLine_circle <- function(q,v,c_c,r) {
   # computes the intersections, if any,
   # between a semi-line and a circumference in the complex plane
   # INPUTS:        q   : anchor point of the semi-line
   #                v   : unit direction vector of the semi-line
   #                c_c : center of the circumference
   #                r   : radius of the circumference
   # OUTPUTS:       I1  : intersection I1 (complex value)
   #                I2  : intersection I2 (complex value)
   #                nI  : indicates which Intersections are relevant
   #                      nI = 0 : I1 and I2 are both irrelevant,
   #                               or there are no intersections
   #                      nI = 1 : only one intersection is relevant,
   #                      nI = 2 : I1 and I2 are both relevant
   # NOTE: Outputs are packed into matrix OUT

   # the first two columns in OUT store the intersections
   # and the 3rd column is the relevance status for intersections
   OUT   <- matrix(0,nrow=length(q),ncol=3)
   nI    <- rep(NA,length(q))   # 3rd column in OUT; we store here
                                # the relevance status for I1 and I2

   # we make sure that v is a unit vector
   v = v/abs(v)

   # vector from q to c_c
   c_q <- c_c - q

   # coordinates of vector c_c - q
   Re_c_q <- Re(c_q)
   Im_c_q <- Im(c_q)

   # coordinates of semi-line’s direction vector
   vr <- Re(v)
   vi <- Im(v)

   # (cc - q) dot v
   c_q_dot_v <- Re_c_q*vr + Im_c_q*vi

   # squared norm of c_c - q
   c_q_dot_c_q <- Re_c_q*Re_c_q + Im_c_q*Im_c_q

   # discriminant
   delta <- c_q_dot_v^2 + r^2 - c_q_dot_c_q

   # parametric values t1, t2 of the semi-line where
   # intersections with the circumference occur
   t1 <- rep(NA,length(q))
   t2 <- rep(NA,length(q))

   i     <- which( delta>=0 )
   t1[i] <- c_q_dot_v[i] - sqrt(delta[i])
   t2[i] <- c_q_dot_v[i] + sqrt(delta[i])
```

```
   # intersections between semi-line and circumference
   I1 <- q + t1*v
   I2 <- q + t2*v

   # RELEVANCE CONDITIONS
   # I1 and I2 are both irrelevant, or there are no intersections
   i0 <- which( (t1<0 & t2<0) | is.na(t1) )

   # only one intersection is relevant (I2)
   i1 <- which( t1<0 & t2>=0 )

   # both intersections are relevant (I1 and I2)
   i2 <- which( t1>=0 & t2>=0 )

   # assigns relevance values
   nI[i0] <- 0
   nI[i1] <- 1
   nI[i2] <- 2

   # prepares output
   OUT[,1] <- I1
   OUT[,2] <- I2
   OUT[,3] <- nI

   return( OUT )
}
```

### Description

This function `intersect_semiLine_circle`, like the function `center_zcircle`, takes advantage of R's capabilities to vectorize operations, in order to process "in parallel" $N$ cases ($N \geq 1$) where it is required to find the intersections between a semi-line and a circumference. The following is a description of how the function `intersect_semiLine_circle` processes each one of the $N$ cases.

The function `intersect_semiLine_circle` takes as input the radius $r$ and center $c = h + ik$ of a circumference, as well as the starting point (anchor point) $q = q_r + iq_i$ and unit direction vector $v = v_r + iv_i = \cos\theta + i\sin\theta$ of a semi-line; $\theta \in [-\pi, \pi)$. The goal of this function of course is to find the intersections $I_1$ and $I_2$ (if they exist at all) between a semi-line and a circumference, both specified by the input parameters for the function. Figure A1.2a illustrates a situation where there are two *relevant* intersections between a semi-line and a circumference.

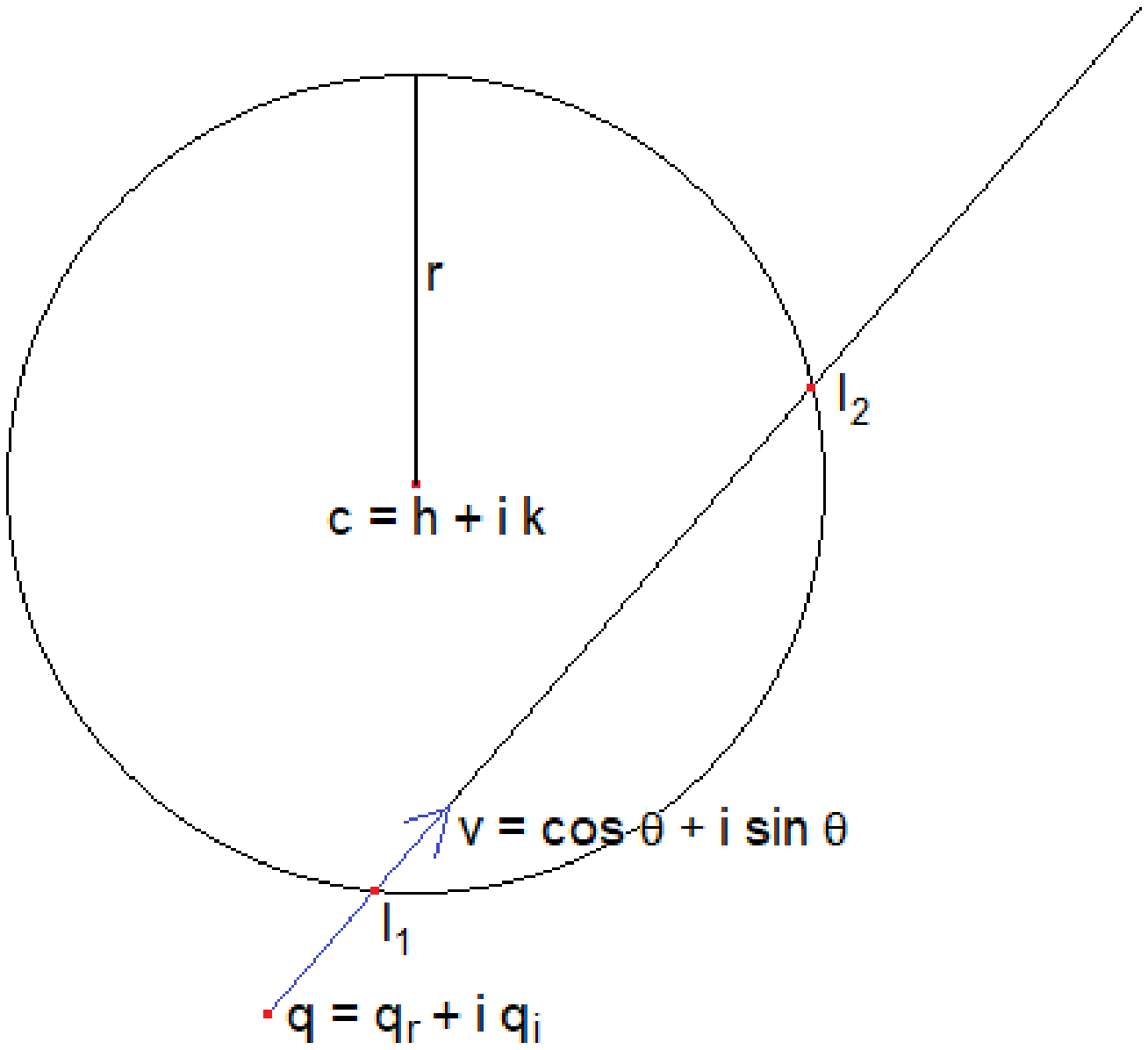

**Fig. A1.2a**. Semi-line $q + tv$, $t \in \mathbb{R}^+ \cup \{0\}$, intersecting at two relevant points $I_1$ and $I_2$ with the circumference whose equation in $\mathbb{R}^2$ is $(x - h)^2 + (y - k)^2 = r^2$.

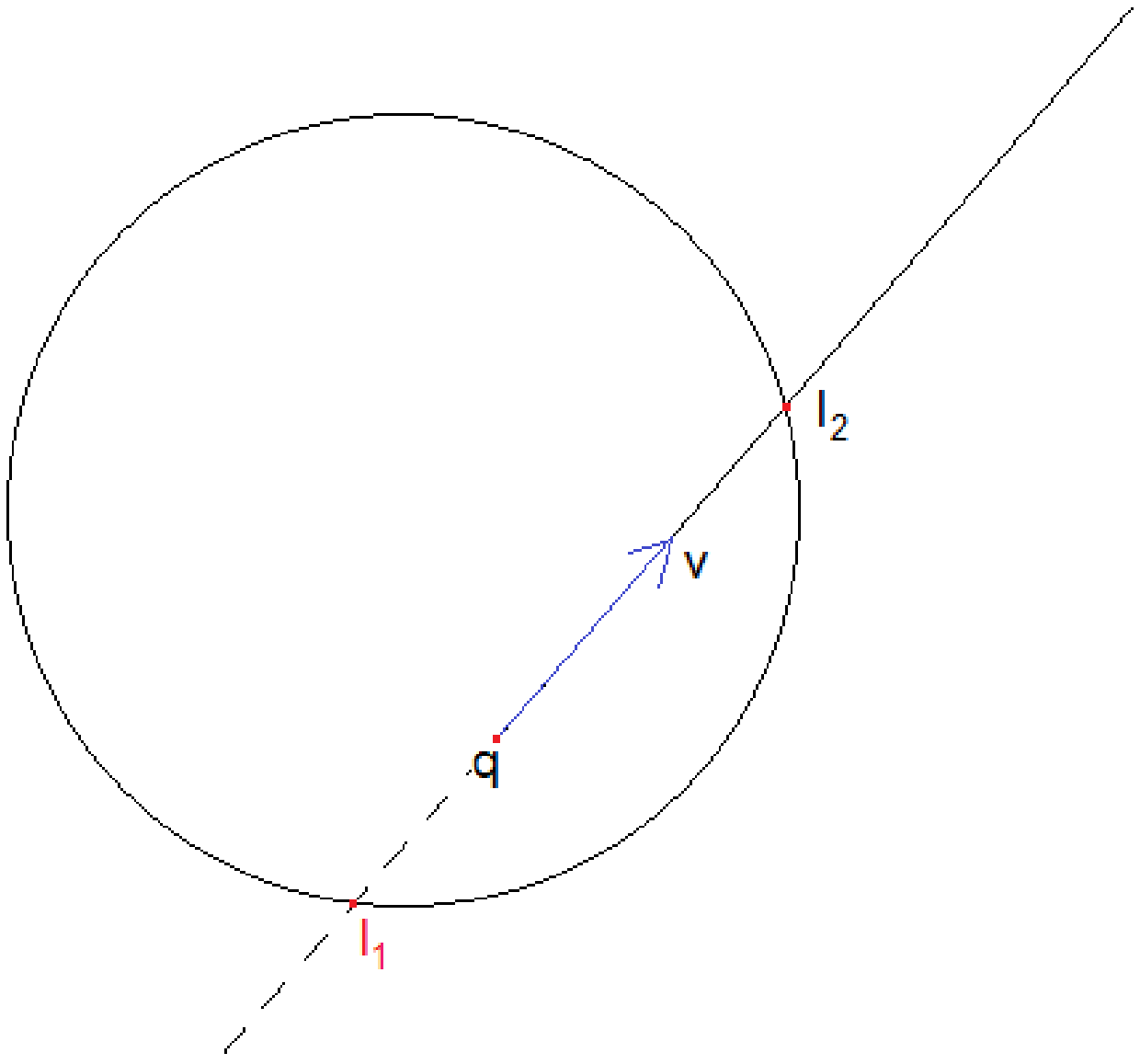

**Fig. A1.2b**. The line parallel to vector $v$ and containing point $q$ intersects at two points with a circumference, but only one intersection is relevant, namely, the intersection between the corresponding semi-line (with anchor point $q$ and direction vector $v$) and the circumference.

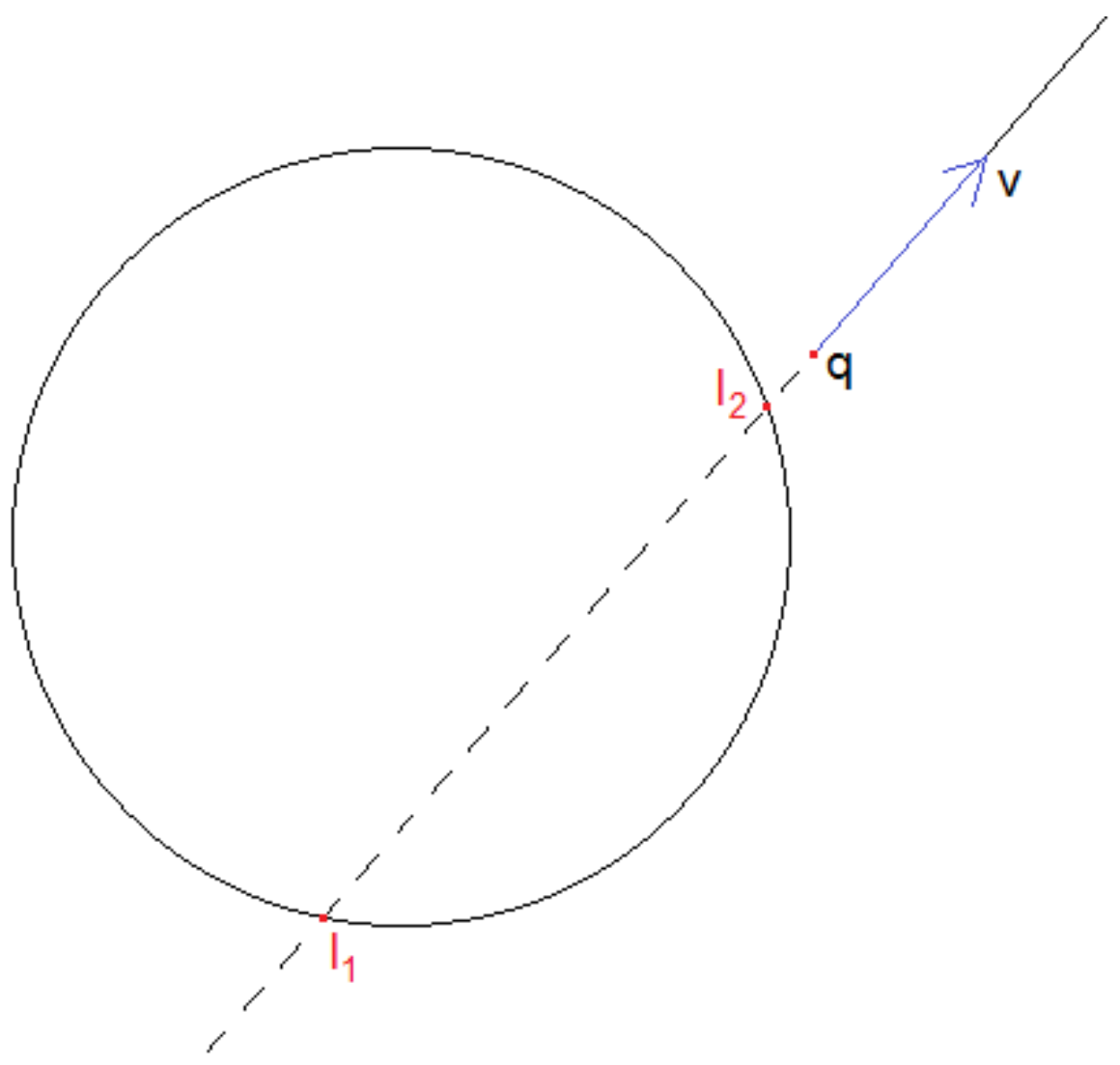

**Fig. A1.2c**. The line parallel to vector $v$ and containing point $q$ intersects at two points with a circumference; however, no intersection is relevant, since the corresponding semi-line (with anchor point $q$ and direction vector $v$) is outside the circumference.

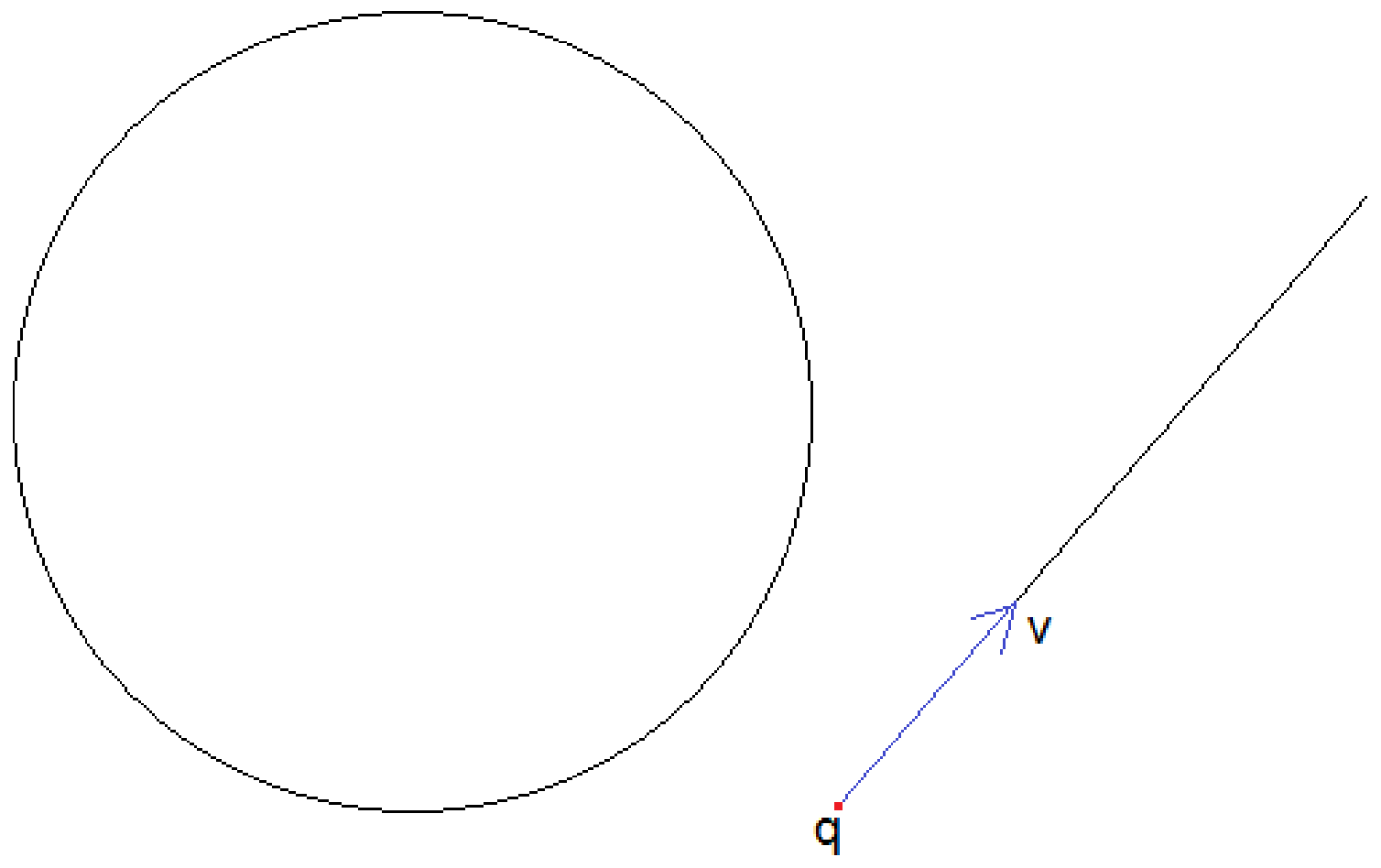

**Fig. A1.2d**. The line parallel to vector $v$ and containing point $q$ does not intersect with the circumference, so the corresponding semi-line (with anchor point $q$ and direction vector $v$) does not intersect with the circumference either.

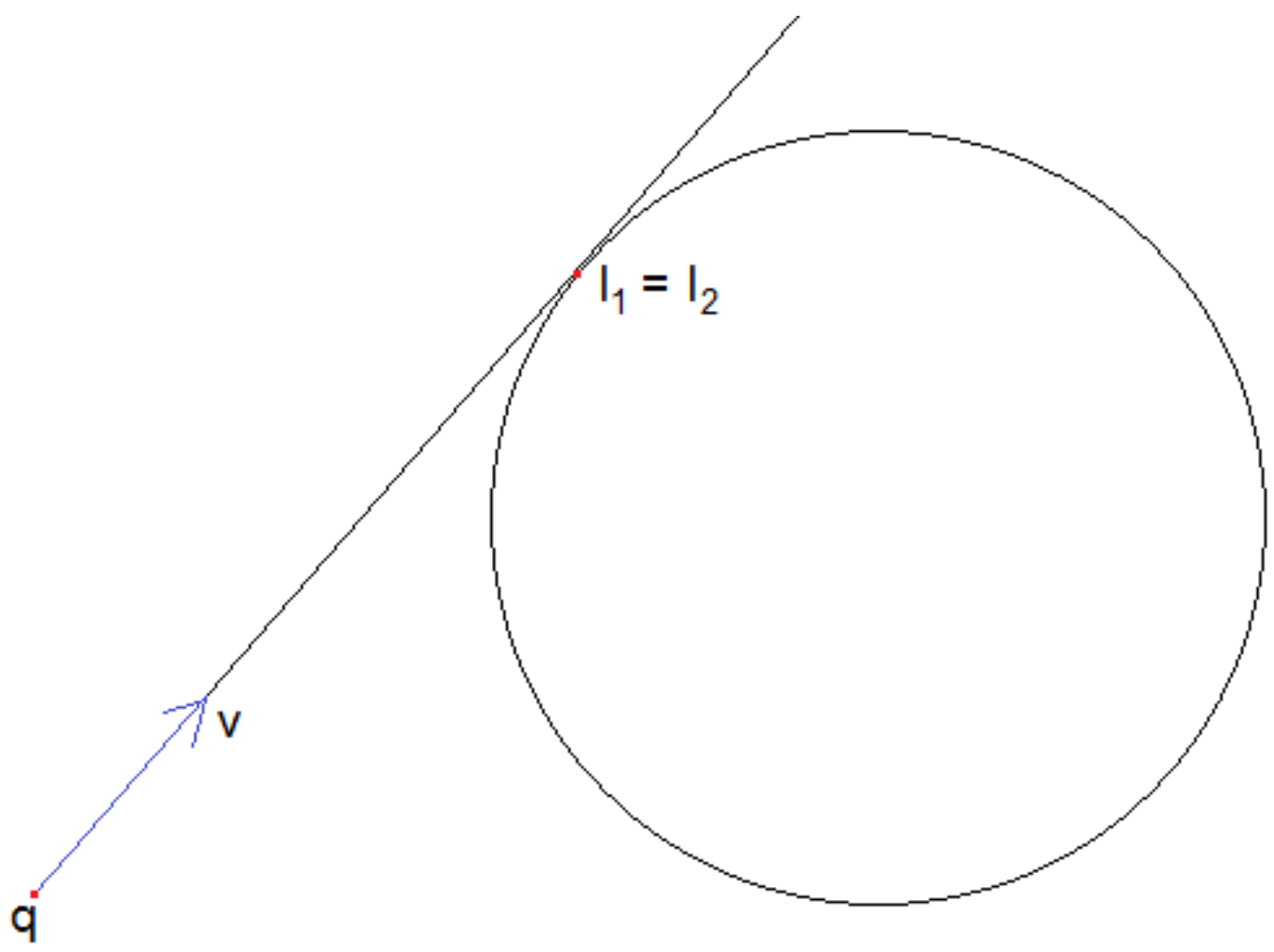

**Fig. A1.2e**. The line parallel to vector $v$ and containing point $q$, intersects tangentially with the circumference at a single point; the tangential intersection between the corresponding semi-line (with anchor point $q$ and direction vector $v$) and the circumference may or may not be relevant (in this illustration, the tangential intersection is relevant).

We say that an intersection $I = q + t_0 v$ between a semi-line and a circumference is relevant if $t_0 \geq 0$. With this in mind, we see that it is possible to have situations similar to those illustrated in figures A1.2b, A1.2c, A1.2d and A1.2e.

In order to find the coordinates of the possible intersections $I_1$ and $I_2$ between a semi-line and a circumference, we first consider the more general problem of finding the intersections between a line $q + tv$, $t \in \mathbb{R}$, and a circumference whose equation in $\mathbb{R}^2$ is $(x - h)^2 + (y - k)^2 = r^2$; once we have solved this more general problem, it will be easy to determine which of the intersections are relevant, with respect to the semi-line with anchor point $q$ and direction vector $v$.

Our goal then is to first find values $t_1$, $t_2$ such that $I_1 = q + t_1 v$, $I_2 = q + t_2 v$ are also in the circumference $(x - h)^2 + (y - k)^2 = r^2$. For this, we plug the parametric expressions of the line into the center-radius equation of the circumference:

$$[(q_r + t\cos\theta) - h]^2 + [(q_i + t\sin\theta) - k]^2 = r^2 \qquad \text{(A1.7)}$$

We solve for parameter $t$ in equation (A1.7):

$$(q_r + t\cos\theta)^2 - 2(q_r + t\cos\theta)h + h^2 + (q_i + t\sin\theta)^2 - 2(q_i + t\sin\theta)k + k^2 = r^2$$

$$q_r^2 + 2t\,q_r\cos\theta + t^2\cos^2\theta - 2hq_r - 2ht\cos\theta + h^2$$
$$+q_i^2 + 2t\,q_i\sin\theta + t^2\sin^2\theta - 2kq_i - 2kt\sin\theta + k^2 = r^2$$

$$(\cos^2\theta + \sin^2\theta)t^2 + 2(q_r\cos\theta + q_i\sin\theta)t - 2(h\cos\theta + k\sin\theta)t$$
$$+h^2 - 2q_r h + q_r^2 + k^2 - 2q_i k + q_i^2 - r^2 = 0$$

$$t^2 + 2[(q_r - h)\cos\theta + (q_i - k)\sin\theta]t + [(h - q_r)^2 + (k - q_i)^2 - r^2] = 0$$
$$t^2 + 2[(q - c)\cdot v]t + (|q - c|^2 - r^2) = 0$$
$$t^2 + 2[(q - c)\cdot v]t + [(q - c)\cdot(q - c) - r^2] = 0 \qquad \text{(A1.8)}$$

Equation (A1.8) is a polynomial equation of degree 2 in the variable $t$ (we would expect to see two real roots in this equation); if we set

$A = 1,$
$B = 2[(q - c) \cdot v],$
$C = [(q - c) \cdot (q - c) - r^2],$

then we have

$$t = (-B \pm \sqrt{B^2 - 4AC})/(2A)$$

$$t = (-2[(q - c) \cdot v] \pm \sqrt{4[(q - c) \cdot v]^2 - 4[(q - c) \cdot (q - c) - r^2]})/2$$

$$t = -[(q - c) \cdot v] \pm \sqrt{[(q - c) \cdot v]^2 - [(q - c) \cdot (q - c) - r^2]}$$

$$t = [(c - q) \cdot v] \pm \sqrt{[(c - q) \cdot v]^2 - (c - q) \cdot (c - q) + r^2} \qquad \text{(A1.9)}$$

We rewrite equation (A1.9) as

$$t_1 = [(c - q) \cdot v] - \sqrt{\Delta} \qquad \text{(A1.10)}$$

$$t_2 = [(c - q) \cdot v] + \sqrt{\Delta} \qquad \text{(A1.11)}$$

where $\Delta = [(c - q) \cdot v]^2 - (c - q) \cdot (c - q) + r^2$ (A1.12)

From expressions (A1.10), (A1.11) and (A1.12) we see that if $\Delta > 0$, then $t_1 < t_2$. In this case, if $t_1 \geq 0$, we have two relevant intersections, as illustrated in figure A1.2a; if $t_1 < 0$ and $t_2 \geq 0$, we have only one relevant intersection, as illustrated in figure A1.2b; if $t_1 < t_2 < 0$, we have no relevant intersections, as illustrated in figure A1.2c.

If $\Delta < 0$ in expression (A1.12), this means that there are no intersections between the semi-line and the circumference, as illustrated in figure A1.2d. On the other hand, if $\Delta = 0$, $t_1 = t_2$; in such a case, if $t_1 \geq 0$, then the semi-line and the circumference coincide at a single relevant tangential intersection, as illustrated in figure A1.2e.

All of these considerations are encoded in the function `intersect_semiLine_circle`, as can be seen in the source code listed above.

## 3. Function to Obtain Weighted Error Associated with Structure $LzC(\theta)$

### Source Code

```
error_tl_zc <- function(p,v,ptl,pzc) {

   # given an intersection between a terminal semi-line and a
   # z_circumference, this function finds the normalized difference
   # (or weighted error) between the intersection's projection via
   # the terminal semi-line and the intersection's projection via
   # the z-circumference onto line L1 in the complex plane
   # INPUTS:      p   : fixed point of line L1
   #              v   : unit direction vector of L1
   #              ptl : intersection's projection via
   #                    terminal semi-line onto L1
   #              pzc : intersection's projection via
   #                    z-circumference onto L1
   # OUTPUT:      dif : normalized difference

   dif <- sign(Re((ptl-pzc)/v))*abs(ptl-pzc)/(abs(ptl-p)+abs(pzc-p))
   return( dif )
}
```

### Description

This function simply carries out the calculations indicated by expression (3.3) of chapter 3, which is completely analogous to expressions (3.4a) and (3.4b) also in chapter 3, as well as to expressions (4.11a) and (4.11b) in chapter 4, and to expressions (5.16a) and (5.16b) in chapter 5. Despite its simplicity, the function `error_tl_zc` is fundamental in the construction of *discrete* angular proximity maps associated with structures $LzC(\theta)$, regardless of the degree of the univariate polynomial in question. The respective input parameters of function `error_tl_zc` are, according to the concepts developed in chapters 3, 4 and 5, complex values $P_1$, $v_\theta$, $i_0$ and $c_0$. In the same way, the output value of this function is the so-called *weighted error* (also called *normalized difference*) $e(\theta)$. As we saw in chapter 3, this output value is always between -1 and 1, because of the way the weighted error $e(\theta)$ is constructed. Note that, like the functions described in the two previous sections of this annex, function `error_tl_zc` is capable of processing multiple cases "in parallel" by means of the vectorization of operations in R; this capability depends on the physical characteristics of the hardware that is executing this function within the R language (amount of RAM, number of processing cores).

The function `error_tl_zc` is useful in obtaining discrete approximations of the weighted error functions $e_A(\theta)$ and $e_B(\theta)$, which are used in the numerical examples of chapters 3, 4 and 5.

# 4. Function to Approximate Theta Roots

## Source Code

```
get_angle_roots <- function( TH, E, TOL ) {

   # Finds the crossings of function E(TH) with horizontal axis y=0.
   # Crossings are approximated via linear interpolation
   # INPUTS :   TH  : array of discrete values for the independent
   #                  variable
   #            E   : image of discrete values TH under the real
   #                  function of real variable E
   #            TOL : maximum vertical distance allowed between
   #                  two discrete consecutive values of E with
   #                  opposite signs, involved in a smooth
   #                  crossing of E(TH) with horizontal axis y=0
   # OUTPUTS :  th  : abscissae of the smooth crossings of E(TH)
   #                  with y=0 obtained by linear interpolation
   #            de  : absolute vertical distances from consecutive
   #                  elements of E involved in smooth crossings;
   #                  these distances are not greater than TOL
   # NOTE: outputs are packed into a list L

   N <- length(TH)
   i <- 1:N

   ic <- which( E[i]*E[i+1]<=0 & abs(E[i]-E[i+1])<=TOL )

   th1 <- TH[ic]
   th2 <- TH[ic+1]
   e1  <- E[ic]
   e2  <- E[ic+1]

   th <- th1 + e1*(th2-th1)/(e1-e2)
   de <- abs( e1-e2 )

   L      <- list()
   L[[1]] <- th
   L[[2]] <- de
   return( L )
}
```

## Description

This function receives as input a discrete mapping $\vec{\boldsymbol{\theta}} \to e(\vec{\boldsymbol{\theta}})$, with $\vec{\boldsymbol{\theta}}, e(\vec{\boldsymbol{\theta}}) \in \mathbb{R}^N$, where $N$ is the number of elements in vectors $\vec{\boldsymbol{\theta}} = \langle \theta_1, \theta_2, \ldots, \theta_N \rangle$, $e(\vec{\boldsymbol{\theta}}) = \langle e_1, e_2, \ldots, e_N \rangle$. From the source code, we easily identify that `TH` is $\vec{\boldsymbol{\theta}}$, and `E` is $e(\vec{\boldsymbol{\theta}})$. From the received discrete mapping, function `get_angle_roots` attempts to approximate the crossings of the corresponding continuous mapping $\theta \to e(\theta)$ with horizontal axis $y = 0$ by means of linear interpolation. Initially, function `get_angle_roots` identifies consecutive items from the input discrete array $e(\vec{\boldsymbol{\theta}})$ that meet the condition

$$e_i e_{i+1} \le 0 \text{ and } |e_i - e_{i+1}| \le \text{TOL}, \quad i = 1,2, \ldots, N-1. \qquad \text{(A1.13)}$$

That is to say, `get_angle_roots` identifies consecutive discrete values from $e(\vec{\boldsymbol{\theta}})$, $e_i$, $e_{i+1}$, for which there is a change of sign (so between them there must be at least a crossing of $e(\theta)$ with horizontal axis $y = 0$) and also for which their absolute difference is not greater than input parameter `TOL`; with this last restriction we try to make an (informal) filtering in which we only consider changes of sign of $e(\theta)$ associated with moderate and not abrupt variations of this function, as we would expect the behavior of a real continuous and differentiable (smooth) function to be.

Indices $i$ that meet condition (A1.13) are stored in vector `ic`. These indices in `ic` are then used to extract items $\theta_i$, $\theta_{i+1}$, $e_i$, $e_{i+1}$ (from input vectors $\vec{\boldsymbol{\theta}}$, $e(\vec{\boldsymbol{\theta}})$) that are in a situation similar to that shown in figure A1.3. With these items, we want to calculate the value $\theta^*$, which is the crossing of the linear segment that goes from point $(\theta_i, e_i)$ to point $(\theta_{i+1}, e_{i+1})$, with horizontal axis $y = 0$ (note that we are implicitly assuming that between $\theta_i$ and $\theta_{i+1}$ there is only one intersection of $e(\theta)$ with horizontal axis $y = 0$; this assumption will be valid as long as the continuous mapping $\theta \rightarrow e(\theta)$ is properly discretized). This can be easily achieved by remembering the point-slope form for the equation of a line:

$$y - y_1 = m(x - x_1). \qquad \text{(A1.14)}$$

If we link the elements shown in figure A1.3 to those of equation (A1.14), we have

$$e(\theta) - e(\theta^*) = \frac{e_{i+1} - e_i}{\theta_{i+1} - \theta_i}(\theta - \theta^*). \qquad \text{(A1.15)}$$

By assumption, we have $e(\theta^*) = 0$, so, if we solve for $\theta^*$ in (A1.15), we obtain

$$\theta^* = \theta + e(\theta)\frac{\theta_{i+1} - \theta_i}{e_i - e_{i+1}}. \qquad \text{(A1.16)}$$

Equation (A1.16) is valid for any point in the rectilinear segment with endpoints $(\theta_i, e_i)$ and $(\theta_{i+1}, e_{i+1})$; for simplicity's sake, we substitute the indeterminate values in equation (A1.16) with the coordinates of the left endpoint of this rectilinear segment, with which we finally obtain

$$\theta^* = \theta_i + e_i\frac{\theta_{i+1} - \theta_i}{e_i - e_{i+1}}. \qquad \text{(A1.17)}$$

Equation (A1.17) corresponds to instruction `th <- th1 + e1*(th2-th1)/(e1-e2)` within function `get_angle_roots`. With this, we obtain values $\theta^*$, which are stored in a list, along with corresponding absolute vertical differences $|e_i - e_{i+1}|$; this list is the output of function `get_angle_roots`.

Note 1: In the graph from figure A1.3, the line segment with endpoints $(\theta_i, e_i)$ and $(\theta_{i+1}, e_{i+1})$ has a positive slope; however, the above analysis is equally valid, regardless of the slope's sign.

Note 2: Like all the functions described in this annex, function `get_angle_roots` takes advantage of R's capabilities to vectorize operations, so that the computations are carried out "in parallel", subject to available hardware capabilities.

Note 3: `get_angle_roots` will have difficulties in approximating the zeros of $e(\theta)$ when this function does not exhibit changes of sign throughout its domain (such as the function $e(\theta) = \theta^2$).

Note 4: Initially, function `get_angle_roots` was designed to find the theta roots of the discretized version of the weighted error function $e(\theta)$. More generally, `get_angle_roots` can be used to approximate the crossings, with horizontal axis $y = 0$, of **any** discretized mapping obtained by sampling some real function of one real variable, preferably continuous and differentiable. Within this project, `get_angle_roots` is also useful for approximating theta roots of discrete mappings $\frac{d}{d\theta} d_\theta^2(t^*(\theta))$ and $\frac{d}{d\theta} t^*(\theta)$. Potentially, function `get_angle_roots` could also be used within this project as an auxiliary tool to find the minimizing arguments $t^*$ of function $d_\theta^2(t)$, by applying `get_angle_roots` to a discrete version of the derivative $\frac{d}{dt} d_\theta^2(t)$ for a given value $\theta$.

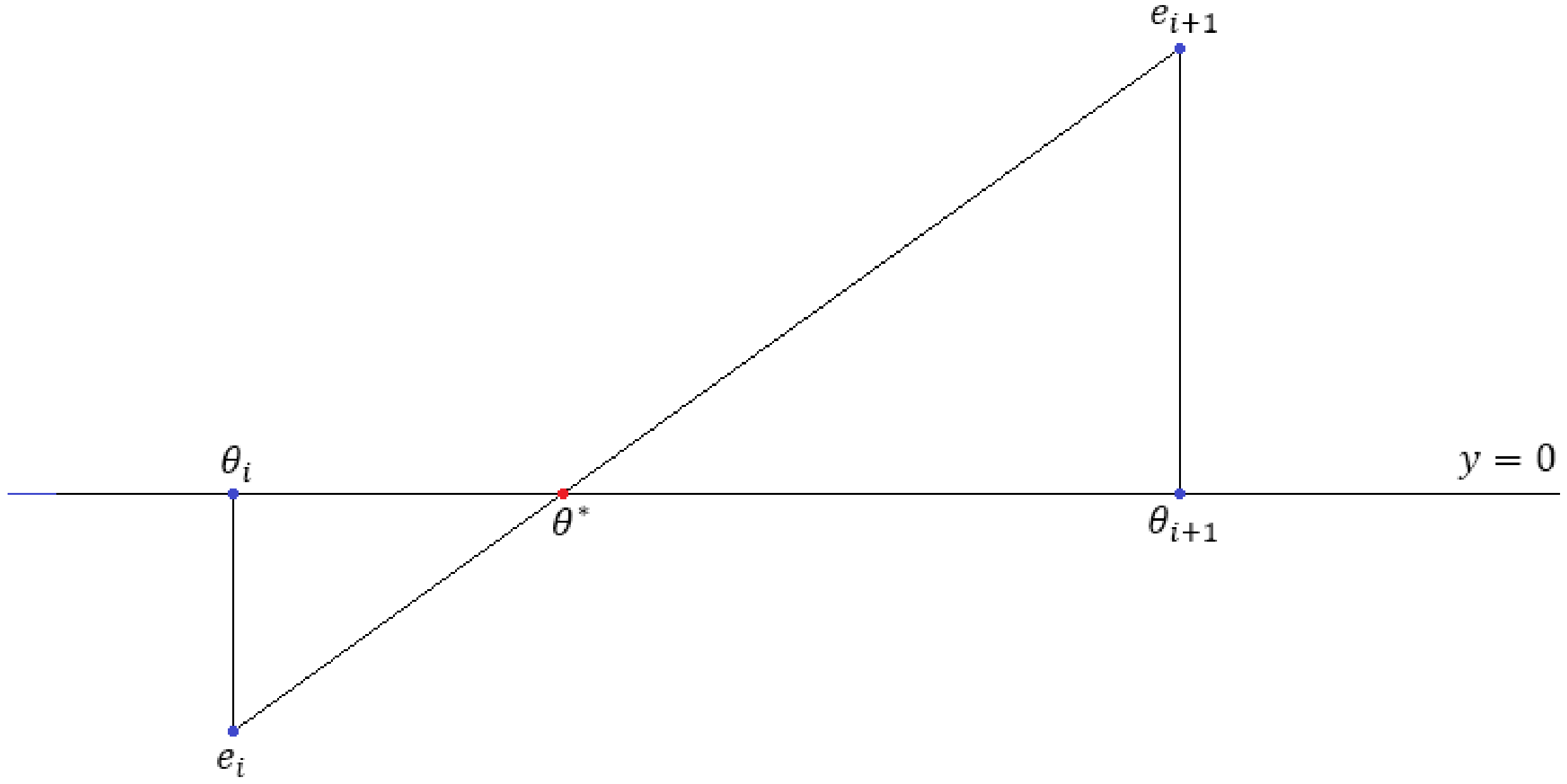

**Fig. A1.3**. A situation where two consecutive discrete elements of $e(\vec{\boldsymbol{\theta}})$, $e_i$ and $e_{i+1}$, have different signs, so the line segment whose endpoints are $(\theta_i, e_i)$ and $(\theta_{i+1}, e_{i+1})$, crosses horizontal axis $y = 0$ at point $(\theta^*, 0)$.

# Annex 2: Central R Functions to Implement the LC Method

## 1. Function to Build Structures *LzC* for Quadratic Polynomials

### Source Code

```
solveP2 <- function(c1,c2) {
   # Computes the roots of a degree 2 univariate polynomial
   # given its complex coefficients, by means of basic LC method
   # INPUTS :       c1,c2           : coefficients of polynomial
   #                                  z^2 + c1*z + c2
   # OUTPUTS:       r1,r2           : estimated polynomial roots
   #                V               : unit direction vector of
   #                                  line L1 containing r1, r2,
   #                                  expressed as complex number.
   #                                  L1's fixed point is -c1/2
   # NOTE:          The output is packed into a vector with
   #                3 components: r1, r2, V

   # fixed points of geometric construction
   P1 <- -c1/2
   P2 <-  c2

   # direction vector of line derived from L1
   n_v2 <- P2 - P1^2

   # angle of L1's direction vector
   arg_v <- atan2(Im(-n_v2),Re(-n_v2))/2

   # L1's direction vector
   V <- cos(arg_v) + sin(arg_v)*1i

   # two points from z-circumference zC
   zc1 <- P2/P1
   zc2 <- P2/(P1+V)

   # center of z-circumference
   czc <- center_zcircle(zc1,zc2)

   # the roots are the intersections between L1 and z-circumference zC
   I  <- intersect_semiLine_circle(P1, V, czc, abs(czc))
   r1 <- I[,1]
   r2 <- I[,2]

   return( c(r1,r2,V) )
}
```

**Remarks on this Function**

Function `solveP2` implements the construction of the unique structure $LzC(C_1, C_2, \theta^*)$ which produces the approximations to the roots $R_1, R_2$ of polynomial $z^2 + C_1 z + C_2$ by means of the two intersections between line $\ell_1(\theta^*)$ and z-circumference $C_2/\ell_1(\theta^*)$. The description of how to construct the structure $LzC(C_1, C_2, \theta^*)$ associated with polynomial $z^2 + C_1 z + C_2$ is found in chapter 2. `solveP2` takes as input coefficients $C_1, C_2$, and returns numerical approximations to the roots $R_1, R_2$, as well as to the unit direction vector $v_{\theta^*}$ of line $\ell_1(\theta^*)$. Function `solveP2` uses the basic functions `center_zcircle`, `intersect_semiLine_circle` from annex 1 for the construction of the structure $LzC(C_1, C_2, \theta^*)$.

## 2. Function to Build Structures $LzC$ for Cubic Polynomials

### Source Code

```
LzC3 <- function(CC, Theta) {

   # builds structures LzC for degree 3 univariate polynomials.
   # INPUTS:        CC    : vector with coefficients C1, C2, C3 of a
   #                        polynomial of the form
   #                        z^3 + C1*z^2 + C2*z + C3
   #                Theta : vector with inclination angles (with respect
   #                        to positive real axis) for lines L1 with
   #                        fixed point P1 = -C1/2
   # OUTPUTS:       E_A   : vector of weighted errors associated with
   #                        intersection I1 between tL and zC; an error
   #                        for each structure LzC
   #                E_B   : vector of weighted errors associated with
   #                        intersection I2 between tL and zC; an error
   #                        for each structure LzC
   #                i0I1  : vector of projections onto L1 of intersections
   #                        I1 between tL(Theta) and zC(Theta) via tL
   #                i0I2  : vector of projections onto L1 of intersections
   #                        I2 between tL(Theta) and zC(Theta) via tL
   #                c0I1  : vector of projections onto L1 of intersections
   #                        I1 between tL(Theta) and zC(Theta) via zC
   #                c0I2  : vector of projections onto L1 of intersections
   #                        I2 between tL(Theta) and zC(Theta) via zC
   #                ii1   : indicates the positions of structures LzC
   #                        (relative to array Theta) where I1 exists
   #                ii2   : indicates the positions of structures LzC
   #                        (relative to array Theta) where I2 exists
   #                iiU   : indicates the positions of structures LzC
   #                        (relative to array Theta) where only I2 exists
   # NOTE: Outputs are packed into a list OUT

   P1  = -CC[1]/2                     # fixed point of lines L1
   V   = cos(Theta) + sin(Theta)*1i   # direction vectors of lines L1
   N   = length(Theta)                # number of elements LzC
   AP  = rep(CC[2]-P1^2,N)            # fixed point of terminal semi-lines
   VTL = V^2               # direction vectors of terminal semi-lines

   # centers and radii of z-circumferences
   PzC1      = rep(-CC[3]/P1,N)    # points 1 for z-circumferences
   PzC2      = -CC[3]/(P1 + V)     # points 2 for z-circumferences
   Centre_zC = center_zcircle(PzC1,PzC2)
   Radius_zC = abs(Centre_zC)

   # intersections between terminal semi-lines and z-circumferences
   I_tL_zC = intersect_semiLine_circle(AP, VTL, Centre_zC, Radius_zC)
   I1      = I_tL_zC[,1]                # intersections I1
   I2      = I_tL_zC[,2]                # intersections I2
   ii1     = which(Re(I_tL_zC[,3])==2)  # indices where I1 exists
   ii2     = which(Re(I_tL_zC[,3])>0)   # indices where I2 exists
   iiU     = which(Re(I_tL_zC[,3])==1)  # indices where only I2 exists

   # projections of intersections between terminal semi-lines and
   # z-circumferences onto L1 via terminal semi-lines
   i0I1 = P1 + sqrt(Re((I1[ii1]-AP[ii1])/VTL[ii1]))*V[ii1]
   i0I2 = P1 + sqrt(Re((I2[ii2]-AP[ii2])/VTL[ii2]))*V[ii2]
```

```
    # projections of intersections between terminal semi-lines and
    # z-circumferences onto L1 via z-circumferences
    c0I1 = -CC[3]/I1[ii1]
    c0I2 = -CC[3]/I2[ii2]

    # initializes weighted error arrays
    E_A = rep(NA,N)
    E_B = rep(NA,N)

    # computes discrete mappings of normalized differences
    E_A[ii1] = error_tl_zc( rep(P1,length(ii1)), V[ii1], i0I1, c0I1 )
    E_B[ii2] = error_tl_zc( rep(P1,length(ii2)), V[ii2], i0I2, c0I2 )

    # prepares outputs
    OUT                  <- list()
    OUT[[length(OUT)+1]] <- E_A
    OUT[[length(OUT)+1]] <- E_B
    OUT[[length(OUT)+1]] <- i0I1
    OUT[[length(OUT)+1]] <- i0I2
    OUT[[length(OUT)+1]] <- c0I1
    OUT[[length(OUT)+1]] <- c0I2
    OUT[[length(OUT)+1]] <- ii1
    OUT[[length(OUT)+1]] <- ii2
    OUT[[length(OUT)+1]] <- iiU
    return(OUT)
}
```

## Remarks on this Function

Function `LzC3` implements the construction of structures $LzC(C_1, C_2, C_3, \theta)$ associated with univariate polynomials of degree 3 of the form $z^3 + C_1 z^2 + C_2 z + C_3$: for this, it uses the basic functions `center_zcircle`, `intersect_semiLine_circle`, and `error_tl_zc` from annex 1; the construction process of such structures $LzC$ is described in detail in chapter 3. Function `LzC3` takes advantage of the vectorization of operations in R to build (in parallel) multiple structures $LzC$, and can be used with an array (vector) `Theta` of $N$ equidistant angular values $\theta_k$, $k = 0,1,2, \dots, N-1$ (representing a regular partition of an interval on the real line $\mathbb{R}$), if a complete discrete proximity map $\{\hat{e}_A(\theta_k), \hat{e}_B(\theta_k)\}$, or a part of it, is to be constructed. Additionally, function `LzC3` can be used with a few input scalar values $\hat{\theta}_i^*$ ($i = 1,2,3$), packed into a vector `Theta` of length 3, in order to obtain approximations $\hat{R}_i$ to the roots $R_i$ of polynomial $z^3 + C_1 z^2 + C_2 z + C_3$; of course, these input values $\hat{\theta}_i^*$ should be good approximations to the true theta roots $\theta_i^*$, in order to obtain reasonable estimates $\hat{R}_i$ from the corresponding structures $LzC(C_1, C_2, C_3, \hat{\theta}_i^*)$.

For the case in which function `LzC3` is used to obtain a discrete proximity map $\{\hat{e}_A(\theta_k), \hat{e}_B(\theta_k)\}$ in the interval $\theta \in [-\pi, \pi)$ or in a region of this interval, output values `E_A`, `E_B` correspond to the weighted errors $e_A(\theta_k)$, $e_B(\theta_k)$ computed according to expressions (3.4a) and (3.4b) of chapter 3. On the other hand, for the case in which function `LzC3` is used in the computation of initial approximations $\hat{R}_i$ to the roots of polynomial $z^3 + C_1 z^2 + C_2 z + C_3$, output values `E_A`, `E_B`, `i0I1`, `i0I2`, `c0I1`, `c0I2`, `ii1`, `ii2`, `iiU` are useful in the implementation of expression (3.5) of chapter 3 for obtaining approximations $\hat{R}_i$. Function `LzC3` is used for both purposes within integrative function `approxP3` (listed in section 5 of this annex), which finally implements strategy 3.1 of chapter 3.

# 3. Function to Build Structures $LzC$ for Polynomials of Degree $n \geq 4$

## Source Code

```
LzC <- function(CC, Theta, MinT) {
   # Builds structures LzC for univariate polynomials of degree > = 4.
   # INPUTS:         CC    : vector with coefficients C1,C2,...,Cn of a
   #                         polynomial of the form
   #                         z^n+C1*z^(n-1)+C2*z^(n-2)+...+C(n-1)*z+Cn; n>=4
   #                 Theta : vector with inclination angles (with respect
   #                         to positive real axis) for lines L1 with
   #                         fixed point P1 = -C1/2
   #                 MinT  : vector with values t* that minimize the
   #                         Dynamic Squared Distance between terminal curve
   #                         tC(Theta) and z-circumference zC(Theta);
   #                         vectors Theta and MinT are of equal length
   # OUTPUTS:        E_A   : vector of weighted errors associated with
   #                         intersection I1 between tL and zC; an error
   #                         for each structure LzC
   #                 E_B   : vector of weighted errors associated with
   #                         intersection I2 between tL and zC; an error
   #                         for each structure LzC
   #                  Rx   : vector of best estimates of a root in terms
   #                         of MinT(Theta); an estimate for each
   #                         structure LzC
   #                  AP   : vector of anchor points for terminal semi-
   #                         lines tL, one anchor point per structure LzC
   #                 VTL   : array of direction vectors (expressed as
   #                         complex numbers) for terminal semi-lines tL;
   #                         one direction vector per structure LzC
   # NOTE: Outputs are packed into a list OUT

   n   = length(CC)                  # polynomial degree
   P1  = -CC[1]/2                    # fixed point of lines L1
   V   = cos(Theta) + sin(Theta)*1i  # direction vectors of lines L1
   Rx  = P1 + MinT*V           # best root estimates in terms of MinT(Theta)
   N   = length(Theta)               # number of elements LzC to process
   i   = 1:N                         # indices to access structures LzC

   # builds anchor points for terminal semi-lines tL;
   # see expressions (5.13) in chapter 5
   AP = apply(as.data.frame(i),1,function(i) {
       ((-1)^(n+1))*(CC[(n-1):3]%*%Rx[i]^(0:(n-4))+(CC[2]-P1^2)*Rx[i]^(n-3))
   })

   # builds direction vectors for terminal semi-lines tL;
   # see expressions (5.13) in chapter 5
   VTL = ((-1)^(n+1))*(Rx^(n-3))*V^2

   # centers and radii of z-circumferences zC
   PzC1      = rep(((-1)^n)*CC[n]/P1,N)   # points 1 for z-circumferences
   PzC2      = ((-1)^n)*CC[n]/(P1 + V)    # points 2 for z-circumferences
   Centre_zC = center_zcircle(PzC1,PzC2)
   Radius_zC = abs(Centre_zC)

   # intersections between terminal semi-lines and z-circumferences
   I_tL_zC = intersect_semiLine_circle(AP, VTL, Centre_zC, Radius_zC)
   I1      = I_tL_zC[,1]                # intersection I1
   I2      = I_tL_zC[,2]                # intersection I2
   ii1     = which(Re(I_tL_zC[,3])==2)  # Indices where I1 exists
   ii2     = which(Re(I_tL_zC[,3])>0)   # Indices where I2 exists
```

```
    # projections of intersections between terminal semi-lines and
    # z-circumferences onto L1 via terminal semi-lines
    i0I1 = P1 + sqrt(Re((I1[ii1]-AP[ii1])/VTL[ii1]))*V[ii1]
    i0I2 = P1 + sqrt(Re((I2[ii2]-AP[ii2])/VTL[ii2]))*V[ii2]

    # projections of intersections between terminal semi-lines and
    # z-circumferences onto L1 via z-circumferences
    c0I1 = ((-1)^n)*CC[n]/I1[ii1]
    c0I2 = ((-1)^n)*CC[n]/I2[ii2]

    #initializes the weighted error arrays
    E_A = rep(NA,N)
    E_B = rep(NA,N)

    # computes discrete mappings of normalized differences
    E_A[ii1] = error_tl_zc( rep(P1,length(ii1)), V[ii1], i0I1, c0I1 )
    E_B[ii2] = error_tl_zc( rep(P1,length(ii2)), V[ii2], i0I2, c0I2 )

    # prepares outputs
    OUT                  <- list()
    OUT[[length(OUT)+1]] <- E_A
    OUT[[length(OUT)+1]] <- E_B
    OUT[[length(OUT)+1]] <- Rx
    OUT[[length(OUT)+1]] <- AP
    OUT[[length(OUT)+1]] <- VTL
    return(OUT)
}
```

**Remarks on this Function**

Function `LzC` implements the construction of structures $LzC(C_1, C_2, \ldots, C_n, \theta)$ associated with univariate polynomials of degree $n \geq 4$ of the form $z^n + C_1 z^{n-1} + \cdots + C_{n-1} z + C_n$: for this, it uses the basic functions `center_zcircle`, `intersect_semiLine_circle`, and `error_tl_zc` from annex 1; the construction process of such structures is described in detail in chapters 4 and 5. Function `LzC` takes advantage of the vectorization of operations in R to build (in parallel) multiple structures $LzC$, and can be used with an array (vector) `Theta` of $N$ equidistant angular values $\theta_k$, $k = 0,1,2, \ldots, N-1$ (representing a regular partition of an interval on the real line $\mathbb{R}$), if a complete discrete proximity map $\{\hat{e}_A(\theta_k), \hat{e}_B(\theta_k)\}$, or a part of it, is to be constructed. Additionally, function `LzC` can be used with arbitrary input scalar values $\theta_s$, packed into a vector `Theta`, in order to analyze specific elements, be they anchor points $aP(\theta_s)$ and direction vectors $v_{tL}(\theta_s)$ of terminal semi-lines $tL(\theta_s)$, or best root estimates $R_x(\theta_s)$.

Note that, with respect to function `LzC3` from section 2 of this annex, function `LzC` requires an additional input, vital for its operation: the vector of values `MinT`, which must coincide in length and correspond, element by element, with input vector `Theta`; vector `MinT` is the output of function `min_D2` described in section 4 of this annex, and it contains approximations to the values $t^*(\theta)$ that minimize the dynamic squared distance function $d_\theta^2(t)$ between terminal curve $t\mathfrak{C}(\theta)$ and z-circumference $zC(\theta)$ (see chapters 4 and 5 for a detailed description of these concepts).

For the case in which function `LzC` is used to obtain a discrete proximity map $\{\hat{e}_A(\theta_k), \hat{e}_B(\theta_k)\}$ in the interval $\theta \in [-\pi, \pi)$ or in a region of this interval, output values `E_A`, `E_B` correspond to the weighted errors $e_A(\theta_k)$, $e_B(\theta_k)$ computed according to expressions (5.16a) and (5.16b) of chapter 5. On the other hand, for the case in which function `LzC` is used with input values $\theta_s$ to obtain best estimates $R_x(\theta_s)$, or anchor points $aP(\theta_s)$ and direction vectors $v_{tL}(\theta_s)$ of terminal semi-lines $tL(\theta_s)$, output values `Rx`, `aP`, `VTL` correspond to expressions (5.9) and (5.13) in chapter 5. Function `LzC` is used within the integrative function `approxLC` (listed in section 6 of this annex) to generate discrete proximity maps $\{\hat{e}_A(\theta_k), \hat{e}_B(\theta_k)\}$ and to implement the strategy 5.1 from chapter 5. Function `LzC` is also used within the script listed in annex 3 section 3, in order to generate plots of trajectories $R_x(\theta)$ and $aP(\theta)$, which are shown in the numerical examples of chapters 4 and 5; these plots are obtained directly from the output vectors `Rx`, `aP` of function `LzC`.

## 4. Routine to Estimate Global Minimum of Dynamic Squared Distance

### Source Code

```
min_D2 <- function(V,C, MaxIt, Temp, TMax) {
   # approximates global minimum of Dynamic Squared Distance (DSD) function
   # between Terminal Curve tC and z-circumference zC
   # INPUTS:         V : array with unit direction vectors for lines L1
   #                 C : array with coefficients of polynomial of the form
   #                     p(z) = z^n + C[1:n] %*% z^((n-1):0); n>3
   # MaxIt, Temp, TMax : parameters for the Simulated Annealing phase of the
   #                     optimization process
   # OUTPUTS:     MinF : stores global minimum value of DSD function
   #                     for each element of V
   #              MinT : stores parametric value t* where the global minimum
   #                     of DSD function for each element of V occurs
   # The output is packed into a matrix, for which its first column is MinF,
   # and its second column is MinT
   N  = length(V)  # number of structures LzC to be processed
   i  = 1:N        # indices to process DSD functions via R function apply
   n  = length(C)  # degree of polynomial
   P1 = -C[1]/2    # fixed point of lines L1
   D2 = vector("list",N) # stores DSD functions, one per structure LzC
   D2 = apply(as.data.frame(i),1,function(i) {
           v = V[i]
           function(t) {
              Z = P1+t*v
              tC = ((-1)^(n+1))*(C[(n-1):2]%*%Z^(0:(n-3))-(P1-t*v)*Z^(n-2))
              zC = ((-1)^n)*C[n]/Z
              return( Re(tC-zC)^2 + Im(tC-zC)^2 )
           }
        })
   # finds the global minimum for each DSD function:
   # phase 1: optimization using simulated annealing
   t_min0 = apply(as.data.frame(i),1,function(i) {
               optim(1,D2[[i]],method="SANN",
                     control=list(maxit=MaxIt,temp=Temp,tmax=TMax))
            })
   # phase 2: refinement of phase 1’s result via BFGS
   t_min = apply(as.data.frame(i),1,function(i) {
               optim(t_min0[[i]]$par,D2[[i]],method="BFGS")
           })
   # retrieves minimum functional value of DSD function
   MinF = rep(NA,N)
   MinF[i] = unlist(apply(as.data.frame(i),1,function(i) {
                if(t_min[[i]]$par>=0) {
                   t_min[[i]]$value
                } else {
                   NA
                }
            }))
   # retrieves parametric value t* where DSD’s global minimum occurs
   MinT = rep(NA,N)
   MinT[i] = unlist(apply(as.data.frame(i),1,function(i) {
                if(t_min[[i]]$par>=0) {
                   t_min[[i]]$par
                } else {
                   NA
                }
             }))
   return(cbind(MinF, MinT))
}
```

**Remarks on this Function**

Function `min_D2` seeks to find value $t^*$ in parametric space $t$ where the global minimum of the dynamic squared distance function $d_\theta^2(t)$ defined by expression (5.8) of chapter 5 is located, keeping fixed value $\theta$; it also obtains the corresponding functional value of the global minimum, $d_\theta^2(t^*)$. It is worth mentioning that all values $t^*(\theta)$ obtained by function `min_D2` are non-negative. `min_D2` is indispensable in the construction of discrete proximity maps associated with univariate polynomials of degree 4 or higher, which is why it is used within functions responsible for constructing such maps (see, for example, function `approxLC` in section 6 of this annex).

`min_D2` resorts to two optimization phases; in the first phase, an initial approximation to the minimizing argument $t^*$ is obtained by means of the *simulated annealing* algorithm, which in theory is able to approximate $t^*$ with arbitrary precision, if given enough time; in the second phase, `min_D2` refines the initial approximation of $t^*$ computed by the first phase, for which it uses the Broyden-Fletcher-Goldfarb-Shanno (BFGS) algorithm. Conceptually, `min_D2` processes in parallel multiple values $\theta_k$ associated with structures $LzC(C_1, C_2, \dots, C_n, \theta_k)$, which is useful in the construction of proximity maps based on the first-order discrete derivative of their outputs, `MinF` and `MinT`, which are arrays with sample values from respective functions (of $\theta$) $d_\theta^2[t^*(\theta)] = \min_t d_\theta^2(t)$ and $t^*(\theta) = \arg\min_t d_\theta^2(t)$.

In the first part of `min_D2`, a dynamic squared distance function $d_\theta^2(t)$ associated with input coefficients `C` is constructed for each discrete value of $\theta$ (encoded within input array `V`, which contains unit direction vectors for lines $\ell_1$) by using R's function `apply`. Functions $d_\theta^2(t)$ built in this way are stored into object `D2`. As we saw in chapters 4 and 5, $d_\theta^2(t)$ represents the dynamic squared distance (the adjective "dynamic" refers to the variation of parameter $t$) between terminal curve $t\mathfrak{C}(\theta, t)$ and z-circumference $zC(\theta, t)$; in addition, we also saw that $d_\theta^2(t)$ is the reflection, with respect to vertical axis $t = 0$, of $d_{\theta+\pi}^2(t)$. This symmetry property for function $d_\theta^2(t)$ in principle can be used to halve the number of dynamic squared distance functions constructed within `min_D2` when obtaining global maps in interval $\theta \in [-\pi, \pi)$; this is not done, however, in order to give `min_D2` greater flexibility: `min_D2` can be used both to construct global or local regions of proximity maps, and to obtain quality measures associated with root estimates $\hat{R}_i$ for polynomial $z^n + C_1 z^{n-1} + \cdots + C_{n-1} z + C_n$, with $n \geq 4$. Details on how to obtain the algebraic expressions for terminal curve $t\mathfrak{C}(\theta, t)$ and z-circumference $zC(\theta, t)$ that are used within `min_D2` can be consulted in chapter 5.

In the second part of `min_D2`, we try to obtain, by means of function `optim` (optimization of functions in R), the global minimum of each function $d_\theta^2(t)$ built in the first part, distributing the "parallel" processing work by means of function `apply`. As we mentioned above, the optimization strategy uses two methods (both available within `optim`): simulated annealing (method `SANN` in `optim`) to obtain an initial approximation of $t^*$, and `BFGS` to refine this initial approximation. Method `SANN` uses the following parameters, which according to [A2.1] and [A2.2] are: `par=1` as the starting point for the simulated annealing algorithm in the search space (in our case variable $t$), `MaxIt` as the maximum number of iterations within the algorithm, `TMax` as the number of function evaluations at each value of the decreasing temperature whose initial value is given by parameter `Temp`, and of course the objective function `fn`, which in our case is $d_\theta^2(t)$, encoded within `min_D2` as `D2[[i]]`. Method `BFGS`, on the other hand, is used with only two input parameters, namely, the output of method `SANN`, `t_min0[[i]]$par`, which serves as starting value `par` for `BFGS`, and the objective function `fn`, which of course is `D2[[i]]`; the rest of the input parameters available to `BFGS` are not specified, so R assigns to them default values; for example, since no value is specified for input parameter `gr`, which provides information about the functional form of the gradient (first derivative of the objective function), by default `gr=NULL`, which means that finite difference approximations are used for the gradient within `BFGS`, according to [A2.2].

The simulated annealing (`SANN`) method is a probabilistic method (also called *metaheuristic* in the field of computer science) used in unconstrained optimization problems; `SANN` seeks to find the global minimum of an objective function (in our case, the dynamic squared distance function $d_\theta^2(t)$) with an arbitrarily large search space. `SANN` is inspired by metallurgical processes, where first a certain metal alloy is subjected to high temperatures, and then its cooling is controlled so that the final product has certain desirable physical properties. During the cooling process, `SANN` moves, with a probability proportional to the current temperature in the process, to solutions (evaluations of the objective function) that may be inferior to previously found solutions; this allows to explore the search space in an exhaustive way and helps prevent the method from getting stuck in local minima. The disadvantage of `SANN` is that it converges slowly to the optimal solution, so it is important to properly select its input parameters so that in a relatively small number of iterations it reaches an approximation to the global minimum perhaps not very accurate, but reasonably close to the true optimal value. For more details on the inner workings of the `SANN` method within R's function `optim`, see [A2.1]. The `BFGS` method, on the other hand, is a quasi-Newton iterative algorithm; this means that it is a method that uses Newton's step (or iteration) to approximate to a minimum value of a scalar field in several variables, but that uses a more economical estimate of the inverse of the Hessian matrix rather than computing it directly at each iteration; the information for computing such an estimate is obtained from gradients of the objective function computed in previous iterations, so `BFGS` can be seen as a generalization of the secant method. `BFGS` obtains a minimum of the objective function with greater precision and speed compared to probabilistic method `SANN`, although in our experience, `BFGS` will almost certainly get stuck in a local minimum if it is not provided with a starting point reasonably close to the minimizing argument; in our strategy, `SANN` provides such a starting point reasonably close to the minimizing argument, while `BFGS` obtains a more refined solution. For a better technical understanding of the BFGS method, refer to the introductory lesson [A2.3].

Returning to our function `min_D2`, what we do in its third and last part is simply recover the optimal values $d_\theta^2[t^*(\theta)]$ and $t^*(\theta)$ obtained by the BFGS method, and store them in the respective vectors `MinF` and `MinT`; for this, we make sure to stick only with values $t^*(\theta)$ greater than or equal to zero; otherwise we store a value `NA` (not available) in the corresponding elements `MinF[i]` and `MinT[i]`. Why do we only consider non-negative values for the optimal argument $t^*$? Because we always assume that an approximation $\hat{R}_i$ to a root of the polynomial $z^n + C_1 z^{n-1} + \cdots + C_{n-1} z + C_n$ must be of the form $\hat{R}_i = -C_1/2 + t^*(\hat{\theta}_i^*) v_{\hat{\theta}_i^*}$, where $v_{\hat{\theta}_i^*}$ is the unit direction vector of line $\ell_1$ with direction $\hat{\theta}_i^*$ such that $v_{\hat{\theta}_i^*}$ starts at fixed point $-C_1/2$ and points towards $\hat{R}_i$.

# 5. Function to Build LC Proximity Maps for Cubic Polynomials

## Source Code

```
approxP3 <- function( c1, c2, c3, N, Th_i, Th_f ) {
   # Obtains approximations to roots of a univariate cubic polynomial
   # given its complex coefficients, by means of the LC method
   # INPUTS :         c1,c2,c3   : coefficients of monic polynomial
   #                               z^3 + c1*z^2 + c2*z + c3
   #                 N           : number of elements LzC used in the
   #                               construction of the proximity map
   #                 Th_i, Th_f : initial and final values in the
   #                               Theta interval of proximity map
   # OUTPUTS:        E_A    : discrete values for weighted errors associated
   #                         with intersections I1 between tL and zC
   #                 E_B    : discrete values for weighted errors associated
   #                         with intersections I2 between tL and zC
   #                 Theta : inclination angles of lines L1 (with respect to
   #                         positive real axis) considered in the discrete
   #                         proximity map; these angles are associated with
   #                         the corresponding components of E_A, E_B;
   #                         together they form the proximity map
   #            APPROX_LC : table of approximations
   #                          Rx    - approximations to roots of polynomial
   #                          Theta - approximations to "theta roots"
   #                          Vert  - absolute vertical distances
   #                                  associated with discrete
   #                                  approximations to theta roots
   # NOTE: outputs are packed into a list OUT
   i  = 1:N        # indices for the discrete proximity map

   # constructs a regular partition of Theta interval
   Theta = Th_i+(Th_f-Th_i)*(i-1)/N    # interval [Th_i, Th_f)

   CC <- c(c1,c2,c3)  # places polynomial coefficients into a vector
   OUTPUT1 <- LzC3(CC,Theta) # generates the proximity map
   E_A     <- OUTPUT1[[1]]
   E_B     <- OUTPUT1[[2]]

   # obtains numerical approximations of "theta roots"
   # via crossings of functions E_A, E_B with horizontal axis y=0
   if (Th_i==-pi & Th_f==pi) {
      # takes advantage of the periodicity of weighted error function
      O1 = get_angle_roots( c(Theta,pi), c(E_A,E_A[1]), 1.0 )
      O2 = get_angle_roots( c(Theta,pi), c(E_B,E_B[1]), 1.0 )
   } else {
      O1 = get_angle_roots( Theta, E_A, 1.0 )
      O2 = get_angle_roots( Theta, E_B, 1.0 )
   }

   # stores approximations to theta roots, as well as the
   # corresponding vertical differences between the
   # discrete elements involved
   theta_roots  = c( O1[[1]], O2[[1]] )
   vertical_dif = c( O1[[2]], O2[[2]] )

   # ADDITIONAL ROUTINE TO RECOVER A POSSIBLE ROOT APPROXIMATION
   # finds minimum vertical distance between E_A and E_B, restricted to
   # the case in which E_A and E_B have opposite signs; this is an attempt
   # to recover a theta root that occurs near a tangential intersection
   # between tL and zC
```

```
   if (length(theta_roots)<3 & (Th_i==-pi & Th_f==pi)) {
      DMIN=min(abs(E_A-E_B)[which(E_A*E_B<0)])
      theta_roots  = c( theta_roots, Theta[which(abs(E_A-E_B)==DMIN)] )
      vertical_dif = c( vertical_dif, DMIN )
   }

   #############################
   # obtains root approximations Rx,
   # corresponding to approximations to theta roots

   # step 1: obtains structures LzC
   # for approximations to theta roots
   OUTPUT2  <- LzC3(CC,theta_roots)
   E_Ar     <- OUTPUT2[[1]]
   E_Br     <- OUTPUT2[[2]]
   i0I1     <- OUTPUT2[[3]]
   i0I2     <- OUTPUT2[[4]]
   c0I1     <- OUTPUT2[[5]]
   c0I2     <- OUTPUT2[[6]]
   ii1      <- OUTPUT2[[7]]
   iiU      <- OUTPUT2[[9]]

   # step 2: implements expression (3.5) from chapter 3 to build Rx
   Rx = rep(NA,length(theta_roots))
   Rx[ii1] = ifelse(abs(E_Ar[ii1])<abs(E_Br[ii1]),
                   (i0I1[ii1]+c0I1[ii1])/2,
                   (i0I2[ii1]+c0I2[ii1])/2)
   Rx[iiU] = (i0I2[iiU]+c0I2[iiU])/2

   # listing of initial approximations of roots, theta roots,
   # and vertical distances between corresponding consecutive functional
   # values of discrete LC map with opposite signs
   APPROX_LC = setNames(data.frame( Rx[order(vertical_dif)],
                                  theta_roots[order(vertical_dif)],
                                  vertical_dif[order(vertical_dif)] ),
                       c("Rx","Theta","Vert"))
   # prepares output
   OUT                  <- list()
   OUT[[length(OUT)+1]] <- E_A
   OUT[[length(OUT)+1]] <- E_B
   OUT[[length(OUT)+1]] <- Theta
   OUT[[length(OUT)+1]] <- APPROX_LC
   return(OUT)
}
```

## Remarks on this Function

The integrative function `approxP3` implements strategy 3.1 of chapter 3, considering of course the weighted error definitions (3.4a) and (3.4b), as well as definition (3.5) that allows us to compute initial approximations to roots. Function `approxP3` also incorporates an additional routine (see example 3.3 of chapter 3) that tries to recover an extra approximation to a root, in case the number of approximations obtained by means of the global LC map is less than 3. `approxP3` uses function `LzC3` from section 2 in this annex, and function `get_angle_roots` from annex 1 section 4; in turn, `approxP3` is used within the script listed in annex 3 section 2, which is responsible for displaying in the R console the final outputs obtained by this function `approxP3`.

# 6. Function to Build LC Proximity Maps for Polynomials of Degree $n \geq 4$

## Source Code

```
approxLC <- function( CC, N, Th_i, Th_f, MAXIT, TEMP, TMAX ) {

   # Obtains initial approximations to roots of a polynomial of degree >= 4
   # given its complex coefficients, by means of the LC method
   #
   # INPUTS:   CC         : vector with coefficients C1,C2,...,Cn
   #                        of monic polynomial
   #                        z^n+C1*z^(n-1)+C2*z^(n-2)+...+C(n-1)*z+Cn; n>=4
   #          N           : number of elements LzC used in the
   #                        construction of the proximity map
   #          Th_i, Th_f  : initial and final values for the
   #                        Theta interval of proximity map
   #    MAXIT, TEMP, TMAX : parameters for the simulated annealing phase
   #                        of the optimization process that finds the
   #                        global minimum of the dynamic squared distance
   #                        (DSD) associated to corresponding element LzC
   # OUTPUTS:         MinF : vector of global minima of DSDs associated
   #                        with elements LzC of proximity map
   #                  MinT : vector of arguments t* where the global minima
   #                        of the DSDs associated to elements LzC occur
   #                  E_A  : discrete values for weighted errors associated
   #                        with intersections I1 between tL and zC
   #                  E_B  : discrete values for weighted errors associated
   #                        with intersections I2 between tL and zC
   #                 Theta : inclination angles of lines L1 (with respect to
   #                        positive real axis) considered in the discrete
   #                        proximity map; these angles are associated with
   #                        the corresponding components of E_A, E_B;
   #                        together they form the proximity map
   #            APPROX_LC : table of approximations
   #                         Rx    - approximations to polynomial roots
   #                         Theta - approximations to "theta roots"
   #                         Vert  - absolute vertical distances
   #                                 associated with discrete
   #                                 approximations to theta roots
   #                         MinF  - global minima of DSDs associated
   #                                 with approximations to theta roots
   # NOTE: outputs are packed into a list OUT

   i  = 1:N           # Indices for the discrete proximity map

   # constructs a regular partition of Theta interval
   Theta = Th_i+(Th_f-Th_i)*(i-1)/N    # discrete interval [Th_i, Th_f)

   # gets global minimum of DSD between terminal curve tC(Theta)
   # and z-circumference zC(Theta)
   V = cos(Theta) + sin(Theta)*1i  # unit direction vectors of lines L1
   OUT1 = min_D2(V, CC, MaxIt=MAXIT, Temp=TEMP, TMax=TMAX)
   MinF = OUT1[,1]
   MinT = OUT1[,2]


   OUT2 <- LzC(CC, Theta, MinT)    # generates LC proximity map
   E_A  <- OUT2[[1]]
   E_B  <- OUT2[[2]]
```

```
   # obtains numerical approximations of "theta roots"
   # via crossings of discrete mappings E_A, E_B with horizontal axis y=0
   if (Th_i==-pi & Th_f==pi) {
      # takes advantage of the periodicity of weighted error function
      O1 = get_angle_roots( c(Theta,pi), c(E_A,E_A[1]), 1.0 )
      O2 = get_angle_roots( c(Theta,pi), c(E_B,E_B[1]), 1.0 )
   } else {
      O1 = get_angle_roots( Theta, E_A, 1.0 )
      O2 = get_angle_roots( Theta, E_B, 1.0 )
   }

   # stores approximations to theta roots, as well as the
   # corresponding vertical differences between the
   # discrete elements involved
   theta_roots  = c( O1[[1]], O2[[1]] )
   vertical_dif = c( O1[[2]], O2[[2]] )

   # obtains global minima of DSDs between terminal curves
   # tC(theta_roots) and corresponding z-circumferences zC(theta_roots),
   # as an additional quality measure for root estimates
   r = 1:(length(theta_roots))
   Vr = cos(theta_roots[r]) + sin(theta_roots[r])*1i
   OUT3 = min_D2(Vr, CC, MaxIt=MAXIT, Temp=TEMP, TMax=TMAX)
   MinFr = OUT3[,1]
   MinTr = OUT3[,2]

   # listing of initial estimates of roots, theta roots,
   # vertical distances between corresponding consecutive functional
   # values of discrete LC map with opposite signs, and global minima
   # of DSDs associated with estimated theta roots
   APPROX_LC = setNames(data.frame( -CC[1]/2 + (MinTr*Vr)[order(MinFr)],
                                    theta_roots[order(MinFr)],
                                    vertical_dif[order(MinFr)],
                                    MinFr[order(MinFr)] ),
                        c("Rx","Theta","Vert","MinF"))

   # prepares output
   OUT                  <- list()
   OUT[[length(OUT)+1]] <- MinF
   OUT[[length(OUT)+1]] <- MinT
   OUT[[length(OUT)+1]] <- E_A
   OUT[[length(OUT)+1]] <- E_B
   OUT[[length(OUT)+1]] <- Theta
   OUT[[length(OUT)+1]] <- APPROX_LC
   return(OUT)
}
```

**Remarks on this Function**

The integrative function `approxLC` can be seen as a generalization of function `approxP3` listed in section 5 of this annex, since it extends the approximation of roots of cubic polynomials, via the LC method, to the case of univariate polynomials of degree 4 or higher. `approxLC` implements strategy 5.1 of chapter 5 (which is a generalization of strategy 4.1 of chapter 4); for this, `approxLC` internally uses 1) function `min_D2` from section 4 in this annex, 2) function `LzC` from section 3 in this annex, and 3) function `get_angle_roots` from annex 1 section 4. `approxLC` generates the results related to LC maps that are shown in the numerical examples of chapters 4 and 5; for this, `approxLC` is used within the script listed in annex 3 section 3, which is responsible for visualizing in the R console all the final results associated with the three types of proximity maps developed in this work.

Note that the calls to function `get_angle_roots` within `approxLC` use the argument `TOL=1.0`; this value works well with the numerical examples in chapters 4 and 5 as an upper bound of what we could consider as a discrete smooth crossing of function $e(\theta)$ with horizontal axis $y = 0$; however, we should not use this particular value for any arbitrary case; this particular value is used only for the illustrative purposes in this work, and of course it would be necessary to add source code in order to compute an argument `TOL` suitable to the norms of the input coefficients `CC`, considering that weighted errors $e(\theta)$ are limited to values from the interval $(-1,1]$.

# 7. Function to Build $dD$ or $dt$ Proximity Maps

## Source Code

```
approxDMin <- function(CC, Min_, Th_i, Th_f, MAXIT, TEMP, TMAX, Tol) {
   # computes a proximity map, as well as the corresponding initial
   # estimates, for the roots of a polynomial of degree greater than
   # or equal to 4, by means of the discrete derivative of function
   # Min_(Theta) with respect to Theta
   #
   # INPUTS :     CC                 : vector with coefficients C1,C2,...,Cn
   #                                   of a polynomial of the form
   #                                   z^n+C1*z^(n-1)+...+C(n-1)*z+Cn; n>=4
   #              Min_               : vector of global minima d^2(t*)
   #                                   (or of minimizing arguments t*)
   #                                   associated with each element LzC
   #              Th_i, Th_f         : initial and final values for Theta
   #                                   interval of proximity map
   #              MAXIT, TEMP, TMAX : parameters for SANN phase of
   #                                   optimization process min_D2 that
   #                                   finds global minima associated
   #                                   with theta roots from proximity map
   #              Tol                : upper threshold to detect smooth
   #                                   crossings of (d/dTheta)Min_(Theta)
   #                                   with horizontal axis y=0 in terms of
   #                                   vertical differences between
   #                                   consecutive functional values of
   #                                   opposite signs
   # OUTPUTS :    DD_Min_            : vector with discrete values of
   #                                   derivative (d/dTheta)Min_(Theta)
   #              Theta_shift        : support values for DD_Min_
   #              APPROX_DDMIN       : table of approximations
   #                    Rx           - approximations to polynomial roots
   #                    Theta        - approximations to theta roots
   #                    Vert         - absolute vertical distances
   #                                   associated with discrete
   #                                   approximations to theta roots
   #                    MinF         - global minima d^2(t*) associated
   #                                   with approximations to theta roots
   # NOTE: outputs are packed into a list OUT

   N = length(Min_)       # number of discrete elements to be processed
   i = 1:N                # index for accessing elements
   j = 2:N                # index for accessing elements of derivative

   Theta = Th_i + (Th_f-Th_i)*(i-1)/N  # regular partition of [Th_i, Th_f)

   delta_theta = Theta[2]-Theta[1]     # differential of Theta

   # discrete derivative of Min_ with respect to Theta
   # (d/dTheta)Min_(Theta), stored into array DD_Min_
   if (Th_i==-pi & Th_f==pi) {
      # takes advantage of the periodicity of global function Min_
      DD_Min_ = c( (Min_[1]-Min_[N])/delta_theta,
                   (Min_[j]-Min_[j-1])/delta_theta )
   } else {
      DD_Min_ = (Min_[j]-Min_[j-1])/delta_theta
   }
```

```
   # computes numerical approximations to theta roots (Th_Z) by means of
   # crossings of DD_Min_ with horizontal axis y=0
   if (Th_i==-pi & Th_f==pi) {
      Th_Z = get_angle_roots(Theta-delta_theta/2, DD_Min_, Tol)
   } else {
      Th_Z = get_angle_roots(Theta[-1]-delta_theta/2, DD_Min_, Tol)
   }

   # computes global minima of dynamic squared distances between
   # terminal curves tC(Th_Z) and corresponding z-circumferences
   # zC(Th_Z) as a quality measure for estimates Th_Z
   k     = 1:(length(Th_Z[[1]]))
   Vk    = cos(Th_Z[[1]][k]) + (sin(Th_Z[[1]][k]))*1i
   OUT   = min_D2(Vk, CC, MaxIt=MAXIT, Temp=TEMP, TMax=TMAX)
   MinFk = OUT[,1]
   MinTk = OUT[,2]

   # listing of initial estimates of roots, theta roots,
   # vertical distances between corresponding consecutive functional
   # values of discrete map DD_Min_ with opposite signs, and global
   # minima of DSDs associated with estimated theta roots Th_Z
   APPROX_DDMIN = setNames(data.frame( -CC[1]/2 + (MinTk*Vk)[order(MinFk)],
                                       Th_Z[[1]][order(MinFk)],
                                       Th_Z[[2]][order(MinFk)],
                                       MinFk[order(MinFk)] ),
                           c("Rx","Theta","Vert","MinF"))
   # prepares output
   OUT                  <- list()
   OUT[[length(OUT)+1]] <- DD_Min_
   if (Th_i==-pi & Th_f==pi) {
      OUT[[length(OUT)+1]] <- Theta-delta_theta/2
   } else {
      OUT[[length(OUT)+1]] <- Theta[-1]-delta_theta/2
   }
   OUT[[length(OUT)+1]] <- APPROX_DDMIN
   return(OUT)
}
```

## Remarks on this Function

Function `approxDMin` is the result of experimental observations arising from the numerical examples in chapter 4, where we realize that it is possible to obtain proximity maps by means of the derivative, with respect to $\theta$, both of function $t^*(\theta)$, and of function $d_\theta^2[t^*(\theta)]$. `approxDMin` follows the same philosophy of function `approxLC` from section 6 in this annex; its main role is to derive any of the discrete versions of input functions $t^*(\theta)$ or $d_\theta^2[t^*(\theta)]$ with respect to $\theta$, and then to approximate crossings of these discrete derivatives with horizontal axis $y = 0$ by means of function `get_angle_roots` listed in annex 1 section 4, thus obtaining initial approximations to theta roots $\theta_i^*$. Both `approxDMin` and `approxLC` are called within the script listed in annex 3 section 3 in order to generate proximity maps associated with univariate polynomials of degree $n \geq 4$. `approxDMin` helps in the implementation of strategies 5.2 and 5.3 from chapter 5.

`approxDMin` computes quotients of first-order differences of the form

```
(Min_[j]-Min_[j-1])/(Theta[j]-Theta[j-1]),
```

which can be regarded as numerical approximations to sample values of the derivative of the input function (whose sample values are the elements in vector `Min_`) with respect to independent variable `Theta`; **each of these quotients approximate the derivative of input function `Min_` evaluated at point `(Theta[j-1]+Theta[j])/2`.**

Care must be taken when using this experimental function `approxDMin`, as the first-order discrete derivative tends to amplify the noise contained in the functional values on which it operates, and this amplifying effect becomes more noticeable if we increase the sampling rate, such as when we want to produce discrete high-resolution proximity maps.

Note that the calls to function `get_angle_roots` within `approxDMin` use the same argument `Tol` that `approxDMin` itself receives as input; this is done in this way because the upper bound of what we could consider as a discrete smooth crossing of function $\frac{d}{d\theta}t^*(\theta)$ or of function $\frac{d}{d\theta}d_\theta^2[t^*(\theta)]$ with horizontal axis $y = 0$ depends on the norms of the polynomial coefficients, and consequently, on the spatial location of the roots in complex plane $\mathbb{C}$. See details in the numerical examples of chapters 4 and 5 to get an idea of what the appropriate value for the argument `Tol` might be. The way in which function `approxDMin` is designed obeys the illustrative purposes in this work, and of course it would be necessary to add source code within this function in order to automate the process for finding a proper value for `Tol`, based on input coefficients `CC`.

# Annex 3: R Scripts for Visualizing Proximity Maps

## 1. R Script to Display Results of the Examples in Chapter 2

### Source Code

```
rm(list=ls(all=TRUE))
# set.seed(2019) # remove comment to reproduce example 2.2

##########################################################################
# insert here code for center_zcircle            (annex 1 section 1)
# insert here code for intersect_semiLine_circle (annex 1 section 2)
# insert here code for solveP2                   (annex 2 section 1}

########################################################################
# randomly generates 2 roots
# in the complex plane region {x+iy : -1<=x<=1, -1<=y<=1}
# by using the uniform distribution
R1 = runif(1,-1,1) + runif(1,-1,1)*1i
R2 = runif(1,-1,1) + runif(1,-1,1)*1i

# generates coefficients of p(z)= z^2 + C1*z + C2 by means of Vieta
C1 <- -(R1 + R2)
C2 <- R1*R2

#########################################################################
# approximates roots of polynomial by means of the basic LC method
OUTPUT <- solveP2(C1,C2)
# recovers items from output
Ra1 <- OUTPUT[1]  # approximation to one of the roots
Ra2 <- OUTPUT[2]  # approximation to one of the roots
V1  <- OUTPUT[3]  # direction vector of L1(Theta*)

#########################################################################
# generates graph where geometric construction LzC(Theta*) associated
# with polynomial z^2 + C1*z + C2 is shown, together with the roots and
# some other fixed points
PrintLabels = TRUE               # print algebraic expressions?
t = seq(from=-50,to=50,by=0.001) # parameter to construct trajectories
P1 = -C1/2                       # fixed point on L1
# trajectories: line L1, Derived line L1D and z-circumference ZC
L1   <- P1 + t*V1
L1D  <- P1^2 - (t^2)*(V1^2)
ZC   <- C2/L1
# draws a region in the complex plane
par(mar=c(2,2,2,2))
plot(0,0,'n',xlim=c(min(Re(ZC),Re(P1^2),Re(C2)),max(Re(ZC),Re(P1^2),Re(C2))),
             ylim=c(min(Im(ZC),Im(P1^2),Im(C2)),max(Im(ZC),Im(P1^2),Im(C2))),
             asp=1)
abline(h=seq(from=-50,to=50,by=0.2),col="light grey")
abline(v=seq(from=-50,to=50,by=0.2),col="light grey")
abline(h=0)
abline(v=0)
```

```
# line L1 (blue)
lines(L1,col="blue")
points(P1,pch=20,col="blue")

# line derived from L1 (red)
lines(L1D,col="red")
points(P1^2,pch=20,col="green")
points(C2,pch=20,col="blue")

# Z-circumference (blue)
lines(ZC,col="blue")
points(C2/P1,pch=20,col="black")

# approximate roots (red dots)
points(Ra1,pch=20,col="red")
points(Ra2,pch=20,col="red")

# line segments (black)
segments(0,0,Re(P1),Im(P1),col="black")
segments(Re(P1),Im(P1),Re(C2/P1),Im(C2/P1),col="black")

if (PrintLabels) {
   text(P1,pos=4,expression(P[1]))
   text(P1^2,pos=4,expression(P[1]^2))
   text(C2,pos=4,expression(C[2]))
   text(C2/P1,pos=4,expression(C[2]/P[1]))
   text(0+0i,pos=4,"0")
   text(Ra1,pos=4,expression(z[1]))
   text(Ra2,pos=4,expression(z[2]))
}

################################################################
# true roots
R1
R2

# coefficients of polynomial
C1
C2

# estimates of roots
Ra1
Ra2

# estimate of L1 direction vector
V1
```

**Remarks on this Script**

This script randomly generates points $R_1$, $R_2$ inside the squared region in $\mathbb{C}$

$\{x + iy : -1 \leq x \leq 1, -1 \leq y \leq 1\}$;

then, it uses Vieta's relations to generate coefficients $C_1, C_2$ from the generated points $R_1$, $R_2$. The script then calls function `solveP2` (listed in annex 2 section 1) using generated coefficients $C_1, C_2$ as inputs, and takes the outputs of `solveP2` to produce a graph of the structure $LzC(C_1, C_2, \theta^*)$ associated with polynomial $z^2 + C_1 z + C_2$. At the end, this script shows, in the R console, the generated points $R_1$, $R_2$, the generated coefficients $C_1, C_2$, and the approximations `Ra1`, `Ra2`, `V1` generated by function `solveP2`.

To reproduce the results from example 2.1 of chapter 2, it is enough to replace the instructions

```
R1 = runif(1,-1,1) + runif(1,-1,1)*1i
R2 = runif(1,-1,1) + runif(1,-1,1)*1i
```

with the instructions

```
R1 = -1+1i
R2 =  0-2i
```

## 2. R Script to Display Results of the Examples in Chapter 3

### Source Code

```
rm(list=ls(all=TRUE))
#set.seed(2020)  # reproduce results from example 3.1
#set.seed(2)     # reproduce results from example 3.3

########################################################################
# insert here code for center_zcircle                 (annex 1 section 1)
# insert here code for intersect_semiLine_circle      (annex 1 section 2)
# insert here code for error_tl_zc                    (annex 1 section 3)
# insert here code for get_angle_roots                (annex 1 section 4)
# insert here code for LzC3                           (annex 2 section 2)
# insert here code for approxP3                       (annex 2 section 5)

##########################################################################
# randomly generates 3 roots in the complex plane region
# {x+iy : -1<=x<=1, -1<=y<=1}, by using the uniform distribution
R1 = runif(1,-1,1) + runif(1,-1,1)*1i
R2 = runif(1,-1,1) + runif(1,-1,1)*1i
R3 = runif(1,-1,1) + runif(1,-1,1)*1i
R = c(R1,R2,R3)

##########################################################################
# constructs coefficients of polynomial z^3 + C1*z^2 + C2*z + C3,
# whose roots are R1, R2, R3, through Vieta's relations
C1 = -(R1 + R2 + R3)
C2 = R1*R2 + R1*R3 + R2*R3
C3 = -R1*R2*R3
c(C1,C2,C3)   # shows coefficients

# true theta root values (to compare them with approximations)
ThetaR = atan2(Im(R+C1/2),Re(R+C1/2))

####### approximates roots of polynomial by means of LC method ########
OUTPUT <- approxP3( C1, C2, C3, N=1000, -pi, pi )
# recovers items from output
E_A        <- OUTPUT[[1]]
E_B        <- OUTPUT[[2]]
Theta      <- OUTPUT[[3]]
ROOTS_LC_A <- OUTPUT[[4]]

# graphs functions of weighted errors.
# Adds points indicating the true theta roots,
# to compare them with graphical approximations
dev.new()
plot(Theta, E_A, 'l', xlim=c(-pi,pi), ylim=c(-1,1), col='blue',
     lty=1, xlab=expression(theta), ylab=expression(e(theta)))
lines(Theta, E_B, 'l', col='green3', lty=3)
abline(h=0, col='red')
points(ThetaR, rep(0,length(ThetaR)), col='red')
legend(x="bottomright",
       legend=c(expression(e[A](theta)),expression(e[B](theta))),
       col=c("blue","green3"),
       lty=c(1,3)
)
```

```
# true roots
R

# true values of theta roots
ThetaR

# table of approximations from LC method
ROOTS_LC_A
```

**Remarks on this Script**

In this script we first randomly generate three points $R_1$, $R_2$, $R_3$ inside square region in $\mathbb{C}$

$\{x + iy : -1 \leq x \leq 1, -1 \leq y \leq 1\}$;

then, we use Vieta's relations to generate coefficients $C_1, C_2, C_3$ from generated points $R_1$, $R_2$, $R_3$. Next, the script calls integrative function `approxP3` (listed in annex 2 section 5) using as inputs generated coefficients $C_1, C_2, C_3$, the number $N$ of discrete elements to generate the proximity map, and the starting and ending limits for the region where the proximity map is to be generated. Objects generated by `approxP3` are retrieved and shown in the R console: proximity map $\{\hat{e}_A(\theta), \hat{e}_B(\theta)\}$ is shown through a graph; generated points $R_1$, $R_2$, $R_3$, generated coefficients $C_1, C_2$, $C_3$, and true theta root values $\theta_1^*$, $\theta_2^*$, $\theta_3^*$ are listed; finally, the script shows a table with the initial approximations to the roots of polynomial $z^3 + C_1 z^2 + C_2 z + C_3$, obtained by means of function `approxP3`.

In order to reproduce the results of example 3.2 for polynomial

$z^3 + (1 + i)z^2 + (2 + 2i)z + (3 + 3i)$,

just replace the instructions

```
R1 = runif(1,-1,1) + runif(1,-1,1)*1i
R2 = runif(1,-1,1) + runif(1,-1,1)*1i
R3 = runif(1,-1,1) + runif(1,-1,1)*1i
```

with instructions

```
Rc2 = polyroot(c(3+3i,2+2i,1+1i,1))
R1 = Rc2[1]
R2 = Rc2[2]
R3 = Rc2[3]
```

the values returned by `polyroot` are used in chapter 3 as reference values for the measurement of errors associated with the generated approximations.

For the examples 3.4 and 3.5, just replace the instructions

```
R1 = runif(1,-1,1) + runif(1,-1,1)*1i
R2 = runif(1,-1,1) + runif(1,-1,1)*1i
R3 = runif(1,-1,1) + runif(1,-1,1)*1i
```

with instructions (example 3.4)

```
R1 = 2 + 5i
R2 = 2 + 5i
R3 = 1 - 2i
```

or with instructions (example 3.5)

```
R1 = 2 + 3i
R2 = 2 + 3i
R3 = 2 + 3i
```

## 3. R Script to Display Results of the Examples of Chapters 4 and 5

### Source Code

```
rm(list=ls(all=TRUE))
# set.seed(1987) # change value according to the example to be reproduced

##########################################################################
# insert here code for center_zcircle                 (annex 1 section 1)
# insert here code for intersect_semiLine_circle      (annex 1 section 2)
# insert here code for error_tl_zc                    (annex 1 section 3)
# insert here code for get_angle_roots                (annex 1 section 4)
# insert here code for min_D2                         (annex 2 section 4)
# insert here code for LzC                            (annex 2 section 3)
# insert here code for approxLC                       (annex 2 section 6)
# insert here code for approxDMin                     (annex 2 section 7)

########################################################################
# INSERT HERE CODE TO GENERATE REFERENCE ROOTS
R = # array with generated roots

###########################
# constructs polynomial coefficients from roots R
# INSERT HERE CODE TO GENERATE COEFFICIENTS
n = length(CC)    # degree of the polynomial
CC                # shows coefficients

# true values of theta roots
# (to compare them with approximations)
ThetaR = atan2(Im(R+CC[1]/2),Re(R+CC[1]/2))

############################################################################
# builds discrete proximity maps
N = 2500    # number of elements LzC to build proximity maps
# LC proximity map by means of intersections between tL and zC
system.time({
   OUT <- approxLC( CC, N, -pi, pi, MAXIT=500, TEMP=10, TMAX=200 )
   OUT[[1]] -> MinF
   OUT[[2]] -> MinT
   OUT[[3]] -> E_A
   OUT[[4]] -> E_B
   OUT[[5]] -> Theta
   OUT[[6]] -> APPROX_LC
})

# proximity map by means of discrete derivative D_d^2
system.time({
   OUT <- approxDMin(CC,MinF,-pi,pi,MAXIT=2000,TEMP=12,TMAX=200,Tol=2.0)
   OUT[[1]] -> DD_MinF
   OUT[[2]] -> Theta_shift_F
   OUT[[3]] -> APPROX_DMINF
})

# proximity map by means of discrete derivative D_t*
system.time({
   OUT <- approxDMin(CC,MinT,-pi,pi,MAXIT=2000,TEMP=12,TMAX=200,Tol=2.0)
   OUT[[1]] -> DD_MinT
   OUT[[2]] -> Theta_shift_T
   OUT[[3]] -> APPROX_DMINT
})
```

```
########################################################################
# graphical results associated with LC map
OUT <- LzC(CC, Theta, MinT)
OUT[[3]] -> Rx    # best root estimate for each structure LzC
OUT[[4]] -> aP    # anchor points of terminal semi-lines in elements LzC

# anchor points of true roots (only as a visual reference)
OUT <- LzC(CC, ThetaR, Re((R+CC[1]/2)/(cos(ThetaR)+sin(ThetaR)*1i)) )
OUT[[4]] -> aPR      # anchor points associated with true roots

# plots trajectory of roots defined by MinT, and trajectory of anchor
# points of terminal semi-lines
dev.new()
plot(Rx,type="l",asp=1, xlim=c(min(Re(Rx),Re(aP),na.rm=TRUE),
                               max(Re(Rx),Re(aP),na.rm=TRUE)),
                        ylim=c(min(Im(Rx),Im(aP),na.rm=TRUE),
                               max(Im(Rx),Im(aP),na.rm=TRUE)),
     xlab=expression(Re(R[x])), ylab=expression(Im(R[x])), mgp=c(2,1,0))
lines(aP, lty=3, col="blue") # trajectory of anchor points
points(R,col="red")          # true roots (for visual reference only)
points(aPR,col="red")        # anchor points associated with true roots
#segments(Re(aPR), Im(aPR), Re(R), Im(R), col="red")
legend(x="bottomleft",
       legend=c(expression(R[x]),expression(aP)),
       lty=c(1,3),
       col=c("black","blue")
)

# plots LC proximity map of weighted errors e(Theta). Adds points
# indicating true theta roots, to compare them with graphical
# approximations (crossings of e(Theta) with horizontal axis y=0)
dev.new()
plot(Theta, E_A, 'l', ylim=c(-1,1), col='blue', mgp=c(2,1,0),
     lty=1, xlab=expression(theta), ylab=expression(e(theta)))
lines(Theta, E_B, 'l', col='green3', lty=3)
abline(h=0, col='red')
points(ThetaR, rep(0,n), col='red')
legend(x="bottomright",
       legend=c(expression(e[A](theta)),expression(e[B](theta))),
       col=c("blue","green3"),
       lty=c(1,3)
)

########################################################################
# graphical results associated with maps D_d^2 and D_t*

# plots global minima d^2(t*(Theta)) vs Theta.
# Adds vertical line segments at the locations of true theta roots,
# to compare them with graphical approximations (zeros of d^2(t*(Theta)))
dev.new()
plot(Theta,MinF,type="l",mgp=c(2,1,0),
     xlab=expression(theta), ylab=expression(d[theta]^2~(t^"*"~(theta))))
segments(x0=ThetaR,y0=rep(0,n),x1=ThetaR,y1=rep(-1,n),col="red")
abline(h=0)

# plots minimizing arguments t* for functions d^2 vs Theta.
# Adds vertical lines at the locations of true theta roots,
# to compare them with graphical approximations (critical points of t*)
dev.new()
plot(Theta,MinT,type="l", mgp=c(2,1,0),
     xlab=expression(theta), ylab=expression(t^"*"~(theta)))
abline(v=ThetaR,col="red")
```

```
# plots discrete derivative D_d^2(t*(Theta)). Adds points at the
# locations of true theta roots, to compare them with graphical
# approximations (crossings of D_d^2(t*(Theta)) with axis y=0)
dev.new()
plot(Theta_shift_F, DD_MinF, type="l", ylim=c(-1,1),mgp=c(1.7,1,0),
        xlab=expression(theta),
        ylab=expression(over(d,d*theta)~d[theta]^2~(t^"*"~(theta))))
abline(h=0, col="red")
points(ThetaR, rep(0,n), col='red')

# plots discrete derivative D_t*(Theta). Adds points at the
# locations of true theta roots, to compare them with graphical
# approximations (crossings of D_t*(theta) with horizontal axis y=0)
dev.new()
plot(Theta_shift_T, DD_MinT, type="l", ylim=c(-1,1),mgp=c(1.7,1,0),
        xlab=expression(theta),
        ylab=expression(over(d,d*theta)~t^"*"~(theta)))
abline(h=0, col="red")
points(ThetaR, rep(0,n), col='red')

###############################################################
# tables of root approximations corresponding to proximity maps
# true roots
R
# true theta roots
ThetaR
# root approximations from map D_d^2(t*(Theta))
APPROX_DMINF

###############################################################
# true roots
R
# true theta roots
ThetaR
# root approximations from map D_t*(Theta)
APPROX_DMINT

###############################################################
# true roots
R
# true theta roots
ThetaR
# root approximations from map e(theta)
APPROX_LC
```

## Remarks on this Script

This script first incorporates the basic and central functions that are necessary for its operation; see annexes 1 and 2 for further details on these functions. Next, the script generates reference roots $R_1$, $R_2$, ..., $R_n$ ($n \geq 4$), and from these, it generates, through Vieta's relations, the coefficients $C_1$, $C_2$, ..., $C_n$ of polynomial $z^n + C_1 z^{n-1} + C_2 z^{n-2} + \cdots + C_{n-2} z^2 + C_{n-1} z + C_n$.

After this, reference theta roots $\theta_1^*$, $\theta_2^*$, ..., $\theta_n^*$ are generated; these are the inclination angles of lines $\ell_1$ with fixed point $-C_1/2$ that contain each of the respective roots $R_1$, $R_2$, ..., $R_n$. Reference values $R_1$, $R_2$, ..., $R_n$ and $\theta_1^*$, $\theta_2^*$, ..., $\theta_n^*$ are used to measure, within the numerical examples in chapters 4 and 5, the degree of accuracy of the approximations generated by this script.

This script generates the following three proximity maps:

1) An LC map given by the weighted error function $e(\theta)$; it consists of discrete arrays `E_A`, `E_B`, with support `Theta`. This LC map is generated by function `approxLC` (annex 2 section 6), which also produces discrete outputs $d_\theta^2(t^*(\theta))$ (array `MinF`) and $t^*(\theta)$ (array `MinT`), which are fundamental in the construction of the other two proximity maps generated by this script.

2) A map given by $\frac{d}{d\theta} d_\theta^2(t^*(\theta))$; it consists of array `DD_MinF`, with support `Theta_shift_F`. This map is generated by a first call to function `approxDMin` (annex 2 section 7), which takes as input the array `MinF` generated by `approxLC` at 1).

3) A map given by $\frac{d}{d\theta} t^*(\theta)$; it consists of array `DD_MinT`, with support `Theta_shift_T`. This map is generated by a second call to function `approxDMin`, which takes as input the array `MinT` generated by `approxLC` at 1).

The final part of the script generates plots for these three proximity maps, for functions $d_\theta^2(t^*(\theta))$ and $t^*(\theta)$, and for trajectories $aP(\theta)$ and $R_x(\theta)$ obtained by function `LzC` (annex 2 section 3); in addition, the script concludes with instructions for displaying (at the R console) numerical tables of initial approximations to roots, constructed by functions `approxLC` and `approxDMin`, and corresponding to the three generated proximity maps. In order to obtain plots similar to those shown in the numerical examples of chapters 4 and 5, it will be necessary to modify some graphical parameters in the final part of the script.

## Specific Conditions for Reproducing Results from the Examples in Chapter 4

In order to generate roots $R_1$, $R_2$, $R_3$, $R_4$ and coefficients $C_1$, $C_2$, $C_3$, $C_4$ in examples 4.1, 4.2 and 4.4, put at the beginning of this script one of the following instructions (initialization of pseudorandom number generator seed):

```
set.seed(1987)        # for example 4.1
set.seed(2020)        # for example 4.2
set.seed(2038)        # for example 4.4
```

Next, right below of the comment line

```
# INSERT HERE CODE TO GENERATE REFERENCE ROOTS
```

put the following instructions:

```
R1 = runif(1,-1,1) + runif(1,-1,1)*1i
R2 = runif(1,-1,1) + runif(1,-1,1)*1i
R3 = runif(1,-1,1) + runif(1,-1,1)*1i
R4 = runif(1,-1,1) + runif(1,-1,1)*1i
```

Also, replace the line

```
R = # array with generated roots
```

with

```
R = c(R1,R2,R3,R4)
```

Finally, right below of the comment line

```
# INSERT HERE CODE TO GENERATE COEFFICIENTS
```

put the instructions

```
C1 = -(R1 + R2 + R3 + R4)
C2 =   R1*R2 + R1*R3 + R1*R4 + R2*R3 + R2*R4 + R3*R4
C3 = -(R1*R2*R3 + R1*R2*R4 + R1*R3*R4 + R2*R3*R4)
C4 =   R1*R2*R3*R4
CC = c(C1,C2,C3,C4)
```

The number of elements $LzC$ used in the construction of proximity maps is specified by means of the variable `N`; see an example in chapter 4 to place in `N` the value specified in that example. To achieve close-ups of maps to specific regions, adjust values `Th_i`, `Th_f` in the calls to functions `approxLC` and `approxDMin`, as specified in the examples; it might also be necessary to adjust parameters `xlim`, `ylim` in the calls to function `plot`.

In Example 4.3, the random number generator seed is initialized by means of the instruction `set.seed(1987)`; the roots are generated by putting right below of the comment line

```
# INSERT HERE CODE TO GENERATE REFERENCE ROOTS
```

the instructions:

```
Rc1 = polyroot( c(4+4i, 3+3i, 2+2i, 1+1i, 1) )
R1 = Rc1[1]
R2 = Rc1[2]
R3 = Rc1[3]
R4 = Rc1[4]
```

Also, replace the line

```
R = # array with generated roots
```

with

```
R = c(R1,R2,R3,R4)
```

Coefficients `C1`, `C2`, `C3`, `C4` in example 4.3 are generated in the same way as in examples 4.1, 4.2 and 4.4; i.e., they are constructed from the initially generated roots.

NOTE: All **global** maps generated in the numerical examples from chapter 4 use the following lists of arguments in the calls to functions `approxLC` and `approxDMin`:

```
# LC proximity map by means of intersections between tL and zC
OUT <- approxLC( CC, N, -pi, pi, MAXIT=500, TEMP=10, TMAX=200 )
⋮

# proximity map by means of discrete derivative D_d^2
OUT <- approxDMin(CC,MinF,-pi,pi,MAXIT=2000,TEMP=12,TMAX=200,Tol=2.0)
⋮

# proximity map by means of discrete derivative D_t*
OUT <- approxDMin(CC,MinT,-pi,pi,MAXIT=2000,TEMP=12,TMAX=200,Tol=2.0)
⋮
```

To reproduce regional maps, only replace values `-pi`, `pi` with the values `Th_i`, `Th_f` specified in the examples, without changing any other argument in these calls to `approxLC` and `approxDMin`.

## Specific Conditions for Reproducing Results from the Examples in Chapter 5

**In example 5.1**, the random number generator seed is initialized by the instruction `set.seed(2020)` at the beginning of the script. We use `N=2500` elements $LzC$ for the construction of global maps. To generate reference roots, put right below of the comment line

```
# INSERT HERE CODE TO GENERATE REFERENCE ROOTS
```

the following instructions:

```
R1 = runif(1,-1,1) + runif(1,-1,1)*1i
R2 = runif(1,-1,1) + runif(1,-1,1)*1i
R3 = runif(1,-1,1) + runif(1,-1,1)*1i
R4 = runif(1,-1,1) + runif(1,-1,1)*1i
R5 = runif(1,-1,1) + runif(1,-1,1)*1i
R6 = runif(1,-1,1) + runif(1,-1,1)*1i
R7 = runif(1,-1,1) + runif(1,-1,1)*1i
```

Also, replace the line

```
R = # array with generated roots
```

with

```
R = c(R1,R2,R3,R4,R5,R6,R7)
```

Finally, right below of the comment line

```
# INSERT HERE CODE TO GENERATE COEFFICIENTS
```

put the instructions

```
C1 = -sum(apply(combn(R,1),2,prod))
C2 =  sum(apply(combn(R,2),2,prod))
C3 = -sum(apply(combn(R,3),2,prod))
C4 =  sum(apply(combn(R,4),2,prod))
C5 = -sum(apply(combn(R,5),2,prod))
C6 =  sum(apply(combn(R,6),2,prod))
C7 = -sum(apply(combn(R,7),2,prod))

CC = c(C1,C2,C3,C4,C5,C6,C7)
```

In this example, proximity maps are generated by using the same operation arguments used in the numerical examples from chapter 4; that is, calls to functions `approxLC` and `approxDMin` use the following lists of arguments:

```
# LC proximity map by means of intersections between tL and zC
OUT <- approxLC( CC, N, -pi, pi, MAXIT=500, TEMP=10, TMAX=200 )
⋮

# proximity map by means of discrete derivative D_d^2
OUT <- approxDMin(CC,MinF,-pi,pi,MAXIT=2000,TEMP=12,TMAX=200,Tol=2.0)
⋮

# proximity map by means of discrete derivative D_t*
OUT <- approxDMin(CC,MinT,-pi,pi,MAXIT=2000,TEMP=12,TMAX=200,Tol=2.0)
⋮
```

**In example 5.2**, the random number generator seed is initialized by the instruction `set.seed(2022)` at the beginning of the script. We use `N=5000` elements $LzC$ for the construction of global maps. To generate reference roots, just replace the line

```
R = # array with generated roots
```

with:

```
R=polyroot(c(10+10i,9+9i,8+8i,7+7i,6+6i,5+5i,4+4i,3+3i,2+2i,1+1i,1))
```

Next, right below of the comment line

```
# INSERT HERE CODE TO GENERATE COEFFICIENTS
```

put the instructions

```
C1 = -sum(apply(combn(R,1),2,prod))
C2 =  sum(apply(combn(R,2),2,prod))
C3 = -sum(apply(combn(R,3),2,prod))
C4 =  sum(apply(combn(R,4),2,prod))
C5 = -sum(apply(combn(R,5),2,prod))
C6 =  sum(apply(combn(R,6),2,prod))
C7 = -sum(apply(combn(R,7),2,prod))
C8 =  sum(apply(combn(R,8),2,prod))
C9 = -sum(apply(combn(R,9),2,prod))
C10=  sum(apply(combn(R,10),2,prod))

CC = c(C1,C2,C3,C4,C5,C6,C7,C8,C9,C10) # coefficient vector
```

Finally, in the calls to functions `approxLC` and `approxDMin`, use the following lines:

```
# LC proximity map by means of intersections between tL and zC
OUT <- approxLC( CC, N, -pi, pi, MAXIT=500, TEMP=500, TMAX=200 )
⋮

# proximity map by means of discrete derivative D_d^2
OUT <- approxDMin(CC,MinF,-pi,pi,MAXIT=2000,TEMP=500,TMAX=200,Tol=1000.0)
⋮

# proximity map by means of discrete derivative D_t*
OUT <- approxDMin(CC,MinT,-pi,pi,MAXIT=2000,TEMP=500,TMAX=200,Tol=1000.0)
⋮
```

Note that what changes in the calls to functions `approxLC` and `approxDMin`, with respect to the previous examples, are arguments `TEMP`; in the call to `approxLC` we use `500` instead of `10`, and in the calls to `approxDMin` we use `500` instead of `12`. In the two calls to function `approxDMin` (listed in annex 2 section 7), the argument `Tol` changes from `2.0` to `1000.0`, in order to raise the threshold for the detection of discrete "smooth" crossings of the proximity maps $\frac{d}{d\theta} d_\theta^2(t^*(\theta))$ and $\frac{d}{d\theta} t^*(\theta)$ with horizontal axis $y = 0$.

**<u>In example 5.3 part 1</u>**, we use the same operating parameters of example 5.2 (same instruction `set.seed(2022)` at the beginning of the script, same number of elements $LzC$ for the construction of global maps, namely, `N=5000`, and the same arguments in the calls to functions `approxLC` and `approxDMin`); what changes are the instructions for generating reference roots. Put right below of the comment line

```
# INSERT HERE CODE TO GENERATE REFERENCE ROOTS
```

the following instructions:

```
R1 = 1+0i
R2 = 2+0i
R3 = 3+0i
R4 = 4+0i
R5 = 5+0i
```

Also, replace the line

```
R = # array with generated roots
```

with

```
R = c(R1,R2,R3,R4,R5)
```

Finally, right below of the comment line

```
# INSERT HERE CODE TO GENERATE COEFFICIENTS
```

put the instructions

```
C1 = -sum(apply(combn(R,1),2,prod))
C2 =  sum(apply(combn(R,2),2,prod))
C3 = -sum(apply(combn(R,3),2,prod))
C4 =  sum(apply(combn(R,4),2,prod))
C5 = -sum(apply(combn(R,5),2,prod))

CC <- c(C1,C2,C3,C4,C5) # coefficient vector
```

**In example 5.3 part 2**, we recalculate the input coefficients according to the result from annex 5 section 2. For this, we insert, right after the instructions

```
n = length(CC)     # degree of the polynomial
CC                 # shows coefficients
```

the following lines:

```
# recalculates coefficients from the coefficients used
# in part 1, with a displacement a = 0 – 2i,
# and using the result from annex 5 section 2
i = 0:n
a = -2i             # displacement of roots
C = c(1,CC)
D1 = sum((((-1)^i)*choose(n-1+i,i)*(a^i))[i<=1]*C[1-i[i<=1]+1])
D2 = sum((((-1)^i)*choose(n-2+i,i)*(a^i))[i<=2]*C[2-i[i<=2]+1])
D3 = sum((((-1)^i)*choose(n-3+i,i)*(a^i))[i<=3]*C[3-i[i<=3]+1])
D4 = sum((((-1)^i)*choose(n-4+i,i)*(a^i))[i<=4]*C[4-i[i<=4]+1])
D5 = sum((((-1)^i)*choose(n-5+i,i)*(a^i))[i<=5]*C[5-i[i<=5]+1])
CC <- c(D1,D2,D3,D4,D5)  # new coefficient vector
```

Then, we apply to the reference roots of part 1 the displacement $a = 0 - 2i$, by means of the instruction

```
R = R+a
```

The random number generator seed is initialized by placing the instruction `set.seed(2022)` at the beginning of the script. We use `N=2500` elements $LzC$ in the construction of proximity maps.

For the calls to functions `approxLC` and `approxDMin`, we use the following lines:

```
# LC proximity map by means of intersections between tL and zC
OUT <- approxLC( CC, N, 0, pi, MAXIT=500, TEMP=10, TMAX=200 )
⋮

# proximity map by means of discrete derivative D_d^2
OUT <- approxDMin(CC,MinF,0,pi,MAXIT=2000,TEMP=12,TMAX=200,Tol=100)
⋮

# proximity map by means of discrete derivative D_t*
OUT <- approxDMin(CC,MinT,0,pi,MAXIT=2000,TEMP=12,TMAX=200,Tol=100)
⋮
```

**In example 5.4**, the random number generator seed is initialized by the instruction `set.seed(2022)` at the beginning of the script. To generate reference roots, just replace the line

```
R = # array with generated roots
```

with:

```
R=polyroot(c(15,14,13,12,11,10,9,8,7,6,5,4,3,2,1,1))
```

Next, right below of the comment line

```
# INSERT HERE CODE TO GENERATE COEFFICIENTS
```

put the instructions

```
C1 = -sum(apply(combn(R,1),2,prod))
C2 =  sum(apply(combn(R,2),2,prod))
C3 = -sum(apply(combn(R,3),2,prod))
C4 =  sum(apply(combn(R,4),2,prod))
C5 = -sum(apply(combn(R,5),2,prod))
C6 =  sum(apply(combn(R,6),2,prod))
C7 = -sum(apply(combn(R,7),2,prod))
C8 =  sum(apply(combn(R,8),2,prod))
C9 = -sum(apply(combn(R,9),2,prod))
C10=  sum(apply(combn(R,10),2,prod))
C11= -sum(apply(combn(R,11),2,prod))
C12=  sum(apply(combn(R,12),2,prod))
C13= -sum(apply(combn(R,13),2,prod))
C14=  sum(apply(combn(R,14),2,prod))
C15= -sum(apply(combn(R,15),2,prod))

# coefficient vector
CC <- c(C1,C2,C3,C4,C5,C6,C7,C8,C9,C10,C11,C12,C13,C14,C15)
```

Finally, in the calls to functions `approxLC` and `approxDMin`, use the following lines:

```
# LC proximity map by means of intersections between tL and zC
OUT <- approxLC( CC, N, -pi, pi, MAXIT=500, TEMP=800, TMAX=500 )
⋮

# proximity map by means of discrete derivative D_d^2
OUT <- approxDMin(CC,MinF,-pi,pi,MAXIT=2000,TEMP=500,TMAX=200,Tol=5000.0)
⋮

# proximity map by means of discrete derivative D_t*
OUT <- approxDMin(CC,MinT,-pi,pi,MAXIT=2000,TEMP=500,TMAX=200,Tol=10.0)
⋮
```

For this example, we use both `N=1000` and `N=50000` elements $LzC$ in the construction of proximity maps.

# Annex 4: Interactive R Scripts to Display Structures $LzC$

## 1. Interactive Script to Display Structures $LzC$ for Cubic Polynomials

**NOTES**:

1. To run this script, it is necessary to install the R package `tkrplot`, which allows to display graphs in R by means of Tcl/Tk widgets (see [A4.1]). For best results, run this script with the following display settings: display resolution = 1920 x 1080 pixels; scale = 100%.

2. When running this script, no text labels for points, lines, or circumferences are displayed on the interactive graph. The names of the graphic elements are specified in the comments in this script.

3. When running this script, the user can vary angle $\theta$ in the range from 0 to $\pi$ radians, which covers a complete turn of line $\ell_1$; this range differs conceptually from the interval where we usually define global angular proximity maps (from $-\pi$ to $\pi$ radians), because the direction vector $v_\theta$ is not shown here.

4. The function `intersect_l_zc` defined in this script is similar to function `intersect_semiLine_circle` described in annex 1 section 2, but the implementations of both functions differ from each other; `intersect_l_zc` is designed (like the rest of the functions defined in this script) to process a single value $\theta$ at a time, instead of taking advantage of vectorization of operations in R for “simultaneous” processing of multiple values $\theta$.

5. In order to observe the evolution of a structure $LzC(\theta)$ associated with a particular univariate polynomial of degree 3, it may be necessary to modify some graphical parameters within the script listed below, such as `hscale` and `vscale` at the call to function `tkrplot` near the end of this script; these parameters control the size of the window where the interactive graph is displayed. Other parameters that may need to be modified are `xlim` and `ylim`, within the call to function `plot`, within the definition of function `graphic`; these parameters control the minimum and maximum values of the real axis and the imaginary axis displayed on interactive graph.

### Source Code

```
rm(list=ls(all=TRUE))
library(tkrplot)
# insert here code for center_zcircle (annex 1 section 1)

################################################################
mb_line <- function(p1,p2) {
   # finds slope m and intercept b of a line
   # given two points p1, p2 on that line
   x1 <- Re(p1)
```

```
   y1 <- Im(p1)
   x2 <- Re(p2)
   y2 <- Im(p2)
   m  <- (y2-y1)/(x2-x1)
   b  <- y1 - m*x1
   return( c(m,b) )
}

sqrtm <- function(a) {
   # computes the two square roots of complex number a
   arg_a <- atan2( Im(a), Re(a) )
   mod_a <- abs( a )
   r1 <- sqrt(mod_a)*(cos(arg_a/2) + 1i*sin(arg_a/2))
   r2 <- sqrt(mod_a)*(cos(arg_a/2+pi) + 1i*sin(arg_a/2+pi))
   return( c(r1,r2) )
}

intersect_l_zc <- function( l1, l2, c1, c2 ) {
   # finds the intersections between a line
   # and a z-circumference.
   # l1, l2 are two points on the line.
   # c1 and c2 are two points on the z-circumference.
   czc <- center_zcircle(c1, c2)
   h   <- Re(czc)
   k   <- Im(czc)

   if ( Re(l1)!=Re(l2) ) {
      mbl <- mb_line(l1, l2)

      m   <- mbl[1]
      b   <- mbl[2]

      A <- (1+m^2)
      B <- 2*(m*b-h-k*m)
      C <- b*(b-2*k)
      D <- B^2 - 4*A*C

      if (D<0) {
         return(-1)
      } else {
         x1 <- (-B-sqrt(D))/(2*A)
         x2 <- (-B+sqrt(D))/(2*A)
         y1 <- m*x1 + b
         y2 <- m*x2 + b
         return( c( x1+y1*1i, x2+y2*1i) )
      }
   } else {
      A <- 1
      B <- -2*k
      C <- Re(l1)^2 - 2*h*Re(l1)
      D <- B^2 - 4*A*C

      if (D<0) {
         return(-1)
      }
      else {
         y1 <- (-B-sqrt(D))/(2*A)
         y2 <- (-B+sqrt(D))/(2*A)
         return( c(Re(l1)+y1*1i, Re(l1)+y2*1i) )
      }

   }
}
```

```
draw_z_circle <- function( center ) {
   # computes a series of points on a z-circumference
   # given the coordinates of its center; this series of points
   # is then used to draw the z-circumference,
   # by means of a procedure external to this function.
   r     <- abs(center)
   i     <- 0:3000
   theta <- -pi + i*(2*pi)/3000
   Xc    <- Re(center)+r*cos(theta)
   Yc    <- Im(center)+r*sin(theta)
   pzc   <- cbind( Xc, Yc )
   return(pzc)
}

#################################################################
# roots of a polynomial of degree 3
# (red points on the graph)
R1 <- runif(1,-1,1) + 1i*runif(1,-1,1)
R2 <- runif(1,-1,1) + 1i*runif(1,-1,1)
R3 <- runif(1,-1,1) + 1i*runif(1,-1,1)

# coefficients of a polynomial of the form x^3 + C1*x^2 + C2*x + C3
C1 <- -( R1 + R2 + R3 )
C2 <- R1*R2 + R1*R3 + R2*R3
C3 <- -R1*R2*R3

# fixed points associated with polynomial coefficients
# (black points on the graph)
P1 <- -C1/2
P2 <- -C3/P1
P3 <- C2 - P1^2

# double products of roots
# (green points on the graph)
Q1 <- R1*R2
Q2 <- R1*R3
Q3 <- R2*R3

######################################################################
# Generates interactive graph, in which the user can vary the
# inclination angle of line L1 through a slider-type control.
# We can observe L1, terminal semi-line C2-LD1, and z-circumference
# zC; additionally, we can appreciate the projections of the
# intersections between C2-LD1 and zC onto line L1 (blue and green
# points: projections associated with C2-LD1; blue and green
# diamonds: projections associated with zC; color blue is associated
# with negative radical of the quadratic formula used to compute the
# intersection; color green is associated with the positive radical;
# see function intersect_l_zc above in this script)

Sliderstart = 0.000
Slidermin   = 0.000
Slidermax   = pi
Sliderstep  = 0.001

tt <- tktoplevel()
SliderValue <- tclVar(Sliderstart)

graphic<-function(...) {

   # L1's inclination angle
   THETA1=as.numeric(tclvalue(SliderValue))
```

```
# L1's direction vector
V <-cos(THETA1) + 1i*sin(THETA1)

# a moving point on L1
PL <- P1 + V

# a moving point on C2-LD1 (terminal semi-line)
PLd <- C2 - (P1-10*V)*(P1+10*V)

# a moving point on z-circumference
PC <- -C3/PL

# slope and intercept of L1
MBL <- mb_line(P1,PL)
m <- MBL[1]
b <- MBL[2]

# gets coordinates of a discrete set
# of points in the z-circumference
ZC <- draw_z_circle( center_zcircle(P2,PC) )

# prepares graphics area
par(mar=c(1.9,1.9,1.9,1.9))
plot(0,0,'n',xlim=c(-1.1,1.1),ylim=c(-1.1,1.1), asp=1)
abline(h=0)
abline(v=0)

# draws line L1 (red color)
abline(b,m,col="red")

# draws z-circumference zC (green color)
lines(ZC, col='green3')

# draws terminal semi-line C2-LD1 (black color)
segments(Re(P3),Im(P3),Re(PLd),Im(PLd))

# draws polynomial roots (red points)
points(R1,pch=20,col="red")
points(R2,pch=20,col="red")
points(R3,pch=20,col="red")

# draws fixed points of L1, C2-LD1, and zC
# (black points)
points(P1,pch=20,col="black")
points(P2,pch=20,col="black")
points(P3,pch=20,col="black")

# draws double products of roots
# (green points)
points(Q1,pch=20,col="green3")
points(Q2,pch=20,col="green3")
points(Q3,pch=20,col="green3")

# computes intersections between terminal semi-line
# and z-circumference
I <- intersect_l_zc(P3,PLd,P2,PC)

# in case there are intersections
# between terminal semi-line and z-circumference
if (length(I)>1) {

   I1 <- I[1]     # negative radical
   I2 <- I[2]     # positive radical
```

```
      # is I1 on terminal semi-line?
      dir1 <- (I1-P3)/(PLd-P3)
      if (Re(dir1)>0) {
         r1 <- I1+P1^2-C2
         roots1 <- sqrtm( r1 )
         r11    <- roots1[1]
         r12    <- roots1[2]

         # Projections of I1 onto L1, via terminal semi-line
         pld11  <- P1 + r11
         pld12  <- P1 + r12

         # Projection of I1 onto L1, via z-circumference
         pzc1   <- -C3/I1

         # draws I1, and their respective projections onto L1
         # (blue points, and blue diamond for the projection
         # associated with z-circumference)
         points(I1,   pch=20,col="blue")
         points(pld11,pch=20,col="blue")
         points(pld12,pch=20,col="blue")
         points(pzc1, pch=5, col="blue")
      }

      # is I2 on terminal semi-line?
      dir2 <- (I2-P3)/(PLd-P3)
      if (Re(dir2)>0) {
         r2     <- I2+P1^2-C2
         roots2 <- sqrtm( r2 )
         r21    <- roots2[1]
         r22    <- roots2[2]

         # Projections of I2 onto L1, via terminal semi-line
         pld21  <- P1 + r21
         pld22  <- P1 + r22

         # Projection of I2 onto L1, vIa z-circumference
         pzc2   <- -C3/I2

         # draws I2, and their respective projections onto L1
         # (green points, and green diamond for the projection
         # associated with z-circumference)
         points(I2,   pch=20,col="green")
         points(pld21,pch=20,col="green")
         points(pld22,pch=20,col="green")
         points(pzc2, pch=5, col="green")
      }
   } # ends if (length(I)>1)
} # ends definition of function graphic

# builds graph of type tcl/tk
img <- tkrplot(tt, graphic, hscale=1.75, vscale=1.75)
showimage <- function(...) tkrreplot(img)

scl = tkscale(tt, from=Slidermin, to=Slidermax,
              showvalue=TRUE, variable=SliderValue, resolution=Sliderstep,
              command=showimage, orient='vertical')

tkpack(scl, side='left')
tkpack(img, side='right')
```

## 2. Interactive Script to Display Structures $\boldsymbol{LzC}$ for Polynomials of Degree $\boldsymbol{n \geq 4}$

**NOTES**:

1. To run this script, it is necessary to install the R package `tkrplot`, which allows to display graphs in R by means of Tcl/Tk widgets (see [A4.1]). For best results, run this script with the following display settings: display resolution = 1920 x 1080 pixels; scale = 100%.

2. When running this script, some text labels for points (roots and fixed points of lines, curves, and circumferences) are displayed on the interactive graph. The names of the graphic elements involved are specified in the comments in this script.

3. When running this script, the user can vary the inclination angle $\theta$ of line $\ell_1$ in the range from $-\pi$ to $\pi$ radians, which allows to observe the complete evolution in a cycle both of function $d_\theta^2(t)$ and of structure $LzC(\theta)$. Note that only half of line $\ell_1(\theta)$ is plotted, for parametric values $t > 0$.

4. In order to observe the evolution of a function $d_\theta^2(t)$ and of a structure $LzC(\theta)$, both associated with a particular univariate polynomial of degree $n \geq 4$, it may be necessary to modify some graphical parameters within the script listed below, such as `hscale` and `vscale` at the call to function `tkrplot` near the end of this script; these parameters control the size of the window where the interactive graph is displayed. Other parameters that may need to be modified are `xlim` and `ylim`, within the calls to function `plot`, within the definition of function `graphic`; these parameters control the minimum and maximum values of the horizontal and vertical axes shown in the interactive graph.

5. The script listed in this section is designed to visualize structures $LzC(\theta)$ and dynamic squared distance functions $d_\theta^2(t)$ associated with univariate polynomials of degree 7; if you want to change the degree `n` of the polynomial, then it is also necessary to change the instructions in the initial part of the script where the roots `R1`, `R2`, ..., `Rn`, the coefficients `C1`, `C2`, ..., `Cn` and the $n-1$ products `Q1`, `Q2`, ..., `Qn` are generated, as well as the part where the points `R1`, `R2`, ..., `Rn` and `Q1`, `Q2`, ..., `Qn` are plotted, within the definition of function `graphic`.

6. In this script, unlike the routines in annex 2 used in the construction of proximity maps associated with univariate polynomials of degree $n \geq 4$, negative values of $t^*(\theta)$ are considered, which allows us to observe the behavior of the elements within structures $LzC(\theta)$ (such as the terminal semi-line $tL$) built with values $t^*(\theta)$ both positive and negative.

## Source Code

```
rm(list=ls(all=TRUE))
library(tkrplot)
set.seed(2020)

# insert here code for center_zcircle (annex 1 section 1)
# insert here code for draw_z_circle (see script in annex 4 section 1)

##################################################################
# generates 7 random roots
# in the region of the complex plane {x+iy : -1<=x<=1, -1<=y<=1}
# by using the uniform distribution
# (red points on the graph)
R1 = runif(1,-1,1) + runif(1,-1,1)*1i
R2 = runif(1,-1,1) + runif(1,-1,1)*1i
R3 = runif(1,-1,1) + runif(1,-1,1)*1i
R4 = runif(1,-1,1) + runif(1,-1,1)*1i
R5 = runif(1,-1,1) + runif(1,-1,1)*1i
R6 = runif(1,-1,1) + runif(1,-1,1)*1i
R7 = runif(1,-1,1) + runif(1,-1,1)*1i
R = c(R1,R2,R3,R4,R5,R6,R7)

###########################
# constructs coefficients of polynomial
# z^7 + C1*z^6 + C2*z^5 + C3*z^4 + C4*z^3 + C5*z^2 + C6*z + C7,
# whose roots are R1, R2, R3, R4, R5, R6, R7

C1 = -sum(apply(combn(R,1),2,prod))
C2 =  sum(apply(combn(R,2),2,prod))
C3 = -sum(apply(combn(R,3),2,prod))
C4 =  sum(apply(combn(R,4),2,prod))
C5 = -sum(apply(combn(R,5),2,prod))
C6 =  sum(apply(combn(R,6),2,prod))
C7 = -sum(apply(combn(R,7),2,prod))

# sextuple products of roots
# (green fixed points on the graph)
Q1 <- R1*R2*R3*R4*R5*R6
Q2 <- R1*R2*R3*R4*R5*R7
Q3 <- R1*R2*R3*R4*R6*R7
Q4 <- R1*R2*R3*R5*R6*R7
Q5 <- R1*R2*R4*R5*R6*R7
Q6 <- R1*R3*R4*R5*R6*R7
Q7 <- R2*R3*R4*R5*R6*R7

CC <- c(C1,C2,C3,C4,C5,C6,C7)

n=length(CC) # degree of the polynomial

# fixed points associated with polynomial coefficients
# (blue fixed points on the graph)
P1 = -C1/2                                             # on L1
PT = ((-1)^(n+1))*(CC[(n-1):2]%*%P1^(0:(n-3))-P1^(n-1)) # on tC
PC = ((-1)^n)*CC[n]/P1                                 # on zC

########################################
N = 2500       # number of equidistant discrete values of Theta
               # in interval [-pi, pi)
```

```
# stores global minimum value of dynamic squared distance (DSD)
# as a function of Theta
MinF = rep(NA,N)

# stores the value of parameter t where DSD's global minimum
# occurs as a function of Theta
MinT = rep(NA,N)

i = 1:N                         # Index

Theta = -pi + 2*pi*(i-1)/N      # equidistant discrete values of
                                # Theta in interval [-pi, pi)

V = cos(Theta) + sin(Theta)*1i  # direction vectors of lines L1

LF = vector("list",N)  # stores functions DSD
                       # a function DSD for each value of Theta
LF = apply(as.data.frame(i),1,function(i){
   v = V[i]
   function(x) {
      z = P1+x*v
      tC = ((-1)^(n+1))*(CC[(n-1):2]%*%z^(0:(n-3))-(P1-x*v)*z^(n-2))
      zC = ((-1)^n)*CC[n]/z
      return( Re(tC-zC)^2 + Im(tC-zC)^2 )
   }
})
# finds the global minimum for each function DSD:
# phase 1: optimization using simulated annealing
system.time({
t_min0 = apply(as.data.frame(i),1,function(i){
               optim(1,LF[[i]],method="SANN",
                     control=list(maxit=2000,temp=15,tmax=200))
})
})
# phase 2: refinement of the result from phase 1 via BFGS
system.time({
t_min = apply(as.data.frame(i),1,function(i){
               optim(t_min0[[i]]$par,LF[[i]],method="BFGS")
})
})

# retrieves DSD's minimal functional value
MinF[i] = unlist(apply(as.data.frame(i),1,function(i){
                 t_min[[i]]$value
}))

# retrieves parametric value t where DSD's global minimum occurs
MinT[i] = unlist(apply(as.data.frame(i),1,function(i){
                  t_min[[i]]$par
}))

######################################################################
# Generates interactive graph, in which the user can vary the
# inclination angle Theta of line L1 through a slider-type control.
# We can observe line L1, terminal curve tC, z-circumference zC, and
# terminal semi-line tL; all these elements as a function of Theta.
# We can also see function DSD(t) for each sample value of Theta.
tt <- tktoplevel()
sliderstart = 1
slidermin   = 1
slidermax   = N
sliderstep  = 1
chkb_actv   = 1
```

```
chkb_deactv = 0
SliderValue <- tclVar(sliderstart)
cb1Value    <- tclVar(chkb_actv)
cb2Value    <- tclVar(chkb_deactv)

graphic<-function(...) {

   # Index associated with slider
   ii=as.numeric(tclvalue(SliderValue))

   dsd <- function(x) {
      tC <- ((-1)^(n+1))*(CC[(n-1):2]%*%(P1+x*V[ii])^(0:(n-3)) -
                        (P1-x*V[ii])*(P1+x*V[ii])^(n-2))
      zC <- ((-1)^n)*CC[n]/(P1+x*V[ii])
      return(Re(tC-zC)^2 + Im(tC-zC)^2)
   }

   # plots function dsd(t) on the left of graph
   # plots structure LzC on the right of graph
   par( mfrow=c(1,2), mar=c(2,2,2,2) )

   # graphs function dsd(t) for a fixed value of Theta
   t=seq(from=-2,to=2,by=0.005)
   plot(t,sapply(t,dsd),type="l",ylim=c(0,0.1))
   points(MinT[ii], MinF[ii], col=ifelse(MinT[ii]>=0,"blue","red"))
   abline(h=0)
   abline(v=0)
   tm = format(MinT[ii], digits=5)
   dsdm = format(MinF[ii], digits=5)
   TH = format(Theta[ii], digits=5)
   text(1,0.08,bquote(theta == .( TH)))
   text(1,0.075,bquote(t^"*"~(theta) == .( tm)))
   text(1,0.07,bquote(d[theta]^2~(t^"*"~(theta)) == .( dsdm)))

   # computes points of line L1
   t = seq(from=0,to=10*max(abs(MinT)),by=0.005)
   L1 = P1+t*V[ii]

   # computes points of terminal curve tC
   ftc = function(x)((-1)^(n+1))*(CC[(n-1):2]%*%(P1+x*V[ii])^(0:(n-3)) -
                        (P1-x*V[ii])*(P1+x*V[ii])^(n-2))
   tC = sapply(t,ftc)

   # computes points of terminal semi-line tL
   Rx = P1 + MinT[ii]*V[ii]
   AP = ((-1)^(n+1))*(CC[(n-1):3]%*%Rx^(0:(n-4))+(CC[2]-P1^2)*Rx^(n-3))
   VTL= ((-1)^(n+1))*(Rx^(n-3))*V[ii]^2
   tL = rep(AP,length(t)) + (t^2)*VTL

   # computes points of z-circumference
   centreZc = center_zcircle(PC, ((-1)^n)*CC[n]/(P1+V[ii]) )
   zCircle = draw_z_circle(centreZc)

   # plots z-circumference
   plot(zCircle,type="l",xlim=c(-1.1,1.1),ylim=c(-1.1,1.1),
        col="green3", asp=1 )
   abline(h=0)
   abline(v=0)
   points(PC,pch=20,col="blue")
   text(0.05+PC,expression(P[C]),cex=1.0)
```

```
   # draws line L1
   lines(L1,col="red")
   points(P1,pch=20,col="blue")
   text(0.05+P1,expression(P[1]),cex=1.0)

   # draws terminal curve tC
   if (as.numeric(tclvalue(cb1Value))==chkb_actv) {
      lines(tC)
      points(PT,pch=20,col="blue")
      text(0.05+PT,expression(P[T]),cex=1.0)
   }

   # draws terminal semi-line tL
   if (as.numeric(tclvalue(cb2Value))==chkb_actv) {
      lines(tL,col="blue")
      points(AP,pch=20,col="blue")
      text(0.05+AP,expression(P[a]),cex=1.0)
   }

   # draws roots of polynomial (red points)
   points(R1,pch=20,col="red")
   points(R2,pch=20,col="red")
   points(R3,pch=20,col="red")
   points(R4,pch=20,col="red")
   points(R5,pch=20,col="red")
   points(R6,pch=20,col="red")
   points(R7,pch=20,col="red")
   text(0.05+R1,expression(R[1]),cex=1.0)
   text(0.05+R2,expression(R[2]),cex=1.0)
   text(0.05+R3,expression(R[3]),cex=1.0)
   text(0.05+R4,expression(R[4]),cex=1.0)
   text(0.05+R5,expression(R[5]),cex=1.0)
   text(0.05+R6,expression(R[6]),cex=1.0)
   text(0.05+R7,expression(R[7]),cex=1.0)

   # draws sextuple products of roots
   # (green points)
   points(Q1,pch=20,col="green3")
   points(Q2,pch=20,col="green3")
   points(Q3,pch=20,col="green3")
   points(Q4,pch=20,col="green3")
   points(Q5,pch=20,col="green3")
   points(Q6,pch=20,col="green3")
   points(Q7,pch=20,col="green3")

   # draws points of L1, zC and tC associated with value t where
   # DSD's global minimum occurs (only if t is non-negative)
   if (MinT[ii] >= 0) {
      points(((-1)^n)*CC[n]/(P1+MinT[ii]*V[ii]),pch=5)
      points(P1+MinT[ii]*V[ii])
      if (as.numeric(tclvalue(cb1Value))==chkb_actv) {
         tmin <- P1+MinT[ii]*V[ii]
         tmc  <- P1-MinT[ii]*V[ii]
         points(((-1)^(n+1))*(CC[(n-1):2]%*%tmin^(0:(n-3))-tmc*tmin^(n-2)))
      }
   }
} # ends definition of function graphic
```

```
# builds graph of type tcl/tk
img <- tkrplot(tt, graphic, hscale=5.00, vscale=2.50)
showimage <- function(...) tkrreplot(img)
scl = tkscale(tt, from=slidermin, to=slidermax,
              showvalue=TRUE, variable=SliderValue, resolution=sliderstep,
              command=showimage, orient='horizontal')
cb1 =ttkcheckbutton(tt, command=showimage, text='Terminal curve',
                    variable=cb1Value)
cb2 =ttkcheckbutton(tt, command=showimage, text='Terminal line',
                    variable=cb2Value)
tkpack(img, side='bottom')
tkpack(cb2, padx=50, side='right')
tkpack(cb1, padx=20, side='right')
tkpack(scl, padx=200, side='right')
```

# Annex 5: A Pair of Theoretical Results on Polynomials

## 1. How to Find the Roots of a Quartic Polynomial Through a Resolvent Cubic

Let $p(x) = x^4 + C_1x^3 + C_2x^2 + C_3x + C_4$ be a degree 4 polynomial with roots $R_1$, $R_2$, $R_3$, $R_4$.

$\Longrightarrow p(x) = 0 \Leftrightarrow x = R_1$ or $x = R_2$ or $x = R_3$ or $x = R_4$.

From Vieta's relations, we have

$C_1 = -(R_1 + R_2 + R_3 + R_4)$,

$C_2 = R_1R_2 + R_1R_3 + R_1R_4 + R_2R_3 + R_2R_4 + R_3R_4$,

$C_3 = -(R_1R_2R_3 + R_1R_2R_4 + R_1R_3R_4 + R_2R_3R_4)$,

$C_4 = R_1R_2R_3R_4$.

Now, let $q(x) = x^3 + D_1x^2 + D_2x + D_3$ be a degree 3 polynomial with roots

$S_1 = R_1R_2 + R_3R_4$, $S_2 = R_1R_3 + R_2R_4$, $S_3 = R_1R_4 + R_2R_3$.

$\Longrightarrow q(x) = 0 \Leftrightarrow x = S_1$ or $x = S_2$ or $x = S_3$.

Let us obtain the expressions of the coefficients of $q(x)$ in terms of the coefficients of $p(x)$. From Vieta's relations, we have

$D_1 = -(S_1 + S_2 + S_3)$

$D_1 = -(R_1R_2 + R_3R_4 + R_1R_3 + R_2R_4 + R_1R_4 + R_2R_3)$

$D_1 = -C_2$

$D_2 = S_1S_2 + S_1S_3 + S_2S_3$

$D_2 = (R_1R_2 + R_3R_4)(R_1R_3 + R_2R_4)$

$+(R_1R_2 + R_3R_4)(R_1R_4 + R_2R_3)$

$+(R_1R_3 + R_2R_4)(R_1R_4 + R_2R_3)$

$D_2 = R_1^2R_2R_3 + R_1R_2^2R_4 + R_1R_3^2R_4 + R_2R_3R_4^2$

$+R_1^2R_2R_4 + R_1R_2^2R_3 + R_1R_3R_4^2 + R_2R_3^2R_4$

$+R_1^2R_3R_4 + R_1R_2R_3^2 + R_1R_2R_4^2 + R_2^2R_3R_4$

$D_2 = R_1(R_1R_2R_3 + R_1R_2R_4 + R_1R_3R_4) + R_1(R_2R_3R_4 - R_2R_3R_4)$

$+R_2(R_1R_2R_4 + R_1R_2R_3 + R_2R_3R_4) + R_2(R_1R_3R_4 - R_1R_3R_4)$

$+R_3(R_1R_3R_4 + R_2R_3R_4 + R_1R_2R_3) + R_3(R_1R_2R_4 - R_1R_2R_4)$

$+R_4(R_2R_3R_4 + R_1R_3R_4 + R_1R_2R_4) + R_4(R_1R_2R_3 - R_1R_2R_3)$

$D_2 = R_1(-C_3) + R_2(-C_3) + R_3(-C_3) + R_4(-C_3) - 4C_4$

$D_2 = (-C_1)(-C_3) - 4C_4$

$D_2 = C_1C_3 - 4C_4$

$D_3 = -S_1S_2S_3$

$D_3 = -(R_1R_2 + R_3R_4)(R_1R_3 + R_2R_4)(R_1R_4 + R_2R_3)$

$D_3 = -(R_1^2R_2R_3 + R_1R_2^2R_4 + R_1R_3^2R_4 + R_2R_3R_4^2)(R_1R_4 + R_2R_3)$

$D_3 = -(R_1^3R_2R_3R_4 + R_1^2R_2^2R_4^2 + R_1^2R_3^2R_4^2 + R_1R_2R_3R_4^3$

$+R_1^2R_2^2R_3^2 + R_1R_2^3R_3R_4 + R_1R_2R_3^3R_4 + R_2^2R_3^2R_4^2)$

$D_3 = -[C_4(R_1^2 + R_2^2 + R_3^2 + R_4^2) + C_4(2C_2 - 2C_2)$

$+(R_1R_2R_3)^2 + (R_1R_2R_4)^2 + (R_1R_3R_4)^2 + (R_2R_3R_4)^2]$

$D_3 = -[C_4C_1^2 - 2C_4C_2 + C_3^2$

$-2R_1^2R_2^2R_3R_4 - 2R_1^2R_2R_3^2R_4 - 2R_1R_2^2R_3^2R_4 - 2R_1^2R_2R_3R_4^2 - 2R_1R_2^2R_3R_4^2 - 2R_1R_2R_3^2R_4^2]$

$D_3 = -[C_4C_1^2 - 2C_4C_2 + C_3^2 - 2C_4(R_1R_2 + R_1R_3 + R_2R_3 + R_1R_4 + R_2R_4 + R_3R_4)]$

$D_3 = -[C_4C_1^2 - 2C_4C_2 + C_3^2 - 2C_4C_2]$

$D_3 = 4C_2C_4 - C_1^2C_4 - C_3^2$

In summary, $q(x)$ can be expressed in terms of the coefficients of $p(x)$ as

$q(x) = x^3 + (-C_2)x^2 + (C_1C_3 - 4C_4)x + (4C_2C_4 - C_1^2C_4 - C_3^2)$.

$q(x)$ is called *resolvent cubic polynomial* of $p(x)$.

To obtain the roots of $p(x) = x^4 + C_1x^3 + C_2x^2 + C_3x + C_4$ by means of the resolvent cubic polynomial $q(x) = x^3 - C_2x^2 + (C_1C_3 - 4C_4)x + (4C_2C_4 - C_1^2C_4 - C_3^2)$:

1. Solve $q(x) = 0$, to obtain $S_1, S_2, S_3$.

2. Solve $x^2 + C_1x + (C_2 - S_1) = 0$, to obtain $T_{1,1} \coloneqq R_1 + R_2, T_{1,2} \coloneqq R_3 + R_4$.

3. Solve $x^2 - S_1x + C_4 = 0$, to obtain $T_{2,1} \coloneqq R_1R_2, T_{2,2} \coloneqq R_3R_4$.

4. In principle, $R_1$, $R_2$ is obtained by solving $x^2 - T_{1,1}x + T_{2,1} = 0$; and $R_3$, $R_4$ is obtained by solving $x^2 - T_{1,2}x + T_{2,2} = 0$.

In practice, to implement step 4, we must verify whether the roots $\rho_1$, $\rho_2$ of equation $x^2 - T_{1,1}x + T_{2,1} = 0$ satisfy equation $p(x) = 0$:

- If $p(\rho_1) \neq 0$, then the four roots of equation $p(x) = 0$ are obtained by solving the equations $x^2 - T_{1,1}x + T_{2,2} = 0$ and $x^2 - T_{1,2}x + T_{2,1} = 0$.

- On the other hand, if $p(\rho_1) = 0$, then we already have two of the roots of $p(x) = 0$, and the other two roots are obtained by solving equation $x^2 - T_{1,2}x + T_{2,2} = 0$.

## 2. Theorem: Change of Variable and Root Shifting

If in equation

$$x^n + C_1x^{n-1} + C_2x^{n-2} + \cdots + C_{n-2}x^2 + C_{n-1}x + C_n = 0 \qquad \text{(A5.1)}$$

with roots $R_1$, $R_2$,..., $R_n$,

we make the change of variable $z = x + a$, with $a \in \mathbb{C}$,

then equation

$$z^n + D_1z^{n-1} + D_2z^{n-2} + \cdots + D_{n-2}z^2 + D_{n-1}z + D_n = 0 \qquad \text{(A5.2)}$$

with coefficients $D_k = \sum_{i=0}^{k} \binom{n-k+i}{i} (-a)^i C_{k-i}; \qquad C_0 = 1, \quad k \in \{1,2,\dots,n\}$,

has roots $R_1 + a$, $R_2 + a$,..., $R_n + a$.

**Demonstration**

Since $x = z - a$, we can rewrite equation (A5.1) as

$$(z-a)^n + C_1(z-a)^{n-1} + C_2(z-a)^{n-2} + \cdots + C_{n-2}(z-a)^2 + C_{n-1}(z-a) + C_n = 0 \qquad \text{(A5.1a)}$$

Expanding the powers of binomial $(z - a)$ in equation (A5.1a), we obtain

$$\begin{array}{l}
z^n + \binom{n}{1} z^{n-1}(-a) + \binom{n}{2} z^{n-2}(-a)^2 + \cdots + \binom{n}{n-2} z^2(-a)^{n-2} + \binom{n}{n-1} z(-a)^{n-1} + (-a)^n \\
\quad + C_1z^{n-1} + C_1\binom{n-1}{1} z^{n-2}(-a) + \cdots + C_1\binom{n-1}{n-3} z^2(-a)^{n-3} + C_1\binom{n-1}{n-2} z(-a)^{n-2} + C_1(-a)^{n-1} \\
\quad\quad + C_2\, z^{n-2} + \cdots + C_2\binom{n-2}{n-4} z^2(-a)^{n-4} + C_2\binom{n-2}{n-3} z(-a)^{n-3} + C_2(-a)^{n-2} \\
\quad\quad\quad \vdots \qquad\qquad \vdots \qquad\qquad \vdots \\
\quad\quad\quad + C_{n-2}z^2 + C_{n-2}\binom{2}{1} z(-a) + C_{n-2}(-a)^2 \\
\quad\quad\quad\quad + C_{n-1}z + C_{n-1}(-a) \\
\quad\quad\quad\quad\quad + C_n = 0
\end{array}$$

Adding like terms of powers of $z$ in this last expression, and incorporating the artificial coefficient $C_0 = 1$, we obtain

$$z^n + \left[\sum_{i=0}^{1} \binom{n-1+i}{i} (-a)^i C_{1-i}\right] z^{n-1} + \left[\sum_{i=0}^{2} \binom{n-2+i}{i} (-a)^i C_{2-i}\right] z^{n-2} + \cdots$$

$$+ \left[\sum_{i=0}^{n-2} \binom{2+i}{i} (-a)^i C_{n-2-i}\right] z^2 + \left[\sum_{i=0}^{n-1} \binom{1+i}{i} (-a)^i C_{n-1-i}\right] z + \left[\sum_{i=0}^{n} \binom{i}{i} (-a)^i C_{n-i}\right] = 0$$

Or more concisely,

$$z^n + \sum_{k=1}^{n} \left[ \sum_{i=0}^{k} \binom{n-k+i}{i} (-a)^i C_{k-i} \right] z^{n-k} = 0$$

$$z^n + \sum_{k=1}^{n} D_k z^{n-k} = 0, \qquad \text{(A5.2a)}$$

with $D_k = \sum_{i=0}^{k} \binom{n-k+i}{i} (-a)^i C_{k-i}; \qquad C_0 = 1, \quad k \in \{1,2,\dots,n\}.$

From here, we see that expression (A5.2a) is identical to expression (A5.2); furthermore, equation (A5.2a) is equivalent to equation (A5.1a), since the steps taken to get from (A5.1a) to (A5.2a) are reversible.

Now, let $S_j = R_j + a, \qquad j \in \{1,2,\dots,n\};$

then, evaluating the left side of (A5.2a) at $z = S_j$, we obtain

$$S_j^n + \sum_{k=1}^{n} D_k S_j^{n-k} = (S_j - a)^n + \sum_{k=1}^{n} C_k (S_j - a)^{n-k} = (R_j)^n + \sum_{k=1}^{n} C_k (R_j)^{n-k}.$$

Equality $S_j^n + \sum_{k=1}^{n} D_k S_j^{n-k} = (S_j - a)^n + \sum_{k=1}^{n} C_k (S_j - a)^{n-k}$ is due to the equivalence between equations (A5.2a) and (A5.1a); expression $(R_j)^n + \sum_{k=1}^{n} C_k (R_j)^{n-k}$ is nothing more than the left side of equation (A5.1) evaluated at one of its roots, so $(R_j)^n + \sum_{k=1}^{n} C_k (R_j)^{n-k} = 0$, and finally, we have

$$S_j^n + \sum_{k=1}^{n} D_k S_j^{n-k} = 0,$$

which means that $S_j = R_j + a$ is a root of equation (A5.2). ■

# Index